\documentclass[12pt]{report}
\usepackage[utf8]{inputenc}

\usepackage[a4paper, tmargin=3cm, bmargin=2cm, rmargin=2cm, lmargin=3cm]{geometry}

\usepackage[main=english, portuguese]{babel}

\usepackage{times}

\usepackage{setspace}
\usepackage{multicol}

\usepackage{booktabs}
\usepackage{amssymb}
\usepackage{amsmath}

\usepackage{graphicx}

\usepackage{enumitem}

\usepackage{mwe}

\usepackage{colortbl}

\usepackage{caption}

\usepackage{mathrsfs}

\usepackage{marginnote}

\usepackage{bbding}

\usepackage{pdfpages}

\usepackage{float}

\usepackage{multirow}

\usepackage{mathtools}

\usepackage{stix}

\usepackage{longtable}

\usepackage[super]{nth}

\usepackage{lscape}

\usepackage{imakeidx}
\makeindex[columns=2, title=Alphabetical Index, intoc, options= -s example_style.ist]

\usepackage[dvipsnames]{xcolor}

\usepackage{hyperref}
\hypersetup{
    colorlinks=true,
    linkcolor=ForestGreen,
    citecolor=red,    
    urlcolor=blue,
    linktocpage=true
}

\usepackage{rotating}

\usepackage{pifont}

\usepackage{pgfplots}
\pgfplotsset{compat=newest}

\usepgfplotslibrary{polar}
\usepgflibrary{shapes.geometric}

\usepackage{shorttoc}
\makeatletter
\renewcommand*{\tableofcontents}{%
  \begingroup
    \let\@input\@gobble
    \@starttoc{toc}%
  \endgroup
}
\makeatother

\usepackage{array}
\newcolumntype{L}[1]{>{\raggedright\let\newline\\\arraybackslash\hspace{0pt}}m{#1}}
\newcolumntype{C}[1]{>{\centering\let\newline\\\arraybackslash\hspace{0pt}}m{#1}}
\newcolumntype{R}[1]{>{\raggedleft\let\newline\\\arraybackslash\hspace{0pt}}m{#1}}

\usepackage{tikz-cd}
\usetikzlibrary{babel}

\usepackage[export]{adjustbox}
\usetikzlibrary {chains}
\usetikzlibrary{arrows.meta}
\usetikzlibrary{positioning}
\usetikzlibrary{decorations}
\usetikzlibrary{calc}

\usepackage{fancyhdr} 
\fancypagestyle{plain}{
\fancyhf{}
 
\rhead{\thepage}}

\usepackage{makecell}
\newcolumntype{x}[1]{>{\centering\arraybackslash}p{#1}}

\usepackage{tikz}
\newcommand\diag[4]{%
  \multicolumn{1}{p{#2}|}{\hskip-\tabcolsep
  $\vcenter{\begin{tikzpicture}[baseline=0,anchor=south west,inner sep=#1]
  \path[use as bounding box] (0,0) rectangle (#2+2\tabcolsep,\baselineskip);
  \node[minimum width={#2+2\tabcolsep},minimum height=\baselineskip+\extrarowheight] (box) {};
  \draw (box.north west) -- (box.south east);
  \node[anchor=south west] at (box.south west) {#3};
  \node[anchor=north east] at (box.north east) {#4};
 \end{tikzpicture}}$\hskip-\tabcolsep}}
 
 \usetikzlibrary{decorations.markings}
\tikzset{negated/.style={
        decoration={markings,
            mark= at position 0.5 with {
                \node[transform shape] (tempnode) {$\backslash$};
            }
        },
        postaction={decorate}
    }
}

\usepackage{prooftrees,turnstile}

\usepackage{amsthm}

\newtheorem{theorem}{Theorem}
\numberwithin{theorem}{section}

\newtheorem{definition}{Definition}
\numberwithin{definition}{section}

\newtheorem{proposition}{Proposition}
\numberwithin{proposition}{section}

\newtheorem{lemma}{Lemma}
\numberwithin{lemma}{section}

\numberwithin{axiom}{section}

\newtheorem{corollary}{Corollary}
\numberwithin{corollary}{section}

\newtheorem{example}{Example}
\numberwithin{example}{section}

\newtheorem{problem}{Problem}
\numberwithin{problem}{section}

\newcommand{\bI}{\textbf{bI}}
\newcommand{\bIpr}{\textbf{bIpr}}
\newcommand{\nbI}{\textbf{nbI}}
\newcommand{\nbIciw}{\textbf{nbIciw}}
\newcommand{\nbIci}{\textbf{nbIci}}
\newcommand{\nbIcl}{\textbf{nbIcl}}
\newcommand{\CILA}{\textbf{Cila}}
\newcommand{\Ci}{\textbf{Ci}}
\newcommand{\mbC}{\textbf{mbC}}
\newcommand{\ca}{\textbf{ca}}
\newcommand{\ci}{\textbf{ci}}
\newcommand{\cl}{\textbf{cl}}
\newcommand{\cf}{\textbf{cf}}
\newcommand{\mbCci}{\textbf{mbCci}}

\newcommand{\mbCcl}{\textbf{mbCcl}}
\newcommand{\mbCciw}{\textbf{mbCciw}}
\newcommand{\LFI}{\textbf{LFI}}
\newcommand{\shorteq}{%
  \settowidth{\@tempdima}{-}
  \resizebox{\@tempdima}{\height}{=}%
}
\newenvironment{dedication}
  {
   \thispagestyle{empty}
   \vspace*{\stretch{5}}
   \itshape             
   \raggedleft          
  }
  {\par 
   \vspace{\stretch{1}} 
   \clearpage           
  }

\renewcommand{\chaptermark}[1]{%
  \ifnum\value{chapter}>0
    \markboth{Chapter \thechapter{}: #1}{}%
  \else
    \markboth{#1}{}%
  \fi}

  \usepackage[backend=biber, style=alphabetic, maxbibnames=9, sorting=nyt]{biblatex}
\usepackage{subfiles} 

\begin{document}

\pagenumbering{gobble}

\begin{titlepage}
	\centering
	\includegraphics[width=0.15\textwidth, left]{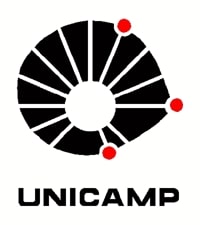}\par\vspace{1cm}
	{\scshape\large Universidade Estadual de Campinas\\ Instituto de Filosofia e Ci\^encias Humanas \par}
	\vspace{2cm}
	{\Large\itshape Guilherme Vicentin de Toledo\par}
	\vspace{1cm}
	{\LARGE\bfseries Multialgebras and Non-Deterministic Semantics applied to Paraconsistent Logics\par}
	\vspace{1cm}
	{\LARGE\bfseries Multi\'algebras e Sem\^anticas N\~ao-Determin\'isticas aplicadas a L\'ogicas Paraconsistentes\par}
	
	\vspace*{\fill}
	{\large Campinas\\ 2022\par}
\end{titlepage}


\begin{titlepage}
\centering
{\itshape Guilherme Vicentin de Toledo\par}
	\vspace{0.5cm}
	{\large\bfseries Multialgebras and Non-Deterministic Semantics applied to Paraconsistent Logics\par}
	\vspace{0.5cm}
	{\large\bfseries Multi\'algebras e Sem\^anticas N\~ao-Determin\'isticas Aplicadas a L\'ogicas Paraconsistentes\par}
\vspace{2.5cm}
\hfill\begin{minipage}[t]{0.5\textwidth}
Tese apresentada ao Instituto de Filosofia
e Ci\^encias Humanas da Universidade Estadual de
Campinas como parte dos requisitos exigidos para a
obten\c{c}\~ao do t\'itulo de Doutor em Filosofia.\par
\vspace{0.5cm}
Thesis presented to the Institute of
Philosophy and Human Sciences of the University of
Cam\-pinas in partial fulfillment of the requirements for
the degree of Doctor, in the area of Philosophy.\par
\end{minipage}
\vspace{1.5cm}

\begin{flushleft}Supervisor/Orientador: Marcelo Esteban Coniglio\\

\vspace{0.5cm}
\begin{minipage}[t]{0.5\textwidth}
ESTE TRABALHO CORRESPONDE \`A
VERS{\~A}O FINAL DA TESE
DEFENDIDA PELO ALUNO GUILHERME VICENTIN DE TOLEDO, E ORIENTADA PELO
PROF. DR. MARCELO ESTEBAN CONIGLIO.\par
\end{minipage}

\end{flushleft}

	\vspace*{\fill}
	{\large Campinas\\ 2022\par}
\end{titlepage}



\includepdf[pages={1}]{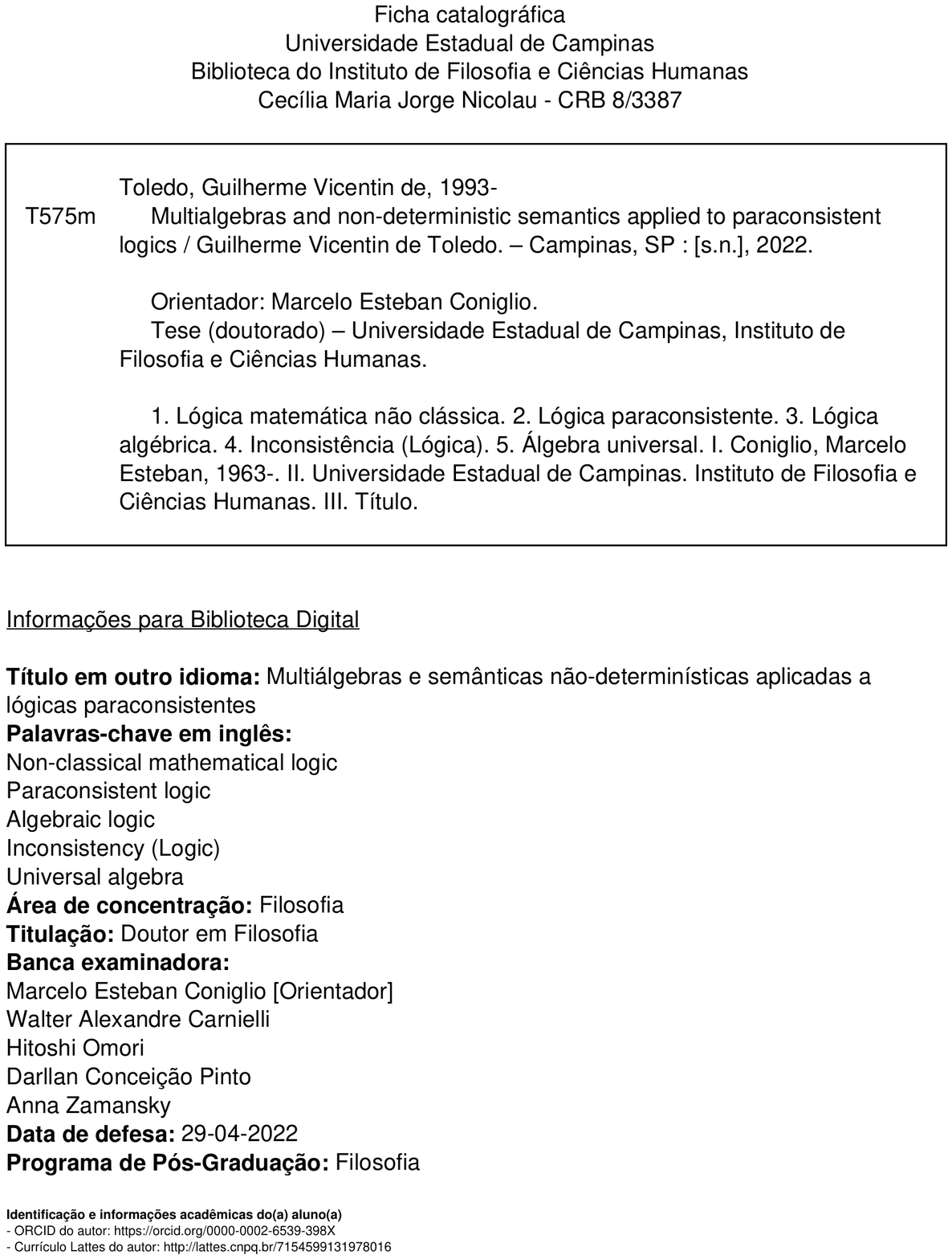}



\begin{titlepage}
\begin{center}
	\includegraphics[width=0.15\textwidth, left]{./logo}\par\vspace{1cm}
	{\scshape\large Universidade Estadual de Campinas\\ Instituto de Filosofia e Ci\^encias Humanas \par}
\end{center}
\vspace{1cm}

A Comiss\~ao Julgadora dos trabalhos de Defesa de Tese de Doutorado, composta pelos Professores Doutores a seguir descritos, em sess\~ao p\'ublica realizada em 29 de abril de 2022, considerou o candidato Guilherme Vicentin de Toledo aprovado.
\vspace{1cm}
\begin{enumerate}
\item[] Prof. Dr. Walter Alexandre Carnielli
\item[] Prof. Dr. Marcelo Esteban Coniglio
\item[] Prof. Dr. Hitoshi Omori
\item[] Prof. Dr. Darllan Concei\c{c}{\~a}o Pinto
\item[] Profa. Dra. Anna Zamansky
\end{enumerate}
\vspace{1.5cm}
\textit{A Ata de Defesa com as respectivas assinaturas dos membros encontra-se no SIGA e na Secretaria do Programa de P\'os-Gradua\c{c}\~ao em Filosofia do Instituto de Filosofia e Ci\^encias Humanas.}
\end{titlepage}


\begin{dedication}
To my father Gustavo, my brother Luis Gustavo, my grandmother Dalva,\\ my aunts Ana and Beatriz, my cousin Marina and my dog Thor.
\end{dedication}

\section*{\centering\LARGE{Acknowledgment}}
\thispagestyle{empty}
I would like to acknowledge, and thank, the \textit{Coordena\c{c}\~ao de Aperfei\c{c}oamento de Pessoal de N\'ivel Superior} (CAPES - Coordination for Improvement of Higher Education Personnel) for the PhD scholarship with process number $88882.329848/2019-01$, given through the \textit{Programa de Excel\^encia Acad\^emica} (PROEX - Academic Excellence Program), without which this work would not have been possible.
\vspace{5mm}

I thank Professor Marcelo E. Coniglio for his advisement, for which I am forever grateful. While at times he only gently steered things in a better direction, be it a definition to be clearer or a theorem to be broader, at others he delved with me into research, when doing so seemed necessary. I firmly believe that, despite the thesis being my own, some developments that it lead to are as his as they are mine, probably more.
\vspace{5mm}

Of utmost importance where the contributions by Professors Walter A. Carnielli and Hugo Mariano, who saw an earlier version of this work in my Ph.D. qualifying exam and suggested several improvements, as well as directions the ongoing research could take. If this thesis has not followed all of these directions, and it hasn't, it was only for lack of more time.
\vspace{5mm}

Many other ameliorations and corrections, all of which I tried to incorporate into the final text, were brought forward in my thesis defense by Professors Walter Carnielli, Itala L. D'Ottaviano, H{\'e}rcules A. Feitosa, Hitoshi Omori, Darllan C. Pinto and Anna Zamansky, who all in addition suggested invaluable ideas for future work.
\vspace{5mm}

And to all those who have otherwise influenced this work, including all students, professors and staff of the \textit{Centro de L\'ogica, Epistemologia e Hist\'oria da Ci\^encia} (CLE - Center for Logic, Epistemology and History of Science), the \textit{Instituto de Filosofia e Ci\^encias Humanas} (IFCH - Institute of Philosophy and Human Sciences), and the \textit{Instituto de Matem\'atica e Estat\'istica da Universidade de S\~ao Paulo} (IME-USP - Institute of Mathematics and Statistics of the University of S\~ao Paulo), I am also deeply grateful.

\newpage
\thispagestyle{empty}
\vspace*{\fill}
\begin{flushright}
\textit{``But man is not made for defeat.\\
A man can be destroyed but not defeated.``\\
(\textit{The Old Man and the Sea}, by Ernest Hemingway)}
\end{flushright}
\newpage

\section*{\centering\LARGE{Abstract}}
\thispagestyle{empty}

This work is divided between two main areas: in the theory of multialgebras, we focus mostly on a new definition of what a freely generated object should be in their category, and on how this category is equivalent to another with partially ordered algebras as objects; we then use non-deterministic semantics, specially those we have named restricted Nmatrices, on paraconsistent logics and some systems dealing with a new presentation of the natural concept of incompatibility, which generalizes inconsistency.

In algebra, we will focus on the non-deterministic ones, also known as multialgebras, whose operations return non-empty subsets of their universes. While the category of algebras over a signature has freely generated objects, which in a sense permit for the unique extension of functions to homomorphisms, the category of multialgebras over a given signature does not have elements with comparable properties. To circumvent this problem, we widen our understanding of algebras of formulas: if a multialgebra generalizes an algebra by having multiple results for a given operation, a multialgebra of formulas should generalize an algebra of formulas by having multiple possibilities for applying a connective to given formulas. Concerning the category of multialgebras itself, we offer an equivalence between it and a category avoiding non-determinism altogether, relying instead on ordered Boolean-like algebras as objects.

On the part devoted to logic, our goals are again roughly twofold: firstly, some logics of formal inconsistency, \textit{exempli gratia} those found in da Costa's hierarchy, cannot be characterized by finite Nmatrices. In what is a very natural development, a restricted Nmatrix (or RNmatrix) restricts those homomorphisms to be taken into consideration when evaluating the validity of a deduction according to an Nmatrix. We show how this distinction gives far greater expressiveness to finite RNmatrices, enough to adequately characterize da Costa's systems and provide decision methods for those logics, both based on row-branching, row-eliminating truth tables, and tableau semantics. In another direction, we generalize logics of formal inconsistency to new systems built around the notion of incompatibility: the \textit{Leitmotiv} being that having two incompatible formulas to simultaneously hold trivializes a deduction, and as a special case, a formula is consistent when it is incompatible with its negation. We show how this notion extends that of inconsistency in a non-trivial way, presenting conservative translations for many simple inconsistent systems into logics of incompatibility; we also provide semantics built on RNmatrices for these new logics, and prove that they can not be characterized by more standard methods.

\vspace*{\fill}
\textbf{Keywords}: Non-classical mathematical logic; Paraconsistent logic; Algebraic logic; Inconsistency (Logic); Universal algebra.

\newpage

\section*{\centering\LARGE{Resumo}}
\thispagestyle{empty}

\foreignlanguage{portuguese}{
Este trabalho est\'a dividido entre duas grandes \'areas: na teoria de multi\'algebras, focamos majoritariamente em uma nova defini\c{c}\~ao do que um objeto livremente gerado deveria ser em sua categoria e em como esta categoria \'e equivalente a outra com \'algebras parcialmente ordenadas como objetos; ent\~ao usamos sem\^anticas n\~ao-determin\'isticas, especialmente aquela que nomeamos Nmatrizes restritas, nas l\'ogicas paraconsistentes e em alguns sistemas lidando com uma nova apresenta\c{c}\~ao do conceito natural de incompatibilidade, que generaliza o conceito de inconsist\^encia.}

\foreignlanguage{portuguese}{
Em \'algebra, nos focaremos nas n\~ao-determin\'isticas, tamb\'em conhecidas como multi\'algebras, cujas opera\c{c}\~oes retornam subconjuntos n\~ao vazios de seus universos. Enquanto a categoria de \'algebras sobre uma assinatura possui objetos livremente gerados, os quais permitem em certo sentido a extens\~ao \'unica de fun\c{c}\~oes a homomorfismos, a categoria de multi\'algebras sobre uma assinatura dada n\~ao possui elementos com propriedades compar\'aveis. Para contornar este problema, estendemos o significado de uma \'algebra de f\'ormulas: se uma multi\'algebra generaliza uma \'algebra ao ter m\'ultiplos resultados para uma dada opera\c{c}\~ao, uma multi\'algebra de f\'ormulas generaliza uma \'algebra de f\'ormulas ao ter m\'ultiplas possibilidade para a aplica\c{c}\~ao de um conectivo a f\'ormulas dadas. Quanto \`a categoria de multi\'algebras propriamente dita, oferecemos uma equival\^encia entre ela e uma categoria livre de n\~ao-determinismo, que alternativamente possui \'algebras ordenadas, semelhantes a \'algebras de Boole, como objetos.}

\foreignlanguage{portuguese}{
Na parte dedicada \`a l\'ogica, nossos objetivos s\~ao novamente dois: primeiramente, algumas l\'ogicas de inconsist\^encia formal, \textit{exempli gratia} aquelas da hierarquia de da Costa, n\~ao podem ser caracterizadas por Nmatrizes finitas. No que \'e um desenvolvimento muito natural, uma Nmatriz restrita, ou RNmatriz, restringe aquelas homomorfismos que devem ser considerados quando testamos a validade de uma dedu\c{c}\~ao por uma Nmatriz. Mostramos como esta distin\c{c}\~ao prov\^e as RNmatrizes finitas com poder expressivo muito superior, suficientente para adequadamente caracterizar os sistemas de da Costa e dar a eles m\'etodos de decis\~ao, tanto baseados em tabelas de verdade quanto em sem\^anticas de tableaux. Em outra dire\c{c}\~ao, generalizamos as l\'ogicas de inconsist\^encia formal a sistemas constru\'idos em torno da no\c{c}\~ao de incompatibilidade: o \textit{Leitmotiv} sendo que duas f\'ormulas incompat\'iveis simultaneamente verdadeiras trivializam uma dedu\c{c}\~ao, e como um caso especial, uma f\'ormula \'e consistente quando \'e incompat\'ivel com sua nega\c{c}\~ao. Mostramos como essa no\c{c}\~ao estende aquela de inconsist\^encia de maneira n\~ao-trivial, apresentando tradu\c{c}\~oes conservativas para muitos dos sistemas inconsistentes mais simples em l\'ogicas de incompatilidade, apresentamos sem\^anticas constru\'idas com RNmatrizes para essas novas l\'ogicas e mostramos que elas n\~ao podem ser caracterizadas por m\'etodos mais usuais.}

\vspace*{\fill}
\foreignlanguage{portuguese}{
\textbf{Palavras-chave}: L\'ogica matem\'atica n\~ao cl\'assica; L\'ogica paraconsistente; L\'ogica alg\'ebrica; Inconsist\^encia (L\'ogica); \'Algebra universal.}

\pagestyle{empty}
\addtocontents{toc}{\protect\thispagestyle{empty}}
\shorttoc{Abridged Table of Contents}{1}
\thispagestyle{empty}
\newpage

\shorttoc{Extended Table of Contents}{3}
\thispagestyle{empty}
\newpage

\tableofcontents
\cleardoublepage

\pagenumbering{arabic}
\setcounter{page}{16}

\newpage
\thispagestyle{plain}
\hspace{0pt}
\vfill
\begin{figure}[h]
\centering
\includegraphics[width=\textwidth]{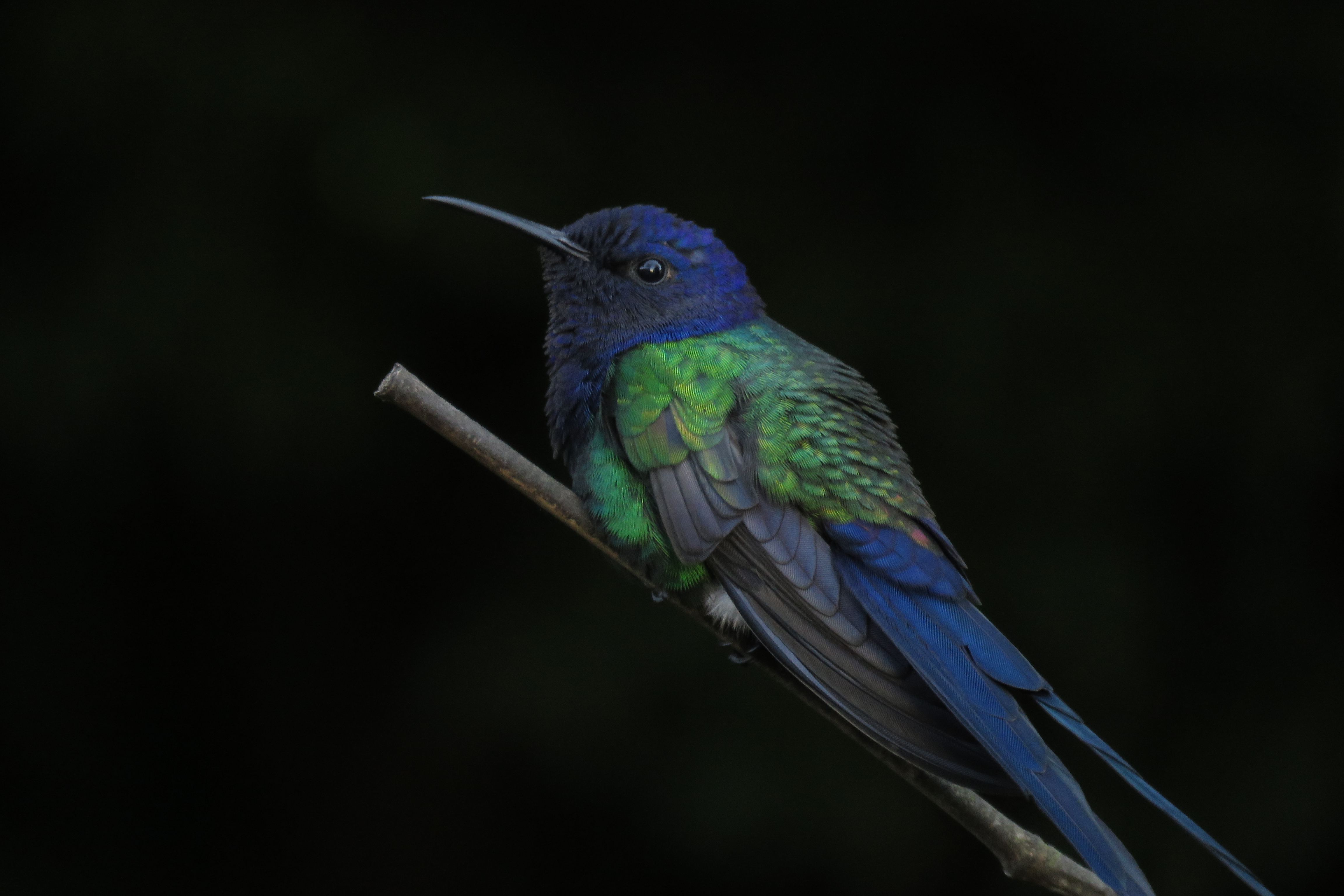}
\caption*{Specimen of \textit{Eupetomena macroura} resting, Monte Alegre do Sul, Brazil.\\Photographed by Guilherme Vicentin de Toledo, all rights reserved.}
\end{figure}
\vfill
\hspace{0pt}

\pagestyle{fancy}

\begin{refsegment}
\defbibfilter{notother}{not segment=\therefsegment}
\chapter*{Introduction}\label{Chapter0}
\chaptermark{Introduction}
\addcontentsline{toc}{chapter}{Introduction}

The research that makes up the bulk of this work may be divided in two grand areas, what explains our division of the thesis in two parts: algebra, corresponding to Part \ref{Part1}, and logic, to Part \ref{Part2}. 

In algebra, we focus on the subject of multialgebras: we present a new concept of \textit{weakly free multialgebras}, which we designed in order to offer a generalization of free algebras; furthermore, we show how the category of multialgebras, as usually defined in the context of non-deterministic semantics, can be presented alternatively as a category of algebras equipped with orders that are compatible with the underlying operations.

In logic, we will delve into paraconsistent logics, together with their generalizations, and non-deterministic semantics: we generalize non-deterministic matrices to \textit{restricted non-deterministic matrices} (also known as RNmatrices), a semantical tool that uses both multialgebras and restrictions over the set of valuations to be considered, and apply this methodology to give an extend analysis of da Costa's hierarchy that includes new decision methods for its logics; we also give a new formalization of a natural concept, that of incompatibility, and treat the resulting logical systems again with RNmatrices.

\subsection*{Multialgebras}

In Part \ref{Part1} we will deal mainly with multialgebras, also known as hyperalgebras or\\ non-deterministic algebras. Chapter \ref{Chapter1} gives a brief introduction to the subject, broaching the definition of multialgebras, as well as those of homomorphisms between multialgebras, submultialgebras, the interpretation of formulas, and so on. We spend a few pages over alternative definitions for a homomorphism of multialgebras that allow for stronger representation theorems, and discuss briefly how these, sadly, are not very practical to use in non-deterministic semantics. Given its frequent use in our work, we also use this chapter to define, for completeness sake, lattices, Boolean algebras and Heyting algebras.

Chapter \ref{Chapter2} offers an alternative solution to a classical problem: the category of algebras, on a given signature, possesses objects satisfying the universal mapping property, namely free algebras; meanwhile, the category of multialgebras over a signature $\Sigma$, denoted by $\textbf{MAlg}(\Sigma)$, does not. Algebras of formulas over a signature $\Sigma$ and variables $\mathcal{V}$, that we denote by $\textbf{F}(\Sigma, \mathcal{V})$, are then extended to structures that admit non-determinism: after all, if a multialgebra generalizes an algebra by allowing multiple results to an operation, a multialgebra of formulas should allow many formulas to be obtained from the application of a connective to given formulas. This new concept is shown to have many desirable characterizations, some inspired by linear algebra, others graph-theoretic in nature, motivating us to coin the nomenclature of weakly free multialgebras for them. They are shown to greatly simplify the proof that $\textbf{MAlg}(\Sigma)$ does not have free objects, and some other categorical considerations are then developed by the end of the chapter.

Concerning the category $\textbf{MAlg}(\Sigma)$ itself, Chapter \ref{Chapter3} proves this category is equivalent to another that avoids non-determinism altogether by considering ordered, Boolean-like algebras. Intuitively, one wants to capture both the operations and order of the natural algebra over the powerset of the universe $A$ of a multialgebra $\mathcal{A}$, where: the order is the usual order for a powerset; and, given non-empty subsets $A_{1}$ trough $A_{n}$ of $A$, an $n$-ary operation $\sigma$ on $(A_{1}, \dotsc , A_{n})$ is given by the union of $\sigma(a_{1}, \dotsc , a_{n})$ for $(a_{1}, \dotsc , a_{n})$ in $A_{1}\times\cdots\times A_{n}$. This presents an alternative for those logicians wishing to stick to deterministic semantics without any loss of expressiveness, and can be generalized to offer equivalences for other categories with non-deterministic algebras, such as partial multialgebras.

\subsection*{Paraconsistent Logic}

Part \ref{Part2} involves paraconsistent logics, as well as non-deterministic semantics and broader systems. It starts with Chapter \ref{Chapter4}, which provides a brief introduction to formal logic, and semantics of logical matrices, non-deterministic matrices (also known as Nmatrices) and restricted matrices. The high point of the chapter, however, is the definition of restricted non-deterministic matrices, also known as restricted Nmatrices or RNmatrices, which very naturally combine restricted and non-deterministic matrices: some theoretical considerations are made about these semantics, as well as a brief analysis of its previous, unrecognized uses in the literature. In essence, an RNmatrix is a triple $(\mathcal{A}, D, \mathcal{F})$, with $\mathcal{A}$ a $\Sigma$-multialgebra, $D$ a subset of its universe and $\mathcal{F}$ a set of homomorphisms $\nu:\textbf{F}(\Sigma, \mathcal{V})\rightarrow\mathcal{A}$. Their motivation is to give finite semantics to logical systems that are not characterizable by finite Nmatrices, such as the logics of formal inconsistency ($\textbf{LFI}'s$) between $\textbf{mbCcl}$ and $\textbf{Cila}$, and the whole hierarchy of da Costa.

The $C$-systems of da Costa, due to their complexity, are treated separately in Chapter \ref{Chapter5}: in it, we start by defining these logics first devised in order to formalize the notion of inconsistency, or paraconsistency. We then proceed to provide RNmatrices of $n+2$ elements capable of characterizing each and every $C_{n}$ of the hierarchy, starting from $C_{2}$ to fix ideas: the intuition is that an RNmatrix for $C_{n}$ must have two classical values, standing from true and false, as well as $n$ inconsistent values, each standing for a different degree of inconsistency achieved in the logic; of course, this suggests that the $n$-th logic in the hierarchy could be regarded as an $n+2$-valued, non-deterministic logic. Furthermore, we show how these finite RNmatrices can be made into decision methods for da Costa's systems through both row-branching, row-eliminating truth tables, and tableau semantics.

Chapter \ref{Chapter6} extends the RNmatrices for $C_{n}$ from the previous chapter to allow for model-theoretic considerations: while we started with $n+2$-valued restricted non-deterministic matrices, these can be shown to be constructed from the two-valued Boolean algebra as swap structures; the next logical step is to construct these swap structures over any Boolean algebras, leading to a class of RNmatrices capable itself of characterizing $C_{n}$. Further results include a brief combinatorial description of the snapshots found in these swap structures, as well as the construction of a category of swap structures for $C_{n}$, which is then proven to be isomorphic to the category of non-trivial Boolean algebras; this has important applications to the model theory of da Costa's hierarchy.

Chapter \ref{Chapter7} introduces logics of incompatibility: a first attempt of formalizing the notion of incompatibility from natural language would be to declare two formulas incompatible if, and only if, together they trivialize a deduction. As it is done in logics of formal inconsistency, we weaken that condition to state that two incompatible formulas, a concept which may be primitive, are capable of trivializing an argument; we will denote by $\alpha\uparrow\beta$ the fact that $\alpha$ and $\beta$ are incompatible. We start by defining some very basic systems of incompatibility, as of yet without negation, characterize then with RNmatrices and provide decision methods: of $\bI^{-}$ we ask nothing, while $\bI$ must have a commutative incompatibility connective and $\bIpr$ propagates incompatibility under some conditions. We finish the chapter by discussing axioms that collapse $\uparrow$ to its classical interpretation, and finally comparing our approach to a preexisting formalization of incompatibility, that of Brandom.

It is natural, once we have a connective for incompatibility, to consider its interplay with negation: the systems of Chapter \ref{Chapter7} are devoid of negation, by design, but Chapter \ref{Chapter8} consider those logics now equipped with negation and some axioms governing its interaction with incompatibility, most of them heavily inspired by the most common axioms for paraconsistent systems: while $\nbI$ only adds a negation to $\bI$ satisfying \textit{tertium non datur}, $\nbIciw$, $\nbIci$ and $\nbIcl$ generalize the logics of paraconsistency $\mbCciw$, $\mbCci$ and $\mbCcl$, respectively. Of course, we then provide these systems with characterizing RNmatrices, as well as decision methods through row-branching, row-eliminating truth tables and tableau calculi, and explain why these semantics, instead of more classical ones, are necessary: not only our basic incompatible systems are not algebraizable by Blok and Pigozzi, they are also not characterizable by either finite Nmatrices or finite restricted matrices. We finish by analyzing the generalizations of logical matrices we have broached here, how do they relate to each other and to other possible generalizations.

Finally, Chapter \ref{Chapter9} studies an equivalence merely implied in previous chapters: in logics of formal inconsistency, a formula $\alpha$ being consistent is recurrently expressed by $\circ\alpha$, where the connective $\circ$ stands precisely for consistency. When dealing with incompatibility, the fact consistency is adequately reintroduced as incompatibility with negation becomes apparent: that is, $\circ\alpha$ may be viewed as $\alpha\uparrow\neg\alpha$. Accordingly, we define a function from the logics of formal inconsistency into the logics of incompatibility that is not only a translation, but a conservative one nonetheless. This seems to imply that our interpretation of incompatibility in logic strictly extends the notion of inconsistency in a non-trivial way.

\subsection*{Publications}

A considerable portion of this thesis was made available as preprints in

\begin{enumerate}

\item \fullcite{AbsFreeHyp};

\item \fullcite{CostaRNmatrix};

\item \fullcite{RestrictedSwap};

\item \fullcite{Frominconsistency},

\end{enumerate}

and published in

\begin{enumerate}

\item \fullcite{WeaklyFreeMultialgebras};

\item \fullcite{TwoDecisionProcedures}.

\end{enumerate}

\subsection*{Presentations}

The work found in this thesis has also lead to the following presentations.

\begin{enumerate}

\item ``\textbf{RNmatrices for da Costa's hierarchy}'', in \textit{I Enc(ue-o)ntro de L\'ogica Brasil\\ Col(o-\^o)mbia} (1st Meeting Brazil-Colombia in Logic), December of 2021.

\item ``\textbf{Uma categoria de \'algebras ordenadas equivalente a uma categoria de hiper\'algebras}'', in \textit{I Encontro Brasileiro em Teoria das Categorias}, January of 2021.

\item ``\textbf{Sem\^anticas n\~ao-determin\'isticas para l\'ogicas n\~ao-cl\'assicas: Uma abordagem da\\ perspectiva de Teoria de Modelos e de \'Algebra Universal}'', in \textit{Encontro Brasileiro de L\'ogica 2019} (19th Brazilian Logic Conference), May of 2019.

\item ``\textbf{Sem\^anticas n\~ao-determin\'isticas para l\'ogicas n\~ao-cl\'assicas: Uma abordagem da\\ perspectiva de Teoria de Modelos e de \'Algebra Universal}'', in \textit{XVIII Encontro Nacional da ANPOF}, October of 2018.

\item ``\textbf{Nondeterministic Semantics for Nonclassical Logics: An Approach from the\\ Perspective of Model Theory and Universal Algebra}'', in \textit{Fourth Workshop CLE-Buenos Aires Logic Group}, April of 2018.
\end{enumerate}

\end{refsegment}

\part{Multialgebras}\label{Part1}

\begin{refsegment}
\defbibfilter{notother}{not segment=\therefsegment}
\setcounter{chapter}{0}
\chapter{$\Sigma$-Multialgebras and lattices}\label{Chapter1}\label{Chapter 1}

A multialgebra, or hyperalgebra\index{Hyperalgebra}, is a generalization of the notion of algebra usually found in the context of universal algebras, see \cite{Burris} for a standard approach to universal algebra. The main objective of such a generalization is to address the possibility that the outcome of an operation may be diffuse, that is, non-deterministic: we may have an idea of what the outcome should be, but not be certain about it. The first known appearance of multialgebras in literature may be found in \cite{Marty}.

A collection of disjoint sets $\Sigma=\{\Sigma_{n}\}_{n\in\mathbb{N}}$ indexed by $\mathbb{N}$ will be called a signature\index{Signature}\label{signature}; the elements of the sets $\Sigma_{n}$ will be called functional symbols of arity $n$, or $n$-ary functional symbols. For simplicity, we will denote $\bigcup_{n\in\mathbb{N}}\Sigma_{n}$ also by $\Sigma$.

A pair $\mathcal{A}=(A, \{\sigma_{\mathcal{A}}\}_{\sigma\in\Sigma})$ is said to be, where $\mathcal{P}(A)$\label{powerset} denotes the powerset (also known as power set) of $A$:
\begin{enumerate}
\item a $\Sigma$-algebra\index{Algebra} if, for every $\sigma\in\Sigma_{n}$, $\sigma_{\mathcal{A}}$ is a function of the form
\[\sigma_{\mathcal{A}}:A^{n}\rightarrow A;\]
\item a $\Sigma$-multialgebra\index{Multialgebra} if $A\neq\emptyset$ and, for every $\sigma\in\Sigma_{n}$, $\sigma_{\mathcal{A}}$ is a function of the form
\[\sigma_{\mathcal{A}}:A^{n}\rightarrow \mathcal{P}(A)\setminus\{\emptyset\};\]
\end{enumerate}

The set $A$ is called the universe\index{Universe} of $\mathcal{A}$. 

Given a $\Sigma$-algebra $\mathcal{A}=(A, \{\sigma_{\mathcal{A}}\}_{\sigma\in\Sigma})$, one can always define the $\Sigma$-multialgebra $\overline{\mathcal{A}}=(A, \{\sigma_{\overline{\mathcal{A}}}\}_{\sigma\in\Sigma})$ such that, for $\sigma\in\Sigma_{n}$ and $a_{1}, \dotsc  , a_{n}\in A$, 
\[\sigma_{\overline{\mathcal{A}}}(a_{1}, \dotsc  , a_{n})=\{\sigma_{\mathcal{A}}(a_{1}, \dotsc  , a_{n})\};\]
it is clear that $\overline{\mathcal{A}}$ carries the same information that $\mathcal{A}$, and so one can see $\Sigma$-algebras as $\Sigma$-multialgebras. For most of our studies here we will focus mainly on multialgebras.


\section{Homomorphisms}


\subsection{Single-valued homomorphisms}

Given $\Sigma$-algebras $\mathcal{A}=(A, \{\sigma_{\mathcal{A}}\}_{\sigma\in\Sigma})$ and $\mathcal{B}=(B, \{\sigma_{\mathcal{B}}\}_{\sigma\in\Sigma})$, a function $\varphi:A\rightarrow B$ is said to be a homomorphism from $\mathcal{A}$ to $\mathcal{B}$ if, for every $\sigma\in\Sigma_{n}$ and $a_{1}, \dotsc  , a_{n}\in A$ we have that
\[\varphi(\sigma_{\mathcal{A}}(a_{1}, \dotsc  , a_{n}))=\sigma_{\mathcal{B}}(\varphi(a_{1}), \dotsc  , \varphi(a_{n})).\]

Given $\Sigma$-multialgebras $\mathcal{A}=(A, \{\sigma_{\mathcal{A}}\}_{\sigma\in\Sigma})$ and $\mathcal{B}=(B, \{\sigma_{\mathcal{B}}\}_{\sigma\in\Sigma})$, we would like to define once again a morphism between the two of them. But in the case of multialgebras, one can come up with many possible straightforward definitions of homomorphisms, all of them starting with a function $\varphi:A\rightarrow B$. The following two definitions are the most useful to our purposes:
\begin{enumerate}
\item if, for all $\sigma\in\Sigma_{n}$ and $a_{1}, \dotsc  , a_{n}\in A$, 
\[\{\varphi(a)\ :\ a\in\sigma_{\mathcal{A}}(a_{1}, \dotsc  , a_{n})\}\subseteq \sigma_{\mathcal{B}}(\varphi(a_{1}), \dotsc  , \varphi(a_{n})),\]
we call $\varphi$ a homomorphism\index{Homomorphism}, or $\Sigma$-homomorphism;
\item we call $\varphi$ a full homomorphism\index{Homomorphism, Full} if the above condition is replaced by 
\[\{\varphi(a)\ :\  a\in\sigma_{\mathcal{A}}(a_{1}, \dotsc  , a_{n})\}=\sigma_{\mathcal{B}}(\varphi(a_{1}), \dotsc  , \varphi(a_{n})).\]
\end{enumerate}
If the function $\varphi:A\rightarrow B$ is a homomorphism from $\mathcal{A}$ to $\mathcal{B}$, we will simply write $\varphi:\mathcal{A}\rightarrow\mathcal{B}$\label{homomorphism}.

\begin{theorem}
The class of all $\Sigma$-multialgebras becomes a category $\textbf{MAlg}(\Sigma)$\label{MAlg} or $\textbf{MAlg}_{=}(\Sigma)$\label{MAlg=} when the set of morphisms between two $\Sigma$-multialgebras $\mathcal{A}$ and $\mathcal{B}$ is, respectively:
\begin{enumerate}
\item the set of all homomorphisms between $\mathcal{A}$ and $\mathcal{B}$;
\item the set of all full homomorphisms between $\mathcal{A}$ and $\mathcal{B}$.
\end{enumerate}
In both cases, the composition of morphisms is the usual composition of functions.
\end{theorem}

\begin{proof}
We must show that all two of those alleged categories have identity morphisms and that their compositions are well-defined, meaning that composing two morphisms returns again a morphism; clearly there is no need to show the associativity of composition, since it is know that the composition of functions is indeed associative.

For every multialgebra $\mathcal{A}$ we consider the morphism $Id_{\mathcal{A}}:\mathcal{A}\rightarrow \mathcal{A}$ given by, for every $a\in A$, $Id_{\mathcal{A}}(a)=a$. This morphism is the desired identity morphism for $\mathcal{A}$ in both categories, since: it is, in fact, a morphism in the two of them, given that it is a full homomorphism and therefore also a homomorphism; it is the identity for the composition of functions and therefore, for any morphisms $\varphi:\mathcal{A}\rightarrow \mathcal{B}$ and $\psi:\mathcal{C}\rightarrow\mathcal{A}$, 
\[\varphi\circ Id_{\mathcal{A}}=\varphi\quad\text{and}\quad Id_{\mathcal{A}}\circ\psi=\psi.\]

Now to prove that the composition is well-defined, fix $\sigma\in\Sigma_{n}$ and $a_{1}, \dotsc  , a_{n}\in A$.
\begin{enumerate}
\item If $\varphi:\mathcal{A}\rightarrow\mathcal{B}$ and $\psi:\mathcal{B}\rightarrow\mathcal{C}$ are simply homomorphisms,
\[\{\varphi(a)\ :\  a\in\sigma_{\mathcal{A}}(a_{1}, \dotsc  , a_{n})\}\subseteq\sigma_{\mathcal{B}}(\varphi(a_{1}), \dotsc  , \varphi(a_{n}))\] 
and therefore 
\[\{\psi\circ\varphi(a)\ :\  a\in\sigma_{\mathcal{A}}(a_{1}, \dotsc  , a_{n})\}\subseteq\{\psi(b)\ :\  b\in \sigma_{\mathcal{B}}(\varphi(a_{1}), \dotsc  , \varphi(a_{n}))\};\]
since $\psi$ is a homomorphism, we have that
\[\{\psi(b)\ :\  b\in \sigma_{\mathcal{B}}(\varphi(a_{1}), \dotsc  , \varphi(a_{n}))\}\subseteq \sigma_{\mathcal{C}}(\psi\circ\varphi(a_{1}), \dotsc  , \psi\circ\varphi(a_{n})),\]
and from that $\psi\circ\varphi$ is also a homomorphism.

\item If $\varphi:\mathcal{A}\rightarrow\mathcal{B}$ and $\psi:\mathcal{B}\rightarrow\mathcal{C}$ are full homomorphisms, is enough to replace all "$\subseteq$" on the proof above by equalities to obtain a proof that $\psi\circ\varphi$ is also a full homomorphism.
\end{enumerate}
\end{proof}

\begin{definition}
A full homomorphism $\varphi:\mathcal{A}\rightarrow\mathcal{B}$ is said to be an isomorphism\index{Isomorphism} if $\varphi:A\rightarrow B$ is a bijection.\footnote{We do not attempt to define isomorphisms for homomorphisms that are not full since the latter class is not closed under inverses: the inverse of a non-full homomorphism is an antihomomorphism, which satisfies $\sigma_{\mathcal{B}}(\varphi(a_{1}), \dotsc, \varphi(a_{n}))\subseteq \{\varphi(a)\ :\ a\in\sigma_{\mathcal{A}}(a_{1}, \dotsc a_{n})\}$; considering antihomomorphisms does not make the underlying theory uninteresting, but it does make the theory much harder.}
\end{definition}

\begin{proposition}
Let $\mathcal{A}$ and $\mathcal{B}$ be $\Sigma$-multialgebras and $\varphi:A\rightarrow B$ be a bijective function with inverse $\psi:B\rightarrow A$: if $\varphi$ is a full homomorphism, so is $\psi$.
\end{proposition}

\begin{proof}
Let $\sigma\in\Sigma_{n}$, $b_{1}, \dotsc  , b_{n}\in B$ and $a_{1}=\psi(b_{1}), \dotsc  , a_{n}=\psi(b_{n})$, so that $b_{1}=\varphi(a_{1}), \dotsc  , b_{n}=\varphi(a_{n})$. We have that
\[ \{\psi(b)\ :\  b\in \sigma_{\mathcal{B}}(b_{1}, \dotsc  , b_{n})\}=\{\psi(b)\ :\  b\in\sigma_{\mathcal{B}}(\varphi(a_{1}), \dotsc  , \varphi(a_{n}))\}=\]
\[\{\psi(b)\ :\  b\in\{\varphi(a)\ :\  a\in\sigma_{\mathcal{A}}(a_{1}, \dotsc  , a_{n})\}\}=\sigma_{\mathcal{A}}(a_{1}, \dotsc  , a_{n})=\sigma_{\mathcal{A}}(\psi(b_{1}), \dotsc  , \psi(b_{n})),\]
and therefore $\psi$ is indeed a full homomorphism.
\end{proof}


\subsection{Multi-valued homomorphisms}\label{Multi-valued homomorphisms}

One natural consideration when defining morphisms between multialgebras is that, if the operations are of a non-deterministic nature, perhaps the morphisms should be as well.

Although there exist a plethora of possible definitions in this case, we will mention only two of them, even avoiding considerations about isomorphisms for they are out of scope for this text. Given a signature $\Sigma$ and $\Sigma$-multialgebras $\mathcal{A}=(A, \{\sigma_{\mathcal{A}}\}_{\sigma\in\Sigma})$ and $\mathcal{B}=(B, \{\sigma_{\mathcal{B}}\}_{\sigma\in\Sigma})$, a function $\varphi:A\rightarrow\mathcal{P}(B)\setminus\{\emptyset\}$ is said to be a multihomomorphism\index{Multihomomorphism} if, for all $\sigma\in\Sigma_{n}$ and $a_{1}, \dotsc  , a_{n}\in A$,
\[\bigcup_{a\in \sigma_{\mathcal{A}}(a_{1}, \dotsc  , a_{n})}\varphi(a)\subseteq\bigcup_{(b_{1}, \dotsc  , b_{n})\in \varphi(a_{1})\times\cdots\times\varphi(a_{n})}\sigma_{\mathcal{B}}(b_{1}, \dotsc  , b_{n});\]
the same function is said to be a full multihomomorphism\index{Multihomomorphism, Full} if this conditions is replaced by 
\[\bigcup_{a\in \sigma_{\mathcal{A}}(a_{1}, \dotsc  , a_{n})}\varphi(a)=\bigcup_{(b_{1}, \dotsc  , b_{n})\in \varphi(a_{1})\times\cdots\times\varphi(a_{n})}\sigma_{\mathcal{B}}(b_{1}, \dotsc  , b_{n}).\]
We will denote a multihomomorphism $\varphi$ between $\mathcal{A}$ and $\mathcal{B}$ simply by $\varphi:\mathcal{A}\rightarrow\mathcal{B}$.

\begin{lemma}\label{lemma about unions}
For an $n\in\mathbb{N}$, take sets $X_{1}, \dotsc  , X_{n}$; for $X=\bigcup_{i=1}^{n}X_{i}$, let $\{Y_{x}\}_{x\in X}$ be a family indexed by $X$; then
\[\bigcup_{(x_{1}, \dotsc  , x_{n})\in X_{1}\times\cdots\times X_{n}} Y_{x_{1}}\times\cdots\times Y_{x_{n}}\subseteq\bigcup_{x_{1}\in X_{1}}Y_{x_{1}}\times\cdots\times\bigcup_{x_{n}\in X_{n}}Y_{x_{n}}.\]
If $n=1$, we have instead an equality.
\end{lemma}

\begin{proof}
Suppose $(y_{1}, \dotsc  , y_{n})$ belongs to the left side of the inequality, and there must exist\\ $(x_{1}, \dotsc  , x_{n})$ in $X_{1}\times\cdots\times X_{n}$ such that $(y_{1}, \dotsc  , y_{n})\in Y_{x_{1}}\times\cdots\times Y_{x_{n}}$.

Now, for any $i\in\{1, \dotsc  , n\}$, since $y_{i}\in Y_{x_{i}}$ for an $x_{i}\in X_{i}$, $y_{i}\in \bigcup_{x\in X_{i}}Y_{x}$. It follows that $(y_{1}, \dotsc  , y_{n})\in\bigcup_{x_{1}\in X_{1}}Y_{x_{1}}\times\cdots\times\bigcup_{x_{n}\in X_{n}}Y_{x_{n}}$.

The equality if $n=1$ is trivial.
\end{proof}

\begin{theorem}
The class of all $\Sigma$-multialgebras becomes a category $\textbf{MMAlg}(\Sigma)$\label{MMAlg} or $\textbf{MMAlg}_{=}(\Sigma)$\label{MMAlg=} when the set of morphisms between two $\Sigma$-multialgebras $\mathcal{A}$ and $\mathcal{B}$ is, respectively:
\begin{enumerate}
\item the set of all multihomomorphisms between $\mathcal{A}$ and $\mathcal{B}$;
\item the set of all full multihomomorphisms between $\mathcal{A}$ and $\mathcal{B}$.
\end{enumerate}
In both cases, the composition $\psi\circ\varphi$ of multihomomorphisms $\varphi:\mathcal{A}\rightarrow\mathcal{B}$ and $\psi:\mathcal{B}\rightarrow\mathcal{C}$ is given by, on an element $a\in A$, $\psi\circ\varphi(a)=\bigcup_{b\in\varphi(a)}\psi(b)$.
\end{theorem}

\begin{proof}
We must show the existence of identity morphisms and that the composition of morphisms is well-defined and associative.

For every $\Sigma$-multialgebra $\mathcal{A}$ we consider the morphism $Id_{\mathcal{A}}:\mathcal{A}\rightarrow\mathcal{A}$ given by, for every $a\in A$, $Id_{\mathcal{A}}(a)=\{a\}$. It is a full multihomomorphism, and therefore also a multihomomorphism, given that, for $\sigma\in\Sigma_{n}$ and $a_{1}, \dotsc  , a_{n}\in A$,
\[\bigcup_{a\in \sigma_{\mathcal{A}}(a_{1}, \dotsc  , a_{n})}Id_{\mathcal{A}}(a)=\bigcup_{a\in\sigma_{\mathcal{A}}(a_{1}, \dotsc  , a_{n})}\{a\}=\sigma_{\mathcal{A}}(a_{1}, \dotsc  , a_{n})=\bigcup_{(b_{1}, \dotsc  , b_{n})\in \{a_{1}\}\times\cdots\times\{a_{n}\}}\sigma_{\mathcal{A}}(b_{1}, \dotsc  , b_{n})=\]
\[\bigcup_{(b_{1}, \dotsc  , b_{n})\in Id_{\mathcal{A}}(a_{1})\times\cdots\times Id_{\mathcal{A}}(a_{n})}\sigma_{\mathcal{A}}(b_{1}, \dotsc  , b_{n}).\]
For any multihomomorphisms $\varphi:\mathcal{A}\rightarrow\mathcal{B}$ and $\psi:\mathcal{C}\rightarrow\mathcal{A}$ and elements $a\in A$ and $c\in C$, we have that
\[\varphi\circ Id_{\mathcal{A}}(a)=\bigcup_{d\in Id_{\mathcal{A}}(a)}\varphi(d)=\bigcup_{d\in \{a\}}\varphi(d)=\varphi(a)\]
and
\[Id_{\mathcal{A}}\circ\psi(c)=\bigcup_{e\in\psi(c)}Id_{\mathcal{A}}(e)=\bigcup_{e\in\psi(c)}\{e\}=\psi(c),\]
so that $\varphi\circ Id_{\mathcal{A}}=\varphi$ and $Id_{\mathcal{A}}\circ\psi=\psi$, and $Id_{\mathcal{A}}$ is indeed an identity.

To see that the composition is well-defined, fix $\sigma\in\Sigma$ and $a_{1}, \dotsc  , a_{n}\in A$.

\begin{enumerate}
\item If $\varphi:\mathcal{A}\rightarrow\mathcal{B}$ and $\psi:\mathcal{B}\rightarrow\mathcal{C}$ are multihomomorphisms, 
\[\bigcup_{a\in \sigma_{\mathcal{A}}(a_{1}, \dotsc  , a_{n})}\psi\circ\varphi(a)=\bigcup_{a\in \sigma_{\mathcal{A}}(a_{1}, \dotsc  , a_{n})}\bigcup_{b\in \varphi(a)}\psi(b);\]
since 
\[\bigcup_{a\in \sigma_{\mathcal{A}}(a_{1}, \dotsc  , a_{n})}\varphi(a)\subseteq \bigcup_{(b_{1}, \dotsc  , b_{n})\in\varphi(a_{1})\times\cdots\times\varphi(a_{n})}\sigma_{\mathcal{B}}(b_{1}, \dotsc  , b_{n}),\]
we have that the last set on the equality is contained in 
\[\bigcup_{(b_{1}, \dotsc  , b_{n})\in \varphi(a_{1})\times\cdots\times\varphi(a_{n})}\bigcup_{b\in \sigma_{\mathcal{B}}(b_{1}, \dotsc  , b_{n})}\psi(b)\subseteq\bigcup_{(b_{1}, \dotsc  , b_{n})\in \varphi(a_{1})\times\cdots\times\varphi(a_{n})}\bigcup_{(c_{1}, \dotsc  , c_{n})\in \psi(b_{1})\times\cdots\times\psi(b_{n})}\sigma_{\mathcal{C}}(c_{1}, \dotsc  , c_{n});\]
by Lemma \ref{lemma about unions}, for $X_{i}=\varphi(a_{i})$ and $Y_{x}=\psi(x)$, this last set is contained in
\[\bigcup_{(c_{1}, \dotsc  , c_{n})\in \bigcup_{b_{1}\in \varphi(a_{1})}\psi(b_{1})\times\cdots\times\bigcup_{b_{n}\in \varphi(a_{n})}\psi(b_{n})}\sigma_{\mathcal{C}}(c_{1}, \dotsc  , c_{n})=\bigcup_{(c_{1}, \dotsc  , c_{n})\in \psi\circ\varphi(a_{1})\times\cdots\times\psi\circ\varphi(a_{n})}\sigma_{\mathcal{C}}(c_{1}, \dotsc  , c_{n}),\]
and from that $\psi\circ\varphi$ is also a multihomomorphism.

\item If $\varphi:\mathcal{A}\rightarrow\mathcal{B}$ and $\psi:\mathcal{B}\rightarrow\mathcal{C}$ are full multihomomorphisms, is enough to replace all occurrences of "$\subseteq$" in the proof above by "$=$" to obtain a proof that $\psi\circ\varphi$ is also a full multihomomorphism.
\end{enumerate}

Finally, it remains to be proved that such a notion of composition is associative: let $\varphi:\mathcal{A}\rightarrow\mathcal{B}$, $\psi:\mathcal{B}\rightarrow\mathcal{C}$ and $\theta:\mathcal{C}\rightarrow\mathcal{D}$ be multihomomorphisms and $a\in A$; then, using Lemma \ref{lemma about unions} again, we have that
\[[\theta\circ(\psi\circ\varphi)](a)=\bigcup_{c\in \psi\circ\varphi(a)}\theta(c)=\bigcup_{c\in \bigcup_{b\in\varphi(a)}\psi(b)}\theta(c)=\bigcup_{b\in\varphi(a)}\bigcup_{c\in\psi(b)}\theta(c)=\]
\[\bigcup_{b\in\varphi(a)}\theta\circ\psi(b)=[(\theta\circ\psi)\circ\varphi](a)\]
\end{proof}

We can arrange the categories that have so far appeared in our study of multialgebras in the following helpful diagram.

\[ \begin{tikzcd}
\textbf{MAlg}_{=}(\Sigma)\arrow[hook]{dd}\arrow{rr}{J_{=}}  & & \textbf{MMAlg}_{=}(\Sigma)\arrow[hook]{dd}\\
&&\\
\textbf{MAlg}(\Sigma) \arrow{rr}{J} & & \textbf{MMAlg}(\Sigma)
\end{tikzcd}
\]

It is clear that $\textbf{MAlg}_{=}(\Sigma)$ is a subcategory of $\textbf{MAlg}(\Sigma)$, and that $\textbf{MMAlg}_{=}(\Sigma)$ is a subcategory of $\textbf{MMAlg}(\Sigma)$, all of these four categories having as objects the class of all $\Sigma$-multialgebras. And, although $\textbf{MAlg}_{=}(\Sigma)$ is not a subcategory of $\textbf{MMAlg}_{=}(\Sigma)$, nor is $\textbf{MAlg}(\Sigma)$ a subcategory of $\textbf{MMAlg}(\Sigma)$, we can easily define functors 
\[J_{=}:\textbf{MAlg}_{=}(\Sigma)\rightarrow \textbf{MMAlg}_{=}(\Sigma)\quad\text{and}\quad J:\textbf{MAlg}(\Sigma)\rightarrow \textbf{MMAlg}(\Sigma)\]
such that the image of $J$ (respectively $J_{=}$) is a subcategory of $\textbf{MMAlg}(\Sigma)$ ($\textbf{MMAlg}_{=}(\Sigma)$) isomorphic to $\textbf{MAlg}(\Sigma)$ ($\textbf{MAlg}_{=}(\Sigma)$); we define $J$ ($J_{=}$) as the identity on objects, and for a (full) homomorphism $\varphi:\mathcal{A}\rightarrow\mathcal{B}$, the (full) multihomomorphism $J\varphi$ (respectively $J_{=}\varphi$) from $\mathcal{A}$ to $\mathcal{B}$, on an element $a$ of $\mathcal{A}$, is just $\{\varphi(a)\}$.


\section{Formulas and how to interpret them}\label{Formulas and how to interpret them}

\begin{definition}
Given a set $\mathcal{V}$, whose elements we will call propositional variables, we define the formulas\index{Formula} over the signature $\Sigma$ on the variables $\mathcal{V}$ by recursion, and only by the following rules:
\begin{enumerate}
\item all elements of $\mathcal{V}$ are formulas;
\item all elements of $\Sigma_{0}$ are formulas;
\item if $\sigma\in\Sigma_{n}$ and $\alpha_{1}, \dotsc  , \alpha_{n}$ are formulas, $\sigma(\alpha_{1}, \dotsc  , \alpha_{n})$ is a formula.
\end{enumerate}
\end{definition}

However, the expression $\sigma(\alpha_{1}, \dotsc  , \alpha_{n})$ is not entirely formally defined: in what would be the correct, formal definition, yet not as clear, a formula over the signature $\Sigma$ on the variables $\mathcal{V}$ is any function $f:I_{k}\rightarrow\mathcal{V}\cup\bigcup_{n\in\mathbb{N}}\Sigma_{n}$, where $I_{k}=\{0, \dotsc  , k\}$ is a initial segment of $\mathbb{N}$, such that:
\begin{enumerate}
\item either $k=0$ and $f(0)\in\mathcal{V}\cup\Sigma_{0}$;
\item or there exist $m\in\mathbb{N}\setminus\{0\}$, $\sigma\in\Sigma_{m}$ and formulas 
\[f_{1}:I_{k_{1}}\rightarrow \mathcal{V}\cup\bigcup_{n\in\mathbb{N}}\Sigma_{n}\quad\text{through}\quad f_{m}:I_{k_{m}}\rightarrow\mathcal{V}\cup\bigcup_{n\in\mathbb{N}}\Sigma_{n}\]
such that $k=1+\sum_{j=1}^{m}k_{j}$, $f(0)=\sigma$ and
\[\text{$f(i-k_{p}+\sum_{j=1}^{p}k_{j})=f_{p}(i)$, for $p\in\{1, \dotsc  , m\}$ and $i\in I_{k_{p}}$.}\]
\end{enumerate}

The set of all formulas over the signature $\Sigma$ on the variables $\mathcal{V}$ is denoted by $F(\Sigma, \mathcal{V})$\label{FSigmaV}. If we define, for a $\sigma\in\Sigma_{n}$, the function $\sigma_{\textbf{F}(\Sigma, \mathcal{V})}:F(\Sigma, \mathcal{V})^{n}\rightarrow F(\Sigma, \mathcal{V})$ to be, for formulas $\alpha_{1}, \dotsc  , \alpha_{n}$,
\[\sigma_{\textbf{F}(\Sigma, \mathcal{V})}(\alpha_{1}, \dotsc  , \alpha_{n})=\sigma(\alpha_{1}, \dotsc  , \alpha_{n}),\]
we obtain a $\Sigma$-algebra 
\[\textbf{F}(\Sigma, \mathcal{V})=(F(\Sigma, \mathcal{V}), \{\sigma_{\textbf{F}(\Sigma, \mathcal{V})}\}_{\sigma\in\Sigma}),\]
\label{tFSigmaV}said to be the free $\Sigma$-algebra on the variables $\mathcal{V}$, or the $\Sigma$-algebra of formulas on the variables $\mathcal{V}$.

We will also use $\textbf{F}(\Sigma, \mathcal{V})$ to denote the corresponding $\Sigma$-multialgebra.

\begin{definition}
We define the order\index{Order of a formula}\label{order}, or complexity\index{Complexity of a formula}, $|\alpha|$ of a formula $\alpha$ in $F(\Sigma, \mathcal{V})$ as:
\begin{enumerate}
\item $0$, if $\alpha$ is an element of $\mathcal{V}$;
\item $0$, if $\alpha$ is an element of $\Sigma_{0}$;
\item if $\alpha$ is of the form $\sigma(\alpha_{1}, \dotsc  , \alpha_{n})$, 
\[|\alpha|=1+\max_{1\leq i\leq n}|\alpha_{i}|.\]
\end{enumerate}
\end{definition}

Now, we have given the simplest definition of formula: when we approach choice-dependent freely generated multialgebras, we will find what one can identify as a broader definition of formula, but with which we shall not deal with until later. 


\subsection{Different notions of interpretation}

Perhaps more important than defining what is a formula is interpreting this formula, and here the non-determinism of multialgebras once again gives us an array of different notions.

\begin{definition}
Given a $\Sigma$-multialgebra $\mathcal{A}$, a homomorphism $\nu:\textbf{F}(\Sigma, \mathcal{V})\rightarrow \mathcal{A}$ is called a legal valuation.\index{Valuation, Legal}
\end{definition}

The notion of legal valuation is one attempt to give an interpretation of a formula $\alpha$: in this case, the formula $\alpha$, under the legal valuation $\nu$, takes the value $\nu(\alpha)$.

A serious problem that arises from the notion of legal valuation is that one interpretation of the variables can lead to several legal valuations, and therefore different interpretations of a formula. This does not occur to $\Sigma$-algebras, where one interpretation of the variables implies an unique interpretation of the formulas.

\begin{example}
Take the signature $\Sigma$ with $\Sigma_{1}=\{s\}$ and $\Sigma_{n}=\emptyset$ for $n\neq 1$; take the $\Sigma$-multialgebras $\textbf{F}(\Sigma,  \mathcal{V})$, for $\mathcal{V}=\{x\}$, and $\mathcal{A}=(A, \{\sigma_{\mathcal{A}}\}_{\sigma\in\Sigma})$ such that $A=\{0,1\}$ and $s_{\mathcal{A}}(0)=s_{\mathcal{A}}(1)=\{0,1\}$.
\begin{figure}[H]
\centering
\begin{minipage}[t]{4cm}
\centering
\begin{tikzcd}[row sep=1cm]
\vdots\\
s^{2}(x)\arrow{u}{s_{\textbf{F}(\Sigma, \mathcal{V})}}\\
s(x)\arrow{u}{s_{\textbf{F}(\Sigma, \mathcal{V})}}\\
x\arrow{u}{s_{\textbf{F}(\Sigma, \mathcal{V})}}
\end{tikzcd}
\caption*{$\textbf{F}(\Sigma, \mathcal{V})$}
\end{minipage}
\hspace{3cm}
\centering
\begin{minipage}[t]{4cm}
\centering
\begin{tikzcd}[column sep=1cm]
0\arrow[loop left, out=150, in=-150, distance=3em, swap]{}{s_{\mathcal{A}}}\arrow[shift right, swap]{r}{s_{\mathcal{A}}} & 1\arrow[loop right, out=30, in=-30, distance=3em]{}{s_{\mathcal{A}}}\arrow[shift right, swap]{l}{s_{\mathcal{A}}}
\end{tikzcd}
\caption*{$\mathcal{A}$}
\end{minipage}
\end{figure}
We state that $\nu_{1}:F(\Sigma, \mathcal{V})\rightarrow A$, given by $\nu_{1}(\alpha)=0$ for every $\alpha\in F(\Sigma, \mathcal{V})$, and $\nu_{2}:F(\Sigma, \mathcal{V})\rightarrow A$, given by $\nu_{2}(\alpha)=1$ for every $\alpha\in F(\Sigma, \mathcal{V})\setminus\{x\}$ and $\nu_{2}(x)=0$, are homomorphisms. In fact, for every $\alpha\in F(\Sigma, \mathcal{V})$, we have that
\[\{\nu_{i}(\beta)\ :\  \beta\in s_{\textbf{F}(\Sigma, \mathcal{V})}(\alpha)\}=\{\nu_{i}(s(\alpha))\}\subseteq \{0,1\}=s_{\mathcal{A}}(\nu_{i}(\alpha)),\]
where $i\in \{1,2\}$. So $\nu_{1}$ and $\nu_{2}$ are different legal valuations, despite being the same over the variables.
\end{example}

We will say that a legal valuation $\nu:\textbf{F}(\Sigma, \mathcal{V})\rightarrow\mathcal{A}$ is associated\index{Valuation, Associated legal} to the function $\chi:\mathcal{V}\rightarrow A$ given by $\chi=\nu|_{\mathcal{V}}$,\footnote{Here, it is important to clarify the notation: for a function $f:X\rightarrow Y$ and a set $Z\subset X$, we will denote by $f|_{Z}$ the restriction of $f$ to $Z$.} which we shall call a valuation or interpretation of the variables; more generally, we say $\nu$ is associated to the function $\chi:\mathcal{V}\rightarrow\mathcal{P}(A)\setminus\{\emptyset\}$ if $\nu(x)\in \chi(x)$ for every $x\in\mathcal{V}$. What we have shown in the previous example is that two distinct legal valuations may be associated to the same interpretation of the variables.

So, if legal valuations are problematic, what other definition should we take? Many attempts to correctly interpret a formula in a multialgebra have been made, each with its own advantages and drawbacks. When studying choice-dependent freely generated multialgebras, we will show that maybe one needs to specify both an interpretation of the variables and what we shall call a collection of choices. For now, we shall offer a few other possibilities.

\begin{definition}
Given a $\Sigma$-multialgebra $\mathcal{A}$, a full multihomomorphism $\nu:\textbf{F}(\Sigma, \mathcal{V})\rightarrow\mathcal{A}$ is called a full valuation\index{Valuation, Full}.
\end{definition}

\begin{proposition}
Two full valuations $\nu_{1}, \nu_{2}:\textbf{F}(\Sigma, \mathcal{V})\rightarrow\mathcal{A}$ such that $\nu_{1}|_{\mathcal{V}}=\nu_{2}|_{\mathcal{V}}$ are the same.
\end{proposition}

\begin{proof}
We will prove that, for any formula $\alpha$, $\nu_{1}(\alpha)=\nu_{2}(\alpha)$ by induction on the order of $\alpha$. If $|\alpha|=0$, either $\alpha$ is an element $x\in\mathcal{V}$ and then $\nu_{1}(x)=\nu_{2}(x)$ since $\nu_{1}|_{\mathcal{V}}=\nu_{2}|_{\mathcal{V}}$, or $\alpha$ is a $\sigma\in\Sigma_{0}$, and then 
\[\nu_{1}(\alpha)=\sigma_{\mathcal{A}}=\nu_{2}(\alpha).\]

Assuming the proposition is true for formulas of order at most $m$, if $\alpha=\sigma(\alpha_{1}, \dotsc  , \alpha_{n})$ is a formula of order $m+1$ then $\alpha_{1}$ through $\alpha_{n}$ are formulas of order at most $m$, and therefore $\nu_{1}(\alpha_{1})=\nu_{2}(\alpha_{2}), \dotsc  , \nu_{1}(\alpha_{n})=\nu_{2}(\alpha_{n})$. It follows that 
\[\nu_{1}(\alpha)=\bigcup_{(a_{1}, \dotsc  , a_{n})\in \nu_{1}(\alpha_{1})\times\cdots\nu_{1}(\alpha_{n})}\sigma_{\mathcal{A}}(a_{1}, \dotsc  , a_{n})=\bigcup_{(a_{1}, \dotsc  , a_{n})\in \nu_{2}(\alpha_{1})\times\cdots\nu_{2}(\alpha_{n})}\sigma_{\mathcal{A}}(a_{1}, \dotsc  , a_{n})=\nu_{2}(\alpha).\]
\end{proof}

So, to every function $\chi:\mathcal{V}\rightarrow \mathcal{P}(A)\setminus\{\emptyset\}$, there corresponds at most one full valuation: but we can also show that there is at least one full valuation associated to $\chi$, which we denote by $\overline{\chi}$. We define it by induction on the order of a formula $\alpha$:
\begin{enumerate}
\item if $|\alpha|=0$, either $\alpha$ is a $x\in\mathcal{V}$, when we define $\overline{\chi}(\alpha)=\chi(x)$, or $\alpha$ is a $\sigma\in\Sigma_{0}$, when $\overline{\chi}(\alpha)=\sigma_{\mathcal{A}}$;
\item supposing $\overline{\chi}$ is defined for formulas of degree at most $m$, if $\alpha=\sigma(\alpha_{1}, \dotsc  , \alpha_{n})$ is of degree $m+1$, and so $\alpha_{1}, \dotsc  , \alpha_{n}$ are of degree at most $m$, we define
\[\overline{\chi}(\alpha)=\bigcup_{(a_{1}, \dotsc  , a_{n})\in\overline{\chi}(\alpha_{1})\times\cdots\times\overline{\chi}(\alpha_{n})}\sigma_{\mathcal{A}}(a_{1}, \dotsc  , a_{n}).\]
\end{enumerate}

We can generalize both legal and full valuations through what we will call simply valuations: any multihomomorphism $\nu:\textbf{F}(\Sigma, \mathcal{V})\rightarrow \mathcal{A}$ will be said to be a valuation on $\mathcal{A}$, associated to any of the functions $\chi:\mathcal{V}\rightarrow \mathcal{P}(A)\setminus\{\emptyset\}$ such that $\nu(x)\subseteq\chi(x)$ for every $x\in\mathcal{V}$.

It is quite clear how this is a generalization of a full valuation, but one can see this is not a direct generalization of a legal valuation, since those are function from $F(\Sigma, \mathcal{V})$ to $A$, instead of $\mathcal{P}(A)\setminus\{\emptyset\}$; however, when one considers for the legal valuation $\nu$ the function $\nu^{*}:F(\Sigma, \mathcal{V})\rightarrow \mathcal{P}(A)\setminus\{\emptyset\}$ such that $\nu^{*}(\alpha)=\{\nu(\alpha)\}$ for every $\alpha\in F(\Sigma, \mathcal{V})$, one sees that $\nu^{*}$ is a valuation.

Without much further ado, we will identify $\nu^{*}$ with $\nu$, and characterize some valuations as legal ones.

\begin{proposition}
Fixed $\textbf{F}(\Sigma,\mathcal{V})$, the set $V(\mathcal{A})$ of all valuations from it to the $\Sigma$-multialgebra $\mathcal{A}$ is a partially ordered set, when for valuations $\nu_{1}, \nu_{2}\in V(\mathcal{A})$ one considers $\nu_{1}\leq \nu_{2}$ if and only if 
\[\nu_{1}(\alpha)\subseteq \nu_{2}(\alpha),\quad \forall \alpha\in F(\Sigma, \mathcal{V}),\]
closed under suprema of non-empty sets, with as minimal elements the set of all legal valuations $LV(\mathcal{A})$  on $\mathcal{A}$.
\end{proposition}

\begin{proof}
First, we prove that the mentioned order is in fact an order.
\begin{enumerate}
\item For any $\nu\in V(\mathcal{A})$ and all $\alpha\in F(\Sigma, \mathcal{V})$, it is true that $\nu(\alpha)\subseteq\nu(\alpha)$, and therefore $\nu\leq \nu$.
\item For any two $\nu_{1}, \nu_{2}\in V(\mathcal{A})$, if $\nu_{1}\leq\nu_{2}$ and $\nu_{2}\leq \nu_{1}$, then for every $\alpha\in F(\Sigma, \mathcal{V})$ we have that $\nu_{1}(\alpha)\subseteq\nu_{2}(\alpha)$ and $\nu_{2}(\alpha)\subseteq\nu_{1}(\alpha)$, and therefore $\nu_{1}(\alpha)=\nu_{2}(\alpha)$; it follows that $\nu_{1}=\nu_{2}$.
\item For any three $\nu_{1}, \nu_{2}, \nu_{3}\in V(\mathcal{A})$, if $\nu_{1}\leq\nu_{2}$ and $\nu_{2}\leq\nu_{3}$, then for every $\alpha\in F(\Sigma, \mathcal{V})$ we see that $\nu_{1}(\alpha)\subseteq\nu_{2}(\alpha)$ and $\nu_{2}(\alpha)\subseteq\nu_{3}(\alpha)$, and therefore $\nu_{1}(\alpha)\subseteq\nu_{3}(\alpha)$; it follows that $\nu_{1}\leq\nu_{3}$.
\end{enumerate}

Now, we state that the supremum of a non-empty set $\Lambda\subseteq V(\mathcal{A})$ is the valuation such that, for every $\alpha\in F(\Sigma, \mathcal{V})$, 
\[\nu(\alpha)=\bigcup_{\lambda\in \Lambda}\lambda(\alpha);\]
it is, in fact, a valuation since it is a multihomomorphism: if $\sigma\in\Sigma_{n}$ and $\alpha_{1}, \dotsc  , \alpha_{n}$ are formulas in $F(\Sigma, \mathcal{V})$, we have that
\[\nu(\sigma(\alpha_{1}, \dotsc  , \alpha_{n}))=\bigcup_{\lambda\in\Lambda}\lambda(\sigma(\alpha_{1}, \dotsc  , \alpha_{n}))\subseteq\bigcup_{\lambda\in\Lambda}\bigcup_{(a_{1}, \dotsc  , a_{n})\in \lambda(\alpha_{1})\times\cdots\times\lambda(\alpha_{n})}\sigma_{\mathcal{A}}(a_{1}, \dotsc  , a_{n})\subseteq\]
\[\bigcup_{\lambda\in\Lambda}\bigcup_{(a_{1}, \dotsc  , a_{n})\in \nu(\alpha_{1})\times\cdots\times\nu(\alpha_{n})}\sigma_{\mathcal{A}}(a_{1}, \dotsc  , a_{n})=\bigcup_{(a_{1}, \dotsc  , a_{n})\in \nu(\alpha_{1})\times\cdots\times\nu(\alpha_{n})}\sigma_{\mathcal{A}}(a_{1}, \dotsc  , a_{n}).\]

Clearly $\nu$ is an upper bound of $\Lambda$, so suppose $\nu^{*}$ is another upper bound: for every formula $\alpha$ and $\lambda\in \Lambda$ we have that $\lambda(\alpha)\subseteq \nu^{*}(\alpha)$, and from that $\nu(\alpha)=\bigcup_{\lambda\in\Lambda}\lambda(\alpha)\subseteq \nu^{*}(\alpha)$, that is, $\nu\leq\nu^{*}$ and so $\nu=\sup\Lambda$.

Since for any legal valuation $\nu\in LV(\mathcal{A})$ and formula $\alpha$ the set $\nu(\alpha)$ is always of cardinality $1$, legal valuations are clearly minimal elements.
\end{proof}

For a function $\chi:\mathcal{V}\rightarrow\mathcal{P}(A)\setminus\{\emptyset\}$, we would like to consider the set 
\[V(\chi)=\{\nu\in V(\mathcal{A})\ :\  \nu\leq\overline{\chi}\}.\]
It is clearly a partially ordered set, containing $\overline{\chi}$ and all valuations, legal or not, associated to $\chi$; it has as minimal elements $LV(\chi)=LV(\mathcal{A})\cap V(\chi)$ and maximum $\overline{\chi}$; it is closed under suprema since for any $\Lambda\subseteq V(\chi)$ we have that $\overline{\chi}$ is an upper bound for $\Lambda$, and therefore $\sup\Lambda\leq \overline{\chi}$, meaning that $\sup\Lambda\in V(\chi)$.

Notice that, alternatively, $\nu\in V(\chi)$ if, and only if, $\nu$ is a valuation for which, for every $x\in\mathcal{V}$, $\nu(x)\subseteq \chi(x)$.

The join-semilattice $V(\chi)$ is then a very interesting object: if one is in doubt of how to interpret a formula given the interpretation of variables $\chi$, using valuations, legal valuations or full valuations, $V(\chi)$ captures the information offered by all of these attempts, despite being considerably more complex.


\subsection{Full valuations and adjointedness}

Take the class of all sets and, as morphisms between any two sets $X$ and $Y$, take the functions $f:X\rightarrow\mathcal{P}(Y)\setminus\{\emptyset\}$, with the composition of two such functions $f:X\rightarrow\mathcal{P}(Y)\setminus\{\emptyset\}$ and $g:Y\rightarrow\mathcal{P}(Z)\setminus\{\emptyset\}$ being defined as the function $g\circ f:X\rightarrow\mathcal{P}(Z)\setminus\{\emptyset\}$ such that, for any $x\in X$, 
\[g\circ f(x)=\bigcup_{y\in f(x)}g(y).\]

As we saw before, such a composition is associative, and has as identity, for a given set $X$, the function $Id_{X}:X\rightarrow\mathcal{P}(X)\setminus\{\emptyset\}$ such that, for every $x\in X$, $Id_{X}(x)=\{x\}$. Therefore, this object is a category, which we will call the category of sets with multifunctions, and denote by $\textbf{MSet}$\label{MSet}; its morphisms will be called multifunctions.

Now, fixed a signature $\Sigma$, take for every set $X$ the $\Sigma$-multialgebra $FX=\textbf{F}(\Sigma, X)$ and, for every morphism $f$ from $X$ to $Y$, that is, function $f:X\rightarrow \mathcal{P}(Y)\setminus\{\emptyset\}$, take the full multihomomorphism $Ff=\overline{f}:\textbf{F}(\Sigma, X)\rightarrow \textbf{F}(\Sigma, Y)$.

Given the identity $Id_{X}:X\rightarrow X$ on $\textbf{MSet}$, we state that $\overline{Id_{X}}:\textbf{F}(\Sigma, X)\rightarrow \textbf{F}(\Sigma, X)$ is exactly the identity $Id_{\textbf{F}(\Sigma, X)}$ of this multialgebra in $\textbf{MMAlg}_{=}(\Sigma)$: one can see this by induction over the order of a formula $\alpha$.
\begin{enumerate}
\item If $\alpha$ is of order $0$, either $\alpha$ is an $x\in X$ and then $\overline{Id_{X}}(x)=Id_{X}(x)=\{x\}$, or $\alpha$ is a $\sigma\in\Sigma_{0}$, when $\overline{Id_{X}}(\sigma)=\sigma_{\textbf{F}(\Sigma, X)}=\{\sigma\}$.
\item If $\alpha=\sigma(\alpha_{1}, \dotsc  , \alpha_{n})$, for $\sigma\in\Sigma_{n}$ and $\alpha_{1}$ through $\alpha_{n}$ of order less than that of $\alpha$, we have that
\[\overline{Id_{X}}(\alpha)=\bigcup_{(\beta_{1}, \dotsc  , \beta_{n})\in \overline{Id_{X}}(\alpha_{1})\times\cdots\times\overline{Id_{X}}(\alpha_{n})}\sigma_{\textbf{F}(\Sigma, X)}(\beta_{1}, \dotsc  , \beta_{n})=\bigcup_{(\beta_{1}, \dotsc  , \beta_{n})\in \{\alpha_{1}\}\times\cdots\times\{\alpha_{n}\}}\{\sigma(\beta_{1}, \dotsc  , \beta_{n})\}=\]
\[\{\sigma(\alpha_{1}, \dotsc  , \alpha_{n})\}=\{\alpha\}.\]
\end{enumerate}

Furthermore, given morphisms $f$ from $X$ to $Y$ and $g$ from $Y$ to $Z$ in $\textbf{MSet}$, we see that for every $x\in X$ 
\[\overline{g}\circ\overline{f}(x)=\bigcup_{y\in \overline{f}(x)}\overline{g}(y)=\bigcup_{y\in f(x)}\overline{g}(y),\]
and since $f(x)\subseteq Y$, for any $y\in f(x)$ we have that $\overline{g}(y)=g(y)$ and therefore
\[\bigcup_{y\in f(x)}\overline{g}(y)=\bigcup_{y\in f(x)}g(y)=g\circ f(x)=\overline{g\circ f}(x),\]
since $g\circ f$ goes from $X$ to $Z$; therefore both $\overline{g}\circ\overline{f}$ and $\overline{g\circ f}$ are full multihomomorphisms equal over the variables, and are therefore equal. 

Since $FId_{X}=Id_{FX}$ and $Fg\circ Ff=F(g\circ f)$, we have proved
\[F:\textbf{MSet}\rightarrow\textbf{MMAlg}_{=}(\Sigma)\]
is a functor. Now, we can also consider the forgetful functor
\[\mathcal{U}:\textbf{MMAlg}_{=}(\Sigma)\rightarrow \textbf{MSet},\]
since full multihomomorphisms are multifunctions; what we shall prove now is that $F$ and $\mathcal{U}$ are actually adjoint. So, we consider the functions, for a set $X$ and a $\Sigma$-multialgebra $\mathcal{A}$,
\[\Phi_{\mathcal{A}, X}:Hom_{\textbf{MSet}}(X, \mathcal{U}\mathcal{A})\rightarrow Hom_{\textbf{MMAlg}_{=}(\Sigma)}(FX, \mathcal{A})\]
associating a map $f:X\rightarrow\mathcal{P}(A)\setminus\{\emptyset\}$ to the full multihomomorphism $\overline{f}:\textbf{F}(\Sigma, X)\rightarrow\mathcal{A}$ extending $f$. Since, for every $f$, there corresponds one and only one $\overline{f}$, the $\Phi_{\mathcal{A}, X}$ are bijections.

Now, given sets $X$ and $Y$, $\Sigma$-multialgebras $\mathcal{A}$ and $\mathcal{B}$, a morphism $f$ from $Y$ to $X$ in $\textbf{MSet}$ and a full multihomomorphism $\varphi:\mathcal{A}\rightarrow\mathcal{B}$, we must only prove that the following diagram commutes.

\[ \begin{tikzcd}[row sep=5em, column sep=3em]
Hom_{\textbf{MSet}}(X, \mathcal{U}\mathcal{A}) \arrow{r}{\Phi_{\mathcal{A}, X}} \arrow{d}{Hom(f, \mathcal{U}\varphi)}& Hom_{\textbf{MMAlg}_{=}(\Sigma)}(FX, \mathcal{A}) \arrow{d}{Hom(Ff, \varphi)} \\%
Hom_{\textbf{MSet}}(Y, \mathcal{U}\mathcal{B}) \arrow{r}{\Phi_{\mathcal{B}, Y}}& Hom_{\textbf{MMAlg}_{=}(\Sigma)}(FY, \mathcal{B})
\end{tikzcd}
\]

So, denoting the universe of $\mathcal{A}$ by $A$ as usual, we take a function $g:X\rightarrow\mathcal{P}(A)\setminus\{\emptyset\}$ in $Hom_{\textbf{MSet}}(X, \mathcal{U}\mathcal{A})$: on the top edge of the diagram we obtain the full multihomomorphism $\Phi_{\mathcal{A}, X}(g)=\overline{g}$ from $\textbf{F}(\Sigma, X)$ to $\mathcal{A}$; and on the right edge we obtain the once again full multihomomorphism $\varphi\circ \overline{g}\circ\overline{f}$, since $Ff=\overline{f}$, from $\textbf{F}(\Sigma, Y)$ to $\mathcal{B}$.

On the left edge, we obtain the multifunction $\varphi\circ g\circ f$ from $Y$ to the universe $B$ of $\mathcal{B}$, since $\mathcal{U}\varphi=\varphi$; applying $\Phi_{\mathcal{B}, Y}$ on the bottom edge of the diagram we obtain $\overline{\varphi\circ g\circ f}$, also a full multihomomorphism from $\textbf{F}(\Sigma, Y)$ to $\mathcal{B}$.

Now, to prove that the diagram commutes or, what is equivalent, that $\varphi\circ\overline{g}\circ\overline{f}=\overline{\varphi\circ g\circ f}$, since both are full valuations on $\textbf{F}(\Sigma, Y)$ is enough to prove that, when restricted to $Y$, both are the same. So, for $y\in Y$,
\[\varphi\circ\overline{g}\circ\overline{f}(y)=\bigcup_{x\in \overline{f}(x)}\varphi\circ\overline{g}(x)=\bigcup_{x\in\overline{f}(x)}\bigcup_{a\in \overline{g}(x)}\varphi(a)=\bigcup_{x\in f(x)}\bigcup_{a\in\overline{g}(x)}\varphi(a)=\bigcup_{x\in f(y)}\bigcup_{a\in g(x)}\varphi(a)=\]
\[\bigcup_{x\in f(y)}\varphi\circ g(x)=\varphi\circ g\circ f(y)=\overline{\varphi\circ g\circ f}(y),\]
what ends the proof that $F$ and $\mathcal{U}$ are adjoint.

\begin{proposition}
Let $\chi:\mathcal{V}\rightarrow\mathcal{P}(A)\setminus\{\emptyset\}$ be a function and $\varphi:\textbf{F}(\Sigma, \mathcal{V})\rightarrow\mathcal{A}$ a multihomomorphism such that, for every $x\in\mathcal{V}$, $\varphi(x)\subseteq\chi(x)$: then $\varphi\leq \overline{\chi}$.
\end{proposition}

\begin{proof}
It is enough to proceed by induction over the order of a formula $\alpha$.

If it is $0$, either $\alpha$ is an element $x\in\mathcal{V}$, when the result is true by hypothesis; or $\alpha$ is a $\sigma\in\Sigma_{0}$, and since $\varphi$ is a multihomomorphism we have that $\varphi(\sigma)\subseteq \sigma_{\mathcal{A}}=\overline{\chi}(\sigma)$.

Now, suppose the result is true for formulas of order smaller than that of $\alpha$, and suppose $\alpha$ is of the form $\sigma(\alpha_{1}, \dotsc  , \alpha_{n})$: then 
\[\varphi(\alpha)=\bigcup_{\beta\in\sigma_{\textbf{F}(\Sigma, \mathcal{V})}(\alpha_{1}, \dotsc  , \alpha_{n})}\varphi(\beta)\subseteq\bigcup_{(\beta_{1}, \dotsc  , \beta_{n})\in \varphi(\alpha_{1})\times\cdots\times\varphi(\alpha_{n})}\sigma_{\mathcal{A}}(\beta_{1}, \dotsc  , \beta_{n})\subseteq\]
\[\bigcup_{(\beta_{1}, \dotsc  , \beta_{n})\in\overline{\chi}(\alpha_{1})\times\cdots\times\overline{\chi}(\alpha_{n})}\sigma_{\mathcal{A}}(\beta_{1}, \dotsc  , \beta_{n})=\bigcup_{\beta\in\sigma_{\textbf{F}(\Sigma, \mathcal{V})}(\alpha_{1}, \dotsc  , \alpha_{n})}\overline{\chi}(\beta)=\overline{\chi}(\alpha).\]
\end{proof}

This way, we see that full valuations have many important properties: they allow, together with the algebras of formulas, to build a functor adjoint to a forgetful functor on a rather natural category ($\textbf{MSet}$); furthermore, they are maximal among valuations associated to a given evaluation of the variables.

And, at the same time, they are not very useful to most applications in logic, one reason being that, from the point of view of full valuations, the free objects of $\textbf{MMAlg}_{=}(\Sigma)$ are the algebras of formulas $\textbf{F}(\Sigma, \mathcal{V})$, meaning we gain no new structures to analyze. A second reason for or disinterest on full valuations is that they are not at all very precise: legal valuations, that are only functions instead of multifunctions, are much more useful and "precise". Unfortunately, to an evaluation of the variables there are many legal valuations associated: one approach in the direction of fixing this problem is to consider collections of choices (Definition \ref{Collection of Choices}), which will allow us at the same time to consider generalizations of the algebras of formulas (Definition \ref{Multialgebra of formulas}).

We will work those objects in Chapter \ref{Chapter2}. Their algebraic structure is rather rich, regardless of how difficult their treatment is trough category theory. And, speaking of category theory, the categories shown up until now can be arranged as in the diagram below. The functor $J_{\textbf{Set}}$ is the identity on sets and takes functions to multifunctions in the obvious way, like $J$ and $J_{=}$; the question of whether the forgetful functor from $\textbf{MAlg}_{=}(\Sigma)$ to $\textbf{Set}$, here represented by a dashed arrow, has a left adjoint is answered negatively at the end of Section \ref{Freely generated multialgebra}.

\[ \begin{tikzcd}
\textbf{MAlg}(\Sigma) \arrow{dd}{J} && \textbf{MAlg}_{=}(\Sigma)\arrow[hook]{ll}\arrow[dashed]{rr} \arrow{dd}{J_{=}} &&\textbf{Set}\arrow{dd}{J_{\textbf{Set}}}  \\
&&&&\\
 \textbf{MMAlg}(\Sigma) && \textbf{MMAlg}_{=}(\Sigma)\arrow[hook]{ll}\arrow[transform canvas={yshift=.7ex}]{rr}{\mathcal{U}} && \textbf{MSet}\arrow[transform canvas={yshift=-.7ex}]{ll}{F}
\end{tikzcd}
\]

\section{Lattices, and Boolean and Heyting algebras}\label{Lattices, and Boolean... }

Consider the signature $\Sigma_{\textbf{Lat}}$\label{SigmaLat} with two binary connectives, $\vee$ and $\wedge$: that is $(\Sigma_{\textbf{Lat}})_{2}=\{\vee, \wedge\}$ and $(\Sigma_{\textbf{Lat}})_{n}=\emptyset$, for $n\neq 2$. 

\begin{definition}
A lattice\index{Lattice} is a $\Sigma_{\textbf{Lat}}$-algebra $\mathcal{L}=(L,\{\sigma_{\mathcal{L}}\}_{\sigma\in\Sigma})$ satisfying, for every $x,y,z\in L$,
\begin{enumerate}
\item $x\vee y=y\vee x$ and $x\wedge y=y\wedge x$ (commutative laws);
\item $x\vee(y\vee z)=(x\vee y)\vee z$ and $x\wedge(y\wedge z)=(x\wedge y)\wedge z$ (associative laws);
\item $x\vee x=x$ and $x\wedge x=x$ (idempotent laws);
\item $x\vee(x\wedge y)=x$ and $x\wedge(x\vee y)=x$ (absorption laws).
\end{enumerate}
\end{definition}
In this definition, for simplicity, we have denoted $\vee_{\mathcal{L}}(x,y)$ and $\wedge_{\mathcal{L}}(x,y)$ by, respectively, $x\vee y$ and $x\wedge y$.

This is often regarded as the algebraic definition of a lattice, an order-theoretic one existing as well. A partial order in a set $X$ is a binary relation $\leq$ on $X$, meaning a subset of $X\times X$, such that, for all $x, y, z\in X$,
\begin{enumerate}
\item $x\leq x$ (reflexivity);
\item $x\leq y$ and $y\leq x$ imply $x=y$ (antisymmetry);
\item $x\leq y$ and $y\leq z$ imply $x\leq z$ (transitivity),
\end{enumerate}
where we will denote the fact that $(x,y)\in\leq$ by $x\leq y$. A set $X$, together with a partial order $\leq$ in it, will be named a partially ordered set, or poset. 

In a poset $(X, \leq)$, we say $x$ is an upper bound of $S\subseteq X$ if $y\leq x$ for every $y\in S$; $x$ is a lower bound for $S$ if $x\leq y$ for every $y\in S$. Then, the supremum of $S\subseteq X$, if it exists, is its least upper bound, meaning an upper bound $\sup S\in X$ for $S$ such that, for any other upper bound $x$ for $S$, $\sup S\leq x$; the infimum of $S$, if it exists, is its greatest lower bound, meaning a lower bound $\inf S$ for $S$ such that, for any other lower bound $x$ for $S$, $x\leq \inf S$.

\begin{definition}
A lattice is a poset $\mathcal{L}=(L,\leq)$ such that, for any elements $x, y\in L$, $\sup\{x,y\}$ and $\inf\{x,y\}$ both exist.
\end{definition}

One can translate between the two definitions: given a lattice $\mathcal{L}$ presented as a $\Sigma_{\textbf{Lat}}$-algebra $(L,\{\sigma_{\mathcal{L}}\}_{\sigma\in\Sigma_{\textbf{Lat}}})$, we may define an order $\leq$ in $L$ by 
\[x\leq y\quad\text{if and only if}\quad x=x\wedge y,\]
and it is easy to prove that in this case $(L,\leq)$ is a lattice, where $x\wedge y$ and $x\vee y$ are, respectively, the infimum and supremum of $x$ and $y$. Reciprocally, given a lattice $\mathcal{L}$ presented as a poset $(L,\leq)$, we may define operations $\vee$ and $\wedge$ on $L$ by
\[x\wedge y=\inf\{x,y\}\quad\text{and}\quad x\vee y=\sup\{x,y\},\]
from what we obtain a $\Sigma_{\textbf{Lat}}$-algebra on which $x\leq y$ if, and only if, $x=x\wedge y$.

An element $0$ of a lattice $\mathcal{L}$ is a minimum, or bottom element, when either: $x\vee 0=x$, for every $x\in L$; $x\wedge 0=0$, for every $x\in L$; or $0\leq x$, for every $x\in X$. Notice that all three conditions are equivalent. An element $1$ of $\mathcal{L}$ is a maximum, or top element, when either: $x\vee 1=1$, for every $x\in L$; $x\wedge 1=x$, for every $x\in L$; or $x\leq 1$, for every $x\in L$. Again, all three conditions are equivalent. 

To accommodate bottom and top elements, we extend the signature $\Sigma_{\textbf{Lat}}$ to $\Sigma_{\textbf{Lat}}^{0}$\label{SigmaLat0}, $\Sigma_{\textbf{Lat}}^{1}$\label{SigmaLat1} and $\Sigma_{\textbf{Lat}}^{0,1}$\label{SigmaLat01} by adding, respectively, a $0$-ary symbol $0$, a $0$-ary symbol $1$ and $0$-ary symbols $0$ and $1$. For simplicity, we will drop the indexes $\mathcal{L}$ from $0_{\mathcal{L}}$ and $1_{\mathcal{L}}$.

\begin{definition}
\begin{enumerate}
\item A lattice with bottom (respectively top) element is a $\Sigma_{\textbf{Lat}}^{0}$-algebra ($\Sigma_{\textbf{Lat}}^{1}$-algebra) $\mathcal{L}=(L,\{\sigma_{\mathcal{L}}\}_{\sigma\in\Sigma_{\textbf{Lat}}^{0}})$ such that $(L,\{\sigma_{\mathcal{L}}\}_{\sigma\in\Sigma_{\textbf{Lat}}})$ is a lattice and $0$ is a bottom element ($1$ is a top element).
\item A bounded lattice\index{Lattice, Bounded} is a $\Sigma_{\textbf{Lat}}^{0,1}$-algebra $\mathcal{L}=(L,\{\sigma_{\mathcal{L}}\}_{\sigma\in\Sigma_{\textbf{Lat}}^{0,1}})$ such that $(L,\{\sigma_{\mathcal{L}}\}_{\sigma\in\Sigma_{\textbf{Lat}}})$ is a lattice, $0$ is a bottom element and $1$ is a top element.
\end{enumerate}
\end{definition}

In a lattice, "$x$ implies $y$" is defined as the element
\[\sup\{z\in L: x\wedge z\leq y\},\]
if it exists, and denoted by $x\rightarrow y$. We add to the signatures $\Sigma_{\textbf{Lat}}^{1}$ and $\Sigma_{\textbf{Lat}}^{0,1}$ the binary symbol $\rightarrow$, obtaining respectively the signatures $\Sigma_{\textbf{Imp}}$\label{SigmaImp} and $\Sigma_{\textbf{Hey}}$\label{SigmaHey}.

\begin{definition}
\begin{enumerate}
\item An implicative lattice\index{Lattice, Implicative} is a $\Sigma_{\textbf{Imp}}$-algebra $\mathcal{L}=(L,\{\sigma_{\mathcal{L}}\}_{\sigma\in\Sigma_{\textbf{Imp}}})$ such that $(L, \{\sigma_{\mathcal{L}}\}_{\sigma\in\Sigma_{\textbf{Lat}}^{1}})$ is a lattice with top element, $\sup\{z\in L: x\wedge z\leq y\}$ exists for any two $x,y\in L$ and 
\[x\rightarrow y=\sup\{z\in L: x\wedge z\leq y\}.\]

\item A Heyting algebra\index{Heyting algebra} is a $\Sigma_{\textbf{Hey}}$-algebra such that $(L,\{\sigma_{\mathcal{L}}\}_{\sigma\in\Sigma_{\textbf{Imp}}})$ is an implicative lattice and $0$ is a bottom element.
\end{enumerate}
\end{definition}

It is important here to point out that various different definitions of implicative lattices and Heyting algebras exist in the literature, most, if not all, equivalent to each other; we chose definitions most befitting of our purposes, taken from \cite{ParLog}.

Another relevant point to make is that Heyting algebras are no strange beasts: classical propositional logic is to Boolean algebras, which we shall soon formally define, as intuitionistic logic is to Heyting algebras, meaning that not only are Heyting algebras models of intuitionistic logic, but their class characterizes intuitionistic logic (unfortunately, the parallel ends there: although the two-valued Boolean algebra is, by itself, capable of characterizing $\textbf{CPL}$, no finite algebra, Heyting or not, can do the same to intuitionistic logic, see \cite{Godel}). One could argue that while intuitionistic logic possesses a negation, Heyting algebras do not, but this can be easily dealt with: it is enough to define
\[\neg x=x\rightarrow 0.\]

With this, we can define Heyting algebras over the signature $\Sigma_{\textbf{Boo}}$\label{SigmaBoo}, obtained from $\Sigma_{\textbf{Hey}}$ by addition of an unary symbol "$\neg$", by merely requesting that $\neg x=x\rightarrow 0$ for all elements $x$; we say $\neg x$ is the partial complement, or (intuitionistic) negation of $x$.

\begin{definition}
A Boolean algebra\index{Boolean algebra} is a Heyting algebra over $\Sigma_{\textbf{Boo}}$ such that, for every $x$ in its universe, we have the law of excluded middle, or \textit{tertium non datur}:
\[x\vee\neg x=1.\]
\end{definition}

The following lemma will need an observation about the definition of implication on a Heyting algebra: notice that, if $x\wedge z\leq y$, then $z\leq x\rightarrow y$, by the very definition of $x\rightarrow y$ as the supremum of $\{z\in L: x\wedge z\leq y\}$. Reciprocally, if $z\leq x\rightarrow y$, $x\wedge z\leq x\wedge(x\rightarrow y)\leq y$,\footnote{Here we are using that $x\leq y$ implies $x\wedge z\leq y\wedge z$, but this is trivial to prove: if $x\leq y$, $x=x\wedge y$, and so $(x\wedge z)\wedge(y\wedge z)=(x\wedge y)\wedge z=x\wedge z$.} meaning that $z\leq x\rightarrow y$ if, and only if, $x\wedge z\leq y$. The following proof, although quite standard, follows closely the one found in \cite{HeytingAlgebras}.

\begin{lemma}
In a Heyting algebra, $\vee$ is distributive over $\wedge$ and vice-versa, meaning that
\[x\vee(y\wedge z)=(x\vee y)\wedge(x\vee z)\quad\text{and}\quad x\wedge(y\vee z)=(x\wedge y)\vee(x\wedge z),\]
for all $x$, $y$ and $z$ in the universe.
\end{lemma}

\begin{proof}
We will prove that $x\wedge(y\vee z)=(x\wedge y)\vee(x\wedge z)$, being the other equality proved analogously. Even more, we need only to prove that $x\wedge(y\vee z)\leq (x\wedge y)\vee(x\wedge z)$, since in any lattice one finds that $y\leq y\vee z$ and $z\leq y\vee z$, meaning that $x\wedge y\leq x\wedge (y\vee z)$ and $x\wedge z\leq x\wedge(y\vee z)$, and therefore $(x\wedge y)\vee(x\wedge z)\leq x\wedge (y\vee z)$.

Let $w$ denote $(x\wedge y)\vee(x\wedge z)$, and since $x\wedge y\leq w$ and $x\wedge z\leq w$, we obtain $y\leq x\rightarrow w$ and $z\leq x\rightarrow w$. So $y\vee z\leq x\rightarrow w$, and we get, as desired,
\[x\wedge(y\vee z)\leq x\wedge(x\rightarrow w)\leq w=(x\wedge y)\vee(x\wedge z).\]
\end{proof}

In the following proposition, we omit several more trivial steps, such as proving that $\neg x\wedge x=0$. The reader is invited to fill them in as they appear.

\begin{proposition}
A Heyting algebra $\mathcal{H}$, with universe $H$, is a Boolean algebra if, and only if, $x\vee(x\rightarrow y)=1$ for all $x, y\in H$.
\end{proposition}

\begin{proof}
Start by assuming that we have $x\vee(x\rightarrow y)=1$ for all $x, y\in H$: then, since in particular $0\in H$, $x\vee\neg x=x\vee(x\rightarrow 0)=1$ for all $x\in H$, meaning $\mathcal{H}$ is a Boolean algebra.

Reciprocally, suppose $\mathcal{H}$ is a Boolean algebra. We state that $x\rightarrow y=\neg x\vee y$, being then necessary to prove that $\neg x\vee y$ is the supremum of $\{z\in H: x\wedge z\leq y\}$: since 
\[x\wedge(\neg x\vee y)=(x\wedge\neg x)\vee(x\wedge y)=0\vee (x\wedge y)=x\wedge y,\]
which is smaller or equal to $y$, $\neg x\vee y$ is in fact an element of this set; furthermore, if $x\wedge z\leq y$, then 
\[z\leq \neg x\vee z=1\wedge (\neg x\vee z)=(\neg x\vee x)\wedge (\neg x\vee z)=\neg x\vee(x\wedge z)\leq \neg x\vee y.\]

Then, since $x\rightarrow y=\neg x\vee y$,
\[x\vee (x\rightarrow y)=x\vee(\neg x\vee y)=(x\vee \neg x)\vee y=1\vee y=1,\]
as we wished to prove.
\end{proof}

Of course, this means we could avoid adding a negation to the signature of Heyting algebras to be able to express Boolean algebras: a Boolean algebra is a $\Sigma_{\textbf{Hey}}$-algebra which is a Heyting algebra and satisfies, for all $x$ and $y$ in its universe, $x\vee(x\rightarrow y)=1$. But we prefer to have a negation at hand when dealing with Boolean algebras, it is just more convenient that way.

Now, it is important to point out that some symbols on the signature $\Sigma_{\textbf{Boo}}$ can be changed, usually because we are inserting Boolean algebras in a context where they are already in use. So $0$ may be replaced with "$\bot$", $1$ with "$\top$" and $\neg$ with "$\sim$". We will make such changes more or less freely, without much ado.

As a final remark, we would like to make clear that lattices, implicative or not, and Heyting algebras will play a very minor role in what is to come: they are, more importantly, milestones in defining Boolean algebras that also appear when dealing with some swap structures, specifically for paraconsistent logics (\cite{ParLog}). Boolean algebras, however, will be used time and time again in our study: in Section \ref{bottomless Boolean algebras} we will use them to search for a category equivalent to that of multialgebras, when we will offer a more order-theoretic approach to the subject; Section \ref{B-valuations} shows a new semantics of valuations for da Costa's logics $C_{n}$ based on Boolean algebras, which is used in Section \ref{RNmatrices RMBCn} to offer yet another semantics for those logics, proven in Section \ref{BA and RSwap} to generate a category of models for $C_{n}$ isomorphic to the category of Boolean algebras itself.

\newpage

\printbibliography[segment=\therefsegment,heading=subbibliography]
\end{refsegment}

\begin{refsegment}
\defbibfilter{notother}{not segment=\therefsegment}
\setcounter{chapter}{1}
\chapter{Weakly free multialgebras}\label{Chapter2}\label{Chapter 2}

In the study of universal algebra, it is well known (\cite{Burris}) that there exist $\Sigma$-algebras $\mathcal{A}$ freely generated\index{Freely generated} (or absolutely freely generated) by subsets $X$ of their universes $A$, meaning that for any $\Sigma$-algebra $\mathcal{B}$, with universe $B$, and function $f:X\rightarrow B$, there exists precisely one homomorphism $\overline{f}: \mathcal{A}\rightarrow\mathcal{B}$ extending $f$, and therefore commuting the following diagram in $\textbf{Set}$, for $j:X\rightarrow A$ the inclusion.
\[
\begin{tikzcd}
 A \arrow{ddrr}{\overline{f}} & & \\
 & & \\
X \arrow{uu}{j} \arrow{rr}{f} & & B
\end{tikzcd}
\]
Moreover, the $\Sigma$-algebra freely generated by $X$ is isomorphic to the $\Sigma$-algebra of formulas over $X$, that is $\textbf{F}(\Sigma, X)$, and therefore unique up to isomorphisms; even more, given $\textbf{F}(\Sigma, X)$ and $\textbf{F}(\Sigma, Y)$ are isomorphic whenever $X$ and $Y$ are of the same cardinality, we discover there is precisely one freely generated $\Sigma$-algebra, up to isomorphisms, for each cardinality. Equivalently, in the language of categories, the forgetful functor 
\[U:\textbf{Alg}(\Sigma)\rightarrow\textbf{Set},\]
from the category of $\Sigma$-algebras (with homomorphisms between them) to the category of sets, which takes a $\Sigma$-algebra and returns its underlying universe, has a left adjoint $F$: it associates to a set $X$ any $\Sigma$-algebra freely generated by $X$, and to a function the only homomorphism extending it.

However, an algebraic structure being absolutely free is a concept that does not extend well to the context of multialgebras (\cite{Marty}), specially when one restricts oneself to the non-partial multialgebras (whose multioperations do not return the empty-set), as we often do here given our interest on non-deterministic semantics, specifically those designed for paraconsistency. It is easy to prove that, first of all, freely generated multialgebras, which generalize freely generated algebras in the most obvious way, do not exist, and second, that the forgetful functor 
\[\mathcal{U}:\textbf{MAlg}(\Sigma)\rightarrow\textbf{Set},\]
from the category of multialgebras over the signature $\Sigma$ to the category of sets, does not have a left adjoint. These results are deeply folkloric, being difficult to pinpoint a proof of them in standard literature, although it seems no one is not aware of their validity.

We here expose our reasons to believe that understanding as formulas, in $\textbf{MAlg}(\Sigma)$, only those elements found in the universe of $\textbf{F}(\Sigma, \mathcal{V})$ disregards other multialgebras with an astoundingly similar behavior, so that we generalize the algebras of formulas to multialgebras of formulas. These structures indeed share many of the properties one expects of the algebras of formulas (or, equivalently, freely generated algebras):

\begin{enumerate}
\item they possess the unique extension property not for functions, but rather pairs of functions and what we will call collections of choices, that ``select'' how a homomorphism will approach indeterminacies; 

in fact, given multialgebras $\mathcal{A}$ and $\mathcal{B}$ and an $n$-ary $\sigma$, for all $n$-tuples $(a_{1},\dotsc , a_{n})\in A^{n}$ and $(b_{1}, \dotsc , b_{n})\in B^{n}$ a collection of choices determines how to map those elements of $\sigma_{\mathcal{A}}(a_{1},\dotsc , a_{n})$ into those of $\sigma_{\mathcal{B}}(b_{1}, \dotsc , b_{n})$;

\item they are somewhat ``free'' of identities, an intuition we formalize through disconnected multialgebras, and they are generated by a set of ``indecomposable'' elements we shall call the ground of the multialgebra, much like variables (which are formulas without proper subformulas and therefore indecomposable in some sense); 

more formally, a multialgebra is disconnected whenever different operations, or the same operation performed on different elements, always return disjoint results, while the ground $G(\mathcal{A})$ of a multialgebra $\mathcal{A}$ is the subset of its universe of all elements $a$ for which there do not exist $\sigma$ (of arity $n$) and elements $a_{1}, \dotsc  , a_{n}$ of $\mathcal{A}$ such that $a\in\sigma_{\mathcal{A}}(a_{1}, \dotsc  , a_{n})$;

\item strengthening the previous point, they are disconnected and have a minimum generating set that behaves quite similarly to a basis, of e.g. a vector space;

notice, however, that while a basis of a vector space is a minimal generating set, we are looking here at minimum generating sets, so we use the terminology of ``strong basis'', which end up being precisely the grounds we have mentioned earlier;

\item a final pair of properties we present is that they are simultaneously disconnected and satisfy that every sequence of further simpler and simpler elements eventually ends on an indecomposable element, condition we call being ``chainless'' and that implies having a strong basis;

essentially, by $a$ being simpler than $b$ we mean that there exist $\sigma$ (of arity $n$) and elements $a_{1}, \dotsc  , a_{n}$ such that $a_{i}=a$, for some $i\in\{1, \dotsc  , n\}$, and $b\in\sigma_{\mathcal{A}}(a_{1}, \dotsc  , a_{n})$.
\end{enumerate}

Furthermore, we prove all these four listed items are equivalent, characterizing exactly the same multialgebras; and with these weakly free multialgebras we here define at hand, we can offer simple proofs that freely generated multialgebras (now in the naive generalization) do not exist, and that $\mathcal{U}$ does not have a left adjoint.

Most of the research presented in this chapter was submitted as a preprint in \cite{AbsFreeHyp} and finally published in \cite{WeaklyFreeMultialgebras}.

\section{Formulas}

Here we briefly recall the basic notions involving formulas found in Section \ref{Formulas and how to interpret them}. Given a set $\mathcal{V}$ of propositional variables and a signature $\Sigma=\{\Sigma_{n}\}_{n\in\mathbb{N}}$, the algebra of formulas freely generated by $\mathcal{V}$ over $\Sigma$ will be denoted $\textbf{F}(\Sigma, \mathcal{V})$, and its universe will be denoted $F(\Sigma, \mathcal{V})$. Intuitively, the set of formulas $F(\Sigma, \mathcal{V})$ is the smallest set containing:
\begin{enumerate}
\item the variables $\mathcal{V}$;
\item the formula $\sigma(\alpha_{1}, \dotsc  , \alpha_{n})$, given a $\sigma\in\Sigma_{n}$ and already defined formulas $\alpha_{1}, \dotsc  , \alpha_{n}$ in $F(\Sigma, \mathcal{V})$.
\end{enumerate}

More formally, $F(\Sigma, \mathcal{V})$ should be smallest set containing:
\begin{enumerate}
\item for every $x\in\mathcal{V}$, the function $f_{x}:\{0\}\rightarrow \mathcal{V}\cup \Sigma$ defined by $f_{x}(0)=x$;
\item the function 
\[f_{\sigma(\alpha_{1}, \dotsc  , \alpha_{n})}:\{0, \dotsc  , 1+\sum_{i=1}^{n}m_{n}\}\rightarrow \mathcal{V}\cup\Sigma,\]
given a $\sigma\in\Sigma_{n}$ and already defined formulas $f_{\alpha_{1}}:\{0, \dotsc  , m_{1}\}\rightarrow\mathcal{V}\cup\Sigma, \dotsc  , f_{\alpha_{n}}:\{0, \dotsc  , m_{n}\}\rightarrow\mathcal{V}\cup\Sigma$, such that $f_{\sigma(\alpha_{1}, \dotsc  , \alpha_{n})}(0)=\sigma$ and, for every $j\in\{1, \dotsc  , n\}$ and $k\in\{0, \dotsc  , m_{j}\}$,
\[f_{\sigma(\alpha_{1}, \dotsc  , \alpha_{n})}(k+\sum_{i=0}^{j-1}m_{i})=f_{\alpha_{j}}(k),\]
where for simplicity we define $m_{0}=1$.
\end{enumerate}
One should notice that $f_{\sigma(\alpha_{1}, \dotsc  , \alpha_{n})}$ is simply the formalization of the polish notation of $\sigma(\alpha_{1}, \dotsc  , \alpha_{n})$.

The set $F(\Sigma, \mathcal{V})$ becomes the $\Sigma-$algebra $\textbf{F}(\Sigma, \mathcal{V})$ when we define, for a $\sigma\in\Sigma_{n}$ and formulas $\alpha_{1}, \dotsc  , \alpha_{n}$ in $F(\Sigma, \mathcal{V})$, 
\[\sigma_{\textbf{F}(\Sigma, \mathcal{V})}(\alpha_{1}, \dotsc  , \alpha_{n})=\sigma(\alpha_{1}, \dotsc  , \alpha_{n}).\]

We define the order, or complexity, of an element of $F(\Sigma, \mathcal{V})$ as: $0$ if the formula is a variable or a constant, that is, a $\sigma\in \Sigma_{0}$; as $1+\max\{p_{1}, \dotsc  , p_{n}\}$ if the formula is of the form $\sigma(\alpha_{1}, \dotsc  , \alpha_{n})$, with $\alpha_{j}$ being of order $p_{j}$.

\subsection{Multialgebra of non-deterministic formulas}

\begin{definition}
Given a signature $\Sigma$ and a cardinal $\kappa>0$, the expanded signature\index{Signature, Expanded}\label{expandedsignature} $\Sigma^{\kappa}=\{\Sigma_{n}^{\kappa}\}_{n\in\mathbb{N}}$ is the signature such that $\Sigma_{n}^{\kappa}=\Sigma_{n}\times \kappa$, where we will denote the pair $(\sigma, \beta)$ by $\sigma^{\beta}$, for $\sigma\in\Sigma$ and $\beta\in\kappa$. 
\end{definition}

We demand that $\kappa$ is greater than zero: hence, if $\Sigma$ is non-empty, so is $\Sigma^{\kappa}$.

\begin{definition}\label{Multialgebra of formulas}
Given a set of  variables $\mathcal{V}$, a signature $\Sigma$ and a cardinal $\kappa>0$, we define the $\Sigma$-multialgebra of non-deterministic formulas, or simply multialgebra of formulas\index{Formulas, Non-deterministic}\index{Multialgebra of non-deterministic formulas}\label{mFSigmakappaV}, as
\[\textbf{mF}(\Sigma, \mathcal{V}, \kappa)=(F(\Sigma^{\kappa}, \mathcal{V}), \{\sigma_{\textbf{mF}(\Sigma, \mathcal{V}, \kappa)}\}_{\sigma\in\Sigma})\]
with universe $F(\Sigma^{\kappa}, \mathcal{V})$ and such that, for $\sigma\in \Sigma_{n}$ and $\alpha_{1}, \dotsc   , \alpha_{n}\in {T}(\Sigma^{\kappa}, \mathcal{V})$,
\[\sigma_{\textbf{mF}(\Sigma, \mathcal{V}, \kappa)}(\alpha_{1}, \dotsc   , \alpha_{n})=\{\sigma^{\beta}(\alpha_{1}, \dotsc  , \alpha_{n}) \ : \  \beta\in\kappa\}.\]
\end{definition}

The intuition behind this definition is that, connecting given formulas $\alpha_{1}$ through $\alpha_{n}$ with a connective $\sigma$ can, in a broader interpretation taking into account non-determinism, return many formulas with the same general shape, that is $\sigma(\alpha_{1}, \dotsc  ,\alpha_{n})$, over which we maintain certain degree of control by counting them, what we achieve by using an index to our connective, $\sigma^{\beta}$.

One can ask why all connectives must return the exact same number of generalized formulas, that is $\kappa$, but this will not be the case: more useful to us shall be the submultialgebras of $\textbf{mF}(\Sigma, \mathcal{V}, \kappa)$, where the cardinality will vary as long as it is bounded by $\kappa$; we have defined the multialgebras of formulas as above since defining its submultialgebras directly is substantially more difficult.

 Here, we will restrict ourselves to the cases where $\Sigma_{0}\neq\emptyset$ or $\mathcal{V}\neq\emptyset$, so that $\textbf{mF}(\Sigma, \mathcal{V}, \kappa)$ is always well defined. 

We will understand as the order of an element $\alpha$ of $\textbf{mF}(\Sigma, \mathcal{V}, \kappa)$ simply its order as an element of $F(\Sigma^{\kappa}, \mathcal{V})$. Notice that, if 
\[\sigma_{\textbf{mF}(\Sigma, \mathcal{V}, \kappa)}(\alpha_{1}, \dotsc   , \alpha_{n})\cap \theta_{\textbf{mF}(\Sigma, \mathcal{V}, \kappa)}(\beta_{1}, \dotsc    ,\beta_{m})\neq\emptyset,\]
then $\sigma=\theta$, $n=m$ and $\alpha_{1}=\beta_{1}, \dotsc , \alpha_{n}=\beta_{m}$, since if the intersection is not empty there are $\beta, \gamma\in\kappa$ such that $\sigma^{\beta}(\alpha_{1}, \dotsc  , \alpha_{n})= \theta^{\gamma}(\beta_{1}, \dotsc  ,\beta_{m})$ and by the structure of $F(\Sigma^{\kappa}, \mathcal{V})$ we find that $\sigma^{\beta}=\theta^{\gamma}$. 

\begin{example}
The $\Sigma$-algebras of formulas $\textbf{F}(\Sigma, \mathcal{V})$, when considered as multialgebras such that $\sigma_{\textbf{F}(\Sigma, \mathcal{V})}(\alpha_{1}, \dotsc   , \alpha_{n})=\{\sigma(\alpha_{1}, \dotsc  ,\alpha_{n})\}$, are multialgebras of formulas, with $\kappa=1$; that is, $\textbf{F}(\Sigma, \mathcal{V})$ and $\textbf{mF}(\Sigma, \mathcal{V}, 1)$ are isomorphic.
\end{example}

From now on, the cardinal of a set $X$ will be denoted by $|X|$.\label{|X|}

\begin{example}\label{s}
A directed graph is a pair $(V, A)$, with $V$ a non-empty set of elements called vertices and $A\subseteq V^{2}$ a set of elements called arrows, where we say that there is an arrow from $u$ to $v$, both in $V$, if $(u,v)\in A$; we say that the $n$-tuple $(v_{1}, \dotsc   , v_{n})$ is a path between $u$ and $v$ if $u=v_{1}$, $v=v_{n}$ and $(v_{i}, v_{i+1})\in A$ for every $i\in\{1, \dotsc   , n-1\}$; we say that a vertex $u\in V$ has a successor if there exists $v\in V$ such that $(u,v)\in A$, and $u$ has a predecessor if there exists $v\in V$ such that $(v,u)\in A$.

A directed graph $F=(V, A)$ is a forest if, for any two vertices $u, v\in V$, there exists at most one path between $u$ and $v$, and a forest is said to have height $\omega$ if every vertex has a successor. Then, we state that forests of height $\omega$ are in bijection with the submultialgebras of the multialgebras of formulas over the signature $\Sigma_{s}$\label{Sigmas} with exactly one operator $s$ of arity $1$.

Essentially, take as $\mathcal{V}$ the set of elements of $F$ that have no predecessor and define, for $u\in V$,
\[s_{\mathcal{A}}(u)=\{v\in V \ : \ (u,v)\in A\},\]
and we have that the $\Sigma_{s}$-multialgebra $\mathcal{A}=(V, \{s_{\mathcal{A}}\})$, submultialgebra of $\textbf{mF}(\Sigma_{s}, \mathcal{V}, |V|)$, carries the same information that $F$.
\end{example}

\begin{example}
More generally, a directed multi-graph\index{Multi-graph} \cite{ConiglioSernadas}, or directed $m$-graph, is a pair $(V, A)$ with $V$ a non-empty set of vertices and $A$ a subset of $V^{+}\times V$, where $V^{+}=\bigcup_{n\in\mathbb{N}\setminus\{0\}}V^{n}$ is the set of finite, non-empty, sequences over $V$. We will say that $(v_{1}, \dotsc   , v_{n})$ is a path between $u$ and $v$ if $u=v_{1}$, $v=v_{n}$ and, for every $i\in\{1, \dotsc   , n-1\}$, there exist $v_{i_{1}}, \dotsc   , v_{i_{m}}$ such that $((v_{i_{1}}, \dotsc   , v_{i_{m}}), v_{i+1})$, with $v_{i}=v_{i_{j}}$ for some $j\in\{1, \dotsc   , m\}$.

Then an $m$-forest is a directed $m$-graph such that any two elements are connected by at most one path; and an $m$-forest is said to have $n$-height $\omega$, for $n\in\mathbb{N}\setminus\{0\}$, if, for any $(u_{1}, \dotsc   , u_{n})\in V^{n}$, there exists $v\in V$ such that $((u_{1}, \dotsc   , u_{n}), v)\in A$. Finally, we see that every $m$-forest $F=(V, A)$ with $n$-height $\omega$, for every $n\in S\subseteq\mathbb{N}\setminus\{\emptyset\}$, is essentially equivalent to the $\Sigma_{S}$-multialgebra $\mathcal{A}=(V, \{\sigma_{\mathcal{A}}\}_{\sigma\in\Sigma_{S}})$, with
\[\sigma_{\mathcal{A}}(u_{1}, \dotsc   , u_{m})=\{v\in V \ : \ ((u_{1}, \dotsc   , u_{m}), v)\in A\},\]
for $\sigma$ of arity $m$, and $\Sigma_{S}$ the signature with exactly one operator of arity $n$, for every $n\in S$. It is not hard to see that $\mathcal{A}$ is a submultialgebra of $\textbf{mF}(\Sigma_{S}, \mathcal{V}, |V|)$, with $\mathcal{V}$ the set of elements $v$ of $V$ such that, for no $(u_{1}, \dotsc   , u_{n})\in V^{+}$, $((u_{1}, \dotsc   , u_{n}), v)\in A$.
\end{example}

\section{Equivalences for being a submultialgebra of $\textbf{mF}(\Sigma, \mathcal{V}, \kappa)$}
\subsection{Being $\textbf{cdf}$-generated}
 
Now, in universal algebra, the algebras of formulas $\textbf{F}(\Sigma, \mathcal{V})$ are absolutely free, also said to be freely generated, also said to be freely generated in the variety of all $\Sigma$-algebras: this means that there exists a set, in their case the set of variables $\mathcal{V}$, such that, for every other $\Sigma$-algebra $\mathcal{B}$ with universe $B$ and function $f:\mathcal{V}\rightarrow B$, there exists an unique homomorphism $\overline{f}:\textbf{F}(\Sigma, \mathcal{V})\rightarrow\mathcal{B}$ extending $f$, essentially defined as:
\begin{enumerate}
\item $\overline{f}(p)=f(p)$, for every $p\in \mathcal{V}$;
\item $\overline{f}(\sigma(\alpha_{1}, \dotsc  , \alpha_{n}))=\sigma_{\mathcal{B}}(\overline{f}(\alpha_{1}), \dotsc  , \overline{f}(\alpha_{n}))$. 
\end{enumerate}
As we mentioned before, this is no longer true when dealing with multialgebras, but we can define a closely related concept with the aid of what we will call collections of choices.

Collections of choices are motivated by legal valuations, notion first defined in Avron and Lev's seminal paper \cite{AvronLev} on non-deterministic logical matrices. A map $\nu$ from $\textbf{F}(\Sigma, \mathcal{V})$ to the universe of a $\Sigma$-multialgebra $\mathcal{A}$ is a legal valuation whenever 
\[\nu(\sigma(\alpha_{1}, \dotsc  ,\alpha_{n}))\in \sigma_{\mathcal{A}}(\nu(\alpha_{1}), \dotsc   , \nu(\alpha_{n}));\]
essentially, at every formula $\sigma(\alpha_{1}, \dotsc   ,\alpha_{n})$, $\nu$ ``chooses'' a value from all the possible values\\ $\sigma_{\mathcal{A}}(\nu(\alpha_{1}), \dotsc   , \nu(\alpha_{n}))$, possible values which depend themselves on previous choices $\nu(\alpha_{1})$ through $\nu(\alpha_{n})$ performed by $\nu$. What a collection of choices does is then to automatize these choices, what justifies its name. 

\begin{definition}\label{Collection of Choices}
Given multialgebras $\mathcal{A}=(A, \{\sigma_{\mathcal{A}}\}_{\sigma\in\Sigma})$ and $\mathcal{B}=(B, \{\sigma_{\mathcal{B}}\}_{\sigma\in\sigma})$ over the signature $\Sigma$, a collection of choices\index{Collection of choices} from $\mathcal{A}$ to $\mathcal{B}$ is a collection $C=\{C_{n}\}_{n\in\mathbb{N}}$ of collections of functions
\[C_{n}=\{C\sigma_{a_{1}, \dotsc  , a_{n}}^{b_{1}, \dotsc  , b_{n}} : \sigma\in\Sigma_{n}, a_{1}, \dotsc  , a_{n}\in A, b_{1}, \dotsc   ,b_{n}\in B\}\]
such that, for $\sigma\in\Sigma_{n}$, $a_{1}, \dotsc  , a_{n}\in A$ and $b_{1}, \dotsc  , b_{n}\in B$, $C\sigma_{a_{1}, \dotsc  , a_{n}}^{b_{1}, \dotsc  , b_{n}}$ is a function of the form 
\[C\sigma_{a_{1}, \dotsc  , a_{n}}^{b_{1}, \dotsc  , b_{n}}:\sigma_{\mathcal{A}}(a_{1}, \dotsc   ,a_{n})\rightarrow\sigma_{\mathcal{B}}(b_{1}, \dotsc  , b_{n}).\]
\end{definition}

\begin{example}
If $\mathcal{B}$ is actually an algebra (meaning all its operations return singletons), there only exists one collection of choices from any $\mathcal{A}$ to $\mathcal{B}$ (since there exists only one function to a set with only one element); this means that in universal algebra, collections of choices are somewhat irrelevant.
\end{example}

\begin{example}
A directed tree is a directed forest where there exists exactly one element without predecessor; we say that $v$ ramifies from $u$ if there exists an arrow from $u$ to $v$. Then, given two directed trees $T_{1}=(V_{1}, A_{1})$ and $T_{2}=(V_{2}, A_{2})$ of height $\omega$, seem as $\Sigma_{s}$-multialgebras, and a collection of choices $C$ from $T_{1}$ to $T_{2}$, for every $v\in V_{1}$ and $u\in V_{2}$ the function $Cs_{v}^{u}$ chooses, for each of the elements that ramify from $v$, one element that ramifies from $u$.
\end{example}

\begin{definition}
Given a signature $\Sigma$, a $\Sigma$-multialgebra $\mathcal{A}=(A, \{\sigma_{\mathcal{A}}\}_{\sigma\in\Sigma})$ is choice-dependent freely generated\index{Choice-dependent freely generated} by $X$ if $X\subseteq A$ and, for all $\Sigma$-multialgebras $\mathcal{B}=(B, \{\sigma_{\mathcal{B}}\}_{\sigma\in\Sigma})$, all functions $f:X\rightarrow B$ and all collections of choices $C$ from $\mathcal{A}$ to $\mathcal{B}$, there is a unique homomorphism $f_{C}:\mathcal{A}\rightarrow \mathcal{B}$ such that:

\begin{enumerate} \item $f_{C}|_{X}=f$;

\item for all $\sigma\in \Sigma_{n}$ and $a_{1}, \dotsc   , a_{n}\in A$, 
\[f_{C}|_{\sigma_{\mathcal{A}}(a_{1}, \dotsc   , a_{n})}=C\sigma_{a_{1}, \dotsc   , a_{n}}^{f_{C}(a_{1}), \dotsc   , f_{C}(a_{n})}.\]
\end{enumerate}
\end{definition}

For simplicity, when $\mathcal{A}$ is choice-dependent freely generated by $X$, we will write that $\mathcal{A}$ is $\textbf{cdf}$-generated\label{cdfgenerated} by $X$, or merely that $\mathcal{A}$ is $\textbf{cdf}$-generated, when the set $X$ is not important.

We now introduce the concept of ground to indicate what elements of a multialgebra are not ``achieved'', ``reached'' by its multioperations; alternatively, while thinking of formulas and their respective algebras, the ground is the set of indecomposable formulas, that is, variables.

\begin{definition}
Given a $\Sigma$-multialgebra $\mathcal{A}=(A, \{\sigma_{\mathcal{A}}\}_{\sigma\in\Sigma})$, we define its build\index{Build}\label{build} as
\[B(\mathcal{A})=\bigcup \big\{\sigma_{\mathcal{A}}(a_{1}, \dotsc   , a_{n}) \ : \  n\in\mathbb{N}, \, \sigma\in\Sigma_{n}, \, a_{1}, \dotsc   , a_{n}\in A \big\}.\]
We define the ground\index{Ground}\label{ground} of $\mathcal{A}$ as
\[G(\mathcal{A})=A\setminus B(\mathcal{A}).\]
\end{definition}

\begin{example}\label{ground of formulas}
$B(\textbf{F}(\Sigma, \mathcal{V}))={T}(\Sigma, \mathcal{V})\setminus\mathcal{V}$ and $G(\textbf{F}(\Sigma, \mathcal{V}))=\mathcal{V}$.
\end{example}

\begin{example}
If $F=(V, A)$ is a directed forest of height $\omega$, thought as a $\Sigma_{s}$-multialgebra, its ground is the set of elements $v$ in $V$ without predecessors.
\end{example}

\begin{proposition}
Let $\mathcal{A}$ and $\mathcal{B}$ be $\Sigma$-multialgebras.
\begin{enumerate} 
\item If $f:\mathcal{A}\rightarrow\mathcal{B}$ is a homomorphism between $\Sigma$-multialgebras, then 
\[B(\mathcal{A})\subseteq f^{-1}(B(\mathcal{B}))\quad\text{and}\quad f^{-1}(G(\mathcal{B}))\subseteq G(\mathcal{A}).\]

\item If $\mathcal{B}$ is a submultialgebra of $\mathcal{A}$, $B(\mathcal{B})\subseteq B(\mathcal{A})$ and $G(\mathcal{A})\subseteq G(\mathcal{B})$.
\end{enumerate}
\end{proposition}

\begin{proof}
\begin{enumerate} \item If $a\in B(\mathcal{A})$, there exist $\sigma\in\Sigma_{n}$ and $a_{1}, \dotsc   , a_{n}\in A$ such that $a\in\sigma_{\mathcal{A}}(a_{1}, \dotsc   , a_{n})$. Since $f(\sigma_{\mathcal{A}}(a_{1}, \dotsc   , a_{n})) \subseteq \sigma_{\mathcal{B}}(f(a_{1}), \dotsc    ,f(a_{n}))$,
we find that $f(a)\in \sigma_{\mathcal{B}}(f(a_{1}), \dotsc   , f(a_{n}))$ and therefore $f(a)\in B(\mathcal{B})$, meaning that $a\in f^{-1}(B(\mathcal{B}))$. Using that $G(\mathcal{A})=A\setminus B(\mathcal{A})$ and $G(\mathcal{B})=B\setminus B(\mathcal{B})$ we obtain the second mentioned inclusion.

\item If $b\in B(\mathcal{B})$, there exist $\sigma\in\Sigma_{n}$ and $b_{1}, \dotsc   , b_{n}\in B$ such that $b\in \sigma_{\mathcal{B}}(b_{1}, \dotsc   , b_{n})$, and given that $\sigma_{\mathcal{B}}(b_{1}, \dotsc   , b_{n})\subseteq \sigma_{\mathcal{A}}(b_{1}, \dotsc   , b_{n})$ we obtain $b\in B(\mathcal{A})$. Using again that $G(\mathcal{A})=A\setminus B(\mathcal{A})$ and $G(\mathcal{B})=B\setminus B(\mathcal{B})$ we finish the proof.
\end{enumerate}
\end{proof}

From this it also follows that if $f:\mathcal{A}\rightarrow\mathcal{B}$ is a homomorphism, $G(\mathcal{B})\cap f(A)$ is contained in $\{f(a)\ : \  a\in G(\mathcal{A})\}$. Indeed,  if $b$ is in $G(\mathcal{B})\cap f(A)$, $a\in A$ such that $f(a)=b$ is in $f^{-1}(G(\mathcal{B}))$ and, by the previous proposition, it is in $G(\mathcal{A})$, and therefore $b$ is in $\{f(a) \ : \  a\in G(\mathcal{A})\}$.

Generalizing Example \ref{ground of formulas}, we have that $G(\textbf{mF}(\Sigma, \mathcal{V}, \kappa))=\mathcal{V}$: we proceed by induction to show that $B(\textbf{mF}(\Sigma, \mathcal{V}, \kappa))={T}(\Sigma^{\kappa}, \mathcal{V})\setminus \mathcal{V}$, what is equivalent. If $\alpha$ is of order $0$, either we have $\alpha=\sigma^{\beta}$, for a $\sigma\in\Sigma_{0}$ and $\beta\in\kappa$, and therefore $\alpha\in B(\textbf{mF}(\Sigma, \mathcal{V}, \kappa))$; or we have that $\alpha=p\in \mathcal{V}$, and if there exists $\sigma\in\Sigma_{m}$ and $\alpha_{1}, \dotsc   , \alpha_{m}\in {T}(\Sigma^{\kappa}, \mathcal{V})$ such that 
\[p\in \sigma_{\textbf{mF}(\Sigma, \mathcal{V}, \kappa)}(\alpha_{1}, \dotsc   , \alpha_{m})\]
we have $p=\sigma^{\beta}(\alpha_{1}, \dotsc  ,\alpha_{m})$, for $\beta\in\kappa$, which is absurd given the structure of $F(\Sigma^{\kappa}, \mathcal{V})$, forcing us to conclude that $x\notin B(\textbf{mF}(\Sigma, \mathcal{V}, \kappa))$. If $\alpha$ is of order $n+1$, we have that $\alpha=\sigma^{\beta}(\alpha_{1}, \dotsc  ,\alpha_{m})$ for a $\sigma\in\Sigma_{m}$, $\beta\in\kappa$ and $\alpha_{1}, \dotsc   , \alpha_{m}$ of order at most $n$, and therefore we have $\alpha$ in $\sigma_{\textbf{mF}(\Sigma, \mathcal{V}, \kappa)}(\alpha_{1}, \dotsc   , \alpha_{m})$, meaning that $\alpha\in B(\textbf{mF}(\Sigma, \mathcal{V}, \kappa))$.

\begin{definition}
Given a $\Sigma$-multialgebra $\mathcal{A}=(A, \{\sigma_{\mathcal{A}}\}_{\sigma\in\Sigma})$ and a set $S\subseteq A$, we define the sets $\langle S\rangle_{m}$\label{mth generated} by induction: $\langle S\rangle _{0}=S\cup\bigcup_{\sigma\in \Sigma_{0}}\sigma_{\mathcal{A}}$; and assuming we have defined $\langle S\rangle_{m}$, we make
\[\langle S\rangle_{m+1}=\langle S\rangle_{m}\cup \bigcup \big\{\sigma_{\mathcal{A}}(a_{1}, \dotsc   , a_{n}) \ : \  n\in\mathbb{N}, \, \sigma\in\Sigma_{n}, \, a_{1}, \dotsc   , a_{n}\in \langle S\rangle_{m} \big\}.\]
The set generated\index{Generated}\label{generated} by $S$, denoted by $\langle S\rangle$, is then defined as $\langle S\rangle=\bigcup_{m\in\mathbb{N}}\langle S\rangle_{m}$.

We say $\mathcal{A}$ is generated by $S$ if $\langle S\rangle=A$.
\end{definition}

\begin{lemma}\label{sub mF is gen ground}
Every submultialgebra $\mathcal{A}$ of $\textbf{mF}(\Sigma, \mathcal{V}, \kappa)$ is generated by $G(\mathcal{A})$.
\end{lemma}

\begin{proof}
Suppose $a$ is an element of $\mathcal{A}$ not contained in $\langle G(\mathcal{A})\rangle$ of minimum order: since $a$ cannot belong to $G(\mathcal{A})\cup\bigcup_{\sigma\in\Sigma_{0}}\sigma_{\mathcal{A}}=\langle G(\mathcal{A})\rangle_{0}$, there exist $n>0$, $\sigma\in\Sigma_{n}$ and $a_{1}, \dotsc   , a_{n}\in A$ such that $a\in\sigma_{\mathcal{A}}(a_{1}, \dotsc   , a_{n})$.

Since $\sigma_{\mathcal{A}}(a_{1}, \dotsc   , a_{n})\subseteq \sigma_{\textbf{mF}(\Sigma, X, \kappa)}(a_{1}, \dotsc   , a_{n})$ we derive that $a_{1}$ to $a_{n}$ are of order less than that of $a$: by our hypothesis, there must exist $m_{1}, \dotsc   , m_{n}$ such that $a_{j}\in \langle G(\mathcal{A})\rangle_{m_{j}}$ for all $j\in\{1, \dotsc   , n\}$; taking $m=\max\{m_{1}, \dotsc   , m_{n}\}$, $a_{1}, \dotsc   , a_{n}\in \langle G(\mathcal{A})\rangle_{m}$, and therefore
\[a\in \sigma_{\mathcal{A}}(a_{1}, \dotsc   , a_{n})\subseteq \langle G(\mathcal{A})\rangle_{m+1},\]
which contradicts our assumption and proves the lemma.
\end{proof}

\begin{theorem}\label{sub mF is cdf-gen ground}
Every submultialgebra $\mathcal{A}$ of $\textbf{mF}(\Sigma, \mathcal{V}, \kappa)$ is $\textbf{cdf}$-generated by $G(\mathcal{A})$.
\end{theorem}

\begin{proof}
Let $\mathcal{A}=(A, \{\sigma_{\mathcal{A}}\}_{\sigma\in\Sigma})$ be a submultialgebra of $\textbf{mF}(\Sigma, \mathcal{V}, \kappa)$, let $\mathcal{B}=(B, \{\sigma_{\mathcal{B}}\}_{\sigma\in\Sigma})$ be any $\Sigma$-multialgebra, let $f:G(\mathcal{A})\rightarrow B$ be a function and $C$ a collection of choices from $\mathcal{A}$ to $\mathcal{B}$. We define $f_{C}:\mathcal{A}\rightarrow\mathcal{B}$ by induction on $\langle G(\mathcal{A})\rangle_{m}$:

\begin{enumerate} \item if $a\in \langle G(\mathcal{A})\rangle_{0}$ and $a\in G(\mathcal{A})$, we define $f_{C}(a)=f(a)$;

\item if $a\in \langle G(\mathcal{A})\rangle_{0}$ and $a\in\sigma_{\mathcal{A}}$, for a $\sigma\in\Sigma_{0}$, we define $f_{C}(a)=C\sigma(a)$;

\item if $f_{C}$ is defined for all elements of $\langle G(\mathcal{A})\rangle_{m}$, $a_{1}, \dotsc   , a_{n}\in \langle G(\mathcal{A})\rangle_{m}$ and $\sigma\in\Sigma_{n}$, for every element $a\in \sigma_{\mathcal{A}}(a_{1}, \dotsc   , a_{n})$ we define 
\[f_{C}(a)=C\sigma_{a_{1}, \dotsc   , a_{n}}^{f_{C}(a_{1}), \dotsc   , f_{C}(a_{n})}(a).\]
\end{enumerate}

First, we must prove that $f_{C}$ is well defined. There are two possibly problematic cases to consider for an element $a\in A$: the one in which $a\in G(\mathcal{A})$ and there are $\sigma\in\Sigma_{n}$ and $a_{1}, \dotsc   , a_{n}\in A$ for which $a\in\sigma_{\mathcal{A}}(a_{1}, \dotsc   , a_{n})$, corresponding to $a$ falling simultaneously in the cases $(1)$ and $(2)$, or $(1)$ and $(3)$ of the definition; and the one where there are $\sigma\in\Sigma_{n}$, $\theta\in\Sigma_{m}$, $a_{1}, \dotsc   , a_{n}\in A$ and $b_{1}, \dotsc   , b_{m}\in A$ such that $a\in\sigma_{\mathcal{A}}(a_{1}, \dotsc   , a_{n})$ and $a\in\theta_{\mathcal{A}}(b_{1}, \dotsc   , b_{m})$, corresponding to the cases $(2)$ and $(3)$, $(2)$ and $(2)$, or $(3)$ and $(3)$ occurring simultaneously.

The first case is not possible, since $G(\mathcal{A})\subseteq A\setminus\sigma_{\mathcal{A}}(a_{1}, \dotsc   , a_{n})$ for every $\sigma\in\Sigma_{n}$ and $a_{1}, \dotsc   , a_{n}\in A$; in the second case, we find that 
\[a\in \sigma_{\mathcal{A}}(a_{1}, \dotsc   , a_{n})\cap\theta_{\mathcal{A}}(b_{1}, \dotsc   , b_{m})\subseteq\sigma_{\textbf{mF}(\Sigma, \mathcal{V}, \kappa)}(a_{1}, \dotsc   , a_{n})\cap \theta_{\textbf{mF}(\Sigma, \mathcal{V}, \kappa)}( b_{1}, \dotsc    , b_{m}),\]
so $n=m$, $\sigma=\theta$ and $a_{1}=b_{1}, \dotsc   , a_{n}=b_{m}$, and therefore $f_{C}(a)$ is well-defined.

Second, we must prove that $f_{C}$ is defined over all of $A$: that is simple, for $f_{C}$ is defined over all of $\langle G(\mathcal{A})\rangle$ and we established in Lemma~\ref{sub mF is gen ground} that $A=\langle G(\mathcal{A})\rangle$.

So $f_{C}:A\rightarrow B$ is a well-defined function: it remains to be shown that it is a homomorphism; given $\sigma\in\Sigma_{n}$ and $a_{1}, \dotsc   , a_{n}$, we see that
\[f_{C}(\sigma_{\mathcal{A}}(a_{1}, \dotsc   , a_{n}))=\big\{C\sigma_{a_{1}, \dotsc   , a_{n}}^{f_{C}(a_{1}), \dotsc   , f_{C}(a_{n})}(a)\ : \  a\in\sigma_{\mathcal{A}}(a_{1}, \dotsc   , a_{n})\big\}\subseteq\sigma_{\mathcal{B}}(f_{C}(a_{1}), \dotsc   , f_{C}(a_{n})),\]
while we also have that $f_{C}$ clearly extends both $f$ and all $C\sigma_{a_{1}, \dotsc   , a_{n}}^{f_{C}(a_{1}, \dotsc   , f_{C}(a_{n})}$. 

To finish the proof, suppose $g:\mathcal{A}\rightarrow\mathcal{B}$ is another homomorphism extending both $f$ and all $C\sigma_{a_{1}, \dotsc   , a_{n}}^{g(a_{1}, \dotsc   , g(a_{n})}$: we will prove that $g=f_{C}$ again by induction on the $m$ of $\langle G(\mathcal{A})\rangle_{m}$. For $m=0$, an element $a\in \langle G(\mathcal{A})\rangle_{0}$ is either in $G(\mathcal{A})$, when we have $g(a)=f(a)=f_{C}(a)$, or in $\sigma_{\mathcal{A}}$ for a $\sigma\in\Sigma_{0}$, when $g(a)=C\sigma(a)=f_{C}(a)$.

Suppose $g$ is equal to $f_{C}$ in $\langle G(\mathcal{A})\rangle_{m}$ and take an $a\in \langle G(\mathcal{A})\rangle_{m+1}\setminus \langle G(\mathcal{A})\rangle_{m}$: we have that there exist $\sigma\in\Sigma_{n}$ and $a_{1}, \dotsc   , a_{n}\in \langle G(\mathcal{A})\rangle_{m}$ such that $a\in\sigma_{\mathcal{A}}(a_{1}, \dotsc   , a_{n})$, and then 
\[g(a)=C\sigma_{a_{1}, \dotsc   , a_{n}}^{g(a_{1}), \dotsc   , g(a_{n})}(a)=C\sigma_{a_{1}, \dotsc   , a_{n}}^{f_{C}(a_{1}), \dotsc   , f_{C}(a_{n})}(a)=f_{C}(a),\]
proving that $g=f_{C}$ and that, in fact, $f_{C}$ is unique. This means that $\mathcal{A}$ is $\textbf{cdf}$-generated by $G(\mathcal{A})$.
\end{proof}

The following lemma may be found in section 2 of \cite{CFG}.

\begin{lemma}\label{image of hom is sub}
If $\mathcal{A}=(A, \{\sigma_{\mathcal{A}}\}_{\sigma\in\Sigma})$ and $\mathcal{B}=(B, \{\sigma_{\mathcal{B}}\}_{\sigma\in\Sigma})$ are $\Sigma$-multialgebras and $f:\mathcal{A}\rightarrow\mathcal{B}$ is a homomorphism, $\mathcal{C}=(f(A), \{\sigma_{\mathcal{C}}\}_{\sigma\in\Sigma})$ such that
\[\sigma_{\mathcal{C}}(c_{1}, \dotsc  , c_{n})=\bigcup\{f(\sigma_{\mathcal{A}}(a_{1}, \dotsc  , a_{n})) : f(a_{1})=c_{1}, \dotsc  , f(a_{n})=c_{n}\}\]
is a $\Sigma$-submultialgebra of $\mathcal{B}$, while $f:\mathcal{A}\rightarrow\mathcal{C}$ is an epimorphism.\footnote{A epimorphism between $\Sigma$-multialgebras $\mathcal{A}=(A, \{\sigma_{\mathcal{A}}\}_{\sigma\in\Sigma})$ and $\mathcal{B}=(B, \{\sigma_{\mathcal{B}}\}_{\sigma\in\Sigma})$ is defined, as usual, as a homomorphism $\varphi:A\rightarrow B$ between $\mathcal{A}$ and $\mathcal{B}$ that is surjective.} The $\Sigma$-multialgebra $\mathcal{C}$ is known as the direct image\index{Direct image} of $\mathcal{A}$ trough $f$.
\end{lemma}

\begin{proof}
First of all, $A$ is not empty, and therefore so is $f(A)$.

Take $c_{1}, \dotsc  , c_{n}\in f(A)$: there must exist $a_{1}, \dotsc  , a_{n}\in A$ such that $f(a_{1})=c_{1}, \dotsc ,$\\ $f(a_{n})=c_{n}$, and since $\mathcal{A}$ is a multialgebra, $\sigma_{\mathcal{A}}(a_{1}, \dotsc  , a_{n})\neq\emptyset$, implying that $f(\sigma_{\mathcal{A}}(a_{1}, \dotsc  , a_{n}))$ is not empty and, therefore, that $\sigma_{\mathcal{C}}(c_{1}, \dotsc  , c_{n})$ is non-empty, given it contains $f(\sigma_{\mathcal{A}}(a_{1}, \dotsc  , a_{n}))$. It is obvious that, as defined, $\sigma_{\mathcal{C}}(c_{1}, \dotsc  , c_{n})$ is a subset of $f(A)$, and so we can deduce that $\mathcal{C}$ is a $\Sigma$-multialgebra.

Now, given $f:\mathcal{A}\rightarrow\mathcal{B}$ is a homomorphism, for all $a_{1}, \dotsc  , a_{n}\in A$ we have $f(\sigma_{\mathcal{A}}(a_{1}, \dotsc  , a_{n}))\subseteq \sigma_{\mathcal{B}}(f(a_{1}), \dotsc  , f(a_{n}))$, meaning
\[\sigma_{\mathcal{C}}(c_{1}, \dotsc  , c_{n})=\bigcup\{f(\sigma_{\mathcal{A}}(a_{1}, \dotsc  , a_{n})) : f(a_{1})=c_{1}, \dotsc  , f(a_{n})=c_{n}\}\subseteq\]
\[\bigcup\{\sigma_{\mathcal{B}}(f(a_{1}), \dotsc  , f(a_{n})) : f(a_{1})=c_{1}, \dotsc  , f(a_{n})=c_{n}\}=\sigma_{\mathcal{B}}(c_{1}, \dotsc  , c_{n}),\]
or what is equivalent, that $\mathcal{C}$ is a submultialgebra of $\mathcal{B}$.

Finally, $f:A\rightarrow f(A)$ is still a well-defined function, obviously surjective: for any $n$-ary $\sigma$ and elements $a_{1}$ through $a_{n}$ of $A$, one has
\[f(\sigma_{\mathcal{A}}(a_{1}, \dotsc  , a_{n}))\subseteq \bigcup\{f(\sigma_{\mathcal{A}}(a'_{1}, \dotsc  , a'_{n})) : f(a'_{1})=f(a_{1}), \dotsc  , f(a'_{n})=f(a_{n})\}=\]
\[\sigma_{\mathcal{C}}(f(a_{1}), \dotsc  , f(a_{n}))\]
and that, in conclusion, $f$ is a homomorphism.
\end{proof}

\begin{theorem}\label{cdf-gen is iso to sub mF}
If the multialgebra $\mathcal{A}=(A, \{\sigma_{\mathcal{A}}\}_{\sigma\in\Sigma})$ over $\Sigma$ is $\textbf{cdf}$-generated by $X$, then $\mathcal{A}$ is isomorphic to a submultialgebra of $\textbf{mF}(\Sigma, X, |A|)$ containing $X$.
\end{theorem}

\begin{proof}
Take $f:X\rightarrow F(\Sigma^{|A|}, X)$ to be the inclusion (such that $f(x)=x$), and take a collection of choices $C$ such that, for $\sigma\in\Sigma_{n}$, $a_{1}, \dotsc   , a_{n}\in A$ and $\alpha_{1}, \dotsc   , \alpha_{n}\in F(\Sigma^{|A|}, X)$, 
\[C\sigma_{a_{1}, \dotsc   , a_{n}}^{\alpha_{1}, \dotsc   , \alpha_{n}}:\sigma_{\mathcal{A}}(a_{1}, \dotsc   , a_{n})\rightarrow\sigma_{\textbf{mF}(\Sigma, X, |A|)}(\alpha_{1}, \dotsc   , \alpha_{n})\]
is an injective function. Such collection of choices exist since $\sigma_{\mathcal{A}}(a_{1}, \dotsc   , a_{n})\subseteq A$ and\\ $\sigma_{\textbf{mF}(\Sigma, X, |A|)}(\alpha_{1}, \dotsc   , \alpha_{n})$ is of cardinality $|A|$. Now, since $\mathcal{A}$ is $\textbf{cdf}$-generated by $X$, there exists a homomorphism $f_{C}:\mathcal{A}\rightarrow\textbf{mF}(\Sigma, X, |A|)$ extending $f$ and each $C\sigma_{a_{1}, \dotsc   , a_{n}}^{f_{C}(a_{1}), \dotsc   , f_{C}(a_{n})}$.

Let $\mathcal{B}=(f_{C}(A), \{\sigma_{\mathcal{B}}\}_{\sigma\in\Sigma})$ be the direct image of $\mathcal{A}$ trough $f_{C}$, so that $f_{C}:\mathcal{A}\rightarrow\mathcal{B}$ is an epimorphism, what is possible given Lemma \ref{image of hom is sub}: notice too that 
\[X=X\cap f_{C}(A)=G(\textbf{mF}(\Sigma, X, |A|))\cap f_{C}(A)\subseteq G(\mathcal{B})\]
because $\mathcal{B}$ is a submultialgebra of $\textbf{mF}(\Sigma, X, |A|)$. Take any $g:G(\mathcal{B})\rightarrow A$ such that $g(x)=x$, for every $x\in X$. And take a collection of choices $D$ from $\mathcal{B}$ to $\mathcal{A}$ such that, for any $\sigma\in\Sigma_{n}$, $b_{1}, \dotsc   , b_{n}\in f_{C}(A)$ and $a_{1}, \dotsc   , a_{n}\in A$, the function
\[D\sigma_{b_{1}, \dotsc   , b_{n}}^{a_{1}, \dotsc   , a_{n}}:\sigma_{\mathcal{B}}(b_{1}, \dotsc   , b_{n})\rightarrow\sigma_{\mathcal{A}}(a_{1}, \dotsc   , a_{n})\]
satisfies that, if $a \in \sigma_{\mathcal{A}}(a_{1}, \dotsc   , a_{n})$ is such that $C\sigma_{a_{1}, \dotsc  , a_{n}}^{b_{1}, \dotsc  , b_{n}}(a) \in \sigma_{\mathcal{B}}(b_{1}, \dotsc  , b_{n})$, then\\ $D\sigma_{b_{1}, \dotsc   , b_{n}}^{a_{1}, \dotsc   , a_{n}}(C\sigma_{a_{1}, \dotsc   , a_{n}}^{b_{1}, \dotsc   , b_{n}}(a))=a$. Given that $C\sigma_{a_{1}, \dotsc   , a_{n}}^{b_{1}, \dotsc   , b_{n}}$ is injective, this condition is well-defined.

Since $\mathcal{B}$ is $\textbf{cdf}$-generated by $G(\mathcal{B})$, we know to exist a homomorphism $g_{D}:\mathcal{B}\rightarrow \mathcal{A}$ extending $g$ and the functions $D\sigma_{b_{1}, \dotsc   , b_{n}}^{g_{D}(b_{1}), \dotsc   , g_{D}(b_{n})}$.

Finally, we take $g_{D}\circ f_{C}:\mathcal{A}\rightarrow\mathcal{A}$: it extends the injection $id=g\circ f:X\rightarrow A$, for which $id(x)=x$; it also extends the collection of choices $E$ defined by
\[E\sigma_{a_{1}, \dotsc  , a_{n}}^{a'_{1}, \dotsc  , a'_{n}}=D\sigma_{f_{C}(a_{1}), \dotsc  , f_{C}(a_{n})}^{a'_{1}, \dotsc  , a'_{n}}\circ C_{a_{1}, \dotsc  , a_{n}}^{f_{C}(a_{1}), \dotsc  , f_{C}(a_{n})}:\sigma_{\mathcal{A}}(a_{1}, \dotsc  , a_{n})\rightarrow\sigma_{\mathcal{A}}(a'_{1}, \dotsc  , a'_{n}),\]
for $\sigma\in\Sigma_{n}$ and $a_{1}, \dotsc  a_{n}, a'_{1}, \dotsc  , a'_{n}\in A$. 
This way, $E\sigma_{a_{1}, \dotsc  , a_{n}}^{a_{1}, \dotsc  , a_{n}}$ is the identity on $\sigma_{\mathcal{A}}(a_{1}, \dotsc  , a_{n})$: indeed, for any $a\in \sigma_{\mathcal{A}}(a_{1}, \dotsc  , a_{n})$, 
\[C\sigma_{a_{1}, \dotsc  , a_{n}}^{f_{C}(a_{1}), \dotsc  , f_{C}(a_{n})}(a)=f_{C}(a)\]
by definition of $f_{C}$, and, given $f_{C}$ is a homomorphism, $f_{C}(a)\in \sigma_{\mathcal{B}}(f_{C}(a_{1}), \dotsc  , f_{C}(a_{n}))$, meaning\\ $C\sigma_{a_{1}, \dotsc  , a_{n}}^{f_{C}(a_{1}), \dotsc  , f_{C}(a_{n})}(a) \in \sigma_{\mathcal{B}}(f_{C}(a_{1}), \dotsc  , f_{C}(a_{n}))$; then
\[E\sigma_{a_{1}, \dotsc  , a_{n}}^{a_{1}, \dotsc  , a_{n}}(a) = D\sigma_{f_{C}(a_{1}), \dotsc  , f_{C}(a_{n})}^{a_{1}, \dotsc  , a_{n}}(C_{a_{1}, \dotsc  , a_{n}}^{f_{C}(a_{1}), \dotsc  , f_{C}(a_{n})}(a)) = a\]
by definition of $D$.

But notice that the identical homomorphism $Id_{\mathcal{A}}:\mathcal{A}\rightarrow \mathcal{A}$ also extends both $id$ and $E$ and, given the unicity of such extensions on the definition of being $\textbf{cdf}$-generated, we obtain that $Id_{\mathcal{A}}=g_{D}\circ f_{C}$. Of course, the fact that $f_{C}:\mathcal{A}\rightarrow\mathcal{B}$ has a left inverse implies that is injective, and by definition of $\mathcal{B}$ it is also surjective, meaning it is a bijective function; moreover, $g_{D}$ is the inverse function of $f_{C}$. Finally, for $\sigma\in \Sigma_{n}$ and $a_{1}, \dotsc  , a_{n}\in A$,
\[f_{C}(\sigma_{\mathcal{A}}(a_{1}, \dotsc  , a_{n}))\subseteq \sigma_{\mathcal{B}}(f_{C}(a_{1}), \dotsc  , f_{C}(a_{n})),\]
since $f_{C}$ is a homomorphism; however, given $g_{D}$ is also a homomorphism.
\[g_{D}(\sigma_{\mathcal{B}}(f_{C}(a_{1}), \dotsc  , f_{C}(a_{n})))\subseteq \sigma_{\mathcal{A}}(g_{D}\circ f_{C}(a_{1}), \dotsc  , g_{D}\circ f_{C}(a_{n}))=\sigma_{\mathcal{A}}(a_{1}, \dotsc  , a_{n}),\]
and by applying $f_{C}$ to both sides, one obtains
\[\sigma_{\mathcal{B}}(f_{C}(a_{1}), \dotsc  , f_{C}(a_{n}))=f_{C}(g_{D}(\sigma_{\mathcal{B}}(f_{C}(a_{1}), \dotsc  , f_{C}(a_{n}))))\subseteq f_{C}(\sigma_{\mathcal{A}}(a_{1}, \dotsc  , a_{n})),\]
what proves $f_{C}$ is a full homomorphism, that is, an isomorphism.
\end{proof}

Notice that, \textit{mutatis mutandis}, the previous proof shows that if $\mathcal{A}=(A, \{\sigma_{\mathcal{A}}\}_{\sigma\in\Sigma})$ is $\textbf{cdf}$-generated by $X$, then $\mathcal{A}$ is isomorphic to a submultialgebra of $\textbf{mF}(\Sigma, X, M(\mathcal{A}))$, where
\[M(\mathcal{A})=\sup \big\{|\sigma_{\mathcal{A}}(a_{1}, \dotsc   , a_{n})| \ : \  n\in\mathbb{N}, \, \sigma\in\Sigma_{n}, \, a_{1}, \dotsc   , a_{n}\in A \big\}.\]
\label{M(A)}Quite obviously, $M(\textbf{mF}(\Sigma, \mathcal{V}, \kappa))=\kappa$ and $M(\mathcal{A})\leq \kappa$ for any submultialgebra of $\textbf{mF}(\Sigma, \mathcal{V}, \kappa)$. The value $M(\mathcal{A})$ has been regarded for some time, in the literature of multialgebras, as one of their fundamental aspects, see for example ~\cite{CuponaMadarasz}; unfortunately, their definition of homomorphism is grossly different from ours, meaning their results are not applicable to or studies.

Notice, furthermore, that written in classical terms, the previous Theorems~\ref{sub mF is cdf-gen ground} and~\ref{cdf-gen is iso to sub mF} state something quite well known: an algebra is absolutely free if, and only if, it is isomorphic to some algebra of formulas over the same signature.

\begin{corollary}\label{cdf-gen is cdf-gen by ground}
Every $\textbf{cdf}$-generated multialgebra $\mathcal{A}$ is generated by its ground $G(\mathcal{A})$.
\end{corollary}

\begin{proof}
Since every $\textbf{cdf}$-generated multialgebra is isomorphic to a submultialgebra of some\\ $\textbf{mF}(\Sigma, X, \kappa)$, from Theorem \ref{cdf-gen is iso to sub mF}, and every submultialgebra of $\textbf{mF}(\Sigma, X, \kappa)$ is generated by its ground, by Lemma \ref{sub mF is gen ground}, the result follows.
\end{proof}

\begin{corollary}
Every $\textbf{cdf}$-generated multialgebra $\mathcal{A}$ is $\textbf{cdf}$-generated by its ground $G(\mathcal{A})$.
\end{corollary}

\begin{definition}
A $\Sigma$-multialgebra $\mathcal{A}=(A, \{\sigma_{\mathcal{A}}\}_{\sigma\in\Sigma})$ is said to be disconnected\index{Disconnected} if, for every $\sigma\in\Sigma_{n}, \theta\in \Sigma_{m}$, $a_{1}, \dotsc    , a_{n}, b_{1}, \dotsc    , b_{m}\in A$,
\[\sigma_{\mathcal{A}}(a_{1}, \dotsc   , a_{n})\cap \theta_{\mathcal{A}}(b_{1}, \dotsc   , b_{m})\neq\emptyset\]
implies that $n=m$, $\sigma=\theta$ and $a_{1}=b_{1}, \dotsc   , a_{n}=b_{m}$.
\end{definition}

\begin{example}
$\textbf{F}(\Sigma, \mathcal{V})$ is disconnected.
\end{example}

\begin{example}
All directed forests of height $\omega$, when considered as $\Sigma_{s}$-multialgebras, are disconnected, given that no two arrows point to the same element.
\end{example}

It is clear that if $\mathcal{B}$ is a submultialgebra of $\mathcal{A}$ and $\mathcal{A}$ is disconnected, then $\mathcal{B}$ is also disconnected, since if $\sigma_{\mathcal{B}}(a_{1}, \dotsc   , a_{n})\cap \theta_{\mathcal{B}}(b_{1}, \dotsc   , b_{m})\neq\emptyset$, for $a_{1}, \dotsc   , a_{n}, b_{1}, \dotsc   , b_{m}\in B$, given that $\sigma_{\mathcal{B}}(a_{1}, \dotsc   , a_{n})\subseteq\sigma_{\mathcal{A}}(a_{1}, \dotsc   , a_{n})$ and $\theta_{\mathcal{B}}(b_{1}, \dotsc   , b_{m})\subseteq \theta_{\mathcal{A}}(b_{1}, \dotsc   , b_{m})$, we find that 
\[\sigma_{\mathcal{A}}(a_{1}, \dotsc   , a_{n})\cap \theta_{\mathcal{A}}(b_{1}, \dotsc   , b_{m})\neq\emptyset\]
and therefore $n=m$, $\sigma=\theta$ and $a_{1}=b_{1}, \dotsc   , a_{n}=b_{m}$.

We noticed before that $\textbf{mF}(\Sigma, \mathcal{V}, \kappa)$ is disconnected, and by Theorem \ref{cdf-gen is iso to sub mF} we obtain that every $\textbf{cdf}-$generated algebra is disconnected. This has a deeper meaning: being disconnected is an attempt to measure how free of identities a multialgebra is. After all, having no two multioperations to coincide, on any elements, is strongly indicative that the multialgebra does not satisfy any identities. The fact that our candidates for an absolutely free multialgebra (the submultialgebras of $\textbf{mF}(\Sigma, \mathcal{V}, \kappa)$) are all disconnected suggests that they are well-deserving of the ``weakly free'' title.

\subsection{Being disconnected and generated by the ground}

We have offered, so far, two characterizations of the multialgebras we chose to call weakly free: first of all, they are submultialgebras of some $\textbf{mF}(\Sigma, \mathcal{V}, \kappa)$, and perhaps this is to be taken as their definition; second, they are $\textbf{cdf}$-generated. Now, we look at other possible characterizations of being weakly free that could lead to possible future definitions of relatively free multialgebras. 

Algebras of formulas satisfy no identities, what would partially correspond here to the concept of being disconnected. However, there is one property that is maybe more representative of our intuition of formulas (which are, up to isomorphism, the elements of all absolutely free algebras): whenever one deals with formulas, one starts by defining them from elements that are as simple as possible (variables), and continues indefinitely by combining them trough operations (connectives). 

The concept of simplest (or indecomposable) element, here, is replaced by that of being an element of the ground, so one would expect that being generated by it plays some sort of role in the objects we have defined so far: a multialgebra which is generated by a set has all of its elements either on the set, or as the result of increasingly more complex (multi)operations performed on that very set.

\begin{lemma}\label{X in ground if cdf-gen by X}
If $\mathcal{A}$ is $\textbf{cdf}$-generated by $X$, then $X\subseteq G(\mathcal{A})$.
\end{lemma}

\begin{proof}
If $\mathcal{A}$ is $\textbf{cdf}$-generated by $X$, then $\mathcal{A}$ is isomorphic to a submultialgebra of $\textbf{mF}(\Sigma, X, |A|)$ containing $X$, from Theorem \ref{cdf-gen is iso to sub mF}: let us assume that $\mathcal{A}$ is equal to this submultialgebra, without loss of generality.

Then we have $X=G(\textbf{mF}(\Sigma, X, |A|))\cap A\subseteq G(\mathcal{A})$.
\end{proof}

\begin{lemma}\label{no proper cdf-gen sets}
If $\mathcal{A}$ is $\textbf{cdf}$-generated by both $X$ and $Y$, with $X\subseteq Y$, then $X=Y$.
\end{lemma}

\begin{proof}
Suppose $X\neq Y$ and let $y\in Y\setminus X$: take a $\Sigma$-multialgebra $\mathcal{B}$ over the same signature as $\mathcal{A}$ such that $|B|\geq 2$, and a collection of choices $C$ from $\mathcal{A}$ to $\mathcal{B}$.

Take also two functions $g, h:Y\rightarrow B$ such that $g|_{X}=h|_{X}$ and $g(y)\neq h(y)$, which is possible since $|B|\geq 2$: given that $\mathcal{A}$ is $\textbf{cdf}$-generated by $Y$, there exist unique homomorphisms $g_{C}$ and $h_{C}$ extending both, respectively, $g$ and $C$ and $h$ and $C$.

However, $g_{C}$ and $h_{C}$ extend both $g|_{X}:X\rightarrow B$ and $C$, and since $\mathcal{A}$ is $\textbf{cdf}$-generated by $X$, we find that $g_{C}=h_{C}$. This is not possible, since $g_{C}(y)\neq h_{C}(y)$, what must imply that $Y\setminus X=\emptyset$ and therefore $X=Y$.
\end{proof}

\begin{theorem}\label{ground is only cdf-gen set}
Every $\textbf{cdf}$-generated multialgebra $\mathcal{A}$ is uniquely $\textbf{cdf}$-generated by its ground.
\end{theorem}

\begin{proof}
From Theorem \ref{cdf-gen is cdf-gen by ground}, $\mathcal{A}$ is $\textbf{cdf}$-generated by $G(\mathcal{A})$, and from Lemma \ref{X in ground if cdf-gen by X}, if $\mathcal{A}$ is also $\textbf{cdf}$-generated by $X$, then $X\subseteq G(\mathcal{A})$. By Lemma \ref{no proper cdf-gen sets}, this implies that $X=G(\mathcal{A})$.
\end{proof}

We have proved so far that, if $\mathcal{A}$ is $\textbf{cdf}$-generated, then $\mathcal{A}$ is generated by its ground and disconnected. We would like to prove that this is enough to characterize a $\textbf{cdf}$-generated multialgebra: that is, if $\mathcal{A}$ is generated by its ground and disconnected, then it is $\textbf{cdf}$-generated, exactly by its ground.

The idea is similar to the one we used to prove, in Theorem \ref{sub mF is cdf-gen ground}, that all submultialgebras of $\textbf{mF}(\Sigma, \mathcal{V}, \kappa)$ are $\textbf{cdf}$-generated: take a multialgebra $\mathcal{A}$ that is both generated by its ground $G(\mathcal{A})$, which we will denote by $X$, and disconnected, and fix a multialgebra $\mathcal{B}$ over the same signature, a function $f:X\rightarrow B$ and a collection of choices $C$ from $\mathcal{A}$ to $\mathcal{B}$.

We define a function $f_{C}:A\rightarrow B$ using induction on the $\langle X\rangle_{m}$: for $m=0$, either we have an element $x\in X$, when we define $f_{C}(x)=f(x)$, or we have an $a\in\sigma_{\mathcal{A}}$, for a $\sigma\in\Sigma_{0}$, when we define $f_{C}(a)=C\sigma(a)$. Notice how, up to this point, there are no contradictions on the definition, given an element cannot belong both to $X$ and to a $\sigma_{\mathcal{A}}$, since $X=G(\mathcal{A})$.

Suppose we have successfully defined $f_{C}$ on $\langle X\rangle_{m}$ and take an $a\in \sigma_{\mathcal{A}}(a_{1}, \dotsc   , a_{n})$ for $a_{1}, \dotsc   , a_{n}\in\langle X\rangle_{m}$. We then define 
\[f_{C}(a)=C\sigma_{a_{1}, \dotsc   , a_{n}}^{f_{C}(a_{1}), \dotsc   , f_{C}(a_{n})}(a).\]
Again the function remains well-defined: $a$ cannot belong to $X$, since $X=G(\mathcal{A})$, and cannot belong to a $\theta_{\mathcal{A}}(b_{1}, \dotsc   , b_{p})$ unless $p=n$, $\theta=\sigma$ and $b_{1}=a_{1}, \dotsc    ,b_{p}=a_{n}$, since $\mathcal{A}$ is disconnected.

Clearly $f_{C}$ is a homomorphism, since the image of $\sigma_{\mathcal{A}}(a_{1}, \dotsc   , a_{n})$ under $f_{C}$ is contained in $\sigma_{\mathcal{B}}(f_{C}(a_{1}), \dotsc   , f_{C}(a_{n}))$, and $f_{C}$ extends both $f$ and $C$.

\begin{lemma}\label{gen by ground and disc implies cdf-gen}
If a multialgebra $\mathcal{A}$ is both generated by its ground $X$ and disconnected, $\mathcal{A}$ is $\textbf{cdf}$-generated by $X$.
\end{lemma}

\begin{proof}
It remains for us to show that $f_{C}$, as defined above, is the only homomorphism extending $f$ and $C$. Suppose $g$ is another such homomorphism and we shall proceed yet again by induction.

On $\langle X\rangle_{0}$, we have that $f_{C}(x)=f(x)=g(x)$ for all $x\in X$; and for $a\in\sigma_{\mathcal{A}}$ and $\sigma\in\Sigma_{0}$ we have that 
\[f_{C}(a)=C\sigma(a)=g(a),\]
hence $f_{C}$ and $g$ coincide on $\langle X\rangle_{0}$. Suppose that $f_{C}$ and $g$ are equal on $\langle X\rangle_{m}$ and take $a\in\sigma_{\mathcal{A}}(a_{1}, \dotsc   , a_{n})$ for $a_{1}, \dotsc   , a_{n}\in \langle X\rangle_{m}$: we have by induction hypothesis that
\[f_{C}(a)=C\sigma_{a_{1}, \dotsc   , a_{n}}^{f_{C}(a_{1}), \dotsc   , f_{C}(a_{n})}(a)=C\sigma_{a_{1}, \dotsc   , a_{n}}^{g(a_{1}), \dotsc   , g(a_{n})}(a)=g(a),\]
which concludes our proof.
\end{proof}

\begin{theorem}\label{cdf-gen iff gen by ground and disc}
A multialgebra $\mathcal{A}$ is $\textbf{cdf}$-generated if, and only if, $\mathcal{A}$ is generated by its ground and disconnected.
\end{theorem}

We have introduced several concepts that, although dissimilar in their definitions, are intrinsically connected; so it becomes important to analyze whether they are indeed distinct: are there multialgebras that are disconnected but not generated by their grounds? Are there multialgebras that are generated by their grounds but not disconnected? Or does being generated by its ground implies being disconnected, or vice-versa? We show below that this is not the case, for we provide examples answering positively both previous questions.

\begin{example}\label{C}
Take the signature $\Sigma_{s}$ with a single unary operator, first defined in Example \ref{s}. Consider the $\Sigma_{s}-$multialgebra $\mathcal{C}=(\{-1,1\}, \{s_{\mathcal{C}}\})$ such that $s_{\mathcal{C}}(-1)=\{1\}$ and $s_{\mathcal{C}}(1)=\{-1\}$ (that is, $s_{\mathcal{C}}(x)=\{-x\}$).

We state that $\mathcal{C}$ is disconnected, but not generated by its ground. $\mathcal{C}$ is clearly disconnected since $s_{\mathcal{C}}(-1)\cap s_{\mathcal{C}}(1)=\emptyset$; now, $B(\mathcal{C})=s_{\mathcal{C}}(-1)\cup s_{\mathcal{C}}(1)=\{-1,1\}$, and so $G(\mathcal{C})=\emptyset$. Since $\Sigma_{0}=\emptyset$, $\bigcup_{\sigma\in\Sigma_{0}}\sigma_{\mathcal{C}}=\emptyset$ and therefore $\langle G(\mathcal{C})\rangle_{n}=\emptyset$ for every $n\in\mathbb{N}$, so that $G(\mathcal{C})$ does not generate $\mathcal{C}$.
\begin{figure}[H]
\centering
\begin{tikzcd}
    -1 \arrow[rr, bend left=50, "s_{\mathcal{C}}"]  && 1 \arrow[ll, bend left=50, "s_{\mathcal{C}}"]
  \end{tikzcd}
\caption*{The $\Sigma_{s}$-multialgebra $\mathcal{C}$}
\end{figure}
\end{example}

\begin{example}\label{B}
Take the signature $\Sigma_{s}$ with a single unary operator. Consider the $\Sigma_{s}-$multialge\-bra $\mathcal{B}=(\{0, 1\}, \{s_{\mathcal{B}}\})$ such that $s_{\mathcal{B}}(0)=\{1\}$ and $s_{\mathcal{B}}(1)=\{1\}$ (that is, $s_{\mathcal{B}}(x)=\{1\}$).

Then $\mathcal{B}$ is clearly not disconnected, since $s_{\mathcal{B}}(0)\cap s_{\mathcal{B}}(1)=\{1\}$, yet $\mathcal{B}$ is generated by its ground: $B(\mathcal{B})=\{1\}$ and so $G(\mathcal{B})=\{0\}$, and we see that $\langle G(\mathcal{B})\rangle_{1}$ is already $\{0,1\}$.
\begin{figure}[H]
\centering
\begin{tikzcd}
    0 \arrow[r, "s_{\mathcal{B}}"]  & 1 \arrow[loop right, out=30, in=-30, distance=3em]{}{s_{\mathcal{B}}}
  \end{tikzcd}
\caption*{The $\Sigma_{s}$-multialgebra $\mathcal{B}$}
\end{figure}

\end{example}

\subsection{Being disconnected and having a strong basis}

Now, we define the notion of a strong basis (a minimum generating set), and prove that, on one hand, being generated by the ground implies having a strong basis, what means that being disconnected and generated by the ground implies being disconnected and having a strong basis; reciprocally, we also prove having a strong basis and being disconnected implies being disconnected and generated by the ground (although having a strong basis does not imply being generated by the ground). This provides a third characterization of our weakly free multialgebras.

Our motivation, when coining the definition of a strong basis, was to be able to weaken that very condition: after all, absolutely free algebras (that is, algebras freely generated on the variety of all algebras on a given signature) are easier to define than the relatively free ones (which are the algebras freely generated on any variety one wants to consider), so it is natural that we start this study with ``absolutely free'' multialgebras. However, we still would like to be able, in the future, to define what should be a relatively free multialgebra (whatever a variety of multialgebras may be); to this end, weakening a strong basis to be a minimal (instead of minimum) generating set makes sense, given many relatively free algebras, on domains such as that of vector spaces, indeed have basis.

\begin{definition}
We say $B\subseteq A$ is a strong basis\index{Strong basis} of the $\Sigma$-multialgebra $\mathcal{A}=(A, \{\sigma_{\mathcal{A}}\}_{\sigma\in\Sigma})$ if it is the minimum of the set $\mathcal{G}=\{S\subseteq A\ : \  \langle S\rangle=A\}$ when ordered by inclusion.
\end{definition}

\begin{example}
The set of variables $\mathcal{V}$ is a strong basis of $\textbf{F}(\Sigma, \mathcal{V})$.
\end{example}

\begin{example}
The set of elements without predecessor of a directed forest of height $\omega$ is a strong basis of the forest, when considered as a $\Sigma_{s}$-multialgebra.
\end{example}

\begin{lemma}\label{ground cap span is in set}
For every subset $S$ of the universe of a $\Sigma$-multialgebra $\mathcal{A}$, $G(\mathcal{A})\cap \langle S\rangle\subseteq S$.
\end{lemma}

\begin{proof}
Suppose $x\in G(\mathcal{A})\cap \langle S\rangle$: if $x\notin S$, we will show that $x$ cannot be in $\langle S\rangle$, which contradicts our assumption. Indeed, if $x\notin S$ then 
\[x\notin\langle S\rangle_{0}=S\cup \bigcup_{\sigma\in\Sigma_{0}}\sigma_{\mathcal{A}},\]
since $x\notin S$, and $x\in G(\mathcal{A})$ implies that 
\[x\in A\setminus B(\mathcal{A}) \subseteq A\setminus \bigcup_{\sigma\in\Sigma_{0}}\sigma_{\mathcal{A}}.\]

Now, for induction hypothesis, suppose that $x\notin\langle S\rangle_{m}$; then
\[x\notin \langle S\rangle_{m+1}=\langle S\rangle_{m}\cup\bigcup \big\{\sigma_{\mathcal{A}}(a_{1}, \dotsc   , a_{n}) \ : \  n\in\mathbb{N}, \, \sigma\in\Sigma_{n}, \, a_{1}, \dotsc   , a_{n}\in \langle S\rangle_{m} \big\}\]
since $x\notin \langle S\rangle_{m}$, and $x\in G(\mathcal{A})$ implies that 
\[x\in A\setminus B(\mathcal{A}) \subseteq A\setminus \bigcup \big\{\sigma_{\mathcal{A}}(a_{1}, \dotsc   , a_{n}) \ : \  n\in\mathbb{N}, \, \sigma\in\Sigma_{n}, \, a_{1}, \dotsc   , a_{n}\in \langle S\rangle_{m} \big\}.\]
\end{proof}

\begin{theorem}\label{ground is in basis}
If the $\Sigma$-multialgebra $\mathcal{A}$ has a strong basis $B$, $G(\mathcal{A})\subseteq B$.
\end{theorem}

\begin{proof}
By lemma \ref{ground cap span is in set}, $G(\mathcal{A})=G(\mathcal{A})\cap A=G(\mathcal{A})\cap\langle B\rangle\subseteq B$.
\end{proof}

Notice lemma \ref{ground cap span is in set} leads us too to the fact that, if $\mathcal{A}$ is generated by its ground, then it has the ground as a strong basis: this is because, if $\langle S\rangle=A$, $G(\mathcal{A})=G(\mathcal{A})\cap\langle S\rangle\subseteq S$, and therefore $G(\mathcal{A})$ becomes a minimum generating set.

\begin{definition}
If $B$ is a strong basis of a disconnected $\Sigma$-multialgebra $\mathcal{A}$, we define the $B$-order\index{Order, $B$-}\label{Border} of an element $a\in A$ as the natural number 
\[o_{B}(a)=\min \big\{k\in\mathbb{N}\ : \  a\in\langle B\rangle_{k}\big\}.\]
\end{definition}

This is a clear generalization of the order, or complexity, of a formula: in fact, the order of a formula in ${T}(\Sigma, \mathcal{V})$ is exactly its $\mathcal{V}$-order.

\begin{proposition}
If $a\in \sigma_{\mathcal{A}}(a_{1}, \dotsc   , a_{n})$ and $o_{B}(a)\geq 1$, then $o_{B}(a_{1}), \dotsc   , o_{B}(a_{n})<o_{B}(a)$.
\end{proposition}

\begin{proof}
Suppose $o_{B}(a)$ equals $m+1$, with $m\in\mathbb{N}$, implying that
\[a\in \langle B\rangle_{m+1}=\langle B\rangle_{m}\cup \bigcup \big\{\sigma_{\mathcal{A}}(a_{1}, \dotsc   , a_{n}) \ : \  n\in\mathbb{N}, \, \sigma\in\Sigma_{n}, \, a_{1}, \dotsc   , a_{n}\in \langle B\rangle_{m} \big\};\]
since $m+1=\min\{k\in\mathbb{N}\ : \  a\in\langle B\rangle_{k}\}$, we have that $a\notin \langle B\rangle_{m}$ and therefore
\[a\in \bigcup \big\{\sigma_{\mathcal{A}}(a_{1}, \dotsc   , a_{n}) \ : \  n\in\mathbb{N}, \, \sigma\in\Sigma_{n}, \, a_{1}, \dotsc   , a_{n}\in \langle B\rangle_{m} \big\}.\]
Finally, we obtain that there exist $p\in\mathbb{N}$, $\theta\in\Sigma_{p}$ and $b_{1}, \dotsc   , b_{p}\in \langle B\rangle_{m}$ such that $a\in \theta_{\mathcal{A}}(b_{1}, \dotsc   , b_{p})$. Since $a\in \sigma_{\mathcal{A}}(a_{1}, \dotsc   , a_{n})$, this implies that  $\sigma_{\mathcal{A}}(a_{1}, \dotsc   , a_{n})\cap\theta_{\mathcal{A}}(b_{1}, \dotsc   , b_{p})\neq\emptyset$, and therefore $p=n$, $\theta=\sigma$ and $b_{1}=a_{1}, \dotsc , b_{p}=a_{n}$, so that $o_{B}(a_{1}), \dotsc   , o_{B}(a_{n})\leq m$.
\end{proof}

But what if $a\in \sigma_{\mathcal{A}}(a_{1}, \dotsc   , a_{n})$ and $o_{B}(a)=0$, implying $a\in B$? We state this case cannot occur, for if it does, 
\[B^{*}=\big(B\cup\{a_{1}, \dotsc   , a_{n}\}\big)\setminus\{a\}\]
generates $A$, while clearly not containing $B$, even in the case where $n=0$. We have that $a\in \langle B^{*}\rangle_{1}$, since $a_{1}, \dotsc   , a_{n}\in B^{*}$ and $a\in\sigma_{\mathcal{A}}(a_{1}, \dotsc   , a_{n})$, and given that $B\setminus\{a\}\subseteq B^{*}$, it follows that $B\subseteq \langle B^{*}\rangle_{1}$, and so $\langle B\rangle_{0}\subseteq\langle B^{*}\rangle_{1}$. 

It is then true that, for every $m\in\mathbb{N}$, $\langle B\rangle_{m}\subseteq \langle B^{*}\rangle_{m+1}$: if this is true for $m$, let $b\in \langle B\rangle_{m+1}$, and then either $b\in \langle B\rangle_{m}$, so that $b\in \langle B^{*}\rangle_{m+1}\subseteq \langle B^{*}\rangle_{m+2}$, or there exist $\theta\in\Sigma_{p}$ and $b_{1}, \dotsc   , b_{p}\in \langle B\rangle_{m}$ such that $b\in\theta_{\mathcal{A}}(b_{1}, \dotsc   , b_{p})$; in this case, since $\langle B\rangle_{m}\subseteq\langle B^{*}\rangle_{m+1}$, we have that
\[b\in \theta_{\mathcal{A}}(b_{1}, \dotsc   , b_{p})\subseteq \bigcup \big\{\sigma_{\mathcal{A}}(a_{1}, \dotsc   , a_{n}) \ : \  n\in\mathbb{N}, \, \sigma\in\Sigma_{n}, \, a_{1}, \dotsc   , a_{n}\in \langle B^{*}\rangle_{m+1} \big\} \subseteq\langle B^{*}\rangle_{m+2},\]
so once again $b\in\langle B^{*}\rangle_{m+2}$. Since $\langle B\rangle=\bigcup_{m\in\mathbb{N}}\langle B\rangle_{m}$ equals $A$, we have that $\langle B^{*}\rangle$ also equals $A$, as we previously stated. This is absurd, since $B$ is the minimum of $\{S\subseteq A\ : \  \langle S\rangle=A\}$ ordered by inclusion and $B\not\subseteq B^{*}$. The conclusion must be that if $a\in\sigma_{\mathcal{A}}(a_{1}, \dotsc   , a_{n})$, then $o_{B}(a_{1}), \dotsc   , o_{B}(a_{n})<o_{B}(a)$, regardless of the value of $o_{B}(a)$.

\begin{lemma}\label{ground is basis if disc}
If $\mathcal{A}$ is disconnected and has a strong basis $B$, $B=G(\mathcal{A})$ and so $\mathcal{A}$ is generated by its ground.
\end{lemma}

\begin{proof}
Suppose $a\in B\setminus G(\mathcal{A})$: since $a$ is in the build of $\mathcal{A}$, there exist $\sigma\in\Sigma_{n}$ and elements $a_{1}, \dotsc   , a_{n}\in A$ such that $a\in\sigma_{\mathcal{A}}(a_{1}, \dotsc   , a_{n})$. If $n>0$, $o_{B}(a)>o_{B}(a_{1})\geq 0$, which contradicts the fact that $a\in B$ and therefore $o_{B}(a)=0$.

If $n=0$, it is clear that $B^{*}=B\setminus \{a\}$ is a generating set smaller than $B$: generating set because, if $a\in \sigma_{\mathcal{A}}$, $a\in \bigcup_{\sigma\in\Sigma_{0}}\sigma_{\mathcal{A}}$ and therefore $B\subseteq \langle B^{*}\rangle_{0}$, so that $\langle B\rangle_{m}\subseteq\langle B^{*}\rangle_{m+1}$ and $\bigcup_{m\in\mathbb{N}}\langle B\rangle_{m}=\bigcup_{m\in\mathbb{N}}\langle B^{*}\rangle_{m}$. This is also a contradiction, since $B$ is a strong basis.
\end{proof}

\begin{theorem}\label{theorem 3}
$\mathcal{A}$ is generated by its ground and disconnected if, and only if, it has a strong basis and it is disconnected.
\end{theorem}

\begin{proof}
We already proved, in Lemma \ref{ground is basis if disc}, that if $\mathcal{A}$ is disconnected and has a strong basis $B$, then it is generated by its ground and disconnected. Reciprocally, if $\mathcal{A}$ is disconnected and generated by its ground, first of all it is clearly disconnected.

Now, if $\langle G(\mathcal{A})\rangle=A$ one has that $G(\mathcal{A})\subseteq S$ for every $S\in \{S\subseteq A\ : \  \langle S\rangle=A\}$, by Lemma \ref{ground cap span is in set}. Therefore, the ground is a strong basis.
\end{proof}

Once again, we ask ourselves whether the concepts we have defined in this section are truly independent: does being disconnected imply having a strong basis? Does having a strong basis imply being disconnected? We show that neither is the case by providing examples of a multialgebra that is disconnected but does not have a strong basis and one of a multialgebra that has a strong basis but is not disconnected.

\begin{example}
Take the $\Sigma_{s}-$multialgebra $\mathcal{C}$ from Example \ref{C}.

We know that $\mathcal{C}$ is disconnected, but we also state that it does not have a strong basis: in fact, we see that the set $\{S\subseteq \{-1,1\}\ :\  \langle S\rangle=\{-1,1\}\}$ is exactly $\{\{-1\}, \{1\}, \{-1,1\}\}$, and this set has no minimum.
\end{example}

\begin{example}
Take the $\Sigma_{s}-$multialgebra $\mathcal{B}$ from Example \ref{B}.

As we saw before, $\mathcal{B}$ is not disconnected, but we state that it has a strong basis: $B=\{0\}$ generates $\mathcal{B}$ and, since $\{1\}$ does not generate the multialgebra, we find that $B$ is a minimum generating set.
\end{example}

In these two examples we presented a multialgebra ($\mathcal{B}$) which has a strong basis and is generated by its ground, and one multialgebra ($\mathcal{C}$) which does not have a strong basis and is not generated by its ground. Clearly being generated by its ground implies having a strong basis, so it is natural to hypothesize that having a strong basis and being generated by its ground could be equivalent notions; however, as we show in the example below, having a basis does not imply being generated by its ground.

\begin{example}
Take the signature $\Sigma_{s}$ with a single unary operator from Example \ref{s}. Consider the $\Sigma_{s}-$multialgebra $\mathcal{M}=(\{-1, 0, 1\}, \{s_{\mathcal{M}}\})$ such that $s_{\mathcal{M}}(0)=\{0\}$, $s_{\mathcal{M}}(1)=\{1\}$ and $s_{\mathcal{M}}(-1)=\{1\}$ (that is, $s_{\mathcal{M}}(x)=\{|x|\}$, where $|X|$ denotes the absolute value of $x$).

We have that $G(\mathcal{M})=\{-1\}$ and that $\langle\{-1\}\rangle=\{-1, 1\}$, so that $\mathcal{M}$ is not generated by its ground. But we state that $\{-1, 0\}$ is a strong basis: first of all, it clearly generates $\mathcal{M}$; furthermore, the generating sets of $\mathcal{M}$ are only $\{-1, 0\}$ and $\{-1, 0, 1\}$, so that $\{-1, 0\}$ is in fact the smallest generating set.

\begin{figure}[H]
\centering
\begin{tikzcd}
    -1 \arrow[rr, bend right=50]{}{s_{\mathcal{M}}}  & 0 \arrow[u, loop]{}{s_{\mathcal{M}}} & 1 \arrow[loop right, out=30, in=-30, distance=3em]{}{s_{\mathcal{M}}}
  \end{tikzcd}
\caption*{The $\Sigma_{s}$-multialgebra $\mathcal{M}$}
\end{figure}
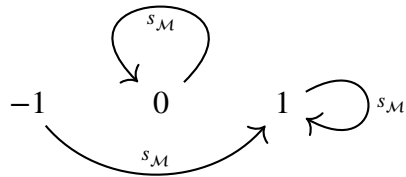
\end{example}

\subsection{Being disconnected and chainless}

The last equivalence to being a submultialgebra of $\textbf{mF}(\Sigma, \mathcal{V}, \kappa)$ we give depends on the notion of being chainless, which is very graph-theoretical in nature. Think of a tree that ramifies ever downward: one can pick any vertex and proceed, against the arrows, upwards until an element without predecessor is reached. More than that, it is not possible to find an infinite path, starting in any one vertex, by always going against the arrows: such a path, if it existed, would be what we shall call a chain. A multialgebra without chains is, very naturally, chainless.

As it was in the case of strong bases, there isn't a parallel concept to being chainless among the theory of universal algebra: it seems that this concept is far more natural when dealing with multioperations, although it can be easily applied to algebras if one wishes to do so. A close, although not equivalent, entity are the branches in the formation trees of formulas: if allowed to grow infinitely, these would became chains.

The main result of this section is, perhaps, the fact that being chainless implies being generated by its ground, which, as we know, implies in turn having a strong basis. Of course, when adding disconnectedness to the equation, all three concepts become equivalent to one another and to the fact that the multialgebra at hand is weakly free. We have, then, the following schematic diagram, that offers an overview of the aforementioned results.

\begin{figure}[H]
\centering
\begin{minipage}[t]{3cm}
\centering
\begin{tikzcd}[arrows=Rightarrow, every arrow/.append style={shift right=3ex}]
\text{$\mathcal{A}$ is chainless} \arrow[Rightarrow]{d}\\
\text{$\mathcal{A}$ is generated by $G(\mathcal{A})$} \arrow[Rightarrow]{d}\arrow[negated]{u}\\
\text{$\mathcal{A}$ has a strong basis} \arrow[negated]{u}
  \end{tikzcd}
\end{minipage}
\hspace{2cm}
\centering
\begin{minipage}[t]{6cm}
\centering
\begin{tikzcd}[arrows=Rightarrow, every arrow/.append style={shift right=3ex}]
\begin{tabular}{c}$\mathcal{A}$ is chainless\\ and disconnected\end{tabular} \arrow[Rightarrow]{d}\\
\begin{tabular}{c}$\mathcal{A}$ is generated by $G(\mathcal{A})$\\ and disconnected\end{tabular} \arrow[Rightarrow]{d}\arrow{u}\\
\begin{tabular}{c}$\mathcal{A}$ has a strong basis\\ and is disconnected\end{tabular} \arrow{u}
  \end{tikzcd}
\end{minipage}
\end{figure}

Given a function $\tau:\{1, \dotsc   , n\}\rightarrow\{1, \dotsc   , n\}$ in $S_{n}$, the group of permutations on $n$ elements (meaning $\tau$ is bijective), the action of $\tau$ in an $n$-tuple $(x_{1}, \dotsc   , x_{n})\in X^{n}$ is given by
\[\tau(x_{1}, \dotsc   , x_{n})=(x_{\tau(1)}, \dotsc   , x_{\tau(n)}).\]
Given $1\leq i, j\leq n$, we define $[i,j]$ to be the permutation such that $[i,j](i)=j$, $[i,j](j)=i$ and, for $k\in\{1, \dotsc   , n\}$ different from $i$ and $j$, $[i,j](k)=k$.

\begin{definition}
Given a $\Sigma$-multialgebra $\mathcal{A}$, a sequence $\{a_{n}\}_{n\in\mathbb{N}}\subseteq A$ is said to be a chain\index{Chain} if, for every $n\in\mathbb{N}$, there exist a natural number $m_{n}\in\mathbb{N}$, a functional symbol $\sigma^{n}\in\Sigma_{m_{n}}$, a permutation $\tau_{n}\in S_{m_{n}}$ and elements $a_{1}^{n}, \dotsc   , a_{m_{n}-1}^{n}\in A$ such that 
\[a_{n}\in \sigma^{n}_{\mathcal{A}}(\tau_{n}(a_{n+1}, a_{1}^{n}, \dotsc   , a_{m_{n}-1}^{n})).\]

A $\Sigma$-multialgebra is said to be chainless\index{Chainless} when it has no chains.
\end{definition}

\begin{example}
Take a directed forest of height $\omega$ and add to it a loop, that is, choose a vertex $v$ and add an arrow from $v$ to $v$: then $\{a_{n}\}_{n\in\mathbb{N}}$, such that $a_{n}=v$ for every $n\in\mathbb{N}$, is a chain.
\end{example}

\begin{example}
$\textbf{F}(\Sigma, \mathcal{V})$ is chainless.
\end{example}

\begin{lemma}\label{chainless implies gen by ground}
If $\mathcal{A}$ is chainless, then it is generated by its ground.
\end{lemma}

\begin{proof}
Suppose that this not hold, so $A\setminus\langle G(\mathcal{A})\rangle$ is not empty, and must therefore contain some element $a_{0}$. We create a chain $\{a_{n}\}_{n\in\mathbb{N}}$ by induction, being the case $n=0$ already done.

So, suppose we have created a finite sequence of elements $a_{0}, \dotsc   , a_{k}\in A\setminus \langle G(\mathcal{A})\rangle$ such that, for each $0\leq n< k$, there exist a positive integer $m_{n}\in\mathbb{N}\setminus\{0\}$, a functional symbol $\sigma^{n}\in\Sigma_{m_{n}}$, a permutation $\tau_{n}\in S_{m_{n}}$ and elements $a_{1}^{n}, \dotsc   , a_{m_{n}-1}^{n}\in A$ such that 
\[a_{n}\in \sigma^{n}_{\mathcal{A}}(\tau_{n}(a_{n+1}, a_{1}^{n}, \dotsc   , a_{m_{n}-1}^{n})).\]

Since $a_{k}\in A\setminus \langle G(\mathcal{A})\rangle$, we have that $a_{k}$ is not an element of the ground; so, there must exist $m_{k}\in\mathbb{N}$, a functional symbol $\sigma^{k}\in\Sigma_{m_{k}}$ and elements $b_{1}^{k}, \dotsc   , b_{m_{k}}^{k}\in A$ such that
\[a_{k}\in\sigma^{k}_{\mathcal{A}}(b_{1}^{k}, \dotsc   , b_{m_{k}}^{k}).\]
Now, if all $b_{1}^{k}, \dotsc   , b_{m_{k}}^{k}$ belonged to $\langle G(\mathcal{A})\rangle$, so would $a_{k}$: there must be an element $a_{k+1}\in\{b_{1}^{k}, \dotsc   , b_{m_{k}}^{k}\}$, say $b_{l}^{k}$, such that $a_{k+1}\in A\setminus \langle G(\mathcal{A})\rangle$. We then define $a_{i}^{k}$ as $b_{j}^{k}$, for $i\in\{1, \dotsc   , m_{k}-1\}$ and $j=\min\{i\leq p\leq m_{k}\ : \  p\neq l\}$, and 
\[\tau_{k}=[l-1, l]\circ\cdots\circ[1,2],\]
and then it is clear that $\{a_{n}\}_{n\in\mathbb{N}}$ becomes a chain, with the extra condition that $\{a_{n}\}_{n\in\mathbb{N}}\subseteq A\setminus \langle G(\mathcal{A})\rangle$. Since $\mathcal{A}$ is chainless, we reach a contradiction, so we must have instead that $A\setminus \langle G(\mathcal{A})\rangle=\emptyset$, and therefore $\mathcal{A}$ is generated by its ground.
\end{proof}

So, a disconnected, chainless multialgebra is, by Lemma \ref{chainless implies gen by ground}, disconnected and generated by its ground. We state, that, in fact, the reciprocal also holds, when we arrive at yet another characterization of being an weakly free multialgebra.

So, suppose $\mathcal{A}$ is disconnected and generated by its ground, and let $\{a_{n}\}_{n\in\mathbb{N}}$ be a chain in $\mathcal{A}$: clearly no $a_{n}$ can belong to the ground, since 
\[a_{n}\in \sigma^{n}_{\mathcal{A}}(\tau_{n}(a_{n+1}, a_{1}^{n}, \dotsc   , a_{m_{n}-1}^{n})),\]
and therefore $o_{G(\mathcal{A})}(a_{n+1})<o_{G(\mathcal{A})}(a_{n})$, that is, the $G(\mathcal{A})$-order of $a_{n+1}$ is less than the $G(\mathcal{A})$-order of $a_{n}$; we reach a contradiction, since if $o_{G(\mathcal{A})}(a_{0})=m$, then $o_{G(\mathcal{A})}(a_{m+1})<0$, what is impossible. $\mathcal{A}$ must then be chainless.

\begin{theorem}\label{generated by ground and disc. iff chainless and disc.}
$\mathcal{A}$ is generated by its ground and disconnected if, and only if, it is chainless and disconnected.
\end{theorem}

Finally, Theorems \ref{sub mF is cdf-gen ground}, \ref{cdf-gen is iso to sub mF}, \ref{cdf-gen iff gen by ground and disc}, \ref{theorem 3} and \ref {generated by ground and disc. iff chainless and disc.} can be summarized as follows; the object we decided to call a weakly free multialgebra is then a multialgebra which satisfies any, and therefore all, of the conditions found in the following theorem.

\begin{theorem}
Are equivalent:

\begin{enumerate} \item $\mathcal{A}$ is a submultialgebra of some $\textbf{mF}(\Sigma, \mathcal{V}, \kappa)$;

\item $\mathcal{A}$ is $\textbf{cdf}$-generated;

\item $\mathcal{A}$ is generated by its ground and disconnected;

\item $\mathcal{A}$ has a strong basis and is disconnected;

\item $\mathcal{A}$ is chainless and disconnected.
\end{enumerate}
\end{theorem}

An important point to stress is that, although not all concepts present in the previous theorem have natural counterparts in universal algebra, by defining them for algebras, presented as multialgebras, we find that all of the conditions in the theorem are valid only for algebras of formulas. This follows easily from the fact that the only $\textbf{cdf}$-generated algebras are the algebras of formulas themselves.

Now, a few examples concerning being chainless, disconnected, having a strong basis and being generated by the ground will be given. Essentially, they show that the notion of being chainless is independent from that of being disconnected, and that despite being chainless implies being generated by the ground (and so having a strong basis), the reciprocal does not hold, and there exist multialgebras generated by the ground which are not chainless.

\begin{example}
Take the signature $\Sigma_{s}$ from Example \ref{s}, and consider the $\Sigma_{s}-$multialgebra $\mathcal{Y}=(\mathbb{N}\cup\{a,b\}, \{s_{\mathcal{Y}}\})$ such that $s_{\mathcal{Y}}(n)=\{n+1\}$, for $n\in\mathbb{N}$, and $s_{\mathcal{Y}}(a)=s_{\mathcal{Y}}(b)=\{0\}$.

We see that $\mathcal{Y}$ is chainless since, given a chain $\{a_{n}\}_{n\in\mathbb{N}}$, it must be contained in the build of $\mathcal{Y}$, that is, $\mathbb{N}$: but then $a_{n+1}=a_{n}-1$, what is a contradiction, since there are only a finite number of elements with index smaller than that of $a_{0}$. At the same time, $\mathcal{Y}$ is not disconnected, since $s_{\mathcal{Y}}(a)=s_{\mathcal{Y}}(b)$.

\begin{figure}[H]
\centering
\begin{tikzcd}
a\arrow[dr]{}{s_{\mathcal{Y}}} & & & \\
    & 0\arrow[r]{}{s_{\mathcal{Y}}} & 1\arrow[r]{}{s_{\mathcal{Y}}} & \cdots\\
b\arrow[ur]{}{s_{\mathcal{Y}}} & & & \\
  \end{tikzcd}
\caption*{The $\Sigma_{s}$-multialgebra $\mathcal{Y}$}
\end{figure}
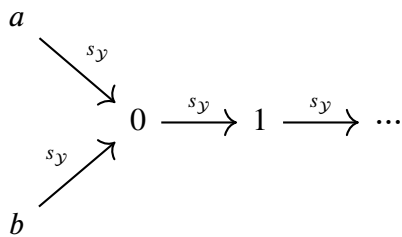
\end{example}

\begin{example}
Take the $\Sigma_{s}-$multialgebra $\mathcal{C}$ from Example \ref{C}.

We know that $\mathcal{C}$ is disconnected, but we state that it is not chainless: in fact,\\ $\{(-1)^{n}\}_{n\in\mathbb{N}}$ and $\{(-1)^{n+1}\}_{n\in\mathbb{N}}$ are chains in $\mathcal{C}$.
\end{example}

\begin{example}
Take the $\Sigma_{s}-$multialgebra $\mathcal{B}$ from Example \ref{B}.

We have already established that $\mathcal{B}$ has a basis and is generated by its ground, $\{0\}$, yet it is not chainless: $\{ 1\}_{n\in\mathbb{N}}$ is a chain in $\mathcal{B}$.
\end{example}

\section{Freely generated multialgebras}\label{Freely generated multialgebra}

So, as we mentioned time and time again, the algebras of formulas, in the study of universal algebras, are specially significant given they are the only absolutely free algebras, meaning they are the only ones with the unique extension property in the variety of all algebras over their signature.

So we turn now to a somewhat folkloric result: the unique extension property, relative to the class of all multialgebras on a given signature, is satisfied by no one; alternatively, this statement can be reworded as to say that the forgetful functor from the category of multialgebras to $\textbf{Set}$ does not have a left adjoint, what implies in particular that this category has no initial objects. 

Of course, such a result is stated in various ways, depending on your definition of homomorphism and, perhaps most importantly, even on your own definition of multialgebra. So, we offer what we believe to be a simplified proof of the result for the category $\textbf{MAlg}(\Sigma)$ as we have defined it. We proceed by contradiction, and start defining the multialgebras which do satisfy the unique extension property.

\begin{definition}
A $\Sigma$-multialgebra $\mathcal{A}=(A, \{\sigma_{\mathcal{A}}\}_{\sigma\in\Sigma})$ is freely generated by $X$, for a $X\subseteq A$, if for every $\Sigma$-multialgebra $\mathcal{B}=(B, \{\sigma_{\mathcal{B}}\}_{\sigma\in\Sigma})$ and map $f:X\rightarrow B$ there exists only one homomorphism $\overline{f}:\mathcal{A}\rightarrow\mathcal{B}$ extending $f$.
\end{definition}

In other words, if $j:X\rightarrow A$ is the inclusion, there exists only one homomorphism $\overline{f}:\mathcal{A} \to \mathcal{B}$ commuting the following diagram in {\bf Set}.
\[
\begin{tikzcd}
 & A \arrow{dr}{\overline{f}} \\
X \arrow{ur}{j} \arrow{rr}{f} && B
\end{tikzcd}
\]

\begin{lemma}\label{bijective homomorphism is full}
If $\varphi$ is a bijective $\Sigma$-homomorphism from $\mathcal{A}=(A, \{\sigma_{\mathcal{A}}\})$ to $\mathcal{B}=(B, \{\sigma_{\mathcal{B}}\})$ whose inverse is also a homomorphism, then $\varphi$ is a isomorphism.
\end{lemma}

\begin{proof}
We only have to prove that $\varphi$ is a full homomorphism, so take $\sigma\in\Sigma_{n}$ and $a_{1}, \dotsc  , a_{n}\in A$: from the fact $\varphi$ is a homomorphism we know that 
\[\{\varphi(a) \ : \  a\in\sigma_{\mathcal{A}}(a_{1}, \dotsc  , a_{n})\}\subseteq\sigma_{\mathcal{B}}(\varphi(a_{1}), \dotsc  , \varphi(a_{n}));\]
and since $\varphi^{-1}$ is also a homomorphism, we have that
\[\{\varphi^{-1}(b) \ : \  b\in\sigma_{\mathcal{B}}(\varphi(a_{1}), \dotsc  , \varphi(a_{n}))\}\subseteq \sigma_{\mathcal{A}}(a_{1}, \dotsc  , a_{n}).\]
Since $\varphi$ is a bijection, we may apply it to the last equation while preserving the inclusion, thus getting 
\[\sigma_{\mathcal{B}}(\varphi(a_{1}), \dotsc  , \varphi(a_{n}))\subseteq \{\varphi(a) \ : \  a\in \sigma_{\mathcal{A}}(a_{1}, \dotsc  , a_{n})\}.\]
\end{proof}

\begin{proposition}
If there exist $\Sigma$-multialgebras $\mathcal{A}$ and $\mathcal{B}$ which are freely generated by, respectively, $X$ and $Y$ such that $|X|=|Y|$, then $\mathcal{A}$ and $\mathcal{B}$ are isomorphic.
\end{proposition}

\begin{proof}
Since $X$ and $Y$ are of the same cardinality, there exists bijective functions $f:X\rightarrow Y$ and $g:Y\rightarrow X$ inverses of each other. Take the extensions $\overline{f}:\mathcal{A}\rightarrow\mathcal{B}$ and $\overline{g}:\mathcal{B}\rightarrow\mathcal{A}$ and we have that $\overline{g}\circ\overline{f}$ is a homomorphism extending $g\circ f=Id_{X}$, the identity on $X$.

Since the identical homomorphism $Id_{\mathcal{A}}:\mathcal{A}\rightarrow \mathcal{A}$ also extends $Id_{X}$, we have that $Id_{\mathcal{A}}=\overline{g}\circ\overline{f}$. In a similar way we have that $Id_{\mathcal{B}}=\overline{f}\circ\overline{g}$, so that $\mathcal{A}$ and $\mathcal{B}$ are isomorphic by use of Lemma \ref{bijective homomorphism is full}.
\end{proof}

This way we can refer ourselves to the unique $\Sigma$-multialgebra freely generated by $X$, up to isomorphisms.

We will denote by $\mathcal{U}:\textbf{MAlg}(\Sigma)\rightarrow \textbf{Set}$ the forgetful functor, which takes a multialgebra to its universe and a homomorphism to itself, seem now merely as a function between sets.

\begin{lemma}\label{adjoint of forget}
If there exist, for every set $X$, a $\Sigma$-multialgebra freely generated by $X$, then the functor $F:\textbf{Set}\rightarrow \textbf{MAlg}(\Sigma)$, associating a set $X$ to a $\Sigma$-multialgebra freely generated by $X$, which we will denote $FX$, and a function $f:X\rightarrow Y$ to the only homomorphism $\overline{f}:FX\rightarrow FY$ extending $f$, is a left adjoint of $\mathcal{U}$.
\end{lemma}

\begin{proof}
For $X$ a set and $\mathcal{A}$ a $\Sigma$-multialgebra with universe $A$ we consider
\[\Phi_{\mathcal{A}, X}:Hom_{\textbf{Set}}(X,\mathcal{U}\mathcal{A})\rightarrow Hom_{\textbf{MAlg}(\Sigma)}(FX, \mathcal{A})\]
associating a map $f:X\rightarrow A$ to the only homomorphism $\overline{f}:FX\rightarrow \mathcal{A}$ extending $f$. Each $\Phi_{\mathcal{A}, X}$ is clearly a bijection given that $FX$ is freely generated by $X$.

Now, given sets $X$ and $Y$, $\Sigma$-multialgebras $\mathcal{A}$ and $\mathcal{B}$, a function $f:Y\rightarrow X$ and a homomorphism $h:\mathcal{A}\rightarrow \mathcal{B}$, we have only left to prove that the following diagram commutes in {\bf Set}.
\[ \begin{tikzcd}[sep=huge]
Hom_{\textbf{Set}}(X, \mathcal{U}\mathcal{A}) \arrow{r}{\Phi_{\mathcal{A}, X}} \arrow{d}{Hom(f, \mathcal{U}\varphi)}& Hom_{\textbf{MAlg}(\Sigma)}(FX, \mathcal{A}) \arrow{d}{Hom(Ff, \varphi)} \\%
Hom_{\textbf{Set}}(Y, \mathcal{U}\mathcal{B}) \arrow{r}{\Phi_{\mathcal{B}, Y}}& Hom_{\textbf{MAlg}(\Sigma)}(FY, \mathcal{B})
\end{tikzcd}
\]

So we take a function $g:X\rightarrow \mathcal{U}\mathcal{A}$: taking the top-right edges of the diagram we have $\Phi_{\mathcal{A}, X}g=\overline{g}$ and $Hom(Ff, h)\overline{g}=h\circ \overline{g}\circ Ff$; on the bottom-left ones, $Hom(f, \mathcal{U}h)g=\mathcal{U}h\circ g\circ f$ and $\Phi_{\mathcal{B}, Y}(\mathcal{U}h\circ g\circ f)=\overline{\mathcal{U}h\circ g\circ f}$.

Now both $h\circ \overline{g}\circ Ff$ and $\overline{\mathcal{U}h\circ g\circ f}$ are homomorphisms from $FY$ to $\mathcal{B}$ extending $\mathcal{U}h\circ g\circ f:Y\rightarrow \mathcal{U}\mathcal{B}$: for the second one this is obvious, for the first we take an element $y\in Y$ and see that
\[h\circ \overline{g}\circ Ff(y)=h\circ\overline{g}\circ f(y)=h\circ g\circ f(y)=\mathcal{U}h\circ g\circ f(y)\]
since, respectively, $Ff=\overline{f}$, which extends $f$; $\overline{g}$ extends $g$, which is defined on $X$ that contains $f(y)$; and $\mathcal{U}h=h$, seen only as a function between sets.

Given that $FY$ is freely generated over $Y$, we have that $h\circ \overline{g}\circ Ff=\overline{\mathcal{U}\varphi\circ g\circ f}$ and the diagram in fact commutes.
\end{proof}

\begin{theorem}\label{no f-gen multi}
Given a non-empty signature $\Sigma$ and a set $X$, there does not exist a $\Sigma$-multialgebra freely generated by $X$.
\end{theorem}

\begin{proof}
Suppose that $\mathcal{A}=(A, \{\sigma_{\mathcal{A}}\}_{\sigma\in\Sigma})$ is freely generated by $X$ and let $\mathcal{V}$ be a set that properly contains $X$, meaning $\mathcal{V}\neq \emptyset$ and therefore that $\textbf{F}(\Sigma, \mathcal{V})$ is well defined; then take the identity function $j: X\rightarrow F(\Sigma, \mathcal{V})$, such that $j(x)=x$ for every $x\in X$, and the homomorphism $\overline{j}:\mathcal{A}\rightarrow\textbf{F}(\Sigma, \mathcal{V})$ extending $j$.

Now, take the identity function $id:\mathcal{V}\rightarrow F(\Sigma^{2}, \mathcal{V})$ and the collections of choices $C$ and $D$ from $\textbf{F}(\Sigma, \mathcal{V})$ to $\textbf{mF}(\Sigma, \mathcal{V}, 2)$ such that, for $\sigma\in\Sigma_{n}$,
\[C\sigma_{\alpha_{1}, \dotsc   , \alpha_{n}}^{\beta_{1}, \dotsc   , \beta_{n}}(\sigma(\alpha_{1}, \dotsc  ,\alpha_{n}))=\sigma^{0}(\beta_{1}, \dotsc  ,\beta_{n})\]
and
\[D\sigma_{\alpha_{1}, \dotsc   , \alpha_{n}}^{\beta_{1}, \dotsc   , \beta_{n}}(\sigma(\alpha_{1}, \dotsc  ,\alpha_{n}))=\sigma^{1}(\beta_{1}, \dotsc  ,\beta_{n}),\]
and consider the only homomorphisms $id_{C}, id_{D}:\textbf{F}(\Sigma, \mathcal{V})\rightarrow\textbf{mF}(\Sigma, \mathcal{V}, 2)$ extending, respectively, $id$ and $C$, and $id$ and $D$, which we know to exist given $\textbf{F}(\Sigma, \mathcal{V})$ is \textbf{cdf}-generated by $\mathcal{V}$. Since $id_{C}\circ \overline{j}, id_{D}\circ\overline{j}:\mathcal{A}\rightarrow \textbf{mF}(\Sigma, \mathcal{V}, 2)$ both extend the function $j':X\rightarrow {T}(\Sigma^{2}, \mathcal{V})$ such that $j'(x)=x$ for every $x\in X$ (remember $\mathcal{V}$ properly contains $X$), we have $id_{C}\circ \overline{j}=id_{D}\circ\overline{j}$.

Now, if $\alpha\in {T}(\Sigma, \mathcal{V})\setminus \mathcal{V}$, we have that there exist $\sigma\in \Sigma_{n}$, for $n\in\mathbb{N}$, and elements $\alpha_{1}, \dotsc   , \alpha_{n}\in F(\Sigma, \mathcal{V})$ such that $\alpha=\sigma(\alpha_{1}, \dotsc  ,\alpha_{n})$; in this case, 
\[id_{C}(\alpha)=\sigma^{0}(id_{C}(\alpha_{1}), \dotsc  ,id_{C}(\alpha_{n}))\neq \sigma^{1}(id_{D}(\alpha_{1}), \dotsc  ,id_{D}(\alpha_{n}))=id_{D}(\alpha),\]
given that the main connectives are distinct, and therefore implying that $id_{C}$ and $id_{D}$ are always different outside of $\mathcal{V}$. 

Since $id_{C}\circ \overline{j}=id_{D}\circ\overline{j}$, we must have that $\overline{j}(A)\subseteq \mathcal{V}$, and this is absurd since we are assuming $\Sigma$ non-empty: if $\Sigma_{0}\neq\emptyset$, for a $\sigma\in\Sigma_{0}$ and $a\in\sigma_{\mathcal{A}}$ we have that $\overline{j}(a)=\sigma$ in $F(\Sigma, \mathcal{V})$, which is not in $X$; if it is another $\Sigma_{n}$ which is not empty, given $a\in A$ (which exists given the universe of multialgebras are assumed to be non-empty) we have that, for $b\in \sigma_{\mathcal{A}}(a, \dotsc   , a)$, is valid that $\overline{j}(b)=\sigma(\overline{j}(a), \dotsc   , \overline{j}(a))$, which is again not in $\mathcal{V}$.

We must conclude that there are no freely generated multialgebras.
\end{proof}

\begin{corollary}\label{no free obj}
The category $\textbf{MAlg}(\Sigma)$ does not have an initial object.
\end{corollary}

\begin{proof}
We state that if $\mathcal{A}$ is an initial object, $\mathcal{A}$ is freely generated by $\emptyset$: in fact, for every $\Sigma$-multialgebra $\mathcal{B}$ and map $f:\emptyset\rightarrow B$, there exists a single homomorphism $!_{\mathcal{B}}:\mathcal{A}\rightarrow \mathcal{B}$ extending $f=\emptyset$, that is, the only homomorphism between $\mathcal{A}$ and $\mathcal{B}$.

By Theorem~\ref{no f-gen multi}, freely generated multialgebras do not exist, what ends the proof.
\end{proof}

\begin{theorem}\label{no left adjoints}
The forgetful functor $\mathcal{U}:\textbf{MAlg}(\Sigma)\rightarrow\textbf{Set}$ does not have a left adjoint.
\end{theorem}

\begin{proof}
For suppose we have a left adjoint $F:\textbf{Set}\rightarrow\textbf{MAlg}(\Sigma)$ of $\mathcal{U}$, so that $F$ has a right adjoint and is therefore cocontinuous. Since $\emptyset$ is the initial object in $\textbf{Set}$, we have that $F\emptyset$ must be an initial object in $\textbf{MAlg}(\Sigma)$, which by Corollary~\ref{no free obj} does not exist.
\end{proof}

Answering a question posed at the end of Section \ref{Formulas and how to interpret them}, we can prove with arguments similar to those found in Theorem \ref{no f-gen multi} and Corollary \ref{no free obj} that the category of $\Sigma$-multialgebras equipped with full homomorphisms, $\textbf{MAlg}_{=}(\Sigma)$, does not have intial objects, and so its forgetful functor into $\textbf{Set}$ does not have a left adjoint.

In fact, suppose $\mathcal{A}$ is an initial object of $\textbf{MAlg}_{=}(\Sigma)$, and let $!$ be the only full homomorphism from $\mathcal{A}$ to $\textbf{mF}(\Sigma, \mathcal{V}, 2)$; we take the identity function $id:\mathcal{V}\rightarrow F(\Sigma^{2}, \mathcal{V})$ and the collections of choices $C$ and $D$ from $\textbf{mF}(\Sigma, \mathcal{V}, 2)$ into itself such that
\[C\sigma_{\alpha_{1}, \dotsc   , \alpha_{n}}^{\beta_{1}, \dotsc   , \beta_{n}}(\sigma^{0}(\alpha_{1}, \dotsc  ,\alpha_{n}))=\sigma^{0}(\beta_{1}, \dotsc  ,\beta_{n}),\quad C\sigma_{\alpha_{1}, \dotsc   , \alpha_{n}}^{\beta_{1}, \dotsc   , \beta_{n}}(\sigma^{1}(\alpha_{1}, \dotsc  ,\alpha_{n}))=\sigma^{1}(\beta_{1}, \dotsc  ,\beta_{n}),\]
\[D\sigma_{\alpha_{1}, \dotsc   , \alpha_{n}}^{\beta_{1}, \dotsc   , \beta_{n}}(\sigma^{0}(\alpha_{1}, \dotsc  ,\alpha_{n}))=\sigma^{1}(\beta_{1}, \dotsc  ,\beta_{n})\quad\text{and}\quad D\sigma_{\alpha_{1}, \dotsc   , \alpha_{n}}^{\beta_{1}, \dotsc   , \beta_{n}}(\sigma^{1}(\alpha_{1}, \dotsc  ,\alpha_{n}))=\sigma^{0}(\beta_{1}, \dotsc  ,\beta_{n}),\]
and it is clear that $id_{C}$ and $id_{D}$ are both full homomorphisms satisfying $id_{C}\circ !=id_{D}\circ !$. Since $id_{C}$ and $id_{D}$ are always distinct over the build of $\textbf{mF}(\Sigma, \mathcal{V}, 2)$, we get that the image of $!$ is contained in $\mathcal{V}$, a contradiction as we saw in the proof of Theorem \ref{no f-gen multi}.


\section{A categorical characterization}

Given the category $\textbf{MAlg}(\Sigma)$ whose objects are all $\Sigma-$multialgebras and morphisms from $\mathcal{A}$ to $\mathcal{B}$ are all homomorphisms from $\mathcal{A}$ to $\mathcal{B}$, we would like to construct a related category $\textbf{MG}(\Sigma)$ whose objects are also exactly all $\Sigma-$multialgebras, but where the $\textbf{cdf}$-generated multialgebras can be construed as the image of an adjoint-having functor. The motivation behind the search for such a $\textbf{MG}(\Sigma)$ is that having such a category, and the associated pair of adjoint functors, would make weakly free multialgebras far more similar to absolutely free algebras than they are now: after all, the latter are images under the so called free functor, right adjoint of the forgetful functor from the category of $\Sigma$-algebras $\textbf{Alg}(\Sigma)$ to the category of sets $\textbf{Set}$. Unfortunately, since, as we showed in Theorem \ref{no left adjoints}, the forgetful functor of $\textbf{MAlg}(\Sigma)$ does not have a left-adjoint, the task is not as simple as in the classical case and will require far more craftsmanship.

\begin{definition}
Given a homomorphism $\varphi:\mathcal{A}\rightarrow\mathcal{B}$ between $\Sigma-$multialgebras $\mathcal{A}$ and $\mathcal{B}$, we say that $\varphi$ is \textit{ground-preserving}, or that it preserves grounds, if $\varphi(G(\mathcal{A}))\subseteq G(\mathcal{B})$.
\end{definition}

\begin{proposition}
Taking as objects all $\Sigma-$multialgebras and as morphisms between the\\ $\Sigma-$multialgebras $\mathcal{A}$ and $\mathcal{B}$ all ground-preserving homomorphisms between $\mathcal{A}$ and $\mathcal{B}$, the resulting structure $\textbf{MAlg}_{G}(\Sigma)$ is a subcategory of $\textbf{MAlg}(\Sigma)$.
\end{proposition}

\begin{proof}
Clearly the identical homomorphism $Id_{\mathcal{A}}$ is ground-preserving, so $Id_{\mathcal{A}}$ is in\\ $Hom_{\textbf{MAlg}_{G}(\Sigma)}(\mathcal{A}, \mathcal{A})$.

And if $\varphi:\mathcal{A}\rightarrow \mathcal{B}$ and $\psi:\mathcal{B}\rightarrow\mathcal{C}$ are ground-preserving, then 
\[\psi\circ\varphi(G(\mathcal{A}))=\psi(\varphi(G(\mathcal{A}))\subseteq\psi(G(\mathcal{B}))\subseteq G(\mathcal{C}),\]
so that $\psi\circ\varphi$ is also ground-preserving and $\textbf{MAlg}_{G}(\Sigma)$ is in fact a subcategory of $\textbf{MAlg}(\Sigma)$.
\end{proof}

Now, for ground-preserving homomorphisms $\varphi, \psi:\mathcal{A}\rightarrow\mathcal{B}$, we define an equivalence relation "${\sim}$" in $Hom_{\textbf{MAlg}(\Sigma)}(\mathcal{A}, \mathcal{B})$ by means of
\[\varphi{\sim}\psi \Leftrightarrow \varphi|_{G(\mathcal{A})}=\psi|_{G(\mathcal{A})}.\]
That is in fact an equivalence relation since, given homomorphisms $\varphi, \psi, \nu:\mathcal{A}\rightarrow\mathcal{B}$ we have that:
\begin{enumerate}
\item $\varphi{\sim}\varphi$, since $\varphi|_{G(\mathcal{A})}=\varphi|_{G(\mathcal{A})}$;
\item if $\varphi{\sim}\psi$, then $\varphi|_{G(\mathcal{A})}=\psi|_{G(\mathcal{A})}$ and therefore $\psi|_{G(\mathcal{A})}=\varphi|_{G(\mathcal{A})}$, and so $\psi{\sim}\varphi$;
\item if $\varphi{\sim}\psi$ and $\psi{\sim}\nu$, we have that $\varphi|_{G(\mathcal{A})}=\psi|_{G(\mathcal{A})}$ and $\psi|_{G(\mathcal{A})}=\nu|_{G(\mathcal{A})}$, and therefore $\varphi|_{G(\mathcal{A})}=\nu|_{G(\mathcal{A})}$ and $\varphi{\sim}\nu$.
\end{enumerate}
We will denote the equivalence class by ${\sim}$ with representative a homomorphism $\varphi$ by $[\varphi]$. Finally, we can define the morphisms from $\mathcal{A}$ to $\mathcal{B}$ in $\textbf{MG}(\Sigma)$\label{MGSigma} to be
\[Hom_{\textbf{MG}(\Sigma)}(\mathcal{A}, \mathcal{B})=Hom_{\textbf{MAlg}_{G}(\Sigma)}(\mathcal{A}, \mathcal{B})/{\sim},\]
and their composition by, given $[\varphi]\in Hom_{\textbf{MG}(\Sigma)}(\mathcal{A}, \mathcal{B})$ and $[\psi]\in Hom_{\textbf{MG}(\Sigma)}(\mathcal{B}, \mathcal{C})$, 
\[[\psi]\circ[\varphi]=[\psi\circ\varphi].\]
We state that this definition of composition is well behaved: if $[\varphi_{1}]=[\varphi_{2}]$ and $[\psi_{1}]=[\psi_{2}]$ for $\varphi_{1}, \varphi_{2}:\mathcal{A}\rightarrow \mathcal{B}$ and $\psi_{1}, \psi_{2}:\mathcal{B}\rightarrow\mathcal{C}$ ground-preserving homomorphisms, we have that for an $a\in G(\mathcal{A})$ it is valid that
\[\psi_{1}\circ\varphi_{1}(a)=\psi_{1}(\varphi_{1}(a))=\psi_{1}(\varphi_{2}(a));\]
since $\varphi_{2}$ is ground-preserving, we have that $\varphi_{2}(a)\in G(\mathcal{B})$, implying
\[\psi_{1}(\varphi_{2}(a))=\psi_{2}(\varphi_{2}(a))=\psi_{2}\circ\varphi_{2}(a),\]
and therefore $\psi_{1}\circ\varphi_{1}|_{G(\mathcal{A})}=\psi_{2}\circ\varphi_{2}|_{G(\mathcal{A})}$, meaning $[\psi_{1}\circ\varphi_{1}]=[\psi_{2}\circ\varphi_{2}]$.

\begin{theorem}
$\textbf{MG}(\Sigma)$ is a category.
\end{theorem}

\begin{proof}
By the above remark $Hom_{\textbf{MG}(\Sigma)}$ has a well-defined composition: to see that this composition is also associative, take ground-preserving homomorphisms $\varphi:\mathcal{A}\rightarrow \mathcal{B}$, $\psi:\mathcal{B}\rightarrow\mathcal{C}$ and $\theta:\mathcal{C}\rightarrow \mathcal{D}$: then
\[[\theta]\circ([\psi]\circ[\varphi])=[\theta]\circ[\psi\circ\varphi]=[\theta\circ(\psi\circ\varphi)]=[(\theta\circ\psi)\circ\varphi]=[\theta\circ\psi]\circ[\varphi]=([\theta]\circ[\psi])\circ[\varphi].\]

Furthermore, $[Id_{\mathcal{A}}]$, for $Id_{\mathcal{A}}\in Hom_{\textbf{MAlg}(\Sigma)}(\mathcal{A}, \mathcal{A})$ the usual identical homomorphism, is the identity on $Hom_{\textbf{MG}(\Sigma)}(\mathcal{A}, \mathcal{A})$, since for any other $\Sigma-$multialgebra $\mathcal{B}$ and ground-preserving homomorphisms $\varphi:\mathcal{A}\rightarrow \mathcal{B}$ and $\psi:\mathcal{B}\rightarrow\mathcal{A}$ we have
\[[\varphi]\circ[Id_{\mathcal{A}}]=[\varphi\circ Id_{\mathcal{A}}]=[\varphi]\]
and
\[[Id_{\mathcal{A}}]\circ[\psi]=[Id_{\mathcal{A}}\circ\psi]=[\psi].\]
\end{proof}

Our goal is to prove that the functor associating a set $X$ to a $\textbf{cdf}$-generated multialgebra over $X$ is the adjoint of a "forgetful functor" on $\textbf{MG}(\Sigma)$: the clear problem is that the morphisms on this category are not functions, but rather equivalence classes of functions.  

So we define the functor $U:\textbf{MG}(\Sigma)\rightarrow \textbf{Set}$, which we shall call forgetful, by $U\mathcal{A}=G(\mathcal{A})$ and, for a morphism $[\varphi]:\mathcal{A}\rightarrow\mathcal{B}$, by $U[\varphi]=\varphi|_{G(\mathcal{A})}$, which is well-defined since $\varphi$ is ground-preserving.

$U$ is indeed a functor since: if $[Id_{\mathcal{A}}]$ is the identity on $\mathcal{A}$, $U[Id_{\mathcal{A}}]=Id_{\mathcal{A}}|_{G(\mathcal{A})}$, which is the identical function on $G(\mathcal{A})$; and given homomorphisms $\varphi:\mathcal{A}\rightarrow\mathcal{B}$ and $\psi:\mathcal{B}\rightarrow\mathcal{C}$, 
\[U[\psi\circ\varphi]=\psi\circ\varphi|_{G(\mathcal{A})}=\psi|_{\varphi(G(\mathcal{A}))}\circ\varphi|_{G(\mathcal{A})}=\psi|_{G(\mathcal{B})}\circ\varphi|_{G(\mathcal{A})}=U[\psi]\circ U[\varphi].\]

\begin{theorem}
The functor 
\[F_{\kappa}:\textbf{Set}\rightarrow\textbf{MG}(\Sigma)\]
associating to a set $X$ the $\Sigma-$multialgebra $\textbf{mF}(\Sigma, X, \kappa)$ and to a function $f:X\rightarrow Y$ the morphism $[f_{C}]$, for any collection of choices $C$ from $\textbf{mF}(\Sigma, X, \kappa)$ to $\textbf{mF}(\Sigma, Y, \kappa)$, is a left-adjoint of the forgetful functor $U:\textbf{MG}(\Sigma)\rightarrow\textbf{Set}$.
\end{theorem}

\begin{proof}
First of all $F$ is well defined over morphisms given that for any two collection of choices $C$ and $D$ from $\textbf{mF}(\Sigma, X, \kappa)$ to $\textbf{mF}(\Sigma, Y, \kappa)$ we have that $f_{C}|_{X}=f=f_{D}|_{X}$, so $[f_{C}]=[f_{D}]$.

It is a functor since:
\begin{enumerate}
\item given the identity between sets $Id_{X}:X\rightarrow X$, the morphism $F_{\kappa}Id_{X}$, when restricted to $X$, equals $Id_{\textbf{mF}(\Sigma, X, \kappa)}$, and therefore $F_{\kappa}Id_{X}=[Id_{\textbf{mF}(\Sigma, X, \kappa)}]$;
\item given functions $f:X\rightarrow Y$ and $g:Y\rightarrow Z$, 
\[F_{\kappa}g\circ F_{\kappa}f=[g_{D}]\circ[f_{C}]=[g_{D}\circ f_{C}],\]
for collections of choices $C$ from $\textbf{mF}(\Sigma, X, \kappa)$ to $\textbf{mF}(\Sigma, Y, \kappa)$ and $D$ from $\textbf{mF}(\Sigma, Y, \kappa)$ to $\textbf{mF}(\Sigma, Z, \kappa)$; for any collection of choices $E$ from $\textbf{mF}(\Sigma, X, \kappa)$ to $\textbf{mF}(\Sigma, Z, \kappa)$ we have that $(g\circ f)_{E}$, when restricted to the ground, is the same as $g_{D}\circ f_{C}$, that is, $g\circ f$, and therefore 
\[F_{\kappa}g\circ F_{\kappa}f=[(g\circ f)_{E}]=F_{\kappa}(g\circ f).\]
\end{enumerate}

We define the bijection $\Phi_{X,\mathcal{A}}:Hom_{\textbf{Set}}(X, U\mathcal{A})\rightarrow Hom_{\textbf{MG}(\Sigma)}(F_{\kappa}X, \mathcal{A})$ by, given a function $f:X\rightarrow U\mathcal{A}$ (where one should remember that $U\mathcal{A}=G(\mathcal{A})$), $\Phi_{X, \mathcal{A}}f=[f_{C}]$ for any collection of choices $C$ from $\textbf{mF}(\Sigma, X, \kappa)$ to $\mathcal{A}$: again, this is clearly well-defined since if $D$ is another collection of choices we have $f_{C}|_{X}=f=f_{D}|_{X}$.

To see that such $\phi_{X \mathcal{A}}$ is really a bijection is easy:
\begin{enumerate}
\item given two functions $f, g:X\rightarrow G(\mathcal{A})$, if $\Phi_{X, \mathcal{A}}f=\Phi_{X, \mathcal{A}}g$, then $[f_{C}]=[g_{D}]$, for any collections of choices $C$ and $D$ from $\textbf{mF}(\Sigma, X, \kappa)$ to $\mathcal{A}$, and therefore $f=f_{C}|_{X}=g_{D}|_{X}=g$;
\item for any homomorphism $\varphi:\textbf{mF}(\Sigma, X, \kappa)\rightarrow\mathcal{A}$ and any collection of choices $C$ from $\textbf{mF}(\Sigma, X,\kappa)$ to $\mathcal{A}$, once we define $f=\varphi|_{X}$ we have that $[\varphi]=[f_{C}]=F_{\kappa}f$.
\end{enumerate}

Now we must prove that for any sets $X$ and $Y$, any $\Sigma-$multialgebras $\mathcal{A}$ and $\mathcal{B}$, any function $f:Y\rightarrow X$ and morphism $[\varphi]:\mathcal{A}\rightarrow\mathcal{B}$ we have that the following diagram commutes.
\[ \begin{tikzcd}[row sep=5em, column sep=3em]
Hom_{\textbf{Set}}(X, U\mathcal{A}) \arrow{r}{\Phi_{X, \mathcal{A}}} \arrow{d}{Hom(f, U[\varphi])}& Hom_{\textbf{MG}(\Sigma)}(F_{\kappa}X, \mathcal{A}) \arrow{d}{Hom(F_{\kappa}f, [\varphi])} \\%
Hom_{\textbf{Set}}(Y, U\mathcal{B}) \arrow{r}{\Phi_{Y, \mathcal{B}}}& Hom_{\textbf{MG}(\Sigma)}(F_{\kappa}Y, \mathcal{B})
\end{tikzcd}
\]
Given a map $g:X\rightarrow U\mathcal{A}$, $\Phi_{X, \mathcal{A}}g=[g_{C}]$ for some collection of choices $C$ from $\textbf{mF}(\Sigma, X, \kappa)$ to $\mathcal{A}$, and 
\[Hom(F_{\kappa}f, [\varphi])\Phi_{X,\mathcal{A}}g=Hom(F_{\kappa}f, [\varphi])[g_{C}]=[\varphi]\circ[g_{C}]\circ [f_{D}]=[\varphi\circ g_{C}\circ f_{D}],\]
for some collection of choices $D$ from $\textbf{mF}(\Sigma, Y, \kappa)$ to $\textbf{mF}(\Sigma, X, \kappa)$; on the other branch of the diagram, we have that $Hom(f, U[\varphi])g=U[\varphi]\circ g\circ f$ and therefore
\[\Phi_{Y, \mathcal{B}} (Hom(f, U[\varphi])g)=\Phi_{Y, \mathcal{B}}(U[\varphi]\circ g\circ f)=[(U[\varphi]\circ g\circ f)_{E}]\]
for some collection of choices $E$ from $\textbf{mF}(\Sigma, Y, \kappa)$ to $\mathcal{B}$. Now, is enough to prove that $\varphi\circ g_{C}\circ f_{D}$ and $(U[\varphi]\circ g\circ f)_{E}$ are equal over the ground of $\textbf{mF}(\Sigma, Y, \kappa)$, that is, $Y$. This is easy, since for an element $y\in Y$ we have that
\[\varphi\circ g_{C}\circ f_{D}(y)=\varphi\circ g_{C}(f(y))=\varphi(g(f(y))\]
and that
\[(U[\varphi]\circ g\circ f)_{E}(y)=U[\varphi]\circ g\circ f(y)=U[\varphi]\circ g(f(y))=U[\varphi](g(f(y)))=\varphi(g(f(y))).\]
\end{proof}

Notice that two multialgebras $\mathcal{A}$ and $\mathcal{B}$ $\textbf{cdf}$-generated respectively by sets $X$ and $Y$ of the same cardinality are isomorphic in $\textbf{MG}(\Sigma)$: for, given a bijective function $f:X\rightarrow Y$ and collections of choices $C$ from $\mathcal{A}$ to $\mathcal{B}$ and $D$ from $\mathcal{B}$ to $\mathcal{A}$, we have that $[f^{-1}_{D}]\circ[f_{C}]=[Id_{\mathcal{A}}]$, since $f^{-1}\circ f$ is the identity on $X$, and $[f_{C}]\circ[f^{-1}_{D}]=[Id_{\mathcal{B}}]$. 

In general, this makes the classes of isomorphisms of multialgebras in $\textbf{MG}(\Sigma)$ too coarse, meaning that this category can not tell apart multialgebras that, intuitively, we consider very much distinct, such as $\textbf{mF}(\Sigma, X, \kappa)$ and $\textbf{mF}(\Sigma, X, \beta)$ for distinct cardinals $\kappa$ and $\beta$; of course, this is not to say that $\textbf{MG}(\Sigma)$ doesn't play a role in attempting to categorify weakly free multialgebras, just that this category is not completely successful at that task. Of course, before anything definitive is said about $\textbf{MG}(\Sigma)$'s success, or lack of it, as an ambient category for further studies on non-deterministic algebraization, a better understanding on the relationship between legal valuations, and grounds and collections of choices is necessary: after all, while the former are already quite well established as the central object of study of non-deterministic semantics, the latter have only recently been formally defined and, despite the results outlined here, are arguably not very well understood.

Now, while still in the topic of categories, we may represent those found in this chapter through the following diagram.

\[ \begin{tikzcd}
\textbf{MG}(\Sigma)\arrow[transform canvas={yshift=.9ex}]{ddrr}{U} & & \textbf{MAlg}_{G}(\Sigma)\arrow[swap]{ll}{Q}\arrow[hook]{rr}\arrow[dashed]{dd} && \textbf{MAlg}(\Sigma) \arrow{ddll}{\mathcal{U}}\\
&&&&\\
&& \textbf{Set}\arrow[transform canvas={yshift=-.9ex}]{uull}{F_{\kappa}}&&
\end{tikzcd}
\]

Clearly $\textbf{MAlg}_{G}(\Sigma)$ is a subcategory of $\textbf{MAlg}(\Sigma)$, sharing all of its objects, while the dashed arrow represents the forgetful functor from $\textbf{MAlg}
_{G}(\Sigma)$ into $\textbf{Set}$, of which we have nothing interesting to say. We can also define the functor $Q:\textbf{MAlg}_{G}(\Sigma)\rightarrow\textbf{MG}(\Sigma)$, that is the identity on objects and takes ground-preserving homomorphism to their equivalence classes, to better understand the relation between $\textbf{MAlg}_{G}(\Sigma)$ and $\textbf{MG}(\Sigma)$: $Q$ is surjective on both objects and morphisms, but is only one-to-one on objects. Finally, from Theorem \ref{no left adjoints}, $\mathcal{U}$ does not have a left adjoint.

\newpage
\printbibliography[segment=\therefsegment,heading=subbibliography]
\end{refsegment}

\begin{refsegment}
\defbibfilter{notother}{not segment=\therefsegment}
\setcounter{chapter}{2}
\chapter{Multialgebras as partially ordered algebras}\label{Chapter3}\label{Chapter 3}

It is a fundamental result (see \cite{Oosten} for a proof) that there exists a (bijective) correspondence between sets and complete, atomic Boolean algebras (\textit{CABA}'s\index{CABA}): on one hand, a set is taken to its powerset (which is certainly a \textit{CABA}), while on the other direction a \textit{CABA} is taken to its set of atomic elements; accordingly, a map $f:X\rightarrow Y$ is taken to the map $\mathcal{P}(f):\mathcal{P}(Y)\rightarrow\mathcal{P}(X)$ such that $\mathcal{P}(f)(B)=\{a\in A\ :\  f(a)\in B\}$, for any $B\subseteq Y$, and a continuous homomorphism of Boolean algebras $\varphi:\mathcal{A}\rightarrow\mathcal{B}$ is taken to the function $A\varphi$, from atoms of $\mathcal{B}$ to atoms of $\mathcal{A}$, such that $A\varphi(b)=a$ is the unique atom of $\mathcal{A}$ satisfying that $b\leq \varphi(a)$, for any atom $b$ of $\mathcal{B}$.

These two assignments can be made into functors, giving rise to an equivalence of $\textbf{Set}^{op}$ and $\textbf{CABA}$\label{CABA}, the category with \textit{CABA}'s as objects and continuous homomorphisms of Boolean algebras as morphisms. This is part of a broader area of study, known by Stone dualities, which studies relationships between partially ordered sets and topological spaces (or rather, their respective categories), and was established by Stone (\cite{Stone}) and his representation theorem, which states that every Boolean algebra is isomorphic to a field of sets, specifically the algebra of open and close subsets of its Stone space, a topological space where points are ultrafilters of the original Boolean algebra: an ultrafilter of a Boolean algebra $\mathcal{A}$ is any non-empty subset $U$ of its universe such that 
\begin{enumerate}
\item $a, b\in U$ imply the existence of $c\in U$ satisfying $c\leq a\wedge b$; 
\item $a\in U$ and $a\leq b$ imply $b\in U$; 
\item and, for every element $a$ of $\mathcal{A}$, either $a\in U$ or $\neg a\in U$. 
\end{enumerate}
The topology of the Stone space $S(\mathcal{A})$ of $\mathcal{A}$ is generated by the sets of ultrafilters sharing an element $a$ of $\mathcal{A}$, that is, $\{U\in S(\mathcal{A}): a\in U\}$. Of course, Stone' theorem corresponds to an equivalence between the category $\textbf{BA}$ of Boolean algebras and that of Stone spaces, with continuous functions between them as morphisms.

In this search of dualities with categories of partially ordered sets, we focus on a more palpable equivalence, closely related to the one between $\textbf{Set}^{op}$ and $\textbf{CABA}$: we shift from $\textbf{Set}$ to the category $\textbf{MAlg}(\Sigma)$ of multialgebras over a given signature $\Sigma$ by adding further structure to $\textbf{Set}$, that is, multioperations; correspondingly, we replace $\textbf{CABA}$ by a category of \textit{CABA}'s equipped with $\Sigma$-operations compatible with their orders.

We make only small adjustments to these categories, in order to better suit our needs (although results for the aforementioned categories do exist): since we are most interested in non-partial multialgebras, we consequently exchange \textit{CABA}'s by posets corresponding to powersets with the empty-set removed (that is, \textit{CABA}'s without minimum elements, an apparently paradoxical concept we make precise further ahead). This way, a multialgebra with universe $A$ is taken to an algebra over the set of non-empty subsets of $A$, with order given by inclusion and operations given by ``accumulating'' the operations of the multialgebra: that is, an $n$-ary operator $\sigma$, when evaluated on non-empty subsets $A_{i}$ of $A$, returns the non-empty set
\[\bigcup_{(a_{1}, \dotsc  , a_{n})\in A_{1}\times \cdots\times A_{n}}\sigma_{\mathcal{A}}(a_{1}, \dotsc  , a_{n});\]
reciprocally, an ``almost'' Boolean' algebra, as we have loosely described, is taken to its set of atomic elements, transformed into a multialgebra by addition of multioperations returning, for $\sigma$ of arity $n$ and atoms $a_{i}$, 
\[\{\text{$a$ is an atom} : a\leq \sigma_{\mathcal{A}}(a_{1}, \dotsc  , a_{n})\}.\]

In non-deterministic semantics (\cite{AvronLev}), all these considerations involving a category equivalent to that of multialgebras offer an alternative: many logicians, mainly due to philosophical objections, are reluctant to use multialgebras in order to characterize a given logic; the equivalence we here present shows one can, if one chooses to, replace those semantics based on multialgebras by, arguably more complex, semantics involving both an algebra and an underlying order, which are however very classically behaved; of course, adding order-theoretic elements to algebraic semantics is nothing new, see \cite{Raftery} for an example, and the already existing theory of order-algebraizable logics may offer new insights into non-deterministic semantics through the paradigm here presented.

\section{Complete, atomic and bottomless Boolean algebras}\label{bottomless Boolean algebras}

We have discussed Boolean algebras in some depth in Section \ref{Lattices, and Boolean... }, but we return to the topic now with a more order-theoretic mindset. A partially ordered set\index{Partially ordered set}, or \textit{poset}\index{Poset}, is a pair $(A, \leq)$, with $A$ a set and $\leq$ a relation on $A$ (meaning it is a subset of $A\times A$) such that:
\begin{enumerate}
\item $\leq$ is reflexive, meaning that for any $a\in A$, $(a, a)$ is in the relation, what we write as $a\leq a$;
\item $\leq$ is anti-symmetric, meaning that if $a\leq b$ and $b\leq a$, then $a=b$, for any $a, b\in A$;
\item $\leq$ is transitive, what means that if $a\leq b$ and $b\leq c$, then $a\leq c$, for any $a, b, c\in A$.
\end{enumerate}
Whenever $a\leq b$, we may also write $b\geq a$; if $a\leq b$ and we know $a\neq b$, one can also write $a<b$, what may also be denoted by $b>a$.

\begin{definition}
Given a poset $(A, \leq)$, an element $a\in A$ is:
\begin{enumerate}
\item a maximum if, for all $b\in A$, $b\leq a$;
\item a minimum if, for all $b\in A$, $a\leq b$;
\item maximal if, for all $b\in A$, $a\leq b$ implies $a=b$;
\item minimal if, for all $b\in A$, $b\leq a$ implies $a=b$.
\end{enumerate}

Furthermore, given a subset $S\subseteq A$, we say $a\in A$ is:
\begin{enumerate}
\item an upper bound for $S$ if, for any $s\in S$, $s\leq a$;
\item a lower bound for $S$ if, for any $s\in S$, $a\leq s$;
\item the supremum of $S$, when we write $a=\sup S$, if it is the minimum of all upper bounds for $S$, meaning that:
\begin{enumerate}
\item for any $s\in S$, $s\leq a$;
\item if $b\in A$ is such that, for any $s\in S$, $s\leq b$, then $a\leq b$;\footnote{Notice that there is indeed only one supremum of a set, as well as only one infimum: if $a$ and $b$ are both minimum upper bounds for $S$, the fact that $a$ is a minimum upper bound and that $b$ is an upper bound gives us $a\leq b$; reciprocally, the fact $b$ is a minimum upper bound and $a$ is an upper bound implies $b\leq a$, and so $a=b$. The same reasoning applies to infima.}
\end{enumerate}
\item the infimum of $S$, when we write $a=\inf S$, if it is the maximum of all lower bounds for $S$, meaning that:
\begin{enumerate}
\item for any $s\in S$, $a\leq s$;
\item if $b\in A$ is such that, for any $s\in S$, $a\leq s$, then $b\leq a$.
\end{enumerate}
\end{enumerate}
\end{definition}

Of course, we may define maxima, minima, maximal and minimal elements for a subset $S\subseteq A$ of $A$, by merely restricting the order of $(A, \leq)$ to $S$.

A Boolean algebra, to which we have already given on Section \ref{Lattices, and Boolean... } a purely algebraic formulation, is a partially-ordered sets $(A, \leq)$ such that: there are a maximum (denoted by $1$) and a minimum ($0$) elements, which we shall assume distinct; for every pair of elements $a, b\in A$, the set $\{a, b\}$ has a supremum, denoted by $a\vee b$, and an infimum, denoted by $a\wedge b$; and every element $a$ has a complement $b$, which satisfies
\[b=\inf\{c\in A:  \sup\{a, c\}=1\}\]
and
\[b=\sup\{c\in A:  \inf\{a, c\}=0\}.\]

A poset $(A, \leq)$ is said to be complete\index{Poset, complete} if every $S\subseteq A$ has a supremum and an infimum. 

\begin{lemma}
\begin{enumerate} 
\item Every Boolean algebra $(A, \leq)$ is distributive, meaning 
\[a\vee (b\wedge c)=(a\vee b)\wedge(a\vee c)\quad\text{and}\quad a\wedge(b\vee c)=(a\wedge b)\vee(a\wedge c),\]
for any $a, b, c\in A$;

\item every complete Boolean algebra $(A, \leq)$ is infinite distributive, meaning that for any $S\cup\{a\}\subseteq A$, 
\[\sup\{\inf\{a, s\}: s\in S\}=\inf\{a, \sup S\}\quad\text{and}\quad\inf\{\sup\{a, s\}: s\in S\}=\sup\{a, \inf S\}.\]
\end{enumerate}
\end{lemma}

\begin{proof}
We only prove the first equation related to distributivity, given that the proof for the other is very similar. Denote $p=a\vee (b\wedge c)$, meaning 
\[p=\sup\{a, \inf\{b, c\}\},\]
and $q=(a\vee b)\wedge (a\vee c)$, that is 
\[q=\inf\{\sup\{a, b\}, \sup\{a, c\}\};\]
since $\inf\{b, c\}\leq b$ and $\inf\{b, c\}\leq c$, $p=\sup\{a, \inf\{b, c\}\}\leq \sup\{a, b\}$ and, analogously, $p\leq \sup\{a, c\}$, implying $p$ is a lower bound for $\{\sup\{a, b\}, \sup\{a, c\}\}$. Since $q$ is the largest of these lower bounds, we get $p\leq q$.

Reciprocally, $p=\sup\{a, \inf\{b, c\}\}$ implies $a\leq p$ and $\inf\{b, c\}\leq p$: if $\inf\{b, c\}=p$, this means $a\leq \inf\{b, c\}$ and therefore $a\leq b$ and $a\leq c$, meaning that $\sup\{a, b\}=b$ and $\sup\{a, c\}=c$ and therefore $q=\inf\{b, c\}=p$; if, otherwise, $\inf\{b, c\}<p$, $p$ is not a lower bound for $\{b,c\}$ and so either $b\leq p$ or $c\leq p$, implying either $\sup\{a, b\}\leq p$ or $\sup\{a, c\}\leq p$ and consequently $q=\inf\{\sup\{a, b\}, \sup\{a, c\}\}\leq p$, from what we get $p=q$.

More generally, we now prove one of the equalities of infinite-distributivity, the other being analogous. We know both 
\[p=\sup\{\inf\{a, s\}: s\in S\}\quad\text{and}\quad q=\inf\{a, \sup S\}\]
exist, given by hypothesis our Boolean algebra is complete. Now, for any $s\in S$, $s\leq \sup S$, and therefore $\inf\{a, s\}\leq \inf\{a, \sup S\}=q$, implying $q$ is an upper bound for $\{\inf\{a, s\}: s\in S\}$ and therefore $p\leq q$, given $p$ is the least upper bound for that set.

Reciprocally, $q=\inf\{a, \sup S\}$ implies $q\leq a$ and $q\leq \sup S$: if $q=\sup S$, this means $\sup S\leq a$, and therefore $s\leq a$ for every $s\in S$, meaning $\inf\{a, s\}=s$, for every $s\in S$ and so $p=\sup\{s: s\in S\}=\sup S$, giving us the desired equality; otherwise, $q<\sup S$, and there must exist $s\in S$ such that $q\leq s$. Since we also have $q\leq a$, $q\leq\inf\{a, s\}$, and since $p$ is an upper bound for $\{\inf\{a, s\}: s\in S\}$, $p\geq q$, what finishes the proof.
\end{proof}

An element $a$ of a Boolean algebra is an atom\index{Atom} if it is minimal in $A\setminus\{0\}$, what means that if $b\leq a$, then either $b=0$ or $b=a$; the set of atoms smaller than $a$ will be denoted by $A_{a}$; and a Boolean algebra is said to be atomic if, for every one of its elements $a$, $a=\sup A_{a}$.

Complete, atomic Boolean algebras are essentially powersets (\cite{Oosten}): if one takes, for a Boolean algebra $\mathcal{A}=(A, \leq)$, the set $A_{1}$ of all of its atoms (what is valid given all atoms are smaller than $1$), one sees $\mathcal{A}$ is equivalent to $\mathcal{P}(A_{1})$, the powerset of $A_{1}$, where an element $a\in A\setminus\{0\}$ is taken, by this isomorphism, to $A_{a}$ (and $0$ to $\emptyset$). Reciprocally, the complete, atomic Boolean algebra associated to a set $X$ is precisely $\mathcal{P}(X)$. For more information, look at Theorem $2.4$ of \cite{Oosten}.

Useful to our intents and purposes are Boolean algebras that are, simultaneously, complete, atomic and bottomless, meaning they lack a minimum element: this may seem a contradiction, given we assumed Boolean algebras to have such elements, but this can be adequately formalized.

\begin{definition}
Given a partially ordered set $\mathcal{A}=(A, \leq_{\mathcal{A}})$, we define\label{A0} 
\[\mathcal{A}_{0}=(A\cup\{0\}, \leq_{\mathcal{A}_{0}}),\]
where we assume $0\notin A$, as the poset such that $a\leq_{\mathcal{A}_{0}} b$ if and only if:

\begin{enumerate} 
\item either $a\leq_{\mathcal{A}} b$;

\item or $a=0$.
\end{enumerate}
\end{definition}

\begin{definition}
The non-empty partially ordered set $\mathcal{A}$ is a complete, atomic and bottomless Boolean algebra\index{Boolean algebra, Complete, atomic and bottomless} if, and only if, $\mathcal{A}_{0}$ is a complete, atomic Boolean algebra.
\end{definition}

Notice that, since $\mathcal{P}(\emptyset)$ only has $\emptyset$ as element, for any complete, atomic and bottomless Boolean algebra $\mathcal{A}$ we can not have $\mathcal{A}_{0}$ equivalent to $\mathcal{P}(\emptyset)$, given $\mathcal{A}$ has at least one element and therefore $\mathcal{A}_{0}$ must have at least two. This means complete, atomic and bottomless Boolean algebras correspond to the powerset of non-empty sets with $\emptyset$ removed.

\begin{proposition}\label{A0 is poset}
If $\mathcal{A}$ is a partially ordered set, so it is $\mathcal{A}_{0}$.
\end{proposition}

\begin{proof}
\begin{enumerate}
\item For any $a\in A\cup\{0\}$, we have that: if $a\in A$, since $\mathcal{A}$ is a partially ordered set, $a\leq_{\mathcal{A}}a$, and therefore $a\leq_{\mathcal{A}_{0}}a$; if $a=0$, we immediately have that $a\leq_{\mathcal{A}_{0}}a$, and therefore $\leq_{\mathcal{A}_{0}}$ is reflexive.
\item Suppose $a\leq_{\mathcal{A}_{0}}b$ and $b\leq_{\mathcal{A}_{0}}a$: if both $a$ and $b$ are in $A$, this means $a\leq_{\mathcal{A}}b$ and $b\leq_{\mathcal{A}}a$, and therefore $a=b$; if $a=0$, $b\leq_{\mathcal{A}_{0}}a$ implies that $b=0$, and therefore $a=b$; if $b=0$, $a\leq_{\mathcal{A}_{0}}b$ implies that $a=0$, and from that $a=b$; if $a=b=0$, there is nothing to be done, following from that that $\leq_{\mathcal{A}_{0}}$ is anti-symmetric.
\item Suppose $a\leq_{\mathcal{A}_{0}}b$ and $b\leq_{\mathcal{A}_{0}}c$: if $c=0$, $b\leq_{\mathcal{A}_{0}}c$ implies $b=0$, and then $a\leq_{\mathcal{A}_{0}}b$ implies $a=0$, and so $a\leq_{\mathcal{A}_{0}}c$; if $b=0$, $a\leq_{\mathcal{A}_{0}}b$ implies $a=0$, and so $a\leq_{\mathcal{A}_{0}}c$; obviously, if $a=0$ we promptly have that $a\leq_{\mathcal{A}_{0}}c$; so suppose $a, b, c\in A$ as our last case, and then $a\leq_{\mathcal{A}_{0}}b$ implies $a\leq_{\mathcal{A}}b$ while $b\leq_{\mathcal{A}_{0}}c$ implies $b\leq_{\mathcal{A}}c$, and since $\mathcal{A}$ is a partially ordered set we get that $a\leq_{\mathcal{A}}c$ and then $a\leq_{\mathcal{A}_{0}}c$, which proves $\leq_{\mathcal{A}_{0}}$ is transitive.
\end{enumerate}
\end{proof}

\begin{lemma}\label{sup of low bounds is low bound}
Given a partially ordered set $(A, \leq)$, for elements $a, b\in A$ we have that the supremum of the lower bounds of $\{a,b\}$, if it exists, is itself a lower bound for $\{a,b\}$.
\end{lemma}

\begin{proof}
Let $S$ be he supremum of the lower bounds of $\{a, b\}$: by definition, this means that for any upper bound $d$ for the set $\{c\in A:  c\leq a, c\leq b\}$ of lower bounds, $S\leq d$; but, since $a$ and $b$ are such upper bounds, we find that $S\leq a$ and $S\leq b$.
\end{proof}

\begin{theorem}\label{list}
A partially ordered set $(A, \leq)$ which satisfies all the following conditions is a complete, atomic and bottomless Boolean algebra.

\begin{enumerate} 
\item It has a maximum element $1$.
\item All non-empty subsets $S$ of $A$ have a supremum.
\item For every $a\in A$ different from $1$ there exists $b\in A$, named the complement of $a$, such that 
\[b=\inf\{c\in A:  \sup\{a, c\}=1\}\]
and
\[b=\sup\{c\in A:  \text{$\inf\{a, c\}$ does not exist}\},\]
property we call being semi-complemented.
\item Denoting by $A_{a}$ the set of minimal elements smaller than $a$, $a=\sup A_{a}$ (whenever a poset satisfies this condition, we will call it atomic).
\end{enumerate}
\end{theorem}

\begin{proof}
Suppose that $\mathcal{A}=(A, \leq_{\mathcal{A}})$ is a partially ordered set satisfying the aforementioned conditions. Since $\mathcal{A}$ is a partially ordered set, so is $\mathcal{A}_{0}$ from Proposition \ref{A0 is poset}. The maximum $1$ of $\mathcal{A}$ remains a maximum in $\mathcal{A}_{0}$, while $0$ becomes a minimum. For non-zero elements $a$ and $b$, the supremum of the two in $\mathcal{A}_{0}$ remains the same as in $\mathcal{A}$, while if $a=0$ or $b=0$ the supremum is simply the largest of the two. 

If $a$ or $b$ equals $0$, the infimum is $0$, while if $a, b\in A$ there are two cases to consider: if $\inf\{a, b\}$ was defined in $\mathcal{A}$, it remains the same in $\mathcal{A}_{0}$; if the infimum was not defined in $\mathcal{A}$, this means that there were no lower bounds for both $a$ and $b$, since otherwise we would have
\[\inf\{a,b\}=\sup\{c\in A:  c\leq_{\mathcal{A}}a\quad\text{and}\quad c\leq_{\mathcal{A}}b\}\]
by Lemma \ref{sup of low bounds is low bound}, and therefore the infimum of both in $\mathcal{A}_{0}$ is $0$. Every element $a\in A\setminus\{1\}$ already has a complement $b$ in $\mathcal{A}$ such that $b=\inf\{c\in A:  \sup\{a, c\}=1\}$ and
\[b=\sup\{c\in A:  \text{$\inf\{a, c\}$ does not exist}\};\]
of course the first equality keeps on holding in $\mathcal{A}_{0}$, while the second becomes, remembering that the non-defined infima in $\mathcal{A}$ become $0$ in $\mathcal{A}_{0}$,
\[b=\sup\{c\in A:  \inf\{a, c\}=0\};\]
the complement of $1$ is clearly $0$ and vice-versa. This, of course, proves $\mathcal{A}_{0}$ is a Boolean algebra.

Since $\mathcal{A}$ is closed under suprema of non-empty sets and $\sup\emptyset=0$ in $\mathcal{A}_{0}$, it is clear that $\mathcal{A}_{0}$ is closed under any suprema; it is clear that, for any set $S$, $\inf S=(\sup\{s^{c}: s\in S\})^{c}$, where $a^{c}$ denotes the complement of $s$. Clearly $\mathcal{A}_{0}$ remains atomic, since $\mathcal{A}$ is atomic, what finishes the proof that the previous list of conditions imply $\mathcal{A}$ is a complete, atomic and bottomless Boolean algebra.
\end{proof}

\begin{theorem}
The reciprocal of Theorem \ref{list} holds, meaning that complete, atomic and bottomless Boolean algebras satisfy the list of conditions found in that theorem.
\end{theorem}

\begin{proof}
Given a partially ordered set $\mathcal{A}$, suppose $\mathcal{A}_{0}$ is a complete, atomic Boolean algebra.

\begin{enumerate} 
\item The maximum $1$ of $\mathcal{A}_{0}$ is still a maximum in $\mathcal{A}$.
\item The supremum of any non-empty set in $\mathcal{A}$ is just its supremum in $\mathcal{A}_{0}$.
\item Given any element $a\neq 1$, its complement $b$ in $\mathcal{A}_{0}$ ends up being also its complement in $\mathcal{A}$: clearly 
\[b=\inf\{c\in A:  \sup\{a,c\}=1\}.\]
Now, $\inf\{a,c\}$ does not exist in $\mathcal{A}$ if, and only if, $\inf\{a,c\}=0$ in $\mathcal{A}_{0}$: we already proved that if $\inf\{a,c\}$ does not exist in $\mathcal{A}$ then $\inf\{a,c\}=0$ in $\mathcal{A}_{0}$, remaining to show the reciprocal; if $\inf\{a, c\}$ existed in $\mathcal{A}$, since $\inf\{a, c\}=0$ in $\mathcal{A}_{0}$ and given the unicity of the infimum one would find that $\inf\{a, c\}=0$ in $\mathcal{A}$, contradicting the fact that $0$ is not in $\mathcal{A}$. This way, we find that in $\mathcal{A}$
\[b=\sup\{c\in A:  \text{$\inf\{a,c\}$ does not exist}\},\]
as necessary.
\item Clearly $\mathcal{A}_{0}$ being atomic implies $\mathcal{A}$ is atomic.
\end{enumerate}
\end{proof}

\begin{proposition}\label{sem-inf-dist}
If $(A, \leq_{\mathcal{A}})$ is a complete, atomic and bottomless Boolean algebra, for any $S\subseteq A$, if 
\[S^{a}=\{s\in S: \text{$\inf\{a, s\}$ exists}\}\neq\emptyset,\]
then
\[\sup\{\inf\{a, s\}: s\in S^{a}\}=\inf\{a, \sup S\};\]
if $S^{a}=\emptyset$, $\inf\{a, \sup S\}$ also does not exist.
\end{proposition}

\begin{proof}
If $S^{a}=\emptyset$ this means that $\inf\{a,s\}=0$ for every $s\in S$ in $\mathcal{A}_{0}$, and therefore\\ $\inf\{a, \sup S\}=0$, so that the same infimum no longer exists in $\mathcal{A}$.

If $S^{a}\neq\emptyset$, all infima and suprema in $\sup\{\inf\{a,s\}: s\in S^{a}\}$ and $\inf\{a, \sup S\}$ exist in $\mathcal{A}$ and are therefore equal to their counterparts in $\mathcal{A}_{0}$; given $\sup\{\inf\{a, s\}: s\in S^{a}\}=\sup\{\inf\{a, s\}: s\in S\}$ in $\mathcal{A}_{0}$, since $s\in S\setminus S^{a}$ implies $\inf\{a, s\}=0$, by the infinite-distributivity of $\mathcal{A}_{0}$ one proves the desired result.
\end{proof}

To summarize the developments on this chapter so far, atomic and bottomless Bool\-ean algebra are powersets (of non-empty sets) with their empty-set removed. The relevance of removing the empty-set from consideration is that the multialgebras we will work with are not partial, meaning the empty-set is never the result of an operation.


\section{A first attempt}\label{A first attempt}

For simplicity, let us denote the set of non-empty subsets of $A$ by $\mathcal{P}^{*}(A)$, that is, $\mathcal{P}(A)\setminus\{\emptyset\}=\mathcal{P}^{*}(A)$.

Consider the categories $\textbf{Alg}(\Sigma)$\label{AlgSigma} of $\Sigma$-algebras, with homomorphisms between $\Sigma$-algebras as morphisms, and $\textbf{MAlg}(\Sigma)$ of $\Sigma$-multialgebras, with homomorphisms between $\Sigma$-multialgebras as morphisms.

Consider the transformation $\mathcal{P}:\textbf{MAlg}(\Sigma)\rightarrow\textbf{Alg}(\Sigma)$ taking:
\begin{enumerate}
\item a $\Sigma$-multialgebra $\mathcal{A}=(A, \{\sigma_{\mathcal{A}}\}_{\sigma\in\Sigma})$ to the $\Sigma$-algebra 
\[\mathcal{P}(\mathcal{A})=(\mathcal{P}^{*}(A), \{\sigma_{\mathcal{P}(\mathcal{A})}\}_{\sigma\in\Sigma}),\]
where $\mathcal{P}(A)$ is the powerset of $A$ and, for a $\sigma\in\Sigma_{n}$ and nonempty $A_{1}, \dotsc  , A_{n}\subseteq A$,
\[\sigma_{\mathcal{P}(\mathcal{A})}(A_{1}, \dotsc  , A_{n})=\bigcup_{(a_{1}, \dotsc  , a_{n})\in A_{1}\times\cdots\times A_{n}}\sigma_{\mathcal{A}}(a_{1}, \dotsc  , a_{n});\]
\item for $\mathcal{A}$ and $\mathcal{B}$ $\Sigma$-multialgebras, a homomorphism $\varphi:\mathcal{A}\rightarrow\mathcal{B}$ to the function $\mathcal{P}(\varphi):\mathcal{P}(A)\rightarrow\mathcal{P}(B)$ such that, for a $\emptyset\neq A'\subseteq A$, 
\[\mathcal{P}(\varphi)(A')=\{\varphi(a)\in B\ :\  a\in A'\}.\]
\end{enumerate}

One could hope that $\mathcal{P}(\varphi)$ is a homomorphism between $\Sigma$-algebras, and perhaps that $\mathcal{P}$ is a functor from $\textbf{MAlg}(\Sigma)$ to $\textbf{Alg}(\Sigma)$, but this does not happen: we have instead an inclusion.

\begin{lemma}\label{homomorphisms lead to almost homomorphisms}
For $\mathcal{A}$ and $\mathcal{B}$ two $\Sigma$-multialgebras and $\varphi:\mathcal{A}\rightarrow\mathcal{B}$ a homomorphism, $\mathcal{P}(\varphi)$ satisfies
\[\mathcal{P}(\varphi)(\sigma_{\mathcal{P}(\mathcal{A})}(A_{1}, \dotsc  , A_{n}))\subseteq \sigma_{\mathcal{P}(\mathcal{B})}(\mathcal{P}(\varphi)(A_{1}), \dotsc  , \mathcal{P}(\varphi)(A_{n}))\]
for all $\sigma\in\Sigma$ and nonempty $A_{1}, \dotsc  , A_{n}\subseteq A$; if $\varphi$ is a full homomorphism, $\mathcal{P}(\varphi)$ is a homomorphism.

\end{lemma}

\begin{proof}
Given $\sigma\in \Sigma_{n}$ and nonempty $A_{1}, \dotsc  , A_{n}\subseteq A$, we have that
\[\sigma_{\mathcal{P}(\mathcal{B})}(\mathcal{P}(\varphi)(A_{1}), \dotsc  , \mathcal{P}(\varphi)(A_{n}))=\bigcup_{(b_{1}, \dotsc  , b_{n})\in \mathcal{P}(\varphi)(A_{1})\times\cdots\times \mathcal{P}(\varphi)(A_{n})}\sigma_{\mathcal{B}}(b_{1}, \dotsc  , b_{n})=\]
\[\bigcup_{(b_{1}, \dotsc  , b_{n})\in \{\varphi(a)\ :\  a\in A_{1}\}\times\cdots\times \{\varphi(a)\ :\  a\in A_{n}\}}\sigma_{\mathcal{B}}(b_{1}, \dotsc  , b_{n})=\bigcup_{(a_{1}, \dotsc  , a_{n})\in A_{1}\times\cdots\times A_{n}}\sigma_{\mathcal{B}}(\varphi(a_{1}), \dotsc  , \varphi(a_{n})),\]
which clearly contains
\[\bigcup_{(a_{1}, \dotsc  , a_{n})\in A_{1}\times\cdots\times A_{n}}\{\varphi(a)\ :\  a\in\sigma_{\mathcal{A}}(a_{1}, \dotsc  , a_{n})\}=\{\varphi(a)\ :\  a\in \bigcup_{(a_{1}, \dotsc  , a_{n})\in A_{1}\times\cdots\times A_{n}}\sigma_{\mathcal{A}}(a_{1}, \dotsc  , a_{n})\}=\]
\[\{\varphi(a)\ :\  a\in \sigma_{\mathcal{P}(\mathcal{A})}(A_{1}, \dotsc  , A_{n})\}=\mathcal{P}(\varphi)(\sigma_{\mathcal{P}(\mathcal{A})}(A_{1}, \dotsc  , A_{n})),\]
so that $\mathcal{P}(\varphi)$ satisfies the required property; if $\varphi$ is a full homomorphism, $\sigma_{\mathcal{B}}(\varphi(a_{1}), \dotsc  , \varphi(a_{n}))=\{\varphi(a)\ :\  a\in\sigma_{\mathcal{A}}(a_{1}, \dotsc  , a_{n})\}$, and the inclusion between the equations above becomes an equality.
\end{proof}

So, let us restrict $\mathcal{P}$ for a moment to the category $\textbf{MAlg}_{=}(\Sigma)$, of $\Sigma$-multialgebras with only full homomorphisms between them as morphisms, and let us call this new transformation $\mathcal{P}_{=}:\textbf{MAlg}_{=}(\Sigma)\rightarrow \textbf{Alg}(\Sigma)$.

\begin{proposition}
$\mathcal{P}_{=}$ is, in fact, a functor.
\end{proposition}

\begin{proof}
\begin{enumerate}
\item Given a $\Sigma$-multialgebra $\mathcal{A}$, let $Id_{\mathcal{A}}:\mathcal{A}\rightarrow\mathcal{A}$ be the identity homomorphism. For any $\emptyset\neq A'\subseteq A$, we have that 
\[\mathcal{P}_{=}(Id_{\mathcal{A}})(A')=\{Id_{\mathcal{A}}(a)\ :\  a\in A'\}=\{a\ :\  a\in A'\}=A',\]
and therefore $\mathcal{P}_{=}(Id_{\mathcal{A}})$ is the identity homomorphism on $\mathcal{P}^{*}(\mathcal{A})$.

\item Let $\mathcal{A}$, $\mathcal{B}$ and $\mathcal{C}$ be $\Sigma$-multialgebras and $\varphi:\mathcal{A}\rightarrow\mathcal{B}$ and $\psi:\mathcal{B}\rightarrow\mathcal{C}$ be full homomorphisms between them. We have that, for every $\emptyset\neq A'\subseteq A$, 
\[\mathcal{P}_{=}(\psi)\circ\mathcal{P}_{=}(\varphi)(A')=\mathcal{P}_{=}(\psi)(\{\varphi(a)\in B\ :\  a\in A'\})=\{\psi(\varphi(a))\in C\ :\  a\in A'\}=\]
\[\mathcal{P}_{=}(\psi\circ\varphi)(A'),\]
and so $\mathcal{P}_{=}(\psi)\circ\mathcal{P}_{=}(\varphi)=\mathcal{P}_{=}(\psi\circ\varphi)$.
\end{enumerate}
\end{proof}

What we actually want is an equivalence of categories, but we can be certain that $\mathcal{P}_{=}$ will not provide this equivalence: unfortunately, $\mathcal{P}_{=}$ is not injective on objects. Take the signature $\Sigma_{s}$ with a single unary operator $s$, and consider the $\Sigma$-multialgebras $\mathcal{A}=(\{0,1\}, \{s_{\mathcal{A}}\})$ and $\mathcal{B}=(\{0,1\}, \{s_{\mathcal{B}}\})$ such that: $s_{\mathcal{A}}(0)=\{1\}$ and $s_{\mathcal{A}}(1)=\{1\}$; and $s_{\mathcal{B}}(0)=s_{\mathcal{B}}(1)=\{0,1\}$.

\begin{figure}[H]
\centering
\begin{minipage}[t]{4cm}
\centering
\begin{tikzcd}
    0 \arrow[r, "s_{\mathcal{A}}"]  & 1 \arrow[loop right, out=30, in=-30, distance=3em]{}{s_{\mathcal{A}}}
  \end{tikzcd}
\caption*{The $\Sigma_{s}$-multialgebra $\mathcal{A}$}
\end{minipage}
\hspace{3cm}
\centering
\begin{minipage}[t]{4cm}
\centering
\begin{tikzcd}
    0 \arrow[rr, bend left=50, "s_{\mathcal{B}}"]\arrow[loop right, out=210, in=150, distance=3em]{}{s_{\mathcal{B}}}  && 1 \arrow[ll, bend left=50, "s_{\mathcal{B}}"]\arrow[loop right, out=30, in=-30, distance=3em]{}{s_{\mathcal{B}}}
  \end{tikzcd}
\caption*{The $\Sigma_{s}$-multialgebra $\mathcal{B}$}
\end{minipage}
\end{figure}

Clearly the two of then are not isomorphic, given that the result of an operation in $\mathcal{A}$ always has cardinality $1$ and in $\mathcal{B}$ always has cardinality $2$.

However, we have that $s_{\mathcal{P}_{=}(\mathcal{A})}(\{0\})$, $s_{\mathcal{P}_{=}(\mathcal{A})}(\{1\})$, and $s_{\mathcal{P}_{=}(\mathcal{A})}(\{0, 1\})$ all equal $\{1\}$, while $s_{\mathcal{P}_{=}(\mathcal{B})}(\{0\})$, $s_{\mathcal{P}_{=}(\mathcal{B})}(\{1\})$, and $s_{\mathcal{P}_{=}(\mathcal{B})}(\{0, 1\})$ all equal $\{0,1\}$.

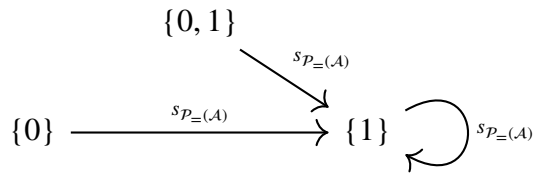
\begin{figure}[H]
\centering
\begin{tikzcd}
  & \{0,1\} \arrow[dr, "s_{\mathcal{P}_{=}(\mathcal{A})}"] &\\
\{0\}\arrow[rr, "s_{\mathcal{P}_{=}(\mathcal{A})}"] & & \{1\}\arrow[loop right, out=30, in=-30, distance=3em]{}{s_{\mathcal{P}_{=}(\mathcal{A})}}
  \end{tikzcd}
\caption*{The $\Sigma_{s}$-multialgebra $\mathcal{P}_{=}(\mathcal{A})$}
\end{figure}

\begin{figure}[H]
\centering
\begin{tikzcd}
  & \{0,1\} \arrow[loop, swap, "s_{\mathcal{P}_{=}(\mathcal{B})}"] &\\
\{0\}\arrow[ur, "s_{\mathcal{P}_{=}(\mathcal{B})}"] & & \{1\}\arrow[ul, swap, "s_{\mathcal{P}_{=}(\mathcal{B})}"]
  \end{tikzcd}
\caption*{The $\Sigma_{s}$-multialgebra $\mathcal{P}_{=}(\mathcal{B})$}
\end{figure}
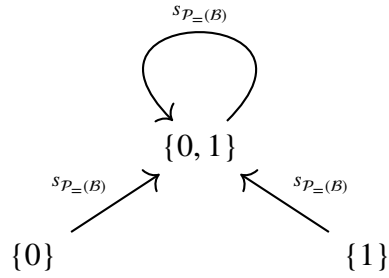

Taking the function $\varphi:\mathcal{P}^{*}(A)\rightarrow\mathcal{P}^{*}(B)$ such that $\varphi(\{0\})=\{0\}$, $\varphi(\{1\})=\{0,1\}$, and $\varphi(\{0,1\})=\{1\}$, we see that it is a bijection and a homomorphism, and therefore $\varphi:\mathcal{P}_{=}(\mathcal{A})\rightarrow\mathcal{P}_{=}(\mathcal{B})$ is an isomorphism.

Curiously, we can also quite easily show that $\mathcal{P}_{=}$ does not have a right adjoint: if it did, that functor would itself posses a left adjoint, and therefore be would be continuous, preserving, among other constructions, terminal objects; and, while $\textbf{Alg}(\Sigma)$ has as terminal objects the $\Sigma$-algebras with one element (where all operations return that very element), $\textbf{MAlg}_{=}(\Sigma)$ does not have them.\footnote{Assuming, what is quite reasonable, that $\Sigma$ is not empty.} In fact, for any $\Sigma$-multialgebra $\mathcal{A}$, by taking a cardinal $\kappa$ greater than the cardinality of the universe of $\mathcal{A}$ we see that there can be no full homomorphism from $\textbf{mF}(\Sigma, \mathcal{V}, \kappa)$ to $\mathcal{A}$.


\section{An improvement}

The problem with our definition of $\mathcal{P}_{=}$ is that it disregards the structure of the universe of $\mathcal{P}(\mathcal{A})$, which admits an order. So, we change our target category to reflect this structure, allowing the objects to carry with them such orderings.

\begin{definition}
Given a signature $\Sigma$, a $(\Sigma, \leq)$-algebra\index{Algebra, $(\Sigma, \leq)$-} $\mathcal{A}$ is a triple 
\[(A, \{\sigma_{\mathcal{A}}\}_{\sigma\in\Sigma}, \leq_{\mathcal{A}})\]
such that:
\begin{enumerate}
\item $(A, \{\sigma_{\mathcal{A}}\}_{\sigma\in\Sigma})$ is a $\Sigma$-algebra;
\item $(A,\leq_{\mathcal{A}})$ is a complete, atomic and bottomless Boolean algebra;
\item if $A_{a}$ is the set of minimal elements of $(A, \leq_{\mathcal{A}})$ (atoms) below $a$, for all $\sigma\in\Sigma_{n}$ and $a_{1}, \dotsc  , a_{n}$ we have that
\[\sigma_{\mathcal{A}}(a_{1}, \dotsc  , a_{n})=\sup\{\sigma_{\mathcal{A}}(b_{1}, \dotsc  , b_{n})\ :\  (b_{1}, \dotsc  , b_{n})\in A_{a_{1}}\times\cdots\times A_{a_{n}}\}.\]
\end{enumerate}
\end{definition}

\begin{proposition}
For $\mathcal{A}$ a $(\Sigma, \leq)$-algebra, $\sigma\in\Sigma_{n}$ and $a_{1}, \dotsc  , a_{n}, b_{1}, \dotsc  , b_{n}\in A$ such that $a_{1}\leq_{\mathcal{A}} b_{1}, \dotsc  , a_{n}\leq_{\mathcal{A}} b_{n}$,
\[\sigma_{\mathcal{A}}(a_{1}, \dotsc  , a_{n})\leq_{\mathcal{A}} \sigma_{\mathcal{A}}(b_{1}, \dotsc  , b_{n}).\]
\end{proposition}

\begin{proof}
Since, for every $i\in\{1, \dotsc  , n\}$, $a_{i}\leq_{\mathcal{A}}b_{i}$, we have that $A_{a_{i}}\subseteq A_{b_{i}}$, and therefore $A_{a_{1}}\times\cdots\times A_{a_{n}}\subseteq A_{b_{1}}\times\cdots\times A_{b_{n}}$; this way,
\[\sigma_{\mathcal{A}}(a_{1}, \dotsc  , a_{n})=\sup\{\sigma_{\mathcal{A}}(c_{1}, \dotsc  , c_{n})\ :\  (c_{1}, \dotsc  , c_{n})\in A_{a_{1}}\times\cdots\times A_{a_{n}}\}\leq_{\mathcal{A}}\]
\[\sup\{\sigma_{\mathcal{A}}(c_{1}, \dotsc  , c_{n})\ :\  (c_{1}, \dotsc  , c_{n})\in A_{b_{1}}\times\cdots\times A_{b_{n}}\}=\sigma_{\mathcal{A}}(b_{1}, \dotsc  , b_{n}).\]
\end{proof}

For a $\Sigma$-multialgebra $\mathcal{A}=(A, \{\sigma_{\mathcal{A}}\}_{\sigma\in\sigma})$, we define $\mathcal{P}_{\mathcal{B}}(\mathcal{A})$\label{Functor PB} as the $(\Sigma, \leq)$-algebra 
\[(\mathcal{P}^{*}(A), \{\sigma_{\mathcal{P}_{\mathcal{B}}(\mathcal{A})}\}_{\sigma\in\Sigma}, \leq_{\mathcal{P}_{\mathcal{B}}(\mathcal{A})})\]
such that $(\mathcal{P}^{*}(A), \{\sigma_{\mathcal{P}_{\mathcal{B}}(\mathcal{A})}\}_{\sigma\in\Sigma})$ is exactly $\mathcal{P}(\mathcal{A})$ (defined in the beginning of Section \ref{A first attempt}) and, for nonempty subsets $A_{1}$ and $A_{2}$ of $A$, $A_{1}\leq_{\mathcal{P}_{\mathcal{B}}(\mathcal{A})} A_{2}$ if and only if $A_{1}\subseteq A_{2}$. Since:
\begin{enumerate}
\item $\mathcal{P}(\mathcal{A})$ is a $\Sigma$-algebra;
\item $(\mathcal{P}^{*}(A), \leq_{\mathcal{P}_{\mathcal{B}}(\mathcal{A})})$ is a complete, atomic and bottomless Boolean algebra (since $\mathcal{P}(A)$ is a complete, atomic Boolean algebra);
\item and, for $\sigma\in\Sigma_{n}$ and $\emptyset\neq A_{1}, \dotsc  , A_{n}\subseteq A$, since the atoms of $A_{i}$ are exactly $A_{A_{i}}=\{\{a\}\ :\  a\in A_{i}\}$,
\[\sigma_{\mathcal{P}_{\mathcal{B}}(\mathcal{A})}(A_{1}, \dotsc  , A_{n})=\bigcup_{(a_{1}, \dotsc  , a_{n})\in A_{1}\times\cdots\times A_{n}}\sigma_{\mathcal{A}}(a_{1}, \dotsc  , a_{n})=\]
\[\bigcup_{(\{a_{1}\}, \dotsc  , \{a_{n}\})\in A_{A_{1}}\times\cdots\times A_{A_{n}}}\sigma_{\mathcal{P}_{\mathcal{B}}(\mathcal{A})}(\{a_{1}\}, \dotsc  , \{a_{n}\});\]
\end{enumerate}
we truly have that $\mathcal{P}_{\mathcal{B}}(\mathcal{A})$ is a $(\Sigma, \leq)$-algebra.

\begin{definition}
Given $(\Sigma, \leq)$-algebras $\mathcal{A}=(A, \{\sigma_{\mathcal{A}}\}_{\sigma\in\Sigma}, \leq_{\mathcal{A}})$ and $\mathcal{B}=(B, \{\sigma_{\mathcal{B}}\}_{\sigma\in\Sigma}, \leq_{\mathcal{B}})$, a function $\varphi:A\rightarrow B$ is said to be a $(\Sigma, \leq)$-homomorphism\index{Homomorphism, $(\sigma, \leq)$-}, in which case we write $\varphi:\mathcal{A}\rightarrow\mathcal{B}$, when:
\begin{enumerate}
\item for all $\sigma\in\Sigma_{n}$ and $a_{1}, \dotsc  , a_{n}\in A$ we have that 
\[\varphi(\sigma_{\mathcal{A}}(a_{1}, \dotsc  , a_{n}))\leq_{\mathcal{B}}\sigma_{\mathcal{B}}(\varphi(a_{1}), \dotsc  , \varphi(a_{n}));\]
\item $\varphi$ is continuous, that is, for every non-empty subset $A'\subseteq A$, 
\[\varphi(\sup A')=\sup\{\varphi(a)\ :\  a\in A'\};\]
\item $\varphi$ maps minimal elements of $(A, \leq_{\mathcal{A}})$ to minimal elements of $(B, \leq_{\mathcal{B}})$.
\end{enumerate}
\end{definition}

Notice that a $(\Sigma, \leq)-$homomorphism is essentially an "almost $\Sigma-$homomorphism" which is also continuous and minimal-elements-preserving; notice too that a $(\Sigma, \leq)$-homomor\-phism is order preserving: if $a\leq_{\mathcal{A}}b$, then $b=\sup\{a,b\}$, and therefore $\varphi(b)=\sup\{\varphi(a), \varphi(b)\}$, meaning that $\varphi(a)\leq_{\mathcal{B}}\varphi(b)$.


\subsection{$\mathcal{P}_{\mathcal{B}}$ is a functor}

\begin{proposition}
When we take as objects all $(\Sigma, \leq)$-algebras and as morphisms all the $(\Sigma, \leq)$-homomorphisms between them, the resulting object is a category, denoted by $\textbf{Alg}_{\mathcal{B}}(\Sigma)$\label{AlgBSigma}.
\end{proposition}

\begin{proof}
It is enough to show the composition of $(\Sigma, \leq)$-homomorphisms returns a $(\Sigma, \leq)$-homo\-morphism, that this composition is associative and has an identity element for every object in the category. Let $\mathcal{A}$, $\mathcal{B}$, $\mathcal{C}$ and $\mathcal{D}$ be $(\Sigma, \leq)-$algebras and $\varphi:\mathcal{A}\rightarrow\mathcal{B}$, $\psi:\mathcal{B}\rightarrow\mathcal{C}$ and $\chi:\mathcal{C}\rightarrow\mathcal{D}$ be $(\Sigma, \leq)$-homomorphisms.

\begin{enumerate}
\item\begin{enumerate}
\item $\psi\circ\varphi$ obviously is a function from $A$ to $C$, so let $\sigma\in\Sigma_{n}$ and $a_{1}, \dotsc  , a_{n}\in A$: we have that, since $\varphi$ is a $(\Sigma, \leq)$-homomorphism, $\varphi(\sigma_{\mathcal{A}}(a_{1}, \dotsc  , a_{n}))\leq_{\mathcal{B}}\sigma_{\mathcal{B}}(\varphi(a_{1}), \dotsc  , \varphi(a_{n}))$, and since $\psi$ is order-preserving
\[\psi\circ\varphi(\sigma_{\mathcal{A}}(a_{1}, \dotsc  , a_{n}))=\psi(\varphi(\sigma_{\mathcal{A}}(a_{1}, \dotsc  , a_{n})))\leq_{\mathcal{C}}\psi(\sigma_{\mathcal{B}}(\varphi(a_{1}), \dotsc  , \varphi(a_{n})));\]
and since $\psi$ is a $(\Sigma, \leq)$-homomorphism,
\[\psi(\sigma_{\mathcal{B}}(\varphi(a_{1}), \dotsc  , \varphi(a_{n})))\leq_{\mathcal{C}}\sigma_{\mathcal{C}}(\psi(\varphi(a_{1})), \dotsc  , \psi(\varphi(a_{n})))=\sigma_{\mathcal{C}}(\psi\circ\varphi(a_{1}), \dotsc  , \psi\circ\varphi(a_{n})).\]

\item Given a non-empty $A'\subseteq A$, we have that $\varphi(\sup A')=\sup\{\varphi(a)\ :\  a\in A'\}$ and, denoting $\{\varphi(a)\ :\  a\in A'\}$ by $B'$, we have that $\psi(\sup B')=\sup\{\psi(b)\ :\  b\in B'\}$; since $\sup B'=\varphi(\sup A')$, we obtain
\[\psi\circ\varphi(\sup A')=\sup\{\psi(b)\ :\  b\in B'\}=\sup\{\psi\circ\varphi(a)\ :\  a\in A'\},\]
which means that $\psi\circ\varphi$ is continuous.

\item Finally, if $a\in A$ is a minimal element, $\varphi(a)\in B$ is a minimal element, since $\varphi$ preserves minimal elements, and for the same reason $\psi\circ\varphi(a)=\psi(\varphi(a))\in C$ remains a minimal element still, and from all of the above $\psi\circ\varphi$ is a $(\Sigma, \leq)$-homomorphism.
\end{enumerate}

\item It is clear that $\chi\circ(\psi\circ\varphi)=(\chi\circ\psi)\circ\varphi$ since, as functions, both are the same.

\item Consider the identity function $Id_{\mathcal{A}}:A\rightarrow A$: it is a $(\Sigma, \leq)$-homomorphism, since for $n\in\Sigma_{n}$ and $a_{1}, \dotsc  , a_{n}\in A$,
\[Id_{\mathcal{A}}(\sigma_{\mathcal{A}}(a_{1}, \dotsc  , a_{n}))=\sigma_{\mathcal{A}}(a_{1}, \dotsc  , a_{n})=\sigma_{\mathcal{A}}(Id_{\mathcal{A}}(a_{1}), \dotsc  , Id_{\mathcal{A}}(a_{n}));\]
for a non-empty $A'\subseteq A$, 
\[Id_{\mathcal{A}}(\sup A')=\sup A'=\sup\{Id_{\mathcal{A}}(a)\ :\  a\in A'\};\]
and for every minimal element $a$ of $\mathcal{A}$, $Id_{\mathcal{A}}(a)=a$ clearly remains a minimal element. Furthermore, $Id_{\mathcal{A}}$ is clearly an identity for the composition, and is therefore an identity morphism.
\end{enumerate}

\end{proof}

So the transformation taking a $\Sigma$-multialgebra $\mathcal{A}$ to $\mathcal{P}_{\mathcal{B}}(\mathcal{A})$ and a homomorphism $\varphi:\mathcal{A}\rightarrow\mathcal{B}$ to the $(\Sigma, \leq)$-homomorphism $\mathcal{P}_{\mathcal{B}}(\varphi):\mathcal{P}_{\mathcal{B}}(\mathcal{A})\rightarrow\mathcal{P}_{\mathcal{B}}(\mathcal{B})$ such that, for an $\emptyset\neq A'\subseteq A$,
\[\mathcal{P}_{\mathcal{B}}(\varphi)(A')=\{\varphi(a)\in B\ :\  a\in A'\},\]
is a functor, of the form $\mathcal{P}_{\mathcal{B}}:\textbf{MAlg}(\Sigma)\rightarrow \textbf{Alg}_{\mathcal{B}}(\Sigma)$. First we must show that $\mathcal{P}_{\mathcal{B}}(\varphi)$ is, in fact, a $(\Sigma, \leq)$-homomorphism: for $\sigma\in\Sigma_{n}$ and $\emptyset\neq A_{1}, \dotsc  , A_{n}\subseteq A$, as we saw in Lemma \ref{homomorphisms lead to almost homomorphisms}, 
\[\mathcal{P}(\varphi)(\sigma_{\mathcal{P}(\mathcal{A})}(A_{1}, \dotsc  , A_{n}))\subseteq \sigma_{\mathcal{P}(\mathcal{B})}(\mathcal{P}(\varphi)(A_{1}), \dotsc  , \mathcal{P}(\varphi)(A_{n})),\]
and since $\mathcal{P}(\varphi)=\mathcal{P}_{\mathcal{B}}(\varphi)$, we have that $\mathcal{P}_{\mathcal{B}}(\varphi)$ satisfies the first condition for being a $(\Sigma, \leq)-$\\homomorphism; and, if $\emptyset \neq A''$ is a subset of $\mathcal{P}(A)$, we have that
\[\mathcal{P}_{\mathcal{B}}(\varphi)(\sup A'')=\{\varphi(a)\ :\  a\in \sup A''\}=\{\varphi(a)\ :\  a\in \bigcup_{A'\in A''} A'\}=\]
\[\bigcup_{A'\in A''}\{\varphi(a)\ :\  a\in A'\}=\bigcup_{A'\in A''}\mathcal{P}_{\mathcal{B}}(\varphi)(A')=\sup \{\mathcal{P}_{\mathcal{B}}(\varphi)(A')\ :\  A'\in A''\},\]
what proves the satisfaction of the second condition; for the third condition, we remember that the minimal elements of $(\mathcal{P}(A)\setminus\emptyset, \subseteq)$ are the singletons, that is, sets of the form $\{a\}$ with $a\in A$, and since $\mathcal{P}_{\mathcal{B}}(\varphi)(\{a\})=\{\varphi(a)\}$, $\mathcal{P}_{\mathcal{B}}(\varphi)$ preserves minimal elements.

Now we show that $\mathcal{P}_{\mathcal{B}}$ is, in fact, a functor.
\begin{enumerate}
\item Let $\mathcal{A}$ be a $\Sigma$-multialgebra and $Id_{\mathcal{A}}$ be its identity homomorphism: we then have that, for every $\emptyset\neq A'\subseteq A$, $\mathcal{P}_{\mathcal{B}}(Id_{\mathcal{A}})(A')$ equals $\{Id_{\mathcal{A}}(a)\ :\  a\in A'\}=A'$, and therefore $\mathcal{P}_{\mathcal{B}}(Id_{\mathcal{A}})$ is the identity $(\Sigma, \leq)$-homomorphism on $\mathcal{P}_{\mathcal{B}}(\mathcal{A})$.

\item Let $\mathcal{A}$, $\mathcal{B}$ and $\mathcal{C}$ be $\Sigma$-multialgebras and $\varphi:\mathcal{A}\rightarrow\mathcal{B}$ and $\psi:\mathcal{B}\rightarrow\mathcal{C}$ be homomorphisms: for every $\emptyset\neq A'\subseteq A$, we have that
\[\mathcal{P}_{\mathcal{B}}(\psi\circ\varphi)(A')=\{\psi\circ\varphi(a)\ :\  a\in A'\}=\mathcal{P}_{\mathcal{B}}(\psi)(\{\varphi(a)\ :\  a\in A'\})=\]
\[\mathcal{P}_{\mathcal{B}}(\psi)(\mathcal{P}_{\mathcal{B}}(\varphi)(A'))=\mathcal{P}_{\mathcal{B}}(\psi)\circ\mathcal{P}_{\mathcal{B}}(\varphi)(A').\]
\end{enumerate}

\begin{theorem}
$\mathcal{P}_{\mathcal{B}}:\textbf{MAlg}(\Sigma)\rightarrow \textbf{Alg}_{\mathcal{B}}(\Sigma)$ is a functor.
\end{theorem}

\subsection{$\mathcal{P}_{\mathcal{B}}$ is part of a monad}

Remember that, for functors $F$ and $G$, both from the category $B$ to the category $C$, a natural transformation $\eta$ from $F$ to $G$ (what we may denote by $\eta:F\rightarrow G$) is a  collection of morphisms $\eta_{U}$ of $C$ indexed by the objects $U$ of $B$ such that: for objects $U$ and $V$ of $B$, $\eta_{U}$ and $\eta_{V}$ are morphisms from, respectively, $F(U)$ and $F(V)$ to, also respectively, $G(U)$ and $G(V)$ such that, for any morphism $f:U\rightarrow V$,
\[\eta_{V}\circ F(f)=G(f)\circ\eta_{U}.\]

For categories $B$, $C$ and $D$, and $\eta:F\rightarrow G$ a natural transformation between functors $F, G:B\rightarrow C$, given a third functor $H:C\rightarrow D$ we denote by $H\eta$ the natural transformation between $H\circ F$ and $H\circ G$ given by, in an object $X$ of $C$, 
\[(H\eta)_{X}=H\eta_{X};\]
given yet a fourth category $E$ and a functor $I:E\rightarrow B$, we also define the natural transformation $\eta I:F\circ I\rightarrow G\circ I$ given by, for an object $Y$ of $E$,
\[(\eta I)_{Y}=\eta_{I(Y)}.\]

\begin{definition}
A monad\index{Monad} (\cite{CWM}), in a category $C$, is an endofunctor $P:C\rightarrow C$ together with natural transformations $\eta:1_{C}\rightarrow P$ and $\epsilon:P\circ P\rightarrow P$, where $1_{C}$ is the identity functor on $C$, such that the following conditions hold:
\begin{enumerate}
\item $\epsilon\circ P\epsilon=\epsilon\circ \epsilon P$;
\item $\epsilon\circ P\eta=\epsilon\circ\eta P=1_{P}$, where $1_{P}$ denotes the identity natural transformation on $P$.
\end{enumerate}
\end{definition}

The conditions required of a monad are equivalent to the commutativity of the following diagrams.

\begin{figure}[H]
\centering
\begin{minipage}[t]{4cm}
\centering
\begin{tikzcd}[row sep=large,column sep=large]
   P^{3} \arrow[rightarrow]{r}{P\epsilon}\arrow[rightarrow]{d}{\epsilon P}  & P^{2} \arrow[rightarrow]{d}{\epsilon}\\
   P^{2}\arrow[rightarrow]{r}{\epsilon} & P
  \end{tikzcd}
\end{minipage}
\hspace{3cm}
\centering
\begin{minipage}[t]{4cm}
\centering
\begin{tikzcd}[row sep=large,column sep=large]
P \arrow[rightarrow]{r}{\eta P}\arrow[rightarrow]{d}{T\eta}\arrow[equal]{dr} & P^{2} \arrow[rightarrow]{d}{\epsilon}\\
P^{2}\arrow[rightarrow]{r}{\epsilon} & P
  \end{tikzcd}
\end{minipage}
\end{figure}

Now, we state that the functors $\mathcal{P}_{\mathcal{B}}$ and $\mathcal{P}$ may be adequately adapted as to be part of a monad. Let $P:\textbf{MAlg}(\Sigma)\rightarrow\textbf{MAlg}(\Sigma)$\label{P} be the functor taking: a $\Sigma$-multialgebra $\mathcal{A}=(A, \{\sigma_{\mathcal{A}}\}_{\sigma\in\Sigma})$ into the $\Sigma$-multialgebra $P\mathcal{A}=(\mathcal{P}(A)\setminus\{\emptyset\}, \{\sigma_{P\mathcal{A}}\}_{\sigma\in\Sigma})$ such that, for non-empty subsets $A_{1}, \dotsc  , A_{n}$ of $A$,
\[\sigma_{P\mathcal{A}}(A_{1}, \dotsc  , A_{n})=\{\{a\}\ :\  a\in \bigcup_{(a_{1}, \dotsc  , a_{n})\in A_{1}\times\cdots\times A_{n}}\sigma_{\mathcal{A}}(a_{1}, \dotsc  , a_{n})\};\]
and, for $\Sigma$-multialgebras $\mathcal{A}$ and $\mathcal{B}$, a homomorphism $\varphi:\mathcal{A}\rightarrow\mathcal{B}$ to the function $P\varphi:\mathcal{P}(A)\setminus\{\emptyset\}\rightarrow\mathcal{P}(B)\setminus\{\emptyset\}$ such that, for a non-empty $A'\subseteq A$, 
\[P\varphi(A')=\{\varphi(a)\in B\ :\  a\in A'\}.\]
Notice, first of all, that $P\varphi$ is, indeed, a homomorphism from $P\mathcal{A}$ to $P\mathcal{B}$, whenever $\varphi$ is a homomorphism from $\mathcal{A}$ to $\mathcal{B}$: for any non-empty subsets $A_{1}$ through $A_{n}$ of $A$ and a $\sigma$ of arity $n$, one has that, for every $(a_{1}, \dotsc  , a_{n})\in A_{1}\times\cdots\times A_{n}$, $\{\varphi(a)\ :\  a\in \sigma_{\mathcal{A}}(a_{1}, \dotsc  , a_{n})\}\subseteq\sigma_{\mathcal{B}}(\varphi(a_{1}), \dotsc  , \varphi(a_{n}))$ since $\varphi$ is a homomorphism, and therefore 
\[\{P\varphi(A')\ :\  A'\in \sigma_{P\mathcal{A}}(A_{1}, \dotsc  , A_{n})\}=\{P\varphi(\{a\})\ :\  a\in\bigcup_{(a_{1}, \dotsc  , a_{n})\in A_{1}\times\cdots\times A_{n}}\sigma_{\mathcal{A}}(a_{1}, \dotsc  , a_{n})\}=\]
\[\{\{\varphi(a)\}\ :\  a\in\bigcup_{(a_{1}, \dotsc  , a_{n})\in A_{1}\times\cdots\times A_{n}}\sigma_{\mathcal{A}}(a_{1}, \dotsc  , a_{n})\}=\]
\[\bigcup_{(a_{1}, \dotsc  , a_{n})\in A_{1}\times\cdots\times A_{n}}\{\{\varphi(a)\} \ :\ a\in \sigma_{\mathcal{A}}(a_{1}, \dotsc  , a_{n})\}\subseteq\]
\[\bigcup_{(a_{1}, \dotsc  , a_{n})\in A_{1}\times\cdots\times A_{n}}\{\{b\}\ :\  b\in \sigma_{\mathcal{B}}(\varphi(a_{1}), \dotsc  , \varphi(a_{n}))\}=\]
\[\bigcup_{(b_{1}, \dotsc  , b_{n})\in \{\varphi(a_{1})\ :\  a_{1}\in A_{1}\}\times\cdots\times\{\varphi(a_{n})\ :\  a_{n}\in A_{n}\}}\{\{b\}\ :\  b\in \sigma_{\mathcal{B}}(b_{1}, \dotsc  , b_{n})\}=\]
\[\bigcup_{(b_{1}, \dotsc  , b_{n})\in P\varphi(A_{1})\times\cdots\times P\varphi(A_{n})}\{\{b\}\ :\  b\in \sigma_{\mathcal{B}}(b_{1}, \dotsc  , b_{n})\}=\]
\[\{\{b\}\ :\  b\in\bigcup_{(b_{1}, \dotsc  , b_{n})\in P\varphi(A_{1})\times\cdots\times P\varphi(A_{n})}\sigma_{\mathcal{B}}(b_{1}, \dotsc  , b_{n})\}=\sigma_{P\mathcal{B}}(P\varphi(A_{1}), \dotsc  , P\varphi(A_{n})).\]
It is easy to see that the identity homomorphism of $\mathcal{A}$ is mapped, by $P$, to the identity homomorphism of $P\mathcal{A}$, and given homomorphisms $\varphi:\mathcal{A}\rightarrow\mathcal{B}$ and $\psi:\mathcal{B}\rightarrow\mathcal{C}$, we see that, for a non-empty subset $A'$ of $A$, one has
\[P(\psi\circ\varphi)(A')=\{\psi\circ\varphi(a)\ :\  a\in A'\}=\{\psi(b)\ :\  b\in P\varphi(A')\}=(P\psi\circ P\varphi)(A'),\]
what finishes demonstrating that $P$ is indeed an endofunctor of $\textbf{MAlg}(\Sigma)$. 

$P$ is essentially equal to both $\mathcal{P}$ and $\mathcal{P}_{\mathcal{B}}$, given that $\Sigma$-algebras are multialgebras whose operations are always singletons and, disregarding its order, a $(\Sigma, \leq)$-algebra is nothing but a $\Sigma$-algebra; the difference here is that, while the operations in $\mathcal{P}(\mathcal{A})$ return sets of elements of $\mathcal{A}$, the operations in $P\mathcal{A}$ returns sets of singletons of elements of $\mathcal{A}$, precisely those elements obtained in the operation as performed in $\mathcal{P}(\mathcal{A})$.

Now, to form a monad we need adequate natural transformations $\eta:1_{\textbf{MAlg}(\Sigma)}\rightarrow P$ and $\epsilon:P\circ P\rightarrow P$, or what is equivalent, homomorphisms $\eta_{\mathcal{A}}:\mathcal{A}\rightarrow P\mathcal{A}$ and $\epsilon_{\mathcal{A}}:PP\mathcal{A}\rightarrow P\mathcal{A}$, for every $\Sigma$-multialgebra $\mathcal{A}$. And we simply take the obvious candidates $\eta_{\mathcal{A}}(a)=\{a\}$ and, for a non-empty set of non-empty subsets $\{A_{i}\}_{i\in I}$ of $A$, $\epsilon_{\mathcal{A}}(\{A_{i}\}_{i\in I})=\bigcup_{i\in I}A_{i}$.

\begin{proposition}
For any $\Sigma$-multialgebra $\mathcal{A}$, $\eta_{\mathcal{A}}$ and $\epsilon_{\mathcal{A}}$, as previously defined, are full homomorphisms.
\end{proposition}

\begin{proof}
\begin{enumerate}
\item Let $a_{1}, \dotsc  , a_{n}$ be elements of $\mathcal{A}$ and $\sigma$ an $n$-ary operator: then 
\[\{\eta_{\mathcal{A}}(a)\ :\  a\in \sigma_{\mathcal{A}}(a_{1}, \dotsc  , a_{n})\}=\{\{a\}\ :\  a\in\sigma_{\mathcal{A}}(a_{1}, \dotsc  , a_{n})\}=\]
\[\sigma_{P\mathcal{A}}(\{a_{1}\}, \dotsc  , \{a_{n}\})=\sigma_{P\mathcal{A}}(\eta_{\mathcal{A}}(a_{1}), \dotsc  , \eta_{\mathcal{A}}(a_{n})).\]

\item Let $\{A_{i}^{1}\}_{i\in I_{1}}$ through $\{A_{i}^{n}\}_{i\in I_{n}}$ be elements of $PP\mathcal{A}$ and $\sigma$ be $n$-ary: we have that 
\[\bigcup_{(i_{1}, \dotsc  , i_{n})\in I_{1}\times\cdots\times I_{n}}\bigcup_{(a_{1}, \dotsc  , a_{n})\in A_{i_{1}}^{1}\times\cdots\times A_{i_{n}}^{n}}\sigma_{\mathcal{A}}(a_{1}, \dotsc  , a_{n})=\bigcup_{(a_{1}, \dotsc  , a_{n})\in \bigcup_{i_{1}\in I_{1}}A_{i_{1}}\times\cdots\times\bigcup_{i_{n}\in I_{n}}A_{i_{n}}}\sigma_{\mathcal{A}}(a_{1}, \dotsc  , a_{n})=\]
\[\bigcup_{(a_{1}, \dotsc  , a_{n})\in \epsilon_{\mathcal{A}}(\{A_{i}^{1}\}_{i\in I_{1}})\times\cdots\times\epsilon_{\mathcal{A}}(\{A_{i}^{n}\}_{i\in I_{n}})}\sigma_{\mathcal{A}}(a_{1}, \dotsc  , a_{n}),\]
implying that 
\[\sigma_{PP\mathcal{A}}(\{A_{i}^{1}\}_{i\in I_{1}}, \dotsc  , \{A_{i}^{n}\}_{i\in I_{n}})=\{\{A'\}\ :\  A'\in \bigcup_{(i_{1}, \dotsc  , i_{n})\in I_{1}\times\cdots\times I_{n}}\sigma_{P\mathcal{A}}(A_{i_{1}}^{1}, \dotsc  , A_{i_{n}}^{n})\}=\]
\[\{\{\{a\}\}\ :\  a\in \bigcup_{(i_{1}, \dotsc  , i_{n})\in I_{1}\times\cdots\times I_{n}}\bigcup_{(a_{1}, \dotsc  , a_{n})\in A_{i_{1}}^{1}\times\cdots\times A_{i_{n}}^{n}}\sigma_{\mathcal{A}}(a_{1}, \dotsc  , a_{n})\}=\]
\[\{\{\{a\}\}\ :\  a\in \bigcup_{(a_{1}, \dotsc  , a_{n})\in \epsilon_{\mathcal{A}}(\{A_{i}^{1}\}_{i\in I_{1}})\times\cdots\times\epsilon_{\mathcal{A}}(\{A_{i}^{n}\}_{i\in I_{n}})}\sigma_{\mathcal{A}}(a_{1}, \dotsc  , a_{n})\};\]
with this, we obtain that 
\[\{\epsilon_{\mathcal{A}}(\{A_{i}\}_{i\in I})\ :\  \{A_{i}\}_{i\in I}\in \sigma_{PP\mathcal{A}}(\{A_{i}^{1}\}_{i\in I_{1}}, \dotsc  , \{A_{i}^{n}\}_{i\in I_{n}})\}=\]
\[\{\epsilon_{\mathcal{A}}(\{\{a\}\})\ :\  a\in \bigcup_{(a_{1}, \dotsc  , a_{n})\in \epsilon_{\mathcal{A}}(\{A_{i}^{1}\}_{i\in I_{1}})\times\cdots\times\epsilon_{\mathcal{A}}(\{A_{i}^{n}\}_{i\in I_{n}})}\sigma_{\mathcal{A}}(a_{1}, \dotsc  , a_{n})\}=\]
\[\{\{a\}\ :\  a\in \bigcup_{(a_{1}, \dotsc  , a_{n})\in \epsilon_{\mathcal{A}}(\{A_{i}^{1}\}_{i\in I_{1}})\times\cdots\times\epsilon_{\mathcal{A}}(\{A_{i}^{n}\}_{i\in I_{n}})}\sigma_{\mathcal{A}}(a_{1}, \dotsc  , a_{n})\}=\]
\[\sigma_{P\mathcal{A}}(\epsilon_{\mathcal{A}}(\{A_{i}^{1}\}_{i\in I_{1}}), \dotsc  ,\epsilon_{\mathcal{A}}(\{A_{i}^{n}\}_{i\in I_{n}})).\]
\end{enumerate}
\end{proof}

\begin{proposition}
For any $\Sigma$-multialgebras $\mathcal{A}$ and $\mathcal{B}$, and homomorphism $\varphi:\mathcal{A}\rightarrow\mathcal{B}$, 
\begin{enumerate}
\item $P\varphi\circ \eta_{\mathcal{A}}=\eta_{\mathcal{B}}\circ \varphi$ and
\item $P\varphi\circ\epsilon_{\mathcal{A}}=\epsilon_{\mathcal{B}}\circ PP\varphi$,
\end{enumerate}
meaning $\eta$ and $\epsilon$ are natural transformations from, respectively, $Id_{\textbf{MAlg}(\Sigma)}$ to $P$, and $P\circ P$ to $P$.
\end{proposition}

\begin{proof}
\begin{enumerate}
\item Let $a$ be an element of $\mathcal{A}$: we have that $P\varphi\circ\eta_{\mathcal{A}}(a)=P\varphi(\eta_{\mathcal{A}}(a))$, and since $\eta_{\mathcal{A}}(a)=\{a\}$, we have this equals $\{\varphi(a)\}$; meanwhile, $\eta_{\mathcal{B}}\circ\varphi(a)=\eta_{\mathcal{A}}(\varphi(a))=\{\varphi(a)\}$, and as stated both expressions end up being equal.

\item Let $\{A_{i}\}_{i\in I}$ be an element of $PP\mathcal{A}$, meaning it is a non-empty set of non-empty subsets of $\mathcal{A}$: $P\varphi\circ\epsilon_{\mathcal{A}}(\{A_{i}\}_{i\in I})=P\varphi(\epsilon_{\mathcal{A}}(\{A_{i}\}_{i\in I}))$, and since $\epsilon_{\mathcal{A}}(\{A_{i}\}_{i\in I})=\bigcup_{i\in I}A_{i}$, the whole expression simplifies to $\{\varphi(a)\ :\  a\in \bigcup_{i\in I}A_{i}\}$; at the other hand, 
\[\epsilon_{\mathcal{B}}\circ PP\varphi(\{A_{i}\}_{i\in I})=\epsilon_{\mathcal{B}}(\{\{\varphi(a)\ :\  a\in A_{i}\}\ :\  i\in I\})=\bigcup_{i\in I}\{\varphi(a)\ :\  a\in A_{i}\}=\]
\[\{\varphi(a)\ :\  a\in \bigcup_{i\in I}A_{i}\},\]
giving us the desired equality.
\end{enumerate}
\end{proof}

\begin{theorem}
The triple formed by $P$, $\eta$ and $\epsilon$ constitutes a monad.
\end{theorem}

\begin{proof}
Let $\mathcal{A}$ be a $\Sigma$-multialgebra.
\begin{enumerate}
\item We must prove that $\epsilon\circ P\epsilon=\epsilon\circ\epsilon P$, what amounts to $\epsilon_{\mathcal{A}}\circ P\epsilon_{\mathcal{A}}=\epsilon_{\mathcal{A}}\circ \epsilon_{P\mathcal{A}}$, as homomorphisms from $P^{3}\mathcal{A}$ to $P\mathcal{A}$. So, let $\{\{A_{i}^{j}\}_{i\in I}\}_{j\in J}$ be an element of $P^{3}\mathcal{A}$, where $I$ and $J$ are non-empty sets of indexes and all $A_{i}^{j}$ are non-empty subsets of $A$: \[\epsilon_{\mathcal{A}}\circ P\epsilon_{\mathcal{A}}(\{\{A_{i}^{j}\}_{i\in I}\}_{j\in J})=\epsilon_{\mathcal{A}}(\{\epsilon_{\mathcal{A}}(\{A_{i}^{j}\ :\  i\in I\})\ :\  j\in J\})=\epsilon_{\mathcal{A}}(\{\bigcup_{i\in I}A_{i}^{j}\ :\  j\in J\})=\]
\[\bigcup_{j\in J}\bigcup_{i\in I}A_{i}^{j},\]
while
\[\epsilon_{\mathcal{A}}\circ\epsilon_{P\mathcal{A}}(\{\{A_{i}^{j}\}_{i\in I}\}_{j\in J})=\epsilon_{\mathcal{A}}(\bigcup_{j\in J}\{A_{i}^{j}\}_{i\in I})=\bigcup_{i\in I}\bigcup_{j\in J}A_{i}^{j},\]
and it is clear both sets are the same.

\item It remains to be proven that $\epsilon\circ P\eta=\epsilon\circ\eta P=1_{P}$, meaning $\epsilon_{\mathcal{A}}\circ \eta_{P\mathcal{A}}=\epsilon_{\mathcal{A}}\circ P\eta_{\mathcal{A}}$, as homomorphisms from $P\mathcal{A}$ to $P\mathcal{A}$, and this equals the identity homomorphism on this multialgebra as well. So, we take a non-empty subset $A'$ of $A$, and we have that
\[\epsilon_{\mathcal{A}}\circ\eta_{P\mathcal{A}}(A')=\epsilon_{\mathcal{A}}(\{A'\})=A',\]
while for the other other expression one derives
\[\epsilon_{\mathcal{A}}\circ P\eta_{\mathcal{A}}(A')=\epsilon_{\mathcal{A}}(\{\eta_{\mathcal{A}}(a)\ :\  a\in A'\})=\epsilon_{\mathcal{A}}(\{\{a\}\ :\  a\in A'\})=\bigcup_{a\in A'}\{a\}=A',\]
what finishes the proof.
\end{enumerate}
\end{proof}

\subsection{Multialgebras of atoms}

Given a $(\Sigma, \leq)$-algebra $\mathcal{A}$, take the set $\mathbb{A}((A, \leq_{\mathcal{A}}))$ of atoms of $(A,\leq_{\mathcal{A}})$, that is, the set of minimal elements of this poset.

For a $\sigma\in\Sigma_{n}$ and atoms $a_{1}, \dotsc  , a_{n}\in \mathbb{A}((A, \leq_{\mathcal{A}}))$, we define 
\[\sigma_{\mathbb{A}(\mathcal{A})}(a_{1}, \dotsc  , a_{n})=\{a\in \mathbb{A}((A, \leq_{\mathcal{A}}))\ :\  a\leq_{\mathcal{A}} \sigma_{\mathcal{A}}(a_{1}, \dotsc  , a_{n})\}=A_{\sigma_{\mathcal{A}}(a_{1}, \dotsc  , a_{n})}\neq\emptyset.\]

This way, $(\mathbb{A}((A, \leq_{\mathcal{A}})), \{\sigma_{\mathbb{A}(\mathcal{A})}\}_{\sigma\in\Sigma})$ becomes a $\Sigma$-multialgebra, that we will denote by $\mathbb{A}(\mathcal{A})$ and call the multialgebra of atoms\index{Multialgebra of atoms}\label{A} of $\mathcal{A}$.

Given $(\Sigma, \leq)$-algebras $\mathcal{A}$ and $\mathcal{B}$ and a $(\Sigma, \leq)$-homomorphism $\varphi:\mathcal{A}\rightarrow \mathcal{B}$, we define $\mathbb{A}(\varphi):\mathbb{A}((A, \leq_{\mathcal{A}}))\rightarrow\mathbb{A}((B, \leq_{\mathcal{B}}))$ as the restriction of $\varphi$ to $\mathbb{A}((A, \leq_{\mathcal{A}}))\subseteq A$: it is well-defined since every $(\Sigma, \leq)$-homomorphism preserves minimal elements, that is, atoms.

For $\sigma\in\Sigma_{n}$ and atoms $a_{1}, \dotsc  , a_{n}\in \mathbb{A}((A, \leq_{\mathcal{A}}))$ we have that 
\[\{\mathbb{A}(\varphi)(a)\ :\  a\in \sigma_{\mathbb{A}(\mathcal{A})}(a_{1}, \dotsc  , a_{n})\}=\{\varphi(a)\ :\  a\in\sigma_{\mathbb{A}(\mathcal{A})}(a_{1}, \dotsc  , a_{n})\}=\]
\[\{\varphi(a)\in \mathbb{A}((B, \leq_{\mathcal{B}}))\ :\  a\leq_{\mathcal{A}} \sigma_{\mathcal{A}}(a_{1}, \dotsc  , a_{n})\}\]
and, since $a\leq_{\mathcal{A}} \sigma_{\mathcal{A}}(a_{1}, \dotsc  , a_{n})$ implies $\varphi(a)\leq_{\mathcal{B}}\varphi(\sigma_{\mathcal{A}}(a_{1}, \dotsc  , a_{n}))$ given $\varphi$ is order preserving, which in turn implies $\varphi(a)\leq_{\mathcal{B}}\sigma_{\mathcal{B}}(\varphi(a_{1}), \dotsc  , \varphi(a_{n}))$ since $\varphi$ is an "almost homomorphism", we get that
\[\{\varphi(a)\in \mathbb{A}((B, \leq_{\mathcal{B}}))\ :\  a\leq_{\mathcal{A}} \sigma_{\mathcal{A}}(a_{1}, \dotsc  , a_{n})\}\subseteq\{b\in \mathbb{A}((B, \leq_{\mathcal{B}}))\ :\  b\leq_{\mathcal{B}}\sigma_{\mathcal{B}}(\varphi(a_{1}), \dotsc  , \varphi(a_{n})\}=\]
\[\sigma_{\mathbb{A}(\mathcal{B})}(\varphi(a_{1}), \dotsc  , \varphi(a_{n}))=\sigma_{\mathbb{A}(\mathcal{B})}(\mathbb{A}(\varphi)(a_{1}), \dotsc  , \mathbb{A}(\varphi)(a_{n})),\]
what proves $\mathbb{A}(\varphi)$ is a homomorphism between $\Sigma$-multialgebras, and we can write $\mathbb{A}(\varphi):\mathbb{A}(\mathcal{A})\rightarrow\mathbb{A}(\mathcal{B})$.

The natural question is whether $\mathbb{A}:\textbf{Alg}_{\mathcal{B}}(\Sigma)\rightarrow\textbf{MAlg}(\Sigma)$ is a functor, to which the answer is yes.

\begin{enumerate}
\item Let $\mathcal{A}$ be a $(\Sigma, \leq)$-algebra and $Id_{\mathcal{A}}$ be its identity $(\Sigma, \leq)$-homomorphism: we then have that, for every $a\in \mathbb{A}((A, \leq_{\mathcal{A}}))$, $\mathbb{A}(Id_{\mathcal{A}})(a)=Id_{\mathcal{A}}(a)=a$, and therefore $\mathbb{A}(I_{\mathcal{A}})$ is the identity $\Sigma$-homomorphism of $\mathbb{A}(\mathcal{A})$.
\item Let $\mathcal{A}$, $\mathcal{B}$ and $\mathcal{C}$ be $(\Sigma, \leq)$-algebras, and $\varphi:\mathcal{A}\rightarrow\mathcal{B}$ and $\psi:\mathcal{B}\rightarrow\mathcal{C}$ be $(\Sigma, \leq)$-homomorphisms: for every $a\in\mathbb{A}((A, \leq_{\mathcal{A}}))$ we have that
\[\mathbb{A}(\psi\circ\varphi)(a)=\psi\circ\varphi(a)=\psi(\varphi(a))=\mathbb{A}(\psi)(\varphi(a))=\mathbb{A}(\psi)(\mathbb{A}(\varphi)(a))=(\mathbb{A}(\psi)\circ\mathbb{A}(\varphi))(a).\]
\end{enumerate}

%

\section{$\textbf{Alg}_{\mathcal{B}}(\Sigma)$ and $\textbf{MAlg}(\Sigma)$ are equivalent}

Now, we aim to prove that $\textbf{Alg}_{\mathcal{B}}(\Sigma)$ and $\textbf{MAlg}(\Sigma)$ are actually equivalent categories, the equivalence being given by the functors $\mathcal{P}_{\mathcal{B}}$ and $\mathbb{A}$: so, let $Id_{\textbf{Alg}_{\mathcal{B}}(\Sigma)}$ and $Id_{\textbf{MAlg}(\Sigma)}$ be the identity functors on $\textbf{Alg}_{\mathcal{B}}(\Sigma)$ and $\textbf{MAlg}(\Sigma)$ respectively, and we should prove that $\mathcal{P}_{\mathcal{B}}\circ\mathbb{A}$ and $Id_{\textbf{Alg}_{\mathcal{B}}(\Sigma)}$, and $\mathbb{A}\circ\mathcal{P}_{\mathcal{B}}$ and $Id_{\textbf{MAlg}(\Sigma)}$ are naturally isomorphic.

Equivalently, $\mathcal{P}_{\mathcal{B}}$ and $\mathbb{A}$ form an equivalence of categories if both are full and faithful and $\mathbb{A}$ is a right adjoint of $\mathcal{P}_{\mathcal{B}}$, what we will prove instead.

%

\subsection{$\mathcal{P}_{\mathcal{B}}$ and $\mathbb{A}$ are full and faithful}

It is easy to see $\mathcal{P}_{\mathcal{B}}$ is faithful: given $\Sigma$-multialgebras $\mathcal{A}$ and $\mathcal{B}$ and homomorphisms $\varphi, \psi:\mathcal{A}\rightarrow\mathcal{B}$, if $\mathcal{P}_{\mathcal{B}}(\varphi)=\mathcal{P}_{\mathcal{B}}(\psi)$, we have that for every $a\in A$
\[\{\varphi(a)\}=\mathcal{P}_{\mathcal{B}}(\varphi)(\{a\})=\mathcal{P}_{\mathcal{B}}(\psi)(\{a\})=\{\psi(a)\},\]
and therefore $\varphi=\psi$.

\begin{proposition}
$\mathbb{A}$ is faithful.
\end{proposition}

\begin{proof}
Given $(\Sigma, \leq)$-algebras $\mathcal{A}$ and $\mathcal{B}$, and $(\Sigma, \leq)$-homomorphisms $\varphi, \psi:\mathcal{A}\rightarrow\mathcal{B}$, suppose we have $\mathbb{A}(\varphi)=\mathbb{A}(\psi)$: then, for every $a\in A$, we can write $a=\sup A_{a}$, since $(A, \leq_{\mathcal{A}})$ is atomic.

Since $\varphi$ and $\psi$ are continuous, $\varphi(a)=\sup\{\varphi(a')\ :\  a'\in A_{a}\}$ and $\psi(a)=\sup\{\psi(a')\ :\  a'\in A_{a}\}$; but, since $\mathbb{A}(\varphi)=\mathbb{A}(\psi)$, $\varphi$ and $\psi$ are the same when restricted to atoms, and therefore $\{\varphi(a')\ :\  a'\in A_{a}\}=\sup\{\psi(a')\ :\  a'\in A_{a}\}$, what means $\varphi(a)=\psi(a)$ and, since $a$ is arbitrary, $\varphi=\psi$.
\end{proof}

Now, given $\Sigma$-multialgebras $\mathcal{A}$ and $\mathcal{B}$ and a $(\Sigma, \leq)$-homomorphism $\varphi:\mathcal{P}_{\mathcal{B}}(\mathcal{A})\rightarrow\mathcal{P}_{\mathcal{B}}(\mathcal{B})$, to prove $\mathcal{P}_{\mathcal{B}}$ is also full we must find a homomorphism $\psi:\mathcal{A}\rightarrow\mathcal{B}$ such that $\mathcal{P}_{\mathcal{B}}(\psi)=\varphi$.

For every $a\in A$, $\{a\}$ is an atom and, since $\varphi$ preserves atoms, $\varphi(\{a\})$ is an atom of $\mathcal{P}_{\mathcal{B}}(B)$, and therefore of the form $\{b_{a}\}$ for some $b_{a}\in B$: we define $\psi:\mathcal{A}\rightarrow\mathcal{B}$ by $\psi(a)=b_{a}$. First of all, we must show $\psi$ is in fact a homomorphism, which is quite analogous to the proof of the same fact for $\mathbb{A}(\varphi)$: given $\sigma\in\Sigma_{n}$ and $a_{1}, \dotsc  , a_{n}\in A$,
\[\{\psi(a)\ :\  a\in \sigma_{\mathcal{A}}(a_{1}, \dotsc  , a_{n})\}=\{b_{a}\ :\  a\in\sigma_{\mathcal{A}}(a_{1}, \dotsc  , a_{n})\}=\sup\{\{b_{a}\}\ :\  a\in\sigma_{\mathcal{A}}(a_{1}, \dotsc  , a_{n})\}=\]
\[\sup\{\varphi(\{a\})\ :\  a\in\sigma_{\mathcal{A}}(a_{1}, \dotsc  , a_{n})\}=\varphi(\sup \{\{a\}\ :\  a\in\sigma_{\mathcal{A}}(a_{1}, \dotsc  , a_{n})\})=\varphi(\sigma_{\mathcal{A}}(a_{1}, \dotsc  , a_{n}))=\]
\[\varphi(\sigma_{\mathcal{P}_{\mathcal{B}}(\mathcal{A})}(\{a_{1}\}, \dotsc  , \{a_{n}\}))\subseteq \sigma_{\mathcal{P}_{\mathcal{B}}(\mathcal{B})}(\varphi(\{a_{1}\}), \dotsc  , \varphi(\{a_{n}\}))= \sigma_{\mathcal{P}_{\mathcal{B}}(\mathcal{B})}(\{b_{a_{1}}\}, \dotsc  , \{b_{a_{n}}\})=\]
\[\sigma_{\mathcal{B}}(b_{a_{1}}, \dotsc  , b_{a_{n}})=\sigma_{\mathcal{B}}(\psi(a_{1}), \dotsc  , \psi(a_{n})).\]
Now, when we consider $\mathcal{P}_{\mathcal{B}}(\psi)$, we see that, for every atom $\{a\}$ of $\mathcal{P}_{\mathcal{B}}(\mathcal{A})$, $\mathcal{P}_{\mathcal{B}}(\psi)(\{a\})=\{b_{a}\}=\varphi(\{a\})$, and so the restrictions of $\varphi$ and $\mathcal{P}_{\mathcal{B}}(\psi)$ to atoms are the same, therefore $\mathbb{A}(\varphi)=\mathbb{A}(\mathcal{P}_{\mathcal{B}}(\psi))$. Since $\mathbb{A}$ is faithful, we discover that $\varphi=\mathcal{P}_{\mathcal{B}}(\psi)$ and, as we stated before, $\mathcal{P}_{\mathcal{B}}$ is full.

Now it remains to be shown that $\mathbb{A}$ is also full: given $(\Sigma, \leq)$-algebras $\mathcal{A}$ and $\mathcal{B}$ and a homomorphism $\varphi:\mathbb{A}(\mathcal{A})\rightarrow\mathbb{A}(\mathcal{B})$; we then define $\psi:\mathcal{A}\rightarrow\mathcal{B}$ by
\[\psi(a)=\sup\{\varphi(c)\ :\  c\in A_{a}\}.\]
First of all, we must prove that $\psi$ is a $(\Sigma, \leq)$-homomorphism, for which we shall need a few lemmas.

\begin{lemma}\label{sup. of the union is sup. of sup.}
Given a family of subsets $\{X_{i}\}_{i\in I}$ of a complete poset $(S, \leq)$, suppose we define $x_{i}=\sup X_{i}$ for $i\in I$ and $X=\bigcup_{i\in I}X_{i}$: then $\sup\{x_{i}\ :\  i\in I\}=\sup X$.
\end{lemma}

\begin{proof}
We define $a=\sup\{x_{i}\ :\  i\in I\}$ and $b=\sup X$: first, we show that $a$ is an upper bound for $X$, so that $a\geq b$. For every $x\in X$, we have that, since $X=\bigcup_{i\in I}X_{i}$, there exists $j\in I$ such that $x\in X_{j}$, and therefore $x_{j}\geq x$; since $a=\sup\{x_{i}\ :\  i\in I\}$, we have that $a\geq x_{j}$, and by transitivity $a\geq x$, and therefore $a$ is indeed an upper bound for $X$.

Now we show that $b$ is an upper bound for $\{x_{i}\ :\  i\in I\}$, and so $b\geq a$ (and $a=b$). For every $i\in I$, we have that $b$ is an upper bound for $X_{i}$  (since $X_{i}\subseteq X$ and $b$ is an upper bound for $X$), and therefore $b\geq x_{i}$, since $x_{i}$ is the smallest upper bound for $X_{i}$; it follows that $b$ is indeed an upper bound for $\{x_{i}\ :\  i\in I\}$, what finishes the proof.
\end{proof}

\begin{lemma}\label{Union of atoms is atoms of sup.}
In a complete, atomic and bottomless Boolean algebra, given a non-empty set $C$, $\bigcup_{c\in C}A_{c}=A_{\sup C}$.
\end{lemma}

\begin{proof}
If $d\in A_{c}$ for a $c\in C$, $d$ is an atom smaller than $c$; since $c\leq_{\mathcal{A}}\sup C$, $d\leq_{\mathcal{A}}\sup C$, and therefore belongs to $A_{\sup C}$. So $\bigcup_{c\in C}A_{c}\subseteq A_{\sup C}$.

Reciprocally, suppose $d\in A_{\sup C}$: then $d$ is an atom which is smaller than $\sup C$, and therefore $\inf\{d, \sup C\}=d$; it follows the subset $C^{d}\subseteq C$ of $c\in C$ such that $\inf\{d,c\}$ exists is not empty, by the infinite-distributivity of $\mathcal{A}$. But if $c\in C^{d}$, $\inf\{d,c\}$ exists, and since $d$ is an atom, we have that $d\in A_{c}\subseteq \bigcup_{c\in C}A_{c}$, and from that $\bigcup_{c\in C}A_{c}=A_{\sup C}$.
\end{proof}

Since $\sigma_{\mathcal{A}}(a_{1}, \dotsc  , a_{n})$ equals the supremum of $\{\sigma_{\mathcal{A}}(c_{1}, \dotsc  , c_{n})\ :\  (c_{1}, \dotsc  , c_{n})\in A_{a_{1}}\times\cdots\times A_{a_{n}}\}$, by Lemma \ref{Union of atoms is atoms of sup.} we have that $A_{\sigma_{\mathcal{A}}(a_{1}, \dotsc  , a_{n})}$ equals 
\[\bigcup_{(c_{1}, \dotsc   ,c_{n})\in A_{a_{1}}\times\cdots\times A_{a_{n}}}A_{\sigma_{\mathcal{A}}(c_{1}, \dotsc  , c_{n})},\]
that is, we have the following lemma.

\begin{lemma}\label{atoms of operation are union of operations over atoms}
For a $\sigma\in\Sigma_{n}$, and $a_{1}, \dotsc  , a_{n}\in A$,
\[A_{\sigma_{\mathcal{A}}(a_{1}, \dotsc  , a_{n})}=\bigcup_{(c_{1}, \dotsc  , c_{n})\in A_{a_{1}}\times\cdots\times A_{a_{n}}}A_{\sigma_{\mathcal{A}}(c_{1}, \dotsc  , c_{n})}.\]
\end{lemma}

\begin{theorem}
$\mathbb{A}$ is full.
\end{theorem}

\begin{proof}
So, first of all, we prove using Lemmas \ref{sup. of the union is sup. of sup.}, \ref{Union of atoms is atoms of sup.} and \ref{atoms of operation are union of operations over atoms} that $\psi$ is a $(\Sigma, \leq)$-homomor\-phism.
\begin{enumerate}
\item First, it is clear that $\psi$ maps atoms into atoms: if $a$ is an atom,
\[\psi(a)=\sup\{\varphi(c)\ :\  c\in A_{a}\}=\sup\{\varphi(a)\}=\varphi(a),\]
which is an atom since $\varphi$ is a map between $\mathbb{A}(\mathcal{A})$ and $\mathbb{A}(\mathcal{B})$.

\item $\psi$ is continuous: for any non-empty set $C\subseteq A$, we remember that $\psi(c)=\sup\{\varphi(d)\ :\  d\in A_{c}\}$, and from Lemma \ref{sup. of the union is sup. of sup.} we get that
\[\sup\{\psi(c)\ :\  c\in C\}=\sup\{\sup\{\varphi(d)\ :\  d\in A_{c}\}\ :\  c\in C\}=\sup \bigcup_{c\in C}\{\varphi(d)\ :\  d\in A_{c}\};\]
from Lemma \ref{Union of atoms is atoms of sup.},
\[\sup \bigcup_{c\in C}\{\varphi(d)\ :\  d\in A_{c}\}=\sup\{\varphi(d)\ :\  d\in A_{\sup C}\}=\psi(\sup C).\]

\item 
Since $\{\varphi(a)\ :\  a\in\sigma_{\mathbb{A}(\mathcal{A})}(a_{1}, \dotsc  , a_{n})\}\subseteq \sigma_{\mathbb{A}(\mathcal{B})}(\varphi(a_{1}), \dotsc  , \varphi(a_{n}))$, it follows from Lemma \ref{atoms of operation are union of operations over atoms} that
\[\psi(\sigma_{\mathcal{A}}(a_{1}, \dotsc  , a_{n}))=\sup\{\varphi(c)\ :\  c\in A_{\sigma_{\mathcal{A}}(a_{1}, \dotsc  , a_{n})}\}=\]
\[\sup\{\varphi(c)\ :\  c\in \bigcup_{(c_{1}, \dotsc  , c_{n})\in A_{a_{1}}\times\cdots\times A_{a_{n}}}A_{\sigma_{\mathcal{A}}(c_{1}, \dotsc  , c_{n})}\};\]
since, for atoms $c_{1}, \dotsc  , c_{n}$ of $\mathcal{A}$, $\sigma_{\mathbb{A}(\mathcal{A})}(c_{1}, \dotsc  , c_{n})=A_{\sigma_{\mathcal{A}}(c_{1}, \dotsc  , c_{n})}$, we have that the previous expression equals
\[\sup\{\varphi(c)\ :\  c\in \bigcup_{(c_{1}, \dotsc  , c_{n})\in A_{a_{1}}\times\cdots\times A_{a_{n}}}\sigma_{\mathbb{A}(\mathcal{A})}(c_{1}, \dotsc  , c_{n})\}=\]
\[\sup\bigcup_{(c_{1}, \dotsc  , c_{n})\in A_{a_{1}}\times\cdots\times A_{a_{n}}}\{\varphi(c)\ :\  c\in \sigma_{\mathbb{A}(\mathcal{A})}(c_{1}, \dotsc  , c_{n})\},\]
which is less than or equal to 
\[\sup \bigcup_{(c_{1}, \dotsc  , c_{n})\in A_{a_{1}}\times\cdots\times A_{a_{n}}}\sigma_{\mathbb{A}(\mathcal{B})}(\varphi(c_{1}), \dotsc  , \varphi(c_{n}))\]
in $\mathcal{B}$; since, for atoms $d_{1}, \dotsc  , d_{n}$ of $\mathcal{B}$, we have that $\sigma_{\mathbb{A}(\mathcal{B})}(d_{1}, \dotsc  , d_{n})=A_{\sigma_{\mathcal{B}}(d_{1}, \dotsc  , d_{n})}$, this equals
\[\sup \bigcup_{(c_{1}, \dotsc  , c_{n})\in A_{a_{1}}\times\cdots\times A_{a_{n}}}A_{\sigma_{\mathcal{B}}(\varphi(c_{1}), \dotsc  , \varphi(c_{n}))}=\sup \bigcup_{(c_{1}, \dotsc  , c_{n})\in A_{a_{1}}\times\cdots\times A_{a_{n}}}A_{\sigma_{\mathcal{B}}(\psi(c_{1}), \dotsc  , \psi(c_{n}))};\]
now, since $\psi$ is continuous, $c_{i}\leq_{\mathcal{A}}a_{i}$, for every $i\in \{1, \dotsc  , n\}$, implies $\psi(c_{i})\leq_{\mathcal{B}}\psi(a_{i})$, and therefore $\sigma_{\mathcal{B}}(\psi(c_{1}), \dotsc  , \psi(c_{n}))\leq_{\mathcal{B}}\sigma_{\mathcal{B}}(\psi(a_{1}), \dotsc  , \psi(a_{n}))$ for $(c_{1}, \dotsc  , c_{n})\in A_{a_{1}}\times\cdots\times A_{a_{n}}$.

It follows that $\bigcup_{(c_{1}, \dotsc  , c_{n})\in A_{a_{1}}\times\cdots\times A_{a_{n}}}A_{\sigma_{\mathcal{B}}(\psi(c_{1}), \dotsc  , \psi(c_{n}))}$ is contained on $A_{\sigma_{\mathcal{B}}(\psi(a_{1}), \dotsc  , \psi(a_{n}))}$, and therefore 
\[\sup \bigcup_{(c_{1}, \dotsc  , c_{n})\in A_{a_{1}}\times\cdots\times A_{a_{n}}}A_{\sigma_{\mathcal{B}}(\psi(c_{1}), \dotsc  , \psi(c_{n}))}\leq_{\mathcal{B}}\sup A_{\sigma_{\mathcal{B}}(\psi(a_{1}), \dotsc  , \psi(a_{n}))}=\sigma_{\mathcal{B}}(\psi(a_{1}), \dotsc  , \psi(a_{n})).\]
\end{enumerate}

Now, for every atom $a$ of $\mathcal{A}$, we have that $\psi(a)=\varphi(a)$, and therefore the restriction of $\psi$ to atoms equals $\varphi$, that is, $\mathbb{A}(\psi)=\varphi$, and since $\varphi$ was taken arbitrarily, $\mathbb{A}$ is full.
\end{proof}

%

\subsection{$\mathcal{P}_{\mathcal{B}}$ and $\mathbb{A}$ are adjoint}

It remains to be shown that $\mathcal{P}_{\mathcal{B}}$ and $\mathbb{A}$ are adjoint: consider the bijections
\[\Phi_{\mathcal{B}, \mathcal{A}}:Hom_{\textbf{MAlg}(\Sigma)}(\mathbb{A}(\mathcal{B}), \mathcal{A})\rightarrow Hom_{\textbf{Alg}_{\mathcal{B}}(\Sigma)}(\mathcal{B}, \mathcal{P}_{\mathcal{B}}(\mathcal{A})),\]
for $\mathcal{A}$ a $\Sigma$-multialgebra and $\mathcal{B}$ a $(\Sigma, \leq)$-algebra, given by, for $\varphi:\mathbb{A}(\mathcal{B})\rightarrow\mathcal{A}$ a homomorphism and $b$ an element of $\mathcal{B}$,
\[\Phi_{\mathcal{B},\mathcal{A}}(\varphi)(b)=\{\varphi(c)\ :\  c\in A_{b}\}.\]
First of all, we must prove $\Phi_{\mathcal{B}, \mathcal{A}}(\varphi)$ is truly a $(\Sigma, \leq)$-homomorphism.

\begin{enumerate}
\item If $b$ is an atom, $A_{b}=\{b\}$, and therefore 
\[\Phi_{\mathcal{B}, \mathcal{A}}(\varphi)(b)=\{\varphi(c)\ :\  c\in A_{b}\}=\{\varphi(b)\},\]
which is a singleton and therefore an atom of $\mathcal{P}_{\mathcal{B}}(\mathcal{A})$.

\item 
Let $B'$ be a non-empty subset of $\mathcal{B}$, and we have that
\[\Phi_{\mathcal{B}, \mathcal{A}}(\varphi)(\sup B')=\{\varphi(c)\ :\  c\in A_{\sup B'}\}=\{\varphi(c)\ :\  c\in \bigcup_{b\in B'}A_{b}\}=\]
\[\bigcup_{b\in B'}\{\varphi(c)\ :\  c\in A_{b}\}=\bigcup_{b\in B'}\Phi_{\mathcal{B}, \mathcal{A}}(\varphi)(b)=\sup\{\Phi_{\mathcal{B}, \mathcal{A}}(\varphi)(b)\ :\  b\in B'\},\]
since $A_{\sup B'}=\bigcup_{b\in B'}A_{d}$, from Lemma \ref{Union of atoms is atoms of sup.}, and the supremum in $\mathcal{P}_{\mathcal{B}}(\mathcal{A})$ is simply the union.

\item For $\sigma\in\Sigma_{n}$ and $b_{1}, \dotsc  , b_{n}$ elements of $\mathcal{B}$, we have that, using Lemma \ref{atoms of operation are union of operations over atoms},
\[\Phi_{\mathcal{B}, \mathcal{A}}(\varphi)(\sigma_{\mathcal{B}}(b_{1}, \dotsc  , b_{n}))=\{\varphi(c)\ :\  c\in A_{\sigma_{\mathcal{B}}(b_{1}, \dotsc  , b_{n})}\}=\]
\[\{\varphi(c)\ :\  c\in \bigcup_{(c_{1}, \dotsc  , c_{n})\in A_{b_{1}}\times\cdots\times A_{b_{n}}}A_{\sigma_{\mathcal{B}}(c_{1}, \dotsc  , c_{n})}\}=\bigcup_{(c_{1}, \dotsc  , c_{n})\in A_{b_{1}}\times\cdots\times A_{b_{n}}}\{\varphi(c)\ :\  c\in A_{\sigma_{\mathcal{B}}(c_{1}, \dotsc  , c_{n})}\},\]
and, since $c_{1}, \dotsc  , c_{n}$ are atoms and $\varphi$ is a homomorphism, this equals
\[\bigcup_{(c_{1}, \dotsc  , c_{n})\in A_{b_{1}}\times\cdots\times A_{b_{n}}}\{\varphi(c)\ :\  c\in \sigma_{\mathbb{A}(\mathcal{B})}(c_{1}, \dotsc  , c_{n})\}\subseteq\]
\[\bigcup_{(c_{1}, \dotsc  , c_{n})\in A_{b_{1}}\times\cdots\times A_{b_{n}}}\sigma_{\mathcal{A}}(\varphi(c_{1}), \dotsc  , \varphi(c_{n}))=\bigcup_{(a_{1}, \dotsc  , a_{n})\in \{\varphi(c)\ :\  c\in A_{b_{1}}\}\times\cdots\times \{\varphi(c)\ :\  c\in A_{b_{n}}\}}\sigma_{\mathcal{A}}(a_{1}, \dotsc  , a_{n})=\]
\[\bigcup_{(a_{1}, \dotsc  , a_{n})\in \Phi_{\mathcal{B}, \mathcal{A}}(\varphi)(b_{1})\times\cdots\times \Phi_{\mathcal{B}, \mathcal{A}}(\varphi)(b_{n})}\sigma_{\mathcal{A}}(a_{1}, \dotsc  , a_{n})=\sigma_{\mathcal{P}_{\mathcal{B}}(\mathcal{A})}(\Phi_{\mathcal{B}, \mathcal{A}}(\varphi)(b_{1}), \dotsc  , \Phi_{\mathcal{B}, \mathcal{A}}(\varphi)(b_{n})).\]
\end{enumerate}

Second, the $\Phi_{\mathcal{B}, \mathcal{A}}$ must be bijections between $Hom_{\textbf{MAlg}(\Sigma)}(\mathbb{A}(\mathcal{B}), \mathcal{A})$ and\\ $Hom_{\textbf{Alg}_{\mathcal{B}}(\Sigma)}(\mathcal{B}, \mathcal{P}_{\mathcal{B}}(\mathcal{A}))$. They are certainly injective: if $\Phi_{\mathcal{B}, \mathcal{A}}(\varphi)=\Phi_{\mathcal{B}, \mathcal{A}}(\psi)$, for every atom $b$ we have that
\[\{\varphi(b)\}=\Phi_{\mathcal{B}, \mathcal{A}}(\varphi)(b)=\Phi_{\mathcal{B}, \mathcal{A}}(\psi)(b)=\{\psi(b)\},\]
and therefore $\varphi=\psi$.

For the surjectivity, take a $(\Sigma, \leq)$-homomorphism $\varphi:\mathcal{B}\rightarrow \mathcal{P}_{\mathcal{B}}(\mathcal{A})$: we then define $\psi:\mathbb{A}(\mathcal{B})\rightarrow \mathcal{A}$ by $\psi(b)=a$, for an atom $b$ in $\mathcal{B}$ with $\varphi(b)=\{a\}$. It is well-defined since a $(\Sigma, \leq)$-homomorphism takes atoms to atoms, and the atoms of $\mathcal{P}_{\mathcal{B}}(\mathcal{A})$ are exactly the singletons.

We must show that $\psi$ is truly a homomorphism: for $\sigma\in\Sigma_{n}$ and atoms $b_{1}, \dotsc  , b_{n}$ in $\mathbb{A}(\mathcal{B})$ such that $\varphi(b_{i})=\{a_{i}\}$ for every $i\in\{1, \dotsc  , n\}$, we have that
\[\varphi(\sigma_{\mathcal{B}}(b_{1}, \dotsc  , b_{n}))\subseteq \sigma_{\mathcal{P}_{\mathcal{B}}(\mathcal{A})}(\varphi(b_{1}), \dotsc  , \varphi(b_{n})),\]
since $\varphi$ is a $(\Sigma, \leq)$-homomorphism, and therefore
\[\{\psi(b)\ :\  b\in\sigma_{\mathbb{A}(\mathcal{B})}(b_{1}, \dotsc  , b_{n})\}=\{\psi(b)\ :\  b\in A_{\sigma_{\mathcal{B}}(b_{1}, \dotsc  , b_{n})}\}=\]
\[\bigcup_{b\in A_{\sigma_{\mathcal{B}}(b_{1}, \dotsc  , b_{n})}}\varphi(b)=\sup\{\varphi(b)\ :\ b\in A_{\sigma_{\mathcal{B}}(b_{1}, \dotsc  , b_{n})}\}=\varphi(\sup A_{\sigma_{\mathcal{B}}(b_{1}, \dotsc  , b_{n})})=\varphi(\sigma_{\mathcal{B}}(b_{1}, \dotsc  , b_{n}))\subseteq\]
\[ \sigma_{\mathcal{P}_{\mathcal{B}}(\mathcal{A})}(\varphi(b_{1}), \dotsc  , \varphi(b_{n}))=\sigma_{\mathcal{P}_{\mathcal{B}}(\mathcal{A})}(\{a_{1}\}, \dotsc  , \{a_{n}\})=\sigma_{\mathcal{A}}(a_{1}, \dotsc  , a_{n})=\sigma_{\mathcal{A}}(\psi(b_{1}), \dotsc  , \psi(b_{n})).\]

Now $\Phi_{\mathcal{B}, \mathcal{A}}(\psi)=\varphi$ since, for any element $b$ in $\mathcal{B}$, we have that
\[\Phi_{\mathcal{B}, \mathcal{A}}(\psi)(b)=\{\psi(c)\ :\  c\in A_{b}\}=\bigcup_{c\in A_{b}}\varphi(c)=\sup\{\varphi(c)\ :\ c\in A_{b}\}=\varphi(\sup A_{b})=\varphi(b),\]
and therefore the $\Phi_{\mathcal{B}, \mathcal{A}}$ are, indeed, bijective.

Finally, for $\mathcal{A}$ and $\mathcal{C}$ $\Sigma$-multialgebras, $\mathcal{B}$ and $\mathcal{D}$ $(\Sigma, \leq)$-algebras, $\varphi:\mathcal{A}\rightarrow\mathcal{C}$ a homomorphism and $\psi:\mathcal{D}\rightarrow\mathcal{B}$ a $(\Sigma, \leq)$-homomorphism, we must now only prove that the following diagram commutes.

\[
\begin{tikzcd}[row sep=5em, column sep=3em]
Hom_{\textbf{MAlg}(\Sigma)}(\mathbb{A}(\mathcal{B}), \mathcal{A}) \arrow{r}{\Phi_{\mathcal{B}, \mathcal{A}}} \arrow{d}{Hom(\mathbb{A}(\psi), \varphi)}& Hom_{\textbf{Alg}_{\mathcal{B}}(\Sigma)}(\mathcal{B}, \mathcal{P}_{\mathcal{B}}(\mathcal{A})) \arrow{d}{Hom(\psi, \mathcal{P}_{\mathcal{B}}(\varphi))} \\%
Hom_{\textbf{MAlg}(\Sigma)}(\mathbb{A}(\mathcal{D}), \mathcal{C}) \arrow{r}{\Phi_{\mathcal{D}, \mathcal{C}}}& Hom_{\textbf{Alg}_{\mathcal{B}}(\Sigma)}(\mathcal{D}, \mathcal{P}_{\mathcal{B}}(\mathcal{C}))
\end{tikzcd}\]

So, we take a homomorphism $\theta:\mathbb{A}(\mathcal{B})\rightarrow\mathcal{A}$ and an element $d$ of $\mathcal{D}$. We have that
\[Hom(\psi, \mathcal{P}_{\mathcal{B}}(\varphi))\Phi_{\mathcal{B}, \mathcal{A}}(\theta)=\mathcal{P}_{\mathcal{B}}(\varphi)\circ\Phi_{\mathcal{B}, \mathcal{A}}(\theta)\circ\psi,\]
and therefore the top-right edges of the diagram give us 
\[\mathcal{P}_{\mathcal{B}}(\varphi)\circ\Phi_{\mathcal{B}, \mathcal{A}}(\theta)\circ\psi(d)=\mathcal{P}_{\mathcal{B}}(\varphi)(\{\theta(b)\ :\  b\in A_{\psi(d)}\})=\{\varphi\circ\theta(b)\ :\  b\in A_{\psi(d)}\}.\]

Meanwhile, the left-bottom edges give us 
\[\Phi_{\mathcal{D}, \mathcal{C}}(\varphi\circ \theta\circ \mathbb{A}(\psi))(d)=\{\varphi\circ \theta\circ \mathbb{A}(\psi)(e)\ :\  e\in A_{d}\}=\{\varphi\circ\theta\circ\psi(e)\ :\  e\in A_{d}\}.\]

If $d$ is an atom, the top-right edges give the singleton containing only $\varphi\circ\theta\circ\psi(d)$, since in this case $A_{d}=\{d\}$ and, given that $\psi$ preserves atoms, $A_{\psi(d)}=\{\psi(d)\}$; the left-bottom edges also give the singleton formed by $\varphi\circ\theta\circ\psi(d)$, since again $A_{d}=\{d\}$. Since a $(\Sigma, \leq)$-homomorphism is determined by its action on atoms, we find that the left and right sides of the diagram are equal, and therefore the diagram commutes.

As observed before, this proves $\mathbb{A}$ and $\mathcal{P}_{\mathcal{B}}$ are adjoint and, therefore, that $\textbf{MAlg}(\Sigma)$ and $\textbf{Alg}_{\mathcal{B}}(\Sigma)$ are equivalent. We then have, more or less, the following functors and categories to consider.

\[ \begin{tikzcd}
\textbf{MAlg}_{=}(\Sigma)\arrow[hook]{rr}\arrow{dd}{\mathcal{P}_{=}} && \textbf{MAlg}(\Sigma)\arrow[looseness=3, out=30, in=330, distance=4em, "P"]\arrow{ddll}{\mathcal{P}}\arrow[transform canvas={xshift=.7ex}]{dd}{\mathcal{P}_{\mathcal{B}}}\\
&&\\
\textbf{Alg}(\Sigma) && \textbf{Alg}_{\mathcal{B}}(\Sigma)\arrow[transform canvas={xshift=-.7ex}]{uu}{\mathbb{A}}\arrow{ll}{\mathbb{U}}
\end{tikzcd}
\]

While $\mathcal{P}$ and $\mathcal{P}_{=}$ are our failed attempts, $\mathcal{P}_{\mathcal{B}}$ is the successful one, and $P$ its presentation as an endofunctor of $\textbf{MAlg}(\Sigma)$ and part of a monad. Here, $\mathbb{U}$ is the forgetful functor from $\textbf{Alg}_{\mathcal{B}}(\Sigma)$ into $\textbf{Alg}(\Sigma)$, that ignores the order of a $(\Sigma, \leq)$-algebra, leaving us simply with a $\Sigma$-algebra.

\section{Some related results}

The result that $\textbf{MAlg}(\Sigma)$ and $\textbf{Alg}_{\mathcal{B}}(\Sigma)$ are equivalent has a few consequences, and related results, we would like to stress. We start by mentioning a consequence, and move to a few related results to which, although we do not offer a proof, providing a demonstration appears straightforward.

So, take the empty signature, with no operators at all: in that case, given all multialgebras have non-empty  universes, $\textbf{MAlg}(\Sigma)$ becomes (or, rather, is equivalent to) the category of non-empty sets $\textbf{Set}^{*}$\label{Set*} (a multialgebra corresponding to its universe), with functions between them as morphisms; this is because, once one disregards the conditions demanded of a homomorphism, what remains is purely a function between the universes of the multialgebras. 

Meanwhile, $\textbf{Alg}_{\mathcal{B}}(\Sigma)$ becomes the category with complete, atomic and bottomless Boolean algebras as objects (given we simply drop the operations from a $(\Sigma, \leq)$-algebra), with continuous, atoms-preserving functions between them as morphisms, category we will call $\textbf{CABo}$\label{CABo} (standing for ``complete, atomic and bottomless''). More explicitly, the equivalence is given by the pair of functors: 
\[F:\textbf{Set}^{*}\rightarrow \textbf{CABo},\]
which takes a non-empty set $X$ to the complete, atomic and bottomless Boolean algebra $\mathcal{P}^{*}(X)$ of its non-empty sets, and a function $f:X\rightarrow Y$ to the continuous, atoms-preserving function $Ff:\mathcal{P}^{*}(X)\rightarrow\mathcal{P}^{*}(Y)$ such that, for a non-empty $A\subseteq X$, $Ff(A)=\{f(x)\in Y : x\in A\}$; and
\[G:\textbf{CABo}\rightarrow\textbf{Set}^{*},\]
which takes a complete, atomic and bottomless Boolean algebra to its set of minimal elements (atoms), and a continuous, atoms-preserving function $\varphi:(A, \leq_{\mathcal{A}})\rightarrow(B, \leq_{\mathcal{B}})$ to the function $G\varphi$ which maps an atom $a$ in $A$ to the, still minimal, element $\varphi(a)$ of $B$.

Notice this is very closely related to the equivalence between $\textbf{CABA}$ and $\textbf{Set}^{op}$: the morphisms on the former are merely continuous functions, so the only extra requirement to the morphisms of $\textbf{CABo}$ we are making is that they should preserve atoms; this, of course, allows one to exchange the opposite category of $\textbf{Set}$ by $\textbf{Set}$ itself (or rather $\textbf{Set}^{*}$); the extra requirement that the maps should preserve atoms also allows one to prove a similar equivalence between an enriched $\textbf{CABA}$ and $\textbf{Set}$.

\subsection{Partial multialgebras}

Now, a generalization of our result is to partial multialgebras\index{Partial multialgebras}: given a signature $\Sigma$, a partial $\Sigma$-multialgebra is a pair $\mathcal{A}=(A, \{\sigma_{\mathcal{A}}\}_{\sigma\in\Sigma})$ (with $A$ possibly empty) such that, if $\sigma\in\Sigma_{n}$, $\sigma_{\mathcal{A}}$ is a function of the form
\[\sigma_{\mathcal{A}}:A^{n}\rightarrow\mathcal{P}(A)\]
(no longer $\mathcal{P}(A)\setminus\{\emptyset\}$ as in the case of multialgebras); that is, a partial multialgebra is a multialgebra where operations may return the empty-set. Given partial $\Sigma$-multialgebras $\mathcal{A}$ and $\mathcal{B}$, a homomorphism between them is a function $\varphi:A\rightarrow B$ such that, as is the case for homomorphisms for multialgebras,
\[\{\varphi(a): a\in\sigma_{\mathcal{A}}(a_{1}, \dotsc  , a_{n})\}\subseteq \sigma_{\mathcal{B}}(\varphi(a_{1}), \dotsc  , \varphi(a_{n})),\]
for all $\sigma\in\Sigma_{n}$ and $a_{1}, \dotsc  , a_{n}\in A$, where one or both of these sets may be empty. The class of all partial $\Sigma$-multialgebras, with these homomorphisms between them as morphisms, becomes a category, which we shall denote by $\textbf{PMAlg}(\Sigma)$\label{PMAlgSigma}.

Correspondingly, we modify the category $\textbf{Alg}_{\mathcal{B}}(\Sigma)$: take, as objects of $\textbf{Alg}_{CABA}(\Sigma)$\label{AlgCABASigma}, triples $\mathcal{A}=(A, \{\sigma_{\mathcal{A}}\}_{\sigma\in\Sigma}, \leq_{\mathcal{A}})$ such that:
\begin{enumerate}
\item $(A, \{\sigma_{\mathcal{A}}\}_{\sigma\in\Sigma})$ is a $\Sigma$-algebra;
\item $(A, \leq_{\mathcal{A}})$ is a complete, atomic Boolean algebra (no longer a bottomless Boolean algebra);
\item for $A_{a}$ the set of atoms smaller than $a$ (equal to $\emptyset$ if $a=0$), $\sigma\in\Sigma_{n}$ and $a_{1}, \dotsc  , a_{n}\in A$, 
\[\sigma_{\mathcal{A}}(a_{1}, \dotsc  , a_{n})=\sup\{\sigma_{\mathcal{A}}(b_{1}, \dotsc  , b_{n}) : (b_{1}, \dotsc  , b_{n})\in A_{a_{1}}\times\cdots\times A_{a_{n}}\}\]
(where, now, the set of which we take the supremum may be empty).
\end{enumerate}
The morphisms in $\textbf{Alg}_{CABA}(\Sigma)$, between two of its objects $\mathcal{A}$ and $\mathcal{B}$, will be functions $\varphi:A\rightarrow B$ satisfying 
\begin{enumerate}
\item for $\sigma\in\Sigma_{n}$ and $a_{1}, \dotsc  , a_{n}\in A$, $\varphi(\sigma_{\mathcal{A}}(a_{1}, \dotsc  , a_{n}))\leq_{\mathcal{B}}\sigma_{\mathcal{B}}(\varphi(a_{1}), \dotsc  , \varphi(a_{n}))$;
\item $\varphi$ is continuous, \textit{id est}, for any $S\subseteq A$ (possibly empty), $\varphi(\sup S)=\sup\{\varphi(a) : a\in S\}$;
\item $\varphi$ maps atoms of $\mathcal{A}$ to atoms of $\mathcal{B}$.
\end{enumerate}

Notice that the category $\textbf{Alg}_{CABA}(\Sigma)$ is really closely related to $\textbf{Alg}_{\mathcal{B}}(\Sigma)$: what we changed is that the underlying objects are not complete, atomic and bottomless Boolean algebras anymore, but rather complete, atomic Boolean algebras, while the morphisms are now required to map $0$ to $0$ (that is, be continuous on the empty set). And, with the categories $\textbf{PMAlg}(\Sigma)$ and $\textbf{Alg}_{CABA}(\Sigma)$ at hand, it becomes easy to show both of them are equivalent.

In one direction, the equivalence is given by the functor\\$F:\textbf{PMAlg}(\Sigma)\rightarrow \textbf{Alg}_{CABA}(\Sigma)$ which takes: a partial multialgebra $\mathcal{A}$ with universe $A$ to the powerset of $A$, equipped with operations given by
\[\sigma_{F\mathcal{A}}(A_{1}, \dotsc  , A_{n})=\bigcup_{(a_{1}, \dotsc  , a_{n})\in A_{1}\times\cdots\times A_{n}}\sigma_{\mathcal{A}}(a_{1}, \dotsc  , a_{n})\]
for possibly empty $A_{1}, \dotsc  , A_{n}\subseteq A$; and a homomorphism $\varphi:\mathcal{A}\rightarrow\mathcal{B}$ between partial $\Sigma$-multialgebras to the function given by, for a possibly empty $A'\subseteq A$, $F\varphi(A')=\{\varphi(a)\in B : a\in A'\}$.

In the other direction, the equivalence is given by $G:\textbf{Alg}_{CABA}(\Sigma)\rightarrow \textbf{PMAlg}(\Sigma)$, taking: an object of $\textbf{Alg}_{CABA}(\Sigma)$ to its set of atomic elements, equipped with multioperations such that $\sigma_{G\mathcal{A}}(a_{1}, \dotsc  , a_{n})$, for atoms $a_{1}$ through $a_{n}$, is the set of atoms smaller than $\sigma_{\mathcal{A}}(a_{1}, \dotsc  , a_{n})$ (now possibly empty, in the case the result of the operation is $0$); and a morphism $\varphi:\mathcal{A}\rightarrow\mathcal{B}$ in $\textbf{Alg}_{CABA}(\Sigma)$ to its restriction to the atoms of $\mathcal{A}$.

The equivalence between $\textbf{PMAlg}(\Sigma)$ and $\textbf{Alg}_{CABA}(\Sigma)$ is obtained directly from the equivalence between $\textbf{MAlg}(\Sigma)$ and $\textbf{Alg}_{\mathcal{B}}(\Sigma)$ by merely allowing operations to return the empty-set, and correcting whatever definitions that require a set to be non-empty accordingly. 

This may seem a long commentary for a result which changes very little, but there is a reason why we focus on $\textbf{PMAlg}(\Sigma)$: a pair $\mathcal{M}=(\mathcal{A}, D)$, such that $\mathcal{A}$ is a partial $\Sigma$-multialgebra with universe $A$ and $D$ is a subset of $A$, is called a PNmatrix, standing for a partial, non-deterministic matrix. Given a set of formulas $\Gamma\cup\{\varphi\}$ in the signature $\Sigma$, we say $\Gamma$ proves $\varphi$ according to $\mathcal{M}$ if, for every homomorphism $\nu:\textbf{F}(\Sigma, \mathcal{V})\rightarrow\mathcal{A}$ of partial multialgebras,
\[\nu(\Gamma)\subseteq D\quad\text{implies}\quad \nu(\varphi)\in D,\]
when we write $\Gamma\vDash_{\mathcal{M}}\varphi$.

It is easy to see that PNmatrices produce a broader semantics than that of Nmatrices: after all, every Nmatrix is a PNmatrix, but not vice-versa; so, in the effort of algebraizing logics, PNmatrices may offer new developments on logics that were, without them, uncharacterizable. PNmatrices were first developed in \cite{Baaz}, although the presentation given here is closer to that of \cite{CalMar}, and are indubitably very powerful; and yet, they are still not very widespread. We believe the reason for this is that they add partiality, a not so desired property for decision methods, to an already disliked by many semantical methodology, that of Nmatrices. So we extract from the logicians who have no philosophical objections against non-determinism those with no objections against partiality to find those willing to use PNmatrices.

From our studies presented here it is clear we are supporters of Nmatrices, and of PNmatrices as well, so our intention here is not to replace them, but rather to give an alternative approach, more classically behaved, to those logicians not willing to appeal to them.

\subsection{Multihomomorphisms}

Finally, consider the modified notions of homomorphism between $\Sigma$-multialgebras $\mathcal{A}=$\\$(A, \{\sigma_{\mathcal{A}}\}_{\sigma\in\Sigma})$ and $\mathcal{B}=(B, \{\sigma_{\mathcal{B}}\}_{\sigma\in\Sigma})$ we briefly met in Section \ref{Multi-valued homomorphisms}: the, perhaps, simplest multihomomorphism one can define is a function $\varphi:A\rightarrow\mathcal{P}(B)\setminus\{\emptyset\}$ satisfying, for all $\sigma\in\Sigma_{n}$ and $a_{1}, \dotsc  , a_{n}\in A$,

\[\bigcup_{a\in \sigma_{\mathcal{A}}(a_{1}, \dotsc  , a_{n})}\varphi(a)\subseteq \bigcup_{(b_{1}, \dotsc  , b_{n})\in \varphi(a_{1})\times\cdots\times\varphi(a_{n})}\sigma_{\mathcal{B}}(b_{1}, \dotsc  , b_{n});\]
then, the category with $\Sigma$-multialgebras as objects, and these multihomomorphisms as morphisms, will be denoted by $\textbf{MMAlg}(\Sigma)$. Notice that the composition of multihomomorphisms $\varphi$, between $\mathcal{A}$ and $\mathcal{B}$, and $\psi$, between $\mathcal{B}$ and $\mathcal{C}$, is not given by their composition as functions, but rather as 
\[\psi\circ\varphi(a)=\bigcup_{b\in \varphi(a)}\psi(b),\]
for any element $a$ of $\mathcal{A}$.

The idea here is quite clear: now, morphisms in our category of multialgebras return non-empty sets, much alike the multioperations in the multialgebras themselves. To find an equivalent category, we must reflect this change, what is actually rather easy: we demanded that the morphisms in $\textbf{Alg}_{\mathcal{B}}(\Sigma)$ should preserve atoms, exactly because we wanted them all to be images, under $\mathcal{P}_{\mathcal{B}}$, of some homomorphism $\varphi$ in $\textbf{MAlg}(\Sigma)$; and $\mathcal{P}_{\mathcal{B}}(\varphi)$ preserves singletons, that is, atoms in the $(\Sigma, \leq)$-algebras $\mathcal{P}_{\mathcal{B}}(\mathcal{A})$.

So, consider the category $\textbf{MAlg}_{\mathcal{B}}(\Sigma)$\label{MAlgBSigma}, with $(\Sigma, \leq)$-algebras as objects and, as morphisms between $(\Sigma, \leq)$-algebras $\mathcal{A}=(A, \{\sigma_{\mathcal{A}}\}_{\sigma\in\Sigma}, \leq_{\mathcal{A}})$ and $\mathcal{B}=(B, \{\sigma_{\mathcal{B}}\}_{\sigma\in\Sigma}, \leq_{\mathcal{B}})$, any function $\varphi:A\rightarrow B$ such that:
\begin{enumerate}
\item for all $\sigma\in\Sigma_{n}$ and $a_{1}, \dotsc  , a_{n}\in A$, 
\[\varphi(\sigma_{\mathcal{A}}(a_{1}, \dotsc  , a_{n}))\leq_{\mathcal{B}}\sigma(\varphi(a_{1}), \dotsc  , \varphi(a_{n}));\]
\item $\varphi$ is continuous on non-empty sets, meaning that for any non-empty $S\subseteq A$,
\[\varphi(\sup S)=\sup\{\varphi(a) : a\in S\}.\]
\end{enumerate}

Of course, slight changes are necessaries on the functors realizing an equivalence between $\textbf{MMAlg}(\Sigma)$ and $\textbf{MAlg}_{\mathcal{B}}(\Sigma)$. The first of them is $F:\textbf{MMAlg}(\Sigma)\rightarrow\textbf{MAlg}_{\mathcal{B}}(\Sigma)$, which takes a multialgebra $\mathcal{A}$ to $\mathcal{P}_{\mathcal{B}}(\mathcal{A})$, but takes a multihomomorphism between $\mathcal{A}$ and $\mathcal{B}$ (with universes $A$ and $B$) to the function $F\varphi:\mathcal{P}(A)\setminus\{\emptyset\}\rightarrow\mathcal{P}(B)\setminus\{\emptyset\}$ such that
\[F\varphi(A')=\bigcup_{a\in A'}\varphi(a),\]
for a non-empty $A'\subseteq A$. Reciprocally, we take the functor $G:\textbf{MAlg}_{\mathcal{B}}(\Sigma)\rightarrow\textbf{MMAlg}(\Sigma)$ mapping a $(\Sigma, \leq)$-algebra $\mathcal{A}$ to the multialgebra $\mathbb{A}(\mathcal{A})$, but a morphism $\varphi:\mathcal{A}\rightarrow\mathcal{B}$ of $\textbf{MAlg}_{\mathcal{B}}(\Sigma)$ to the multihomomorphism between $\mathbb{A}(\mathcal{A})$ and $\mathbb{A}(\mathcal{B})$ such that, for an atom $a$ of $\mathcal{A}$,
\[G\varphi(a)=A_{\varphi(a)},\]
that is, $G\varphi(a)$ is the set of atoms, in $\mathcal{B}$, below $\varphi(a)$.

\newpage
\printbibliography[segment=\therefsegment,heading=subbibliography]
\end{refsegment}

\newpage
\thispagestyle{plain}
\hspace{0pt}
\vfill
\begin{figure}[h]
\centering
\includegraphics[width=\textwidth]{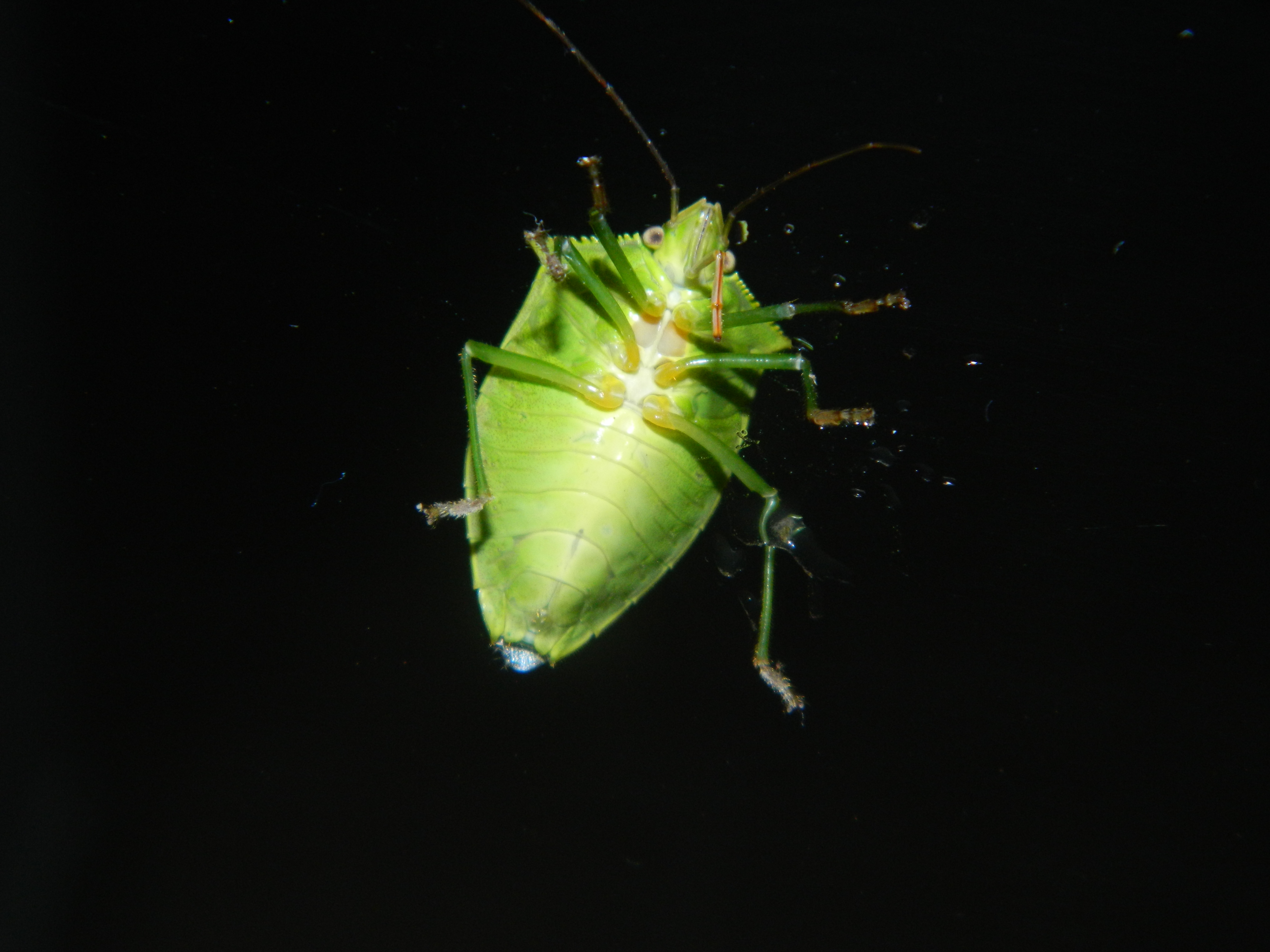}
\caption*{Specimen of \textit{Nezara viridula} against a window, Campinas, Brazil.\\Photographed by Guilherme Vicentin de Toledo, all rights reserved.}
\end{figure}
\vfill
\hspace{0pt}

\part{Paraconsistent Logic}\label{Part2}

\begin{refsegment}
\defbibfilter{notother}{not segment=\therefsegment}
\setcounter{chapter}{3}
\chapter{Logics and Nmatrices}\label{Chapter4}\label{Chapter 4}

\section{Basic notions in logic}

Given a signature $\Sigma$ and a countable set $\mathcal{V}$, to which we will refer as a set of propositional variables and whose elements we will denote by $p_{i}$ for $i\in\mathbb{N}$, the language\index{Language}, or propositional language, $\mathcal{L}_{\Sigma}$\label{language} generated by $\Sigma$ from $\mathcal{V}$, which we may denote simply by $\mathcal{L}$ if $\Sigma$ is fixed, is the set of formulas $F(\Sigma, \mathcal{V})$.

A consequence relation\index{Consequence relation}\label{vdash} $\vdash$ on a language $\mathcal{L}$ is a binary relation on $\mathcal{P}(\mathcal{L})$, and if $(\Gamma, \Delta)$ is in $\vdash$, we will simply write $\Gamma\vdash\Delta$. A logic\index{Logic}\label{logic} $\mathscr{L}$ over a signature $\Sigma$ is then a pair $(\mathcal{L}, \vdash_{\mathscr{L}})$, where $\vdash_{\mathscr{L}}$ is a consequence relation on $\mathcal{L}$ which will, at times, be denoted simply by $\vdash$.

When $\Gamma\vdash_{\mathscr{L}}\Delta$, we will say that $\Gamma$ proves $\Delta$; and for sets of formulas $\Gamma_{1}, \dotsc  , \Gamma_{n}$ and $\Delta_{1}, \dotsc  , \Delta_{m}$ and formulas $\alpha_{1}, \dotsc  , \alpha_{k}$ and $\beta_{1}, \dotsc  , \beta_{l}$, we will denote
\[\Gamma_{1}\cup\cdots\cup\Gamma_{n}\cup\{\alpha_{1}, \dotsc  , \alpha_{k}\}\vdash_{\mathscr{L}}\Delta_{1}\cup\cdots\cup\Delta_{m}\cup\{\beta_{1}, \dotsc  , \beta_{l}\}\]
simply by 
\[\Gamma_{1}, \dotsc  , \Gamma_{n}, \alpha_{1}, \dotsc  , \alpha_{k}\vdash_{\mathscr{L}}\Delta_{1}, \dotsc  , \Delta_{m}, \beta_{1}, \dotsc  , \beta_{l},\]
meaning we will drop unions and curly brackets to simplify our notation. If $\emptyset\vdash_{\mathscr{L}}\varphi$, we will simply write $\vdash_{\mathscr{L}}\varphi$, and call $\varphi$ a tautology\index{Tautology}.

\begin{definition}
A logic $\mathscr{L}$ is said to be a tarskian logic\index{Logic, Tarskian} if, for any set $\Gamma\cup\Delta\cup\{\varphi\}$ of formulas in the language of $\mathscr{L}$, the following holds:
\begin{enumerate}
\item if $\varphi\in \Gamma$, $\Gamma\vdash_{\mathscr{L}}\varphi$;
\item if $\Gamma\vdash_{\mathscr{L}}\varphi$ and $\Gamma\subseteq\Delta$, then $\Delta\vdash_{\mathscr{L}}\varphi$;
\item if $\Gamma\vdash_{\mathscr{L}}\varphi$ and $\Delta\vdash_{\mathscr{L}}\gamma$ for every $\gamma\in\Gamma$, then $\Delta\vdash_{\mathscr{L}}\varphi$.
\end{enumerate}
\end{definition}

A logic $\mathscr{L}$ is said to be finitary\index{Logic, Finitary} if, for all sets of formulas $\Gamma\cup\{\varphi\}$ in its language, if $\Gamma\vdash_{\mathscr{L}}\varphi$ then there exists a finite set $\Gamma_{0}\subseteq\Gamma$ such that $\Gamma_{0}\vdash_{\mathscr{L}}\varphi$.

A logic $\mathscr{L}$ over the signature $\Sigma$ with language $\mathcal{L}=F(\Sigma, \mathcal{V})$ is said to be structural\index{Logic, Structural} when, for every $\Gamma\cup\{\varphi\}\subseteq\mathcal{L}$ and every $\Sigma$-homomorphism $\sigma:\textbf{F}(\Sigma, \mathcal{V})\rightarrow\textbf{F}(\Sigma, \mathcal{V})$, if $\Gamma\vdash_{\mathscr{L}}\varphi$ then $\sigma(\Gamma)\vdash_{\mathscr{L}}\sigma(\varphi)$; we call any $\Sigma$-homomorphism $\sigma:\textbf{F}(\Sigma, \mathcal{V})\rightarrow \textbf{F}(\Sigma, \mathcal{V})$ a substitution.

\begin{definition}
A set of formulas $\Gamma$ in a logic $\mathscr{L}$ is said to be closed\index{Closed set of formulas} if, anytime $\Gamma\vdash_{\mathscr{L}}\varphi$, we have that $\varphi\in\Gamma$.
\end{definition}

The following theorem is due to Lindenbaum and \L oz, and will be used several times in our study.

\begin{theorem}
If $\mathscr{L}$ is a tarskian, finitary logic, and $\Gamma\cup\{\varphi\}$ is a set of formulas such that $\Gamma\not\vdash_{\mathcal{L}}\varphi$, there exists a non-trivial, closed set of formulas $\Delta$ containing $\Gamma$ and maximal with respect to not proving $\varphi$.
\end{theorem}

\begin{proof}
A proof can be found in \cite{ParLog}, in Theorem $2.2.6$.
\end{proof}

\subsection{Hilbert systems}

Perhaps the most useful way for us to define a logic is through a set of axiom schemata plus a set of rules of inference, in what is called a Hilbert system\index{Hilbert system}: given a language $\mathcal{L}=F(\Sigma, \mathcal{V})$ over the signature $\Sigma$, an $n$-ary rule of inference\index{Rule of inference}, for $n\in\mathbb{N}$, is an element of $\mathcal{L}^{n}\times\mathcal{L}$. A rule of inference $((\alpha_{1}, \dotsc  , \alpha_{n}), \alpha)$ will usually be denoted by $\alpha_{1}, \dotsc  , \alpha_{n}|\alpha$\label{ruleofinference} or
\[\frac{\alpha_{1}, \dotsc  , \alpha_{n}}{\alpha}.\]
 An axiom scheme\index{Axiom scheme} is any formula $\alpha$, and its instances\index{Axiom scheme, Instance of an} are $\sigma(\alpha)$, for any $\Sigma$-homomorphism $\sigma:\textbf{F}(\Sigma, \mathcal{V})\rightarrow\textbf{F}(\Sigma, \mathcal{V})$: notice that, technically speaking, an axiom is a $0$-ary rule of inference; regardless, we will still treat then differently, given that historically this is the most common approach. More generally, given a rule of inference $\alpha_{1}, \dotsc  , \alpha_{n}|\alpha$ and a homomorphism $\sigma:\textbf{F}(\Sigma, \mathcal{V})\rightarrow\textbf{F}(\Sigma, \mathcal{V})$, $\sigma(\alpha_{1}), \dotsc  , \sigma(\alpha_{n})|\sigma(\alpha)$ is an instance\index{Rule of inference, Instance of a} of the aforementioned rule.

So, given a set of axiom schemata $\mathfrak{A}$\label{frakA} and a set of rules of inference $\mathfrak{R}$\label{frakR}, we can create a logic $\mathscr{L}$ dependent on $\mathfrak{A}$ and $\mathfrak{R}$ such that, given formulas $\Gamma\cup\{\varphi\}$ in $\mathcal{L}$, $\Gamma\vdash_{\mathscr{L}}\varphi$ if and only if there exists a sequence of formulas $\varphi_{1}, \dotsc  , \varphi_{m}$, said to be a demonstration\index{Demonstration} of $\varphi$ from $\Gamma$, such that $\varphi_{m}=\varphi$ and $\varphi_{i}$, for $i\in\{1, \dotsc  , m\}$, is either:
\begin{enumerate}
\item an element of $\Gamma$, when we call it a premise\index{Premise};
\item an instance of axiom $\sigma(\alpha)$, for an $\alpha\in\mathfrak{A}$ and a $\Sigma$-homomorphism $\sigma:\textbf{F}(\Sigma, \mathcal{V})\rightarrow\textbf{F}(\Sigma, \mathcal{V})$;
\item a formula $\sigma(\alpha)$, for a $\Sigma$-homomorphism $\sigma:\textbf{F}(\Sigma, \mathcal{V})\rightarrow\textbf{F}(\Sigma, \mathcal{V})$ and $\alpha_{1}, \dotsc  , \alpha_{n}|\alpha$ an $n$-ary rule of inference in $\mathfrak{R}$ such that there exist $1\leq i_{1}, \dotsc  , i_{n}<i$ with $\varphi_{i_{1}}=\sigma(\alpha_{1}), \dotsc , \varphi_{i_{n}}=\sigma(\alpha_{n})$.
\end{enumerate}

\begin{theorem}
If $\mathscr{L}$ is the logic obtained from the set of axiom schemata $\mathfrak{A}$ and the set of rules of inference $\mathfrak{R}$, $\mathscr{L}$ is a finitary, tarskian logic.
\end{theorem}

\begin{proof}
Let $\Gamma\cup\Delta\cup\{\varphi\}$ be a set of formulas in $\mathcal{L}$.

\begin{enumerate}
\item If $\varphi\in\Gamma$, $\varphi_{1}$ is a demonstration of $\varphi$ from $\Gamma$, with $\varphi_{1}=\varphi$: it is, in fact, a demonstration since $\varphi_{1}$ is an element of $\Gamma$. This way, $\Gamma\vdash_{\mathscr{L}}\varphi$.

\item Suppose $\Gamma\vdash_{\mathscr{L}}\varphi$ and $\Gamma\subseteq\Delta$: given $\varphi_{1}, \dotsc  , \varphi_{m}$, a demonstration of $\varphi$ from $\Gamma$, we state that $\varphi_{1}, \dotsc  , \varphi_{m}$ is also a demonstration of $\varphi$ from $\Delta$. Clearly, first of all, we have that $\varphi_{m}=\varphi$. Then $\varphi_{i}$ is either:
\begin{enumerate}
\item an element of $\Gamma$, and since $\Gamma\subseteq \Delta$, an element of $\Delta$;
\item an instance of axiom $\sigma(\alpha)$, for $\alpha$ an axiom scheme and $\sigma$ a substitution;
\item a formula $\sigma(\alpha)$, for $\sigma$ a substitution and $\alpha_{1}, \dotsc  , \alpha_{n}|\alpha$ a rule of inference for which there exist $1\leq i_{1}, \dotsc  , i_{n}<i$ with $\varphi_{i_{1}}=\sigma(\alpha_{1}), \dotsc  , \varphi_{i_{n}}=\sigma(\alpha_{n})$.
\end{enumerate}
This proves $\varphi_{1}, \dotsc  , \varphi_{m}$ is a demonstration of $\varphi$ from $\Delta$, and therefore $\Delta\vdash_{\mathscr{L}}\varphi$.

\item Suppose $\Gamma\vdash_{\mathscr{L}}\varphi$, and let $\varphi_{1}, \dotsc  , \varphi_{m}$ be a demonstration of $\varphi$ from $\Gamma$. We take $\Gamma_{0}=\Gamma\cap\{\varphi_{1}, \dotsc  , \varphi_{m}\}$: we state then that $\varphi_{1}, \dotsc  , \varphi_{m}$ is a demonstration of $\varphi$ from $\Gamma_{0}\subseteq\Gamma$, which is a finite set.

First of all, $\varphi_{n}=\varphi$. Then $\varphi_{i}$ is either:
\begin{enumerate}
\item an element of $\Gamma$, and then $\varphi_{i}\in\Gamma\cap\{\varphi_{1}, \dotsc  , \varphi_{n}\}$, implying $\varphi_{i}$ is in $\Gamma_{0}$;
\item an instance $\sigma(\alpha)$ of the axiom scheme $\alpha$;
\item a formula $\sigma(\alpha)$, for a rule of inference $\alpha_{1}, \dotsc  , \alpha_{n}|\alpha$ and $1\leq i_{1}, \dotsc  , i_{n}< i$ such that $\varphi_{i_{j}}=\sigma(\alpha_{j})$, for $j\in\{1, \dotsc  , n\}$.
\end{enumerate}
 So, to summarize, we have that $\Gamma_{0}\vdash_{\mathscr{L}}\varphi$, and therefore $\mathscr{L}$ is finitary.

\item Finally, suppose $\Gamma\vdash_{\mathscr{L}}\varphi$ and that $\Delta\vdash_{\mathscr{L}}\gamma$ for every $\gamma\in\Gamma$. By the item above, there exists a finite $\Gamma_{0}\subseteq\Gamma$ such that $\Gamma_{0}\vdash_{\mathscr{L}}\varphi$, so let $\varphi_{1}, \dotsc  , \varphi_{M}$ be a demonstration of $\varphi$ from $\Gamma_{0}$, and let $\Gamma_{0}=\{\gamma_{1}, \dotsc  , \gamma_{l}\}$: we will prove that we can drop one of the premises of $\Gamma_{0}$, e.g. $\gamma_{1}$, by adding as premises $\Delta$, meaning
\[\Delta\cup\{\gamma_{1}, \dotsc  , \gamma_{l}\}\setminus\{\gamma_{1}\}\vdash_{\mathscr{L}}\varphi,\]
and by an inductive argument, we will have that $\Delta\vdash_{\mathscr{L}}\varphi$.

So, suppose $\gamma_{1}\neq\varphi$ (otherwise we already have that $\Delta\vdash_{\mathscr{L}}\varphi$), and consider the subsequence $\varphi'_{1}, \dotsc  , \varphi'_{m}$ of $\varphi_{1}, \dotsc  , \varphi_{M}$ of formulas different from $\gamma_{1}$; we then take a demonstration $\psi_{1}, \dotsc  , \psi_{k}$ of $\gamma_{1}$ from $\Delta$, and state that the sequence $\phi_{1}, \dotsc  , \phi_{k+m}$ equal to
\[\psi_{1}, \dotsc  , \psi_{k}, \varphi'_{1}, \dotsc  , \varphi'_{m}\]
is a demonstration of $\varphi$ from $\Delta\cup\Gamma_{0}\setminus\{\gamma_{1}\}$.

First of all, since $\varphi'_{1}, \dotsc  , \varphi'_{m}$ is the subsequence of $\varphi_{1}, \dotsc  , \varphi_{M}$ of formulas different from $\gamma_{1}$ and, by hypothesis, $\gamma_{1}\neq\varphi$, we have that $\varphi'_{m}=\varphi_{M}=\varphi$. Now, if $i\in \{1, \dotsc  , k\}$, $\phi_{i}=\psi_{i}$, and then:
\begin{enumerate}
\item if $\psi_{i}$ is a premise in $\Delta$, clearly $\phi_{i}$ is a premise in the larger set $\Delta\cup\Gamma_{0}\setminus\{\gamma_{1}\}$;
\item if $\psi_{i}$ is an instance of axiom $\sigma(\alpha)$, clearly $\phi_{i}$ is also an instance of $\alpha$;
\item if $\psi_{i}$ is the formula $\sigma(\alpha)$, for $\alpha_{1}, \dotsc  , \alpha_{n}|\alpha$ an $n$-ary inference rule and $1\leq i_{1}, \dotsc  , i_{n}<i$ such that $\psi_{i_{1}}=\sigma(\alpha_{1}), \dotsc , \psi_{i_{n}}=\sigma(\alpha_{n})$, then, since 
\[\psi_{i_{1}}=\phi_{i_{1}}, \dotsc  , \psi_{i_{n}}=\phi_{i_{n}},\]
we have that $\phi_{i}$ is also $\sigma(\alpha)$, for $\alpha_{1}, \dotsc  , \alpha_{n}|\alpha$ an inference rule and $1\leq i_{1}, \dotsc  , i_{n}<i$ such that $\phi_{i_{1}}=\sigma(\alpha_{1}), \dotsc  , \phi_{i_{n}}=\sigma(\alpha_{n})$.
\end{enumerate}

If $i\in\{k+1, \dotsc  , k+m\}$, $\phi_{i}=\varphi'_{i-k}$, and then
\begin{enumerate}
\item if $\varphi'_{i-k}$ is a premise of the demonstration $\varphi_{1}, \dotsc  , \varphi_{M}$, and therefore in $\Gamma_{0}$, since $\varphi'_{1}, \dotsc  ,$\\$\varphi'_{m}$ is the subsequence of elements different from $\gamma_{1}$ we find $\varphi'_{i-k}\in\Gamma_{0}\setminus\{\gamma_{1}\}$, and $\phi_{i}$ is therefore in $\Delta\cup\Gamma_{0}\setminus\{\gamma_{1}\}$;
\item if $\varphi'_{i-k}$ is an instance of axiom $\sigma(\alpha)$, evidently $\phi_{i}$ is also an instance of the same axiom;
\item if $\varphi'_{i-k}=\varphi_{j}$ is $\sigma(\alpha)$, for an inference rule $\alpha_{1}, \dotsc  , \alpha_{n}|\alpha$ and $1\leq i_{1}, .. , i_{n}<j$ such that $\varphi_{i_{1}}=\sigma(\alpha_{1}), \dotsc  , \varphi_{i_{n}}=\sigma(\alpha_{n})$, there are two cases to consider for a $\varphi_{i_{s}}$: 
\begin{enumerate}
\item either $\varphi_{i_{s}}\neq\gamma_{1}$, and so there exists $k+1\leq j_{s}<k+i$ such that $\varphi'_{j_{s}-k}=\varphi_{i_{s}}$;
\item or $\varphi_{i_{s}}=\gamma_{1}$, when we make $j_{s}=k$ and therefore $\phi_{j_{s}}=\gamma_{1}$;
\end{enumerate}
either way, we find that $\phi_{i}$ is $\sigma(\alpha)$, for $\alpha_{1}, \dotsc  , \alpha_{n}|\alpha$ an $n$-ary inference rule and $1\leq j_{1}, \dotsc  , j_{n}< i$ such that $\phi_{j_{1}}=\sigma(\alpha_{1}), \dotsc  , \phi_{j_{n}}=\sigma(\alpha_{n})$.
\end{enumerate}
\end{enumerate}

\end{proof}

Again, let $\mathscr{L}$ be the logic with axiom schemata $\mathfrak{A}$ and rules of inference $\mathfrak{R}$: beyond being a finitary, tarskian logic, we can prove that is also structural. Assume $\Gamma\vdash_{\mathscr{L}}\varphi$ and let $\varphi_{1}, \dotsc  , \varphi_{m}$ be a demonstration of $\varphi$ from $\Gamma$ and $\sigma:\textbf{F}(\Sigma, \mathcal{V})\rightarrow\textbf{F}(\Sigma, \mathcal{V})$ a $\Sigma$-homomorphism.

We state then that $\sigma(\varphi_{1}), \dotsc  , \sigma(\varphi_{m})$ is a demonstration of $\sigma(\varphi)$ from $\sigma(\Gamma)$. First of all, we have that, since $\varphi_{m}=\varphi$, $\sigma(\varphi_{m})=\sigma(\varphi)$. 
\begin{enumerate}
\item If $\varphi_{i}$ is a premise in $\Gamma$, then $\sigma(\varphi_{i})$ is a premise in $\sigma(\Gamma)$.
\item If $\varphi_{i}$ is an instance $\tau(\alpha)$ of an axiom $\alpha$, since $\sigma\circ\tau$ is still a $\Sigma$-homomorphism, then $\sigma(\varphi_{i})$ is an instance $\sigma\circ\tau(\alpha)$ of the same axiom $\alpha$.
\item If $\varphi_{i}$ equals $\tau(\alpha)$, for a $\Sigma$-homomorphism $\tau$, a rule of inference $\alpha_{1}, \dotsc  , \alpha_{n}|\alpha$ and $1\leq i_{1}, \dotsc  , i_{n}<i$ such that $\varphi_{i_{1}}=\tau(\alpha_{1}), \dotsc  , \varphi_{i_{n}}=\tau(\alpha_{n})$, then $\sigma\circ\tau$ is still a $\Sigma$-homomorphism, and $\sigma(\varphi_{i})$ equals $\sigma\circ\tau(\alpha)$ for a rule of inference $\alpha_{1}, \dotsc  , \alpha_{n}|\alpha$ and $1\leq i_{1}, \dotsc  , i_{n}<i$ such that $\sigma(\varphi_{i_{1}})=\sigma\circ\tau(\alpha_{1}), \dotsc  , \sigma(\varphi_{i_{n}})=\sigma\circ\tau(\alpha_{n})$.
\end{enumerate}

One important thing to notice is that, when presenting an axiom or rule of inference, is common to use a not so precise notation: suppose our signature has a binary symbol $\uparrow$, and that $p_{1}\uparrow p_{2}$ is an axiom; we will usually write it as $\alpha\uparrow\beta$, using $\alpha$ and $\beta$ as meta-variables for formulas, although, in principle, this notation is incorrect. The same will be done for rules of inference, where, to name one example, a rule such as
\[\frac{p_{1}, p_{1}\uparrow p_{2}}{{\sim}  p_{2}}\]
will be denoted by $\alpha, \alpha\uparrow\beta|{\sim} \beta$.

\subsection{Paraconsistent logics and $\textbf{LFI}'s$}

We will understand as classical propositional logic\index{Logic, Classical propositional} the logic we will denote by $\textbf{CPL}$\label{CPL} over the signature $\Sigma^{\textbf{CPL}}$\label{SigmaCPL}, such that $\Sigma^{\textbf{CPL}}_{0}=\{\bot, \top\}$, $\Sigma^{\textbf{CPL}}_{1}=\{{\sim} \}$, $\Sigma^{\textbf{CPL}}_{2}=\{\vee, \wedge, \rightarrow\}$ and $\Sigma^{\textbf{CPL}}_{n}=\emptyset$ for $n>2$, with axiom schemata
\begin{enumerate}
\item[\textbf{Ax\: 1}] $\alpha\rightarrow(\beta\rightarrow\alpha)$;
\item[\textbf{Ax\: 2}] $\big(\alpha\rightarrow (\beta\rightarrow \gamma)\big)\rightarrow\big((\alpha\rightarrow\beta)\rightarrow(\alpha\rightarrow\gamma)\big)$;
\item[\textbf{Ax\: 3}] $\alpha\rightarrow\big(\beta\rightarrow(\alpha\wedge\beta)\big)$;
\item[\textbf{Ax\: 4}] $(\alpha\wedge\beta)\rightarrow \alpha$;
\item[\textbf{Ax\: 5}] $(\alpha\wedge\beta)\rightarrow \beta$;
\item[\textbf{Ax\: 6}] $\alpha\rightarrow(\alpha\vee\beta)$;
\item[\textbf{Ax\: 7}] $\beta\rightarrow(\alpha\vee\beta)$;
\item[\textbf{Ax\: 8}] $(\alpha\rightarrow\gamma)\rightarrow\Big((\beta\rightarrow\gamma)\rightarrow \big((\alpha\vee\beta)\rightarrow\gamma\big)\Big)$;
\item[\textbf{Ax\: 9}] $(\alpha\rightarrow\beta)\rightarrow\big((\alpha\rightarrow{\sim} \beta)\rightarrow{\sim} \alpha\big)$;
\item[\textbf{Ax\: 10}] $\alpha\rightarrow({\sim} \alpha\rightarrow\beta)$;
\item[\textbf{Ax\: 11}] $\alpha\vee{\sim} \alpha$;
\item[\textbf{Ax\: 12}] $\bot\rightarrow\alpha$;
\item[\textbf{Ax\: 13}] $\alpha\rightarrow\top$;
\end{enumerate}
and following the inference rule of Modus Ponens
\[\frac{\alpha\quad\alpha\rightarrow\beta}{\beta}.\]

As opposed to $\textbf{CPL}$, a tarskian logic $\mathscr{L}$ will be said to be paraconsistent\index{Logic, Paraconsistent} when it possesses an unary symbol "$\neg $", that we shall refer to as a negation, such that there exist formulas $\alpha$ and $\beta$ of $\mathscr{L}$ satisfying $\alpha, \neg\alpha\not\vdash_{\mathscr{L}}\beta$. A good, intuitive way to look at paraconsistent logics is through the paradigm introduced by da Costa that the formulas of such a logic can, at times, be divided into well-behaved and badly-behaved (nowadays more commonly referred to as inconsistent) ones. Think of a scientist performing an experiment: from her school years, said scientist knows that opposite poles of a magnet attract each other; so, if during the experiment two poles she thought to be opposite are instead repealing each other, she can be certain that she was wrong, and the two poles are instead the same, both north or both south. Her reasoning works because the sentence ``opposite poles attract'' is well-behaved (that is, it can not coexist with its negation) and true, while what she had apparently observed was that ``opposite poles repeal'', a well-behaved although false sentence.

Now, suppose that our scientist has discovered that all monopoles (conjectured magnetic particles that, instead of having both north and south poles, have only one pole, hence the name) have exactly the same pole and is now testing two hypotheses: ``all monopoles are north poles'' and ``all monopoles are south poles''. Again, we have contradictory sentences, as was the case with ``opposite poles attract'' and ``opposite poles repeal'', but this time our scientist can not reach the conclusion that, at some prior step, she made a mistake: why is that?\footnote{Notice that we are presenting this logical problem in a somewhat convoluted way in order to avoid explicit explosivity of the theory at hand: scientists, and most mathematicians, appear to prefer thinking about these topics in terms of non-contradiction (what is eventually equivalent).} Simply put: none of the hypotheses to be tested, ``all monopoles are north poles'' and ``all monopoles are south poles'', is well-behaved, meaning that we can not discard their negations as necessarily false.

What we are doing is dividing those sentences found in science between true-or-false sentences (such as those involving attracting or repealing poles) and hypothetical sentences (such as those involving all monopoles). But this is not exclusive to scientists: mathematicians also have true-or-false sentences and conjectural sentences: after all, some mathematicians believe that, e.g., the Riemann hypothesis is true, while others believe it to be false, and that doesn't make mathematics as a whole trivial. In daily life, conjectural and hypothetical sentences may be replaced with rumors, that can also be contradictory without making logical reasoning unfeasible. da Costa's hierarchy, which deals with these distinctions plus higher degrees of consistency, will be studied algebraically in Chapters \ref{Chapter5} and \ref{Chapter6}. 

If our logic $\mathscr{L}$ also possesses a binary symbol "$\rightarrow$" satisfying the deduction meta-theorem, i.e., $\Gamma, \psi\vdash_{\mathscr{L}}\varphi$ if and only if $\Gamma\vdash_{\mathscr{L}}\psi\rightarrow\varphi$, then the condition that there exist formulas $\alpha$ and $\beta$ such that $\alpha, \neg\alpha\not\vdash_{\mathscr{L}}\beta$ is equivalent to the existence of formulas $\alpha$ and $\beta$ such that $\not\vdash_{\mathscr{L}}\alpha\rightarrow(\neg\alpha\rightarrow\beta)$, which is know as the failure of the explosion law\index{Explosion law}; if this result were to be true for any $\alpha$ and $\beta$, that would mean the explosion law would be valid, which indeed happens in $\textbf{CPL}$ as one can see by its $\textbf{Ax\: 10}$.

If a formula $\psi$ has all its propositional variables among the set $\{p_{1}, \dotsc  , p_{n}\}$, we may write $\psi(p_{1}, \dotsc  , p_{n})$\label{varformula}: to define logics of formal inconsistency, we shall need a set of formulas $\bigcirc(p)$, each of them dependent exactly on the propositional variable $p$, that is, each one of them contains $p$ and solely $p$ as a variable.

And if we apply to a formula $\psi(p_{1}, \dotsc  , p_{n})$ a homomorphism $\sigma$ sending $p_{i}$ to the formula $\alpha_{i}$, for $i\in\{1, \dotsc  , n\}$, we shall write $\psi(\alpha_{1}, \dotsc  , \alpha_{n})$ for $\sigma(\psi)$: this way, $\bigcirc(\alpha)$ will be the set of formulas obtained from $\bigcirc(p)$ after we apply to each of its elements the homomorphism taking $p$ to $\alpha$, that is,
\[\bigcirc(\alpha)=\{\psi(\alpha)\ :\  \psi(p)\in\bigcirc(p)\}.\]

\begin{definition}\label{Definition of LFI}
A tarskian, finitary and structural logic $\mathscr{L}$ containing a negation and a set of formulas $\bigcirc(p)\neq\emptyset$ depending exactly on the propositional variable $p$ is a logic of formal inconsistency ($\textbf{LFI}$)\index{Logic of formal inconsistency}\label{LFI} when:
\begin{enumerate}
\item there exist formulas $\phi$ and $\theta$ in $\mathscr{L}$ such that $\phi, \neg\phi\not\vdash_{\mathscr{L}}\theta$;
\item there exist formulas $\alpha$ and $\beta$ in $\mathscr{L}$ such that
\begin{enumerate}
\item $\bigcirc(\alpha), \alpha\not\vdash_{\mathscr{L}}\beta$ and
\item $\bigcirc(\alpha), \neg\alpha\not\vdash_{\mathscr{L}}\beta$;
\end{enumerate}
\item for all formulas $\varphi$ and $\psi$ in $\mathscr{L}$,
\[\bigcirc(\varphi), \varphi, \neg\varphi\vdash_{\mathscr{L}}\psi.\]
\end{enumerate}
\end{definition}

The logic $\mathscr{L}$ will be called a weak $\textbf{LFI}$\index{Weak $\textbf{LFI}$} if the second condition of Definition \ref{Definition of LFI} is replaced by there existing formulas $\alpha_{1}$ and $\beta_{1}$ such that
\[\bigcirc(\alpha_{1}), \alpha_{1}\not\vdash_{\mathscr{L}}\beta_{1}\]
and $\alpha_{2}$ and $\beta_{2}$, possibly different from respectively $\alpha_{1}$ and $\beta_{1}$, such that
\[\bigcirc(\alpha_{2}), \neg\alpha_{2}\not\vdash_{\mathscr{L}}\beta_{2}.\]

The logic $\mathscr{L}$ will be called a strong $\textbf{LFI}$\index{Strong $\textbf{LFI}$} if the first and second conditions of Definition \ref{Definition of LFI} are replaced by there existing a single pair of formulas $\alpha$ and $\beta$ satisfying simultaneously
\begin{enumerate}
\item $\alpha, \neg\alpha\not\vdash_{\mathscr{L}}\beta$,
\item $\bigcirc(\alpha), \alpha\not\vdash_{\mathscr{L}}\beta$ and
\item $\bigcirc(\alpha), \neg\alpha\not\vdash_{\mathscr{L}}\beta$.
\end{enumerate}
All these definitions are the Definitions $2.1.7$, $2.1.8$ and $2.1.9$ of \cite{ParLog}.

When the set $\bigcirc(p)$ contains a single formula, this formula will be denoted by $\circ p$\label{circ}, being the consistency of $p$. More often than not, we will want "$\circ$" to be a primitive connective, as well as consistency a primitive notion.

The $\textbf{LFI}$ we will most frequently work with is one that, beyond having a primitive consistency, emulates only the positive fragment of classical propositional logic and excluded middle: this way, $\textbf{mbC}$\label{mbC} is often regarded as the simplest logic of formal inconsistency. It has the signature we will denote by $\Sigma_{\textbf{LFI}}$\label{SigmaLFI}, with $(\Sigma_{\textbf{LFI}})_{0}=\emptyset$, $(\Sigma_{\textbf{LFI}})_{1}=\{\neg , \circ\}$, $(\Sigma_{\textbf{LFI}})_{2}=\{\vee, \wedge, \rightarrow\}$ and $(\Sigma_{\textbf{LFI}})_{n}=\emptyset$ for $n>2$; it has as axiom schemata
\begin{enumerate}
\item[\textbf{Ax\: 1}] $\alpha\rightarrow(\beta\rightarrow\alpha)$;
\item[\textbf{Ax\: 2}] $\big(\alpha\rightarrow (\beta\rightarrow \gamma)\big)\rightarrow\big((\alpha\rightarrow\beta)\rightarrow(\alpha\rightarrow\gamma)\big)$;
\item[\textbf{Ax\: 3}] $\alpha\rightarrow\big(\beta\rightarrow(\alpha\wedge\beta)\big)$;
\item[\textbf{Ax\: 4}] $(\alpha\wedge\beta)\rightarrow \alpha$;
\item[\textbf{Ax\: 5}] $(\alpha\wedge\beta)\rightarrow \beta$;
\item[\textbf{Ax\: 6}] $\alpha\rightarrow(\alpha\vee\beta)$;
\item[\textbf{Ax\: 7}] $\beta\rightarrow(\alpha\vee\beta)$;
\item[\textbf{Ax\: 8}] $(\alpha\rightarrow\gamma)\rightarrow\Big((\beta\rightarrow\gamma)\rightarrow \big((\alpha\vee\beta)\rightarrow\gamma\big)\Big)$;
\item[$\textbf{Ax\: 9}^{*}$] $(\alpha\rightarrow \beta)\vee\alpha$;
\item[$\textbf{Ax\: 11}^{*}$] $\alpha\vee\neg \alpha$,
\end{enumerate}
plus 
\[\tag{\textbf{bc1}}\circ\alpha\rightarrow(\alpha\rightarrow(\neg \alpha\rightarrow \beta)),\]
and as inference rules that of Modus Ponens,
\[\frac{\alpha\quad\alpha\rightarrow\beta}{\beta}.\]

Other logics of formal incompatibility we will make use of are:
\begin{enumerate}
\item $\textbf{mbCciw}$\label{mbCciw}, obtained from the Hilbert system for $\textbf{mbC}$ by adding the axiom schema\label{ciw}
\[\tag{\textbf{ciw}} \circ\alpha\vee(\alpha\wedge\neg \alpha);\]
\item $\textbf{mbCci}$\label{mbCci}, obtained from $\textbf{mbC}$ by adding\label{ci}
\[\tag{\textbf{ci}}\neg \circ\alpha\rightarrow(\alpha\wedge\neg \alpha);\]
\item $\textbf{mbCcl}$\label{mbCcl}, obtained from $\textbf{mbC}$ by adding\label{cl}
\[\tag{\textbf{cl}}\neg(\alpha\wedge\neg \alpha)\rightarrow\circ\alpha.\]
\end{enumerate}


\section{Matrices and generalizations}

A logical matrix\index{Matrix, Logical} over a signature $\Sigma$ is a pair $\mathcal{M}=(\mathcal{A}, D)$ such that $\mathcal{A}=(A, \{\sigma_{\mathcal{A}}\}_{\sigma\in\Sigma})$ is a $\Sigma$-algebra and $D\subseteq A$ is said to be the set of designated elements\index{Designated elements} of the matrix. Given formulas $\Gamma\cup\{\varphi\}$ over the signature $\Sigma$, we say $\Gamma$ semantically proves $\varphi$ according to $\mathcal{M}$ (and write $\Gamma\vDash_{\mathcal{M}}\varphi$\label{semant.proves}) if, for every homomorphism $\nu:\textbf{F}(\Sigma, \mathcal{V})\rightarrow\mathcal{A}$ such that $\nu(\gamma)\in D$, for every $\gamma\in\Gamma$, one has that $\nu(\varphi)\in D$.

A logic $\mathcal{L}$, over the signature $\Sigma$, is said to be characterized\index{Characterized} by $\mathcal{M}$ if, for every set of formulas $\Gamma\cup\{\varphi\}$ over $\Sigma$, $\Gamma\vdash_{\mathcal{L}}\varphi$ if and only if $\Gamma\vDash_{\mathcal{M}}\varphi$. 

Quite analogously, given a class $\mathbb{M}$ of matrices $\mathcal{M}$ over the same signature, we write $\Gamma\vDash_{\mathbb{M}}\varphi$\label{vDashM} if $\Gamma\vDash_{\mathcal{M}}\varphi$ for every $\mathcal{M}\in\mathbb{M}$, and a logic $\mathcal{L}$ is said to be characterized by $\mathbb{M}$ when $\Gamma\vdash_{\mathcal{L}}\varphi$ if and only if $\Gamma\vDash_{\mathbb{M}}\varphi$.

It is a well know result by W\'ojcicki\index{W\'ojciki} that every tarskian logic can be characterized by a suitable class of logical matrices, see \cite{Woj}. And, although one could see such a result as settling the matter of logical matrices, many times the suitable class of logical matrices obtained for a logic is not efficient, as in, to name one example, its algebras can be too large or complex, or the class may be infinite. So alternatives have been offered to classes of logical matrices, as in W\'ojcicki \cite{Woj2}.

Another approach is the one we will call that of restricted matrices\index{Matrix, Restricted}, or Rmatrices\index{Rmatrix}, (although that was not its original nomenclature, Piochi\index{Piochi} named then $\mathcal{E}$-matrices instead), see \cite{Piochi} or, for an approach focusing more on structurality of the related closure operators, \cite{Piochi2} and \cite{Piochi3}. A restricted matrix over a signature $\Sigma$ is a triple $\mathcal{M}=(\mathcal{A}, D, \mathcal{F})$ such that:
\begin{enumerate}
\item $\mathcal{A}=(A,\{\sigma_{\mathcal{A}}\}_{\sigma\in\Sigma})$ is a $\Sigma$-algebra;
\item $D$ is a subset of $A$;
\item $\mathcal{F}$ is a set of homomorphisms from $\textbf{F}(\Sigma, \mathcal{V})$ to $\mathcal{A}$.
\end{enumerate}

The set $\mathcal{F}$ will be called the set of restrictions\index{Restrictions}. Given formulas $\Gamma\cup\{\varphi\}$ over the signature $\Sigma$, we say $\Gamma$ proves $\varphi$ according to a restricted matrix $\mathcal{M}$ (also over $\Sigma$) and write $\Gamma\vDash_{\mathcal{M}}\varphi$ if, for every homomorphism $\nu\in\mathcal{F}$, $\nu(\gamma)\in D$, for every $\gamma\in\Gamma$, implies $\nu(\varphi)\in D$; we say a logic $\mathcal{L}$ is characterized by a restricted matrix $\mathcal{M}$ when $\Gamma\vdash_{\mathcal{L}}\varphi$ if and only if $\Gamma\vDash_{\mathcal{M}}\varphi$, and it is possible to prove that every tarskian logic can be characterized by a, potentially infinite, restricted matrix (see Piochi's \cite{Piochi} and \cite{Piochi2}).

And yet, given sometimes this is the most efficient approach, we may define $\vDash_{\mathbb{M}}$ for $\mathbb{M}$ a class of restricted matrices in much the same way we did for a class of matrices. 

Many of these generalizations have one simple objective, which is to offer a reasonable decision method for a logic; of course, methods that depend on infinite matrices or infinite classes of matrices are not really decision methods, and so finite matrices, or finite sets of finite matrices, are the preferable outcomes of the algebraization process of a logic. Of course, such outcomes are not always possible, as many results on uncharacterizability of logics show: probably the earliest one is G{\"o}del's\index{G{\"o}del} proof that propositional, intuitionistic logic is not characterizable by a single, finite logical matrix, found in \cite{Godel}.

Inspired by G{\"o}del's proof, Dugundji (\cite{Dugundji})\index{Dugundji} proved that no system lying between the modal logics $\textbf{S1}$ and $\textbf{S5}$ admits a single, finite logic matrix which characterizes it as well, and in many places we may refer to results regarding uncharacterizability of logics as ``Dugundji-like'' theorems. As we enter the domain of paraconsistent logics, these results abound, given the intrinsic complexity of many of those systems: one example would be Avron's\index{Avron} proof that many logics, including da Costa's $C_{1}$, do not possess a characterizing finite Nmatrix, or even a characterizing finite set of finite Nmatrices (\cite{Avron, Avron3}).

Most of the work shown in Sections \ref{Restricted Nmatrices}, \ref{A brief history of RNmatrices}, \ref{Structurality} and \ref{Examples: characterizing some logics with RNmatrices} was submitted to an online repository in \cite{CostaRNmatrix}, and then finally published in \cite{TwoDecisionProcedures}.


\subsection{Restricted Nmatrices}\label{Restricted Nmatrices}

As we explained before, although the problem of finding semantics for a given logic is somewhat solved in the case the logic at hand is tarskian, the corresponding class of matrices or even restricted matrix associated to the logic may not be sufficiently efficient, and in the realm of paraconsistency and, furthermore, in the presence of incompatibility, which we will further study ahead, non-deterministic matrices have been proven fruitful. The first approach to non-deterministic matrices is found in the work of Rescher\index{Rescher} (see \cite{Rescher}) and Ivlev\index{Ivlev} (see \cite{Lev}, \cite{Lev2}, \cite{Lev3} and \cite{Lev4}), although our reasoning will be closer to that of \cite{Avron}.

\begin{definition}
Given a signature $\Sigma$, a non-deterministic matrix over $\Sigma$, or Nmatrix\index{Nmatrix}\index{Matrix, Non-deterministic}, is a pair $\mathcal{M}=(\mathcal{A}, D)$ such that $\mathcal{A}=(A, \{\sigma_{\mathcal{A}}\}_{\Sigma})$ is a $\Sigma-$multialgebra and $D$ is a subset of $A$, said to be its set of designated elements. 
\end{definition}

One may check definition $6.3.1$ of \cite{ParLog} and \cite{Avron} for equivalent definitions, with slightly different emphases. 

Given formulas $\Gamma\cup\{\varphi\}$ over the signature $\Sigma$ and an Nmatrix $\mathcal{M}$, we say that $\Gamma$ proves $\varphi$ according to $\mathcal{M}$ if, for every homomorphism of multialgebras $\nu:\textbf{F}(\Sigma, \mathcal{V})\rightarrow \mathcal{A}$ such that $\nu(\gamma)\in D$ for every $\gamma\in\Gamma$, one has $\nu(\varphi)\in D$, when we then write $\Gamma\vDash_{\mathcal{M}}\varphi$: notice how close this definition is of that for $\vDash_{\mathcal{M}}$ when $\mathcal{M}$ is simply a matrix. 

One, perhaps very important, observation is that Nmatrices, as defined, are slightly redundant, being enough to add a symbol to our signature in order to work only their underlying multialgebra: consider an Nmatrix $\mathcal{M}=(\mathcal{A}, D)$ over the signature $\Sigma$, add a symbol $\top$ to $\Sigma_{0}$ therefore producing the signature $\Sigma^{\top}$\label{Sigmatop} and define the $\Sigma^{\top}$-multialgebra $\mathcal{A}^{\top}$ such that
\[\sigma_{\mathcal{A}^{\top}}=\sigma_{\mathcal{A}}, \quad\text{for every $\sigma$ in $\Sigma$},\]
and $\top_{\mathcal{A}^{\top}}=D$. Then, given formulas $\Gamma\cup\{\varphi\}$ in the signature $\Sigma$, we say that $\Gamma$ proves $\varphi$ according to $\mathcal{A}^{\top}$, and write $\Gamma\vDash_{\mathcal{A}^{\top}}\varphi$, if for every homomorphism of $\Sigma$-multialgebras $\nu:\textbf{F}(\Sigma, \mathcal{V})\rightarrow\mathcal{A}^{\top}$,
\[\nu(\Gamma)\subseteq \top_{\mathcal{A}^{\top}}\quad\text{implies}\quad \nu(\varphi)\in \top_{\mathcal{A}^{\top}},\]
what is possible given that $\Sigma\subseteq \Sigma^{\top}$ implies $\Gamma\cup\{\varphi\}\subseteq F(\Sigma, \mathcal{V})\subseteq F(\Sigma^{\top}, \mathcal{V})$.
Of course, such a simplification is not possible when dealing with matrices, unless the set of designated values is a singleton.

Given a class $\mathbb{M}$ of Nmatrices over the same signature $\Sigma$, we say $\Gamma\vDash_{\mathbb{M}}\varphi$ if, for every $\mathcal{M}\in \mathbb{M}$, $\Gamma\vDash_{\mathcal{M}}\varphi$. Given every matrix is an Nmatrix, we have that every tarskian logic may be characterized by a class of Nmatrices.

\begin{definition}
A restricted non-deterministic matrix, or restricted Nmatrix or RNmatrix\index{RNmatrix}\index{Matrix, Restricted non-deterministic}, over a signature $\Sigma$ is a triple 
\[\mathcal{M}=(\mathcal{A}, D, \mathcal{F})\]
such that:
\begin{enumerate}
\item $(\mathcal{A}, D)$ is a non-deterministic matrix over $\Sigma$;
\item $\mathcal{F}$ is a subset of the set of all homomorphisms from $\textbf{F}(\Sigma, \mathcal{V})$ to $\mathcal{A}$.
\end{enumerate}
\end{definition}

We define $\vDash_{\mathcal{M}}$ as in the case that $\mathcal{M}$ is a restricted matrix. Notice that every restricted matrix is a restricted Nmatrix, and then every tarskian logic may be characterized by a restricted Nmatrix, but we can achieve a more powerful result, similar to the one commonly known as Suszko`s Thesis, found in \cite{Suszko}; but, unlike Suszko, we do not focus on bivaluations, and we do not wish to advocate that all logics are two-valued.

\begin{theorem}\label{2-valued RNmatrix}
Every tarskian logic is characterizable by a two-valued RNmatrix.\footnote{By an $n$-valued RNmatrix we understand an RNmatrix $\mathcal{M}=(\mathcal{A}, D, \mathcal{F})$ where the universe of $\mathcal{A}$ has $n$ elements.}
\end{theorem}

\begin{proof}
Let $\mathfrak{L}=(\mathcal{L}, \vdash)$ be a tarskian logic over the signature $\Sigma$. Consider then the $\Sigma$-multialgebra $\textbf{2}(\Sigma)$\label{2Sigma} with universe $\{0,1\}$ and, for an $n$-ary $\sigma\in\Sigma$, operations defined by
\[\sigma_{\textbf{2}(\Sigma)}(x_{1}, \dotsc  , x_{n})=\{0,1\}, \forall x_{1}, \dotsc  , x_{n}\in \{0,1\}.\]
We then define the set $\mathcal{F}_{\mathfrak{L}}$ of valuations $\nu:\textbf{F}(\Sigma, \mathcal{V})\rightarrow\textbf{2}(\Sigma)$ such that there exists a closed set of formulas $\Gamma$ over $\Sigma$ for which $\nu(\gamma)=1$ if, and only if, $\gamma\in\Gamma$; notice that all functions from $\textbf{F}(\Sigma, \mathcal{V})$ to $\textbf{2}(\Sigma)$ are homomorphisms, and therefore no further restrictions are necessary.

We then define the RNmatrix 
\[\textbf{2}(\mathfrak{L})=(\textbf{2}(\Sigma), \{1\}, \mathcal{F}_{\mathfrak{L}}),\]
\label{2L}and we state that $\Gamma\vdash\varphi$ if, and only if, $\Gamma\vDash_{\textbf{2}(\mathfrak{L})}\varphi$. 

Regarding the first direction, if $\Gamma\vdash\varphi$, for any valuation $\nu\in\mathcal{F}_{\mathfrak{L}}$ such that $\nu(\Gamma)\subseteq\{1\}$ there must exist a closed set $\Delta$ of formulas over $\Sigma$ such that $\nu(\delta)=1$ if, and only if, $\delta\in\Delta$, and since $\nu(\Gamma)\subseteq\{1\}$ we have $\Gamma\subseteq \Delta$; now, given $\Gamma\vdash\varphi$ and $\Gamma\subseteq \Delta$, it follows that $\Delta\vdash\varphi$ and, since $\Delta$ is closed, $\varphi\in \Delta$, meaning $\nu(\varphi)=1$ and, therefore, $\Gamma\vDash_{\textbf{2}(\mathfrak{L})}\varphi$.

For the second direction, assume $\Gamma\vDash_{\textbf{2}(\mathfrak{L})}\varphi$: if $\Gamma\not\vdash\varphi$, by Lindenbaum-\L oz there exists a non-trivial extension $\Delta$ of $\Gamma$ such that $\varphi\notin\Delta$; if $\nu$ is the valuation of $\mathcal{F}_{\mathfrak{L}}$ such that $\nu(\delta)=1$ if, and only if, $\delta\in\Delta$, this means that $\nu(\Delta)\subseteq\{1\}$, and therefore $\nu(\Gamma)\subseteq\{1\}$ since $\Gamma\subseteq\Delta$, and $\nu(\varphi)=0$, contradicting the fact that $\Gamma\vDash_{\textbf{2}(\mathfrak{L})}\varphi$. This proves $\Gamma\vDash_{\textbf{2}(\mathfrak{L})}\varphi$ implies $\Gamma\vdash\varphi$, and the theorem is proved.
\end{proof}

For a restricted Nmatrix $\mathcal{M}$ we will want to consider, for a subset $\Gamma$ of $F(\Sigma, \mathcal{V})$, the closure $K_{\mathcal{M}}(\Gamma)$\label{KMGamma}, that is, the set of formulas $\varphi\in F(\Sigma, \mathcal{V})$ such that $\Gamma\vDash_{\mathcal{M}}\varphi$.

\begin{definition}
Given a signature $\Sigma$, an operator $K:\mathcal{P}(F(\Sigma, \mathcal{V}))\rightarrow \mathcal{P}(F(\Sigma, \mathcal{V}))$ is said to be a tarskian\index{Operator, Tarskian} (or closure) one if, for all $\Gamma, \Theta\subseteq F(\Sigma, \mathcal{V})$, it satisfies:
\begin{enumerate}
\item $\Gamma\subseteq K(\Gamma)$;
\item if $\Theta\subseteq \Gamma$, $K(\Theta)\subseteq K(\Gamma)$;
\item if $\Theta=K(\Gamma)$, $K(\Theta)=\Theta$.
\end{enumerate}
\end{definition}

Notice how this generalizes the notion of a logic being tarskian: the consequence relation in a tarskian logic is a tarskian operator; for such a reason, one may often call the pair $(\textbf{F}(\Sigma, \mathcal{V}), K)$, for $K$ a tarskian operator on $F(\Sigma, \mathcal{V})$, a sentential logic itself. Even more, we may define a logic $\mathcal{L}$ over the signature $\Sigma$ by $\Gamma\vdash_{\mathcal{L}}\varphi$ if, and only if, $\varphi\in K(\Gamma)$, and it is clear how $\mathcal{L}$ is tarskian: to every tarskian operator there corresponds a tarskian logic and vice-versa.

\begin{proposition}
Given a restricted Nmatrix $\mathcal{M}=(\mathcal{A}, D, \mathcal{F})$, the operator $\Gamma\mapsto K_{\mathcal{M}}(\Gamma)$ is a tarskian one.
\end{proposition}

\begin{proof}
\begin{enumerate}
\item Take $\varphi\in \Gamma$, and suppose $\nu\in\mathcal{F}$ is such that $\nu(\gamma)\in D$ for all $\gamma\in\Gamma$; then $\nu(\varphi)\in D$ and therefore $\varphi\in K_{\mathcal{M}}(\Gamma)$, so that $\Gamma\subseteq K_{\mathcal{M}}(\Gamma)$.

\item Suppose $\Theta\subseteq\Gamma$ and that $\varphi\in K_{\mathcal{M}}(\Theta)$: then, if $\nu\in\mathcal{F}$ satisfies that $\nu(\theta)\in D$, for every $\theta\in\Theta$, one has $\nu(\varphi)\in D$.

Now, take a restricted homomorphism $\nu\in\mathcal{F}$ such that $\nu(\gamma)\in D$, for every $\gamma\in\Gamma$: since $\Theta\subseteq \Gamma$, in this case $\nu(\theta)\in D$, $\forall \theta\in \Theta$, and therefore $\nu(\varphi)\in D$, implying that $\varphi\in K_{\mathcal{M}}(\Gamma)$. It follows that $K_{\mathcal{M}}(\Theta)\subseteq K_{\mathcal{M}}(\Gamma)$.

\item By the first item above, $\Theta\subseteq K_{\mathcal{M}}(\Theta)$, so it only remains to be shown that $K_{\mathcal{M}}(\Theta)\subseteq \Theta$. Given one $\varphi\in K_{\mathcal{M}}(\Theta)$, for any $\nu\in\mathcal{F}$ such that $\nu(\theta)\in D$, for all $\theta\in\Theta$, $\nu(\varphi)\in D$.

Now suppose $\nu\in\mathcal{F}$ satisfies that $\nu(\gamma)\in D$, for every $\gamma\in\Gamma$: then, since $\Theta=K_{\mathcal{M}}(\Gamma)$, $\nu(\theta)\in D$ for all $\theta\in\Theta$, and therefore $\nu(\varphi)\in D$. It follows that $\varphi\in K_{\mathcal{M}}(\Gamma)=\Theta$ and $K_{\mathcal{M}}(\Theta)=\Theta$.
\end{enumerate}
\end{proof}

Given a class $\mathbb{M}$ of RNmatrices, we can then consider the closure of a set $\Gamma\subseteq F(\Sigma, X)$ under $\mathbb{M}$, that is, $\varphi\in K_{\mathbb{M}}(\Gamma)$\label{KMMGamma} if and only if $\Gamma\vDash_{\mathbb{M}}\varphi$. Clearly 
\[K_{\mathbb{M}}(\Gamma)=\bigcap_{\mathcal{M}\in\mathbb{M}}K_{\mathcal{M}}(\Gamma).\]

\begin{lemma}
If $K_{\lambda}$ is a tarskian operator for every $\lambda\in\Lambda$, $K$ defined as 
\[K(\Gamma)=\bigcap_{\lambda\in\Lambda}K_{\lambda}(\Gamma)\]
is also tarskian.
\end{lemma}

\begin{proof}
\begin{enumerate}
\item Since, for every $\lambda\in\Lambda$ we have that $\Gamma\subseteq K_{\lambda}(\Gamma)$, we have that $\Gamma\subseteq \bigcap_{\lambda\in\Lambda}K_{\lambda}(\Gamma)=K(\Gamma)$.

\item If $\Theta\subseteq \Gamma$, we have $K_{\lambda}(\Theta)\subseteq K_{\lambda}(\Gamma)$ for every $\lambda\in\Lambda$, so 
\[K(\Theta)=\bigcap_{\lambda\in\Lambda}K_{\lambda}(\Theta)\subseteq\bigcap_{\lambda\in\Lambda}K_{\lambda}(\Gamma)=K(\Gamma).\]

\item Finally, if $\Theta=K(\Gamma)$, then $\bigcap_{\lambda\in\Lambda}K_{\lambda}(\Gamma)=\Theta$, so that $\Theta\subseteq K_{\lambda}(\Gamma)$ for every $\lambda\in\Lambda$: clearly $\Theta\subseteq K(\Theta)$, so it remains for us to show that $K(\Theta)\subseteq \Theta$.

Let us denote $K_{\lambda}(\Gamma)$ by $\Theta_{\lambda}$, and since $\Theta\subseteq K_{\lambda}(\Gamma)$, $K_{\lambda}(\Theta)\subseteq K_{\lambda}(\Theta_{\lambda})=\Theta_{\lambda}$, for every $\lambda\in\Lambda$; notice that $\Theta=\bigcap_{\lambda\in\Lambda}\Theta_{\lambda}$.

Then, 
\[K(\Theta)=\bigcap_{\lambda\in\Lambda}K_{\lambda}(\Theta)\subseteq \bigcap_{\lambda\in\Lambda}\Theta_{\lambda}=\Theta,\]
what finishes the proof.
\end{enumerate}
\end{proof}

\begin{theorem}
Given a class $\mathbb{M}$ of RNmatrices, $K_{\mathbb{M}}$ is a tarskian operator.
\end{theorem}

In other words, classes of restricted Nmatrices can, at most, describe tarskian logics, and from Theorem \ref{2-valued RNmatrix} they indeed characterize all of these. So no non-tarskian logic can be characterized by either RNmatrices or their classes,\footnote{Although slight generalizations of RNmatrices could possibly change this.} but the fact is we do not see that as a real problem: we are not looking exclusively for expressive power in our semantics, being efficiency a more desirable property. That is, the point of restricted Nmatrices will not be what they can express, but how easily will be to define and use them. Most importantly, however, is that RNmatrices will be able to provide decision methods through what are, essentially, truth-tables, where none were available before: one example, we will stress repeatedly, is that of $C_{1}$; although decidable, this system is not only not characterizable by finite matrices, but neither by finite sets of finite matrices, finite Nmatrices or finite sets of finite Nmatrices (\cite{Avron}). And, despite all of this, this logic, and in fact all of da Costa's hierarchy, may be characterized by finite RNmatrices, what we will achieve in Chapter \ref{Chapter5}.

\subsection{A brief history of RNmatrices}\label{A brief history of RNmatrices}

We first developed RNmatrices by studying semantics for the logic $\textbf{mbC}$: when analyzing Fidel structures for that very system, we hoped to simplify the ($3$-valued) Nmatrix for $\textbf{mbC}$ to a matrix, by slightly altering its signature. By replacing the unary connective, standing for consistency, for a binary connective ($\uparrow$) for incompatibility, and therefore creating the logics of incompatibility we study in Chapters \ref{Chapter7}, \ref{Chapter8} and \ref{Chapter9}, we created a system (within $\textbf{nbI}^{-}$, that is, $\textbf{nbI}$ without the commutativity of $\uparrow$) equivalent to $\textbf{mbC}$, with simpler semantics, but which was also, unfortunately, not characterizable by finite matrices.

Well, a very natural system related to $\textbf{nbI}^{-}$ is $\textbf{bI}$, which adds the commutativity of the binary incompatibility connective and subtracts the negation; rather unfortunate was then our discovery that this logic is not characterizable by, not only finite matrices, but rather finite Nmatrices as well. But we had a semantic for it almost ready, simply by adapting the decision method (a two-valued Nmatrix) for $\textbf{bI}^{-}$: it was enough to demand that every valuation to be taken into consideration should satisfy $\nu(\alpha\uparrow\beta)=\nu(\beta\uparrow\alpha)$. 

Of course, within a semantics of matrices, or even Nmatrices, this is not permissible; but, inspired by Piochi's work on $\mathcal{E}$-matrices (which do restrict valuations) and our previous knowledge of Nmatrices, we have coined what we begun to call restricted non-deterministic matrices, or RNmatrices. We quickly grew to realize that RNmatrices are actually recurrent in the history of non-classical logics, although not defined as such but rather as a specific solution to a difficult system: examples abound, including \index{Bivaluation}bivaluations, Fidel structures, Kearn's $4$-valued semantics for modal logics, static semantics, PNmatrices, among others.

To our great surprise, we found out that RNmatrices semantics have quite recently started to resurface: Pawlowski and Urbaniak, in \cite{Pawlowski, PawlowskUrbaniak}, had noticed, shortly before us, previous uses of the semantics of RNmatrces (although their nomenclature is, naturally, different), specially in the areas of informal provability and modal logic. We add to the topic a more comprehensive analysis of the literature and a more in depth study of the expressiveness of these semantics, in addition to the first explicit use of RNmatrices to paraconsistent logics.

And it is important, here, to explain what we mean by the expressiveness of RNmatrix semantics: usually, the expressiveness of semantics refers to which logics can these semantics characterize, a concept that is distinctively easy to understand when our semantics are matricial in nature. However, according to Theorem \ref{2-valued RNmatrix}, RNmatrices are, in a way, as expressive as possible, meaning that any tarskian logic can be characterized by a finite RNmatrix. Hence, we are actually more concerned about, first of all, the expressiveness of decidable RNmatrices, \textit{i. e.} RNmatrices where the problem of identifying which valuations are restricted valuations is decidable, see Section \ref{Decision methods} of Chapter \ref{Chapter5} for more details; second, the expressiveness of other, related semantics, such as PNmatrices of Example \ref{example of PNMatrices} below, which we know to be no more expressive than RNmatrices (that is, all logics characterized by PNmatrices can also be characterized by RNmatrices), but are not sure whether they are as expressive as RNmatrices (meaning, if PNmatrices characterize all tarskian logics, as RNmatrices do).

\begin{example}\label{Example Kearns}
As we have already discussed, Dugundji proved that the logics between the modal systems $\textbf{S1}$ and $\textbf{S5}$ can not be characterized by single, finite matrices; to circumvent that, J. Kearns\index{Kearns} (\cite{Kearns}) proposed a (four-valued) Nmatrix semantics for the modal $\textbf{T}$, $\textbf{S4}$ and $\textbf{S5}$ for which just some specific valuations can be considered, what clearly characterizes these as RNmatrices semantics. Kearns idea proved itself to be very popular, and generalizations can be found in, \textit{exempli gratia}, \cite{CCP} (and \cite{CCPErr}) and \cite{OmoriSkurt}.

Using a more modern approach to Kearns technique, as in \cite{CCP}, we proceed as follows:\label{Semantics, Kearns} consider the universe $\{T, t, f, F\}$; operations are defined as necessary for each logic $\textbf{L}\in\{\textbf{T}, \textbf{S4}, \textbf{S5}\}$. Consider the set $Val^{\textbf{L}}$ of all valuations for $\textbf{L}$; we define $Val_{k}^{\textbf{L}}\subseteq Val^{\textbf{L}}$ by induction as $Val_{0}^{\textbf{L}}=Val^{\textbf{L}}$ and, for $k\in\mathbb{N}$,\[Val_{k+1}^{\textbf{L}}=\{\nu\in Val_{k}^{\textbf{L}} : \text{for every formula $\alpha$, $Val_{k}^{\textbf{L}}(\alpha)\subseteq D$ implies $\nu(\alpha)=T$}\},\]
where $Val_{k}^{\textbf{L}}(\alpha)=\{\mu(\alpha) : \mu\in Val_{k}^{\textbf{L}}\}$. We then define the restricted valuations as $\mathcal{F}_{\textbf{L}}=\bigcap_{k\in\mathbb{N}}Val_{k}^{\textbf{L}}$, and Kearns proved that $\vdash_{\textbf{L}}\alpha$ if, and only if, $\nu(\alpha)=T$ for every $\nu\in\mathcal{F}_{\textbf{L}}$, making his semantics for $\textbf{L}$ correspond to the RNmatrix $\mathcal{K}_{\textbf{L}}=(\mathcal{A}_{\textbf{L}}, \{T\}, \mathcal{F}_{\textbf{L}})$, for $\mathcal{A}_{\textbf{L}}$ the multialgebra for $\textbf{L}$, with universe $\{T, t, f, F\}$, whose precise definition we omitted.

Now, we can also prove $\mathcal{K}_{\textbf{L}}$\label{KL} is structural: through an induction on $k$, one shows that, for every valuation $\nu\in Val_{k}^{\textbf{L}}$ and every substitution $\rho$, $\nu\circ\rho\in Val_{k}^{\textbf{L}}$, being the case $k=0$ trivially true; so, suppose the result holds for a certain $k\in\mathbb{N}$. Given $\nu\in Val_{k+1}^{\textbf{L}}$ and a substitution $\rho$, by definition of $Val_{k+1}^{\textbf{L}}$ we have $\nu\in Val_{k}^{\textbf{L}}$ and, by induction hypothesis, $\nu\circ\rho\in Val_{k}^{\textbf{L}}$. Let $\alpha$ be a formula satisfying $Val_{k}^{\textbf{L}}(\alpha)\subseteq D$: for every $\mu\in Val_{k}^{\textbf{L}}$, $\mu\circ\rho\in Val_{k}^{\textbf{L}}$ (again by induction hypothesis), meaning that $\mu\circ\rho(\alpha)=\mu(\rho(\alpha))\in D$, and therefore $Val_{k}^{\textbf{L}}(\rho(\alpha))\subseteq D$.

This implies, since $\nu\in Val_{k+1}^{\textbf{L}}$, that $\nu\circ\rho(\alpha)=\nu(\rho(\alpha))=T$, and henceforth $\nu\circ\rho\in Val_{k+1}^{\textbf{L}}$. Obviously this proves that $\nu\in Val_{k}^{\textbf{L}}$ implies $\nu\circ\rho\in Val_{k}^{\textbf{L}}$ for all $k\in\mathbb{N}$: so, if $\nu\in \mathcal{F}_{\textbf{L}}$ and $\rho$ is a substitution, for every $k\in \mathbb{N}$, $\nu\in Val_{k}^{\textbf{L}}$ implies $\nu\circ\rho\in Val_{k}^{\textbf{L}}$ (again, for every $k\in\mathbb{N}$), and so $\nu\circ\rho\in \bigcap_{k\in\mathbb{N}}Val_{k}^{\textbf{L}}=\mathcal{F}_{\textbf{L}}$.
\end{example}

\begin{example}
The first proof of the decidability of da Costa's calculi $C_{n}$ (see Chapter \ref{Chapter5} for a definition of these logics) is due to Fidel, which in 1977 created a new class of structures, both algebraic and relational, to provide precisely such a demonstration in \cite{Fidel3}: these objects are now known as \textit{Fidel structures}\index{Fidel structure}. Essentially, a \index{Fidel structure for $C_{n}$}Fidel structure for $C_{n}$, which in this presentation are equipped with a connective $(n)$, is a triple $\mathcal{N}=(\mathcal{A}, \{N_{a}\}_{a\in A}, \{N_{a}^{(n)}\}_{a\in A})$ such that:
\begin{enumerate}
\item $\mathcal{A}$ is a Boolean algebra with universe $A$;
\item for every $a\in A$, $N_{a}$\label{Na} and $N_{a}^{(n)}$ are subsets of $A$ (and therefore unary relations indexed by $A$) satisfying certain desired properties.
\end{enumerate}
A valuation over a Fidel structure $\mathcal{N}$ for $C_{n}$ is a function $\nu$ from the formulas of $C_{n}$ to $A$ satisfying, among other conditions, that $\nu(\alpha\#\beta)=\nu(\alpha)\#\nu(\beta)$, for $\#\in\{\vee, \wedge, \rightarrow\}$, 
\[\nu(\neg\alpha)\in N_{\nu(\alpha)}\quad\text{and}\quad\nu(\alpha^{(n)})\in N_{\nu(\alpha)}^{(n)}.\]
If, for every element $a$ of $\mathcal{A}$, we define $\neg a=N_{a}$ and $n(a)=N_{a}^{(n)}$, we have enriched $\mathcal{A}$ to become a multialgebra, which we denote by $\mathcal{A}_{\mathcal{N}}^{n}$; it is not difficult to see, then, that the consequence operator induced by the Fidel structure $\mathcal{N}$ equals the one produced by the RNmatrix $\mathcal{M}_{\mathcal{N}}^{n}=(\mathcal{A}_{\mathcal{N}}^{n}, D, \mathcal{F}_{\mathcal{N}}^{n})$, for $D=\{1\}$ and $\mathcal{F}_{\mathcal{N}}^{n}$ the set of valuations for $C_{n}$ over $\mathcal{N}$.

Thus, we understand Fidel structures as one of the earliest, and best understood, applications of RNmatrices; more important, however, is its success in surviving the test of time, being a methodology that is still relevantly applied today. See, for an example, the characterizations of $\textbf{LFI}$'s such as $\textbf{mbC}$, $\textbf{mbCcl}$ and $\textbf{CILA}$ in \cite{ParLog} through Fidel structures, all of which may be recast as RNmatrices, one may add.

Now, to be more precise, \cite{Fidel3} uses a Fidel structure called $\textbf{C}$ over the two-element Boolean algebra $\textbf{2}$ as a decision procedure for all of $C_{n}$: that way, for any $n\geq 1$, the multialgebras $\mathcal{A}_{\textbf{C}}^{n}$ are the same and satisfy $\neg 0=n(0)=\{1\}$ and $\neg1=n(1)=\{0,1\}$; this forces most of the semantical power of these structures to lie in their valuations, making of $\mathcal{F}_{\textbf{C}}^{n}$ rather complicated sets (look, for an example of a similar nature, to $\textbf{2}(\mathfrak{L})$ in Theorem \ref{2-valued RNmatrix}). In a certain sense, when we present our own RNmatrices for $C_{n}$ in Chapter \ref{Chapter5} we make a compromise: we allow for the underlying multialgebras to have larger universes (with $n+2$ elements for $C_{n}$), but in return we have far simpler valuations.
\end{example}

\begin{example}
In the work of Avron and Konikowska \cite{AK:05}, valuations over Nmatrices induce semantics which they called dynamic semantics\index{Semantics, Dynamic} over Nmatrices; as opposed to these, they have also considered semantics which restrict the usual valuations, what they have baptized as static semantics\index{Semantics, Static}. Essentially, given an Nmatrix $\mathcal{M}=(\mathcal{A}, D)$, its static semantics is given by $\mathcal{M}^{s}=(\mathcal{A}, D, \mathcal{F}^{s}_{\mathcal{M}})$\label{Ms}, for $\mathcal{F}^{s}_{\mathcal{M}}$ the set of valuations $\nu$ over $\mathcal{M}$ satisfying that: for an $n$-ary connective $\sigma$, and formulas $\sigma(\alpha_{1}, \dotsc  , \alpha_{n})$ and $\sigma(\beta_{1}, \dotsc  , \beta_{n})$,  $\nu(\alpha_{i})=\nu(\beta_{i})$, for every $i\in\{1, \dotsc  , n\}$, implies
\[\nu(\sigma(\alpha_{1}, \dotsc  , \alpha_{n}))=\nu(\sigma(\beta_{1}, \dotsc  , \beta_{n})).\]
It is easy to prove that the RNmatrices $\mathcal{M}^{s}$ are all structural.
\end{example}

\begin{example}\label{example of PNMatrices}
PNmatrices\index{PNmatrix}\index{Matrix, Partial non-deterministic}, whose first appearance in the literature was in \cite{Baaz}, generalize Nmatrices by replacing multialgebras for partial multialgebras (that is, algebras of relations). Our brief presentation here, however, is closer to that of \cite{CalMar}: a partial $\Sigma$-multialgebra is a pair $\mathcal{A}=(A, \{\sigma_{\mathcal{A}}\}_{\sigma\in\Sigma})$ such that, if $\sigma\in\Sigma_{n}$, $\sigma_{\mathcal{A}}$ is a function of the form $\sigma_{\mathcal{A}}:A^{n}\rightarrow\mathcal{P}(A)$ (no longer $\mathcal{P}(A)\setminus\{\emptyset\}$ as in the case of multialgebras). A valuation for a partial $\Sigma$-multialgebra $\mathcal{A}$, as described, is then a map $\nu: F(\Sigma, \mathcal{V})\rightarrow A$ such that
\[\nu(\sigma(\alpha_{1}, \dotsc  , \alpha_{n}))\in\sigma_{\mathcal{A}}(\nu(\alpha_{1}), \dotsc  , \nu(\alpha_{n})),\]
for every $\sigma\in\Sigma_{n}$ and formulas $\alpha_{1}, \dotsc  , \alpha_{n}$; notice that, if $\sigma_{\mathcal{A}}(a_{1}, \dotsc  , a_{n})=\emptyset$, then there cannot exist a valuation $\nu$ and formulas $\alpha_{1}$ trough $\alpha_{n}$ such that $\nu(\alpha_{i})=a_{i}$, for $1\leq i\leq n$, given $\nu$ is supposed to be a total function. Finally, given a partial $\Sigma$-multialgebra $\mathcal{A}$ and a subset $D$ of its universe, $\mathcal{M}=(\mathcal{A}, D)$ is a PNmatrix: given formulas $\Gamma\cup\{\varphi\}$ on $\Sigma$, we say $\Gamma$ proves $\varphi$ according to $\mathcal{M}$, and write $\Gamma\vDash_{\mathcal{M}}\varphi$, whenever, for every valuation $\nu$ for $\mathcal{A}$, $\nu(\Gamma)\subseteq D$ implies $\nu(\varphi)\in D$.

What we can prove, however, is that the semantics induced by PNmatrices can be also induced by RNmatrices, making of the latter a more expressive technique. Given a PNmatrix $\mathcal{M}=(\mathcal{A}, D)$, over $\Sigma$, we consider the $\Sigma$-multialgebra $\mathcal{A}^{\emptyset}=(A\cup\{o\}, \{\sigma_{\mathcal{A}^{\emptyset}}\}_{\sigma\in\Sigma})$, for an element $o\notin A$, such that
\[\sigma_{\mathcal{A}^{\emptyset}}(a_{1}, \dotsc  , a_{n})=
\begin{cases*}
\sigma_{\mathcal{A}}(a_{1}, \dotsc  , a_{n}) & if $a_{1}, \dotsc  , a_{n}\in A$ and $\sigma_{\mathcal{A}}(a_{1}, \dotsc  , a_{n})\neq\emptyset$,\\
\{o\} & otherwise,
\end{cases*}\]
for all $a_{1}, \dotsc  , a_{n}\in A\cup\{o\}$. By defining
\[\mathcal{F}=\{\nu:\textbf{F}(\Sigma, \mathcal{V})\rightarrow\mathcal{A}^{\emptyset}\ :\ o\notin\nu(F(\Sigma, \mathcal{V}))\},\]
where $\nu(F(\Sigma, \mathcal{V}))=\{\nu(\alpha) : \alpha\in F(\Sigma, \mathcal{V})\}$, one sees that the valuations for $\mathcal{A}$ correspond, more or less, to $\mathcal{F}$, that is, the set of valuations for $\mathcal{A}^{\emptyset}$ without $o$ in their range; then $\mathcal{M}^{\emptyset}=(\mathcal{A}^{\emptyset}, D, \mathcal{F})$\label{Mo} has the same deduction operator as $\mathcal{A}$.
\end{example}

\begin{example}\label{PTS}
Possible-translations semantics were first defined in \cite{PTS-first}, but here we use an approach similar to that found in \cite{PTS-defined}, the difference being that we only focus on zeroth order logics: take signatures $\Sigma$ and $\Sigma^{*}$ and its respective languages, $\mathcal{L}_{\Sigma}$ and $\mathcal{L}_{\Sigma^{*}}$; a translation between logics $\mathscr{L}=(\mathcal{L}_{\Sigma}, \vdash_{\mathscr{L}})$ and $\mathscr{L}^{*}=(\mathcal{L}_{\Sigma^{*}}, \vdash_{\mathscr{L}^{*}})$ is any function $t:\mathcal{L}_{\Sigma}\rightarrow \mathcal{L}_{\Sigma^{*}}$ such that, for any set of formulas $\Gamma\cup\{\varphi\}$ of $\mathcal{L}_{\Sigma}$,
\[\text{$\Gamma\vdash_{\mathscr{L}}\varphi$\quad implies\quad $t(\Gamma)\vdash_{\mathscr{L}^{*}}t(\varphi)$,}\]
where $t(\Gamma)=\{t(\gamma)\ :\ \gamma\in\Gamma\}$. Fixed a system $\mathscr{L}$ over $\Sigma$, a possible-translations semantics\index{Possible-translations semantics} for $\mathscr{L}$ is then a pair $\mathcal{PT}=(\textbf{Log}, \textbf{Tr})$ where $\textbf{Log}=\{\mathscr{L}_{i}\}_{i\in I}$ are logics over possibly distinct signatures $\Sigma^{i}$, and $\textbf{Tr}=\{t_{i}\}_{i\in I}$ are translations from $\mathscr{L}$ to $\mathscr{L}_{i}$. The possible-translations semantics leads to a consequence operator $\Vdash_{\mathcal{PT}}$ on the formulas of the signature $\Sigma$ whenever we define, for a set of formulas $\Gamma\cup\{\varphi\}$ over $\Sigma$,
\[\text{$\Gamma\Vdash_{\mathcal{PT}}\varphi$\quad if and only if\quad $t_{i}(\Gamma_{i})\vdash_{\mathscr{L}_{i}}t_{i}(\varphi)$,\quad for all $i\in I$.}\]
Of course, $\mathscr{L}$ is characterized by $\mathcal{PT}$ when $\Gamma\vdash_{\mathscr{L}}\varphi$ iff $\Gamma\Vdash_{\mathcal{PT}}\varphi$; intuitively, a logic is characterized by a given possible-translations semantics whenever $\mathscr{L}$ can be defined by specific properties of the $\mathscr{L}_{i}$, which the translations $t_{i}$ combine together to obtain $\Vdash_{\mathcal{PT}}$. 

Of course, RNmatrices cannot possibly be as expressive as even this weaker definition of possible-translations semantics, given that they can only deal with one signature at a time,\footnote{It seems, however, possible to combine RNmatrices over different signatures as long as one has at their disposal adequate maps between signatures, sometimes also known as translations.} but a connection between the two semantics is self-evident: both of them shift the expressive power of semantics, from structures to maps themselves; specifically, translations in the case of possible-translations semantics, and restricted valuations in the case of RNmatrices.

However, consider a specific case of possible-translations semantics: suppose all logics of $\textbf{Log}$ have the same signature $\Sigma$ as $\mathscr{L}$, and that all of them are characterized by RNmatrices $\mathcal{M}_{i}$ or, what is equivalent, that all of them are tarskian (and, therefore, so is $\mathscr{L}$); then $\mathscr{L}$ is characterized by a class of of RNmatrices in a straightforward way. In fact, if $\mathcal{M}_{i}=(\mathcal{A}_{i}, D_{i}, \mathcal{F}_{i})$, then 
\[\text{$\Gamma\vdash_{\mathscr{L}_{i}}\varphi$\quad iff \quad $\nu(\Gamma)\subseteq D_{i}$\quad implies\quad $\nu(\varphi)\in D_{i}$\quad for every $\nu\in\mathcal{F}_{i}$;}\]
by defining $\mathcal{F}_{i}^{t}=\{\nu\circ t_{i}\ :\ \nu\in\mathcal{F}_{i}\}$ and $\mathcal{M}_{i}^{t}=(\mathcal{A}_{i}, D_{i}, \mathcal{F}_{i}^{t})$, it is easy to see that $t_{i}(\Gamma)\vdash_{\mathscr{L}_{i}}t_{i}(\varphi)$ if, and only if, $\Gamma\vDash_{\mathcal{M}_{i}^{t}}\varphi$. By making $\mathbb{M}=\{\mathcal{M}_{i}^{t}\}_{i\in I}$, we see that $\Gamma\Vdash_{\mathcal{PT}}\varphi$ is equivalent to $\Gamma\vDash_{\mathbb{M}}\varphi$, and so $\mathscr{L}$ is characterized by $\mathbb{M}$.

\end{example}

\subsection{Structurality}\label{Structurality}

\begin{definition}
For any subsemigroup $\mathcal{E}$ of the semigroup $\textit{End}(F(\Sigma, \mathcal{V}))$ of endomorphisms on $F(\Sigma, \mathcal{V})$, we say an operator $K:\mathcal{P}(F(\Sigma, \mathcal{V}))\rightarrow\mathcal{P}(F(\Sigma, \mathcal{V}))$ is $\mathcal{E}$-structural\index{Structural, $\mathcal{E}$-} when 
\[\{\sigma(\varphi)\ :\  \varphi\in K(\Gamma)\}\subseteq K(\{\sigma(\varphi)\ :\  \varphi\in\Gamma\})\]
for every $\sigma\in\mathcal{E}$, or as we shall write it, $\sigma K(\Gamma)\subseteq K(\sigma\Gamma)$; the operator is said to be structural when it is $\textit{End}(F(\Sigma, \mathcal{V}))$-structural. 
\end{definition}

Structurality is rather important when dealing with logics arising from a Hilbert calculus given that instances of axiom schemata and rules of inference are obtained by application of endomorphisms; so, when studying a deduction operator, it is important to understand how structural it is. While most of the systems we study here are structural, not only $\mathcal{E}$-structural for a proper subsemigroup $\mathcal{E}$, developing a general theory of structurality for RNmatrices, as Piochi did for Rmatrices in \cite{Piochi}, \cite{Piochi2} and \cite{Piochi3}, can have useful applications on non-structural logics.

\begin{proposition}
If $\mathcal{M}=(\mathcal{A}, D, \mathcal{F})$ is an RNmatrix such that, for every $\nu\in\mathcal{F}$ and $\sigma\in\mathcal{E}$, $\nu\circ\sigma\in\mathcal{F}$, then $K_{\mathcal{M}}$ is $\mathcal{E}$-structural.
\end{proposition}

\begin{proof}
Let $\Gamma$ be a subset of $F(\Sigma, \mathcal{V})$ and $\varphi\in K_{\mathcal{M}}(\Gamma)$: we must prove that, for a given $\sigma\in \mathcal{E}$, $\sigma(\varphi)$ is in $K_{\mathcal{M}}(\sigma\Gamma)$, or what is equivalent, that $\sigma\Gamma\vDash_{\mathcal{M}}\sigma(\varphi)$. So, let $\nu\in\mathcal{F}$ be a valuation satisfying $\nu(\sigma(\gamma))=\nu\circ\sigma(\gamma)\in D$, $\forall\gamma\in\Gamma$, and we must prove that $\nu(\sigma(\varphi))=\nu\circ\sigma(\varphi)\in D$.

By our hypothesis, $\nu\circ\sigma\in\mathcal{F}$, and since $\nu\circ\sigma(\gamma)\in D$ for every $\gamma\in\Gamma$ and $\varphi\in K_{\mathcal{M}}(\Gamma)$ (meaning $\Gamma\vDash_{\mathcal{M}}\varphi$), we find that $\nu\circ\sigma(\varphi)\in D$, what ends the proof.
\end{proof}

\begin{lemma}
If all operators $K_{\lambda}$, for $\lambda\in\Lambda$, are $\mathcal{E}$-structural, then so it is $K$ defined by 
\[K(\Gamma)=\bigcap_{\lambda\in\Lambda}K_{\lambda}(\Gamma),\]
for every $\Gamma\subseteq F(\Sigma, \mathcal{V})$.
\end{lemma}

\begin{proof}
Take $\sigma\in\mathcal{E}$ and $\varphi\in K(\Gamma)$: we have that $\varphi\in K_{\lambda}(\Gamma)$ for every $\lambda\in\Lambda$, and then $\sigma(\varphi)\in K_{\lambda}(\sigma\Gamma)$, since each $K_{\lambda}$ is structural.

It follows that $\sigma(\varphi)\in\bigcap_{\lambda\in\Lambda}K_{\lambda}(\sigma\Gamma)=K(\sigma\Gamma)$, and so $\sigma K(\Gamma)\subseteq K(\sigma\Gamma)$.
\end{proof}

\begin{theorem}
Given a class $\mathbb{M}$ of RNmatrices $\mathcal{M}=(\mathcal{A}, D, \mathcal{F})$ such that, for every $\nu\in\mathcal{F}$ and $\sigma\in\mathcal{E}$, $\nu\circ\sigma\in\mathcal{F}$, $K_{\mathbb{M}}$ is a $\mathcal{E}$-structural operator.
\end{theorem}

Now, for any tarskian logic $\mathfrak{L}$, we remember for an instant the two-valued RNmatrix $\textbf{2}(\mathfrak{L})$ which characterizes $\mathfrak{L}$ from Theorem \ref{2-valued RNmatrix}: notice that, if $\mathfrak{L}$ is $\mathcal{E}$-structural, so is the deduction operator induced by $\textbf{2}(\mathfrak{L})$. To see that, take an endomorphism $\sigma$ in $\mathcal{E}$; if $\Gamma\vDash_{\textbf{2}(\mathfrak{L})}\varphi$, $\Gamma\vdash\varphi$, and from the fact that $\mathfrak{L}$ is $\mathcal{E}$-structural, $\{\sigma(\gamma)\ :\  \gamma\in\Gamma\}\vdash\sigma(\varphi)$; this means, of course, that 
\[\{\sigma(\gamma)\ :\ \gamma\in\Gamma\}\vDash_{\textbf{2}(\mathfrak{L})}\sigma(\varphi),\]
what proves the result. Of course, if $\mathfrak{L}$ is structural, so is the operator of $\textbf{2}(\mathfrak{L})$.

\section{Examples: characterizing some logics with RNmatrices}\label{Examples: characterizing some logics with RNmatrices}

We would like to know that our approach using restricted Nmatrices has something to offer that other approaches don't. It is a classical result by Avron, found in \cite{Avron}, that no finite Nmatrix can characterize the logics between $\mbCcl$ and $\CILA$, both included: the problem arises specifically from axiom $\cl$. This is an important Dugundji-like uncharacterizability theorem, that limits greatly the expressive power of Nmatrices. We will show, by offering examples of such finite RNmatrices, that no such restriction exists for finite RNmatrices over those logics.

\subsection{$\mbCcl$}\label{RNmatrix for mbCcl}

This logic, together with the closely related $\mbCci$, were first defined, under these names, in \cite{ParLog}; however, the system $\textbf{B}[\{\textbf{i1}, \textbf{i2}\}]$\label{Bi1i2}, equivalent to $\mbCci$, appeared already in \cite{Avron2}, while $\textbf{Bi}$\label{Bi} and $\textbf{Bl}$\label{Bl}, equivalent to respectively $\mbCci$ and $\mbCcl$, had also been defined in \cite{Avron}.

Consider the $\Sigma_{\textbf{LFI}}$-multialgebra $\mathcal{A}_{\mbCciw}$ with universe $\{F, t, T\}$ and operations given by the tables below.
\begin{figure}[H]
\centering
\begin{minipage}[t]{4cm}
\centering
\begin{tabular}{|l|c|c|r|}
\hline
$\vee$ & $F$ & $t$ & $T$ \\ \hline
$F$ & $\{F\}$ & $\{t, T\}$ & $\{t, T\}$ \\ \hline
$t$ & $\{t, T\}$ & $\{t, T\}$ & $\{t, T\}$ \\ \hline
$T$ & $\{t, T\}$ & $\{t, T\}$ & $\{t, T\}$ \\ \hline
\end{tabular}
\caption*{Table for Disjunction}
\end{minipage}
\hspace{3cm}
\centering
\begin{minipage}[t]{4cm}
\centering
\begin{tabular}{|l|c|c|r|}
\hline
$\wedge$ & $F$ & $t$ & $T$ \\ \hline
$F$ & $\{F\}$ & $\{F\}$ & $\{F\}$ \\ \hline
$t$ & $\{F\}$ & $\{t, T\}$ & $\{t, T\}$ \\ \hline
$T$ & $\{F\}$ & $\{t, T\}$ & $\{t, T\}$ \\ \hline
\end{tabular}
\caption*{Table for Conjunction}
\end{minipage}
\end{figure}

\begin{figure}[H]
\centering
\begin{minipage}[t]{3cm}
\centering
\begin{tabular}{|l|r|}
\hline
 & $\neg$ \\ \hline
$F$ & $\{t, T\}$\\ \hline
$t$ & $\{t, T\}$\\ \hline
$T$ & $\{F\}$ \\ \hline
\end{tabular}
\caption*{Table for $\neg$}
\end{minipage}
\begin{minipage}[t]{5cm}
\centering
\begin{tabular}{|l|c|c|r|}
\hline
$\rightarrow$ & $F$ & $t$ & $T$ \\ \hline
$F$ & $\{t, T\}$ & $\{t, T\}$ & $\{t, T\}$ \\ \hline
$t$ & $\{F\}$ & $\{t, T\}$ & $\{t, T\}$ \\ \hline
$T$ & $\{F\}$ & $\{t, T\}$ & $\{t, T\}$ \\ \hline
\end{tabular}
\caption*{Table for Implication}
\end{minipage}
\begin{minipage}[t]{3cm}
\centering
\begin{tabular}{|l|r|}
\hline
 & $\circ$ \\ \hline
$F$ & $\{t, T\}$\\ \hline
$t$ & $\{F\}$\\ \hline
$T$ & $\{t,T\}$ \\ \hline
\end{tabular}
\caption*{Table for $\circ$}
\end{minipage}
\end{figure}

When making $D=\{t, T\}$, is it is shown in Corollary $6.5.5$ of \cite{ParLog} that the Nmatrix $\mathcal{M}_{\mbCciw}=(\mathcal{A}_{\mbCciw}, D)$\label{MmbCciw} is adequate for $\mbCciw$, that is, for any set of formulas $\Gamma\cup\{\varphi\}$ over the signature $\Sigma_{\textbf{LFI}}$ we have $\Gamma\vdash_{\mbCciw}\varphi$ if, and only if, $\Gamma\vDash_{\mathcal{M}_{\mbCciw}}\varphi$.

Now consider the restricted Nmatrix \label{MmbCcl}
\[\mathcal{M}_{\mbCcl}=(\mathcal{A}_{\mbCciw}, D, \mathcal{F}_{\mbCcl})\]
such that $\mathcal{F}_{\mbCcl}$ is the set of homomorphisms $\nu:\textbf{F}(\Sigma_{\textbf{LFI}},\mathcal{V})\rightarrow\mathcal{A}_{\mbCcl}$ satisfying that, if $\nu(\alpha)=t$, then $\nu(\alpha\wedge\neg\alpha)=T$. Clearly such an RNmatrix is structural, since for any endomorphism $\sigma$ of $\textbf{F}(\Sigma_{\textbf{LFI}},\mathcal{V})$ we have that if $\nu\circ\sigma(\alpha)=t$, then $\nu(\sigma(\alpha))=t$, meaning $\nu(\sigma(\alpha)\wedge\neg\sigma(\alpha))=T$ and, therefore, $\nu\circ\sigma(\alpha\wedge\neg\alpha)=T$.

It is easy to see $\mathcal{M}_{\mbCcl}$ models the axiom schemata and rules of inference of\\ $\mbCciw$, but it is also true that, given an instance $\psi=\neg(\alpha\wedge\neg\alpha)\rightarrow\circ\alpha$ of $\cl$, we have $\vDash_{\mathcal{M}_{\mbCcl}}\psi$: assume that, for $\nu\in\mathcal{F}_{\mbCcl}$, $\nu(\psi)\notin D$, meaning that $\nu(\neg(\alpha\wedge\neg\alpha))\in D$ and $\nu(\circ\alpha)=F$, which in turn implies $\nu(\alpha)=t$.

Since $\nu\in\mathcal{F}_{\mbCcl}$, $\nu(\alpha)=t$ implies $\nu(\alpha\wedge\neg\alpha)=T$, and therefore $\nu(\neg(\alpha\wedge\neg\alpha))=F$, reaching a contradiction. We have, then, that $\nu(\psi)\in D$, for any $\nu\in\mathcal{F}_{\mbCcl}$.

\begin{theorem}
Given formulas $\Gamma\cup\{\varphi\}$ of $\mbCcl$, if $\Gamma\vdash_{\mbCcl}\varphi$ then $\Gamma\vDash_{\mathcal{M}_{\mbCcl}}\varphi$.
\end{theorem}

\begin{proof}
If $\Gamma\vdash_{\mbCcl}\varphi$, there exists a demonstration $\alpha_{1}, \dotsc  , \alpha_{n}$ of $\varphi$ from $\Gamma$, with $\alpha_{n}=\varphi$.

Let $\nu\in\mathcal{F}_{\mbCcl}$ be a valuation satisfying that that $\nu(\Gamma)\subseteq D$: we want to prove that, in this case, $\nu(\varphi)\in D$; so we prove, by induction, that $\alpha_{1}$ through $\alpha_{n}$ have image in $D$ under $\nu$, and therefore $\nu(\varphi)=\nu(\alpha_{n})=1$.

The formula $\alpha_{1}$ is either an instance of an axiom, when $\nu(\alpha_{1})\in D$ since all instances of axioms have image in $D$ through any elements of $\mathcal{F}_{\mbCcl}$, or $\alpha_{1}$ is a premise, that is, an element of $\Gamma$, and since $\nu(\Gamma)\subseteq D$ we have that $\nu(\alpha_{1})\in D$.

Suppose then that $\nu(\alpha_{1}), \dotsc  , \nu(\alpha_{i-1})\in D$, and we have three cases to consider:
\begin{enumerate}
\item if $\alpha_{i}$ is an instance of an axiom, as mentioned above $\nu(\alpha_{i})\in D$;
\item if $\alpha_{i}$ is a premise, $\nu(\alpha_{i})\in D$ since $\alpha_{i}\in\Gamma$ and $\nu(\Gamma)\subseteq D$;
\item if there are $\alpha_{j}$ and $\alpha_{k}$ with $j,k<i$ such that $\alpha_{j}=\alpha_{k}\rightarrow\alpha_{i}$ or $\alpha_{k}=\alpha_{j}\rightarrow\alpha_{i}$, since $\nu(\alpha_{j}), \nu(\alpha_{k})\in D$ we find in both cases that $\nu(\alpha_{i})\in D$, what ends the proof.
\end{enumerate}
\end{proof}

Reciprocally, consider that a map $\nu:F(\Sigma_{\textbf{LFI}},\mathcal{V})\rightarrow\{0,1\}$ is said to be a \index{Bivaluation for $\mbCciw$}bivaluation for $\mbCciw$ when it satisfies
\begin{enumerate}
\item $\nu(\alpha\vee\beta)=1$ if and only if $\nu(\alpha)=1$ or $\nu(\beta)=1$;
\item $\nu(\alpha\wedge\beta)=1$ if and only if $\nu(\alpha)=\nu(\beta)=1$;
\item $\nu(\alpha\rightarrow\beta)=1$ if and only if $\nu(\alpha)=0$ or $\nu(\beta)=1$;
\item $\nu(\alpha)=0$ implies $\nu(\neg\alpha)=1$;
\item $\nu(\circ\alpha)=1$ if and only if $\nu(\alpha)\neq\nu(\neg\alpha)$;
\end{enumerate}

A \index{Bivaluation for $\mbCcl$}bivaluation for $\mbCcl$ is simply a bivaluation for $\mbCciw$ such that, in addition to the previous conditions,
\[\nu(\circ\alpha)=0\quad\text{implies}\quad\nu(\neg(\alpha\wedge\neg\alpha))=0.\]

 If we define, for a set of formulas $\Gamma\cup\{\varphi\}$ over the signature $\Sigma_{\textbf{LFI}}$, that $\Gamma\vDash_{\mbCcl}\varphi$ whenever, for every bivaluation for $\mbCcl$ such that $\nu(\Gamma)\subseteq\{1\}$, $\nu(\varphi)=1$, it is proved in Theorem $3.3.28$ of \cite{ParLog} that $\Gamma\vdash_{\mbCcl}\varphi$ if and only if $\Gamma\vDash_{\mbCcl}\varphi$.

We want to show $\Gamma\vDash_{\mathcal{M}_{\mbCcl}}\varphi$ implies $\Gamma\vdash_{\mbCcl}\varphi$, so we will show instead that, given a bivaluation $\nu$ for $\mbCcl$, there exists a homomorphism $\overline{\nu}:\textbf{F}(\Sigma_{\textbf{LFI}},\mathcal{V})\rightarrow\mathcal{A}_{\mbCciw}$ which lies in $\mathcal{F}_{\mbCcl}$ and satisfies that $\nu(\alpha)=1$ if and only if $\overline{\nu}(\alpha)\in D$. 

This way, when we assume $\Gamma\vDash_{\mathcal{M}_{\mbCcl}}\varphi$, if $\nu(\Gamma)\subseteq\{1\}$ we have $\overline{\nu}(\Gamma)\subseteq D$, and therefore $\overline{\nu}(\varphi)\in D$, what means that $\nu(\varphi)=1$, thus proving $\Gamma\vDash_{\mbCcl}\varphi$ or, what is equivalent, $\Gamma\vdash_{\mbCcl}\varphi$.

So, suppose $\Gamma\vDash_{\mathcal{M}_{\mbCcl}}\varphi$ and let $\nu:F(\Sigma_{\textbf{LFI}},\mathcal{V})\rightarrow\{0,1\}$ be a bivaluation for $\mbCcl$: we then consider the map $\overline{\nu}:F(\Sigma_{\textbf{LFI}},\mathcal{V})\rightarrow\{F,t,T\}$ such that:
\begin{enumerate}
\item $\overline{\nu}(\alpha)=F$ if and only if $\nu(\alpha)=0$ and $\nu(\neg\alpha)=1$;
\item $\overline{\nu}(\alpha)=t$ if and only if $\nu(\alpha)=1$ and $\nu(\neg\alpha)=1$;
\item $\overline{\nu}(\alpha)=T$ if and only if $\nu(\alpha)=1$ and $\nu(\neg\alpha)=0$.
\end{enumerate}

Notice $\overline{\nu}$ is well defined since $\nu(\alpha)=0$ implies $\nu(\neg\alpha)=1$, and we therefore can not have $\nu(\alpha)=\nu(\neg\alpha)=0$.
The proof of the following results is long and tedious, but rather straightforward, so we will prefer to wait to show similar proofs in the case of $\CILA$, in Section \ref{RNmatrix for CILA}.

\begin{theorem}
$\overline{\nu}$ is a $\Sigma_{\textbf{LFI}}$-homomorphism between $\textbf{F}(\Sigma_{\textbf{LFI}}, \mathcal{V})$ and $\mathcal{A}_{\mbCcl}$.
\end{theorem}

\begin{theorem}
$\overline{\nu}$ is in $\mathcal{F}_{\mbCcl}$.
\end{theorem}

\subsection{$\CILA$}\label{RNmatrix for CILA}

\label{CILA}This logic was defined in \cite{CM}, and already in this study it was shown to be equivalent to da Costa's $C_{1}$; it is obtained over the signature $\Sigma_{\textbf{LFI}}$ by adding to the Hilbert calculus of $\textbf{mbC}$ the axiom schemata\label{cf}
\[\tag{\textbf{ci}}\neg\circ\alpha\rightarrow(\alpha\wedge\neg\alpha);\]
\[\tag{\textbf{cl}}\neg(\alpha\wedge\neg\alpha)\rightarrow\circ\alpha;\]
\[\tag{$\textbf{ca}_{\wedge}$}(\circ\alpha\wedge\circ\beta)\rightarrow\circ(\alpha\wedge\beta);\]
\[\tag{$\textbf{ca}_{\vee}$}(\circ\alpha\wedge\circ\beta)\rightarrow\circ(\alpha\vee\beta);\]
\[\tag{$\textbf{ca}_{\rightarrow}$}(\circ\alpha\wedge\circ\beta)\rightarrow\circ(\alpha\rightarrow\beta);\]
\[\tag{\textbf{cf}}\neg\neg\alpha\rightarrow\alpha.\]

Consider the $\Sigma_{\textbf{LFI}}$-multialgebra $\mathcal{A}_{\textbf{Cila}}$ given by the tables below.
\begin{figure}[H]
\centering
\begin{minipage}[t]{4cm}
\centering
\begin{tabular}{|l|c|c|r|}
\hline
$\vee$ & $F$ & $t$ & $T$ \\ \hline
$F$ & $\{F\}$ & $\{t, T\}$ & $\{T\}$ \\ \hline
$t$ & $\{t, T\}$ & $\{t, T\}$ & $\{t, T\}$ \\ \hline
$T$ & $\{T\}$ & $\{t, T\}$ & $\{T\}$ \\ \hline
\end{tabular}
\caption*{Disjunction}
\end{minipage}
\hspace{3cm}
\centering
\begin{minipage}[t]{4cm}
\centering
\begin{tabular}{|l|c|c|r|}
\hline
$\wedge$ & $F$ & $t$ & $T$ \\ \hline
$F$ & $\{F\}$ & $\{F\}$ & $\{F\}$ \\ \hline
$t$ & $\{F\}$ & $\{t, T\}$ & $\{t, T\}$ \\ \hline
$T$ & $\{F\}$ & $\{t, T\}$ & $\{T\}$ \\ \hline
\end{tabular}
\caption*{Conjunction}
\end{minipage}
\end{figure}

\begin{figure}[H]
\centering
\begin{minipage}[t]{3cm}
\centering
\begin{tabular}{|l|r|}
\hline
 & $\neg$ \\ \hline
$F$ & $\{T\}$\\ \hline
$t$ & $\{t, T\}$\\ \hline
$T$ & $\{F\}$ \\ \hline
\end{tabular}
\caption*{Negation}
\end{minipage}
\begin{minipage}[t]{5cm}
\centering
\begin{tabular}{|l|c|c|r|}
\hline
$\rightarrow$ & $F$ & $t$ & $T$ \\ \hline
$F$ & $\{T\}$ & $\{t, T\}$ & $\{T\}$ \\ \hline
$t$ & $\{F\}$ & $\{t, T\}$ & $\{t, T\}$ \\ \hline
$T$ & $\{F\}$ & $\{t, T\}$ & $\{T\}$ \\ \hline
\end{tabular}
\caption*{Implication}
\end{minipage}
\begin{minipage}[t]{3cm}
\centering
\begin{tabular}{|l|r|}
\hline
 & $\circ$ \\ \hline
$F$ & $\{T\}$\\ \hline
$t$ & $\{F\}$\\ \hline
$T$ & $\{T\}$ \\ \hline
\end{tabular}
\caption*{Consistency}
\end{minipage}
\end{figure}

We then define the restricted Nmatrix $\mathcal{M}_{\textbf{Cila}}=(\mathcal{A}_{\textbf{Cila}}, D, \mathcal{F}_{\textbf{Cila}})$\label{MCILA} such that $D=\{t, T\}$ and $\mathcal{F}_{\textbf{Cila}}$ is the set of homomorphisms $\nu:\textbf{F}(\Sigma_{\textbf{LFI}}, \mathcal{V})\rightarrow\mathcal{A}_{\CILA}$ with the property that, if $\nu(\alpha)=t$, then $\nu(\alpha\wedge\neg\alpha)=T$.

To prove the soundness of our RNmatrix, we start by defining the logic $\Ci$\label{Ci}, obtained from $\mbC$ by addition of the axiom schemata $\ci$ and $\cf$: by Theorem $6.5.24$ of \cite{ParLog}, $\Ci$ is characterized by the Nmatrix $\mathcal{M}_{\Ci}=(\mathcal{A}_{\Ci}, D)$, with $\mathcal{A}_{\Ci}$ the $\Sigma_{\textbf{LFI}}$-multialgebra given by the tables below.
\begin{figure}[H]
\centering
\begin{minipage}[t]{4cm}
\centering
\begin{tabular}{|l|c|c|r|}
\hline
$\vee$ & $F$ & $t$ & $T$ \\ \hline
$F$ & $\{F\}$ & $\{t, T\}$ & $\{t, T\}$ \\ \hline
$t$ & $\{t, T\}$ & $\{t, T\}$ & $\{t, T\}$ \\ \hline
$T$ & $\{t, T\}$ & $\{t, T\}$ & $\{t, T\}$ \\ \hline
\end{tabular}
\caption*{Disjunction}
\end{minipage}
\hspace{3cm}
\centering
\begin{minipage}[t]{4cm}
\centering
\begin{tabular}{|l|c|c|r|}
\hline
$\wedge$ & $F$ & $t$ & $T$ \\ \hline
$F$ & $\{F\}$ & $\{F\}$ & $\{F\}$ \\ \hline
$t$ & $\{F\}$ & $\{t, T\}$ & $\{t, T\}$ \\ \hline
$T$ & $\{F\}$ & $\{t, T\}$ & $\{t, T\}$ \\ \hline
\end{tabular}
\caption*{Conjunction}
\end{minipage}
\end{figure}

\begin{figure}[H]
\centering
\begin{minipage}[t]{3cm}
\centering
\begin{tabular}{|l|r|}
\hline
 & $\neg$ \\ \hline
$F$ & $\{T\}$\\ \hline
$t$ & $\{t, T\}$\\ \hline
$T$ & $\{F\}$ \\ \hline
\end{tabular}
\caption*{Negation}
\end{minipage}
\begin{minipage}[t]{5cm}
\centering
\begin{tabular}{|l|c|c|r|}
\hline
$\rightarrow$ & $F$ & $t$ & $T$ \\ \hline
$F$ & $\{t, T\}$ & $\{t, T\}$ & $\{t, T\}$ \\ \hline
$t$ & $\{F\}$ & $\{t, T\}$ & $\{t, T\}$ \\ \hline
$T$ & $\{F\}$ & $\{t, T\}$ & $\{t, T\}$ \\ \hline
\end{tabular}
\caption*{Implication}
\end{minipage}
\begin{minipage}[t]{3cm}
\centering
\begin{tabular}{|l|r|}
\hline
 & $\circ$ \\ \hline
$F$ & $\{T\}$\\ \hline
$t$ & $\{F\}$\\ \hline
$T$ & $\{T\}$ \\ \hline
\end{tabular}
\caption*{Consistency}
\end{minipage}
\end{figure}

Since $\mathcal{A}_{\CILA}$ is a submultialgebra of $\mathcal{A}_{\Ci}$, it is clear $\mathcal{M}_{\CILA}$ models Modus Ponens and those axiom schemata of $\Ci$: it remains to be shown that this RNmatrix also models the axioms concerning propagation of consistency and $\cl$.

\begin{enumerate}
\item So, take an instance $\psi=(\circ\alpha\wedge\circ\beta)\rightarrow\circ(\alpha\#\beta)$ of $\ca_{\#}$, for $\#\in\{\vee, \wedge, \rightarrow\}$: we have that, for a $\nu\in\mathcal{F}_{\CILA}$, $\nu(\psi)\notin D$ if and only if $\nu(\circ\alpha\wedge\circ\beta)\in D$ and $\nu(\circ(\alpha\#\beta))=F$, which imply $\nu(\circ\alpha), \nu(\circ\beta)\in D$ and $\nu(\alpha\#\beta)=t$. 

Since $\circ$ always gives classical values (that is, either $F$ or $T$), we have $\nu(\circ\alpha)=\nu(\circ\beta)=T$, and therefore $\nu(\alpha), \nu(\beta)\in \{F, T\}$, where we finally reach the desired contradiction: if both $\nu(\alpha)$ and $\nu(\beta)$ are classically valued, so is $\nu(\alpha\#\beta)$ from the tables for $\mathcal{A}_{\CILA}$, and therefore we must have $\nu(\psi)\in D$.

\item Take an instance $\psi=\neg(\alpha\wedge\neg\alpha)\rightarrow\circ\alpha$ of $\cl$, and assume that, for $\nu\in\mathcal{F}_{\CILA}$, $\nu(\psi)\notin D$, meaning that $\nu(\neg(\alpha\wedge\neg\alpha))\in D$ and $\nu(\circ\alpha)=F$, which implies $\nu(\alpha)=t$.

However, since $\nu\in\mathcal{F}_{\CILA}$, $\nu(\alpha)=t$ in turn implies $\nu(\alpha\wedge\neg\alpha)=T$, and therefore $\nu(\neg(\alpha\wedge\neg\alpha))=F$, a contradiction. We have, then, that $\nu(\psi)\in D$.
\end{enumerate}

\begin{theorem}
Given formulas $\Gamma\cup\{\varphi\}$ of $\CILA$, if $\Gamma\vdash_{\CILA}\varphi$ then $\Gamma\vDash_{\mathcal{M}_{\CILA}}\varphi$.
\end{theorem}

To prove the completeness of our RNmatrix, we will detour through bivaluations.

\begin{definition}
A \index{Bivaluation for $\CILA$}bivaluation for $\CILA$ is a bivaluation for $\mbCcl$ satisfying additionally:
\begin{enumerate}
\item $\nu(\alpha)=0$ implies $\nu(\neg\neg\alpha)=0$;
\item $\nu(\circ\alpha)=\nu(\circ\beta)=1$ implies $\nu(\circ(\alpha\#\beta))=1$, for $\#\in\{\vee, \wedge, \rightarrow\}$;
\item $\nu(\neg\circ\alpha)=1$ implies $\nu(\circ\alpha)=0$.
\end{enumerate}

We say $\Gamma$ semantically proves $\varphi$ in $\CILA$ if, for every bivaluation $\nu$ for $\CILA$ such that $\nu(\gamma)=1$, for every $\gamma\in\Gamma$, one has that $\nu(\varphi)=1$; in this case we write $\Gamma\vDash_{\CILA}\varphi$.
\end{definition}

As it is discussed in \cite{ParLog} and \cite{Handbook}, $\Gamma\vdash_{\CILA}\varphi$ if and only if $\Gamma\vDash_{\CILA}\varphi$. We now want to prove that, if $\Gamma\vDash_{\mathcal{M}_{\CILA}}\varphi$, then $\Gamma\vdash_{\CILA}\varphi$, so what we will do instead is prove that, if $\Gamma\vDash_{\mathcal{M}_{\CILA}}\varphi$, then $\Gamma\vDash_{\CILA}\varphi$. Again, we suppose $\Gamma\vDash_{\mathcal{M}_{\CILA}}\varphi$ and let $\nu:F(\Sigma_{\textbf{LFI}},\mathcal{V})\rightarrow\{0,1\}$ be a bivaluation for $\CILA$, to then consider the map $\overline{\nu}:F(\Sigma_{\textbf{LFI}},\mathcal{V})\rightarrow\{F,t,T\}$ such that:
\begin{enumerate}
\item $\overline{\nu}(\alpha)=F$ if and only if $\nu(\alpha)=0$ and $\nu(\neg\alpha)=1$;
\item $\overline{\nu}(\alpha)=t$ if and only if $\nu(\alpha)=1$ and $\nu(\neg\alpha)=1$;
\item $\overline{\nu}(\alpha)=T$ if and only if $\nu(\alpha)=1$ and $\nu(\neg\alpha)=0$.
\end{enumerate}

\begin{theorem}
$\overline{\nu}$ is a $\Sigma_{\textbf{LFI}}$-homomorphism between $\textbf{F}(\Sigma_{\textbf{LFI}}, \mathcal{V})$ and $\mathcal{A}_{\CILA}$.
\end{theorem}

\begin{proof}
\begin{enumerate}
\item If $\overline{\nu}(\alpha)=t$, $\nu(\alpha)=\nu(\neg\alpha)=1$ and therefore $\nu(\alpha\vee\beta)=1$, meaning 
\[\overline{\nu}(\alpha\vee\beta)\in\{t, T\}=\overline{\nu}(\alpha)\vee\overline{\nu}(\beta).\]
We can do the same if $\overline{\nu}(\beta)=t$.

If $\overline{\nu}(\alpha)=T$ and $\overline{\nu}(\beta)=F$, we have $\nu(\alpha)=\nu(\neg\beta)=1$ and $\nu(\neg\alpha)=\nu(\beta)=0$, which imply $\nu(\circ\alpha)=\nu(\circ\beta)=1$ and therefore $\nu(\circ(\alpha\vee\beta))=1$; since $\nu(\alpha\vee\beta)=1$, $\nu(\neg(\alpha\vee\beta))=0$, and therefore $\overline{\nu}(\alpha\vee\beta)=T\in \overline{\nu}(\alpha)\vee\overline{\nu}(\beta)$, the same happening if $\overline{\nu}(\alpha)=F$ and $\overline{\nu}(\beta)=T$.

If $\overline{\nu}(\alpha)=F$ and $\overline{\nu}(\beta)=F$, we have that $\nu(\alpha\vee\beta)=0$ and, since $\nu(\circ(\alpha\vee\beta))=1$, one has $\nu(\neg(\alpha\vee\beta))=1$; so $\overline{\nu}(\alpha\vee\beta)=F\in \overline{\nu}(\alpha)\vee\overline{\nu}(\beta)$. If $\overline{\nu}(\alpha)=T$ and $\overline{\nu}(\beta)=T$, $\nu(\alpha\vee\beta)=1$ and $\nu(\circ(\alpha\vee\beta))=1$, that is, $\nu(\neg(\alpha\vee\beta))=0$: with this, $\overline{\nu}(\alpha\vee\beta)=T\in \overline{\nu}(\alpha)\vee\overline{\nu}(\beta)$.

This finishes proving that, for any values of $\overline{\nu}(\alpha)$ and $\overline{\nu}(\beta)$, $\overline{\nu}(\alpha\vee\beta)\in \overline{\nu}(\alpha)\vee\overline{\nu}(\beta)$.

\item If $\overline{\nu}(\alpha)=F$, $\nu(\alpha\wedge\beta)=0$ and therefore $\overline{\nu}(\alpha\wedge\beta)=F\in\overline{\nu}(\alpha)\wedge\overline{\nu}(\beta)$. The same can be done if $\overline{\nu}(\beta)=F$.

If $\overline{\nu}(\alpha)=t$ and $\overline{\nu}(\beta)$ is either $t$ or $T$, we have $\nu(\alpha)=\nu(\beta)=1$ and therefore $\nu(\alpha\wedge\beta)=1$, meaning 
\[\overline{\nu}(\alpha\wedge\beta)\in \{t, T\}=\overline{\nu}(\alpha)\wedge\overline{\nu}(\beta).\]
The same can be done if the conditions are reversed, $\overline{\nu}(\alpha)$ is either $t$ and $T$, and $\overline{\nu}(\beta)=t$.

Finally, if $\overline{\nu}(\alpha)=\overline{\nu}(\beta)=T$, $\nu(\alpha\wedge\beta)=1$ and, since $\nu(\circ(\alpha\wedge\beta))=1$, $\nu(\neg(\alpha\wedge\beta))=0$, what means $\overline{\nu}(\alpha\wedge\beta)=T\in \overline{\nu}(\alpha)\wedge\overline{\nu}(\beta)$.

\item First of all, suppose $\overline{\nu}(\beta)=t$, what signifies $\nu(\beta)=1$ and allows to find $\nu(\alpha\rightarrow\beta)=1$, \textit{id est}, $\overline{\nu}(\alpha\rightarrow\beta)\in\{t, T\}=\overline{\nu}(\alpha)\rightarrow\overline{\nu}(\beta)$.

If $\overline{\nu}(\alpha)=F$ and $\overline{\nu}(\beta)$ is either $F$ or $T$, we clearly have $\nu(\alpha\rightarrow\beta)=1$; and since $\nu(\circ\alpha)=\nu(\circ\beta)=1$, $\nu(\circ(\alpha\rightarrow\beta))=1$, meaning $\nu(\neg(\alpha\rightarrow\beta))=0$ and therefore $\overline{\nu}(\alpha\rightarrow\beta)=T\in \overline{\nu}(\alpha)\rightarrow\overline{\nu}(\beta)$.

If $\overline{\nu}(\beta)=F$ and $\overline{\nu}(\alpha)\in \{t, T\}$, $\nu(\alpha)=1$ but $\nu(\beta)=0$, implying that $\nu(\alpha\rightarrow\beta)=0$, $\overline{\nu}(\alpha\rightarrow\beta)=F\in \overline{\nu}(\alpha)\rightarrow\overline{\nu}(\beta)$.

If $\overline{\nu}(\alpha)=t$ and $\overline{\nu}(\beta)=T$, $\nu(\alpha)=1$ and $\nu(\beta)=1$, what means that $\nu(\alpha\rightarrow\beta)=1$ and 
\[\overline{\nu}(\alpha\rightarrow\beta)\in\{t, T\}=\overline{\nu}(\alpha)\rightarrow\overline{\nu}(\beta).\]
If $\overline{\nu}(\alpha)=T$ and $\overline{\nu}(\beta)=T$, $\nu(\alpha\rightarrow\beta)=1$ and $\nu(\circ\alpha)=\nu(\circ\beta)=1$, from what we derive that $\nu(\circ(\alpha\rightarrow\beta))=1$ and $\nu(\neg(\alpha\rightarrow\beta))=0$, that is, $\overline{\nu}(\alpha\rightarrow\beta)=T\in \overline{\nu}(\alpha)\rightarrow\overline{\nu}(\beta)$.

\item If $\overline{\nu}(\alpha)=F$, $\nu(\alpha)=0$, $\nu(\neg\alpha)=1$ and, therefore, $\nu(\neg\neg\alpha)=0$, meaning $\overline{\nu}(\neg\alpha)=T\in \neg\overline{\nu}(\alpha)$.

If $\overline{\nu}(\alpha)=t$, $\nu(\neg\alpha)=1$, and it follows that $\overline{\nu}(\neg\alpha)\in\{t, T\}=\neg\overline{\nu}(\alpha)$.

If $\overline{\nu}(\alpha)=T$, one has $\nu(\alpha)=1$ and $\nu(\neg\alpha)=0$ and, therefore, $\nu(\neg\neg\alpha)=1$, implying $\overline{\nu}(\neg\alpha)=F\in\neg\overline{\nu}(\alpha)$.

\item If $\overline{\nu}(\alpha)$ equals $F$ or $T$, $\nu(\circ\alpha)=1$ and $\nu(\neg\circ\alpha)=0$, since otherwise we would be forced to have $\nu(\circ\alpha)=0$; this means $\overline{\nu}(\circ\alpha)=T\in\circ\overline{\nu}(\alpha)$.

If $\overline{\nu}(\alpha)=t$, $\nu(\circ\alpha)=0$ and therefore $\overline{\nu}(\circ\alpha)=F\in \circ\overline{\nu}(\alpha)$.
\end{enumerate}
\end{proof}

\begin{theorem}
$\overline{\nu}$ is in $\mathcal{F}_{\CILA}$.
\end{theorem}

\begin{proof}
If $\overline{\nu}(\alpha)=t$, $\nu(\alpha)=\nu(\neg\alpha)=1$ and therefore $\nu(\alpha\wedge\neg\alpha)=1$; meanwhile, $\nu(\alpha)=\nu(\neg\alpha)=1$ imply $\nu(\circ\alpha)=0$ and therefore $\nu(\neg(\alpha\wedge\neg\alpha))=0$, meaning $\overline{\nu}(\alpha\wedge\neg\alpha)=T$.
\end{proof}

From the facts that $\overline{\nu}$ is in $\mathcal{F}_{\CILA}$, $\Gamma\vDash_{\mathcal{M}_{\CILA}}\varphi$ and $\nu(\gamma)=1$ for every $\gamma\in\Gamma$, implying $\overline{\nu}(\gamma)\in D$ for every $\gamma\in\Gamma$, we obtain that $\overline{\nu}(\varphi)\in D$, meaning $\nu(\varphi)=1$. With this, $\Gamma\vDash_{\CILA}\varphi$, and we have that $\Gamma\vdash_{\CILA}\varphi$, what finishes proving $\mathcal{M}_{\CILA}$ characterizes $\CILA$.

\begin{example}
With the described RNmatrix for $\CILA$, we have obtained a simple, elegant and rather efficient decision procedure for $\CILA$ and, given their equivalence, $C_{1}$ as well.\footnote{For how this decision procedure works in practice, and the proof that it actually works, look at Section \ref{Row-branching,row-eliminating} ahead.} So, let us stop for a moment to actually apply the method. Consider
\[\neg(\alpha\vee\beta)\rightarrow(\neg\alpha\wedge\neg\beta),\]
which is not a theorem of $\CILA$ (nor of $C_{1}$). Consider the homomorphism $\nu:\textbf{F}(\Sigma_{\textbf{LFI}}, \mathcal{V})\rightarrow \mathcal{A}_{\CILA}$, lying in $\mathcal{F}_{\CILA}$, such that 
\[\nu(\alpha)=t,\quad\nu(\beta)=T,\quad\nu(\alpha\vee\beta)=t\quad\text{and}\quad\nu(\neg(\alpha\vee\beta))=t.\]
Then we have $\nu(\neg\beta)=F$, what means $\nu(\neg\alpha\wedge\neg\beta)=F$ (regardless of the value of $\neg\alpha$) and therefore $\nu(\neg(\alpha\vee\beta)\rightarrow(\neg\alpha\wedge\neg\beta))=F$, what indeed shows that the formula is not a theorem.
\end{example}

\newpage
\printbibliography[segment=\therefsegment,heading=subbibliography]
\end{refsegment}

\begin{refsegment}
\defbibfilter{notother}{not segment=\therefsegment}
\setcounter{chapter}{4}
\chapter{RNmatrices for da Costa's hierarchy}\label{Chapter5}\label{Chapter 5}

The most important result of this chapter will be a full treatment of da Costa's hierarchy. In 1963, Newton C. A. da Costa\index{da Costa} presented his \textit{Tese de C{\'a}tedra} (the Brazilian equivalent of  an habilitation thesis) under the title ``Sistemas Formais Inconsistentes'' (\textit{inconsistent formal systems} in Portuguese, \cite{Costa3}). This work is one of the founding studies in paraconsistent logic, together with Stanislaw Ja{\'s}kowski discussive (or discursive) logic $\textbf{D2}$, based on suggestions made to him by Jan Lukaziewicz (\cite{Jaskowski, Jaskowski2}). da Costa's systems $C_{n}$, however, war far better scrutinized by him than $\textbf{D2}$ by Ja{\'s}kowski, what may explain their far greater success in the long run.

One of da Costa's innovations was the separation of the formulas for a given $C_{n}$ in two: those that, indeed, do not trivialize an argument when presented together with their negation (and are, in a sense, ``badly-behaved''), and those that do, usually called ``well-behaved'' or ``classically-behaved''. To facilitate this division, he introduced a connective $\circ_{n}$, for each $C_{n}$, which asserts the well-behavior of a formula when analyzed in light of the explosion law: this way, in $C_{n}$, $\{\alpha, \neg\alpha\}$ does not necessarily have a trivial deductive closure (meaning its negation is not explosive), while $\{\alpha, \neg\alpha, \circ_{n}\alpha\}$ does have it, implying
\[\alpha, \neg\alpha, \circ_{n}\alpha\vdash_{C_{n}}\beta\]
for any formulas $\alpha$ and $\beta$ of $C_{n}$. This approach was further generalized, later on, by Carnielli\index{Carnielli} and Marcos \cite{CM}, through the notion of logics of formal inconsistency ($\textbf{LFI}$'s) which adds the connective $\circ$, the consistency operator, as a primitive one.

Many logicians will agree that da Costa's hierarchy $C_{n}$ is astoundingly interesting, but very few will deny that it is also of exceedingly difficult treatment: by Avron's aforementioned results, $C_{1}$, which is equivalent in a different signature to $\CILA$, can not be characterized by a finite Nmatrix; even more, it is not characterizable by a finite set of finite Nmatrices, and is not algebraizable by \index{Blok and Pigozzi}Blok and Pigozzi's methodology (\cite{Avron, Mortensen80, Mortensen89, Lewin}). It is, fortunately, decidable, and decision methods vary from bivaluations (\cite{Loparic}) to Fidel structures (\cite{Fidel3}). 

We here provide other decision procedures, now based on RNmatrices, for the entire hierarchy, which we believe to be easier to apply, given the very algebraic nature of RNmatrices, than many other such procedures found in the literature; in particular, they seem to suggest that the $n$th logic of da Costa should be seem as an $n+2$-valued logic, something other approaches cannot do. In Chapter \ref{Chapter6} we proceed to generalize these RNmatrices for $C_{n}$ to a more general, and categorical, semantics through restricted swap structures (\cite{ParLog}).

Most of the work we present in this chapter can be found submitted to an online repository in \cite{CostaRNmatrix}, and published in \cite{TwoDecisionProcedures}.

\section{Definition of da Costa's hierarchy}

Consider the signature $\Sigma_{\textbf{C}}$\label{SigmaC} with $(\Sigma_{\textbf{C}})_{0}=\emptyset$, $(\Sigma_{\textbf{C}})_{1}=\{\neg\}$, $(\Sigma_{\textbf{C}})_{2}=\{\vee, \wedge, \rightarrow\}$ and $(\Sigma_{\textbf{C}})_{n}=\emptyset$ for $n>2$.

For simplicity, we define $\alpha^{0}=\alpha$\label{alpha^} and
\[\alpha^{n+1}=\neg(\alpha^{n}\wedge\neg(\alpha^{n}))\]
for $n\in\mathbb{N}$; we also define $\alpha^{(0)}=\alpha$, $\alpha^{(1)}=\alpha^{1}$\label{alpha^()} and 
\[\alpha^{(n+1)}=\alpha^{(n)}\wedge\alpha^{n+1}\]
for $n\in\mathbb{N}\setminus\{0\}$. Again for simplicity, $\alpha^{1}$ may be denoted by $\alpha^{\circ}$\label{alphacirc}, since this formula will play a role equivalent to $\circ\alpha$'s role in, for example, $\mbC$.

The Hilbert calculus for the positive fragment of intuitionistic logic\index{Logic, Intuitionistic} is composed of the following axiom schemata
\begin{enumerate}
\item[\textbf{Ax\: 1}] $\alpha\rightarrow(\beta\rightarrow\alpha)$;
\item[\textbf{Ax\: 2}] $\big(\alpha\rightarrow (\beta\rightarrow \gamma)\big)\rightarrow\big((\alpha\rightarrow\beta)\rightarrow(\alpha\rightarrow\gamma)\big)$;
\item[\textbf{Ax\: 3}] $\alpha\rightarrow\big(\beta\rightarrow(\alpha\wedge\beta)\big)$;
\item[\textbf{Ax\: 4}] $(\alpha\wedge\beta)\rightarrow \alpha$;
\item[\textbf{Ax\: 5}] $(\alpha\wedge\beta)\rightarrow \beta$;
\item[\textbf{Ax\: 6}] $\alpha\rightarrow(\alpha\vee\beta)$;
\item[\textbf{Ax\: 7}] $\beta\rightarrow(\alpha\vee\beta)$;
\item[\textbf{Ax\: 8}] $(\alpha\rightarrow\gamma)\rightarrow\Big((\beta\rightarrow\gamma)\rightarrow \big((\alpha\vee\beta)\rightarrow\gamma\big)\Big)$;
\end{enumerate}
plus Modus Ponens as inference rule
\[\frac{\alpha\quad\alpha\rightarrow\beta}{\beta}.\]

We define the da costa's system $C_{\omega}$\label{Comega}, over $\Sigma_{\textbf{C}}$, by adding to the Hilbert calculus for the positive fragment of intuitionistic logic the axiom schemata
\begin{enumerate}
\item[\textbf{Ax\: 10}] $\alpha\vee\neg\alpha$;
\item[\textbf{cf}] $\neg\neg\alpha\rightarrow\alpha$.
\end{enumerate}

\begin{definition}
For $n\in\mathbb{N}\setminus\{0\}$, we define the da Costa's\index{Hierarchy, da Costa's} system $C_{n}$\label{Cn}, over the signature $\Sigma_{\textbf{C}}$, as the logic obtained from $C_{\omega}$ by addition of the axiom schemata\label{bcn}\label{phashn}
\[\tag{$\textbf{bc}_{n}$}\alpha^{(n)}\rightarrow\big(\alpha\rightarrow(\neg\alpha\rightarrow\beta)\big);\]
\[\tag{$\textbf{p}\vee_{n}$}(\alpha^{(n)}\wedge\beta^{(n)})\rightarrow(\alpha\vee\beta)^{(n)};\]
\[\tag{$\textbf{p}\wedge_{n}$}(\alpha^{(n)}\wedge\beta^{(n)})\rightarrow(\alpha\wedge\beta)^{(n)};\]
\[\tag{$\textbf{p}\rightarrow_{n}$}(\alpha^{(n)}\wedge\beta^{(n)})\rightarrow(\alpha\rightarrow\beta)^{(n)}.\]
\end{definition}

These systems were originally developed by Newton da Costa, in his seminal work \cite{Costa3}; originally, the axiom $\textbf{bc}_{n}$ was presented as
\[\alpha^{(n)}\rightarrow\Big((\beta\rightarrow\alpha)\rightarrow\big((\beta\rightarrow\neg\alpha)\rightarrow\neg\beta\big)\Big),\]
which deals with non-contradiction instead of explosivity, but it is clear how both approaches are equivalent. For further comments on this finer distinction, refer back to \cite{ParLog}.

\begin{definition}\label{bival-def}
A bivaluation for $C_{n}$, also known as a $C_{n}$-bivaluation\index{Bivaluation for $C_{n}$}, is a function\\ $\mathsf{b}:\textbf{F}(\Sigma_{\textbf{C}}, \mathcal{V})\rightarrow\{0,1\}$ satisfying:
\begin{enumerate}
\item[$(B1)$] $\mathsf{b}(\alpha\wedge\beta)=1$ if and only if $\mathsf{b}(\alpha)=1$ and $\mathsf{b}(\beta)=1$;
\item[$(B2)$] $\mathsf{b}(\alpha\vee\beta)=1$ if and only if $\mathsf{b}(\alpha)=1$ or $\mathsf{b}(\beta)=1$;
\item[$(B3)$] $\mathsf{b}(\alpha\rightarrow\beta)=1$ if and only if $\mathsf{b}(\alpha)=0$ or $\mathsf{b}(\beta)=1$;
\item[$(B4)$] $\mathsf{b}(\alpha)=0$ implies $\mathsf{b}(\neg\alpha)=1$;
\item[$(B5)$] $\mathsf{b}(\neg\neg\alpha)=1$ implies $\mathsf{b}(\alpha)=1$;
\item[$(B6)_{n}$] $\mathsf{b}(\alpha^{n-1})=\mathsf{b}(\neg(\alpha^{n-1}))$ if and only if $\mathsf{b}(\alpha^{n})=0$;
\item[$(B7)$] $\mathsf{b}(\alpha)=\mathsf{b}(\neg\alpha)$ if and only if $\mathsf{b}(\neg(\alpha^{1}))=1$;
\item[$(B8)$] for any $\#\in\{\vee, \wedge, \rightarrow\}$, $\mathsf{b}(\alpha)\neq\mathsf{b}(\neg\alpha)$ and $\mathsf{b}(\beta)\neq\mathsf{b}(\neg\beta)$ imply, together, that $\mathsf{b}(\alpha\#\beta)\neq\mathsf{b}(\neg(\alpha\#\beta))$.
\end{enumerate}
\end{definition}

The previous notion of bivaluations for $C_{n}$ has its origin in \cite{Loparic}. If we denote the fact that, for every $C_{n}$-bivaluation $\mathsf{b}$, $\mathsf{b}(\Gamma)\subseteq\{1\}$ implies $\mathsf{b}(\varphi)=1$ by $\Gamma\vDash_{C_{n}}\varphi$, it is proved in the same paper that 
\[\Gamma\vdash_{C_{n}}\varphi\quad\text{if and only if}\quad\Gamma\vDash_{C_{n}}\varphi,\]
for any set of formulas $\Gamma\cup\{\varphi\}$ over the signature $\Sigma_{\textbf{C}}$.

We shall provide, in the next sections, a semantical approach to da Costa's hierarchy trough finite restricted Nmatrices.

\section{$C_{2}$}

We will start from the simpler case that is $C_{2}$. For $\mathsf{b}$ a $C_{2}$-bivaluation, $(B6)_{2}$ implies that $\mathsf{b}(\alpha^{1})=\mathsf{b}(\neg(\alpha^{1}))$ if and only if $\mathsf{b}(\alpha^{2})=0$; so $\mathsf{b}(\alpha)=\mathsf{b}(\neg\alpha)=\mathsf{b}(\alpha^{1})=1$ implies, from $(B7)$, that $\mathsf{b}(\alpha^{2})=0$. Furthermore, $\mathsf{b}(\alpha)=0$ leads us to $\mathsf{b}(\neg\neg\alpha)=0$ from $(B5)$, and so the following four scenarios are possible.

\begin{figure}[H]
\centering
\begin{tabular}{|l|c|c|c|c|c|c|r|}
\hline
$\alpha$ & $\neg\alpha$ & $\alpha\wedge\neg\alpha$ & $\alpha^{1}$ & $\neg(\alpha^{1})$ & $\alpha^{1}\wedge\neg(\alpha^{1})$ & $\alpha^{2}$ & $\alpha^{(2)}$\\ \hline
\multirow{3}{*}{$1$} & \multirow{2}{*}{$1$} & \multirow{2}{*}{$1$} & $1$ & $1$ & $1$ & $0$ & $0$\\\cline{4-8}
& & & $0$ & $1$ & $0$ & $1$ & $0$\\\cline{2-8}
& $0$ & $0$ & $1$ & $0$ & $0$ & $1$ & $1$\\\hline
$0$ & $1$ & $0$ & $1$ & $0$ & $0$ & $1$ & $1$\\\hline
\end{tabular}
\caption*{Table for the scenarios on $C_{2}$}
\end{figure}

Now, for a formula $\alpha$ and a $C_{2}$-bivaluation, we would like to consider the triple $(\mathsf{b}(\alpha), \mathsf{b}(\neg\alpha), \mathsf{b}(\alpha^{1}))$ in $\{0,1\}^{3}$; from the previous table, we see that there are only $4$ possibilities, that we list below:
\[T_{2}=(1,0,1),\quad t^{2}_{0}=(1,1,0),\quad t^{2}_{1}=(1,1,1,)\quad\text{and}\quad F_{2}=(0,1,1).\]
We call the set of those elements $B_{2}$, and it is clear that 
\[B_{2}=\{z\in \{0,1\}^{3}\ :\  z_{1}\vee z_{2}=1\quad\text{and}\quad (z_{1}\wedge z_{2})\vee z_{3}=1\},\]
where $z_{i}$ will denote the $i$th coordinate of an element $z\in\{0,1\}^{3}$. We then define the $\Sigma_{\textbf{C}}$-multialgebra $\mathcal{A}_{C_{2}}$ with universe $B_{2}$ and operations given by
\[\tilde{\neg}z=\{w\in B_{2}\ :\  w_{1}=z_{2}\quad\text{and}\quad w_{2}\leq z_{1}\}\]
and
\[\text{for}\quad \#\in\{\vee, \wedge, \rightarrow\},\quad z\tilde{\#}w=\begin{cases*}
\{u\in B_{2}\ :\  u_{1}=z_{1}\# w_{1}\}\cap Boo_{2} & if $z, w\in Boo_{2}$\\
\{u\in B_{2}\ :\  u_{1}=z_{1}\# w_{1}\} & otherwise
\end{cases*},\]
where: we will denote an operation $\sigma_{\mathcal{A}_{C_{2}}}$ on $\mathcal{A}_{C_{2}}$ simply by $\tilde{\sigma}$; the operations $\#\in\{\vee, \wedge, \rightarrow\}$ are defined as usual in $\{0,1\}$; and $Boo_{2}=\{F_{2}, T_{2}\}$ is the set of classically-behaving elements.

One could argue that a more natural definition would be that $z\tilde{\#}w$ always equals $\{u\in B_{2}\ :\  u_{1}=z_{1}\# w_{1}\}$, but that would be problematic. Condition $(B8)$ for being a bivaluation for $C_{n}$ states that $\mathsf{b}(\alpha)\neq\mathsf{b}(\neg\alpha)$ and $\mathsf{b}(\beta\neq\mathsf{b}(\neg\beta)$ imply, together, that $\mathsf{b}(\alpha\#\beta)\neq\mathsf{b}(\neg(\alpha\#\beta))$; and one easily sees that the elements $(\mathsf{b}(\alpha), \mathsf{b}(\neg\alpha), \mathsf{b}(\alpha^{1}))$ of $B_{2}$ satisfying $\mathsf{b}(\alpha)\neq\mathsf{b}(\neg\alpha)$ are precisely $F_{2}$ and $T_{2}$. This all means that $(B8)$ translates to demanding that, if $x,y\in Boo_{2}$, then one must necessarily have $x\tilde{\#}y\in Boo_{2}$. Notice that this also justifies the nomenclature of ``classical'' for the elements $F_{2}$ and $T_{2}$: they correspond to formulas such that their negation, and the formula itself, are not simultaneously true.

If we denote the set $\{t_{1}^{2}, t_{0}^{2}, T_{2}\}$ by $D_{2}$, the tables for $\mathcal{A}_{C_{2}}$ are as below. We will also use the notation $\{t_{0}^{2}, t_{1}^{2}\}=I_{2}$ when necessary.
\begin{figure}[H]
\centering
\begin{minipage}[t]{4cm}
\centering
\begin{tabular}{|l|c|c|c|r|}
\hline
$\vee$ & $F_{2}$ & $t_{1}^{2}$ & $t_{0}^{2}$ & $T_{2}$ \\ \hline
$F_{2}$ & $\{F_{2}\}$ & $D_{2}$ & $D_{2}$ & $\{T_{2}\}$ \\ \hline
$t_{1}^{2}$ & $D_{2}$ & $D_{2}$ & $D_{2}$ & $D_{2}$ \\ \hline
$t_{0}^{2}$ & $D_{2}$ & $D_{2}$ & $D_{2}$ & $D_{2}$ \\ \hline
$T_{2}$ & $\{T_{2}\}$ & $D_{2}$ & $D_{2}$ & $\{T_{2}\}$ \\ \hline
\end{tabular}
\caption*{Table for Disjunction}
\end{minipage}
\hspace{3cm}
\centering
\begin{minipage}[t]{4cm}
\centering
\begin{tabular}{|l|c|c|c|r|}
\hline
$\wedge$ & $F_{2}$ & $t_{1}^{2}$ & $t_{0}^{2}$ & $T_{2}$ \\ \hline
$F_{2}$ & $\{F_{2}\}$ & $\{F_{2}\}$ & $\{F_{2}\}$ & $\{F_{2}\}$ \\ \hline
$t_{1}^{2}$ & $\{F_{2}\}$ & $D_{2}$ & $D_{2}$ & $D_{2}$ \\ \hline
$t_{0}^{2}$ & $\{F_{2}\}$ & $D_{2}$ & $D_{2}$ & $D_{2}$ \\ \hline
$T_{2}$ & $\{F_{2}\}$ & $D_{2}$ & $D_{2}$ & $\{T_{2}\}$ \\ \hline
\end{tabular}
\caption*{Table for Conjunction}
\end{minipage}
\end{figure}

\begin{figure}[H]
\centering
\begin{minipage}[t]{4cm}
\centering
\begin{tabular}{|l|r|}
\hline
 & $\neg$ \\ \hline
$F_{2}$ & $\{T_{2}\}$\\ \hline
$t_{1}^{2}$ & $D_{2}$\\ \hline
$t_{0}^{2}$ & $D_{2}$\\ \hline
$T_{2}$ & $\{F_{2}\}$ \\ \hline
\end{tabular}
\caption*{Table for negation}
\end{minipage}
\hspace{3cm}
\centering
\begin{minipage}[t]{5cm}
\centering
\begin{tabular}{|l|c|c|c|r|}
\hline
$\rightarrow$ & $F_{2}$ & $t_{1}^{2}$ & $t_{0}^{2}$ & $T_{2}$ \\ \hline
$F_{2}$ & $\{T_{2}\}$ & $D_{2}$ & $D_{2}$ & $\{T_{2}\}$ \\ \hline
$t_{1}^{2}$ & $\{F_{2}\}$ & $D_{2}$ & $D_{2}$ & $D_{2}$ \\ \hline
$t_{0}^{2}$ & $\{F_{2}\}$ & $D_{2}$ & $D_{2}$ & $D_{2}$ \\ \hline
$T_{2}$ & $\{F_{2}\}$ & $D_{2}$ & $D_{2}$ & $\{T_{2}\}$ \\ \hline
\end{tabular}
\caption*{Table for Implication}
\end{minipage}
\end{figure}

If we consider the Nmatrix $\mathcal{M}_{C_{2}}=(\mathcal{A}_{C_{2}}, D_{2})$, we cannot hope to characterize $C_{2}$ with it given that $C_{2}$ is not characterizable by a single finite Nmatrix (\cite{Avron}). What we will do instead, and which will be a successful endeavor, is to restrict the set of valuations for this Nmatrix, creating an RNmatrix, in order to characterize $C_{2}$.

\begin{definition}\label{FC2}
Let $\mathcal{F}_{C_{2}}$ be the set of homomorphisms $\nu:\textbf{F}(\Sigma_{\textbf{C}}, \mathcal{V})\rightarrow \mathcal{A}_{C_{2}}$ (which are called valuations over $\mathcal{A}_{C_{2}}$) such that:
\begin{enumerate}
\item if $\nu(\alpha)=t_{0}^{2}$, then $\nu(\alpha\wedge\neg\alpha)=T_{2}$;
\item if $\nu(\alpha)=t_{1}^{2}$, then $\nu(\alpha\wedge\neg\alpha)\in I_{2}$ and $\nu(\alpha^{\circ})=t_{0}^{2}$.
\end{enumerate}
We will denote the restricted Nmatrix $(\mathcal{A}_{C_{2}}, D_{2}, \mathcal{F}_{C_{2}})$ by $\mathcal{RM}_{C_{2}}$.
\end{definition}

Suppose $\nu\in\mathcal{F}_{C_{2}}$: notice that, if $\nu(\alpha)=t_{0}^{2}$, then $\nu(\alpha\wedge\neg\alpha)=T_{2}$ and so $\nu(\alpha^{\circ})=F_{2}$; if $\nu(\alpha)=t_{1}^{2}$, $\nu(\alpha^{\circ})=t_{0}^{2}$ and therefore $\nu(\alpha^{\circ\circ})=F_{2}$. Let $\bot_{\alpha}$ denote $(\alpha\wedge\neg\alpha)\wedge\alpha^{(2)}$ and ${\sim}\alpha$ denote $\alpha\rightarrow\bot_{\alpha}$ (the strong negation definable in $C_{2}$), and we arrive at the following table, where an asterisk signifies that a certain value would be different if the table were constructed using the Nmatrix $(\mathcal{A}_{C_{2}}, D_{2})$.

\begin{figure}[H]
\centering
\begin{tabular}{|l|c|c|c|c|c|c|c|c|r|}
\hline
$\alpha$ & $\neg\alpha$ & $\alpha\wedge\neg\alpha$ & $\alpha^{\circ}$ & $\neg(\alpha^{\circ})$ & $\alpha^{\circ}\wedge\neg(\alpha^{\circ})$ & $\alpha^{\circ\circ}$ & $\alpha^{(2)}$ & $\bot_{\alpha}$ & ${\sim}\alpha$ \\ \hline
$T_{2}$ & $F_{2}$ & $F_{2}$ & $T_{2}$ & $F_{2}$ & $F_{2}$ & $T_{2}$ & $T_{2}$ & $F_{2}$ & $F_{2}$\\ \hline
$t_{0}^{2}$ & $D_{2}$ & $T_{2}^{*}$ & $F_{2}$ & $T_{2}$ & $F_{2}$ & $T_{2}$ & $F_{2}$ & $F_{2}$ & $F_{2}$\\ \hline
$t_{1}^{2}$ & $D_{2}$ & $I_{2}^{*}$ & $t_{0}^{2 *}$ & $D_{2}$ & $T_{2}^{*}$ & $F_{2}$ & $F_{2}$ & $F_{2}$ & $F_{2}$ \\ \hline
$F_{2}$ & $T_{2}$ & $F_{2}$ & $T_{2}$ & $F_{2}$ & $F_{2}$ & $T_{2}$ & $T_{2}$ & $F_{2}$ & $T_{2}$ \\ \hline
\end{tabular}
\caption*{Table for the scenarios in $\mathcal{RM}_{C_{2}}$}
\end{figure}

\begin{proposition}
For a homomorphism $\nu$ in $\mathcal{F}_{C_{2}}$ and an endomorphism $\sigma:\textbf{F}(\Sigma_{\textbf{C}}, \mathcal{V})\rightarrow\textbf{F}(\Sigma_{\textbf{C}}, \mathcal{V})$, $\nu\circ\sigma\in \mathcal{F}_{C_{2}}$.
\end{proposition}

\begin{proof}
Of course $\nu\circ\sigma:\textbf{F}(\Sigma_{\textbf{C}}, \mathcal{V})\rightarrow\mathcal{A}_{C_{2}}$ remains a homomorphism.
\begin{enumerate}
\item If $\nu\circ\sigma(\alpha)=t_{0}^{2}$, given $\nu$ is in $\mathcal{F}_{C_{2}}$ we derive that $\nu(\sigma(\alpha)\wedge\neg\sigma(\alpha))=T_{2}$; since $\sigma$ is an endomorphism of $\textbf{F}(\Sigma_{\textbf{C}}, \mathcal{V})$, $\sigma(\alpha)\wedge\neg\sigma(\alpha)=\sigma(\alpha\wedge\neg\alpha)$, and so $\nu\circ\sigma(\alpha\wedge\neg\alpha)=T_{2}$.

\item If $\nu\circ\sigma(\alpha)=t_{1}^{2}$, since $\nu\in\mathcal{F}_{C_{2}}$ we get $\nu(\sigma(\alpha)\wedge\neg\sigma(\alpha))\in I_{2}$ and 
\[\nu(\sigma(\alpha)^{\circ})=\nu\Big(\neg(\sigma(\alpha)\wedge\neg\sigma(\alpha))\Big)=t_{0}^{2};\]
given $\sigma(\alpha)\wedge\neg\sigma(\alpha)=\sigma(\alpha\wedge\neg\alpha)$ and $\neg(\sigma(\alpha)\wedge\neg\sigma(\alpha))=\sigma(\neg(\alpha\wedge\neg\alpha))$, we obtain that $\nu\circ\sigma(\alpha\wedge\neg\alpha)\in I_{2}$ and $\nu\circ\sigma(\neg(\alpha\wedge\neg\alpha))=\nu\circ\sigma(\alpha^{\circ})=t_{0}^{2}$, what finishes the proof.
\end{enumerate}
\end{proof}

The previous proposition implies $\mathcal{RM}_{C_{2}}$ is structural. 

The following technical lemmas are necessary in order to prove the desired Theorem \ref{Char C2}.

\begin{lemma}\label{Sound C2}
For $\nu$ a valuation in $\mathcal{F}_{C_{2}}$, the mapping $\mathsf{b}: F(\Sigma_{\textbf{C}}, \mathcal{V})\rightarrow \{0,1\}$, such that $\mathsf{b}(\alpha)=1$ if and only if $\nu(\alpha)\in D_{2}$, is a $C_{2}$-bivaluation.
\end{lemma}

\begin{proof}
Let us see that $\mathsf{b}$ satisfies the clauses of Definition \ref{bival-def} for $n=2$. Clauses $(B1)$ through $(B3)$ are clearly satisfied: since 
\[\mathsf{b}(\alpha \# \beta) = \nu(\alpha \# \beta)_1 = \nu(\alpha)_1 \# \nu(\beta)_1\]
for $\# \in \{\vee, \wedge, \rightarrow\}$. 

On the other hand, $\mathsf{b}(\neg\alpha) = \nu(\neg\alpha)_1 = \nu(\alpha)_2$, hence $\mathsf{b}(\alpha) \vee \mathsf{b}(\neg\alpha)=1$, by definition of $B_{2}$. From this, clause $(B2)$ is satisfied. Since $\mathsf{b}(\neg\neg\alpha)=\nu(\neg\neg\alpha)_1=\nu(\neg\alpha)_2\leq \nu(\alpha)_1 = \mathsf{b}(\alpha)$ it follows that $\mathsf{b}$ satisfies $(B5)$. 

Concerning $(B6)_2$, suppose that $\mathsf{b}(\alpha^\circ)=\mathsf{b}(\neg(\alpha^\circ))=1$: this means that $\nu(\alpha^\circ)_1=\nu(\neg(\alpha^\circ))_1=\nu(\alpha^\circ)_2=1$. That is, $\nu(\alpha^\circ) \in \{t^2_0,t^2_1\}$. From the possible scenarios in $\mathcal{RM}_{C_{2}}$ it follows that $\nu(\alpha^\circ) = t^2_1$ and so $\nu(\alpha^2)=F$, what implies that $\mathsf{b}(\alpha^2)=\nu(\alpha^2)_1=0$. Conversely, $\mathsf{b}(\alpha^2)=\nu(\alpha^2)_1=0$ means that $\nu(\alpha^2)=F$, which implies that $\nu(\alpha^\circ) = t^2_1$, and so $\mathsf{b}(\alpha^\circ)=\mathsf{b}(\neg(\alpha^\circ))=1$. From this, clause $(B6)_2$ is satisfied. 

Now,  $\mathsf{b}(\alpha)=\mathsf{b}(\neg\alpha)$ if and only if $\nu(\alpha)_1 = \nu(\alpha)_2=1$, what in turn happens if and only if $\nu(\alpha) \in \{t^2_0,t^2_1\}$. This last condition is equivalent to the fact that $\nu(\neg(\alpha^\circ)) \in D_2$, by the possible scenarios in $\mathcal{RM}_{C_{2}}$, which is equivalent to $\mathsf{b}(\neg(\alpha^\circ))=1$. Then, $(B7)$ is satisfied. 

Finally, suppose that $\mathsf{b}(\alpha) \neq \mathsf{b}(\neg\alpha)$ and $v(\beta) \neq \mathsf{b}(\neg\beta)$. This means that $\nu(\alpha),\nu(\beta) \in \{T_2,F_2\}$ and so $\nu(\alpha\# \beta) \in \{T_2,F_2\}$, by definition of  $\mathcal{A}_{C_2}$. Hence clause $(B8)$ is fulfilled, and the proof is complete.
\end{proof}

\begin{lemma}\label{Comp C2}
For $\mathsf{b}$ a $C_{2}$-bivaluation, the mapping $\nu:\textbf{F}(\Sigma_{\textbf{C}}, \mathcal{V})\rightarrow \mathcal{A}_{C_{2}}$, such that $\nu(\alpha)=(\mathsf{b}(\alpha), \mathsf{b}(\neg\alpha), \mathsf{b}(\alpha^{\circ}))$, is a valuation which lies in $\mathcal{F}_{C_{2}}$ and satisfies that $\mathsf{b}(\alpha)=1$ if, and only if, $\nu(\alpha)\in D_{2}$.
\end{lemma}

\begin{proof}
By definition, $\nu(\neg\alpha)=(\mathsf{b}(\neg\alpha),\mathsf{b}(\neg\neg\alpha),\mathsf{b}((\neg\alpha)^\circ))$, and clearly this belongs to $\tilde{\neg} \nu(\alpha)$ according to the definition of $\tilde{\neg}$, by the property $(B5)$ of $\mathsf{b}$. Analogously, by definition of $\tilde{\#}$ it follows that $\nu(\alpha\# \beta) \in \nu(\alpha) \tilde{\#} \nu(\beta)$, with use of the properties $(B1)$ through $(B3)$ of $\mathsf{b}$. This shows that $\nu$ is a valuation for $\mathcal{A}_{C_2}$. It remains to be shown that $\nu$ satisfies conditions $1$ and $2$ of Definition \ref{FC2}.

Regarding the first of these conditions, assume that $\nu(\alpha) = t^2_0$. This means that $\mathsf{b}(\alpha)=\mathsf{b}(\neg\alpha)=1$ and $\mathsf{b}(\alpha^\circ)=0$. Let $\beta=\alpha\wedge\neg\alpha$ (thus $\neg\beta=\alpha^\circ$). By $(B1)$, $\mathsf{b}(\beta)=1$. Since, by hypothesis, $\mathsf{b}(\beta) \neq \mathsf{b}(\neg\beta)$, it follows that $\mathsf{b}(\neg(\beta^\circ))=0$, by~$( B7)$, and so $\mathsf{b}(\beta^\circ)=1$, by~$(B4)$. This shows that 
\[\nu(\alpha\wedge\neg\alpha) = (1,0,1)=T_2,\]
as we wished to prove. 

Regarding the second condition, suppose that $\nu(\alpha) = t^2_1$. Then, $\mathsf{b}(\alpha)=\mathsf{b}(\neg\alpha)=\mathsf{b}(\alpha^\circ)=1$. Consider again $\beta=\alpha \land \neg\alpha$. From this, $\mathsf{b}(\beta)=\mathsf{b}(\neg\beta)=1$ and so 
\[\nu(\alpha\wedge\neg\alpha) = (\mathsf{b}(\beta), \mathsf{b}(\neg\beta), \mathsf{b}(\beta^\circ)) \in \{t^2_1,t^2_1\}.\]
Since $\mathsf{b}(\alpha)=\mathsf{b}(\neg\alpha)$ it follows by $(B7)$ that  $\mathsf{b}(\neg(\alpha^\circ))=1$. Hence, since $\mathsf{b}(\alpha^1)=\mathsf{b}(\neg(\alpha^1))$, we infer that $\mathsf{b}(\alpha^2)=0$, by $(B6)_2$. That is, $\nu(\alpha^\circ) = (1,1,0)=t^2_0$. This, of course, finishes proving that $\nu$ is in $\mathcal{F}_{C_{2}}$.
\end{proof}

\begin{theorem}\label{Char C2}
Given formulas $\Gamma\cup\{\varphi\}$ of $C_{2}$, $\Gamma\vdash_{C_{2}}\varphi$ if, and only if, $\Gamma\vDash_{\mathcal{RM}_{C_{2}}}\varphi$.
\end{theorem}

\begin{proof}
Suppose first that $\Gamma \vdash_{C_{2}} \varphi$, and let $\nu$ be  a valuation in $\mathcal{F}_{C_2}$ such that $\nu(\gamma) \in D_{2}$ for every $\gamma \in \Gamma$. By Lemma \ref{Sound C2}, $\mathsf{b}$ defined as $\mathsf{b}(\alpha)=\nu(\alpha)_1$, for every $\alpha$, is a $C_2$-bivaluation such that, by definition, $\mathsf{b}(\gamma) =1$ for every $\gamma \in \Gamma$. By  hypothesis and by soundness of $C_2$ with respect to bivaluations, it follows that $\mathsf{b}(\varphi)=1$. That is, $\nu(\varphi) \in D_2$. This shows that $\Gamma\vDash_{\mathcal{RM}_{C_{2}}}\varphi$.

Now assume that $\Gamma\vDash_{\mathcal{RM}_{C_{2}}}\varphi$, and let $\mathsf{b}$ be a $C_2$-bivaluation such that $\mathsf{b}(\gamma) =1$ for every $\gamma \in \Gamma$. By Lemma \ref{Comp C2}, the function $\nu$ defined as $\nu(\alpha)=(\mathsf{b}(\alpha),\mathsf{b}(\neg\alpha),\mathsf{b}(\alpha^\circ))$, for every $\alpha$, is a valuation in $\mathcal{F}_{C_2}$ such that, by definition, $\nu(\gamma) \in D_2$ for every $\gamma \in \Gamma$. By hypothesis, $\nu(\varphi) \in D_2$, which means that $\mathsf{b}(\varphi)=1$. By completeness of $C_2$  with respect to bivaluations, it follows that $\Gamma \vdash_{C_2} \varphi$.
\end{proof}

\section{The general case}\label{The general case}

\begin{lemma}\label{Values of i-consistency}
Given a $C_{n}$-bivaluation $\mathsf{b}$, if $\mathsf{b}(\alpha^{i})=1$ and $\mathsf{b}(\neg(\alpha^{i}))=0$ for some $i\in\mathbb{N}$, then $\mathsf{b}(\alpha^{j})=1$ for all $j\geq i$.
\end{lemma}

\begin{proof}
We prove that, for all $j\geq i$, $\mathsf{b}(\alpha^{j})=1$ and $\mathsf{b}(\neg(\alpha^{j}))=0$ by induction, being the base case done. So, suppose $\mathsf{b}(\alpha^{j})=1$ and $\mathsf{b}(\neg(\alpha^{j}))=0$: we have $\mathsf{b}(\alpha^{j}\wedge\neg(\alpha^{j}))=0$ from $(B1)$, and $\mathsf{b}(\alpha^{j+1})=\mathsf{b}(\neg(\alpha^{j}\wedge\neg(\alpha^{j})))=1$ from $(B4)$; and from $(B7)$, $\mathsf{b}(\neg(\alpha^{j+1}))=\mathsf{b}(\neg((\alpha^{j})^{1}))=0$, since $\mathsf{b}(\alpha^{j})\neq\mathsf{b}(\neg(\alpha^{j}))$, what finishes the proof.
\end{proof}

\begin{proposition}
If $\mathsf{b}$ is a $C_{n}$-bivaluation, at most one of the elements $\mathsf{b}(\alpha)$, $\mathsf{b}(\neg\alpha)$, $\mathsf{b}(\alpha^{1}), \dotsc  ,
$\\$\mathsf{b}(\alpha^{n-1})$ equals $0$.
\end{proposition}

\begin{proof}
Suppose $\mathsf{b}(\alpha)=0$: from condition $(B4)$, this means $\mathsf{b}(\neg\alpha)=1$. Then $\mathsf{b}(\alpha\wedge\neg\alpha)=0$ (from $(B1)$), meaning $\mathsf{b}(\alpha^{1})=\mathsf{b}(\neg(\alpha\wedge\neg\alpha))=1$, again from $(B4)$; but notice, furthermore, that, from the converse of $(B5)$, $\mathsf{b}(\alpha\wedge\neg\alpha)=0$ implies $\mathsf{b}(\neg(\alpha^{1}))=\mathsf{b}(\neg\neg(\alpha\wedge\neg\alpha))=0$. From Lemma \ref{Values of i-consistency}, we find $\mathsf{b}(\alpha^{2})=\cdots=\mathsf{b}(\alpha^{n-1})=1$, what finishes proving that, if $\mathsf{b}(\alpha)=0$, then $\mathsf{b}(\neg\alpha)=\mathsf{b}(\alpha^{1})=\cdots=\mathsf{b}(\alpha^{n-1})=1$.

Now, suppose $\mathsf{b}(\neg\alpha)=0$: we already have $\mathsf{b}(\alpha)=1$, since otherwise one could derive $\mathsf{b}(\neg\alpha)=1$ from the previous remarks. So $\mathsf{b}(\alpha\wedge\neg\alpha)=0$ and therefore $\mathsf{b}(\alpha^{1})=1$. Again by the converse of $(B5)$, one finds $\mathsf{b}(\neg(\alpha^{1}))=0$, and again from Lemma \ref{Values of i-consistency}, we obtain $\mathsf{b}(\alpha^{1})=\cdots=\mathsf{b}(\alpha^{n-1})=1$.

Finally, suppose that for some $1\leq i<j< n-1$ one has $\mathsf{b}(\alpha^{i})=\mathsf{b}(\alpha^{j})=0$: we find $\mathsf{b}(\alpha^{i}\wedge\neg(\alpha^{i}))=0$, $\mathsf{b}(\alpha^{i+1})=1$ and, from $(B5)$'s converse, $\mathsf{b}(\neg(\alpha^{i+1}))=0$; once again using Lemma \ref{Values of i-consistency}, we obtain $\mathsf{b}(\alpha^{i+1})=\cdots=\mathsf{b}(\alpha^{n-1})=1$, contradicting $\mathsf{b}(\alpha^{j})=0$. The obvious conclusion is that at most one of $\mathsf{b}(\alpha^{1})$ trough $\mathsf{b}(\alpha^{n-1})$ can equal $0$, in which case one also has $\mathsf{b}(\alpha)=\mathsf{b}(\neg\alpha)=1$.
\end{proof}

We consider $(n+1)$-tuples $z=(z_{1}, \dotsc  , z_{n+1})$ in $\textbf{2}^{n+1}$, corresponding to $\mathsf{b}(\alpha)$, $\mathsf{b}(\neg\alpha)$ and $\mathsf{b}(\alpha^{1})$ trough $\mathsf{b}(\alpha^{n-1})$ under a given $C_{n}$-bivaluation $\mathsf{b}$, which we will call snapshots\index{Snapshots} for $C_{n}$; given the restrictions established by the proposition, there are precisely $n+2$ of such tuples:\label{Tn}\label{Fn}
\[T_{n}=(1, 0, 1, \dotsc  , 1), \quad t_{0}^{n}=(1, 1, 0, 1, \dotsc  , 1),\quad\cdots\quad,\quad t_{n-2}^{n}=(1, \dotsc  , 1, 0), \quad t_{n-1}^{n}=(1, \dotsc  , 1)\]
\[\text{and}\quad F_{n}=(0, 1, \dotsc  , 1).\]
The set $B_{n}$\label{Bn} of all snapshots for $C_{n}$ is clearly the set of all elements of $\{0,1\}^{n+1}$ with at most one coordinate equal to $0$, and may be defined trough the following equality.
\[B_{n}=\{z\in \textbf{2}^{n+1}\ :\  (\bigwedge_{i=1}^{k}z_{i})\vee z^{k+1}=1,\quad\text{for every $1\leq k\leq n$}\}.\]

The set of designated elements, that will soon be part of appropriate restricted Nmatrices for $C_{n}$ with universe $B_{n}$, is $D_{n}=B_{n}\setminus\{F_{n}\}$\label{Dn}. Other useful sets are: 
\begin{enumerate}
\item the set of undesignated elements, $U_{n}=\{F_{n}\}$;
\item the set of Boolean elements\index{Boolean elements}, $Boo_{n}=\{F_{n}, T_{n}\}$, also equal to $\{z\in B_{n}\ :\  z_{1}\wedge z_{2}=0\}$\label{Boon}; 
\item and the set of inconsistent elements, $I_{n}=B_{n}\setminus Boo_{n}$.
\end{enumerate}

Notice that an element $z$ of $B_{n}$ is in $Boo_{n}$ if, and only if, it equals $(a, {\sim}a, 1, \dotsc  , 1)$, for some element $a$ of $\textbf{2}$ and ${\sim}a$ its complement in $\textbf{2}$; this is easy to see given that $z_{1}\vee z_{2}=1$, from the definition of $B_{n}$, and $z_{1}\wedge z_{2}=0$, from the definition of $Boo_{n}$.

Now, we will define the multialgebra $\mathcal{A}_{C_{n}}$\label{ACn}, which happens to be a swap structure\index{Swap structure} in the sense of \cite{ParLog}; unlike the swap structures in this reference, however, we will need to restrict the valuations taken into consideration, making of the semantics with underlying $\mathcal{A}_{C_{n}}$ an RNmatrix.

\begin{definition}\label{Define operations in Cn}
We define the multioperations on the $\Sigma_{\textbf{C}}$-multialgebra $\mathcal{A}_{C_{n}}$\label{ACn}, with universe $B_{n}$, by means of the following equations:
\[\tilde{\neg}z=\{w\in B_{n}\ :\  w_{1}=z_{2}\quad\text{and}\quad w_{2}\leq z_{1}\}\]
and
\[\text{for}\quad \#\in\{\vee, \wedge, \rightarrow\},\quad z\tilde{\#}w=\begin{cases*}
\{u\in B_{n}\ :\  u_{1}=z_{1}\# w_{1}\}\cap Boo_{n} & if $z, w\in Boo_{n}$\\
\{u\in B_{n}\ :\  u_{1}=z_{1}\# w_{1}\} & otherwise
\end{cases*}.\]
\end{definition}

The operations on this multialgebra may be illustrated, concisely, by the following tables.
\begin{figure}[H]
\centering
\begin{minipage}[t]{4cm}
\centering
\begin{tabular}{|l|c|c|r|}
\hline
$\vee$ & $F_{n}$ & $I_{n}$ & $T_{n}$ \\ \hline
$F_{n}$ & $\{F_{n}\}$ &  $D_{n}$ & $\{T_{n}\}$ \\ \hline
$I_{n}$ & $D_{n}$ & $D_{n}$ & $D_{n}$ \\ \hline
$T_{n}$ & $\{T_{n}\}$ & $D_{n}$ & $\{T_{n}\}$ \\ \hline
\end{tabular}
\caption*{Table for Disjunction}
\end{minipage}
\hspace{3cm}
\centering
\begin{minipage}[t]{4cm}
\centering
\begin{tabular}{|l|c|c|r|}
\hline
$\wedge$ & $F_{n}$ & $I_{n}$ & $T_{n}$ \\ \hline
$F_{n}$ & $\{F_{n}\}$ & $\{F_{n}\}$ & $\{F_{n}\}$ \\ \hline
$I_{n}$ & $\{F_{n}\}$  & $D_{n}$ & $D_{n}$ \\ \hline
$T_{n}$ & $\{F_{n}\}$ & $D_{n}$ & $\{T_{n}\}$ \\ \hline
\end{tabular}
\caption*{Table for Conjunction}
\end{minipage}
\end{figure}

\begin{figure}[H]
\centering
\begin{minipage}[t]{4cm}
\centering
\begin{tabular}{|l|r|}
\hline
 & $\neg$ \\ \hline
$F_{n}$ & $\{T_{n}\}$\\ \hline
$I_{n}$ & $D_{n}$\\ \hline
$T_{n}$ & $\{F_{n}\}$ \\ \hline
\end{tabular}
\caption*{Table for negation}
\end{minipage}
\hspace{3cm}
\centering
\begin{minipage}[t]{5cm}
\centering
\begin{tabular}{|l|c|c|r|}
\hline
$\rightarrow$ & $F_{n}$ & $I_{n}$ & $T_{n}$ \\ \hline
$F_{n}$ & $\{T_{n}\}$ & $D_{n}$ & $\{T_{n}\}$ \\ \hline
$I_{n}$ & $\{F_{n}\}$ & $D_{n}$ & $D_{n}$ \\ \hline
$T_{n}$ & $\{F_{n}\}$ & $D_{n}$ & $\{T_{n}\}$ \\ \hline
\end{tabular}
\caption*{Table for Implication}
\end{minipage}
\end{figure}

\begin{definition}
Let $\mathcal{F}_{C_{n}}$\label{FCn} be the set of valuations $\nu$ over $\mathcal{A}_{C_{n}}$ satisfying that:
\begin{enumerate}
\item $\nu(\alpha)=t_{0}^{n}$ implies $\nu(\alpha\wedge\neg\alpha)=T_{n}$;
\item for every $2\leq k\leq n$, $\nu(\alpha)=t_{k-1}^{n}$ implies $\nu(\alpha\wedge\neg\alpha)\in I_{n}$ and $\nu(\alpha^{1})=t_{k-2}^{n}$.
\end{enumerate}
We denote the RNmatrix $(\mathcal{A}_{C_{n}}, D_{n}, \mathcal{F}_{C_{n}})$ by $\mathcal{RM}_{C_{n}}$\label{RMCn}.
\end{definition}

\begin{proposition}
For a homomorphism $\nu$ in $\mathcal{F}_{C_{n}}$ and an endomorphism $\sigma:\textbf{F}(\Sigma_{\textbf{C}}, \mathcal{V})\rightarrow\textbf{F}(\Sigma_{\textbf{C}}, \mathcal{V})$, $\nu\circ\sigma\in \mathcal{F}_{C_{n}}$.
\end{proposition}

\begin{proof}
Given the composition of homomorphisms returns homomorphisms, $\nu\circ\sigma:\textbf{F}(\Sigma_{\textbf{C}}, \mathcal{V})\rightarrow\mathcal{A}_{C_{n}}$ is certainly still a homomorphism.

\begin{enumerate}
\item Suppose $\nu\circ\sigma(\alpha)=t_{0}^{n}$: this means $\nu(\sigma(\alpha))=t_{0}^{n}$, and since $\nu\in\mathcal{F}_{C_{n}}$, $\nu(\sigma(\alpha)\wedge\neg\sigma(\alpha))=T_{n}$. Given $\alpha$ is an endomorphism of $\textbf{F}(\Sigma_{\textbf{C}}, \mathcal{V})$, 
\[\sigma(\alpha)\wedge\neg\sigma(\alpha)=\sigma(\alpha\wedge\neg\alpha),\]
and so $\nu\circ\sigma(\alpha\wedge\neg\alpha)=\nu(\sigma(\alpha\wedge\neg\alpha))=T_{n}$.

\item For a $2\leq k\leq n$, suppose $\nu\circ\sigma(\alpha)=t_{k-1}^{n}$. That is, $\nu(\sigma(\alpha))=t_{k-1}^{n}$, and since $\nu$ is in $\mathcal{F}_{C_{n}}$, we have that $\nu(\sigma(\alpha)\wedge\neg\sigma(\alpha))\in I_{n}$ and $\nu(\sigma(\alpha)^{1})=\nu(\neg(\sigma(\alpha)\wedge\neg\sigma(\alpha)))=t_{k-2}^{n}$. Again, given $\sigma$ is an endomorphism of $\textbf{F}(\Sigma_{\textbf{C}}, \mathcal{V})$,
\[\sigma(\alpha)\wedge\neg\sigma(\alpha)=\sigma(\alpha\wedge\neg\alpha)\quad\text{and}\quad \neg(\sigma(\alpha)\wedge\neg\sigma(\alpha))=\sigma(\neg(\alpha\wedge\neg\alpha)),\]
what implies $\nu\circ\sigma(\alpha\wedge\neg\alpha)=\nu(\sigma(\alpha\wedge\neg\alpha))\in I_{n}$ and $\nu\circ\sigma(\alpha^{1}))=\nu(\sigma(\alpha^{1}))=t_{k-2}^{n}$.
\end{enumerate}
\end{proof}

The previous result proves $\mathcal{RM}_{C_{n}}$ is structural.

Notice that, for $\nu\in\mathcal{F}_{C_{n}}$, $\nu(\alpha)=t_{0}^{n}$ implies $\nu(\alpha\wedge\neg\alpha)=T_{n}$ and, therefore, $\nu(\alpha^{1})=F_{n}$; if $\nu(\alpha)=t_{1}^{n}$, $\nu(\alpha^{1})=t_{0}^{n}$ and, by the previous remark, $\nu(\alpha^{2})=F_{n}$. Proceeding inductively, we discover that if $\nu(\alpha)=t_{i}^{n}$, $\nu(\alpha^{j})\in D_{n}$, for $1\leq j\leq i$, and $\nu(\alpha^{i+1})=F_{n}$. With the aid of these observations, we obtain the following possible scenarios in $\mathcal{RM}_{C_{n}}$, where we will write $X^{*}$ to mean that the value $X$ is obtained from the fact a valuation is in $\mathcal{F}_{C_{n}}$.

\begin{figure}[H]
\centering
\begin{tabular}{|l|c|c|c|c|c|c|c|c|c|r|}
\hline
$\alpha$ & $\alpha\wedge\neg\alpha$ & $\alpha^{1}$ & $\alpha^{1}\wedge\neg(\alpha^{1})$ & $\alpha^{2}$ & $\alpha^{2}\wedge\neg(\alpha^{2})$ & $\cdots$ & $\alpha^{n-1}$ & $\alpha^{n-1}\wedge\neg(\alpha^{n-1})$ & $\alpha^{n}$ & $\alpha^{(n)}$\\ \hline
$T_{n}$ & $F_{n}$ & $T_{n}$ & $F_{n}$ & $T_{n}$ & $F_{n}$ & $\cdots$ & $T_{n}$ & $F_{n}$ & $T_{n}$ & $T_{n}$ \\ \hline
$t_{0}^{n}$ & $T_{n}^{*}$ & $F_{n}$ & $F_{n}$ & $T_{n}$ & $F_{n}$ & $\cdots$ & $T_{n}$ & $F_{n}$ & $T_{n}$ & $F_{n}$ \\ \hline
$t_{1}^{n}$ & $I_{n}^{*}$ & $t_{0}^{n *}$ & $T_{n}^{*}$ & $F_{n}$ & $F_{n}$ & $\cdots$ & $T_{n}$ & $F_{n}$ & $T_{n}$ & $F_{n}$ \\ \hline
$t_{2}^{n}$ & $I_{n}^{*}$ & $t_{1}^{n *}$ & $I_{n}^{*}$ & $t_{0}^{n *}$ & $T_{n}^{*}$ & $\cdots$ & $T_{n}$ & $F_{n}$ & $T_{n}$ & $F_{n}$ \\ \hline
$\vdots$ & $\vdots$ & $\vdots$ & $\vdots$ & $\vdots$ & $\vdots$ & $\ddots$ & $\vdots$ & $\vdots$ & $\vdots$ & $\vdots$ \\ \hline
$t_{n-3}^{n}$ & $I_{n}^{*}$ & $t_{n-4}^{n *}$ & $I_{n}^{*}$ & $t_{n-5}^{n *}$ & $I_{n}^{*}$ & $\cdots$ & $T_{n}$ & $F_{n}$ & $T_{n}$ & $F_{n}$ \\ \hline
$t_{n-2}^{n}$ & $I_{n}^{*}$ & $t_{n-3}^{n *}$ & $I_{n}^{*}$ & $t_{n-4}^{n *}$ & $I_{n}^{*}$ & $\cdots$ & $F_{n}$ & $F_{n}$ & $T_{n}$ & $F_{n}$ \\ \hline
$t_{n-1}^{n}$ & $I_{n}^{*}$ & $t_{n-2}^{n *}$ & $I_{n}^{*}$ & $t_{n-3}^{n *}$ & $I_{n}^{*}$ & $\cdots$ & $t_{0}^{n *}$ & $T_{n}^{*}$ & $F_{n}$ & $F_{n}$ \\ \hline
$F_{n}$ & $F_{n}$ & $T_{n}$ & $F_{n}$ & $T_{n}$ & $F_{n}$ & $\cdots$ &$T_{n}$ & $F_{n}$ & $T_{n}$ & $T_{n}$ \\ \hline
\end{tabular}
\caption*{Table for the scenarios in $\mathcal{RM}_{C_{n}}$}
\end{figure}

One sees that, if $\nu(\alpha)=t_{i}^{n}$, then the fact that $\nu\in\mathcal{F}_{C_{n}}$ restricts the possible values of $2i+1$ formulas of the form $\alpha^{j}$ or $\alpha^{k}\wedge\neg(\alpha^{k})$: $\alpha\wedge\neg\alpha$, $\alpha^{1}$, $\alpha^{1}\wedge\neg(\alpha^{1}), \dotsc  , \alpha^{i}$ and $\alpha^{i}\wedge\neg(\alpha^{i})$.

Also worthy of note is the fact that, if we define $\bot_{\alpha}=(\alpha\wedge\neg\alpha)\wedge\alpha^{(n)}$ (the definable bottom of $C_{n}$), for all $\nu\in\mathcal{F}_{C_{n}}$ one has $\nu(\bot_{\alpha})=F_{n}$; and if we define ${\sim}\alpha=\alpha\rightarrow\bot_{\alpha}$, for $\nu\in\mathcal{F}_{C_{n}}$ (the definable strong negation on $C_{n}$) we find that $\nu(\alpha)\in D_{n}$ implies $\nu({\sim}\alpha)=F_{n}$, while $\nu(\alpha)=F_{n}$ implies $\nu({\sim}\alpha)=T_{n}$.

\begin{lemma}\label{Scenarions for RMCn}
Denote, for simplicity, $T_{n}$ by $t^{n}_{-1}$; for any formula $\alpha$ on $C_{n}$, valuation $\nu$ on $\mathcal{F}_{C_{n}}$ and integer $1\leq k\leq n$, one has that:
\begin{enumerate}
\item if $\nu(\alpha)=t^{n}_{i}$, for $-1\leq i\leq k-2$, $\nu(\alpha^{k})=T_{n}$;
\item if $\nu(\alpha)=t^{n}_{k-1}$, $\nu(\alpha^{k})=F_{n}$;
\item if $\nu(\alpha)=t^{n}_{i}$, for $k\leq i\leq n-1$, $\nu(\alpha^{k})=t^{n}_{i-k}$;
\item if $\nu(\alpha)=F_{n}$, $\nu(\alpha^{k})=T_{n}$.
\end{enumerate}
\end{lemma}

\begin{proof}
We will prove the lemma by induction on $k$, starting with the case $k=1$.
\begin{enumerate}
\item If $\nu(\alpha)=t^{n}_{-1}$, meaning $\nu(\alpha)=T_{n}$, $\nu(\neg\alpha)=F_{n}$ and therefore $\nu(\alpha\wedge\neg\alpha)=F_{n}$, implying $\nu(\alpha^{1})=T_{n}$.
\item If $\nu(\alpha)=t^{n}_{0}$, from the fact that $\nu$ is in $\mathcal{F}_{C_{n}}$ we derive that $\nu(\alpha\wedge\neg\alpha)=T_{n}$, and therefore $\nu(\alpha^{1})=F_{n}$.
\item If $\nu(\alpha)=t^{n}_{i}$, for $1\leq i\leq n-1$, again from the fact that $\nu\in\mathcal{F}_{C_{n}}$ we find $\nu(\alpha^{1})=t^{n}_{i-1}$.
\item Finally, if $\nu(\alpha)=F_{n}$, $\nu(\alpha\wedge\neg\alpha)=F_{n}$ and therefore $\nu(\alpha^{1})=T_{n}$.
\end{enumerate}

Now, suppose the result holds for $k-1$.

\begin{enumerate}
\item If $\nu(\alpha)=t^{n}_{i}$, for $-1\leq i\leq k-2$, there are two cases to consider: if $i\leq k-3$, by induction hypothesis $\nu(\alpha^{k-1})=T_{n}$, meaning $\nu(\neg(\alpha^{k-1}))=F_{n}$, $\nu(\alpha^{k-1}\wedge\neg(\alpha^{k-1}))=F_{n}$ and $\nu(\alpha^{k})=T_{n}$; if $i=k-2$, again by induction hypothesis $\nu(\alpha^{k-2})=F_{n}$, implying $\nu(\alpha^{k-1}\wedge\neg(\alpha^{k-1}))=F_{n}$ and therefore $\nu(\alpha^{k})=T_{n}$.

\item If $\nu(\alpha)=t^{n}_{k-1}$, by induction hypothesis $\nu(\alpha^{k-1})=t^{n}_{0}$, and from the fact $\nu$ is in $\mathcal{F}_{C_{n}}$ we get $\nu(\alpha^{k-1}\wedge\neg(\alpha^{k-1}))=T_{n}$, meaning $\nu(\alpha^{k})=F_{n}$.

\item If $\nu(\alpha)=t^{n}_{i}$, for $k\leq i\leq n-1$, by induction hypothesis we find that $\nu(\alpha^{k-1})=t^{n}_{i-k+1}$. Since $k\leq i\leq n-1$ and $1\leq k\leq n$, it follows that  $1\leq i-k+1\leq n-1$, and from the fact that $\nu$ is in $\mathcal{F}_{C_{n}}$ we obtain $\nu(\alpha^{k})=t^{n}_{i-k}$.

\item Finally, if $\nu(\alpha)=F_{n}$, again by induction hypothesis $\nu(\alpha^{k-1})=F_{n}$, meaning $\nu(\alpha^{k-1}\wedge\neg(\alpha^{k-1}))=F_{n}$ and therefore $\nu(\alpha^{k})=T_{n}$.
\end{enumerate}
\end{proof}

\begin{lemma}\label{restricted valuations and swap part 1}
Let $\nu$ be a valuation in $\mathcal{F}_{C_{n}}$; then the mapping $\mathsf{b}:\textbf{F}(\Sigma, \mathcal{V})\rightarrow\textbf{2}$ given by $\mathsf{b}(\alpha):=\nu(\alpha)_{1}$ (that is, $\mathsf{b}(\alpha)=1$ if and only if $\nu(\alpha)\in D_{n}$) is a $C_{n}$-bivaluation.
\end{lemma}

\begin{proof}
\begin{enumerate}
\item For any $\#\in\{\vee, \wedge, \rightarrow\}$, $\mathsf{b}(\alpha\#\beta)=1$ if and only if $\nu(\alpha\#\beta)_{1}=1$; since $\nu(\alpha\#\beta)_{1}=\nu(\alpha)_{1}\#\nu(\beta)_{1}$ by the definition of $\tilde{\#}$, we have that:
\begin{enumerate}
\item $\mathsf{b}(\alpha\vee\beta)=1$ if and only if either $\mathsf{b}(\alpha)=1$ or $\mathsf{b}(\beta)=1$ (clause $(B2)$ for being a $C_{n}$-bivaluation);
\item $\mathsf{b}(\alpha\wedge\beta)=1$ if and only if $\mathsf{b}(\alpha)=1$ and $\mathsf{b}(\beta)=1$ (clause $(B1)$);
\item $\mathsf{b}(\alpha\rightarrow\beta)=1$ if and only if $\mathsf{b}(\alpha)=0$ or $\mathsf{b}(\beta)=1$ (clause $(B3)$).
\end{enumerate}

\item If $\mathsf{b}(\alpha)=0$, this means $\nu(\alpha)_{1}=0$; since $\nu(\alpha)_{1}\vee\nu(\alpha)_{2}=1$, by definition of $B_{n}$, and $\nu(\neg\alpha)_{1}=\nu(\alpha)_{2}$, by definition of $\tilde{\neg}$, we find that $\mathsf{b}(\neg\alpha)=\nu(\neg\alpha)_{1}=1$, satisfying clause $(B4)$.

\item If $\mathsf{b}(\neg\neg\alpha)=1$, we have $\nu(\neg\neg\alpha)_{1}=1$; we have $\nu(\neg\neg\alpha)_{1}=\nu(\neg\alpha)_{2}$ and $\nu(\neg\alpha)_{2}\leq \nu(\alpha)_{1}$ from the definition of $\tilde{\neg}$, and so $\nu(\alpha)_{1}\geq 1$, what means that $\mathsf{b}(\alpha)=\nu(\alpha)_{1}=1$ and that clause $(B5)$ is satisfied.

\item If $\mathsf{b}(\alpha^{n-1})=\mathsf{b}(\neg(\alpha^{n-1}))$, then $\nu(\alpha^{n-1})_{1}=\nu(\neg(\alpha^{n-1}))_{1}=\nu(\alpha^{n-1})_{2}$; from looking at the definition of $B_{n}$, we see this implies $\nu(\alpha^{n-1})\in I_{n}$. From Lemma \ref{Scenarions for RMCn}, we obtain $\nu(\alpha)=t^{n}_{n-1}$, and in this case $\nu(\alpha^{n})=F_{n}$, which implies $\mathsf{b}(\alpha^{n})=\nu(\alpha^{n})_{1}=0$.

Reciprocally, if $\mathsf{b}(\alpha^{n})=0$, $\nu(\alpha^{n})=F_{n}$, and from Lemma \ref{Scenarions for RMCn} once again $\nu(\alpha)=t^{n}_{n-1}$: so $\nu(\alpha^{n-1})=t^{n}_{0}$, $\nu(\neg(\alpha^{n-1}))\in D_{n}$ and therefore 
\[\mathsf{b}(\alpha^{n-1})=\nu(\alpha^{n-1})_{1}=1=\nu(\neg(\alpha^{n-1}))_{1}=\mathsf{b}(\neg(\alpha^{n-1})).\]
We then find clause $(B6)_{n}$ is satisfied.

\item If $\mathsf{b}(\alpha)=\mathsf{b}(\neg\alpha)$, $\nu(\alpha)_{1}=\nu(\neg\alpha)_{1}=\nu(\alpha)_{2}=1$, what means $\nu(\alpha)\in I_{n}$; from Lemma \ref{Scenarions for RMCn}, in this case we have that $\nu(\alpha^{1})\in B_{n}\setminus\{T_{n}\}$, and from the definition of $\tilde{\neg}$ one obtains that $\nu(\neg(\alpha^{1}))\in D_{n}$, and so $\mathsf{b}(\neg(\alpha^{1}))=\nu(\neg(\alpha^{1}))_{1}=1$.

Reciprocally, if $\mathsf{b}(\neg(\alpha^{1}))=1$, $\nu(\neg(\alpha^{1}))_{1}=1$ and therefore $\nu(\alpha^{1})_{2}=1$, meaning $\nu(\alpha^{1})\in B_{n}\setminus\{T_{n}\}$; this implies $\nu(\alpha)\in I_{n}$, and by the definition of $\tilde{\neg}$, $\nu(\neg\alpha)\in D_{n}$. So $\mathsf{b}(\alpha)=\nu(\alpha)_{1}$ and $\mathsf{b}(\neg\alpha)=\nu(\neg\alpha)_{1}$ are both $1$, and we have handled clause $(B7)$.

\item If $\mathsf{b}(\alpha)\neq\mathsf{b}(\neg\alpha)$ and $\mathsf{b}(\beta)\neq\mathsf{b}(\neg\beta)$, $\nu(\alpha)_{1}\neq\nu(\neg\alpha)_{1}=\nu(\alpha)_{2}$ and $\nu(\beta)_{1}\neq\nu(\neg\beta)_{1}=\nu(\beta)_{2}$, meaning $\nu(\alpha), \nu(\beta)\in \{F_{n}, T_{n}\}$. From the tables for $\#\in\{\vee, \wedge, \rightarrow\}$ we see that $\nu(\alpha\#\beta)\in\{F_{n}, T_{n}\}$, and again $\mathsf{b}(\alpha\#\beta)=\nu(\alpha\#\beta)_{1}$ differs from $\mathsf{b}(\neg(\alpha\#\beta))=\nu(\neg(\alpha\#\beta))_{1}$, what constitutes clause $(B8)$ and finishes the proof.
\end{enumerate}
\end{proof}

\begin{lemma}\label{restricted valuations and swap part 2}
Let $\mathsf{b}$ be a $C_{n}$-bivaluation; then, the mapping $\nu:\textbf{F}(\Sigma, \mathcal{V})\rightarrow\mathcal{A}_{C_{n}}$, given by $\nu(\alpha)=(\mathsf{b}(\alpha), \mathsf{b}(\neg\alpha), \mathsf{b}(\alpha^{1}), \dotsc  , \mathsf{b}(\alpha^{n-1}))$ is a valuation in $\mathcal{F}_{C_{n}}$ such that $\mathsf{b}(\alpha)=1$ if and only if $\nu(\alpha)\in D_{n}$ for every formula $\alpha$.
\end{lemma}

\begin{proof}
First, we show $\nu$ is a homomorphism.
\begin{enumerate}
\item We have that $\nu(\neg\alpha)=(\mathsf{b}(\neg\alpha), \mathsf{b}(\neg\neg\alpha), \mathsf{b}((\neg\alpha)^{1}), \dotsc  ,  \mathsf{b}((\neg\alpha)^{n-1}))$, and so there are three cases to consider: if $\nu(\alpha)=T_{n}$, $\mathsf{b}(\neg\alpha)=0$ and therefore $\nu(\neg\alpha)=F_{n}$; if $\nu(\alpha)\in I_{n}$, $\mathsf{b}(\neg\alpha)=1$, implying $\nu(\neg\alpha)\in D_{n}$; finally, if $\nu(\alpha)=F_{n}$, $\mathsf{b}(\alpha)=0$ and $\mathsf{b}(\neg\alpha)=1$, and from clause $(B5)$ for being a $C_{n}$-bivaluation we get $\mathsf{b}(\neg\neg\alpha)=0$, and so $\nu(\neg\alpha)=T_{n}$. Regardless of the case, we find $\nu(\neg\alpha)\in \tilde{\neg}\nu(\alpha)$.

\item From the definition of $\nu$, $\nu(\alpha\vee\beta)=(\mathsf{b}(\alpha\vee\beta), \mathsf{b}(\neg(\alpha\vee\beta)), \mathsf{b}((\alpha\vee\beta)^{1}), \dotsc  , \mathsf{b}((\alpha\vee\beta)^{n-1})$. 
\begin{enumerate}
\item If $\nu(\alpha)$ or $\nu(\beta)$ equals $T_{n}$ and both are boolean valued, either $\mathsf{b}(\alpha)$ or $\mathsf{b}(\beta)$ equals $1$, meaning $\mathsf{b}(\alpha\vee\beta)=1$ from clause $(B2)$; and both $\mathsf{b}(\alpha)\neq\mathsf{b}(\neg\alpha)$ and $\mathsf{b}(\beta)\neq\mathsf{b}(\neg\beta)$, which implies from clause $(B8)$ that $\mathsf{b}(\alpha\vee\beta)\neq\mathsf{b}(\neg(\alpha\vee\beta))$ and therefore $\nu(\alpha\vee\beta)=T_{n}$.
\item If $\nu(\alpha)$ or $\nu(\beta)$ is in $I_{n}$, we have that either $\mathsf{b}(\alpha)=1$ or $\mathsf{b}(\beta)=1$ and from clause $(B2)$ one gets $\mathsf{b}(\alpha\vee\beta)=1$, meaning $\nu(\alpha\vee\beta)\in D_{n}$.
\item If $\nu(\alpha)$ and $\nu(\beta)$ are both $F_{n}$, $\mathsf{b}(\alpha)=\mathsf{b}(\beta)=0$ and so $\mathsf{b}(\alpha\vee\beta)=0$ from $(B2)$, meaning that $\nu(\alpha\vee\beta)=F_{n}$.
\end{enumerate}
Either way, $\nu(\alpha\vee\beta)$ is in $\nu(\alpha)\tilde{\vee}\nu(\beta)$.

\item We have $\nu(\alpha\wedge\beta)=(\mathsf{b}(\alpha\wedge\beta), \mathsf{b}(\neg(\alpha\wedge\beta)), \mathsf{b}((\alpha\wedge\beta)^{1}), \dotsc  , \mathsf{b}((\alpha\wedge\beta)^{n-1})$. 
\begin{enumerate}
\item If $\nu(\alpha)=\nu(\beta)=T_{n}$, $\mathsf{b}(\alpha)=\mathsf{b}(\beta)=1$ and $\mathsf{b}(\neg\alpha)=\mathsf{b}(\neg\beta)=0$; from clause $(B1)$, $\mathsf{b}(\alpha\wedge\beta)=1$, and from $(B8)$ we get $\mathsf{b}(\neg(\alpha\wedge\beta))=0$, meaning $\nu(\alpha\wedge\beta)=T_{n}$.
\item If $\nu(\alpha)$ or $\nu(\beta)$ equals $F_{n}$, either $\mathsf{b}(\alpha)$ or $\mathsf{b}(\beta)$ equals $0$, and therefore $\mathsf{b}(\alpha\wedge\beta)=0$ (from clause $(B1)$) and so $\nu(\alpha\wedge\beta)=F_{n}$.
\item In the remaining cases, when either $\nu(\alpha)$ or $\nu(\beta)$ is in $I_{n}$ but none of them equals $F_{n}$, one sees that $\mathsf{b}(\alpha)=\mathsf{b}(\beta)=1$ and therefore $\mathsf{b}(\alpha\wedge\beta)=1$, meaning $\nu(\alpha\wedge\beta)\in D_{n}$.
\end{enumerate}
We have just proved $\nu(\alpha\wedge\beta)\in \nu(\alpha)\tilde{\wedge}\nu(\beta)$.

\item Clearly $\nu(\alpha\rightarrow\beta)=(\mathsf{b}(\alpha\rightarrow\beta), \mathsf{b}(\neg(\alpha\rightarrow\beta)), \mathsf{b}((\alpha\rightarrow\beta)^{1}), \dotsc  , \mathsf{b}((\alpha\rightarrow\beta)^{n-1})$. 
\begin{enumerate}
\item If $\nu(\alpha)=F_{n}$ or $\nu(\beta)=T_{n}$, and both are boolean valued, $\mathsf{b}(\alpha)=0$ or $\mathsf{b}(\beta)=1$, and both $\mathsf{b}(\alpha)\neq\mathsf{b}(\neg\alpha)$ and $\mathsf{b}(\beta)\neq\mathsf{b}(\neg\beta)$. By clauses $(B3)$ and $(B8)$, this all means $\mathsf{b}(\alpha\rightarrow\beta)=1$ and $\mathsf{b}(\neg(\alpha\rightarrow\beta))=0$, and so $\nu(\alpha\rightarrow\beta)=T_{n}$.
\item When $\nu(\beta)\in I_{n}$, $\mathsf{b}(\beta)=1$ and so $\mathsf{b}(\alpha\rightarrow\beta)=1$ from clause $(B3)$, meaning $\nu(\alpha\rightarrow\beta)\in D_{n}$.
\item If $\nu(\beta)=F_{n}$ and $\nu(\alpha)\in D_{n}$, $\mathsf{b}(\beta)=0$ and $\mathsf{b}(\alpha)=1$, implying from clause $(B3)$ that $\mathsf{b}(\alpha\rightarrow\beta)=0$ and so $\nu(\alpha\rightarrow\beta)=F_{n}$.
\item Finally, if $\nu(\alpha)\in I_{n}$ and $\nu(\beta)=T_{n}$, $\mathsf{b}(\alpha)=1$ and $\mathsf{b}(\beta)=1$; from clause $(B3)$ once again, $\mathsf{b}(\alpha\rightarrow\beta)=1$ and therefore $\nu(\alpha\rightarrow\beta)\in D_{n}$.
\end{enumerate}
In all cases, $\nu(\alpha\rightarrow\beta)\in\nu(\alpha)\tilde{\rightarrow}\nu(\beta)$.
\end{enumerate}

Now, it remains to be shown that $\nu$ is in $\mathcal{F}_{C_{n}}$.
\begin{enumerate}
\item If $\nu(\alpha)=t^{n}_{0}$, we have $\mathsf{b}(\alpha)=\mathsf{b}(\neg\alpha)=1$ (and so $\mathsf{b}(\alpha\wedge\neg\alpha)=1$) and $\mathsf{b}(\neg(\alpha\wedge\neg\alpha))=\mathsf{b}(\alpha^{1})=0$, from what one derives $\nu(\alpha\wedge\neg\alpha)=T_{n}$.
\item If $\nu(\alpha)=t^{n}_{k-1}$, for $2\leq k\leq n$, we have $\mathsf{b}(\alpha)=\mathsf{b}(\neg\alpha)=1$ (meaning $\mathsf{b}(\alpha\wedge\neg\alpha)=1$) and $\mathsf{b}(\neg(\alpha\wedge\neg\alpha))=\mathsf{b}(\alpha^{1})=1$, and so $\nu(\alpha\wedge\neg\alpha)\in I_{n}$.

Furthermore, $\mathsf{b}((\alpha^{1})^{k-1})=\mathsf{b}(\alpha^{k})=0$, from the fact that $\nu(\alpha)=t^{n}_{k-1}$, and since $\mathsf{b}((\alpha^{1})^{k-1})=0$ we find $\nu(\alpha^{1})=t^{n}_{k-2}$, what ends the proof.
\end{enumerate}
\end{proof}

\begin{theorem}
Given formulas $\Gamma\cup\{\varphi\}$ of $C_{n}$, $\Gamma\vdash_{C_{n}}\varphi$ if, and only if, $\Gamma\vDash_{\mathcal{RM}_{C_{n}}}\varphi$.
\end{theorem}

\begin{proof}
First, suppose $\Gamma\vdash_{C_{n}}\varphi$, and take a valuation $\nu\in\mathcal{F}_{C_{n}}$ for which $\nu(\Gamma)\subseteq D_{n}$: from Lemma \ref{restricted valuations and swap part 1}, the function $\mathsf{b}:F(\Sigma_{\textbf{C}}, \mathcal{V})\rightarrow\textbf{2}$, defined by $\mathsf{b}(\alpha)=\nu(\alpha)_{1}$, is a bivaluation which, by hypothesis, satisfies $\mathsf{b}(\Gamma)\subseteq\{1\}$. Given the soundness of $C_{n}$ with respect to bivaluations and the fact that $\Gamma\vdash_{C_{n}}\varphi$, it follows that $\mathsf{b}(\varphi)=1$, and therefore $\nu(\varphi)\in D_{n}$. This, of course, shows that $\Gamma\vDash_{\mathcal{RM}_{C_{n}}}\varphi$.

Reciprocally, suppose $\Gamma\vDash_{\mathcal{RM}_{C_{n}}}\varphi$ and let $\mathsf{b}$ be a $C_{n}$-bivaluation satisfying $\mathsf{b}(\Gamma)\subseteq\{1\}$: by Lemma \ref{restricted valuations and swap part 2}, $\nu:\textbf{F}(\Sigma_{\textbf{C}}, \mathcal{V})\rightarrow \mathcal{A}_{C_{n}}$ defined by 
\[\nu(\alpha)=(\mathsf{b}(\alpha), \mathsf{b}(\neg\alpha), \mathsf{b}(\alpha^{1}), \dotsc  , \mathsf{b}(\alpha^{n-1})),\]
is a valuation in $\mathcal{F}_{C_{n}}$ for which, in addition, $\nu(\Gamma)\subseteq D_{n}$. Given $\Gamma\vDash_{\mathcal{RM}_{C_{n}}}\varphi$, we obtain that $\nu(\varphi)\in D_{n}$, meaning $\mathsf{b}(\varphi)=1$ and, by completeness of $C_{n}$ with respect to bivaluations, $\Gamma\vdash_{C_{n}}\varphi$.
\end{proof}

\section{Decision methods}\label{Decision methods}

A finite Nmatrix $\mathcal{M}=(\mathcal{A}, D)$ which characterizes a logic $\mathfrak{L}$ always leads to a decision method for when a formula $\varphi$ is a tautology of $\mathfrak{L}$, or when a finite set of premises $\Gamma$ deduces $\varphi$ in $\mathfrak{L}$ as long as this logic satisfies the deduction meta-theorem (if $\Gamma=\{\gamma_{1}, \dotsc  , \gamma_{n}\}$, this is equivalent to testing whether $\gamma_{1}\rightarrow(\gamma_{2}\rightarrow\cdots(\gamma_{n}\rightarrow\varphi)\cdots)$ is a tautology). The method is straightforward: one constructs a row-branching truth table for $\varphi$, based on $\mathcal{A}$, and $\varphi$ is a tautology iff all rows contain a designated value on the column corresponding to $\varphi$. Even more, if a logic is characterized by a finite class of finite Nmatrices, one still derives a decision method.

Things are not as simple when dealing with RNmatrices: even the logic characterized by a single, finite RNmatrix may not have a decision method, given that the problem of finding out whether a certain homomorphism belongs to the set of restricted homomorphisms may not be decidable; an example of something like this happening is found by looking at the RNmatrix semantics obtained by Kearns for $S4$ (\cite{Kearns}) in Example \ref{Example Kearns}. However this is not always the case: some finite RNmatrices are indeed quite efficient in inducing decision methods, and this happens for all $\mathcal{RM}_{C_{n}}$; this is due to the fact that finding those homomorphisms not relevant for a deduction in $\mathcal{RM}_{C_{n}}$ is actually easy. We start by generalizing row-branching truth-tables to row-branching, row-eliminating truth-tables, and proceed to show how tableau semantics are defined from our RNmatrices characterizing $C_{n}$.

\subsection{Row-branching, row-eliminating truth tables}\label{Row-branching,row-eliminating}

Take a formula $\varphi$ in the signature of $C_{n}$, and let $\varphi_{1}, \dotsc  , \varphi_{k}=\varphi$ be a sequence of all subformulas, proper or not, of $\varphi$ ordered by complexity.

\begin{figure}[H]
\begin{center}
\scalebox{0.9}{
\begin{tabular}{|p{1cm}|p{1cm}|p{1cm}|p{1cm}|}
\hline
$\varphi_{1}$ & $\varphi_{2}$ & $\cdots$ & $\varphi_{p}$\\\hline
\multirow{13}{*}{$T_{n}$} & \multirow{4}{*}{$T_{n}$} & \multirow{4}{*}{$\cdots$} & $T_{n}$\\\cline{4-4}
& & & $t^{n}_{0}$\\\cline{4-4}
& & & $\vdots$\\\cline{4-4}
& & & $F_{n}$\\\cline{2-4}
& \multirow{4}{*}{$t^{n}_{0}$} & \multirow{4}{*}{$\cdots$} & $T_{n}$\\\cline{4-4}
& & & $t^{n}_{0}$\\\cline{4-4}
& & & $\vdots$\\\cline{4-4}
& & & $F_{n}$\\\cline{2-4}
& $\vdots$ & $\ddots$ & $\vdots$ \\\cline{2-4}
& \multirow{4}{*}{$F_{n}$} & \multirow{4}{*}{$\cdots$} & $T_{n}$\\\cline{4-4}
& & & $t^{n}_{0}$\\\cline{4-4}
& & & $\vdots$\\\cline{4-4}
& & & $F_{n}$\\\hline
\multirow{13}{*}{$t^{n}_{0}$} & \multirow{4}{*}{$T_{n}$} & \multirow{4}{*}{$\cdots$} & $T_{n}$\\\cline{4-4}
& & & $t^{n}_{0}$\\\cline{4-4}
& & & $\vdots$\\\cline{4-4}
& & & $F_{n}$\\\cline{2-4}
& \multirow{4}{*}{$t^{n}_{0}$} & \multirow{4}{*}{$\cdots$} & $T_{n}$\\\cline{4-4}
& & & $t^{n}_{0}$\\\cline{4-4}
& & & $\vdots$\\\cline{4-4}
& & & $F_{n}$\\\cline{2-4}
& $\vdots$ & $\ddots$ & $\vdots$ \\\cline{2-4}
& \multirow{4}{*}{$F_{n}$} & \multirow{4}{*}{$\cdots$} & $T_{n}$\\\cline{4-4}
& & & $t^{n}_{0}$\\\cline{4-4}
& & & $\vdots$\\\cline{4-4}
& & & $F_{n}$\\\hline
$\vdots$ & $\vdots$ & $\ddots$ & $\vdots$\\\hline
\multirow{13}{*}{$F_{n}$} & \multirow{4}{*}{$T_{n}$} & \multirow{4}{*}{$\cdots$} & $T_{n}$\\\cline{4-4}
& & & $t^{n}_{0}$\\\cline{4-4}
& & & $\vdots$\\\cline{4-4}
& & & $F_{n}$\\\cline{2-4}
& \multirow{4}{*}{$t^{n}_{0}$} & \multirow{4}{*}{$\cdots$} & $T_{n}$\\\cline{4-4}
& & & $t^{n}_{0}$\\\cline{4-4}
& & & $\vdots$\\\cline{4-4}
& & & $F_{n}$\\\cline{2-4}
& $\vdots$ & $\ddots$ & $\vdots$ \\\cline{2-4}
& \multirow{4}{*}{$F_{n}$} & \multirow{4}{*}{$\cdots$} & $T_{n}$\\\cline{4-4}
& & & $t^{n}_{0}$\\\cline{4-4}
& & & $\vdots$\\\cline{4-4}
& & & $F_{n}$\\\hline
\end{tabular}}
\end{center}
\caption{Beginning a truth table: only propositional variables}
\label{Truth table propositional variables}
\end{figure}
\vspace{5mm}

In order to make things as clear as possible, we remember the complexity (or order) of a formula $\varphi$ is $0$ iff the formula is a propositional variable or a nullary connective, and it is $r+1$ iff $\varphi=\sigma(\alpha_{1}, \dotsc  , \alpha_{m})$, $\sigma$ is a $m$-ary connective and the maximum of the orders of $\alpha_{1}$ trough $\alpha_{m}$ is $r$. We say ``a sequence'', rather than ``the sequence'', since different subformulas of $\varphi$ may have the same order, and therefore may be interchanged without neglecting that the whole sequence remains ordered.

Of course, given that the signature of $C_{n}$ has no nullary connectives, we can be certain that the first elements of the sequence are necessarily propositional variables, and so we start to build a truth table. For each $\varphi_{i}$ which is a propositional variable we list all possible values it may assume in $\mathcal{RM}_{C_{n}}$ given the values taken by $\varphi_{1}$ trough $\varphi_{i-1}$; since the values taken by propositional variables are independent of each other, we find that, for any combination of values for $\varphi_{0}$ up to $\varphi_{i-1}$, $\varphi_{i}$ can assume $n+2$ values, $T_{n}, t^{n}_{0}, t^{n}_{1}, \dotsc  , t^{n}_{n-1}$ and $F_{n}$. If $p\geq 1$ propositional variables appear in $\varphi$, this means we start by writing down $n+2$ rows for the column headed by $\varphi_{1}$, each containing a value of $B_{n}$; for the column headed by $\varphi_{2}$, assuming $p\geq 2$, each of these $n+2$ rows is further subdivided into $n+2$ new rows, each filled with a value of $B_{n}$, to a total now of $(n+2)^{2}$ rows; and inductively, we have a table with $(n+2)^{p}$ rows, and $p$ columns (so far).

We summarize the situation so far with the table in Figure \ref{Truth table propositional variables}, showing how the first part of the process, dealing only with propositional variables, goes. For simplicity, on all row-branching truth-tables we will always merge consecutive rows on a given column with repeated values; this is very helpful when truly trying to interpret those objects as ``row-branching''.

Suppose we have written the table up to a column headed by $\varphi_{q}$, for $q\geq p$, and let us look at $\varphi_{q+1}$: it must be of either the form ~(1)~ $\varphi_{q+1}=\varphi_{i}\#\varphi_{j}$, for some $\#\in\{\vee, \wedge, \rightarrow\}$, or ~(2)~ $\varphi_{q+1}=\neg\varphi_{i}$, for $1\leq i, j\leq q$, since the subformulas of $\varphi_{q+1}$ are also subformulas of $\varphi$ of order strictly smaller. In a normal truth-table, the new column, headed by $\varphi_{q+1}$, would be filled on a given row by either $a\tilde{\#}b$, corresponding to case $(1)$, or $\tilde{\neg}a$, corresponding to case $(2)$, where: $a$ is the value taken by $\varphi_{i}$ on the aforementioned row, $b$ is the value taken by $\varphi_{j}$ and $\tilde{\#}$ is the (deterministic) operation corresponding to the connective $\#$. In a truth-table arising from an Nmatrix, we must further subdivide the row in question in as many values are found in $a\tilde{\#}b$ and fill each new row with one of these values, whenever we find ourselves in case $(1)$, or branch the row in as many rows as there are elements in $\tilde{\neg}a$ and fill each sub-row with one of these, corresponding to case $(2)$, where now $\tilde{\#}$ is a multioperation.

The case with the RNmatrix $\mathcal{RM}_{C_{n}}$ is similar, but has extra subtleties. Fix a given row for what follows.
\begin{enumerate}
\item If $\varphi_{q+1}=\varphi_{i}\wedge\neg\varphi_{i}$ and $\varphi_{i}$ takes the value $t^{n}_{0}$, then the row is not branched in the column headed by $\varphi_{q+1}$ and assumes the value $T_{n}$.
\begin{center}
\begin{tabular}{|p{1cm}|p{1cm}|p{1cm}|p{1cm}|p{1cm}|p{2cm}|p{1cm}|}
\hline
$\cdots$ & $\varphi_{i}$ & $\cdots$ & $\neg\varphi_{i}$ & $\cdots$ & $\varphi_{i}\wedge\neg\varphi_{i}$ & $\cdots$ \\\hline
$\ddots$ & $\vdots$ & $\ddots$ & $\vdots$ & $\ddots$ & $\vdots$ & $\ddots$\\\hline
\multirow{4}{*}{$\cdots$} & \multirow{4}{*}{$t^{n}_{0}$} & \multirow{4}{*}{$\cdots$} & $T_{n}$ & $\cdots$ & $T_{n}$ & $\cdots$ \\\cline{4-7}
& & &$t^{n}_{0}$ & $\cdots$ & $T_{n}$ & $\cdots$ \\\cline{4-7}
& & &$\vdots$ & $\cdots$ & $\vdots$ & $\cdots$ \\\cline{4-7}
& & &$t^{n}_{n-1}$ & $\cdots$ & $T_{n}$ & $\cdots$ \\\hline
$\ddots$ & $\vdots$ & $\ddots$ & $\vdots$ & $\ddots$ & $\vdots$ & $\ddots$\\\hline
\end{tabular}
\end{center}

\item If $\varphi_{q+1}=\varphi_{i}\wedge\neg\varphi_{i}$ and $\varphi_{i}$ takes the value $t^{n}_{l}$, for any $1\leq l\leq n-1$, the row is subdivided into $n$ new rows, each assigned a value of $I_{n}=\{t_{0}^{n}, \dotsc  , t_{n-1}^{n}\}$.

\begin{center}
\begin{tabular}{|p{1cm}|p{1cm}|p{1cm}|p{1cm}|p{1cm}|p{1.5cm}|p{1cm}|}
\hline
$\cdots$ & $\varphi_{i}$ & $\cdots$ & $\neg\varphi_{i}$ & $\cdots$ & $\varphi_{i}\wedge\neg\varphi_{i}$ & $\cdots$ \\\hline
$\ddots$ & $\vdots$ & $\ddots$ & $\vdots$ & $\ddots$ & $\vdots$ & $\ddots$ \\\hline
\multirow{13}{*}{$\cdots$} & \multirow{13}{*}{$t^{n}_{l}$} & \multirow{13}{*}{$\cdots$} & \multirow{4}{*}{$T_{n}$} & \multirow{4}{*}{$\cdots$} & $t^{n}_{0}$ & $\cdots$ \\\cline{6-7}
&  &   &   &  & $\vdots$ & $\ddots$ \\\cline{6-7}
&  &   &   & & $t^{n}_{n-1}$ & $\cdots$ \\\cline{4-7}
&  &   &  \multirow{4}{*}{$t^{n}_{0}$} & \multirow{4}{*}{$\cdots$} & $t^{n}_{0}$ & $\cdots$ \\\cline{6-7}
&  &   &   &  & $\vdots$ & $\ddots$ \\\cline{6-7}
&  &   &   &  & $t^{n}_{n-1}$ & $\cdots$ \\\cline{4-7}
&  &   & $\vdots$  & $\ddots$ & $\vdots$ & $\ddots$ \\\cline{4-7}
&  &   &  \multirow{4}{*}{$t^{n}_{n-1}$} & \multirow{4}{*}{$\cdots$} & $t^{n}_{0}$ & $\cdots$ \\\cline{6-7}
&  &   &   &  & $\vdots$ & $\ddots$ \\\cline{6-7}
&  &   &   & & $t^{n}_{n-1}$ & $\cdots$ \\\hline
$\ddots$ & $\vdots$ & $\ddots$ & $\vdots$ & $\ddots$ & $\vdots$ & $\ddots$ \\\hline
\end{tabular}
\end{center}

\item If $\varphi_{q+1}=\neg(\varphi_{i}\wedge\neg\varphi_{i})$ and $\varphi_{i}$ takes the value $t^{n}_{l}$, for any $1\leq l\leq n-1$, the row is not branched any further in the column headed by $\varphi_{q+1}$ and assumes the value $t^{n}_{l-1}$ on said column (look at Figure \ref{Eliminate rows neg wedge neg}).
\end{enumerate}

Notice all of these steps are easily algorithmically performed: after all, they boil down to identifying the form of a formula (namely $\varphi_{q+1}$) and searching the value taken by another formula ($\varphi_{i}$) on a given row; the intuition behind these procedures involves the fact that each row of a completed row-branching, row-eliminating truth-table should correspond to a possible homomorphism, and the previous three steps eliminate all homomorphisms which are not restricted.

If we are not in the cases described above, then one proceeds as if we were in an Nmatrix: assume that, on a fixed row, $\varphi_{i}$ takes the value $a$ and $\varphi_{j}$ takes the value $b$. Then, if $\varphi_{q+1}=\varphi_{i}\#\varphi_{j}$, for $\#\in\{\vee, \wedge, \rightarrow\}$, we divide our row into $|a\tilde{\#}b|$ new rows and assign each of the elements of the set $a\tilde{\#}b$ to one of these new rows, where $\tilde{\#}$ is the operation corresponding to the connective $\#$ in the multialgebra $\mathcal{A}_{C_{n}}$; notice that $|a\tilde{\#}b|$ may equal either $n+1$, when $a\tilde{\#}b=D_{n}=\{T_{n}, t^{n}_{0}, \dotsc  , t^{n}_{n-1}\}$, or $1$, when $a\tilde{\#}b=\{T_{n}\}$ or $a\tilde{\#}b=\{F_{n}\}$. If $\varphi_{q+1}=\neg\varphi_{i}$, the row is subdivided $|\tilde{\neg}a|$ times, and one element of $\tilde{\neg}a$ is placed on each of these; again we can have $|\tilde{\neg}a|=n+1$, when $\tilde{\neg}a=D_{n}$, or $|\tilde{\neg}a|=1$, when $\tilde{\neg}a=\{T_{n}\}$ or $\tilde{\neg}a=\{F_{n}\}$. For completeness sake, compare what the previous rules dictated against what would be expected if we were working with a mere Nmatrix.

\begin{figure}[H]
\begin{center}
\begin{tabular}{|p{1cm}|p{1cm}|p{1cm}|p{1cm}|p{1cm}|p{1.5cm}|p{1cm}|p{2cm}|p{1cm}|}
\hline
$\cdots$ & $\varphi_{i}$ & $\cdots$ & $\neg\varphi_{i}$ & $\cdots$ & $\varphi_{i}\wedge\neg\varphi_{i}$ & $\cdots$ & $\neg(\varphi_{i}\wedge\neg\varphi_{i})$ & $\cdots$\\\hline
$\ddots$ & $\vdots$ & $\ddots$ & $\vdots$ & $\ddots$ & $\vdots$ & $\ddots$ & $\vdots$ & $\ddots$\\\hline
\multirow{13}{*}{$\cdots$} & \multirow{13}{*}{$t^{n}_{l}$} & \multirow{13}{*}{$\cdots$} & \multirow{4}{*}{$T_{n}$} & \multirow{4}{*}{$\cdots$} & $t^{n}_{0}$ & $\cdots$ & $t^{n}_{l-1}$ & $\cdots$\\\cline{6-9}
&  &   &   &  & $\vdots$ & $\ddots$ &  $\vdots$ & $\ddots$\\\cline{6-9}
&  &   &   & & $t^{n}_{n-1}$ & $\cdots$ &  $t^{n}_{l-1}$ & $\cdots$\\\cline{4-9}
&  &   &  \multirow{4}{*}{$t^{n}_{0}$} & \multirow{4}{*}{$\cdots$} & $t^{n}_{0}$ & $\cdots$ &  $t^{n}_{l-1}$ & $\cdots$\\\cline{6-9}
&  &   &   &  & $\vdots$ & $\ddots$ &  $\vdots$ & $\ddots$\\\cline{6-9}
&  &   &   &  & $t^{n}_{n-1}$ & $\cdots$ &  $t^{n}_{l-1}$ & $\cdots$\\\cline{4-9}
&  &   & $\vdots$  & $\ddots$ & $\vdots$ & $\ddots$ & $\vdots$ & $\ddots$\\\cline{4-9}
&  &   &  \multirow{4}{*}{$t^{n}_{n-1}$} & \multirow{4}{*}{$\cdots$} & $t^{n}_{0}$ & $\cdots$ &  $t^{n}_{l-1}$ & $\cdots$\\\cline{6-9}
&  &   &   &  & $\vdots$ & $\ddots$ &  $\vdots$ & $\ddots$\\\cline{6-9}
&  &   &   & & $t^{n}_{n-1}$ & $\cdots$ &  $t^{n}_{l-1}$ & $\cdots$\\\hline
$\ddots$ & $\vdots$ & $\ddots$ & $\vdots$ & $\ddots$ & $\vdots$ & $\ddots$ & $\vdots$ & $\ddots$\\\hline
\end{tabular}
\end{center}
\caption{Eliminating row in the case of a formula $\neg(\varphi_{i}\wedge\neg\varphi_{i})$}
\label{Eliminate rows neg wedge neg}
\end{figure}

\begin{enumerate}
\item If $\varphi_{q+1}=\varphi_{i}\wedge\neg\varphi_{i}$ and $\varphi_{i}$ is $t^{n}_{0}$, $\tilde{\neg}t^{n}_{0}=D_{n}$ and so, for any value $b\in \tilde{\neg}t^{n}_{0}$ that could be given to $\neg\varphi_{i}$, $t^{n}_{0}\tilde{\wedge}b=D_{n}$; so, only looking at $\mathcal{A}_{C_{n}}$ would branch the row $n+1$ times, giving the values $T_{n}, t^{n}_{0}, \dotsc  , t^{n}_{n-1}$ to $\varphi_{q+1}$, instead of only $T_{n}$.
\item If $\varphi_{q+1}=\varphi_{i}\wedge\neg\varphi_{i}$ and $\varphi_{i}$ is $t^{n}_{l}$, for $1\leq l\leq n-1$, one has $\tilde{\neg}t^{n}_{l}=D_{n}$ and, for any value $b\in\tilde{\neg}t^{n}_{l}$ the formula $\neg\varphi_{i}$ could assume, $t^{n}_{l}\tilde{\wedge}b=D_{n}$; this means that only considering $\mathcal{A}_{C_{n}}$ would divide the row $n+1$ times, instead of $n$ times, and give $\varphi_{q+1}$ the values $T_{n}, t^{n}_{0}, \dotsc  , t^{n}_{n-1}$, instead of $t^{n}_{0}, \dotsc  , t^{n}_{n-1}$.
\item If $\varphi_{q+1}=\neg(\varphi_{i}\wedge\neg\varphi_{i})$ and $\varphi_{i}$ is $t^{n}_{l}$, for $1\leq l\leq n-1$, we obtain $\tilde{\neg}t^{n}_{l}=D_{n}$, for any value $b\in\tilde{\neg}t^{n}_{l}$ that $\neg\varphi_{i}$ could assume $t^{n}_{l}\tilde{\wedge}b=D_{n}$ and, for any value $c\in t^{n}_{l}\tilde{\wedge} b$ which $\varphi_{i}\wedge\neg\varphi_{i}$ could assume, $\tilde{\neg}c=\{F_{n}\}$ or $\tilde{\neg}c=D_{n}$; this means $\varphi_{q+1}$ could, in principle, be assigned any value of $B_{n}$, and therefore subdivide its row into $n+2$ new ones. Instead, we demand that it does not branch its row any further and receives the value $t^{n}_{l-1}$. 
\end{enumerate}

So, we are now in position to proceed writing our row-branching, row-eliminating truth-table as long as there are elements on the sequence $\varphi_{1}, \dotsc  , \varphi_{k}=\varphi$, giving us a table with $k$ columns and no more than $(n+2)^{k}$ rows. Of course, if at the end the column headed by $\varphi_{k}$ is filled with nothing but designated values $D_{n}=\{T_{n}, t^{n}_{0}, \dotsc  , t^{n}_{n-1}\}$, then $\varphi$ is a tautology of $C_{n}$, being the converse also true: if $\varphi$ is a tautology, by writing down its correspondent row-branching, row-eliminating truth-table, we obtain the column of $\varphi$ has only designated values. We want, from here on forward, to prove that this is, in fact, the case, that is, that our truth-tables are sound and complete.

\begin{proposition}\label{Rows of tables are homomorphisms}
Let $\Gamma$ be a finite set of formulas of $C_{n}$ closed by subformulas, that is, if $\alpha\in\Gamma$ and $\beta$ is a subformula of $\alpha$, then $\beta\in\Gamma$. Let $\nu:\Gamma\rightarrow B_{n}$ be a function satisfying:
\begin{enumerate}
\item if $\neg\alpha\in\Gamma$, $\nu(\neg\alpha)\in\tilde{\neg}\nu(\alpha)$;
\item if $\alpha\#\beta\in\Gamma$, for $\#\in\{\vee, \wedge, \rightarrow\}$, then $\nu(\alpha\#\beta)\in\nu(\alpha)\tilde{\#}\nu(\beta)$;
\item if $\alpha\wedge\neg\alpha\in\Gamma$ and $\nu(\alpha)=t^{n}_{0}$, $\nu(\alpha\wedge\neg\alpha)=T_{n}$;
\item if $\alpha\wedge\neg\alpha\in\Gamma$ and $\nu(\alpha)=t^{n}_{k}$, for some $1\leq k\leq n-1$, then $\nu(\alpha\wedge\neg\alpha)\in I_{n}=\{t^{n}_{0}, \dotsc  , t^{n}_{n-1}\}$;
\item if $\alpha^{1}=\neg(\alpha\wedge\neg\alpha)\in\Gamma$ and $\nu(\alpha)=t^{n}_{k}$, for some $1\leq k\leq n-1$, then $\nu(\alpha^{1})=t^{n}_{k-1}$.
\end{enumerate}
Then there exists a homomorphism $\overline{\nu}\in\mathcal{F}_{C_{n}}$ extending $\nu$, \textit{id est}, a homomorphism for which $\overline{\nu}(\alpha)=\nu(\alpha)$, for every $\alpha\in\Gamma$.
\end{proposition}

\begin{proof}
We define $\overline{\nu}$ by induction on the order of a formula, but whenever we define $\overline{\nu}(\alpha)$ we will, simultaneously, define $\overline{\nu}(\neg\alpha)$, $\overline{\nu}(\alpha\wedge\neg\alpha)$ and $\overline{\nu}(\alpha^{1})=\overline{\nu}(\neg(\alpha\wedge\neg\alpha))$ as well.

First, notice that, if $\alpha\not\in\Gamma$, then for no formula $\beta\in\Gamma$ we have that $\alpha$ is a subformula of $\beta$, given $\Gamma$ is closed by subformulas; of course this means that if $\alpha\not\in \Gamma$, then neither $\neg\alpha$, $\alpha\wedge\neg\alpha$ or $\alpha^{1}$ are in $\Gamma$. This is important since, whenever we define $\overline{\nu}(\alpha)$ for a formula $\alpha\not\in \Gamma$, we need not worry if this definition may interfere with $\nu(\beta)$ for some $\beta\in\Gamma$ of which $\alpha$ is a subformula: such a $\beta$ can not exist.

For a propositional variable $p$, $\overline{\nu}(p)$ is defined as $\nu(p)$, whenever $p\in\Gamma$, and arbitrarily otherwise; we have then defined $\overline{\nu}(\alpha)$, for every formula $\alpha$ of order $0$, given $C_{n}$ has no nullary connectives, so it remains to define $\overline{\nu}(\neg\alpha)$, $\overline{\nu}(\alpha\wedge\neg\alpha)$ and $\overline{\nu}(\alpha^{1})$. 
\begin{enumerate}
\item If $\neg p\in \Gamma$, $\overline{\nu}(\neg p)$ is equal to $\nu(\neg p)$, and if $\neg p\not\in \Gamma$ we may take $\overline{\nu}(\neg p)$ as equal to any value of the set $\tilde{\neg}\overline{\nu}(p)$.
\item If $p\wedge\neg p\in\Gamma$, we define $\overline{\nu}(p\wedge\neg p)$ as $\nu(p\wedge\neg p)$; if $p\wedge\neg p\not\in\Gamma$ and $\overline{\nu}(p)=t^{n}_{0}$, we must define $\overline{\nu}(p\wedge\neg p)$ as $T_{n}$, but if $p\wedge\neg p\not\in \Gamma$ and $\overline{\nu}(p)=t^{n}_{k}$, for some $1\leq k\leq n-1$, we can make $\overline{\nu}(p\wedge\neg p)$ equal to an arbitrary value in $I_{n}=\{t^{n}_{0}, \dotsc   , t^{n}_{n-1}\}$; if none of the previous cases apply, we define $\overline{\nu}(p\wedge \neg p)$ as an arbitrary element of $\overline{\nu}(p)\tilde{\wedge}\overline{\nu}(\neg p)$.
\item If $p^{1}=\neg(p\wedge\neg p)\in \Gamma$, we make $\overline{\nu}(p^{1})=\nu(p^{1})$; if $p^{1}\not\in\Gamma$ and $\overline{\nu}(p)=t^{n}_{k}$, for some $1\leq k\leq n-1$, we must make $\overline{\nu}(p^{1})$ equal to $t^{n}_{k-1}$; if none of the previous cases apply, we choose $\overline{\nu}(p^{1})$ as any element of $\tilde{\neg}\overline{\nu}(p\wedge \neg p)$.
\end{enumerate}
This finishes the base step of order $0$.

Now we suppose that, for every formula $\alpha$ of order at most $m$, for a $m\geq 0$, all of $\overline{\nu}(\alpha)$, $\overline{\nu}(\neg\alpha)$, $\overline{\nu}(\alpha\wedge\neg\alpha)$ and $\overline{\nu}(\alpha^{1})$ are defined and satisfy the expected properties that would make of $\overline{\nu}$ a homomorphism, our induction hypothesis, and take a formula $\alpha$ of order $m+1$; there are then two cases to consider.
\begin{enumerate}
\item If $\alpha=\neg\beta$, $\overline{\nu}(\alpha)=\overline{\nu}(\neg\beta)$ has been already defined and satisfies: 
\begin{enumerate}
\item $\overline{\nu}(\neg\beta)=\nu(\neg\beta)$, if $\neg\beta\in\Gamma$; 
\item if $\neg\beta\not\in\Gamma$ and $\beta=\gamma\wedge\neg\gamma$, $\overline{\nu}(\gamma^{1})$ equals $t^{n}_{k-1}$, in the case that $\overline{\nu}(\gamma)=t^{n}_{k}$ (for a $1\leq k\leq n-1$), and any value in $\tilde{\neg}\overline{\nu}(\gamma\wedge\neg\gamma)$ otherwise; 
\item and, if none of these is the case, one still has $\overline{\nu}(\neg\beta)\in\tilde{\neg}\overline{\nu}(\beta)$.
\end{enumerate}
Now we define $\overline{\nu}(\neg\alpha)$, $\overline{\nu}(\alpha\wedge\neg\alpha)$ and $\overline{\nu}(\alpha^{1})$.
\begin{enumerate}
\item If $\neg\alpha\in \Gamma$, we make $\overline{\nu}(\neg\alpha)=\nu(\neg\alpha)$, otherwise is enough to require that $\overline{\nu}(\neg\alpha)\in \tilde{\neg}\overline{\nu}(\alpha)$. 

\item If $\alpha\wedge\neg\alpha\in\Gamma$, $\overline{\nu}(\alpha\wedge\neg\alpha)=\nu(\alpha\wedge\neg\alpha)$; if $\alpha\wedge\neg\alpha\not\in\Gamma$ and $\overline{\nu}(\alpha)=t^{n}_{0}$, $\overline{\nu}(\alpha\wedge\neg\alpha)=T_{n}$, and if $\alpha\wedge\neg\alpha\not\in\Gamma$ and $\overline{\nu}(\alpha)=t^{n}_{k}$ (for $1\leq k\leq n-1$), $\overline{\nu}(\alpha\wedge\neg\alpha)\in I_{n}$; finally, if none of these applies, $\overline{\nu}(\alpha\wedge\neg\alpha)\in \overline{\nu}(\alpha)\tilde{\wedge}\overline{\nu}(\neg\alpha)$.

\item If $\alpha^{1}\in\Gamma$, obviously $\overline{\nu}(\alpha^{1})=\nu(\alpha^{1})$; if $\alpha^{1}\not\in\Gamma$ and $\overline{\nu}(\alpha)=t^{n}_{k}$ (for $1\leq k\leq n-1$), $\overline{\nu}(\alpha^{1})=t^{n}_{k-1}$; and if $\alpha^{1}\not\in\Gamma$ and $\overline{\nu}(\alpha)\neq t^{n}_{k}$, we define arbitrarily $\overline{\nu}(\alpha^{1})\in\tilde{\neg}\overline{\nu}(\alpha\wedge\neg\alpha)$.
\end{enumerate}

\item If $\alpha=\beta\#\gamma$, for some $\#\in\{\vee, \wedge, \rightarrow\}$, we have that either: 
\begin{enumerate}
\item $\beta\#\gamma\in\Gamma$, when we simply define $\overline{\nu}(\beta\#\gamma)$ as $\nu(\beta\#\gamma)$;
\item $\beta\#\gamma\not\in\Gamma$ but $\#=\wedge$ and $\gamma=\neg\beta$, meaning $\overline{\nu}(\beta\wedge\neg\beta)$ has already been defined and is equal to $\nu(\beta\wedge\neg\beta)$, if $\beta\wedge\neg\beta\in\Gamma$; or is equal to $T_{n}$, if $\overline{\nu}(\beta)=t^{n}_{0}$; or lies in $I_{n}$, if $\overline{\nu}(\beta)=t^{n}_{k}$ (for $1\leq k\leq n-1$); or lies in $\overline{\nu}(\beta)\tilde{\wedge}\overline{\nu}(\neg\beta)$, if none of the previous cases applies;
\item $\beta\#\gamma\not\in\Gamma$, and $\#\neq\wedge$ or $\gamma\neq\neg\beta$, when we can define $\overline{\nu}(\beta\#\gamma)$ as any value of $\overline{\nu}(\beta)\tilde{\#}\overline{\nu}(\gamma)$.
\end{enumerate}
We are left now with the task of defining $\overline{\nu}(\neg\alpha)$, $\overline{\nu}(\alpha\wedge\neg\alpha)$ and $\overline{\nu}(\alpha^{1})$.
\begin{enumerate}
\item If $\neg(\beta\#\gamma)\in\Gamma$, we make $\overline{\nu}(\neg(\beta\#\gamma))$ equal to $\nu(\neg(\beta\#\gamma))$; if $\neg(\beta\#\gamma)\not\in\Gamma$, $\#=\wedge$, $\gamma=\neg\beta$, and $\overline{\nu}(\beta)=t^{n}_{k}$ (for $1\leq k\leq n-1$), we have already defined $\overline{\nu}(\neg(\beta\#\gamma))=\overline{\nu}(\beta^{1})$ to be equal to $t^{n}_{k-1}$; if we are in none of these cases, is sufficient to take $\overline{\nu}(\neg(\beta\#\gamma))$ as equal to any value of $\tilde{\neg}\overline{\nu}(\beta\#\gamma)$.

\item If $\alpha\wedge\neg\alpha\in\Gamma$, we merely define $\overline{\nu}(\alpha\wedge\neg\alpha)=\nu(\alpha\wedge\neg\alpha)$; if $\alpha\wedge\neg\alpha\notin\Gamma$ and $\overline{\nu}(\alpha)=t^{n}_{0}$, $\overline{\nu}(\alpha\wedge\neg\alpha)=T_{n}$, and if $\alpha\wedge\neg\alpha\not\in\Gamma$ and $\overline{\nu}(\alpha)=t^{n}_{k}$, for $1\leq k\leq n-1$, $\overline{\nu}(\alpha\wedge\neg\alpha)$ may take any value on the set $I_{n}$; if none of these applies, we can assign any value of $\overline{\nu}(\alpha)\tilde{\wedge}\overline{\nu}(\neg\alpha)$ to $\overline{\nu}(\alpha\wedge\neg\alpha)$.

\item Finally, if $\alpha^{1}\in\Gamma$, $\overline{\nu}(\alpha^{1})=\nu(\alpha^{1})$; if $\alpha^{1}\not\in\Gamma$ and $\overline{\nu}(\alpha)=t^{n}_{k}$, for any $1\leq k\leq n-1$, $\overline{\nu}(\alpha^{1})=t^{n}_{k-1}$; and if $\alpha^{1}\not\in\Gamma$ and $\overline{\nu}(\alpha)\neq t^{n}_{k}$ (for $1\leq k\leq n-1$), we simply make $\overline{\nu}(\alpha^{1})$ equal to any value of $\tilde{\neg}\overline{\nu}(\alpha\wedge\neg\alpha)$.
\end{enumerate}
\end{enumerate}
\end{proof}

So, given a formula $\varphi$ of $C_{n}$, construct its row-branching, row-eliminating truth-table as described before: the resulting table is finite and can be constructed in a finite number of steps; take $\Gamma$ as the set formed by $\varphi$ and all its subformulas and it is clear that it is closed by subformulas. By the procedure we used for constructing the table for $\varphi$, its columns correspond to the elements of $\Gamma$, while its rows correspond to functions $\nu:\Gamma\rightarrow B_{n}$ satisfying the hypothesis of Proposition \ref{Rows of tables are homomorphisms}. Of course, for every row and its $\nu$, by Proposition \ref{Rows of tables are homomorphisms} there exists a homomorphism $\overline{\nu}$ in $\mathcal{F}_{C_{n}}$ extending $\nu$; reciprocally, to every homomorphism $\overline{\nu}$ of $\mathcal{F}_{C_{n}}$ one can find a row in correspondence with a $\nu$ which is extended by $\overline{\nu}$ (it is sufficient to take $\nu$ as the restriction of $\overline{\nu}$ to $\Gamma$).

So, assume $\vdash_{C_{n}}\varphi$, and therefore $\vDash_{\mathcal{RM}_{C_{n}}}\varphi$; then, for any row of the table for $\varphi$ and its $\nu$, since $\overline{\nu}\in\mathcal{F}_{C_{n}}$ we have $\overline{\nu}(\varphi)\in D_{n}$, and therefore $\nu(\varphi)\in D_{n}$, meaning all entries on the column headed by $\varphi$ are designated and, therefore, $\varphi$ is proved to be a tautology according to row-branching, row-eliminating truth-tables. Reciprocally, assume $\varphi$ is proved according to row-branching, row-eliminating truth-tables: for any homomorphism $\overline{\nu}\in\mathcal{F}_{C_{n}}$, $\nu=\overline{\nu}|_{\Gamma}$ corresponds to a row of the table for $\varphi$ and therefore satisfies $\nu(\varphi)\in D_{n}$, meaning therefore that $\overline{\nu}(\varphi)\in D_{n}$; this implies $\vDash_{\mathcal{RM}_{C_{n}}}\varphi$, and so $\vdash_{C_{n}}\varphi$, what gives us soundness and completeness.

We can, in much the same way, define row-branching, row-eliminating truth-tables for $\textbf{mbCcl}$ and $\textbf{Cila}$ based on the RNmatrices $\mathcal{M}_{\textbf{mbCcl}}$ (see Section \ref{MmbCcl}) and $\mathcal{M}_{\textbf{CILA}}$ (see Section \ref{MCILA}); in both cases tables are proven to be decision methods for these logics, which are known to be not characterizable by finite Nmatrices yet are decidable by three-valued RNmatrices.

Now, we present some examples; first of all, we stress a slight simplification: in the case that $\Gamma=\{\gamma_{1}, \dotsc  , \gamma_{m}\}$ is a finite set of premises, instead of writing down the table for $\gamma_{1}\rightarrow(\cdots\rightarrow(\gamma_{m}\rightarrow\varphi)\cdots)$ to verify the deduction $\Gamma\vdash_{C_{n}}\varphi$, is enough to list the columns headed by the formulas $\varphi$ and $\gamma_{1}$ through $\gamma_{m}$, and their subformulas and verify that, in every row where all premises assume designated values, so does $\varphi$. This is due to the fact that all $C_{n}$ obey the deduction meta-theorem, meaning that $\Gamma, \psi\vdash_{C_{n}}\varphi$ if, and only if, $\Gamma\vdash_{C_{n}}\psi\rightarrow\varphi$, and this may shorten significantly our truth-tables.

Take a propositional variable $p$ and the deduction $p, \neg p, \neg(p\wedge\neg p)\vdash_{C_{1}}\neg\neg p$. We obtain then the following series of tables that culminate in the completed row-branching, row-eliminating truth-table for that very argument.

\begin{figure}[H]
\centering
\begin{minipage}[t]{2cm}
\centering
\begin{tabular}{|p{1cm}|p{1cm}|}
\hline
$p$ & Row\\ \hline
$T_{1}$ & \nth{1}\\ \hline
$t^{1}_{0}$ & \nth{2}\\ \hline
$F_{1}$ & \nth{3}\\ \hline
\end{tabular}
\caption*{First stage}
\end{minipage}
\hspace{3cm}
\begin{minipage}[t]{3cm}
\centering
\begin{tabular}{|p{1cm}|p{1cm}|p{1cm}|}
\hline
$p$ & $\neg p$ & Row\\ \hline
$T_{1}$ & $F_{1}$ & \nth{1}\\ \hline
\multirow{2}{*}{$t^{1}_{0}$} & $T_{1}$ & \nth{2}\\ \cline{2-3}
& $t^{1}_{0}$ & \nth{3}\\ \hline
$F_{1}$ & $T_{1}$ & \nth{4}\\ \hline
\end{tabular}
\caption*{Second stage}
\end{minipage}
\end{figure}

\begin{figure}[H]
\centering
\begin{minipage}[t]{4cm}
\centering
\begin{tabular}{|p{1cm}|p{1cm}|p{1cm}|p{1cm}|}
\hline
$p$ & $\neg p$ & $\neg\neg p$ & Row\\ \hline
$T_{1}$ & $F_{1}$ & $T_{1}$ & \nth{1}\\ \hline
\multirow{3}{*}{$t^{1}_{0}$} & $T_{1}$ & $F_{1}$ & \nth{2}\\ \cline{2-4}
& \multirow{2}{*}{$t^{1}_{0}$} & $T_{1}$ & \nth{3}\\ \cline{3-4}
& & $t^{1}_{0}$ & \nth{4}\\ \hline
$F_{1}$ & $T_{1}$ & $F_{1}$ & \nth{5}\\ \hline
\end{tabular}
\caption*{Third stage}
\end{minipage}
\hspace{2.5cm}
\begin{minipage}[t]{5cm}
\centering
\begin{tabular}{|p{1cm}|p{1cm}|p{1cm}|p{1cm}|p{1cm}|}
\hline
$p$ & $\neg p$ & $\neg\neg p$ & $p\wedge\neg p$ & Row\\ \hline
$T_{1}$ & $F_{1}$ & $T_{1}$ & $F_{1}$ & \nth{1}\\ \hline
\multirow{3}{*}{$t^{1}_{0}$} & $T_{1}$ & $F_{1}$ & $T_{1}$ & \nth{2}\\ \cline{2-5}
& \multirow{2}{*}{$t^{1}_{0}$} & $T_{1}$ & $T_{1}$ & \nth{3}\\ \cline{3-5}
& & $t^{1}_{0}$ & $T_{1}$ & \nth{4}\\ \hline
$F_{1}$ & $T_{1}$ & $F_{1}$ & $F_{1}$ & \nth{5}\\ \hline
\end{tabular}
\caption*{Fourth stage}
\end{minipage}
\end{figure}

\begin{figure}[H]
\centering
\begin{minipage}[t]{6.5cm}
\centering
\begin{tabular}{|p{1cm}|p{1cm}|p{1cm}|p{1cm}|p{1.5cm}|p{1cm}|}
\hline
$p$ & $\neg p$ & $\neg\neg p$ & $p\wedge\neg p$ & $\neg(p\wedge\neg p)$ & Row\\\hline
$T_{1}$ & $F_{1}$ & $T_{1}$ & $F_{1}$ & $T_{1}$ & \nth{1}\\\hline
\multirow{3}{*}{$t^{1}_{0}$} & $T_{1}$ & $F_{1}$ & $T_{1}$ & $F_{1}$ & \nth{2}\\\cline{2-6}
& \multirow{2}{*}{$t^{1}_{0}$} & $T_{1}$ & $T_{1}$ & $F_{1}$ & \nth{3}\\\cline{3-6}
& & $t^{1}_{0}$ & $T_{1}$ & $F_{1}$ & \nth{4}\\\hline
$F_{1}$ & $T_{1}$ & $F_{1}$ & $F_{1}$ & $T_{1}$ & \nth{5}\\\hline
\end{tabular}
\caption*{Fifth and final stage}
\end{minipage}
\end{figure}

We see that the deduction $p, \neg p, \neg(p\wedge\neg p)\vdash_{C_{1}}\neg\neg p$ is valid: since in no rows one has $p$, $\neg p$ and $\neg(p\wedge\neg p)$ simultaneously assigned designated values, the deduction is vacuously true.

Now, a final remark regarding our truth-tables: while we have presented row-bran\-ching, row-eliminating truth-tables on which rows are discarded as soon as they violate some condition necessary for the corresponding homomorphism to lie in $\mathcal{F}_{C_{n}}$, one could first write down the row-branching truth-table of the Nmatrix $(\mathcal{A}_{C_{n}}, D_{n})$ to only then delete the rows corresponding to homomorphism outside $\mathcal{F}_{C_{n}}$. The order in which this is done is not relevant: although eliminating the rows as soon as they are discovered not to be relevant is slightly more efficient, the alternative way is, arguably, more intuitive. To show one example of how this goes, consider again the deduction $p, \neg p, \neg(p\wedge\neg p)\vdash_{C_{1}}\neg\neg p$; using only the Nmatrix $(\mathcal{A}_{C_{1}}, D_{1})$ one obtains the following table.

\begin{figure}[H]
\centering
\begin{minipage}[t]{6.5cm}
\centering
\begin{center}
\begin{tabular}{|p{1cm}|p{1cm}|p{1cm}|p{1cm}|p{1.5cm}|c|}
\hline
$p$ & $\neg p$ & $\neg\neg p$ & $p\wedge\neg p$ & $\neg(p\wedge\neg p)$ & Row\\\hline
$T_{1}$ & $F_{1}$ & $T_{1}$ & $F_{1}$ & $T_{1}$ & \nth{1}\\\hline
\multirow{9}{*}{$t^{1}_{0}$} & \multirow{3}{*}{$T_{1}$} & \multirow{3}{*}{$F_{1}$} & $T_{1}$ & $F_{1}$ & \nth{2}\\\cline{4-6}
& & & \multirow{2}{*}{$t^{1}_{0}$}& $T_{1}$ & \nth{3}\\\cline{5-6}
& & & & $t^{1}_{0}$ & \nth{4}\\\cline{2-6}
& \multirow{6}{*}{$t^{1}_{0}$} & \multirow{3}{*}{$T_{1}$} & $T_{1}$ & $F_{1}$ & \nth{5}\\\cline{4-6}
& & & \multirow{2}{*}{$t^{1}_{0}$} & $T_{1}$ & \nth{6}\\\cline{5-6}
& & & & $t^{1}_{0}$ & \nth{7}\\\cline{3-6}
& & \multirow{3}{*}{$t^{1}_{0}$} & $T_{1}$ & $F_{1}$ & \nth{8}\\\cline{4-6}
& & & \multirow{2}{*}{$t^{1}_{0}$} & $T_{1}$ & \nth{9}\\\cline{5-6}
& & & & $t^{1}_{0}$ & \nth{10}\\\hline
$F_{1}$ & $T_{1}$ & $F_{1}$ & $F_{1}$ & $T_{1}$ & \nth{11}\\\hline
\end{tabular}
\end{center}
\end{minipage}
\end{figure}

If one tries to evaluate the validity of the argument only from this table, it appears false: rows $3$ and $4$ have all premises designated but not the conclusion. But the rows, in the general case of $C_{n}$, to be eliminated are:
\begin{enumerate}
\item the ones where $\alpha$ is assigned $t^{n}_{0}$ and $\alpha\wedge\neg\alpha$ is not assigned $T_{n}$;
\item the ones where $\alpha$ is assigned the value $t^{n}_{k}$, for $1\leq k\leq n-1$, but $\alpha\wedge\neg\alpha$ is assigned the value $T_{n}$;
\item the ones where $\alpha$ is $t^{n}_{k}$, for $1\leq k\leq n-1$, but $\alpha^{1}$ is not $t^{n}_{k-1}$.
\end{enumerate}
In our specific case, rows $3$, $4$, $6$, $7$, $9$ and $10$ are eliminated, and we retrieve the previous table we had written for the deduction, which makes it valid.

\subsection{Tableau semantics}\label{Tableaux semantics for Cn}

We have presented, in Sections \ref{RNmatrix for mbCcl}, \ref{RNmatrix for CILA} and \ref{The general case} shown RNmatrices that, through the row-branching, row-eliminating truth tables outlined in Section \ref{Row-branching,row-eliminating}, become decision methods for the whole hierarchy of da Costa plus $\textbf{Cila}$ and $\textbf{mbCcl}$, although none of these logics is characterized by finite Nmatrices. However the corresponding truth tables can grow rather rapidly and implementing such a method may not be at all feasible.

So we offer a second decision procedure, based too on the RNmatrices $\mathcal{RM}_{C_{n}}$, by means of tableau semantics\index{Tableaux}. Although recent developments provided general approaches to obtaining tableau-like proof systems from finite Nmatrices (\cite{Pawlowski}), the present study will make use more extensively of \cite{con:far:per:21}, which appears to be more useful when dealing with RNmatrices defined by swap structures. Many of the ideas used to motivate tableaux from RNmatrices may be found in any classical account on the subject, e.g. \cite{Smullyan}.

Our tableaux will have as nodes labeled formulas, that is, pairs of labels $\textsf{L}\in \textsf{B}_{n}=\{\textsf{T}_{n}, \textsf{t}^{n}_{0}, \dotsc  , \textsf{t}^{n}_{n-1}, \textsf{F}_{n}\}$ and formulas $\varphi$ of $C_{n}$. The labels correspond to the snapshots of $B_{n}$, but to the make the distinction between truth value and label clearer we use different fonts; furthermore, the pair formed by $\textsf{L}$ and $\varphi$ will be denoted by $\textsf{L}(\varphi)$. The labeled tableau system for $C_{n}$ will be denoted by $\mathbb{T}_{n}$\label{mathbbTn} and have the following $n+2$ rules of elimination for every connective (that is, $\neg$, $\wedge$, $\vee$ and $\rightarrow$, totaling $4n+8$ rules), where the pictorial rule
\[\begin{array}{c|c|c}\multicolumn{3}{c}{\textsf{L}(\varphi)}\\\hline\textsf{L}_{11}(\varphi_{11}) & \multirow{3}{*}{$\cdots$} & \textsf{L}_{m1}(\varphi_{m1})\\ \vdots& & \vdots\\ \textsf{L}_{1n_{1}}(\varphi_{1n_{1}})& & \textsf{L}_{mn_{m}}(\varphi_{mn_{m}})\end{array}\]
indicates that a branch of a tableau containing $\textsf{L}(\varphi)$ may be forked into $m$ new branches, each containing all formulas already present in the original branch and, additionally, $\textsf{L}_{k1}(\varphi_{k1})$ trough $\textsf{L}_{kn_{k}}(\varphi_{kn_{k}})$, for some $1\leq k\leq m$.

\begin{table}[H]
\centering
$\begin{array}{c}

\begin{array}{c|c|c|c|c|c|c|c}\multicolumn{8}{c}{\textsf{T}_{n}(\varphi\vee\psi)} \\\hline \textsf{T}_{n}(\varphi) & \textsf{t}^{n}_{0}(\varphi) & \cdots & \textsf{t}^{n}_{n-1}(\varphi) & \textsf{T}_{n}(\psi) & \textsf{t}^{n}_{n-1}(\psi) & \cdots & \textsf{t}^{n}_{n-1}(\psi)\end{array}\quad (E_{\vee}\textsf{T}_{n}) \\
\\
\begin{array}{c|c|c|c|c|c}\multicolumn{6}{c}{\textsf{t}^{n}_{0}(\varphi\vee\psi)} \\\hline \textsf{t}^{n}_{0}(\varphi) & \cdots & \textsf{t}^{n}_{n-1}(\varphi) & \textsf{t}^{n}_{0}(\psi) & \cdots & \textsf{t}^{n}_{n-1}(\psi)\end{array}\quad (E_{\vee}\textsf{t}^{n}_{0})\\
\\
\vdots
\\
\begin{array}{c|c|c|c|c|c}\multicolumn{6}{c}{\textsf{t}^{n}_{n-1}(\varphi\vee\psi)} \\\hline \textsf{t}^{n}_{0}(\varphi) & \cdots & \textsf{t}^{n}_{n-1}(\varphi) & \textsf{t}^{n}_{0}(\psi) & \cdots & \textsf{t}^{n}_{n-1}(\psi)\end{array}\quad (E_{\vee}\textsf{t}^{n}_{n-1})\\
\\

\begin{array}{c} \textsf{F}_{n}(\varphi\vee\psi)\\\hline \textsf{F}_{n}(\varphi)\\ \textsf{F}_{n}(\psi)\end{array} \quad (E_{\vee}\textsf{F}_{n})
\end{array}$
\vspace{0.5cm}
\caption*{The $n+2$ rules for eliminating disjunction}
\end{table}

\begin{table}[H]
\centering
$\begin{array}{ccccc}
\begin{array}{c|c|c|c}
\multicolumn{4}{c}{\textsf{T}_{n}(\neg\varphi)}\\ \hline \textsf{t}^{n}_{0}(\varphi) & \cdots & \textsf{t}^{n}_{n-1}(\varphi) & \textsf{F}_{n}(\varphi)
\end{array} \quad (E_{\neg}\textsf{T}_{n}) &\quad\quad\quad &
\begin{array}{c|c|c}\multicolumn{3}{c}{\textsf{t}^{n}_{0}(\neg\varphi)}\\ \hline \textsf{t}^{n}_{0}(\varphi) & \cdots & \textsf{t}^{n}_{n-1}(\varphi)\\\end{array} \quad (E_{\neg}\textsf{t}^{n}_{0}) & \quad\quad & \cdots\\
\end{array}$
\end{table}

\begin{table}[H]
\centering
$\begin{array}{ccccc}
\begin{array}{c|c|c}\multicolumn{3}{c}{\textsf{t}^{n}_{n-1}(\neg\varphi)}\\ \hline \textsf{t}^{n}_{0}(\varphi) & \cdots & \textsf{t}^{n}_{n-1}(\varphi)\\\end{array} \quad (E_{\neg}\textsf{t}^{n}_{n-1}) &\quad\quad\quad &
\begin{array}{c}\textsf{F}_{n}(\neg\varphi)\\\hline\textsf{T}_{n}(\varphi)\\\end{array} \quad (E_{\neg}\textsf{F}_{n}) &  & \\
\end{array}$
\vspace{0.5cm}
\caption*{The $n+2$ rules for eliminating negation}
\end{table}

\begin{landscape}

\begin{table}[H]
\centering
$\begin{array}{c}
\begin{array}{c|c|c|c|c|c|c|c|c|c|c|c|c}\multicolumn{13}{c}{\textsf{T}_{n}(\varphi\wedge\psi)}\\ \hline\textsf{T}_{n}(\varphi) & \textsf{T}_{n}(\varphi) & \multirow{2}{*}{$\cdots$} & \textsf{T}_{n}(\varphi) & \textsf{t}^{n}_{0}(\varphi) & \textsf{t}^{n}_{0}(\varphi) & \multirow{2}{*}{$\cdots$} & \textsf{t}^{n}_{0}(\varphi) & \multirow{2}{*}{$\cdots$} & \textsf{t}^{n}_{n-1}(\varphi) & \textsf{t}^{n}_{n-1}(\varphi) & \multirow{2}{*}{$\cdots$} & \textsf{t}^{n}_{n-1}(\varphi)\\ \textsf{T}_{n}(\psi) & \textsf{t}^{n}_{0}(\psi) &  & \textsf{t}^{n}_{n-1}(\psi) & \textsf{T}_{n}(\psi) & \textsf{t}^{n}_{0}(\psi) &  & \textsf{t}^{n}_{n-1}(\psi) &  & \textsf{T}_{n}(\psi) & \textsf{t}^{n}_{0}(\psi) &  & \textsf{t}^{n}_{n-1}(\psi) \end{array} \quad (E_{\wedge}\textsf{T}_{n}) \\
\\

\begin{array}{c|c|c|c|c|c|c|c|c|c|c|c|c}\multicolumn{13}{c}{\textsf{t}^{n}_{0}(\varphi\wedge\psi)}\\ \hline \textsf{T}_{n}(\varphi) & \multirow{2}{*}{$\cdots$} & \textsf{T}_{n}(\varphi)  & \textsf{t}^{n}_{0}(\varphi) & \multirow{2}{*}{$\cdots$} & \textsf{t}^{n}_{0}(\varphi) & \multirow{2}{*}{$\cdots$}  & \textsf{t}^{n}_{n-1}(\varphi) & \multirow{2}{*}{$\cdots$} & \textsf{t}^{n}_{n-1}(\varphi) & \textsf{t}^{n}_{0}(\varphi)& \multirow{2}{*}{$\cdots$} & \textsf{t}^{n}_{n-1}(\varphi)\\ 
\textsf{t}^{n}_{0}(\psi) &  & \textsf{t}^{n}_{n-1}(\psi) & \textsf{t}^{n}_{0}(\psi) &  & \textsf{t}^{n}_{n-1}(\psi) &  &  \textsf{t}^{n}_{0}(\psi) &  & \textsf{t}^{n}_{n-1}(\psi) & \textsf{T}_{n}(\psi) & & \textsf{T}_{n}(\psi) \end{array} \quad (E_{\wedge}\textsf{t}^{n}_{0})\\
\\
\vdots\\
\\

\begin{array}{c|c|c|c|c|c|c|c|c|c|c|c|c}\multicolumn{13}{c}{\textsf{t}^{n}_{n-1}(\varphi\wedge\psi)}\\ \hline \textsf{T}_{n}(\varphi) & \multirow{2}{*}{$\cdots$} & \textsf{T}_{n}(\varphi)  & \textsf{t}^{n}_{0}(\varphi) & \multirow{2}{*}{$\cdots$} & \textsf{t}^{n}_{0}(\varphi) & \multirow{2}{*}{$\cdots$}  & \textsf{t}^{n}_{n-1}(\varphi) & \multirow{2}{*}{$\cdots$} & \textsf{t}^{n}_{n-1}(\varphi) & \textsf{t}^{n}_{0}(\varphi)& \multirow{2}{*}{$\cdots$} & \textsf{t}^{n}_{n-1}(\varphi)\\ 
\textsf{t}^{n}_{0}(\psi) &  & \textsf{t}^{n}_{n-1}(\psi) & \textsf{t}^{n}_{0}(\psi) &  & \textsf{t}^{n}_{n-1}(\psi) &  &  \textsf{t}^{n}_{0}(\psi) &  & \textsf{t}^{n}_{n-1}(\psi) & \textsf{T}_{n}(\psi) & & \textsf{T}_{n}(\psi) \end{array} \quad (E_{\wedge}\textsf{t}^{n}_{n-1})\\
\\

\begin{array}{c|c}\multicolumn{2}{c}{\textsf{F}_{n}(\varphi\wedge\psi)}\\\hline \textsf{F}_{n}(\varphi) & \textsf{F}_{n}(\psi)\\\end{array} \quad (E_{\wedge}\textsf{F}_{n})
\end{array}$
\vspace{0.5cm}
\caption*{The $n+2$ rules for eliminating conjunction}
\end{table}

\vspace{1cm}

\begin{table}[H]
\centering
$\begin{array}{c}
\begin{array}{c|c|c|c|c}\multicolumn{5}{c}{\textsf{T}_{n}(\varphi\rightarrow\psi)} \\\hline \textsf{F}_{n}(\varphi) & \textsf{T}_{n}(\psi) & \textsf{t}^{n}_{n-1}(\psi) & \cdots & \textsf{t}^{n}_{n-1}(\psi)\end{array}\quad (E_{\rightarrow}\textsf{T}_{n}) \\
\end{array}$
\end{table}

\end{landscape}

\begin{table}[H]
\centering
$\begin{array}{c}
\begin{array}{c|c|c|c|c|c}\multicolumn{6}{c}{\textsf{t}^{n}_{0}(\varphi\rightarrow\psi)}\\\hline \textsf{t}^{n}_{0}(\varphi) & \multirow{2}{*}{$\cdots$} & \textsf{t}^{n}_{n-1}(\varphi) & \multirow{2}{*}{$\textsf{t}^{n}_{0}(\psi)$} & \multirow{2}{*}{$\cdots$} & \multirow{2}{*}{$\textsf{t}^{n}_{n-1}(\psi)$}\\ \textsf{T}_{n}(\psi) & & \textsf{T}_{n}(\psi) & & &\end{array} \quad (E_{\rightarrow}\textsf{t}^{n}_{0})\\
\\
\vdots\\
\\
\begin{array}{c|c|c|c|c|c}\multicolumn{6}{c}{\textsf{t}^{n}_{n-1}(\varphi\rightarrow\psi)}\\\hline \textsf{t}^{n}_{0}(\varphi) & \multirow{2}{*}{$\cdots$} & \textsf{t}^{n}_{n-1}(\varphi) & \multirow{2}{*}{$\textsf{t}^{n}_{0}(\psi)$} & \multirow{2}{*}{$\cdots$} & \multirow{2}{*}{$\textsf{t}^{n}_{n-1}(\psi)$}\\ \textsf{T}_{n}(\psi) & & \textsf{T}_{n}(\psi) & & &\end{array} \quad (E_{\rightarrow}\textsf{t}^{n}_{n-1})\\
\\
\begin{array}{c|c|c|c}\multicolumn{4}{c}{\textsf{F}_{n}(\varphi\rightarrow\psi)}\\\hline \textsf{T}_{n}(\varphi) & \textsf{t}^{n}_{0}(\varphi) & \multirow{2}{*}{$\cdots$} & \textsf{t}^{n}_{n-1}(\varphi)\\ \textsf{F}_{n}(\psi) & \textsf{F}_{n}(\psi) & & \textsf{F}_{n}(\psi)\end{array}\quad (E_{\rightarrow}\textsf{F}_{n})
\end{array}$
\vspace{0.5cm}
\caption*{The $n+2$ rules for eliminating implication}
\end{table}

One sees that the tableau rules for $C_{n}$ are rather difficult to write down given their sizes, so we give an alternative presentation that uses a more compact notation. Let $\textsf{X}$ and $\textsf{Y}$ be sets of labels.
\[\begin{array}{ccccc}
\begin{array}{c|c}\multicolumn{2}{c}{\textsf{L}(\varphi)}\\\hline \textsf{L}^{*}(\psi) & \textsf{X}(\psi)\end{array}&\quad\quad\quad\quad&
\begin{array}{c}\textsf{L}(\varphi)\\\hline \textsf{L}^{*}(\gamma)\\\textsf{X}(\psi)\end{array}&\quad\quad\quad\quad&
\begin{array}{c|c|c}\multicolumn{3}{c}{\textsf{L}(\varphi)}\\\hline \textsf{L}^{*}(\gamma) & \textsf{L}^{**}(\psi) & \textsf{X}(\gamma)\\ \textsf{X}(\psi) & \textsf{Y}(\gamma) & \textsf{Y}(\psi)\end{array}
\end{array}\]

\begin{enumerate}
\item The leftmost rule states that a closed branch containing $\textsf{L}(\varphi)$ will also contain $\textsf{L}^{*}(\psi)$ or $\textsf{x}(\psi)$, for some label $\textsf{x}\in\textsf{X}$; if $\textsf{X}$ has $p$ elements, such a rule forks a branch into $p+1$ new ones.

\item For the rule in the middle, a closed branch containing $\textsf{L}(\varphi)$ must also contain both $\textsf{L}^{*}(\gamma)$ and $\textsf{x}(\psi)$, for some $\textsf{x}\in\textsf{X}$. Notice that, since we do not have the vertical bar, it would appear the branch is not forked at all; but since $\textsf{X}$ is a set, we still have multiple cases to consider, being the branch forked into $p$ new branches.

\item Finally, if $\textsf{L}(\varphi)$ is in a closed branch, the rightmost rule implies one of the following lies in the branch as well: $\textsf{L}^{*}(\gamma)$ and $\textsf{x}(\psi)$, for some $\textsf{x}\in\textsf{X}$; or $\textsf{L}^{**}(\psi)$ and $\textsf{y}(\gamma)$, for some $\textsf{y}\in\textsf{Y}$;\footnote{Notice that, in this rule, we change the order of the formulas ($\gamma$ above and $\psi$ below) that was previously established in order to make the idea of a tableau as a tree growing downwards clearer, leaving the larger set of labels on the bottom.} or $\textsf{x}(\gamma)$ and $\textsf{y}(\psi)$, for $\textsf{x}\in\textsf{X}$ and $\textsf{y}\in\textsf{Y}$. If $\textsf{Y}$ has $q$ elements, this rules forks a branch $pq+p+q$ times.
\end{enumerate}

Using these conventions and defining $\textsf{I}_{n}=\{\textsf{t}^{n}_{0}, \dotsc  , \textsf{t}^{n}_{n-1}\}$ and $\textsf{D}_{n}=\textsf{I}_{n}\cup\{\textsf{T}_{n}\}$, the rules of $\mathbb{T}_{n}$ may be more compactly presented by the following, where $i$ takes any value in $\{0, 1, \dotsc  , n-1\}$.

\begin{table}[H]
\centering
$\begin{array}{ccccc}
\begin{array}{c|c}\multicolumn{2}{c}{\textsf{T}_{n}(\neg\varphi)}\\\hline\textsf{I}_{n}(\varphi) & \textsf{F}_{n}(\varphi)\end{array}\quad (E_{\neg}\textsf{T}_{n}) & \quad &
\begin{array}{c}\textsf{t}^{n}_{i}(\neg\varphi)\\\hline\textsf{I}_{n}(\varphi)\end{array} \quad (E_{\neg}\textsf{t}^{n}_{i})& \quad &
\begin{array}{c}\textsf{F}_{n}(\neg\varphi)\\\hline\textsf{T}_{n}(\varphi)\end{array}\quad (E_{\neg}\textsf{F}_{n})\\
&&&&\\
\begin{array}{c}\textsf{T}_{n}(\varphi\wedge\psi)\\\hline\textsf{D}_{n}(\varphi)\\\textsf{D}_{n}(\psi)\end{array}\quad (E_{\wedge}\textsf{T}_{n}) & & \begin{array}{c|c|c}\multicolumn{3}{c}{\textsf{t}^{n}_{i}(\varphi\wedge\psi)}\\\hline\textsf{T}_{n}(\varphi)&\textsf{I}_{n}(\varphi)&\textsf{T}_{n}(\psi)\\\textsf{I}_{n}(\psi)&\textsf{I}_{n}(\psi)&\textsf{I}_{n}(\varphi)\end{array}\quad (E_{\wedge}\textsf{t}^{n}_{i}) & & 
\begin{array}{c|c}\multicolumn{2}{c}{\textsf{F}_{n}(\varphi\wedge\psi)}\\\hline\textsf{F}_{n}(\varphi)&\textsf{F}_{n}(\psi)\end{array}\quad(E_{\wedge}\textsf{F}_{n})\\
\end{array}$
\end{table}

\begin{table}[H]
\centering
$\begin{array}{ccccc}
\begin{array}{c|c}\multicolumn{2}{c}{\textsf{T}_{n}(\varphi\vee\psi)}\\\hline\textsf{D}_{n}(\varphi)&\textsf{D}_{n}(\psi)\end{array}\quad(E_{\vee}\textsf{T}_{n}) & & 
\begin{array}{c|c}\multicolumn{2}{c}{\textsf{t}^{n}_{i}(\varphi\vee\psi)}\\\hline\textsf{I}_{n}(\varphi)&\textsf{I}_{n}(\psi)\end{array}\quad (E_{\vee}\textsf{t}^{n}_{i})& &
\begin{array}{c}\textsf{F}_{n}(\varphi\vee\psi)\\\hline\textsf{F}_{n}(\varphi)\\\textsf{F}_{n}(\psi)\end{array}\quad(E_{\vee}\textsf{F}_{n})\\
&&&&\\
\begin{array}{c|c}\multicolumn{2}{c}{\textsf{T}_{n}(\varphi\rightarrow\psi)}\\\hline\textsf{F}_{n}(\varphi)&\textsf{D}_{n}(\psi)\end{array}\quad(E_{\rightarrow}\textsf{T}_{n})& &
\begin{array}{c|c}\multicolumn{2}{c}{\textsf{t}^{n}_{i}(\varphi\rightarrow\psi)}\\\hline\textsf{T}_{n}(\psi)&\multirow{2}{*}{$\textsf{I}_{n}(\psi)$}\\\textsf{I}_{n}(\varphi)& \end{array}\quad(E_{\rightarrow}\textsf{t}^{n}_{i})& &
\begin{array}{c}\textsf{F}_{n}(\varphi\rightarrow\psi)\\\hline\textsf{F}_{n}(\psi)\\\textsf{D}_{n}(\varphi)\end{array}\quad(E_{\rightarrow}\textsf{F}_{n})
\end{array}$
\end{table}

\begin{definition}\label{Definitions about tableaux}
A branch of a tableau in $\mathbb{T}_{n}$ is closed\index{Branch, Closed} whenever:
\begin{enumerate}
\item it contains $\textsf{L}(\varphi)$ and $\textsf{L}^{*}(\varphi)$ for labels $\textsf{L}\neq\textsf{L}^{*}$;
\item it contains $\textsf{t}^{n}_{0}(\varphi)$ and $\textsf{t}^{n}_{i}(\varphi\wedge\neg\varphi)$, for $0\leq i\leq n-1$;
\item it contains $\textsf{t}^{n}_{k}(\varphi)$, and either $\textsf{T}_{n}(\varphi\wedge\neg\varphi)$ or $\textsf{L}(\varphi^{1})$, for $1\leq k\leq n-1$ and $\textsf{L}\neq \textsf{t}^{n}_{k-1}$.
\end{enumerate}

A branch $\theta$ is said to be complete\index{Branch, Complete} if, for every signed formula $\textsf{L}(\varphi)$ appearing in it with $\varphi$ not a propositional variable, $\theta$ also contains all the signed formulas of one of the branches of the only rule applicable to $\textsf{L}(\varphi)$. A complete branch is open if it is not closed.

A tableau of $\mathbb{T}_{n}$ is:
\begin{enumerate}
\item closed\index{Tableau, Closed} if all of its branches are closed; 
\item complete\index{Tableau, Complete} if all of its branches are either closed or complete;
\item open\index{Tableau, Open} if it is complete but not closed.
\end{enumerate}
\end{definition}

The rules of $\mathbb{T}_{n}$ are analytic, meaning that if $\textsf{L'}(\psi)$ appears as a consequence of $\textsf{L}(\varphi)$ in some rule of $\mathbb{T}_{n}$, the $\psi$ is a proper subformula of $\varphi$. In other words, the tableaux of $\mathbb{T}_{n}$ can always be completed in a finite number of steps: in fact, if the complexity of $\varphi$ is $m$, any completed tableau starting with $\textsf{L}(\varphi)$ is guaranteed to have at most $(n+1)^{2m}$ branches, of length at most $m$.

We set out now to prove that the method of tableaux here presented for $C_{n}$ is both sound and complete, using arguments very similar to those found in \cite{Smullyan}; combined with the fact any such tableau may be completed in finite time, we will therefore have a different decision procedure for these systems.

\begin{definition}
A formula $\varphi$ of $C_{n}$ is provable according to tableaux of $\mathbb{T}_{n}$, when we write $\vdash_{\mathbb{T}_{n}}\varphi$, if there exists a closed tableau of $\mathbb{T}_{n}$ started by $\textsf{F}_{n}(\varphi)$.

Given a finite set $\Gamma=\{\gamma_{1}, \dotsc  , \gamma_{m}\}$ of formulas of $C_{n}$, $\varphi$ is said to be provable from $\Gamma$ according to tableaux of $\mathbb{T}_{n}$, when we write $\Gamma\vdash_{\mathbb{T}_{n}}\varphi$\label{vdashmathbbTn}, if 
\[\vdash_{\mathbb{T}_{n}}\bigwedge_{i=1}^{m}\gamma_{i}\rightarrow\varphi.\]
\end{definition}

Given a valuation $\nu$ in $\mathcal{F}_{C_{n}}$, a signed formula $\textsf{L}(\varphi)$ is true under $\nu$, also said to be satisfied by $\nu$, if $\nu(\varphi)=L$, where $L$ is, of course, the element of $B_{n}$ corresponding to $\textsf{L}$; if $\textsf{L}(\varphi)$ is not true under $\nu$, we say it is false.

Then a branch $\theta$ of a tableau $\mathcal{T}$ of $\mathbb{T}_{n}$ is true under $\nu$, or is satisfied by $\nu$, if all of its signed formulas are true under $\nu$; $\mathcal{T}$ is true under $\nu$, or satisfied by $\nu$, if at least one of its branches is true under $\nu$. We can see that a closed branch $\theta$, and accordingly a closed tableau, is not satisfied under any $\nu$ in $\mathcal{F}_{C_{n}}$:
\begin{enumerate}
\item if $\theta$ contains $\textsf{L}(\varphi)$ and $\textsf{L}^{*}(\varphi)$, for $\textsf{L}\neq \textsf{L}^{*}$, $\theta$ being true under $\nu$ would mean $\nu(\varphi)=L$ and $\nu(\varphi)=L^{*}$, what is absurd;
\item if $\theta$ contains $\textsf{t}_{0}^{n}(\varphi)$ and $\textsf{t}_{i}^{n}(\varphi\wedge\neg\varphi)$, for $0\leq i\leq n-1$, and $\theta$ is true under $\nu$, from the definition of $\mathcal{F}_{C_{n}}$ one has that $\nu(\varphi)=t_{0}^{n}$ would imply $\nu(\varphi\wedge\neg\varphi)=T_{n}$, contradicting the fact that $\nu(\varphi\wedge\neg\varphi)=t_{i}^{n}$;
\item and if $\theta$ contains $\textsf{t}^{n}_{k}(\varphi)$, and either $\textsf{T}_{n}(\varphi\wedge\neg\varphi)$ or $\textsf{L}(\varphi^{1})$, for $1\leq k\leq n-1$ and $\textsf{L}\neq \textsf{t}^{n}_{k-1}$, then $\theta$ being true under $\nu$ would imply that $\nu(\varphi)=t_{k}^{n}$ and, from the definition of $\mathcal{F}_{C_{n}}$, $\nu(\varphi\wedge\neg\varphi)\in\{t_{0}^{n}, t_{1}^{n}, \dotsc  , t_{n-1}^{n}\}$ and $\nu(\varphi^{1})=t_{k-1}^{n}$, what in any of the two cases leads to a contradiction.
\end{enumerate}

\begin{lemma}\label{validates consequence rule}
For a signed formula $\textsf{L}(\varphi)$, with $\varphi$ not a propositional variable, and $\nu$ a valuation in $\mathcal{F}_{C_{n}}$ under which $\textsf{L}(\varphi)$ is true, all the formulas of at least one of the branches obtained from the application of a rule of $\mathbb{T}_{n}$ to $\textsf{L}(\varphi)$ are also true under $\nu$.
\end{lemma}

\begin{proof}
This is true from the definition of $\mathcal{A}_{C_{n}}$ due to the form of the rules in $\mathbb{T}_{n}$ and the very definition of being true under a valuation. For an example, look at the rule $E_{\wedge}\textsf{F}_{n}$, given by
\[\begin{array}{c|c}\multicolumn{2}{c}{\textsf{F}_{n}(\varphi\wedge\psi)}\\\hline\textsf{F}_{n}(\varphi)&\textsf{F}_{n}(\psi)\end{array}.\]
If $\nu$ satisfies $\textsf{F}_{n}(\varphi\wedge\psi)$, this means $\nu(\varphi\wedge\psi)=F_{n}$; looking at the table for conjunction in $\mathcal{A}_{C_{n}}$, presented below, we see that this implies either $\nu(\varphi)=F_{n}$ or $\nu(\psi)=F_{n}$, and so $\nu$ satisfies either $\textsf{F}_{n}(\varphi)$ or $\textsf{F}_{n}(\psi)$, that is, precisely all formulas of one of the branches obtained from applying $E_{\wedge}\textsf{F}_{n}$ to $\textsf{F}_{n}(\varphi\wedge\psi)$.
\begin{center}
\begin{tabular}{|l|c|c|r|}
\hline
$\wedge$ & $F_{n}$ & $I_{n}$ & $T_{n}$ \\ \hline
$F_{n}$ & $\{F_{n}\}$ & $\{F_{n}\}$ & $\{F_{n}\}$ \\ \hline
$I_{n}$ & $\{F_{n}\}$  & $D_{n}$ & $D_{n}$ \\ \hline
$T_{n}$ & $\{F_{n}\}$ & $D_{n}$ & $\{T_{n}\}$ \\ \hline
\end{tabular}
\end{center}
Of course, a similar reasoning works for all the other rules of $\mathbb{T}_{n}$.
\end{proof}

\begin{theorem}\label{Correct tableaux}
For $\Gamma\cup\{\varphi\}$ a finite set of formulas of $C_{n}$, if $\Gamma\vdash_{\mathbb{T}_{n}}\varphi$, then $\Gamma\vDash_{\mathcal{RM}_{C_{n}}}\varphi$.
\end{theorem}

\begin{proof}
Given $\Gamma=\{\gamma_{1}, \dotsc  , \gamma_{m}\}$, it is clear that $\bigwedge_{i=1}^{m}\gamma_{i}\vdash_{C_{n}}\gamma_{j}$, for every $1\leq j\leq m$, and $\Gamma\vdash_{C_{n}}\bigwedge_{i=1}^{m}\gamma_{i}$; furthermore, $C_{n}$ satisfies the deduction meta-theorem, and by soundness and completeness of this system with respect to $\mathcal{RM}_{C_{n}}$, it is enough to prove the theorem under the assumption that $\Gamma=\emptyset$.

One constructs a completed tableau for $\textsf{F}_{n}(\varphi)$ through a finite sequence of tableaux $\mathcal{T}_{0}, \dotsc  , \mathcal{T}_{k}=\mathcal{T}$, where $\mathcal{T}_{0}$ contains only $\textsf{F}_{n}(\varphi)$ and each $\mathcal{T}_{j+1}$, for $0\leq j\leq k-1$, is obtained from $\mathcal{T}_{j}$ by application of a rule $R_{j}$ of $\mathbb{T}_{n}$ to a, up to this moment, unused signed formula $\textsf{L}_{j}(\psi_{j})$ of a branch $\theta_{j}$ of $\mathcal{T}_{j}$.

Given a valuation $\nu$ of $\mathcal{F}_{C_{n}}$, we can then prove that if $\mathcal{T}_{j}$ is true under $\nu$, then so is $\mathcal{T}_{j+1}$. To see this, suppose $\nu$ indeed validates $\mathcal{T}_{j}$, and therefore validates some branch $\theta$ of $\mathcal{T}_{j}$:
\begin{enumerate}
\item if $\theta=\theta_{j}$, that is, $\mathcal{T}_{j+1}$ is obtained from $\mathcal{T}_{j}$ by extension of precisely the branch satisfied by $\nu$, by Lemma  \ref{validates consequence rule} we find that all new formulas of one of the branches created by applying $R_{j}$ to $\textsf{L}_{j}(\psi_{j})$ are validated by $\nu$; this branch, which extends $\theta_{j}$, is then a branch of $\mathcal{T}_{j+1}$ true under $\nu$;
\item if $\theta\neq\theta_{j}$, $\theta$ is still a branch of $\mathcal{T}_{j+1}$, and in particular it is true under $\nu$.
\end{enumerate}
In any case, the conclusion must be that $\nu$ validates $\mathcal{T}_{j+1}$.

So suppose $\varphi$ is not a tautology according to $\mathcal{RM}_{C_{n}}$, and therefore there exists a valuation $\nu\in\mathcal{F}_{C_{n}}$ with $\nu(\varphi)=F_{n}$. Since $\nu$ satisfies $\textsf{F}_{n}(\varphi)$, and therefore $\mathcal{T}_{0}$, by induction we prove that any completed tableau for $\textsf{F}_{n}(\varphi)$ is true under $\nu$, and so it cannot be closed. Since every completed tableau for $\textsf{F}_{n}(\varphi)$ is open, $\not\vdash_{\mathbb{T}_{n}}\varphi$, and by the contrapositive we have the theorem.
\end{proof}

To prove completeness, we must define what are Hintikka sets for $\mathbb{T}_{n}$.

\begin{definition}\label{Hintikkaset}
A non-empty set $\Gamma$ of labeled formulas $\textsf{L}(\varphi)$, with $\textsf{L}$ in $\textsf{B}_{n}$ and $\varphi$ a formula of $C_{n}$, is a Hintikka set for $\mathbb{T}_{n}$ if all of the following properties are simultaneously satisfied, where $0\leq i\leq n-1$ and $0\leq k\leq n-2$.
\begin{enumerate}
\item If $\textsf{L}(\varphi), \textsf{L}^{*}(\varphi)\in \Gamma$, then $\textsf{L}=\textsf{L}^{*}$.
\item If $\textsf{t}_{0}^{n}(\varphi)$, no $\textsf{t}_{j}^{n}(\varphi\wedge\neg\varphi)$, for $0\leq j\leq n-1$, is in $\Gamma$.
\item If $\textsf{t}_{k+1}^{n}(\varphi)\in\Gamma$, then $\textsf{T}_{n}(\varphi\wedge\neg\varphi)$ is not in $\Gamma$.
\item If $\textsf{t}_{k+1}^{n}(\varphi), \textsf{L}(\varphi^{1})\in\Gamma$, then $\textsf{L}=\textsf{t}_{k}^{n}$.
\item If $\textsf{T}_{n}(\neg\varphi)\in\Gamma$, at least one of $\textsf{t}_{0}^{n}(\varphi), \dotsc  , \textsf{t}_{n-1}^{n}(\varphi), \textsf{F}_{n}(\varphi)$ is also in $\Gamma$.
\item If $\textsf{t}_{i}^{n}(\neg\varphi)\in\Gamma$, at least one of $\textsf{t}_{0}^{n}(\varphi), \dotsc  , \textsf{t}_{n-1}^{n}(\varphi)$ is also in $\Gamma$.
\item If $\textsf{F}_{n}(\neg\varphi)\in\Gamma$, so is $\textsf{T}_{n}(\varphi)$.
\item If $\textsf{T}_{n}(\varphi\wedge\psi)\in\Gamma$, then:
\begin{enumerate}
\item either $\textsf{T}_{n}(\varphi), \textsf{T}_{n}(\psi)\in\Gamma$;
\item or $\textsf{T}_{n}(\varphi), \textsf{t}_{j}^{n}(\psi)\in\Gamma$, for some $0\leq j\leq n-1$;
\item or $\textsf{t}_{j}^{n}(\varphi), \textsf{T}_{n}(\psi)\in \Gamma$, for a $0\leq j\leq n-1$;
\item or $\textsf{t}_{j}^{n}(\varphi), \textsf{t}_{l}^{n}(\psi)\in\Gamma$, for $0\leq j, l\leq n-1$.
\end{enumerate}
\item If $\textsf{t}_{i}^{n}(\varphi\wedge\psi)\in\Gamma$, then either:
\begin{enumerate}
\item $\textsf{T}_{n}(\varphi), \textsf{t}_{j}^{n}(\psi)\in\Gamma$, for a $0\leq j\leq n-1$;
\item or $\textsf{t}_{j}^{n}(\varphi), \textsf{T}_{n}(\psi)\in\Gamma$, for some $0\leq j\leq n-1$;
\item or $\textsf{t}_{j}^{n}(\varphi), \textsf{t}_{l}^{n}(\psi)\in\Gamma$, for $0\leq j, l\leq n-1$.
\end{enumerate}
\item If $\textsf{F}_{n}(\varphi\wedge\psi)\in\Gamma$, either $\textsf{F}_{n}(\varphi)$ or $\textsf{F}_{n}(\psi)$ is in $\Gamma$.
\item If $\textsf{T}_{n}(\varphi\vee\psi)$ is in $\Gamma$, at least one of $\textsf{T}_{n}(\varphi)$, $\textsf{T}_{n}(\psi)$, $\textsf{t}_{j}^{n}(\varphi)$ or $\textsf{t}_{l}^{n}(\psi)$ is also in $\Gamma$, for $0\leq j, l\leq n-1$.
\item If $\textsf{t}_{i}^{n}(\varphi\vee\psi)$ is in $\Gamma$, either $\textsf{t}_{j}^{n}(\varphi)$ or $\textsf{t}_{l}^{n}(\psi)$, for $0\leq j,l\leq n-1$, is also in $\Gamma$.
\item If $\textsf{F}_{n}(\varphi\vee\psi)\in\Gamma$, then both $\textsf{F}_{n}(\varphi)$ and $\textsf{F}_{n}(\psi)$ are in $\Gamma$.
\item If $\textsf{T}_{n}(\varphi\rightarrow\psi)\in\Gamma$, at least one of $\textsf{F}_{n}(\varphi)$, $\textsf{T}_{n}(\psi)$ or $\textsf{t}_{j}^{n}(\psi)$, for $0\leq j\leq n-1$, is also in $\Gamma$.
\item If $\textsf{t}_{i}^{n}(\varphi\rightarrow\psi)$ is in $\Gamma$, either:
\begin{enumerate}
\item $\textsf{t}_{j}^{n}(\varphi)$ and $\textsf{T}_{n}(\psi)$, for a $0\leq j\leq n-1$, are both in $\Gamma$;
\item or $\textsf{t}_{j}^{n}(\psi)$, for a $0\leq j\leq n-1$, is in $\Gamma$.
\end{enumerate}
\item If $\textsf{F}_{n}(\varphi\rightarrow\psi)$ is in $\Gamma$, then:
\begin{enumerate}
\item either $\textsf{T}_{n}(\varphi)$ and $\textsf{F}_{n}(\psi)$ are both in $\Gamma$;
\item or $\textsf{t}_{j}^{n}(\varphi)$ and $\textsf{F}_{n}(\psi)$, for a $0\leq j\leq n-1$, are both in $\Gamma$.
\end{enumerate}
\end{enumerate}
\end{definition}

We must now prove that a Hintikka set for $\mathbb{T}_{n}$ is satisfiable in $\mathcal{F}_{C_{n}}$, that is, that there exists a $\nu\in\mathcal{F}_{C_{n}}$ such that $\textsf{L}(\varphi)\in\Gamma$ implies $\nu(\varphi)=L$. We start by defining, for a Hintikka set $\Gamma$ for $\mathbb{T}_{n}$, $\Gamma_{0}$ as the set of formulas of $C_{n}$ such that $\textsf{L}(\varphi)\in\Gamma$, for some label $\textsf{L}$ of $\textsf{B}_{n}$; then, we take $\nu_{0}:\Gamma_{0}\rightarrow B_{n}$ defined by $\nu_{0}(\varphi)=L$ iff $\textsf{L}(\varphi)\in \Gamma$. That $\nu_{0}$ is indeed a function one can prove by noticing clause $(1)$ of Definition \ref{Hintikkaset}: if one had both $\nu_{0}(\varphi)=L$ and $\nu_{0}(\varphi)=L^{*}$, for $L\neq L^{*}$, then it would follow that $\textsf{L}(\varphi), \textsf{L}^{*}(\varphi)\in \Gamma$, for $\textsf{L}\neq\textsf{L}^{*}$, contradicting clause $(1)$.

We extend $\nu_{0}$ to a function $\nu:F(\Sigma_{\textbf{C}}, \mathcal{V})\rightarrow B_{n}$ by induction on the complexity of a formula in much the same way we proceeded in Proposition \ref{Rows of tables are homomorphisms}. That is, for a propositional variable $p$, we define:
\begin{enumerate}
\item $\nu(p)=\nu_{0}(p)$, if $p\in\Gamma_{0}$, and arbitrarily otherwise;
\item $\nu(\neg p)=\nu_{0}(\neg p)$, if $\neg p\in\Gamma_{0}$, and as any value in $\tilde{\neg}\nu(p)$ otherwise;
\item $\nu(p\wedge\neg p)=\nu_{0}(p\wedge\neg p)$, in the case that $p\wedge\neg p\in\Gamma_{0}$; if $p\wedge\neg p\not\in \Gamma_{0}$ and $\nu(p)=t_{0}^{n}$, we make $\nu(p\wedge\neg p)=T_{n}$, but if $p\wedge\neg p\not\in\Gamma_{0}$ and $\nu(p)=t_{k}^{n}$, for $1\leq k\leq n-1$, we define $\nu(p\wedge\neg p)$ as an arbitrary value of $I_{n}$; if none of these is the case, $\nu(p\wedge\neg p)$ may be given any value in $\nu(p)\tilde{\wedge}\nu(\neg p)$;
\item $\nu(p^{1})=\nu_{0}(p^{1})$, whenever $p^{1}\in \Gamma_{0}$; in the case that $p^{1}\not\in\Gamma_{0}$, if $\nu(p)=t_{k}^{n}$, for a $1\leq k\leq n-1$, we make $\nu(p^{1})=t_{k-1}^{n}$, and otherwise $\nu(p^{1})$ can equal any value of $\tilde{\neg}\nu(p\wedge\neg p)$.
\end{enumerate}

If $\nu(\alpha)$ and $\nu(\beta)$, as well as $\nu(\neg\alpha)$, $\nu(\alpha\wedge\neg\alpha)$, $\nu(\alpha^{1})$, $\nu(\neg\beta)$, $\nu(\beta\wedge\neg\beta)$, and $\nu(\beta^{1})$, have been already defined, we may also define $\nu(\varphi)$ and $\nu(\psi)$, as well as other correlated values of $\nu$, for $\varphi=\neg\alpha$ and $\psi=\alpha\#\beta$ (where $\#\in\{\vee, \wedge, \rightarrow\}$). For $\varphi$ we have the following.
\begin{enumerate}
\item $\nu(\varphi)=\nu(\neg\alpha)$ has already been defined.
\item If $\neg\varphi\in\Gamma_{0}$, make $\nu(\neg\varphi)$ equal to $\nu_{0}(\neg\varphi)$, otherwise it is enough to make $\nu(\varphi)$ equal to some value in $\tilde{\neg}\nu(\varphi)$.
\item If $\varphi\wedge\neg\varphi\in\Gamma_{0}$, one makes $\nu(\varphi\wedge\neg\varphi)$ equal to $\nu_{0}(\varphi\wedge\neg\varphi)$; if $\varphi\wedge\neg\varphi$ is not in $\Gamma_{0}$ and $\nu(\varphi)\in I_{n}$, we either make $\nu(\varphi\wedge\neg\varphi)$ equal $T_{n}$, if $\nu(\varphi)=t^{n}_{0}$, or let it take any value in $I_{n}$, if $\nu(\varphi)=t_{k}^{n}$ for $1\leq k\leq n-1$; if $\varphi\wedge\neg\varphi\not\in\Gamma_{0}$ and $\nu(\varphi)\not\in I_{n}$, it is enough for $\nu(\varphi\wedge\neg\varphi)$ to lie in $\nu(\varphi)\tilde{\wedge}\nu(\neg\varphi)$.
\item We make $\nu(\varphi^{1})$ equal to: $\nu_{0}(\varphi^{1})$ if $\varphi^{1}\in\Gamma_{0}$; $t_{k-1}^{n}$ if $\varphi^{1}\not\in\Gamma_{0}$ and $\nu(\varphi)=t_{k}^{n}$, for $1\leq k\leq n-1$; and any value of $\tilde{\neg}\nu(\varphi\wedge\neg\varphi)$ if $\varphi^{1}\not\in\Gamma_{0}$ and $\nu(\varphi)\in\{T_{n}, t_{0}^{n}, F_{n}\}$.
\end{enumerate}
For $\psi$ the conditions are as bellow.
\begin{enumerate}
\item If $\psi\in\Gamma_{0}$, it is necessary to make $\nu(\psi)=\nu_{0}(\psi)$. In the case that $\psi\not\in\Gamma_{0}$, if $\#=\wedge$ and $\beta=\neg\alpha$, $\nu(\psi)$ has already been defined; otherwise, $\nu(\psi)$ is given any value among those found in $\nu(\alpha)\tilde{\#}\nu(\beta)$.
\item To define $\nu(\neg\psi)$, we again start by demanding that, in the case that $\neg\psi\in\Gamma_{0}$, it must equal $\nu_{0}(\neg\psi)$. If $\neg\psi\not\in\Gamma_{0}$, and $\#=\wedge$ and $\beta=\neg\alpha$, $\nu(\neg\psi)=\nu(\alpha^{1})$ has already been defined; otherwise $\nu(\neg\psi)$ can be given any value found in $\tilde{\neg}\nu(\psi)$.
\item In the case that $\psi\wedge\neg\psi\in\Gamma_{0}$, $\nu(\psi\wedge\neg\psi)=\nu_{0}(\psi\wedge\neg\psi)$. When $\psi\wedge\neg\psi\not\in\Gamma_{0}$: if $\nu(\psi)=t_{0}^{n}$, $\nu(\psi\wedge\neg\psi)$ is defined as $T_{n}$; if $\nu(\psi)=t_{k}^{n}$ (where $1\leq k\leq n-1$), $\nu(\psi\wedge\neg\psi)$ equals an arbitrary value of $I_{n}$; otherwise, $\nu(\psi\wedge\neg\psi)$ can take any value on $\nu(\psi)\tilde{\wedge}\nu(\neg\psi)$.
\item If $\psi^{1}\in\Gamma_{0}$, $\nu(\psi^{1})$ is $\nu_{0}(\psi^{1})$. If $\psi^{1}\not\in\Gamma_{0}$ and $\nu(\psi)=t_{k}^{n}$, for a $1\leq k\leq n-1$, we make $\nu(\psi^{1})$ equal to $t_{k-1}^{n}$. And if $\psi^{1}\not\in\Gamma_{0}$ and $\nu(\psi)\not\in\{t_{1}^{n}, \dotsc  , t_{n-1}^{n}\}$, we can give to $\nu(\psi^{1})$ any value from $\tilde{\neg}\nu(\psi\wedge\neg\psi)$.
\end{enumerate}

It is easy to see both that $\nu$ is a well defined function from $F(\Sigma_{C}, \mathcal{V})$ to $B_{n}$ and that, if $\nu$ is a homomorphism, then it certainly lies in $\mathcal{F}_{C_{n}}$. The difficulty here is convincing oneself that $\nu$ is indeed a homomorphism: the problem here is that, while in Proposition \ref{Rows of tables are homomorphisms} the relevant set of formulas was closed under subformulas, the set $\Gamma_{0}$ is not; it would seem possible for a formula $\beta\in\Gamma_{0}$ to have a subformula $\alpha\not\in\Gamma_{0}$ such that $\nu(\alpha)$, as defined above, is incompatible with $\nu(\beta)=\nu_{0}(\beta)$, in the sense that the function $\nu$ cannot be a homomorphism. We will argue that this does not happen, although the complete picture will be left for the reader to fill out.

Only a few clauses of Definition \ref{Hintikkaset} do not presuppose that the immediate subformulas of $\beta$ are also in $\Gamma_{0}$, specifically clauses $10$, $11$, $12$, $14$ and $15$ (item $(b)$); regarding the remaining clauses, the proof that $\varphi$ behaves like a homomorphism should runs smoothly. So, take one of these difficult clauses to analyze more thoroughly, let us say $10$: so, there are formulas $\varphi$ and $\psi$ such that $\textsf{F}_{n}(\varphi\wedge\psi)\in \Gamma$, meaning then that $\varphi\wedge\psi$ is in $\Gamma_{0}$, and $\textsf{F}_{n}(\varphi)\in \Gamma$, and thus $\varphi\in\Gamma_{0}$; the case in which $\textsf{F}_{n}(\varphi)$ is not in $\Gamma$ but $\textsf{F}_{n}(\psi)$ is can be treated in an analogous way.

By the definition of $\nu$ we have given, we have that $\nu(\varphi\wedge\psi)=\nu(\varphi)=F_{n}$, and $\nu(\psi)$ takes any value as long as $\nu$ is still a well-defined function that satisfies all conditions to be in $\mathcal{F}_{C_{n}}$, except perhaps being a homomorphism. We ask ourselves: is there any value we can give to $\nu(\psi)$ that would lead to $\nu(\varphi\wedge\psi)\not\in\nu(\varphi)\tilde{\wedge}\nu(\psi)$ and therefore force $\nu$ not to be a homomorphism? The answer is no: since $\nu(\varphi)=F_{n}$, looking at the table for $\tilde{\wedge}$ just after Definition \ref{Define operations in Cn}, we see that, for any value possibly taken by $\nu(\psi)$, $\nu(\varphi)\tilde{\wedge}\nu(\psi)=F_{n}\tilde{\wedge}\nu(\psi)=\{F_{n}\}$, which contains $\nu(\varphi\wedge\psi)$. Briefly put, the definition of Hintikka sets, despite not being strong enough to guarantee that those sets are closed under subformulas, carries enough information to make sure the $\nu$, as defined above, is a homomorphism in $\mathcal{F}_{C_{n}}$.

For another example, look now at clause $15(b)$: there are formulas $\varphi$ and $\psi$ and indices $0\leq i, j\leq n-1$ for which $\textsf{t}_{i}^{n}(\varphi\rightarrow\psi)$ and $\textsf{t}_{j}^{n}(\psi)$ are both in $\Gamma$, and so $\varphi\rightarrow\psi, \psi\in \Gamma_{0}$: is it possible to find a value of $\nu(\varphi)$ that makes of $\nu$ not a homomorphism? No: by looking at the table for $\tilde{\rightarrow}$ in $\mathcal{A}_{C_{n}}$ below, we see that $\nu(\varphi\rightarrow\psi)=t_{i}^{n}$ is in $D_{n}=\nu(\varphi)\tilde{\rightarrow}t_{j}^{n}=\nu(\varphi)\tilde{\rightarrow}\nu(\psi)$, for any $\nu(\varphi)\in B_{n}$.

\begin{figure}[H]
\centering
\begin{tabular}{|l|c|c|r|}
\hline
$\rightarrow$ & $F_{n}$ & $I_{n}$ & $T_{n}$ \\ \hline
$F_{n}$ & $\{T_{n}\}$ & $D_{n}$ & $\{T_{n}\}$ \\ \hline
$I_{n}$ & $\{F_{n}\}$ & $D_{n}$ & $D_{n}$ \\ \hline
$T_{n}$ & $\{F_{n}\}$ & $D_{n}$ & $\{T_{n}\}$ \\ \hline
\end{tabular}
\caption*{Table for Implication}
\end{figure}

Summarizing what we have done just above, one gets the following result.

\begin{proposition}\label{Hintikka's pre-lemma}
For any Hintikka set $\Gamma$ for $\mathbb{T}_{n}$, we define $\Gamma_{0}=\{\varphi\in\textbf{F}(\Sigma_{\textbf{C}}, \mathcal{V}): \textsf{L}(\varphi)\in\Gamma\}$ and the function $\nu_{0}:\Gamma_{0}\rightarrow B_{n}$ by $\nu_{0}(\varphi)=L$ iff $\textsf{L}(\varphi)\in\Gamma$. In this case, there is a homomorphism $\nu\in\mathcal{F}_{C_{n}}$ extending $\nu_{0}$.
\end{proposition}

We are now in position to trivially prove the equivalent to Hintikka's lemma in $\mathbb{T}_{n}$.

\begin{theorem}\label{Hintikka's lemma}
For any Hintikka set $\Gamma$ for $\mathbb{T}_{n}$, there exists a $\nu\in\mathcal{F}_{C_{n}}$ satisfying that $\textsf{L}(\varphi)\in\Gamma$ implies $\nu(\varphi)=L$.
\end{theorem}

\begin{proof}
Take $\Gamma_{0}$ as the set of formulas $\varphi$ for which there exists a label $\textsf{L}$ with $\textsf{L}(\varphi)\in\Gamma$, and define the function $\nu:\Gamma_{0}\rightarrow B_{n}$ by $\nu_{0}(\varphi)=L$ iff $\textsf{L}(\varphi)\in\Gamma$, as it would be expected. By Proposition \ref{Hintikka's pre-lemma}, there exists $\nu\in\mathcal{F}_{C_{n}}$ extending $\nu_{0}$ and, of course, if $\textsf{L}(\varphi)\in\Gamma$, $\nu(\varphi)=\nu_{0}(\varphi)=L$.
\end{proof}

Hintikka sets for $\mathbb{T}_{n}$, as it would be expected, correspond to open branches $\theta$ of a complete tableau $\mathcal{T}$ of $\mathbb{T}_{n}$: in fact, define as $\Gamma$ the set of labeled formulas appearing in $\theta$. Since $\theta$ is open, from Definition \ref{Definitions about tableaux} it follows that if $\textsf{L}(\varphi)$ and $\textsf{L}^{*}(\varphi)$ are both in $\Gamma$, meaning that both are in $\theta$, then $\textsf{L}=\textsf{L}^{*}$, what amounts to clause $1$ in Definition \ref{Hintikkaset}.

If $\textsf{t}_{0}^{n}(\varphi)$ and some $\textsf{t}_{j}^{n}(\varphi\wedge\neg\varphi)$, for $0\leq j\leq n-1$, were both in $\Gamma$, and so in $\theta$, this would make of $\theta$ a closed branch, contradicting our supposition that $\theta$ is actually open; of course, we get from this that clause $2$ of Definition \ref{Hintikkaset} is also satisfied. And if $\textsf{t}^{n}_{k+1}(\varphi)\in\Gamma$, for $0\leq k\leq n-2$, then neither $\textsf{T}_{n}(\varphi\wedge\neg\varphi)$ nor $\textsf{L}(\varphi)$, for $\textsf{L}\neq \textsf{t}_{k}^{n}$, can be in $\Gamma$, corresponding to clauses $3$ and $4$. Finally, for $\textsf{L}(\neg\varphi)$ or $\textsf{L}(\varphi\#\psi)$ (with $\#\in\{\vee, \wedge, \rightarrow\}$) in $\Gamma$, one notices that $\mathcal{T}$ being complete forces $\theta$ to contain all labeled formulas of one branch of the rule in $\mathbb{T}_{n}$ headed by $\textsf{L}(\neg\varphi)$ or respectively $\textsf{L}(\varphi\#\psi)$, what leads to clauses $5$ through $16$ of Definition \ref{Hintikkaset} being satisfied. 

These observations, combined with Theorem \ref{Hintikka's lemma}, allow one to prove the following result.

\begin{corollary}\label{corollary to Hintika's lemma}
For $\textsf{L}(\varphi)$ a signed formula of $C_{n}$, let $\theta$ be an open branch of a completed tableau $\mathcal{T}$ in $\mathbb{T}_{n}$ starting with $\textsf{L}(\varphi)$, and $\Gamma$ be the set of labeled formulas appearing in $\theta$. There exists a homomorphism $\nu\in\mathcal{F}_{C_{n}}$ such that $\nu(\psi)=L$ iff $\textsf{L}(\psi)\in\Gamma$.
\end{corollary}

\begin{theorem}\label{completeness tableaux}
For $\Gamma\cup\{\varphi\}$ a finite set of formulas of $C_{n}$, if $\Gamma\vDash_{\mathcal{RM}_{C_{n}}}\varphi$, then $\Gamma\vdash_{\mathbb{T}_{n}}\varphi$.
\end{theorem}

\begin{proof}
As it was done in Theorem \ref{Correct tableaux}, we can assume that $\Gamma=\emptyset$: this is because, if $\Gamma=\{\gamma_{1}, \dotsc  , \gamma_{m}\}$, since $\bigwedge_{i}^{m}\gamma_{i}\vdash_{C_{n}}\gamma_{j}$, for every $1\leq j\leq m$, and $\Gamma\vdash_{C_{n}}\bigwedge_{i=1}^{m}\gamma_{i}$, we get that $\Gamma\vdash_{C_{n}}\varphi$ if and only if $\vdash_{C_{n}}\bigwedge_{i=1}^{m}\gamma_{i}\rightarrow\varphi$. Given that $\Gamma\vdash_{C_{n}}\varphi$ iff $\Gamma\vDash_{\mathcal{RM}_{C_{n}}}\varphi$, and $\Gamma\vdash_{\mathbb{T}_{n}}\varphi$ iff $\vdash_{\mathbb{T}_{n}}\bigwedge_{i=1}^{n}\gamma_{i}\rightarrow\varphi$ by definition, it is indeed enough to take $\Gamma=\emptyset$.

So we suppose that $\vDash_{\mathcal{RM}_{C_{n}}}\varphi$: if $\mathcal{T}$ is a completed tableau in $\mathbb{T}_{n}$ starting with $\textsf{F}_{n}(\varphi)$ with an open branch $\theta$, by Corollary \ref{corollary to Hintika's lemma} there exists a homomorphism $\nu\in\mathcal{F}_{C_{n}}$ such that $\nu(\psi)=L$ iff $\textsf{L}(\psi)$ appears in $\theta$. Of course, since $\textsf{F}_{n}(\varphi)$ is in $\theta$, $\nu(\varphi)=F_{n}$ and therefore $\varphi$ is not a tautology, meaning that $\not\vDash_{\mathcal{RM}_{C_{n}}}\varphi$ and leading to a contradiction.

This means that, if $\vDash_{\mathcal{RM}_{C_{n}}}\varphi$, every completed tableau in $\mathbb{T}_{n}$ starting with $\textsf{F}_{n}(\varphi)$ must be closed, if there exist any; since it is always possible to construct a complete tableau in $\mathbb{T}_{n}$ starting from any given labeled formula in a finite number of steps, it follows that $\vdash_{\mathbb{T}_{n}}\varphi$.
\end{proof}

Of course, Theorems \ref{Correct tableaux} and \ref{completeness tableaux} prove that the calculi $\mathbb{T}_{n}$ are correct and complete in relation to the respective logics $C_{n}$. They are also decision methods, given that our tableaux can always be constructed in a finite number of steps, inspired by the RNmatrices $\mathcal{RM}_{C_{n}}$ as our row-branching, row-eliminating truth-tables also were: it is easy to see that, if $\varphi$ is valid in $C_{n}$, every complete tableau starting with $\textsf{F}_{n}(\varphi)$ in $\mathbb{T}_{n}$ will be closed; and if $\varphi$ is not valid in $C_{n}$, every complete tableau starting with $\textsf{F}_{n}(\varphi)$ in $\mathbb{T}_{n}$ will be open, and each of its open branches produces a $\nu\in\mathcal{F}_{C_{n}}$ with $\nu(\varphi)=F_{n}$.

\subsubsection{Examples and derived rules}

In this section we present some examples of actually using the calculi $\mathbb{T}_{n}$ to test the validity of arguments.

For the argument $\neg\neg\alpha\vdash_{C_{1}}\alpha$, which is indeed valid as shown below, keep in mind that we are listing, at the right side of the tableau, the row (or rows), as well as their label and lead connective, that lead to the inclusion of the present row; this way, one can check how the tableau rules of $\mathbb{T}_{1}$ are being applied.

\begin{center}
\begin{prooftree}
{
}
[\textsf{F}_{1}(\neg\neg\alpha\rightarrow\alpha)
	[\textsf{T}_{1}(\neg\neg\alpha), just=$\textsf{F}_{1}\rightarrow$:!u
		[\textsf{F}_{1}(\alpha), just=$\textsf{F}_{1}\rightarrow$:!uu
			[\textsf{F}_{1}(\neg\alpha), just=$\textsf{T}_{1}\neg$:!uu
				[\textsf{T}_{1}(\alpha), just=$\textsf{F}_{1}\neg$:!u, close={:!uu, !c}]]
			[\textsf{t}_{0}^{1}(\neg\alpha)
				[\textsf{t}_{0}^{1}(\alpha), close={:!uu, !c}]]]]
	[\textsf{t}_{0}^{1}(\neg\neg\alpha)
		[\textsf{F}_{1}(\alpha)
			[\textsf{t}_{0}^{1}(\neg\alpha), just=$\textsf{t}_{0}^{1}\neg$:!uu
				[\textsf{t}_{0}^{1}(\alpha), just=$\textsf{t}_{0}^{1}\neg$:!u, close={:!uu, !c}]]]]]
\end{prooftree}
\end{center}

Testing the same argument but in $C_{2}$, that is $\neg\neg\alpha\vdash_{C_{2}}\alpha$, in the tableau of Figure \ref{Test in C2} shows that it remains valid, but the complexity of the corresponding tableau grows substantially. 

\begin{landscape}

This suggests, correctly, that the main problem found in dealing with the calculi $\mathbb{T}_{n}$ is how fast their tableaux grow horizontally, meaning how ramified they become.

To simplify the task of writing tableaux in $\mathbb{T}_{n}$, we can construct derived tableau rules: much like, when proving something axiomatically, we make use of intermediary results, we may obtain new tableaux rules derived from the ones already present in $\mathbb{T}_{n}$. We begin by defining some sets of labels: for $0\leq i\leq n-1$, 
\[\textsf{D}_{n}^{\leq i}=\{\textsf{T}_{n}, \textsf{t}^{n}_{0}, \dotsc  , \textsf{t}^{n}_{i}\},\quad \textsf{D}_{n}^{\geq i}=\{\textsf{t}^{n}_{i}, \dotsc  , \textsf{t}^{n}_{n-1}\}\quad\text{and}\quad\textsf{D}_{n}^{\leq -1}=\{\textsf{T}_{n}\}\]
(what somewhat reflects the heuristic that $T_{n}$ behaves like a $t^{n}_{-1}$ should).

\begin{figure}[H]
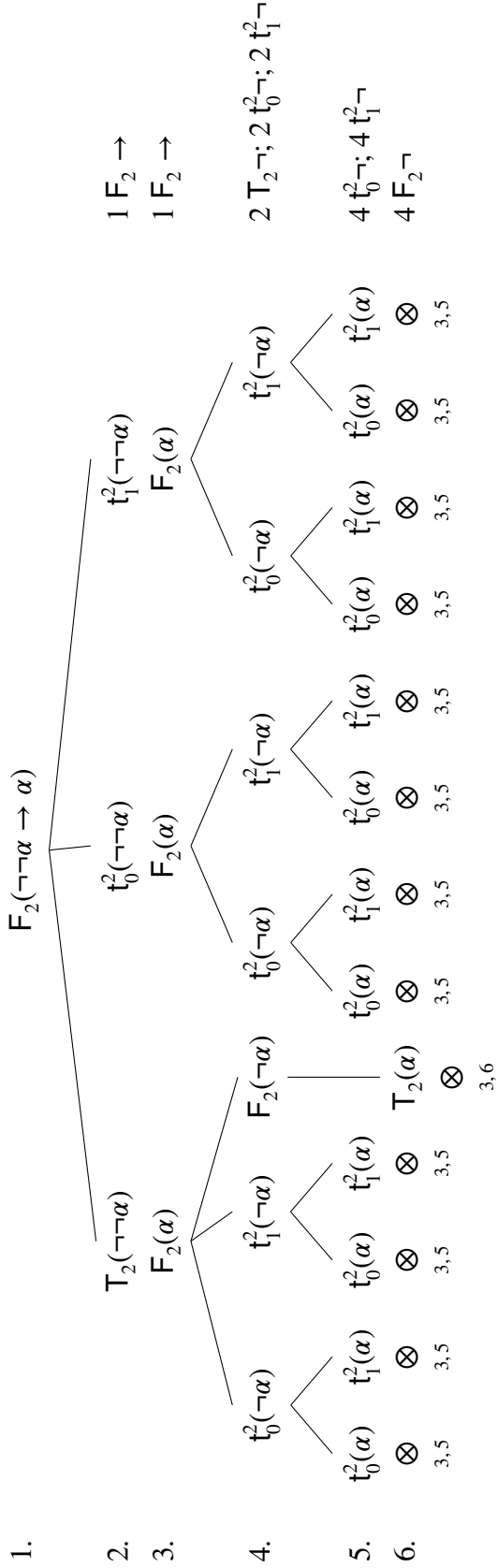

\centering
\begin{prooftree}
{
}
[\textsf{F}_{2}(\neg\neg\alpha\rightarrow\alpha)
	[\textsf{T}_{2}(\neg\neg\alpha), just=$\textsf{F}_{2}\rightarrow$:!u
		[\textsf{F}_{2}(\alpha), just=$\textsf{F}_{2}\rightarrow$:!uu
			[\textsf{t}_{0}^{2}(\neg\alpha), just=$\textsf{T}_{2}\neg$:!uu
				[\textsf{t}_{0}^{2}(\alpha), just=$\textsf{t}_{0}^{2}\neg$:!u, close={:!uu, !c}]
				[\textsf{t}_{1}^{2}(\alpha), close={:!uu, !c}]]
			[\textsf{t}_{1}^{2}(\neg\alpha)
				[\textsf{t}_{0}^{2}(\alpha), just=$\textsf{t}_{1}^{2}\neg$:!u, close={:!uu, !c}]
				[\textsf{t}_{1}^{2}(\alpha), close={:!uu, !c}]]
			[\textsf{F}_{2}(\neg\alpha)
				[\textsf{T}_{2}(\alpha), just=$\textsf{F}_{2}\neg$:!u, close={:!uu, !c}]]]]
	[\textsf{t}_{0}^{2}(\neg\neg\alpha)
		[\textsf{F}_{2}(\alpha)
			[\textsf{t}_{0}^{2}(\neg\alpha), just=$\textsf{t}_{0}^{2}\neg$:!uu
				[\textsf{t}_{0}^{2}(\alpha), close={:!uu, !c}]
				[\textsf{t}_{1}^{2}(\alpha), close={:!uu, !c}]]
			[\textsf{t}_{1}^{2}(\neg\alpha)
				[\textsf{t}_{0}^{2}(\alpha), close={:!uu, !c}]
				[\textsf{t}_{1}^{2}(\alpha), close={:!uu, !c}]]]]
	[\textsf{t}_{1}^{2}(\neg\neg\alpha)
		[\textsf{F}_{2}(\alpha)
			[\textsf{t}_{0}^{2}(\neg\alpha), just=$\textsf{t}_{1}^{2}\neg$:!uu
				[\textsf{t}_{0}^{2}(\alpha), close={:!uu, !c}]
				[\textsf{t}_{1}^{2}(\alpha), close={:!uu, !c}]]
			[\textsf{t}_{1}^{2}(\neg\alpha)
				[\textsf{t}_{0}^{2}(\alpha), close={:!uu, !c}]
				[\textsf{t}_{1}^{2}(\alpha), close={:!uu, !c}]]]]]
\end{prooftree}
\caption{Testing whether $\neg\neg\alpha\vdash_{C_{2}}\alpha$}
\label{Test in C2}
\end{figure}

Taking indices $0\leq i\leq n-1$, $0\leq j\leq n-2$, $k\geq n$ and $l\geq n-1$, we have the following derived rules for formulas of the form $\varphi^{m}\wedge\neg \varphi^{m}$.

\end{landscape}

\begin{table}[H]
\centering
$\begin{array}{ccccc}
\begin{array}{c} \textsf{T}_{n}(\varphi^{i}\wedge\neg\varphi^{i})\\\hline \textsf{t}^{n}_{i}(\varphi)\end{array} & \quad\quad\quad & \begin{array}{c} \textsf{t}^{n}_{i}(\varphi^{j}\wedge\neg\varphi^{j})\\\hline \textsf{D}_{n}^{\geq j+1}(\varphi)\end{array} & \quad\quad\quad & \begin{array}{c|c}\multicolumn{2}{c}{\textsf{F}_{n}(\varphi^{i}\wedge\neg\varphi^{i})}\\\hline \textsf{D}_{n}^{\leq i-1}(\varphi) & \textsf{F}_{n}(\varphi)\end{array}\\
& & & &\\
\begin{array}{c} \textsf{T}_{n}(\varphi^{k}\wedge\neg\varphi^{k})\\\hline \otimes\end{array}& \quad\quad\quad & \begin{array}{c} \textsf{t}^{n}_{i}(\varphi^{l}\wedge\neg\varphi^{l})\\\hline \otimes\end{array}& &
\end{array}$
\end{table}

Here, the symbol $\otimes$ at the bottom of a rule means that any branch where the labeled formula at the top of said rule appears is closed.

Proving these derived rules are valid is long and tedious, but rather trivial: to give one example, take the rule headed by $\textsf{T}_{n}(\varphi^{i}\wedge\neg\varphi^{i})$, and $n=1$ and $i=0$; more complex cases should be treated with inductive arguments. Then we have the following tableau, which indeed only leaves us $\textsf{t}^{1}_{0}(\varphi)$ as the rule stated.

\begin{center}
\begin{prooftree}
{
}
[\textsf{T}_{1}(\varphi\wedge\neg\varphi)
	[\textsf{T}_{1}(\varphi), just=$\textsf{T}_{1}\wedge$:!u
		[\textsf{T}_{1}(\neg\varphi),  just=$\textsf{T}_{1}\wedge$:!uu
			[\textsf{t}^{1}_{0}(\varphi), just=$\textsf{T}_{1}\neg$:!u, close={:!uu, !c}]
			[\textsf{F}_{1}(\varphi), close={:!uu, !c}]]]
	[\textsf{T}_{1}(\varphi)
		[\textsf{t}^{1}_{0}(\neg\varphi)
			[\textsf{t}^{1}_{0}(\varphi), just=$\textsf{t}^{1}_{0}\neg$:!u, close={:!uu, !c}]]]
	[\textsf{t}^{1}_{0}(\varphi)
		[\textsf{T}_{1}(\neg\varphi)
			[\textsf{t}^{1}_{0}(\varphi)]
			[\textsf{F}_{1}(\varphi), close={:!uu, !c}]]]
	[\textsf{t}^{1}_{0}(\varphi)
		[\textsf{t}^{1}_{0}(\neg\varphi)
			[\textsf{t}^{1}_{0}(\varphi)]]]]
\end{prooftree}
\end{center}

We can also easily obtain derived rules for formulas of the form $\varphi^{m}$, as we show below: take now indices $1\leq r\leq n$, $1\leq s\leq n-1$, $0\leq t\leq n-s-1$ and $n-s\leq u\leq n-1$.

\begin{table}[H]
\centering
$\begin{array}{ccccc}
\begin{array}{c|c}\multicolumn{2}{c}{\textsf{T}_{n}(\varphi^{r})}\\\hline \textsf{D}_{n}^{\leq r-2}(\varphi) & \textsf{F}_{n}(\varphi)\end{array} & \quad\quad\quad & \begin{array}{c}\textsf{t}^{n}_{t}(\varphi^{s})\\\hline \textsf{t}^{n}_{s+t}(\varphi)\end{array} & \quad\quad\quad & \begin{array}{c}\textsf{F}_{n}(\varphi^{r})\\\hline \textsf{t}^{n}_{r-1}(\varphi)\end{array}\\
 & & & &\\
& & \begin{array}{c}\textsf{t}^{n}_{u}(\varphi^{s})\\\hline \otimes\end{array} & & 
\end{array}$
\end{table}

Also important are derived rules for labeled formulas of the form $\neg\varphi^{m}$, where the indices are again $0\leq i\leq n-1$, $1\leq r\leq n$ and $1\leq s\leq n-1$.

\begin{table}[H]
\centering
$\begin{array}{ccccc}
\begin{array}{c}\textsf{T}_{n}(\neg\varphi^{r})\\\hline \textsf{D}_{n}^{\geq r-1}(\varphi)\end{array} & \quad\quad\quad & \begin{array}{c}\textsf{t}^{n}_{i}(\neg\varphi^{s})\\\hline \textsf{D}_{n}^{\geq s}(\varphi)\end{array} & \quad\quad\quad & \begin{array}{c|c}\multicolumn{2}{c}{\textsf{F}_{n}(\neg\varphi^{r})}\\\hline \textsf{D}_{n}^{\leq r-2}(\varphi) & \textsf{F}_{n}(\varphi)\end{array}
\end{array}$
\end{table}

All of these derived rules allow us to give more convoluted examples of tableaux, where the $DR$ to the right of a row will signify that a derived rule was used to obtain the said row.

\begin{landscape}

We start by showing that $\alpha, \neg\alpha, \alpha^{1}\vdash_{C_{1}}\neg\neg\alpha$, or what is equivalent, that $\vdash_{C_{1}}(\alpha\wedge\neg\alpha)\wedge\neg(\alpha\wedge\neg\alpha)\rightarrow\neg\neg\alpha$.

\begin{center}
\begin{prooftree}
{
}
[\textsf{F}_{1}(((\alpha\wedge\neg \alpha)\wedge\neg(\alpha\wedge\neg \alpha))\rightarrow\neg\neg \alpha)
	[\textsf{T}_{1}((\alpha\wedge\neg \alpha)\wedge\neg(\alpha\wedge\neg \alpha)), just=$\textsf{F}_{1}\rightarrow$:!u
		[\textsf{F}_{1}(\neg\neg \alpha), just=$\textsf{F}_{1}\rightarrow$:!uu
			[\textsf{t}^{1}_{0}(\alpha\wedge\neg \alpha), just=$DR$:!uu, close={:!c}]]]
	[\textsf{t}^{1}_{0}((\alpha\wedge\neg \alpha)\wedge\neg(\alpha\wedge\neg\alpha))
		[\textsf{F}_{1}(\neg\neg \alpha), close={:!u}]]]
\end{prooftree}
\end{center}

\begin{figure}[H]
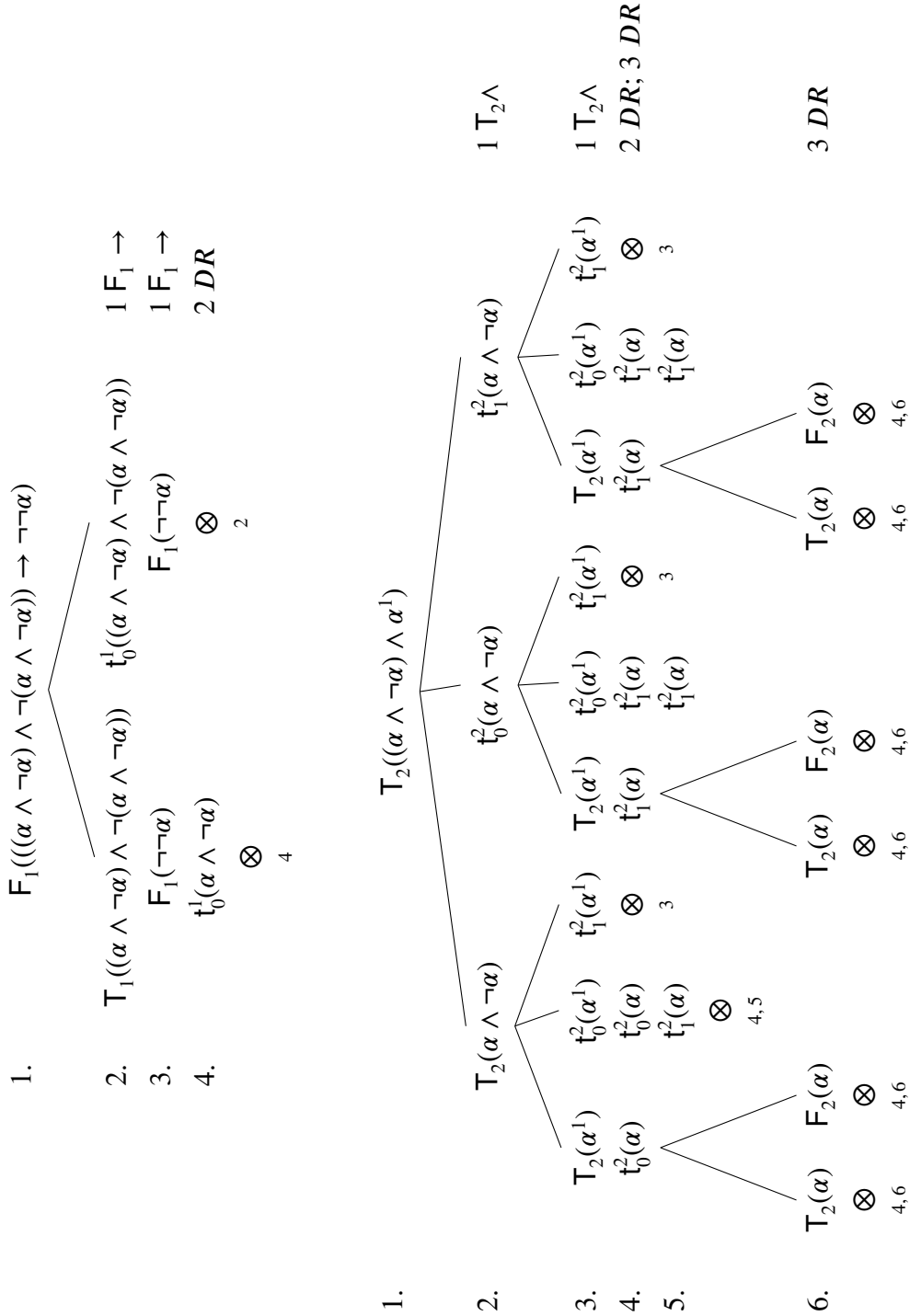

\centering
\begin{center}
\begin{prooftree}
{
}
[\textsf{T}_{2}((\alpha\wedge\neg\alpha)\wedge\alpha^{1})
	[\textsf{T}_{2}(\alpha\wedge\neg\alpha), just=$\textsf{T}_{2}\wedge$:!u
		[\textsf{T}_{2}(\alpha^{1}), just=$\textsf{T}_{2}\wedge$:!uu
			[\textsf{t}^{2}_{0}(\alpha), just=$DR$:!uu
				[\textsf{T}_{2}(\alpha), close={:!u, !c}, just=$DR$:!uu]
				[\textsf{F}_{2}(\alpha), close={:!u, !c}]]]
		[\textsf{t}^{2}_{0}(\alpha^{1})
			[\textsf{t}^{2}_{0}(\alpha), just=$DR$:!u
				[\textsf{t}^{2}_{1}(\alpha), close={:!u, !c}]]]
		[\textsf{t}^{2}_{1}(\alpha^{1}), close={:!c}]]
	[\textsf{t}^{2}_{0}(\alpha\wedge\neg\alpha)
		[\textsf{T}_{2}(\alpha^{1})
			[\textsf{t}^{2}_{1}(\alpha)
				[\textsf{T}_{2}(\alpha), close={:!u, !c}]
				[\textsf{F}_{2}(\alpha), close={:!u, !c}]]]
		[\textsf{t}^{2}_{0}(\alpha^{1})
			[\textsf{t}^{2}_{1}(\alpha)
				[\textsf{t}^{2}_{1}(\alpha)]]]
		[\textsf{t}^{2}_{1}(\alpha^{1}), close={:!c}]]
	[\textsf{t}^{2}_{1}(\alpha\wedge\neg\alpha)
		[\textsf{T}_{2}(\alpha^{1})
			[\textsf{t}^{2}_{1}(\alpha)
				[\textsf{T}_{2}(\alpha), close={:!u, !c}]
				[\textsf{F}_{2}(\alpha), close={:!u, !c}]]]
		[\textsf{t}^{2}_{0}(\alpha^{1})
			[\textsf{t}^{2}_{1}(\alpha)
				[\textsf{t}^{2}_{1}(\alpha)]]]
		[\textsf{t}^{2}_{1}(\alpha^{1}), close={:!c}]]]
\end{prooftree}
\end{center}
\caption{Testing whether $\alpha, \neg\alpha, \alpha^{1}\vdash_{C_{2}}\beta$}
\label{Second test in C2}
\end{figure}

\end{landscape}

Now, still in $C_{1}$, we show that $\alpha^{1}, \beta^{1}\vdash_{C_{1}}(\alpha\vee\beta)^{1}$ or, what is equivalent, $\vdash_{C_{1}}(\alpha^{1}\wedge\beta^{1})\rightarrow(\alpha\vee\beta)^{1}$.

\begin{center}
\begin{prooftree}
{
}
[\textsf{F}_{1}((\alpha^{1}\wedge\beta^{1})\rightarrow(\alpha\vee\beta)^{1})
	[\textsf{T}_{1}(\alpha^{1}\wedge\beta^{1}), just=$\textsf{F}_{1}\rightarrow$:!u
		[\textsf{F}_{1}((\alpha\vee\beta)^{1}), just=$\textsf{F}_{1}\rightarrow$:!uu
			[\textsf{t}^{1}_{0}(\alpha\vee\beta), just=$DR$:!u
				[\textsf{T}_{1}(\alpha), just=$DR$:!uuu
					[\textsf{T}_{1}(\beta), just=$DR$:!uuuu
						[\textsf{t}^{1}_{0}(\alpha), just=$\textsf{t}^{1}_{0}\vee$:!uuu, close={:!uu, !c}]
						[\textsf{t}^{1}_{0}(\beta), close={:!u, !c}]]]
				[\textsf{T}_{1}(\alpha)
					[\textsf{F}_{1}(\beta)
						[\textsf{t}^{1}_{0}(\alpha), close={:!uu, !c}]
						[\textsf{t}^{1}_{0}(\beta), close={:!u, !c}]]]
				[\textsf{F}_{1}(\alpha)
					[\textsf{T}_{1}(\beta)
						[\textsf{t}^{1}_{0}(\alpha), close={:!uu, !c}]
						[\textsf{t}^{1}_{0}(\beta), close={:!u, !c}]]]
				[\textsf{F}_{1}(\alpha)
					[\textsf{F}_{1}(\beta)
						[\textsf{t}^{1}_{0}(\alpha), just=$\textsf{t}^{1}_{0}\vee$:!uuu, close={:!uu, !c}]
						[\textsf{t}^{1}_{0}(\beta), close={:!u, !c}]]]]]]
	[\textsf{t}^{1}_{0}(\alpha^{1}\wedge\beta^{1})
		[\textsf{F}_{1}((\alpha\vee\beta)^{1}), close={:!u}]]]
\end{prooftree}
\end{center}

Finally, in Figure \ref{Second test in C2}, we show a negative result, now in $C_{2}$: it is true that $\alpha, \neg\alpha, \alpha^{1}$\\$\vdash_{C_{2}}\beta$ or, equivalently, $\vdash_{C_{2}}((\alpha\wedge\neg\alpha)\wedge\alpha^{1})\rightarrow\beta$? The answer, obviously, is no, yet for simplicity we prove something slightly different, but equivalent: since $\beta$ can assume any value, for this deduction to be true we must have that $(\alpha\wedge \neg\alpha)\wedge\alpha^{1}$ always equals $F_{2}$, what we show not to be the case; that is, there are homomorphisms on $\mathcal{F}_{C_{2}}$, given by the open branches of the tableau in Figure \ref{Second test in C2}, where $(\alpha\wedge\neg\alpha)\wedge\alpha^{1}$ indeed takes the value $T_{2}$.

\subsubsection{$\mbCcl$ and $\CILA$}

Before showing finite RNmatrices capable of characterizing the logics $C_{n}$, we had done the same for the logics $\mbCcl$ and $\CILA$ in Sections, respectively, \ref{RNmatrix for mbCcl} and \ref{RNmatrix for CILA} and, in exactly the same way we have motivated the introduction of the calculi $\mathbb{T}_{n}$, we may motivate both $\mathbb{T}_{\mbCcl}$ and $\mathbb{T}_{\CILA}$. Starting with $\mbCcl$, of course now our formulas are over the signature $\Sigma_{\textbf{LFI}}$ and our labels are $\textsf{T}$, $\textsf{t}$ or $\textsf{F}$, but the procedure is the same: we look at the tables of the RNmatrix for $\mbCcl$, obtaining the following rules.

\begin{table}[H]
\centering
$\begin{array}{ccccc}
\begin{array}{c|c|c|c}\multicolumn{4}{c}{\textsf{T}(\varphi\vee\psi)}\\\hline \textsf{T}(\varphi) & \textsf{t}(\varphi) & \textsf{T}(\psi) & \textsf{t}(\psi)\end{array} & \quad\quad\quad & \begin{array}{c|c|c|c}\multicolumn{4}{c}{\textsf{t}(\varphi\vee\psi)}\\\hline \textsf{T}(\varphi) & \textsf{t}(\varphi) & \textsf{T}(\psi) & \textsf{t}(\psi)\end{array} & \quad\quad\quad & \begin{array}{c}\textsf{F}(\varphi\vee\psi)\\\hline \textsf{F}(\varphi)\\ \textsf{F}(\psi)\end{array}\\
&&&&\\
\begin{array}{c|c|c|c}\multicolumn{4}{c}{\textsf{T}(\varphi\wedge\psi)}\\\hline \textsf{T}(\varphi) & \textsf{T}(\varphi) & \textsf{t}(\varphi) & \textsf{t}(\varphi)\\ \textsf{T}(\psi) & \textsf{t}(\psi) & \textsf{T}(\psi) & \textsf{t}(\psi) \end{array} & \quad\quad\quad & \begin{array}{c|c|c|c}\multicolumn{4}{c}{\textsf{t}(\varphi\wedge\psi)}\\\hline \textsf{T}(\varphi) & \textsf{T}(\varphi) & \textsf{t}(\varphi) & \textsf{t}(\varphi)\\ \textsf{T}(\psi) & \textsf{t}(\psi) & \textsf{T}(\psi) & \textsf{t}(\psi) \end{array} & \quad\quad\quad & \begin{array}{c|c}\multicolumn{2}{c}{\textsf{F}(\varphi\wedge\psi)}\\\hline \textsf{F}(\varphi)& \textsf{F}(\psi)\end{array}
\end{array}$
\end{table}

\begin{table}[H]
\centering
$\begin{array}{ccccc}
\begin{array}{c|c|c}\multicolumn{3}{c}{\textsf{T}(\varphi\rightarrow\psi)}\\\hline \textsf{F}(\varphi) & \textsf{T}(\psi) & \textsf{t}(\psi)\end{array} & \quad\quad\quad & \begin{array}{c|c|c}\multicolumn{3}{c}{\textsf{t}(\varphi\rightarrow\psi)}\\\hline \textsf{F}(\varphi) & \textsf{T}(\psi) & \textsf{t}(\psi)\end{array} & \quad\quad\quad & \begin{array}{c|c}\multicolumn{2}{c}{\textsf{F}(\varphi\rightarrow\psi)}\\\hline \textsf{T}(\varphi) & \textsf{t}(\varphi)\\ \textsf{F}(\psi) & \textsf{F}(\psi)\end{array}\\
&&&&\\
\begin{array}{c|c}\multicolumn{2}{c}{\textsf{T}(\neg\varphi)}\\\hline \textsf{t}(\varphi) & \textsf{F}(\varphi)\end{array} & \quad\quad\quad & \begin{array}{c|c}\multicolumn{2}{c}{\textsf{t}(\neg\varphi)}\\\hline \textsf{t}(\varphi) & \textsf{F}(\varphi)\end{array} & \quad\quad\quad & \begin{array}{c}\textsf{F}(\neg\varphi)\\\hline \textsf{T}(\varphi)\end{array}\\
&&&&\\
\begin{array}{c|c}\multicolumn{2}{c}{\textsf{T}(\circ\varphi)}\\\hline \textsf{T}(\varphi) & \textsf{F}(\varphi)\end{array} & \quad\quad\quad & \begin{array}{c|c}\multicolumn{2}{c}{\textsf{t}(\circ\varphi)}\\\hline \textsf{T}(\varphi) & \textsf{F}(\varphi)\end{array} & \quad\quad\quad & \begin{array}{c}\textsf{F}(\circ\varphi)\\\hline \textsf{t}(\varphi)\end{array}
\end{array}$
\end{table}

The closure conditions are the adequate translations of the ones in $\mathbb{T}_{1}$: a branch is closed if it contains $\textsf{L}(\varphi)$ and $\textsf{L}^{*}(\varphi)$ with $\textsf{L}\neq\textsf{L}^{*}$ or if it contains $\textsf{t}(\psi\wedge\neg\psi)$.

For $\CILA$, the situation is a bit simpler: one takes the calculus $\mathbb{T}_{1}$ now with formulas over the signature $\Sigma_{\textbf{LFI}}$; exchange the labels $\textsf{T}_{1}$, $\textsf{t}^{1}_{0}$ and $\textsf{F}_{1}$ by, respectively, $\textsf{T}$, $\textsf{t}$ and $\textsf{F}$; and adds the following rules governing the connective $\circ$.

\begin{table}[H]
\centering
$\begin{array}{ccc}
\begin{array}{c|c}\multicolumn{2}{c}{\textsf{T}(\circ\varphi)}\\\hline \textsf{T}(\varphi) & \textsf{F}(\varphi)\end{array} & \quad\quad\quad & \begin{array}{c}\textsf{F}(\circ\varphi)\\\hline \textsf{t}(\varphi)\end{array}
\end{array}$
\end{table}

The closure conditions are the same as those of $\mathbb{T}_{\mbCcl}$, with one addition: a branch is also closed if it contains a labeled formula $\textsf{t}(\circ\psi)$. We will not prove it, but $\mathbb{T}_{\mbCcl}$ and $\mathbb{T}_{\CILA}$ are decision methods for their respective logics.

\newpage
\printbibliography[segment=\therefsegment,heading=subbibliography]
\end{refsegment}

\begin{refsegment}
\defbibfilter{notother}{not segment=\therefsegment}
\setcounter{chapter}{5}
\chapter{Restricted swap structures for da Costa's hierarchy}\label{Chapter6}\label{Chapter 6}

The RNmatrices we defined for $C_{n}$ are swap structures: roughly speaking, they are semantics built by assigning, simultaneously, truth values to a formula $\alpha$ and several other, related, formulas; in our case, $\neg\alpha$, $\alpha^{1}, \dotsc , \alpha^{n-1}$. Most frequently, these truth values are in the two-valued Boolean algebra $\textbf{2}$, but it is possible, in many cases, to instead consider these values in an arbitrary Boolean algebra. One example would be of the $\textbf{LFI}$ known as $\textbf{mbC}$, whose swap structures may be recast over any non-trivial Boolean algebra (\cite{ParLog}); this is made easier by the fact that, in $\textbf{mbC}$'s case, we have only a single, finite Nmatrix. The end result of this process is a class of Nmatrices parameterized by Boolean algebras, which we call a swap structures semantics. This is done in order to produce a wider class of models, what allows one to study a given logic by adapting tools from algebraic logic and model theory to the context of non-deterministic algebras.

However, swap structures have been presented, up to now, only as Nmatrices; given that the logics between $\textbf{mbCcl}$ and $\textbf{Cila}$ are not characterizable by finite Nmatrices, and that swap structures were coined to be decision methods (and henceforth finite), it is easy to see that there haven't been defined swap structures semantics for these logics so far. In Chapter \ref{Chapter5} we changed this by finally defining, not classical swap structures, but their versions as RNmatrices over $\textbf{2}$ for $C_{1}$ and the whole of da Costa's hierarchy; we continue by generalizing, in much the same way swap structures semantics generalize for arbitrary Boolean algebras the finite Nmatrix characterization of $\textbf{mbC}$, the RNmatrices $\mathcal{RM}_{C_{n}}$ to structures defined over any given, non-trivial, Boolean algebra $\mathcal{B}$. As mentioned, this can have important applications in the model theory of $C_{n}$, being the class of RNmatrices obtained a very well-behaved category; but, more importantly, the methodology here outlined offers a new outlook of swap structures altogether.

The first step needed to make this generalization possible is to consider, instead of bivaluations, the concept of $\mathcal{B}$-valuations; the rest of the process is actually rather straightforward, as we have already defined objects such as $B_{n}$, $D_{n}$, $\mathcal{A}_{C_{n}}$ and $\mathcal{F}_{C_{n}}$ in terms that facilitate exchanging $\textbf{2}$ by an arbitrary $\mathcal{B}$: that is, they were only defined in terms of the Boolean operators and the constants $0$ and $1$, present in any Boolean algebra. The only restriction that we make is to non-trivial, also known as non-degenerate, Boolean algebras, where $0\neq 1$.

The research found in this chapter has been submitted as a preprint in \cite{RestrictedSwap}.

\section{$\mathcal{B}$-Valuations generalize bivaluations}\label{B-valuations}

\begin{definition}\label{B-valuation}
Given a Boolean algebra $\mathcal{B}$ with universe $B$, a $\mathcal{B}$-valuation\index{Valuation, $\mathcal{B}$-} for $C_{n}$ is a function $\mathsf{b}:\textbf{F}(\Sigma_{\textbf{C}}, \mathcal{V})\rightarrow B$ satisfying:
\begin{enumerate}
\item[$(V1)$] $\mathsf{b}(\alpha\#\beta)=\mathsf{b}(\alpha)\#\mathsf{b}(\beta)$, for any $\#\in\{\vee, \wedge, \rightarrow\}$;
\item[$(V2)$] ${\sim}\mathsf{b}(\alpha)\leq\mathsf{b}(\neg\alpha)$;
\item[$(V3)$] $\mathsf{b}(\neg\neg\alpha)\leq\mathsf{b}(\alpha)$;
\item[$(V4)_{n}$] $\mathsf{b}(\alpha^{n})={\sim}(\mathsf{b}(\alpha^{n-1})\wedge\mathsf{b}(\neg(\alpha^{n-1})))$;
\item[$(V5)$] $\mathsf{b}(\neg(\alpha^{1}))=\mathsf{b}(\alpha)\wedge\mathsf{b}(\neg\alpha)$;
\item[$(V6)_{n}$] $\mathsf{b}(\alpha^{(n)})\wedge\mathsf{b}(\beta^{(n)})\leq\mathsf{b}((\alpha\#\beta)^{(n)})$, for any $\#\in\{\vee, \wedge, \rightarrow\}$.
\end{enumerate}
\end{definition}

\begin{proposition}\label{Some properties of B-valuations}
Let $\mathcal{B}$ be a Boolean algebra and $\mathsf{b}$ a $\mathcal{B}$-valuation for $C_{n}$: let us denote, for a formula $\alpha$, the $n+1$-tuple $(\mathsf{b}(\alpha), \mathsf{b}(\neg\alpha), \mathsf{b}(\alpha^{1}), \dotsc  , \mathsf{b}(\alpha^{n-1}))$ in $B^{n+1}$ by $z=(z_{1}, \dotsc  , z_{n+1})$.
\begin{enumerate}
\item For any $1\leq k\leq n-1$, 
\[\mathsf{b}(\neg(\alpha^{k}))=\mathsf{b}(\alpha)\wedge\mathsf{b}(\neg\alpha)\wedge\bigwedge_{i=1}^{k-1}\mathsf{b}(\alpha^{i}),\]
meaning at most one among $z_{1}, \dotsc  , z_{n}$ and $z_{n+1}$ equals $0$; furthermore, if $z_{k}=0$, then $z_{i}=1$ for every $k+1\leq i\leq n+1$.

\item For $1\leq k\leq n-1$, $\mathsf{b}(\alpha^{(k)})=\bigwedge_{i=1}^{k}\mathsf{b}(\alpha^{i})$ and
\[\mathsf{b}(\alpha^{n})={\sim}(\mathsf{b}(\alpha)\wedge\mathsf{b}(\neg\alpha)\wedge\bigwedge_{i=1}^{n-1}\mathsf{b}(\alpha^{i})).\]

\item One has $\mathsf{b}(\alpha^{(n)})={\sim}(\mathsf{b}(\alpha)\wedge\mathsf{b}(\neg\alpha))\wedge\bigwedge_{i=1}^{n-1}\mathsf{b}(\alpha^{i})$.

\item We have $\mathsf{b}(\alpha)\vee\mathsf{b}(\neg\alpha)=1$, $(\mathsf{b}(\alpha)\wedge\mathsf{b}(\neg\alpha))\vee\mathsf{b}(\alpha^{1})=1$ and, for every $1<k\leq n-1$,
\[(\mathsf{b}(\alpha)\wedge\mathsf{b}(\neg\alpha)\wedge\bigwedge_{i=1}^{k-1}\mathsf{b}(\alpha^{i}))\vee\mathsf{b}(\alpha^{k})=1.\]
\end{enumerate}
\end{proposition}

\begin{proof}

\begin{enumerate}
\item A simple argument by induction is sufficient: the case $k=1$ is done, from clause $(V5)$ of Definition \ref{B-valuation}; so, suppose for induction hypothesis that $\mathsf{b}(\alpha^{k-1})=\mathsf{b}(\alpha)\wedge\mathsf{b}(\neg\alpha)\wedge\bigwedge_{i=1}^{k-2}\mathsf{b}(\alpha^{i})$, and we have that, again from clause $(V5)$,
\[\mathsf{b}(\neg(\alpha^{k}))=\mathsf{b}(\neg((\alpha^{k-1})^{1}))=\mathsf{b}(\alpha^{k-1})\wedge\mathsf{b}(\neg(\alpha^{k-1}))=\mathsf{b}(\alpha)\wedge\mathsf{b}(\neg\alpha)\wedge\bigwedge_{i=1}^{k-1}\mathsf{b}(\alpha^{i}).\]
If $\mathsf{b}(\alpha)=0$, from $(V2)$ we must have $\mathsf{b}(\neg\alpha)=1$; and from the previous formula for $\mathsf{b}(\neg(\alpha^{k}))$, we find $\mathsf{b}(\neg(\alpha^{1}))=\cdots=\mathsf{b}(\neg(\alpha^{n-1}))=0$, meaning, again from $(V2)$, that $\mathsf{b}(\alpha^{1})=\cdots=\mathsf{b}(\alpha^{n-1})=1$. If $\mathsf{b}(\neg\alpha)=0$, by the same argument one obtains $\mathsf{b}(\neg(\alpha^{1}))=\cdots=\mathsf{b}(\neg(\alpha^{n-1}))=0$ and the desired result.

Finally, if $\mathsf{b}(\alpha^{k})=0$, for a $k+1\leq j\leq n-1$, one obtains $\mathsf{b}(\neg(\alpha^{j}))=\mathsf{b}(\alpha)\wedge\mathsf{b}(\neg\alpha)\wedge\bigwedge_{i=1}^{j-1}\mathsf{b}(\alpha^{i})=0$, and therefore $\mathsf{b}(\alpha^{j})=1$.

\item For the first part, we proceed again by induction: for $k=1$, the result is trivial given the definition $\alpha^{(1)}=\alpha^{1}$; assuming $\mathsf{b}(\alpha^{(k-1)})=\bigwedge_{i=1}^{k-1}\mathsf{b}(\alpha^{i})$, we get 
\[\mathsf{b}(\alpha^{(k)})=\mathsf{b}(\alpha^{k}\wedge\alpha^{(k-1)})=\mathsf{b}(\alpha^{k})\wedge\mathsf{b}(\alpha^{(k-1}))=\bigwedge_{i=1}^{k}\mathsf{b}(\alpha^{i}).\]

Now, remembering $(V4)_{n}$ and the previous item of the proposition,, 
\[\mathsf{b}(\alpha^{n})={\sim}(\mathsf{b}(\alpha^{n-1})\wedge\mathsf{b}(\neg(\alpha^{n-1})))={\sim}(\mathsf{b}(\alpha)\wedge\mathsf{b}(\neg\alpha)\wedge\bigwedge_{i=1}^{n-1}\mathsf{b}(\alpha^{i})).\]

\item From the previous item,
\[\mathsf{b}(\alpha^{(n)})=\mathsf{b}(\alpha^{n}\wedge\alpha^{(n-1)})=\mathsf{b}(\alpha^{n})\wedge\mathsf{b}(\alpha^{(n-1)})={\sim}\big[\mathsf{b}(\alpha)\wedge\mathsf{b}(\neg\alpha)\wedge\bigwedge_{i=1}^{n-1}\mathsf{b}(\alpha^{i})\big]\wedge\bigwedge_{i=1}^{n-1}\mathsf{b}(\alpha^{i})=\]
\[\big[{\sim}(\mathsf{b}(\alpha)\wedge\mathsf{b}(\neg\alpha))\vee{\sim}\bigwedge_{i=1}^{n-1}\mathsf{b}(\alpha^{i})\big]\wedge\bigwedge_{i=1}^{n-1}\mathsf{b}(\alpha^{i})=\]
\[\big[{\sim}(\mathsf{b}(\alpha)\wedge\mathsf{b}(\neg\alpha))\wedge\bigwedge_{i=1}^{n-1}\mathsf{b}(\alpha^{i})\big]\vee\big[{\sim}\bigwedge_{i=1}^{n-1}\mathsf{b}(\alpha^{i})\wedge\bigwedge_{i=1}^{n-1}\mathsf{b}(\alpha^{i})\big]=\]
\[({\sim}(\mathsf{b}(\alpha)\wedge\mathsf{b}(\neg\alpha))\wedge\bigwedge_{i=1}^{n-1}\mathsf{b}(\alpha^{i}))\vee 0={\sim}(\mathsf{b}(\alpha)\wedge\mathsf{b}(\neg\alpha))\wedge\bigwedge_{i=1}^{n-1}\mathsf{b}(\alpha^{i}).\]

\item From $(V2)$, ${\sim}\mathsf{b}(\alpha)\leq\mathsf{b}(\neg\alpha)$, and since $\mathsf{b}(\alpha)\vee{\sim}\mathsf{b}(\alpha)=1$ we find $\mathsf{b}(\alpha)\vee\mathsf{b}(\neg\alpha)=1$; 
from the first item of this proposition, $\mathsf{b}(\neg(\alpha^{1}))=\mathsf{b}(\alpha)\wedge\mathsf{b}(\neg\alpha)$, and since $\mathsf{b}(\alpha)\vee\mathsf{b}(\neg\alpha)=1$ works for any formula $\alpha$, by replacing in the latter equation $\alpha$ with $\alpha^{1}$ one finds 
\[\mathsf{b}(\alpha^{1})\vee(\mathsf{b}(\alpha)\wedge\mathsf{b}(\neg\alpha))=\mathsf{b}(\alpha^{1})\vee\mathsf{b}(\neg(\alpha^{1}))=1;\]
finally, again from the first item of the proposition, for every $1<k\leq n-1$, $\mathsf{b}(\neg(\alpha^{k}))=\mathsf{b}(\alpha)\wedge\mathsf{b}(\neg\alpha)\wedge\bigwedge_{i=1}^{k-1}\mathsf{b}(\alpha^{i})$, and therefore
\[\mathsf{b}(\alpha^{k})\vee(\mathsf{b}(\alpha)\wedge\mathsf{b}(\neg\alpha)\wedge\bigwedge_{i=1}^{k-1}\mathsf{b}(\alpha^{i}))=\mathsf{b}(\alpha^{k})\vee\mathsf{b}(\neg(\alpha^{k}))=1.\]
\end{enumerate}
\end{proof}

\begin{proposition}\label{Negations of B-valuations}
Given a Boolean algebra $\mathcal{B}$, let $\mathsf{b}$ be a $\mathcal{B}$-valuation for $C_{n}$.
\begin{enumerate}
\item For any $\#\in\{\vee, \wedge, \rightarrow\}$, if $\mathsf{b}(\neg\alpha)={\sim}\mathsf{b}(\alpha)$ and $\mathsf{b}(\neg\beta)={\sim}\mathsf{b}(\beta)$, then $\mathsf{b}(\neg(\alpha\#\beta))={\sim}\mathsf{b}(\alpha\#\beta)$;
\item if $\mathsf{b}(\neg\alpha)={\sim}\mathsf{b}(\alpha)$, then $\mathsf{b}(\neg\neg\alpha)={\sim}\mathsf{b}(\neg\alpha)$.
\end{enumerate}
\end{proposition}

\begin{proof}
\begin{enumerate}
\item Assume $\mathsf{b}$ is a $\mathcal{B}$-valuation for $C_{n}$ for which $\mathsf{b}(\neg\alpha)={\sim}\mathsf{b}(\alpha)$ and $\mathsf{b}(\neg\beta)={\sim}\mathsf{b}(\beta)$. Then $\mathsf{b}(\alpha)\wedge\mathsf{b}(\neg\alpha)=\mathsf{b}(\beta)\wedge\mathsf{b}(\neg\beta)=0$ and therefore ${\sim}(\mathsf{b}(\alpha)\wedge\mathsf{b}(\neg\alpha))={\sim}(\mathsf{b}(\beta)\wedge\mathsf{b}(\neg\beta))=1$.

Since, from Proposition \ref{Some properties of B-valuations}, $\mathsf{b}(\neg(\alpha^{k}))=\mathsf{b}(\alpha)\wedge\mathsf{b}(\neg\alpha)\wedge\bigwedge_{i=1}^{k-1}\mathsf{b}(\alpha^{i})$ for every $1\leq k\leq n-1$, we find $\mathsf{b}(\neg(\alpha^{k}))=\mathsf{b}(\neg(\beta^{k}))=0$ for every $1\leq k\leq n-1$, implying from $(V2)$ (in Definition \ref{B-valuation}) that $\mathsf{b}(\alpha^{k})=\mathsf{b}(\beta^{k})=1$ for all $1\leq k\leq n-1$. Again from Proposition \ref{Some properties of B-valuations}, $\mathsf{b}(\alpha^{(n)})={\sim}(\mathsf{b}(\alpha)\wedge\mathsf{b}(\neg\alpha))\wedge\bigwedge_{i=1}^{n-1}\mathsf{b}(\alpha^{i})$, and so $\mathsf{b}(\alpha^{(n)})=\mathsf{b}(\beta^{(n)})=1$, meaning from $(V6)_{n}$ that $\mathsf{b}((\alpha\#\beta)^{(n)})=1$.

But, again from the fact that, for any formula $\alpha$, $\mathsf{b}(\alpha^{(n)})={\sim}(\mathsf{b}(\alpha)\wedge\mathsf{b}(\neg\alpha))\wedge\bigwedge_{i=1}^{n-1}\mathsf{b}(\alpha^{i})$, we obtain that $\mathsf{b}((\alpha\#\beta)^{(n)})\leq{\sim}[\mathsf{b}(\alpha\#\beta)\wedge\mathsf{b}(\neg(\alpha\#\beta))]$, and thus $\mathsf{b}(\alpha\#\beta)\wedge\mathsf{b}(\neg(\alpha\#\beta))=0$, which means $\mathsf{b}(\neg(\alpha\#\beta))={\sim}\mathsf{b}(\alpha\#\beta)$ from $(V2)$.

\item Now suppose $\mathsf{b}$ is a $\mathcal{B}$-valuation for $C_{n}$ for which $\mathsf{b}(\neg\alpha)={\sim}\mathsf{b}(\alpha)$: by $(V3)$, 
\[\mathsf{b}(\neg\neg\alpha)\wedge\mathsf{b}(\neg\alpha)\leq\mathsf{b}(\alpha)\wedge\mathsf{b}(\neg\alpha)=0,\]
and since $\mathsf{b}(\neg\neg\alpha)\vee\mathsf{b}(\neg\alpha)=1$ from $(V2)$, we find $\mathsf{b}(\neg\neg\alpha)={\sim}\mathsf{b}(\neg\alpha)$.

\end{enumerate}
\end{proof}

One notices that, as a corollary to this last proposition, whenever $\mathcal{V}_{0}$ is a non-empty subset of $\mathcal{V}$ and $\mathsf{b}$ is a $\mathcal{B}$-valuation for $C_{n}$ satisfying that $\mathsf{b}(\neg p)={\sim}\mathsf{b}(p)$, for every $p\in\mathcal{V}_{0}$, one has $\mathsf{b}(\neg\alpha)={\sim}\mathsf{b}(\alpha)$ for every formula $\alpha$ in $F(\Sigma_{\textbf{C}}, \mathcal{V}_{0})$.

A simple argument by structural induction is sufficient, being the base case one of the hypothesis: if the result holds for $\alpha$ and $\beta$, from Proposition \ref{Negations of B-valuations} $\mathsf{b}(\neg(\alpha\#\beta))={\sim}\mathsf{b}(\alpha\#\beta)$, for any $\#\in\{\vee, \wedge, \rightarrow\}$; and, again by the proposition, $\mathsf{b}(\neg\neg\alpha)={\sim}\mathsf{b}(\neg\alpha)$. Since $\{\alpha\in F(\Sigma_{\textbf{C}}, \mathcal{V}_{0})\ :\  \mathsf{b}(\neg\alpha)={\sim}\mathsf{b}(\alpha)\}$ contains $\mathcal{V}_{0}$ and is closed under the connectives in $\{\neg, \vee, \wedge, \rightarrow\}$ (all of $\Sigma_{\textbf{C}}$), we obtain this set is the whole of $F(\Sigma_{\textbf{C}}, \mathcal{V}_{0})$.

 Additionally, under those hypothesis, we have that for every $\alpha\in F(\Sigma_{\textbf{C}}, \mathcal{V}_{0})$ we have $\mathsf{b}(\alpha^{(n)})=1$: we can see that from the equality $\mathsf{b}(\alpha^{(n)})={\sim}(\mathsf{b}(\alpha)\wedge\mathsf{b}(\neg\alpha))\wedge\bigwedge_{i=1}^{n-1}\mathsf{b}(\alpha^{i})$. Since $\mathsf{b}(\neg\alpha)={\sim}\mathsf{b}(\alpha)$, $\mathsf{b}(\alpha)\wedge\mathsf{b}(\neg\alpha)=0$ and therefore ${\sim}(\mathsf{b}(\alpha)\wedge\mathsf{b}(\neg\alpha))=1$; furthermore, for any $1\leq k\leq n-1$, $\mathsf{b}(\neg(\alpha^{k}))={\sim}\mathsf{b}(\alpha^{k})$ and $\mathsf{b}(\neg(\alpha^{k}))=\mathsf{b}(\alpha)\wedge\mathsf{b}(\neg\alpha)\wedge\bigwedge_{i=1}^{k-1}\mathsf{b}(\alpha^{i})$, thus we find
 \[\mathsf{b}(\alpha^{k})={\sim}(\mathsf{b}(\alpha)\wedge\mathsf{b}(\neg\alpha)\wedge\bigwedge_{i=1}^{k-1}\mathsf{b}(\alpha^{i})),\]
 and in our current position this means $\mathsf{b}(\alpha^{k})=1$ (for all $1\leq k\leq n-1$), what gives the desired result $\mathsf{b}(\alpha^{(n)})=1$.

 \begin{lemma}\label{B-valuations on n-consistency}
Given a bivaluation $\mathsf{b}$ for $C_{n}$, $\mathsf{b}(\alpha^{(n)})=1$ if and only if $\mathsf{b}(\alpha)\neq\mathsf{b}(\neg\alpha)$.
 \end{lemma}
 
 \begin{proof}
 We first prove that, if $\mathsf{b}(\alpha^{(n)})=1$, then $\mathsf{b}(\alpha)\neq\mathsf{b}(\neg\alpha)$. Start by noticing that if $\mathsf{b}$ is a bivaluation for $C_{n}$, then $\mathsf{b}(\alpha^{(k)})=\bigwedge_{i=1}^{k}\mathsf{b}(\alpha^{i})$ for $k\geq 1$: this is obviously true for $k=1$, given $\alpha^{(1)}=\alpha^{1}$; supposing this is true for $k-1$, we obtain that
 \[\mathsf{b}(\alpha^{(k)})=\mathsf{b}(\alpha^{(k-1)}\wedge\alpha^{k})=\mathsf{b}(\alpha^{(k-1)})\wedge\mathsf{b}(\alpha^{k})=\bigwedge_{i=1}^{k}\mathsf{b}(\alpha^{i}).\]
So, since $\mathsf{b}(\alpha^{(n)})=\bigwedge_{i=1}^{n}\mathsf{b}(\alpha^{i})=1$, $\mathsf{b}(\alpha^{1})=\cdots=\mathsf{b}(\alpha^{n})=1$; from $(B6)_{n}$ of Definition \ref{B-valuation}, $\mathsf{b}(\alpha^{n})=1$ implies $\mathsf{b}(\neg(\alpha^{n-1}))=0$.

Now, from $(B7)$, since $\mathsf{b}(\alpha^{n-1})=\mathsf{b}((\alpha^{n-2})^{1})=0$ one finds out that $\mathsf{b}(\neg(\alpha^{n-2}))=0$; inductively, we obtain $\mathsf{b}(\neg(\alpha^{1}))=\cdots=\mathsf{b}(\neg(\alpha^{n-1}))=0$, and with one final application of $(B7)$, $\mathsf{b}(\neg(\alpha^{1}))=0$ implies $\mathsf{b}(\alpha)\neq\mathsf{b}(\neg\alpha)$.

Reciprocally, $\mathsf{b}(\alpha)\neq\mathsf{b}(\neg\alpha)$ implies $\mathsf{b}(\alpha)\wedge\mathsf{b}(\neg\alpha)=0$ and so $\mathsf{b}(\alpha^{1})=1$, while by $(B7)$ we obtain that $\mathsf{b}(\neg(\alpha^{1}))=0$; this, of course, implies that $\mathsf{b}(\alpha^{2})=\cdots=\mathsf{b}(\alpha^{n})=1$, and therefore $\mathsf{b}(\alpha^{(n)})=\bigwedge_{i=1}^{n}\mathsf{b}(\alpha^{i})=1$.
\end{proof}

\begin{theorem}
Given the Boolean algebra of two elements $\textbf{2}$, a function $\mathsf{b}:F(\Sigma_{\textbf{C}}, \mathcal{V})\rightarrow\{0,1\}$ is a bivaluation for $C_{n}$ if, and only if, it is a $\textbf{2}$-valuation for $C_{n}$.
\end{theorem}

\begin{proof}
Suppose, first of all, that $\mathsf{b}$ is a bivaluation. Condition $(V1)$ follows from $(B1)$, $(B2)$ and $(B3)$. If $\mathsf{b}(\alpha)=0$, from $(B4)$ we find $\mathsf{b}(\neg\alpha)=1$ and therefore ${\sim}\mathsf{b}(\alpha)=1\leq 1=\mathsf{b}(\neg\alpha)$; if $\mathsf{b}(\alpha)=1$, regardless of the value of $\mathsf{b}(\neg\alpha)$ we obtain ${\sim}\mathsf{b}(\alpha)=0\leq \mathsf{b}(\neg\alpha)$, which means that, in all cases, ${\sim}\mathsf{b}(\alpha)\leq\mathsf{b}(\neg\alpha)$, corresponding to clause $(V2)$.

If $\mathsf{b}(\neg\neg\alpha)=1$, condition $(B5)$ implies that $\mathsf{b}(\alpha)=1$, and therefore $\mathsf{b}(\neg\neg\alpha)\leq\mathsf{b}(\alpha)$; if $\mathsf{b}(\neg\neg\alpha)=0$, $\mathsf{b}(\neg\neg\alpha)\leq\mathsf{b}(\alpha)$ is true regardless of the value of $\mathsf{b}(\alpha)$, meaning condition $(V3)$ always holds. Now,  from condition $(B4)$ we can not have $\mathsf{b}(\alpha^{n-1})=\mathsf{b}(\neg(\alpha^{n-1}))=0$, so there remain three cases to consider: if $\mathsf{b}(\alpha^{n-1})=0$ and $\mathsf{b}(\neg(\alpha^{n-1}))=1$, or vice-versa, from $(B6)_{n}$ we obtain $\mathsf{b}(\alpha^{n})=1$ and therefore 
\[\mathsf{b}(\alpha^{n})=1={\sim}(\mathsf{b}(\alpha^{n-1})\wedge\mathsf{b}(\neg(\alpha^{n-1})));\]
if $\mathsf{b}(\alpha^{n-1})=\mathsf{b}(\neg(\alpha^{n-1}))=1$, again from $(B6)_{n}$ we obtain $\mathsf{b}(\alpha^{n})=0$, and so
\[\mathsf{b}(\alpha^{n})=0={\sim}(\mathsf{b}(\alpha^{n-1})\wedge\mathsf{b}(\neg(\alpha^{n-1}))),\]
proving that $(V4)_{n}$ is valid. As mentioned earlier, we can not have $\mathsf{b}(\alpha)=\mathsf{b}(\neg\alpha)=0$, according to $(B4)$, and so if $\mathsf{b}(\alpha)=\mathsf{b}(\neg\alpha)=1$, from $(B7)$ we get $\mathsf{b}(\neg(\alpha^{1}))=1$ and therefore $\mathsf{b}(\neg(\alpha^{1}))=\mathsf{b}(\alpha)\wedge\mathsf{b}(\neg\alpha)$; if $\mathsf{b}(\alpha)=0$ or $\mathsf{b}(\neg\alpha)=0$, again from $(B7)$ we find $\mathsf{b}(\neg(\alpha^{1}))=0$, meaning once more that $\mathsf{b}(\neg(\alpha^{1}))=\mathsf{b}(\alpha)\wedge\mathsf{b}(\neg\alpha)$ and that $(V5)$ holds. 

Finally, notice that if $\mathsf{b}(\alpha^{(n)})=0$ or $\mathsf{b}(\beta^{(n)})=0$, $(V6)_{n}$ already holds true, so let us assume that $\mathsf{b}(\alpha^{(n)})=1$ and $\mathsf{b}(\beta^{(n)})=1$: this means, by Lemma \ref{B-valuations on n-consistency}, that $\mathsf{b}(\alpha)\neq \mathsf{b}(\neg\alpha)$ and $\mathsf{b}(\beta)\neq\mathsf{b}(\neg\beta)$; from $(B8)$, this translates to $\mathsf{b}(\alpha\#\beta)\neq\mathsf{b}(\neg(\alpha\#\beta))$ for any $\#\in\{\vee, \wedge, \rightarrow\}$, and again by Lemma \ref{B-valuations on n-consistency} this becomes $\mathsf{b}((\alpha\#\beta)^{(n)})=1$, which verifies $(V6)_{n}$ and finishes proving $\mathsf{b}$ is a $\textbf{2}$-valuation.

Reciprocally, suppose $\mathsf{b}$ is a $\textbf{2}$-valuation. $(B1)$, $(B2)$ and $(B3)$ clearly follow from $(V1)$. If $\mathsf{b}(\alpha)=0$, from $(V2)$ we obtain $\mathsf{b}(\neg\alpha)\geq{\sim}\mathsf{b}(\alpha)=1$, and therefore $\mathsf{b}(\neg\alpha)=1$, thus satisfying $(B4)$. If $\mathsf{b}(\neg\neg\alpha)=1$, $\mathsf{b}(\alpha)\geq\mathsf{b}(\neg\neg\alpha)=1$ from $(V3)$, and so $\mathsf{b}(\alpha)=1$, meaning $(B5)$ is also validated.

Now, if $\mathsf{b}(\alpha^{n-1})=\mathsf{b}(\neg(\alpha^{n-1}))$, since both can not be $0$ given $(V2)$, both equal $1$ and from $(V4)_{n}$ we get $\mathsf{b}(\alpha^{n})={\sim}(\mathsf{b}(\alpha^{n-1})\wedge\mathsf{b}(\neg(\alpha^{n-1})))=0$, and one direction of $(B6)_{n}$ holds. Reciprocally, if $\mathsf{b}(\alpha^{n})=0$, given $\mathsf{b}(\alpha^{n})={\sim}(\mathsf{b}(\alpha^{n-1})\wedge\mathsf{b}(\neg(\alpha^{n-1})))$ from $(V4)_{n}$ we get $\mathsf{b}(\alpha^{n-1})\wedge\mathsf{b}(\neg(\alpha^{n-1}))=1$, and therefore $\mathsf{b}(\alpha^{n-1})=\mathsf{b}(\neg(\alpha^{n-1}))=1$, meaning that both directions of $(B6)_{n}$ hold.

To analyze $(B7)$, there are again two directions to consider: if $\mathsf{b}(\alpha)=\mathsf{b}(\neg\alpha)$, this means that both are equal to $1$ (given $(V2)$), and so $\mathsf{b}(\neg(\alpha^{1}))=\mathsf{b}(\alpha)\wedge\mathsf{b}(\neg\alpha)=1$, in accordance to $(V5)$; if $\mathsf{b}(\neg(\alpha^{1}))=1$, again from $(V5)$ we obtain $1=\mathsf{b}(\neg(\alpha^{1}))=\mathsf{b}(\alpha)\wedge\mathsf{b}(\neg\alpha)$, meaning $\mathsf{b}(\alpha)=\mathsf{b}(\neg\alpha)=1$ and that $(B7)$ holds true.

Finally, suppose $\mathsf{b}(\alpha)\neq\mathsf{b}(\neg\alpha)$ and $\mathsf{b}(\beta)\neq\neg\mathsf{b}(\beta)$, meaning that $\mathsf{b}(\neg\alpha)={\sim}\mathsf{b}(\alpha)$ and $\mathsf{b}(\neg\beta)={\sim}\mathsf{b}(\beta)$: from Proposition \ref{Negations of B-valuations}, $\mathsf{b}(\neg(\alpha\#\beta))={\sim}\mathsf{b}(\alpha\#\beta)$, for any $\#\in\{\vee, \wedge, \rightarrow\}$, and so $\mathsf{b}(\neg(\alpha\#\beta))\neq \mathsf{b}(\alpha\#\beta)$. This implies $(B8)$ is true, and that $\mathsf{b}$ is then a bivaluation.
\end{proof}

Fix a Boolean algebra $\mathcal{B}$ and take formulas $\Gamma\cup\{\varphi\}\subseteq F(\Sigma_{\textbf{C}}, \mathcal{V})$: we say $\Gamma$ proves $\varphi$ according to $\mathcal{B}$-valuations for $C_{n}$, and write $\Gamma\vDash_{C_{n}}^{\mathcal{B}}\varphi$\label{vDashBCn}, if every $\mathcal{B}$-valuation $\mathsf{b}$ satisfying $\mathsf{b}(\Gamma)\subseteq\{1\}$ also has the property that $\mathsf{b}(\varphi)=1$.

\begin{theorem}
Let $n\in\mathbb{N}\setminus\{0\}$ and $\Gamma\cup\{\varphi\}$ be a set of formulas in the signature $\Sigma_{\textbf{C}}$: then $\Gamma\vdash_{C_{n}}\varphi$ if and only if $\Gamma\vDash_{C_{n}}^{\mathcal{B}}\varphi$ for every Boolean algebra $\mathcal{B}$.
\end{theorem}

\begin{proof}
First assume that $\Gamma\vdash_{C_{n}}\varphi$, and let us fix a Boolean algebra $\mathcal{B}$. It is a straightforward exercise to prove that, if $\varphi$ is an instance of an axiom of $C_{n}$, then $\mathsf{b}(\varphi)=1$ for any $\mathcal{B}$-valuation $\mathsf{b}$ for $C_{n}$; now, again for a $\mathcal{B}$-valuation $\mathsf{b}$ for $C_{n}$, if $\mathsf{b}(\alpha)=\mathsf{b}(\alpha\rightarrow\beta)=1$, since $\mathsf{b}(\alpha\rightarrow\beta)=\mathsf{b}(\alpha)\rightarrow\mathsf{b}(\beta)$ from $(V1)$, one derives that $\mathsf{b}(\beta)=1$, meaning ``$\vDash_{C_{n}}^{\mathcal{B}}$'' models all axiom schemata and rules of inference of $C_{n}$.

So, if $\varphi_{1}, \dotsc  , \varphi_{n}$, with $\varphi_{n}=\varphi$, is a derivation of $\varphi$ from $\Gamma$, we prove by induction that, if $\mathsf{b}(\Gamma)\subseteq\{1\}$, then $\mathsf{b}(\varphi_{1})=\cdots=\mathsf{b}(\varphi_{n})=1$, and therefore if $\mathsf{b}$ is a $\mathcal{B}$-valuation for $C_{n}$, $\mathsf{b}(\Gamma)\subseteq\{1\}$ implies $\mathsf{b}(\varphi)=1$, meaning $\Gamma\vDash_{C_{n}}^{\mathcal{B}}\varphi$. This is rather easy: $\varphi_{1}$ is either in $\Gamma$, where the fact that $\mathsf{b}(\Gamma)\subseteq\{1\}$ implies $\mathsf{b}(\varphi_{1})=1$, or is an instance of an axiom schema, meaning from the previous observations that $\mathsf{b}(\varphi_{1})=1$. So, assume that $\mathsf{b}(\varphi_{1})=\cdots=\mathsf{b}(\varphi_{i-1})=1$, and there are three cases to consider: if $\varphi_{i}$ is in $\Gamma$ or is an instance of an axiom schema, proceed as before; if there exist $\varphi_{j}$ and $\varphi_{k}$, with $1\leq j, k<i$, such that $\varphi_{j}=\varphi_{k}\rightarrow\varphi_{i}$ or $\varphi_{k}=\varphi_{j}\rightarrow\varphi_{i}$, since $\mathsf{b}(\varphi_{j})=\mathsf{b}(\varphi_{k})=1$ we find $\mathsf{b}(\varphi_{i})=1$, and we are done.

Reciprocally, assume $\Gamma\vDash_{C_{n}}^{\mathcal{B}}\varphi$ for every Boolean algebra $\mathcal{B}$; in particular, $\Gamma\vDash_{C_{n}}^{\textbf{2}}\varphi$, and since bivaluations characterize $C_{n}$, we derive that $\Gamma\vdash_{C_{n}}\varphi$.
\end{proof}

\section{The RNmatrices $\mathcal{RM}^{\mathcal{B}}_{C_{n}}$}\label{RNmatrices RMBCn}

Now that we have the appropriate generalization of a bivaluation to arbitrary Boolean algebras, our goal is to generalize the RNmatrix $\mathcal{RM}_{C_{n}}=(\mathcal{A}_{C_{n}}, D_{n}. \mathcal{F}_{C_{n}})$ to a new restricted Nmatrix built around $\mathcal{B}$ by using the notion of snapshots in a broader context, in much the same way it would be done while dealing with swap structures. This will be achieved by generalizing both $\mathcal{A}_{C_{n}}$ and $D_{n}$ to, respectively, the $\Sigma_{\textbf{C}}$-multialgebra $\mathcal{A}_{C_{n}}^{\mathcal{B}}$ and the subset of its universe $D_{n}^{\mathcal{B}}$; this, of course, takes the Nmatrix $\mathcal{M}_{C_{n}}=(\mathcal{A}_{C_{n}}, D_{n})$ to the, still, Nmatrix $\mathcal{M}_{C_{n}}^{\mathcal{B}}=(\mathcal{A}_{C_{n}}^{\mathcal{B}}, D_{n}^{\mathcal{B}})$. The decisive step, however, is generalizing the set of restricted valuations from $\mathcal{F}_{C_{n}}$ to $\mathcal{F}_{C_{n}}^{\mathcal{B}}$, what will allow to trade the RNmatrix $\mathcal{RM}_{C_{n}}=(\mathcal{A}_{C_{n}}, D_{n}, \mathcal{F}_{C_{n}})$ for $\mathcal{RM}_{C_{n}}^{\mathcal{B}}=(\mathcal{A}_{C_{n}}^{\mathcal{B}}, D_{n}^{\mathcal{B}}, \mathcal{F}_{C_{n}}^{\mathcal{B}})$.

For a Boolean algebra $\mathcal{B}$, consider a $\mathcal{B}$-valuation $\mathsf{b}$ for $C_{n}$ and the $(n+1)$-tuple $z=(z_{1}, \dotsc  , z_{n+1})$ in $\mathcal{B}^{n+1}$ such that $z_{1}=\mathsf{b}(\alpha)$, $z_{2}=\mathsf{b}(\neg\alpha)$, $z_{3}=\mathsf{b}(\alpha^{1}), \dotsc  , z_{n+1}=\mathsf{b}(\alpha^{n-1})$, for a formula $\alpha$ of $C_{n}$; these will be the aforementioned generalized snapshots over which we shall work. We can be sure that the following definition works by referring back to Proposition \ref{Some properties of B-valuations}.

\begin{definition}
Given a Boolean algebra $\mathcal{B}$, the set of $\mathcal{B}$-snapshots for $C_{n}$ is\label{BnB}
\[B_{n}^{\mathcal{B}}=\{z\in B^{n+1}\ :\  (\bigwedge_{i=1}^{k}z_{i})\vee z_{k+1}=1,\quad \text{for every} \quad 1\leq k\leq n\}.\]
\end{definition}

Other sets that will be useful to us are:
\begin{enumerate}
\item $D_{n}^{\mathcal{B}}=\{z\in B_{n}^{\mathcal{B}}\ :\  z_{1}=1\}$, the set of designated values;\label{DnB}
\item $Boo_{n}^{\mathcal{B}}=\{z\in B_{n}^{\mathcal{B}}\ :\  z_{1}\wedge z_{2}=0\}$, the set of Boolean values;\label{BoonB}
\item $U_{n}^{\mathcal{B}}=B_{n}^{\mathcal{B}}\setminus D_{n}^{\mathcal{B}}$, the set of undesignated values.
\end{enumerate}

Although somewhat obvious, it is still relevant to point out the set of designated elements correspond to formulas which are true under a $\mathcal{B}$-valuation: that is, if $z=(z_{1}, \dotsc  , z_{n+1})$ satisfies $z_{1}=1$, and $z=(\mathsf{b}(\alpha), \mathsf{b}(\neg\alpha), \mathsf{b}(\alpha^{1}), \dotsc  , \mathsf{b}(\alpha^{n-1}))$, then $\mathsf{b}(\alpha)=1$, that is, $\alpha$ is true according to $\mathsf{b}$. For simplicity, and to maintain the notation we had already adopted, we will denote the $i$-th coordinate of $z\in B^{n+1}$ by $z_{i}$.

Notice that, from the definition of $B_{n}^{\mathcal{B}}$, for any snapshot $z$ one has $z_{1}\vee z_{2}=1$; so, if $z\in Boo_{n}^{\mathcal{B}}$, $z_{1}\wedge z_{2}=0$ and therefore $z_{2}={\sim}z_{1}$; then, from the fact that
\[1=(\bigwedge_{i=1}^{k}z_{i})\vee z_{k+1}=(z_{1}\wedge z_{2}\wedge\bigwedge_{i=3}^{k}z_{i})\vee z_{k+1}=(0\wedge\bigwedge_{i=3}^{k}z_{i})\vee z_{k+1}=0\vee z_{k+1}=z_{k+1},\]
for every $3\leq k\leq n$, we discover that $z\in Boo_{n}^{\mathcal{B}}$ if, and only if, there exists $a$ in $\mathcal{B}$ such that $z=(a, {\sim}a, 1, \dotsc  , 1)$, what establishes a bijection between $\mathcal{B}$ and $Boo_{n}^{\mathcal{B}}$.

\begin{definition}
Given a Boolean algebra $\mathcal{B}$, the full swap structure for $C_{n}$ over $\mathcal{B}$ is the $\Sigma_{\textbf{C}}$-multialgebra $\mathcal{A}_{C_{n}}^{\mathcal{B}}$\label{ACnB}, with universe $B_{n}^{\mathcal{B}}$, with operations defined as, for any two $z, w\in B_{n}^{\mathcal{B}}$,
\[\tilde{\neg}z=\{w\in B_{n}^{\mathcal{B}}\ :\  w_{1}=z_{2}\quad\text{and}\quad w_{2}\leq z_{1}\}\]
and
\[\text{for}\quad \#\in\{\vee, \wedge, \rightarrow\},\quad z\tilde{\#}w=\begin{cases*}
\{u\in B_{n}^{\mathcal{B}}\ :\  u_{1}=z_{1}\# w_{1}\}\cap Boo_{n}^{\mathcal{B}} & if $z, w\in Boo_{n}^{\mathcal{B}}$\\
\{u\in B_{n}^{\mathcal{B}}\ :\  u_{1}=z_{1}\# w_{1}\} & otherwise
\end{cases*}.\]
\end{definition}

Notice that this endows $Boo_{n}^{\mathcal{B}}$ with a structure of a Boolean algebra, top element $(1, 0, 1, \dotsc  , 1)$ and bottom $(0, 1, 1, \dotsc  , 1)$, isomorphic to that of $\mathcal{B}$, giving us a copy of $\mathcal{B}$ inside of $\mathcal{A}_{C_{n}}^{\mathcal{B}}$: in all of this, ``$\tilde{\neg}$'' plays the role of the negation of $Boo_{n}^{\mathcal{B}}$ as a Boolean algebra, meaning $\tilde{\neg}z=\{({\sim}a, a, 1, \dotsc  , 1)\}$, if $z=(a, {\sim}a, 1, \dotsc  , 1)$. This is true since, appealing to its definition, 
\[\tilde{\neg}z=\{w\in B_{n}^{\mathcal{B}}\ :\  w_{1}={\sim}a\quad\text{and}\quad w_{2}\leq a\},\]
what means that $w_{1}\wedge w_{2}\leq {\sim}a\wedge a=0$, and therefore $w\in Boo_{n}^{\mathcal{B}}$ is such that $w_{1}={\sim}a$, and there is only one element satisfying this property, namely $w=({\sim}a, a, 1, \dotsc  , 1)$.

Furthermore, one sees that $B_{n}\subseteq B_{n}^{\mathcal{B}}$, given that $\textbf{2}$ is necessarily a subalgebra of $\mathcal{B}$, and $D_{n}^{\mathcal{B}}\cap B_{n}=D_{n}$ and $Boo_{n}^{\mathcal{B}}\cap B_{n}=Boo_{n}$; the multioperations of $\mathcal{A}_{C_{n}}^{\mathcal{B}}$, when restricted to $B_{n}$, become those of $\mathcal{A}_{C_{n}}$, and it is clear then that $\mathcal{A}_{C_{n}}$ is a submultialgebra of $\mathcal{A}_{C_{n}}^{\mathcal{B}}$.

\begin{definition}\label{restricted valuations for RMBCn}
Given a Boolean algebra $\mathcal{B}$, the restricted Nmatrix for $C_{n}$ over $\mathcal{B}$ is the RNmatrix\label{RMCnB} 
\[\mathcal{RM}_{C_{n}}^{\mathcal{B}}=(\mathcal{A}_{C_{n}}^{\mathcal{B}}. D_{n}^{\mathcal{B}}, \mathcal{F}_{C_{n}}^{\mathcal{B}}),\]
where $\mathcal{F}_{C_{n}}^{\mathcal{B}}$\label{FCnB} is the set of valuations $\nu$ for $\mathcal{A}_{C_{n}}^{\mathcal{B}}$, meaning homomorphisms $\nu:\textbf{F}(\Sigma_{\textbf{C}}, \mathcal{V})\rightarrow \mathcal{A}_{C_{n}}^{\mathcal{B}}$, satisfying, for any two formulas $\alpha$ and $\beta$ of $C_{n}$:
\begin{enumerate}
\item $\nu(\alpha\wedge\neg\alpha)\in\{z\in \nu(\alpha)\tilde{\wedge}\nu(\neg\alpha)\ :\  z_{2}=\nu(\alpha)_{3}\}$;
\item $\nu(\alpha^{1})=(\nu(\alpha)_{3}, \nu(\alpha)_{1}\wedge\nu(\alpha)_{2}, \nu(\alpha)_{4}, \dotsc  , \nu(\alpha)_{n+1}, {\sim}(\bigwedge_{i=1}^{n+1}\nu(\alpha)_{i}))$;
\item for any $\#\in\{\vee, \wedge, \rightarrow\}$, $\nu((\alpha^{(n)}\wedge\beta^{(n)})\rightarrow(\alpha\#\beta)^{(n)})\in D_{n}^{\mathcal{B}}$.
\end{enumerate}
\end{definition}

The following is a somewhat obvious result, but we prove it for completeness sake.

\begin{lemma}\label{index formula}
For every formula $\alpha$ of $C_{n}$, $j, k\in\mathbb{N}$ and endomorphism $\sigma:\textbf{F}(\Sigma_{\textbf{C}}, \mathcal{V})\rightarrow\textbf{F}(\Sigma_{\textbf{C}}, \mathcal{V})$,
\[\sigma(\alpha^{j})=\sigma(\alpha)^{j}\quad\text{and}\quad\sigma(\alpha^{(k)})=\sigma(\alpha)^{(k)}.\]
\end{lemma}

\begin{proof}
One has $\sigma(\alpha^{0})=\sigma(\alpha)=\sigma(\alpha)^{0}$, so assume that $\sigma(\alpha^{j})=\sigma(\alpha)^{j}$:
\[\sigma(\alpha^{j+1})=\sigma(\neg(\alpha^{j}\wedge\neg\alpha^{j}))=\neg(\sigma(\alpha^{j})\wedge\neg\sigma(\alpha^{j}))=\neg(\sigma(\alpha)^{j}\wedge\neg\sigma(\alpha)^{j})=\sigma(\alpha)^{j+1}.\]
Now $\sigma(\alpha^{(0)})=\sigma(\alpha)=\sigma(\alpha)^{(0)}$ and $\sigma(\alpha^{(1)})=\sigma(\alpha^{1})=\sigma(\alpha)^{1}=\sigma(\alpha)^{(0)}$, so assume $\sigma(\alpha^{(k)})=\sigma(\alpha)^{(k)}$:
\[\sigma(\alpha^{(k+1)})=\sigma(\alpha^{(k)}\wedge\alpha^{k+1})=\sigma(\alpha^{(k)})\wedge\sigma(\alpha^{k+1})=\sigma(\alpha)^{(k)}\wedge\sigma(\alpha)^{k+1}=\sigma(\alpha)^{(k+1)}.\]
\end{proof}

\begin{proposition}
For any homomorphism $\nu$ in $\mathcal{F}_{C_{n}}^{\mathcal{B}}$ and any endomorphism $\sigma:\textbf{F}(\Sigma_{\textbf{C}}, \mathcal{V})\rightarrow\textbf{F}(\Sigma_{\textbf{C}}, \mathcal{V})$, $\nu\circ\sigma$ lies in $\mathcal{F}_{C_{n}}^{\mathcal{B}}$.
\end{proposition}

\begin{proof}
Of course $\nu\circ\sigma:\textbf{F}(\Sigma_{\textbf{C}}, \mathcal{V})\rightarrow\mathcal{A}_{C_{n}}^{\mathcal{B}}$ is still a homomorphism of $\Sigma_{\textbf{C}}$-multialgebras.

\begin{enumerate}
\item That $\nu\circ\sigma(\alpha\wedge\neg\alpha)$ is in $\nu\circ\sigma(\alpha)\tilde{\wedge}\nu\circ\sigma(\neg\alpha)$ is a given, since $\nu\circ\sigma$ is a homomorphism. Now, given $\sigma$ is a homomorphism of $\textbf{F}(\Sigma_{\textbf{C}}, \mathcal{V})$, $\nu\circ\sigma(\alpha\wedge\neg\alpha)=\nu(\sigma(\alpha\wedge\neg\alpha))$ equals $\nu(\sigma(\alpha)\wedge\neg\sigma(\alpha))$; since $\nu$ lies in $\mathcal{F}_{C_{n}}^{\mathcal{B}}$, 
\[\nu(\sigma(\alpha)\wedge\neg\sigma(\alpha))_{2}=\nu(\sigma(\alpha))_{3}=\nu\circ\sigma(\alpha)_{3},\]
as we needed to prove.

\item We use, again, that $\sigma$ is a change of variables: from Lemma \ref{index formula}, $\sigma(\alpha^{1})=\sigma(\alpha)^{1}$; keeping in mind that $\nu\circ\sigma(\alpha)_{i}=\nu(\sigma(\alpha))_{i}$,
\[\nu\circ\sigma(\alpha^{1})=\nu(\sigma(\alpha)^{1})=\]
\[\Big(\nu(\sigma(\alpha))_{3}, \nu(\sigma(\alpha))_{1}\wedge\nu(\sigma(\alpha))_{2}, \nu(\sigma(\alpha))_{4}, \dotsc  , \nu(\sigma(\alpha))_{n+1}, {\sim}(\bigwedge_{i=1}^{n+1}\nu(\sigma(\alpha))_{i})\Big)=\]
\[\Big(\nu\circ\sigma(\alpha)_{3}, \nu\circ\sigma(\alpha)_{1}\wedge\nu\circ\sigma(\alpha)_{2}, \nu\circ\sigma(\alpha)_{4}, \dotsc  , \nu\circ\sigma(\alpha)_{n+1}, {\sim}(\bigwedge_{i=1}^{n+1}\nu\circ\sigma(\alpha)_{i})\Big).\]

\item Finally, for any $\#\in\{\vee, \wedge, \rightarrow\}$, 
\[\nu\circ\sigma\Big((\alpha^{(n)}\wedge\beta^{(n)})\rightarrow(\alpha\#\beta)^{(n)}\Big)=\nu\Big((\sigma(\alpha^{(n)})\wedge\nu(\beta^{(n)}))\rightarrow\sigma((\alpha\#\beta)^{(n)})\Big);\]
from Lemma \ref{index formula}, this equals 
\[\nu\Big((\sigma(\alpha)^{(n)}\wedge\sigma(\beta)^{(n)})\rightarrow\sigma(\alpha\#\beta)^{(n)}\Big)=\nu\Big((\sigma(\alpha)^{(n)}\wedge\sigma(\beta)^{(n)})\rightarrow(\sigma(\alpha)\#\sigma(\beta))^{(n)}\Big),\]
which is in $D_{n}^{\mathcal{B}}$ since $\nu$ is in $\mathcal{F}_{C_{n}}^{\mathcal{B}}$ (and $\sigma(\alpha)$ and $\sigma(\beta)$ are still formulas of $C_{n}$).
\end{enumerate}
\end{proof}

It follows, quite easily, that $\mathcal{RM}_{C_{n}}^{\mathcal{B}}$ is structural.

\begin{definition}
To signify that $\varphi$ follows from $\Gamma$ according to the semantical consequence relation with respect to $\mathcal{RM}_{C_{n}}^{\mathcal{B}}$ we write $\Gamma\vDash_{\mathcal{RM}_{C_{n}}^{\mathcal{B}}}^{RN}\varphi$\label{vDashRNRM}.

The semantical consequence relation with respect to the restricted swap structures semantics for $C_{n}$, i.e. the class $\mathcal{RS}_{C_{n}}$ of $\mathcal{RM}_{C_{n}}^{\mathcal{B}}$ for all Boolean algebras $\mathcal{B}$, will be similarly represented by $\Gamma\vDash_{\mathcal{RS}_{C_{n}}}^{RN}\varphi$: this means that\label{vDashRNRS} 
\[\Gamma\vDash_{\mathcal{RS}_{C_{n}}}^{RN}\varphi\quad\text{if, and only if,}\quad \Gamma\vDash_{\mathcal{RM}_{C_{n}}^{\mathcal{B}}}^{RN}\varphi\quad\text{for all Boolean algebras $\mathcal{B}$}.\]
\end{definition}

\begin{lemma}\label{Finding the coordinates}
For $2\leq k\leq n-2$, we have that, for a homomorphism $\nu\in\mathcal{F}_{C_{n}}^{\mathcal{B}}$,
\[\nu(\alpha^{k})=(\nu(\alpha)_{k+2}, \bigwedge_{i=1}^{k+1}\nu(\alpha)_{i}, \nu(\alpha)_{k+3}, \dotsc  ,  \nu(\alpha)_{n+1}, {\sim}\bigwedge_{i=1}^{n+1}\nu(\alpha)_{i}, 1, \dotsc  , 1);\]
of course, we then have 
\[\nu(\alpha^{n-1})=(\nu(\alpha)_{n+1}, \bigwedge_{i=1}^{n}\nu(\alpha)_{i}, {\sim}\bigwedge_{i=1}^{n+1}\nu(\alpha)_{i}, 1, \dotsc  , 1)\]
and 
\[\nu(\alpha^{n})=({\sim}\bigwedge_{i=1}^{n+1}\nu(\alpha)_{i}, \bigwedge_{i=1}^{n+1}\nu(\alpha)_{i}, 1, \dotsc  , 1).\]
\end{lemma}

\begin{proof}
For $k=2$, the result follows from the definition of $\mathcal{F}_{C_{n}}^{\mathcal{B}}$ and the fact that $\alpha^{2}=(\alpha^{1})^{1}$:
\begin{enumerate}
\item $\nu(\alpha^{2})_{1}=\nu(\alpha^{1})_{3}=\nu(\alpha)_{4}$; 
\item $\nu(\alpha^{2})_{2}=\nu(\alpha^{1})_{1}\wedge\nu(\alpha^{1})_{2}=\nu(\alpha)_{1}\wedge\nu(\alpha)_{2}\wedge\nu(\alpha)_{3}$;
\item for $3\leq j\leq n-1$, $\nu(\alpha^{2})_{j}=\nu(\alpha^{1})_{j+1}=\nu(\alpha)_{j+2}$;
\item $\nu(\alpha^{2})_{n}=\nu(\alpha^{1})_{n+1}={\sim}\bigwedge_{i=1}^{n+1}\nu(\alpha)_{i}$;
\item $\nu(\alpha^{2})_{n+1}={\sim}\bigwedge_{i=1}^{n+1}\nu(\alpha^{1})_{i}={\sim}(\bigwedge_{i=1}^{n+1}\nu(\alpha)_{i}\wedge{\sim}\bigwedge_{i=1}^{n+1}\nu(\alpha)_{i})=1$.
\end{enumerate}

For the general case, suppose the result holds for $k-1$, and remember $\alpha^{k}=(\alpha^{k-1})^{1}$.
\begin{enumerate}
\item $\nu(\alpha^{k})_{1}=\nu(\alpha^{k-1})_{3}=\nu(\alpha)_{k+3}$;
\item $\nu(\alpha^{k})_{2}=\nu(\alpha^{k-1})_{1}\wedge\nu(\alpha^{k-1})_{2}=\nu(\alpha)_{k+1}\wedge\bigwedge_{i=1}^{k}\nu(\alpha)_{i}=\bigwedge_{i=1}^{k+1}\nu(\alpha)_{i}$;
\item for $3\leq j\leq n-k+1$, $\nu(\alpha^{k})_{j}=\nu(\alpha^{k-1})_{j+1}=\nu(\alpha)_{j+k}$;
\item $\nu(\alpha^{k})_{n-k+2}=\nu(\alpha^{k-1})_{n-k+3}={\sim}\bigwedge_{i=1}^{k+1}\nu(\alpha)_{i}$;
\item since $\bigwedge_{i=1}^{n-k+2}=\bigwedge_{i=1}^{k+1}\nu(\alpha)_{i}\wedge{\sim}\bigwedge_{i=1}^{k+1}\nu(\alpha)_{i}=0$, all following coordinates must equal $1$.
\end{enumerate}

The cases $\nu(\alpha^{n-1})$ and $\nu(\alpha^{n})$ follow analogously.
\end{proof}

Using the previous lemma, specifically the fact that $\nu(\alpha^{k})_{1}=\nu(\alpha)_{k+2}$, is not difficult to prove that, for $1\leq k\leq n$, $\nu(\alpha^{(k)})_{1}=\bigwedge_{i=3}^{k+2}\nu(\alpha)_{i}$: the base case is rather obvious, since $\nu(\alpha^{(1)})_{1}=\nu(\alpha^{1})_{1}=\nu(\alpha)_{3}$; assuming the result holds for an arbitrary $1\leq k\leq n-1$, 
\[\nu(\alpha^{(k+1)})_{1}=\nu(\alpha^{(k)})_{1}\wedge\nu(\alpha^{k+1})_{1}=\bigwedge_{i=3}^{k+2}\nu(\alpha)_{i}\wedge\nu(\alpha)_{k+3}=\bigwedge_{i=3}^{k+3}\nu(\alpha)_{i}.\]

Then, one has, again from Lemma \ref{Finding the coordinates},
\[\nu(\alpha^{(n)})_{1}=\nu(\alpha^{(n-1)})_{1}\wedge\nu(\alpha^{n})_{1}=\bigwedge_{i=3}^{n+1}\nu(\alpha)_{i}\wedge{\sim}\bigwedge_{i=1}^{n+1}\nu(\alpha)_{i}=\]
\[\bigwedge_{i=3}^{n+1}\nu(\alpha)_{i}\wedge{\sim}\Big(\nu(\alpha)_{1}\wedge\nu(\alpha)_{2}\wedge\bigwedge_{i=3}^{n+1}\nu(\alpha)_{i}\Big)=\bigwedge_{i=3}^{n+1}\nu(\alpha)_{i}\wedge\Big({\sim}(\nu(\alpha)_{1}\wedge\nu(\alpha)_{2})\vee{\sim}\bigwedge_{i=3}^{n+1}\nu(\alpha)_{i}\Big)=\]
\[\Big(\bigwedge_{i=3}^{n+1}\nu(\alpha)_{i}\wedge{\sim}(\nu(\alpha)_{1}\wedge\nu(\alpha)_{2})\Big)\vee \Big(\bigwedge_{i=3}^{n+1}\nu(\alpha)_{i}\wedge{\sim}\bigwedge_{i=3}^{n+1}\nu(\alpha)_{i}\Big)=\]
\[{\sim}(\nu(\alpha)_{1}\wedge\nu(\alpha)_{2})\wedge\bigwedge_{i=3}^{n+1}\nu(\alpha)_{i}={\sim}(\nu(\alpha)_{1}\wedge\nu(\alpha)_{2})\wedge\nu(\alpha^{(n-1)})_{1}\]

\begin{proposition}\label{Value of strong negation}
Let ${\sim}\alpha=\neg\alpha\wedge\alpha^{(n)}$ be the definable strong negation in $C_{n}$; for every $\nu\in\mathcal{F}_{C_{n}}^{\mathcal{B}}$, one has:
\begin{enumerate}
\item $\nu(\alpha\wedge{\sim}\alpha)=F_{n}$;
\item $\nu(\alpha\vee{\sim\alpha})\in D_{n}^{\mathcal{B}}$.
\end{enumerate}
\end{proposition}

\begin{proof}

\begin{enumerate}
\item One has $\nu(\neg\alpha)_{1}=\nu(\alpha)_{2}$ and $\nu(\alpha^{(n)})_{1}={\sim}(\nu(\alpha)_{1}\wedge\nu(\alpha)_{2})\wedge\nu(\alpha^{(n-1)})_{1}$, and therefore
\[\nu(\neg\alpha\wedge\alpha^{(n)})_{1}=\nu(\alpha)_{2}\wedge{\sim}(\nu(\alpha)_{1}\wedge\nu(\alpha)_{2})\wedge\nu(\alpha^{(n-1)})_{1}=\nu(\alpha)_{2}\wedge({\sim}\nu(\alpha)_{1}\vee{\sim}\nu(\alpha)_{2})\wedge\nu(\alpha^{(n-1)})_{1}=\]
\[\Big(\nu(\alpha)_{2}\wedge{\sim}\nu(\alpha)_{1}\wedge\nu(\alpha^{(n-1)})_{1}\Big)\vee\Big(\nu(\alpha)_{2}\wedge{\sim}\nu(\alpha)_{2}\wedge\nu(\alpha^{(n-1)})_{1}\Big)=\nu(\alpha)_{2}\wedge{\sim}\nu(\alpha)_{1}\wedge\nu(\alpha^{(n-1)})_{1}.\]
So, 
\[\nu(\alpha\wedge{\sim}\alpha)_{1}=\nu(\alpha)_{1}\wedge\Big(\nu(\alpha)_{2}\wedge{\sim}\nu(\alpha)_{1}\wedge\nu(\alpha^{(n-1)})_{1}\Big)=0,\]
what means that $\nu(\alpha\wedge{\sim}\alpha)=F_{n}$.

\item Since $\nu(\alpha)$ is in $B_{n}^{\mathcal{B}}$, $\nu(\alpha)_{1}\vee\nu(\alpha)_{2}=1$ and, for any $1\leq k\leq n$, $(\bigwedge_{i=1}^{k}\nu(\alpha)_{i})\vee\nu(\alpha)_{k+1}=1$, meaning 
\[{\sim}\nu(\alpha)_{k+1}\leq \bigwedge_{i=1}^{k}\nu(\alpha)_{i}\leq \nu(\alpha)_{1}.\]
This implies that, for all $3\leq k\leq n+1$, ${\sim}\nu(\alpha)_{k}\leq \nu(\alpha)_{1}$, and so ${\sim}\bigwedge_{i=3}^{n+1}\nu(\alpha)_{i}=\bigvee_{i=3}^{n+1}{\sim}\nu(\alpha)_{i}\leq \nu(\alpha)_{1}$, what implies in turn that
\[\nu(\alpha\vee{\sim}\alpha)_{1}=\nu(\alpha)_{1}\vee\Big(\nu(\alpha)_{2}\wedge{\sim}\nu(\alpha)_{1}\wedge\nu(\alpha^{(n-1)})_{1}\Big)=\]
\[(\nu(\alpha)_{1}\vee\nu(\alpha)_{2})\wedge(\nu(\alpha)_{1}\vee{\sim}\nu(\alpha)_{1})\wedge(\nu(\alpha)_{1}\vee\bigwedge_{i=3}^{n+1}\nu(\alpha)_{i})=1.\]
Of course, this finishes proving that $\nu(\alpha\vee{\sim}\alpha)\in D_{n}^{\mathcal{B}}$.

\end{enumerate}
\end{proof}

\begin{lemma}\label{sound arb}
Let $\nu$ be a valuation in $\mathcal{F}_{C_{n}}^{\mathcal{B}}$; then, the mapping $\mathsf{b}:\textbf{F}(\Sigma_{\textbf{C}}, \mathcal{V})\rightarrow |\mathcal{B}|$\label{|B|}, given by $\mathsf{b}(\alpha)=\nu(\alpha)_{1}$, is a $\mathcal{B}$-valuation for $C_{n}$ such that $\mathsf{b}(\alpha)=1$ if, and only if, $\nu(\alpha)\in D_{n}^{\mathcal{B}}$ for every formula $\alpha$.
\end{lemma}

\begin{proof}
\begin{enumerate}
\item Clause $(V1)$ is quite obvious: since $\nu$ is a homomorphism, $\nu(\alpha\#\beta)\in \nu(\alpha)\tilde{\#}\nu(\beta)$, for any $\#\in\{\vee, \wedge, \rightarrow\}$, and given that $u\in z\tilde{\#}w$ if, and only if, $u_{1}=z_{1}\# w_{1}$, we obtain that $\nu(\alpha\#\beta)_{1}=\nu(\alpha)_{1}\#\nu(\beta)_{1}$ and therefore $\mathsf{b}(\alpha\#\beta)=\mathsf{b}(\alpha)\#\mathsf{b}(\beta)$.

\item We have $\mathsf{b}(\neg\alpha)=\nu(\neg\alpha)_{1}$: since $\nu(\neg\alpha)\in \tilde{\neg}\nu(\alpha)$, given $\nu$ is a homomorphism, and $w\in \tilde{\neg}z$ if, and only if $w_{1}=z_{2}$ and $w_{2}\leq z_{1}$, we obtain $\nu(\neg\alpha)_{1}=\nu(\alpha)_{2}$; from the definition of $B_{n}^{\mathcal{B}}$, $\nu(\alpha)_{1}\vee\nu(\alpha)_{2}=1$, implying $\nu(\alpha)_{2}\geq {\sim}\nu(\alpha)_{1}$. Of course, this means $\mathsf{b}(\neg\alpha)\geq{\sim}\mathsf{b}(\alpha)$, which corresponds to clause $(V2)$.

\item Again from the definition of $\tilde{\neg}$, $\nu(\neg\neg\alpha)_{1}=\nu(\neg\alpha)_{2}$, and $\nu(\neg\alpha)_{2}\leq \nu(\alpha)_{1}$, meaning $\mathsf{b}(\neg\neg\alpha)=\nu(\neg\neg\alpha)_{1}\leq \nu(\alpha)_{1}=\mathsf{b}(\alpha)$. This corresponds to clause $(V3)$.

\item From Lemma \ref{Finding the coordinates}, we have $\nu(\alpha^{n})_{1}={\sim}\bigwedge_{i=1}^{n+1}\nu(\alpha)_{i}$, $\nu(\alpha^{n-1})_{1}=\nu(\alpha)_{n+1}$ and $\nu(\neg(\alpha^{n-1}))_{1}=\nu(\alpha^{n-1})_{2}=\bigwedge_{i=1}^{n}\nu(\alpha)_{i}$, meaning therefore that 
\[\mathsf{b}(\alpha^{n})={\sim}\bigwedge_{i=1}^{n+1}\nu(\alpha)_{i}={\sim}(\nu(\alpha)_{n+1}\wedge \bigwedge_{i=1}^{n}\nu(\alpha)_{i})={\sim}(\mathsf{b}(\alpha^{n-1})\wedge\mathsf{b}(\neg(\alpha^{n-1}))),\]
that is, clause $(V4)_{n}$.

\item We have $\mathsf{b}(\neg(\alpha^{1}))=\nu(\neg(\alpha^{1}))_{1}=\nu(\alpha^{1})_{2}=\nu(\alpha)_{1}\wedge\nu(\alpha)_{2}=\nu(\alpha)_{1}\wedge\nu(\neg\alpha)_{1}=\mathsf{b}(\alpha)\wedge\mathsf{b}(\neg\alpha)$. That validates clause $(V5)$.

\item From the third condition of Definition \ref{restricted valuations for RMBCn}, for any $\#\in\{\vee, \wedge, \rightarrow\}$ one has $\nu((\alpha^{(n)}\wedge\beta^{(n)})\rightarrow(\alpha\#\beta)^{(n)})\in D_{n}^{\mathcal{B}}$, meaning $(\nu(\alpha^{(n)})_{1}\wedge\nu(\beta^{(n)})_{1})\rightarrow\nu((\alpha\#\beta)^{(n)})_{1}=\nu((\alpha^{(n)}\wedge\beta^{(n)})\rightarrow(\alpha\#\beta)^{(n)})_{1}=1$. Of course, this implies $\mathsf{b}(\alpha^{(n)})\wedge\mathsf{b}(\beta^{(n)})=\nu(\alpha^{(n)})_{1}\wedge\nu(\beta^{(n)})_{1}\leq \nu((\alpha\#\beta)^{(n)})_{1}=\mathsf{b}((\alpha\#\beta)^{(n)})$. This shows condition $(V6)_{n}$ is also validated.
\end{enumerate}

Of course, $\mathsf{b}(\alpha)=1$ if, and only if, $\nu(\alpha)_{1}=1$, which is equivalent to $\nu(\alpha)\in D_{n}^{\mathcal{B}}$.
\end{proof}

\begin{lemma}\label{comp arb}
Let $\mathsf{b}$ be a $\mathcal{B}$-valuation for $C_{n}$; then, the mapping $\nu:\textbf{F}(\Sigma_{\textbf{C}}, \mathcal{V})\rightarrow B_{n}^{\mathcal{B}}$ given by $\nu(\alpha)=(\mathsf{b}(\alpha), \mathsf{b}(\neg\alpha), \mathsf{b}(\alpha^{1}), \dotsc  , \mathsf{b}(\alpha^{n-1}))$ is a valuation in $\mathcal{F}_{C_{n}}^{\mathcal{B}}$ such that $\mathsf{b}(\alpha)=1$ if, and only if, $\nu(\alpha)\in D_{n}^{\mathcal{B}}$ for every formula $\alpha$.
\end{lemma}

\begin{proof}
First of all, we prove $\nu$ is a homomorphism.
\begin{enumerate}
\item One sees that, by definition of $\nu$, $\nu(\neg\alpha)_{1}=\mathsf{b}(\neg\alpha)=\nu(\alpha)_{2}$ and, by $(V3)$, $\nu(\neg\alpha)_{2}=\mathsf{b}(\neg\neg\alpha)\leq \mathsf{b}(\alpha)=\nu(\alpha)_{1}$, proving that $\nu(\neg\alpha)\in \tilde{\#}\nu(\alpha)$.

\item For $\#\in\{\vee, \wedge, \rightarrow\}$, from condition $(V1)$ one gets $\nu(\alpha\#\beta)_{1}=\mathsf{b}(\alpha\#\beta)=\mathsf{b}(\alpha)\#\mathsf{b}(\beta)=\nu(\alpha)_{1}\#\nu(\beta)_{1}$. Furthermore, $\nu(\alpha), \nu(\beta)\in Boo_{n}^{\mathcal{B}}$ if and only if $\nu(\alpha)_{1}={\sim}\nu(\alpha)_{2}$ and $\nu(\beta)_{1}={\sim}\nu(\beta)_{2}$, or equivalently, $\mathsf{b}(\neg\alpha)={\sim}\mathsf{b}(\alpha)$ and $\mathsf{b}(\neg\beta)={\sim}\mathsf{b}(\beta)$; from Proposition \ref{Negations of B-valuations}, this implies $\mathsf{b}(\neg(\alpha\#\beta))={\sim}\mathsf{b}(\alpha\#\beta)$, that is, $\nu(\alpha\#\beta)\in Boo_{n}^{\mathcal{B}}$. With all of this, we find that, regardless of the values of $\nu(\alpha)$ and $\nu(\beta)$, $\nu(\alpha\#\beta)\in \nu(\alpha)\tilde{\#}\nu(\beta)$.
\end{enumerate}

Now, we need only to prove that $\nu$ is in $\mathcal{F}_{C_{n}}^{\mathcal{B}}$.

\begin{enumerate}
\item From the fact $\nu$ is a homomorphism, $\nu(\alpha\wedge\neg\alpha)\in \nu(\alpha)\tilde{\wedge}\nu(\neg\alpha)$; moreover, $\nu(\alpha\wedge\neg\alpha)_{2}=\mathsf{b}(\neg(\alpha\wedge\neg\alpha))=\mathsf{b}(\alpha^{1})=\nu(\alpha)_{3}$, what proves the first condition for being in $\mathcal{F}_{C_{n}}^{\mathcal{B}}$ of Definition \ref{restricted valuations for RMBCn} is satisfied.

\item $\nu(\alpha^{1})_{1}=\mathsf{b}(\alpha^{1})=\nu(\alpha)_{3}$, from the definition of $\nu$; from property $(V5)$, $\nu(\alpha^{1})_{2}=\mathsf{b}(\neg(\alpha^{1}))=\mathsf{b}(\alpha)\wedge\mathsf{b}(\neg\alpha)=\nu(\alpha)_{1}\wedge\nu(\alpha)_{2}$; for $3\leq k\leq n$, $\nu(\alpha^{1})_{k}=\mathsf{b}((\alpha^{1})^{k-2})=\mathsf{b}(\alpha^{k-1})=\nu(\alpha)_{k+1}$.

Finally, we have from $(V4)_{n}$ that $\nu(\alpha^{1})_{n+1}=\mathsf{b}(\alpha^{n})={\sim}(\mathsf{b}(\alpha^{n-1})\wedge\mathsf{b}(\neg(\alpha^{n-1})))$; from $(V5)$, $\mathsf{b}(\neg(\alpha^{n-1}))=\mathsf{b}(\alpha^{n-2})\wedge\mathsf{b}(\neg(\alpha^{n-2}))$, and proceeding recursively, one obtains $\nu(\alpha^{1})_{n+1}=(\mathsf{b}(\alpha)\wedge\mathsf{b}(\neg\alpha)\wedge{\sim}\bigwedge_{i=1}^{n-1}\mathsf{b}(\alpha^{i}))={\sim}\bigwedge_{i=1}^{n+1}\nu(\alpha)_{i}$.

\item For any $\#\in\{\vee, \wedge, \rightarrow\}$, from $(V6)_{n}$ we find that $\mathsf{b}(\alpha^{(n)})\wedge\mathsf{b}(\beta^{(n)})\leq\mathsf{b}((\alpha\#\beta)^{(n)})$, meaning $\nu(\alpha^{(n)}\wedge\beta^{(n)})_{1}=\nu(\alpha^{(n)})_{1}\wedge\nu(\beta^{(n)})_{1}\leq \nu((\alpha\#\beta)^{(n)})_{1}$, and therefore $\nu(\alpha^{(n)}\wedge\beta^{(n)})_{1}\rightarrow\nu((\alpha\#\beta)^{(n)})_{1}=\nu((\alpha^{(n)}\wedge\beta^{(n)})\rightarrow(\alpha\#\beta)^{(n)})_{1}=1$, which is equivalent to $\nu((\alpha^{(n)}\wedge\beta^{(n)})\rightarrow(\alpha\#\beta)^{(n)})\in D_{n}^{\mathcal{B}}$.
\end{enumerate}

Clearly, $\mathsf{b}(\alpha)=1$ if, and only if, $\nu(\alpha)_{1}=1$, which is in turn equivalent to $\nu(\alpha)\in D_{n}^{\mathcal{B}}$.
\end{proof}

\begin{theorem}
Given formulas $\Gamma\cup\{\varphi\}$ of $C_{n}$, $\Gamma\vdash_{C_{n}}\varphi$ if, and only if, $\Gamma\vDash_{\mathcal{RS}_{C_{n}}}^{RN}\varphi$.
\end{theorem}

\begin{proof}
We start by assuming that $\Gamma\vdash_{C_{n}}\varphi$, and taking a valuation $\nu$ lying in $\mathcal{F}_{C_{n}}^{\mathcal{B}}$ and satisfying $\nu(\Gamma)\subseteq D_{n}^{\mathcal{B}}$: from Lemma \ref{sound arb}, the function 
\[\mathsf{b}:F(\Sigma_{\textbf{C}}, \mathcal{V})\rightarrow\mathcal{B},\]
such that $\mathsf{b}(\alpha)=\nu(\alpha)_{1}$ for any formula $\alpha$ of $C_{n}$, is a $\mathcal{B}$-valuation for which, by hypothesis, $\mathsf{b}(\Gamma)\subseteq\{1\}$. We know the method of $\mathcal{B}$-valuations to be sound for $C_{n}$, so $\Gamma\vdash_{C_{n}}\varphi$ implies that $\mathsf{b}(\varphi)=1$, and therefore $\nu(\varphi)\in D_{n}^{\mathcal{B}}$. This proves $\Gamma\vDash_{\mathcal{RM}_{C_{n}}}\varphi$, and hence the semantics of restricted swap structures is sound.

Reciprocally, assume $\Gamma\vDash_{\mathcal{RS}_{C_{n}}}^{RN}\varphi$ and suppose $\mathsf{b}$ is a $\mathcal{B}$-valuation with $\mathsf{b}(\Gamma)\subseteq\{1\}$. From Lemma \ref{comp arb}, 
\[\nu:\textbf{F}(\Sigma_{\textbf{C}}, \mathcal{V})\rightarrow \mathcal{A}_{C_{n}}^{\mathcal{B}}\] 
defined as
\[\nu(\alpha)=(\mathsf{b}(\alpha), \mathsf{b}(\neg\alpha), \mathsf{b}(\alpha^{1}), \dotsc  , \mathsf{b}(\alpha^{n-1})),\]
is a valuation in $\mathcal{F}_{C_{n}}^{\mathcal{B}}$ with the additional property that $\nu(\Gamma)\subseteq D_{n}^{\mathcal{B}}$. Since $\Gamma\vDash_{\mathcal{RS}_{C_{n}}}^{RN}\varphi$, it follows that $\nu(\varphi)\in D_{n}^{\mathcal{B}}$, what means $\mathsf{b}(\varphi)=1$ and, by completeness of $C_{n}$ with respect to $\mathcal{B}$-valuations, $\Gamma\vdash_{C_{n}}\varphi$.
\end{proof}

\subsection{$\mathcal{RM}_{C_{n}}^{\mathcal{B}}$ is not trivial}\label{RM is not trivial}

A natural question is then if $\mathcal{RM}_{C_{n}}^{\mathcal{B}}$ adds anything to the already defined $\mathcal{RM}_{C_{n}}$: that is, is there a valuation for $\mathcal{RM}_{C_{n}}^{\mathcal{B}}$ which is not a valuation for $\mathcal{RM}_{C_{n}}$? This is mostly equivalent to the question of whether there exists a (non-trivial) Boolean algebra $\mathcal{B}$ and a $\mathcal{B}$-valuation which is not a bivaluation.

So, take a non-trivial Boolean algebra $\mathcal{B}$, and remember Definition \ref{B-valuation} of a $\mathcal{B}$-valuation . We will construct by structural induction a $\mathcal{B}$-valuation $\mathsf{b}$ for $C_{n}$, extending any given function $b$ from $\mathcal{V}$ into the universe of $\mathcal{B}$, and by taking a $b$ for which there exists a $p\in\mathcal{V}$ with $b(\neg p)\neq {\sim}b(p)$ we will obtain that $\mathsf{b}$ satisfies:
\begin{enumerate}
\item that there exists a formula $\alpha$ for $C_{n}$ with $\mathsf{b}(\neg\alpha)\neq{\sim}\mathsf{b}(\alpha)$;
\item that the image $\{\mathsf{b}(\alpha) : \alpha\in F(\Sigma_{\textbf{C}}, \mathcal{V})\}$ of $\mathsf{b}$ is not contained in $\{0,1\}$.
\end{enumerate}

Furthermore, by taking the homomorphism $\nu(\alpha)=(\mathsf{b}(\alpha), \mathsf{b}(\neg\alpha), \mathsf{b}(\alpha^{1}), \dotsc  , \mathsf{b}(\alpha^{n-1}))$, we find a valuation in $\mathcal{F}_{C_{n}}^{\mathcal{B}}$ whose image is not contained in $B_{n}$ or $Boo_{n}^{\mathcal{B}}$.

Now, concerning the task of defining $\mathsf{b}$, satisfying clauses $(V1)$ through $(V5)$ of Definition \ref{B-valuation} is easy, but in order to guarantee that clause $(V6)_{n}$ is respected we must be very careful. We have that $\mathsf{b}$ must be equal to $b$ over formulas of complexity $0$, namely variables; we must also define, for our recursion to work, $\mathsf{b}(\neg p)$, $\mathsf{b}(p^{i})$ and $\mathsf{b}(\neg p^{i})$, for all $i\in\mathbb{N}$ and $p\in\mathcal{V}$, but we will show how this is done in the more general setting where $p$ is replaced by an arbitrary formula $\alpha$ of bounded complexity. So, assume for a moment that $\mathsf{b}(\alpha)$, $\mathsf{b}(\neg\alpha)$, $\mathsf{b}(\alpha^{i})$ and $\mathsf{b}(\neg\alpha^{i})$, for any $i\in\mathbb{N}$, have been defined for all formulas $\alpha$ of complexity at most $m$: then a formula of complexity $m+1$ is either: $\alpha\#\beta$, for $\alpha$ and $\beta$ of complexity less or equal to $m$ and $\#\in\{\vee, \wedge, \rightarrow\}$, when we define, in order to satisfy $(V1)$,
\[\mathsf{b}(\alpha\#\beta)=\mathsf{b}(\alpha)\#\mathsf{b}(\beta);\]
or $\neg\alpha$, for an $\alpha$ of complexity equal to $m$, when $\mathsf{b}(\neg\alpha)$ has already been defined. So we must now specify, for $\alpha$ and $\beta$ of complexity at most $m$, $\mathsf{b}(\neg\neg \alpha)$, $\mathsf{b}((\neg\alpha)^{i})$, $\mathsf{b}(\neg(\neg\alpha)^{i})$, $\mathsf{b}(\neg(\alpha\#\beta))$, $\mathsf{b}((\alpha\#\beta)^{j})$ and $\mathsf{b}(\neg(\alpha\#\beta)^{j}))$, for $1\leq i,j\leq n$.\footnote{We do not need to worry about defining $\mathsf{b}(\alpha^{j})$ or $\mathsf{b}(\neg\alpha^{j})$ for $j> n$, for any formula $\alpha$, since they must always take the values, respectively, $1$ and $0$.} We start by looking at $\mathsf{b}(\neg\neg\alpha)$, which we make equal to
\begin{enumerate}
\item $\mathsf{b}(\gamma)\wedge\mathsf{b}(\neg\gamma)$, if $\neg\alpha=\gamma^{1}$ (thus validating $(V5)$);
\item any value between $\mathsf{b}(\alpha)$ and ${\sim}\mathsf{b}(\neg\alpha)$, both included, if the previous case does not apply (validating $(V2)$ and $(V3)$).
\end{enumerate}
Regarding $\mathsf{b}((\neg\alpha)^{i})$ and $\mathsf{b}(\neg(\neg\alpha)^{i})$, for $1\leq i\leq n-1$, there isn't a lot we must do: supposing we have defined $\mathsf{b}((\neg\alpha)^{i-1})$ and $\mathsf{b}(\neg(\neg\alpha)^{i-1})$, it is sufficient to demand that $\mathsf{b}((\neg\alpha)^{i})$ is greater or equal to 
\[{\sim}[\mathsf{b}((\neg\alpha)^{i-1}\wedge\neg(\neg\alpha)^{i-1})],\]
in order to satisfy clause $(V2)$, and that $\mathsf{b}(\neg(\neg\alpha)^{i})$ equals $\mathsf{b}((\neg\alpha)^{i-1})\wedge\mathsf{b}(\neg(\neg\alpha)^{i-1})$, in order to validate $(V5)$ (unless $\neg\alpha=\gamma^{k}$ for some $1\leq k\leq n-1$ such that $k+i=n$, in which case $\mathsf{b}((\neg\alpha)^{i})=\mathsf{b}(\gamma^{n})$ is defined according to $(V4)_{n}$). Of course, $\mathsf{b}((\neg\alpha)^{n})$ and $\mathsf{b}(\neg(\neg\alpha)^{n})$ are given values according to clauses, respectively, $(V4)_{n}$ and $(V5)$, being expressed in terms of the already defined $\mathsf{b}((\neg\alpha)^{n-1})$ and $\mathsf{b}(\neg(\neg\alpha)^{n-1})$.

Difficulties finally arise when defining $\mathsf{b}(\neg(\alpha\#\beta))$, $\mathsf{b}((\alpha\#\beta)^{i})$ and $\mathsf{b}(\neg(\alpha\#\beta)^{i})$, because we must be mindful of clause $(V6)_{n}$: denote by $a_{0}$ the value $\mathsf{b}(\alpha_{0})=\mathsf{b}(\alpha)\#\mathsf{b}(\beta)$, and by $a_{\#}$ the value $\mathsf{b}(\alpha^{(n)})\wedge\mathsf{b}(\beta^{(n)})=\bigwedge_{j=1}^{n}\mathsf{b}(\alpha^{j})\wedge\mathsf{b}(\beta^{j})$ (all of $\mathsf{b}(\alpha^{j})$ and $\mathsf{b}(\beta^{j})$, by hypothesis, already defined).

\begin{enumerate}[align=left]
\item[\underline{$\mathsf{b}(\neg(\alpha\#\beta))$} - ] The value for $\mathsf{b}(\neg(\alpha\#\beta))$ will be $a_{\neg}$: on the off chance that $\#=\wedge$, $\alpha=\gamma^{n-1}$ and $\beta=\neg\gamma^{n-1}$, we make $a_{\neg}$ equal to ${\sim}a_{0}={\sim}[\mathsf{b}(\gamma^{n-1})\wedge\mathsf{b}(\neg(\gamma^{n-1}))]$ (as required by $(V4)_{n}$); otherwise, $a_{\neg}$ can be any value such that 
\[a_{\neg}\geq {\sim}a_{0}\quad\text{and}\quad a_{\neg}\wedge a_{0}\leq{\sim}a_{\#}\]
(this is possible, one example of such an element being $a_{\neg}={\sim}a_{0}$); this last case of course implies $\mathsf{b}(\neg(\alpha\#\beta))=a_{\neg}\geq{\sim}a_{0}={\sim}\mathsf{b}(\alpha\#\beta)$, and thus also validates $(V2)$. 

\item[\underline{$\mathsf{b}((\alpha\#\beta)^{1})$} - ] Having defined $a_{\neg}$, we define $\mathsf{b}((\alpha\#\beta)^{1})$ to be any value $a_{1}$ such that $a_{1}\geq {\sim}(a_{0}\wedge a_{\neg})$ (and therefore we have
\[\mathsf{b}\big(\neg((\alpha\#\beta)\wedge\neg(\alpha\#\beta))\big)=\mathsf{b}((\alpha\#\beta)^{1})=a_{1}\geq {\sim}(a_{0}\wedge a_{\neg})={\sim}\mathsf{b}\big((\alpha\#\beta)\wedge\neg(\alpha\#\beta)\big),\]
and thus $(V2)$ remains valid) and ${\sim}(a_{0}\wedge a_{\neg})\wedge a_{1}\geq a_{\#}$ (what would still be possible, being enough to choose, for one, $a_{1}={\sim}(a_{0}\wedge a_{\neg})$).

And we make $\mathsf{b}(\neg(\alpha\#\beta)^{1})=\mathsf{b}(\alpha\#\beta)\wedge\mathsf{b}(\neg(\alpha\#\beta))=a_{0}\wedge a_{\neg}$, thus respecting clause $(V5)$ (and incidentally also $(V2)$ and $(V3)$).

\item[\underline{Inductive step} - ] Supposing we have defined $a_{1}=\mathsf{b}((\alpha\#\beta)^{1})$ trough $a_{k}=\mathsf{b}((\alpha\#\beta)^{k})$, for $1\leq k\leq n-2$, we make $\mathsf{b}((\alpha\#\beta)^{k+1})$ equal to any value $a_{k+1}$ satisfying 
\[a_{k+1}\geq{\sim}(a_{0}\wedge a_{\neg}\wedge \bigwedge_{j=1}^{k}a_{j})\quad\text{and}\quad{\sim}(a_{0}\wedge a_{\neg})\wedge\bigwedge_{j=2}^{k}a_{j}\geq a_{\#};\]

Inductively, we get $\mathsf{b}(\neg(\alpha\#\beta)^{k})=a_{0}\wedge a_{\neg}\wedge \bigwedge_{j=1}^{k-1}a_{j}$, and so: first of all, this means that $\mathsf{b}((\alpha\#\beta)^{k})\wedge\mathsf{b}(\neg(\alpha\#\beta)^{k})=a_{0}\wedge a_{\neg}\wedge \bigwedge_{j=1}^{k}a_{j}$ and therefore 
\[\mathsf{b}((\alpha\#\beta)^{k+1})=\mathsf{b}(\neg((\alpha\#\beta)^{k}\wedge\neg(\alpha\#\beta)^{k})\geq {\sim}\big(\mathsf{b}((\alpha\#\beta)^{k})\wedge\mathsf{b}(\neg(\alpha\#\beta)^{k})\big),\]
meaning $(V2)$ is respected; second, by making $\mathsf{b}(\neg(\alpha\#\beta)^{k+1})=\mathsf{b}((\alpha\#\beta)^{k})\wedge \mathsf{b}(\neg(\alpha\#\beta)^{k})$, so as to respect $(V5)$, we obtain that $\mathsf{b}(\neg(\alpha\#\beta)^{k+1})=a_{0}\wedge a_{\neg}\wedge \bigwedge_{j=1}^{k}a_{j}$ and the pattern remains.

\item[\underline{$\mathsf{b}((\alpha\#\beta)^{n})$} - ] Completing the process outlined in the previous items, we make $\mathsf{b}((\alpha\#\beta)^{n})$ equal to 
\[a_{n}={\sim}(\mathsf{b}((\alpha\#\beta)^{n-1})\wedge \mathsf{b}(\neg(\alpha\#\beta)^{n-1}))={\sim}(a_{0}\wedge a_{\neg}\wedge \bigwedge_{j=1}^{n-1}a_{j}),\]
given clause $(V4)_{n}$; and define $\mathsf{b}(\neg(\alpha\#\beta)^{n})$ to be simply $\mathsf{b}((\alpha\#\beta)^{n-1})\wedge\mathsf{b}(\neg(\alpha\#\beta)^{n-1})$, that is, $a_{0}\wedge a_{\neg}\wedge\bigwedge_{j=1}^{n-1}a_{j}$.
\end{enumerate}

Finally, now that $a_{1}$ trough $a_{n}$ are defined, given that $(\alpha\#\beta)^{(n)}=\bigwedge_{j=1}^{n}(\alpha\#\beta)^{j}$, one obtains 
\[\mathsf{b}((\alpha\#\beta)^{(n)})=\bigwedge_{j=1}^{n}\mathsf{b}((\alpha\#\beta)^{j})=a_{n}\wedge\bigwedge_{j=1}^{n-1}a_{j};\]
so
\[\mathsf{b}((\alpha\#\beta)^{(n)})=[{\sim}(a_{0}\wedge a_{\neg})\vee{\sim}\bigwedge_{j=1}^{n-1}a_{j})]\wedge\bigwedge_{j=1}^{n-1}a_{j}={\sim}(a_{0}\wedge a_{\neg})\wedge\bigwedge_{j=1}^{n-1}a_{j}\geq a_{\#}=\mathsf{b}(\alpha^{(n)})\wedge\mathsf{b}(\beta^{(n)}),\]
satisfying $(V6)_{n}$ and proving, with some difficulty, $\mathsf{b}$ is a $\mathcal{B}$-valuation, as we wanted to show.

\section{Counting snapshots}\label{Counting snapshots}

We would like to prove, for completeness sake, that if $\mathcal{B}$ is a finite Boolean algebra with $2^{m}$ elements, then $B_{n}^{\mathcal{B}}$ has $(n+2)^{m}$ elements. To this end, we begin with a few observations intended to make our job easier. We notice that the very useful equality
\[B_{n+1}^{\mathcal{B}}=\{(a_{1}, \dotsc  , a_{n+2})\in |\mathcal{B}|^{n+2}\ :\  (a_{1}, \dotsc  , a_{n+1})\in B_{n}^{\mathcal{B}}\quad\text{and}\quad a_{n+2}\vee\bigwedge_{i=1}^{n+1}a_{i}=1\}\]
holds for any $n\in\mathbb{N}\setminus\{0\}$.

Now, given all finite Boolean algebras are isomorphic to the field of subsets of some set, we will simply assume all finite Boolean algebras here involved are, in fact, fields of subsets: in that case, the Boolean algebra with $2^{m}$ elements is $\mathcal{P}(X_{m})$, for $X_{m}=\{x_{1}, \dotsc  , x_{m}\}$\label{Xm} a canonical set with $m$ elements. We will say that an element $a$ of a finite Boolean algebra $\mathcal{P}(X_{m})$ is of order $k$ if it is a subset of $X_{m}$ with $k$ elements; it is clear that there are $\binom{m}{0}=1$ elements of order $0$ (namely, $\emptyset$, the bottom), $\binom{m}{1}$ elements of order $1$, \dots , and $\binom{m}{m}=1$ elements of order $m$ (namely, $X_{m}$, the top). For what will follow, it is important to remember that the binomial coefficient is given by
\[\binom{m}{n}=\frac{m!}{n!(m-n)!},\]
for $m, n\in\mathbb{N}$ and $n\leq m$, and that the binomial theorem states 
\[(x+y)^{m}=\sum_{k=0}^{m}\binom{m}{k}x^{k}y^{m-k},\]
for $m\in\mathbb{N}$, and $x$ and $y$ any elements of a commutative ring with unity.

\begin{lemma}\label{finding solutions of a given order}
For an element $a$ of $\mathcal{P}(X_{m})$ of order $k$, there are $\binom{k}{p}2^{m-k}$ elements $b$ such that $a\wedge b$ has order $p\leq k$, and $\binom{m-k}{q}2^{k}$ elements $c$ such that $a\vee c$ has order $q\geq k$.
\end{lemma}

\begin{proof}
If $a\wedge b$ has order $p$, this means $b$ has $p$ elements in common with $a$, and since there are $k$ elements in $a$, we obtain $\binom{k}{p}$ possible values for $a\cap b$; but $a^{c}\cap b$ can be any subset of $a^{c}$, of which there exist $2^{m-k}$ given $a^{c}$ has $m-k$ elements. Combining the two values, we obtain $\binom{k}{p}2^{m-k}$ solutions.

Now, if $a\vee c$ has order $q$, $a^{c}\cap c$ has $q-k$ elements that $a$ doesn't, and since there are $m-k$ elements in $a^{c}$, this gives us $\binom{m-k}{q}$ values for $a^{c}\cap c$; but $a\cap c$ can be any subset of $a$, giving us $2^{k}$ values for $a\cap c$ since $a$ has $k$ elements. This leads to the total $\binom{m-k}{q}2^{k}$ solutions.
\end{proof}

Notice that, by adapting the proof of the previous lemma, we may show that, for an $a$ of order $k$, there are $\binom{k}{p}$ elements $b$ such that $a\vee b=1$ and $a\wedge b$ has order $p\leq k$. This is because: $a^{c}\cap b$ must equal $a^{c}$, so that $a\vee b=1$; and $a\cap b$ must have $p$ elements, and since $a$ has $k$ elements, this gives us $\binom{k}{p}$ solutions.

\begin{lemma}\label{a case of the binomial theorem}
For $m\in\mathbb{N}$ and $p\leq m$,
\[\sum_{j=p}^{m}\binom{j}{p}\binom{m}{j}x^{m-j}=\binom{m}{p}(x+1)^{m-p}.\]
\end{lemma}

\begin{proof}
We see that, from the definition of the binomial coefficients and the binomial theorem,
\[\sum_{j=p}^{m}\binom{j}{p}\binom{m}{j}x^{m-j}=\sum_{j=p}^{m}\frac{j!}{p!(j-p)!}\frac{m!}{j!(m-j)!}x^{m-j}=\]
\[\sum_{j=p}^{m}\frac{m!}{p!(j-p)!(m-j)!}x^{m-j}=\sum_{j=p}^{m}\frac{1}{p!}\frac{m!}{(j-p)!(m-j)!}x^{m-j}=\]
\[\sum_{j=p}^{m}\frac{m!}{p!(m-p)!}\frac{(m-p)!}{(j-p)!(m-j)!}x^{m-j}=\]
\[\binom{m}{p}\sum_{i=0}^{m-p}\frac{(m-p)!}{((i+p)-p)!(m-(i+p))!}x^{m-(i+p)}=\]
\[\binom{m}{p}\sum_{i=0}^{m-p}\frac{(m-p)!}{i!((m-p)-i)!}x^{(m-p)-i}=\binom{m}{p}\sum_{i=0}^{m-p}\binom{m-p}{i}x^{(m-p)-i}1^{i}=\]
\[\binom{m}{p}(x+1)^{m-p}.\]
\end{proof}

\begin{lemma}\label{counting certain snapshots}
For any $n\geq 0$, if $\mathcal{B}$ has $2^{m}$ elements then $B_{n}^{\mathcal{B}}$ has precisely $\binom{m}{p}(n+1)^{m-p}$ elements $(a_{1}, \dotsc  , a_{n+1})$ such that $\bigwedge_{i=1}^{n+1}a_{i}$ has order $p\leq m$.
\end{lemma}

\begin{proof}
Let $(a,b)$ be a pair in $B_{1}^{\mathcal{B}}$, meaning $a\vee b=1$, such that $a\wedge b$ has order $p$. Then, by the commentary after Lemma \ref{finding solutions of a given order}, if $a$ has order $k\geq p$, there are $\binom{k}{p}$ possible values for $b$, and if $k\leq p$ there are none. Since there are $\binom{m}{k}$ elements $a$ with order $k$, this gives us the number of pairs satisfying $a\vee b=1$ and $a\wedge b$ having order $p$ to be
\[\binom{p}{p}\binom{m}{p}+\binom{p+1}{p}\binom{m}{p+1}+\cdots+\binom{m}{p}\binom{m}{m},\]
which equals $\binom{m}{p}2^{m-p}$ from Lemma \ref{a case of the binomial theorem}, once we use $x=1$.

Now, suppose the described conditions hold for $B_{n}^{\mathcal{B}}$: there are, by induction hypothesis, $\binom{m}{k}(n+1)^{m-k}$ elements $(a_{1}, \dotsc  , a_{n+1})$ such that $\bigwedge_{i=1}^{n+1}a_{i}$ is of order $k$, and then there are $\binom{k}{p}$ possible values for $a_{n+3}$,  such that $(a_{1}, \dotsc  , a_{n+2})$ is in $B_{n+1}^{\mathcal{B}}$ (meaning $a_{n+2}\vee\bigwedge_{i=1}^{n+1}a_{i}=1$) and $\bigwedge_{i=1}^{n+2}a_{i}$ has order $p$, what gives the total
\[\binom{p}{p}\binom{m}{p}(n+1)^{m-p}+\cdots+\binom{m}{p}\binom{m}{m}(n+1)^{m-m},\]
which adds to $\binom{m}{p}(n+2)^{m-p}$ once we apply $x=n+1$ to Lemma \ref{a case of the binomial theorem}, ending our proof.
\end{proof}

\begin{theorem}
If $\mathcal{B}$ has $2^{m}$ elements, $B_{n}^{\mathcal{B}}$ has $(n+2)^{m}$ elements.
\end{theorem}

\begin{proof}
By Lemma \ref{counting certain snapshots}, there are $\binom{m}{0}(n+1)^{m-0}$ elements $(a_{1}, \dotsc  , a_{n+1})$ in $B_{n}^{\mathcal{B}}$ such that $\bigwedge_{i=1}^{n+1}a_{i}$ has order $0$, $\binom{m}{1}(n+1)^{m-1}$ elements such that $\bigwedge_{i=1}^{n+1}a_{i}$ has order $1$ and so on, adding up to a total of
\[\binom{m}{0}(n+1)^{m-0}+\binom{m}{1}(n+1)^{m-1}+\cdots+\binom{m}{m}(n+1)^{m-m}=\]
\[\sum_{p=0}^{m}\binom{m}{p}(n+1)^{m-p}1^{p}=((n+1)+1)^{m}=(n+2)^{m}.\]
\end{proof}

Notice, furthermore, that $D_{1}^{\mathcal{B}}$ has $2^{m}$ elements whenever $\mathcal{B}$ has $2^{m}$ elements: this happens since $(a, b)\in D_{1}^{\mathcal{B}}$ whenever $a=1$ and $a\vee b=1$, meaning $b$ can assume any value in $\mathcal{B}$; and when $n\geq 1$,
\[D_{n+1}^{\mathcal{B}}=\{(1, a_{1}, \dotsc  , a_{n+1})\in B_{n+1}^{\mathcal{B}}: (a_{1}, \dotsc  , a_{n+1})\in B_{n}^{\mathcal{B}}\},\]
implying $D_{n+1}^{\mathcal{B}}$ has as many elements as $B_{n}^{\mathcal{B}}$.

\begin{theorem}
If $\mathcal{B}$ has $2^{m}$ elements:
\begin{enumerate}
\item $D_{n}^{\mathcal{B}}$ has $(n+1)^{m}$ elements;
\item $Boo_{n}^{\mathcal{B}}$ has $2^{m}$ elements.
\end{enumerate}
\end{theorem}

A final, relevant, note is on the case that $\mathcal{B}$ is infinite, and of cardinality $\kappa$: since $B_{n}^{\mathcal{B}}$ contains the set of Boolean elements $Boo_{n}^{\mathcal{B}}$, isomorphic to $\mathcal{B}$ itself, we can be sure that $B_{n}^{\mathcal{B}}$ has at least cardinality $\kappa$; however, since $B_{n}^{\mathcal{B}}$ is a subset of $|\mathcal{B}|^{n+1}$, which is too of cardinality $\kappa$ given the assumption this is an infinite cardinal, we obtain that $B_{n}^{\mathcal{B}}$ has precisely $\kappa$ elements.

To provide one example of $B_{n}^{\mathcal{B}}$'s complexity, let us take the four-valued Boolean algebra $\mathcal{B}_{4}$ as the field of subsets of $\{a,b\}$, containing the elements $\emptyset$ (which we shall denote by $0$), $\{a\}$, $\{b\}$, and $\{a,b\}$ (which we will denote by $1$). Then, $B_{1}^{\mathcal{B}_{4}}$ has $9$ snapshots.
\begin{enumerate}
\item Designated and Boolean: $(1, 0)$.
\item Designated and not Boolean: $(1, \{a\})$, $(1, \{b\})$ and $(1, 1)$.
\item Undesignated and Boolean: $(0,1)$, $(\{a\}, \{b\})$ and $(\{b\}, \{a\})$.
\item Undesignated and not Boolean: $(\{a\}, 1)$ and $(\{b\}, 1)$.
\end{enumerate}

Still on the Boolean algebra with four elements, $B_{2}^{\mathcal{B}_{4}}$ has $16$ snapshots.
\begin{enumerate}
\item Designated and Boolean: $(1, 0, 1)$.
\item Designated and not Boolean: $(1, 1, 0)$, $(1, 1, \{a\})$, $(1, 1, \{b\})$, $(1, 1, 1)$, $(1, \{a\}, \{b\})$,\\ $(1, \{b\}, \{a\})$, $(1, \{a\}, 1)$ and $(1, \{b\}, 1)$.
\item Undesignated and Boolean: $(0, 1, 1)$, $(\{a\}, \{b\}, 1)$ and $(\{b\}, \{a\}, 1)$.
\item Undesignated and not Boolean: $(\{a\}, 1, 1)$, $(\{a\}, 1, \{b\})$, $(\{b\}, 1, 1)$ and $(\{b\}, 1, \{a\})$.

\end{enumerate}

Now, consider the eight-valued Boolean algebra $\mathcal{B}_{8}$, which we will take to be the Boolean algebra over the powerset of $\{a, b, c\}$. In that case, $B_{1}^{\mathcal{B}_{8}}$ has $27$ elements:

\begin{enumerate}
\item Designated and Boolean: $(1, 0)$.
\item Designated and not Boolean: $(1, \{a\})$, $(1, \{b\})$, $(1, \{c\})$, $(1, \{b, c\})$, $(1, \{a, c\})$, $(1, \{a, b\})$ and $(1, 1)$.
\item Undesignated and Boolean: $(0, 1)$, $(\{a\}, \{b, c\})$, $(\{b\}, \{a, c\})$, $(\{c\}, \{a, b\})$, $(\{b, c\}, \{a\})$, $(\{a, c\}, \{b\})$ and $(\{a, b\}, \{c\})$.
\item Undesignated and not Boolean: $(\{a\}, 1)$, $(\{b\}, 1)$, $(\{c\}, 1)$, $(\{b, c\}, \{a, c\})$, $(\{b, c\}, \{a, b\})$, $(\{b, c\}, 1)$, $(\{a, c\}, \{b, c\})$, $(\{a, c\}, \{a, b\})$, $(\{a, c\}, 1)$, $(\{a, b\}, \{b, c\})$, $(\{a, b\}, \{a, c\})$\\ and $(\{a, b\}, 1)$.
\end{enumerate}

\section{Category of restricted swap structures for $C_{n}$}

Given a class $C$ of RNmatrices, how to make it into a category $\mathcal{C}$? We believe there may exist more than one natural definition of what a morphism on $\mathcal{C}$ should be, the applications desired for such an object dictating the best ones, but at least one definition seems to be more or less universal. To start defining it, remember: an RNmatrix is a triple $\mathcal{M}=(\mathcal{A}, D, \mathcal{F})$, for $\mathcal{A}$ a $\Sigma$-multialgebra, $D$ a subset of its universe and $\mathcal{F}$ a set of homomorphisms (of multialgebras) $\nu:\textbf{F}(\Sigma, \mathcal{V})\rightarrow\mathcal{A}$. So a morphism between RNmatrices $\mathcal{M}=(\mathcal{A}, D, \mathcal{F})$ and $\mathcal{M}^{\prime}=(\mathcal{A}^{\prime}, D^{\prime}, \mathcal{F}^{\prime})$, over the same signature $\Sigma$, should be, at a minimum:
\begin{enumerate}
\item a homomorphism of multialgebras $h:\mathcal{A}\rightarrow\mathcal{A}^{\prime}$;
\item which preserves designated elements, meaning that $h(D)\subseteq D^{\prime}$;
\item and is absorbed by restricted valuations, meaning that for every $\nu\in\mathcal{F}$, $h\circ\nu\in\mathcal{F}^{\prime}$.
\begin{figure}[H]
\centering
\begin{tikzcd}
    \mathcal{M} \arrow{rr}{h}  &  & \mathcal{M}^{\prime}\\
    & \textbf{F}(\Sigma, \mathcal{V}) \arrow{ul}{\nu} \arrow{ur}{h\circ\nu} & 
  \end{tikzcd}
  \caption*{If $\nu\in\mathcal{F}$, $h\circ\nu$ must be in $\mathcal{F}^{\prime}$}
\end{figure}
\end{enumerate}

\begin{theorem}\label{category of RNmatrices}
A class of RNmatrices $C$, endowed with the morphisms defined above, becomes a category $\mathcal{C}$.
\end{theorem}

\begin{proof}
Take morphisms $g:\mathcal{M}\rightarrow\mathcal{M}^{\prime}$ and $h:\mathcal{M}^{\prime}\rightarrow\mathcal{M}^{\prime\prime}$ as described above, and first of all we will show that their composition, as a composition of functions, returns again such a morphism.
\begin{enumerate}
\item Since $g$ and $h$ are both $\Sigma$-homomorphisms of multialgebras, respectively from $\mathcal{A}$ to $\mathcal{A}^{\prime}$ and from $\mathcal{A}^{\prime}$ to $\mathcal{A}^{\prime\prime}$, it is clear that $h\circ g$ is a $\Sigma$-homomorphism from $\mathcal{A}$ to $\mathcal{A}^{\prime\prime}$.
\item Given that $g(D)\subseteq D^{\prime}$ and $h(D^{\prime})\subseteq D^{\prime\prime}$, we have that $h\circ g(D)=h(g(D))\subseteq h(D^{\prime})\subseteq D^{\prime\prime}$, meaning that $h\circ g$ preserves designated elements.
\item Finally, we have that for any $\nu\in\mathcal{F}$ and any $\nu^{\prime}\in\mathcal{F}^{\prime}$, $g\circ\nu\in\mathcal{F}^{\prime}$ and $h\circ\nu^{\prime}\in\mathcal{F}^{\prime\prime}$; thus, for any $\nu\in\mathcal{F}$, $g\circ \nu\in\mathcal{F}^{\prime}$ and so $(h\circ g)\circ \nu=h\circ(g\circ \nu)\in\mathcal{F}^{\prime\prime}$, proving $h\circ g$ is absorbed by restricted valuations.
\end{enumerate}

The associativity of the composition of these morphisms comes from the associativity of the composition of functions. The identity morphisms are precisely the identity functions, meaning that given $\mathcal{M}=(\mathcal{A}, D, \mathcal{F})$, the identity morphism in $\mathcal{M}$ is the identity function on the universe of $\mathcal{A}$, easily seem to be a homomorphism of multialgebras from $\mathcal{A}$ to itself; that preserves designated elements; and is absorbed by restricted valuations.
\end{proof}

Now, for a fixed da Costa's logic $C_{n}$, we take the class of RNmatrices $\mathcal{RM}^{\mathcal{B}}_{C_{n}}$ for $\mathcal{B}$ a non-trivial Boolean algebra and construct its corresponding category $\textbf{RSwap}_{C_{n}}$\label{RSwapCn} as it was done in Theorem \ref{category of RNmatrices} just above. To spell out our definition in this very important case, $\textbf{RSwap}_{C_{n}}$ is the category: 
\begin{enumerate}
\item with the class of (full) restricted swap structures $\mathcal{RM}_{C_{n}}^{\mathcal{B}}$, for all Boolean algebras $\mathcal{B}$, as objects; 
\item as morphisms between $\mathcal{RM}_{C_{n}}^{\mathcal{B}_{1}}$ and $\mathcal{RM}_{C_{n}}^{\mathcal{B}_{2}}$, all functions $\varphi:B_{n}^{\mathcal{B}_{1}}\rightarrow B_{n}^{\mathcal{B}_{2}}$ such that
\begin{enumerate}
\item $\varphi$ is a $\Sigma_{C}$-homomorphism between $\mathcal{A}_{C_{n}}^{\mathcal{B}_{1}}$ and $\mathcal{A}_{C_{n}}^{\mathcal{B}_{2}}$, seem as $\Sigma_{C}$-multialgebras;
\item for all $d\in D_{n}^{\mathcal{B}_{1}}$, $\varphi(d)$ is in $D_{n}^{\mathcal{B}_{2}}$;
\item for any $\nu:\textbf{F}(\Sigma_{C}, \mathcal{V})\rightarrow\mathcal{A}_{C_{n}}^{\mathcal{B}_{1}}$ in $\mathcal{F}_{C_{n}}^{\mathcal{B}_{1}}$, $\varphi\circ\nu:\textbf{F}(\Sigma_{C}, \mathcal{V})\rightarrow\mathcal{A}_{C_{n}}^{\mathcal{B}_{2}}$ is in $\mathcal{F}_{C_{n}}^{\mathcal{B}_{2}}$.
\end{enumerate}
\end{enumerate}

\begin{figure}[H]
\centering
\begin{tikzcd}
    \mathcal{A}_{C_{n}}^{\mathcal{B}_{1}} \arrow{rr}{\varphi}  &  & \mathcal{A}_{C_{n}}^{\mathcal{B}_{2}}\\
    & \textbf{F}(\Sigma_{C}, \mathcal{V}) \arrow{ul}{\nu} \arrow{ur}{\varphi\circ\nu} & 
  \end{tikzcd}
  \caption*{If $\nu\in\mathcal{F}_{C_{n}}^{\mathcal{B}_{1}}$, $\varphi\circ\nu$ must be in $\mathcal{F}_{C_{n}}^{\mathcal{B}_{2}}$}
\end{figure}

One important point to make is that demanding that a morphism of $\textbf{RSwap}_{C_{n}}$ preserves designated elements is actually superfluous: we could prove this property from the fact that the morphisms in this category are absorbed by restricted valuations. However, the proof of this fact is somewhat involved, so it is easier to assume the preservation of designated elements as one of the defining characteristics of our morphisms.

A first question regarding the definition of $\textbf{RSwap}_{C_{n}}$ could be whether there are truly morphisms $\varphi$ in it other than the identity ones.

\begin{proposition}\label{homomorphisms of boolean algebras lead to morphisms}
For Boolean algebras $\mathcal{B}_{1}$ and $\mathcal{B}_{2}$ and a homomorphism $\psi:\mathcal{B}_{1}\rightarrow\mathcal{B}_{2}$ of Boolean algebras, $\varphi:B_{n}^{\mathcal{B}_{1}}\rightarrow B_{n}^{\mathcal{B}_{2}}$ such that, for every $z=(z_{1}, \dotsc  , z_{n+1})\in B_{n}^{\mathcal{B}_{1}}$,
\[\varphi(z)_{i}=\psi(z_{i}),\quad\text{for every}\quad 1\leq i\leq n+1,\]
what we may write as $\varphi(z)=(\psi(z_{1}), \dotsc  , \psi(z_{n+1}))$, is a morphism between $\mathcal{RM}_{C_{n}}^{\mathcal{B}_{1}}$ and $\mathcal{RM}_{C_{n}}^{\mathcal{B}_{2}}$ in $\textbf{RSwap}_{C_{n}}$.
\end{proposition}

\begin{proof}
Notice that if $z$ is a Boolean element of $\mathcal{A}_{C_{n}}^{\mathcal{B}_{1}}$, then it is of the form $(a, {\sim}a, 1, \dotsc  , 1)$, and $\varphi(z)=(\psi(a), \psi({\sim}a), \psi(1), \dotsc  , \psi(1))=(\psi(a), {\sim}\psi(a), 1, \dotsc  , 1)$, given $\psi$ is a homomorphism; and since $\psi(a)$ is an element of $\mathcal{B}_{2}$, we see $\varphi(z)$ is a Boolean element of $\mathcal{A}_{C_{n}}^{\mathcal{B}_{2}}$, meaning $\varphi$ preserves Boolean elements. 

Take elements $w$ and $z$ in $B_{n}^{\mathcal{B}_{1}}$. First, suppose $w$ and $z$ are not both Boolean elements and $u\in w\tilde{\#}z$, meaning that $u_{1}=w_{1}\# z_{1}$: then 
\[\varphi(u)=(\psi(u_{1}), \psi(u_{2}), \dotsc  , \psi(u_{n+1}))\quad\text{equals}\quad(\psi(w_{1})\#\psi(z_{1}), \psi(u_{2}), \dotsc  , \psi(u_{n+1})),\]
since $\psi$ is a homomorphism, and since $\varphi(w)_{1}=\psi(w_{1})$ and $\varphi(z)_{1}=\psi(z_{1})$, we obtain that $\varphi(u)\in \varphi(w)\tilde{\#}\varphi(z)$. Now, if $w$ and $z$ are both Boolean elements and $u\in w\tilde{\#}z$, then $u_{1}=w_{1}\# z_{1}$ and $u$ is a Boolean element, and from what we saw above we find that $\varphi(u)_{1}=\varphi(w)_{1}\#\varphi(z)_{1}$ and that $\varphi(w)$, $\varphi(z)$ and $\varphi(u)$ are all Boolean elements, meaning $\varphi(u)\in\varphi(w)\tilde{\#}\varphi(z)$.

Now, for any $z\in B_{n}^{\mathcal{B}_{1}}$ and $w\in\tilde{\neg}z$, $w_{1}=z_{2}$ and $w_{2}\leq z_{1}$. Then $\varphi(w)=(\psi(w_{1}), \dotsc  , \psi(w_{n+1}))$, and analogously for $\varphi(z)$, meaning $\varphi(w)_{1}=\psi(w_{1})=\psi(z_{2})=\varphi(z)_{2}$ and $\varphi(w)_{2}=\psi(w_{2})\leq \psi(z_{1})=\varphi(z)_{1}$,\footnote{One should remember that homomorphisms of Boolean algebras preserve order: first of all, in a Boolean algebra $a\leq b$ iff $a=a\wedge b$ iff $b=a\vee b$; given a homomorphism $\psi:\mathcal{B}_{1}\rightarrow \mathcal{B}_{2}$, if $a\leq b$, $b=a\vee b$ and so $\psi(b)=\psi(a\vee b)=\psi(a)\vee\psi(b)$, meaning that $\psi(a)\leq \psi(b)$.}and therefore $\varphi(w)\in\tilde{\neg}\varphi(z)$. With this, $\varphi$ is a $\Sigma_{\textbf{C}}$-homomorphism between $\mathcal{A}_{C_{n}}^{\mathcal{B}_{1}}$ and $\mathcal{A}_{C_{n}}^{\mathcal{B}_{2}}$.

Given a designated element $z=(1, z_{2}, \dotsc  , z_{n+1})$ of $\mathcal{RM}_{C_{n}}^{\mathcal{B}_{1}}$, we get that $\varphi(z)=(\psi(1), \psi(z_{2}), \dotsc  , \psi(z_{n+1}))=(1, \psi(z_{2}), \dotsc  , \psi(z_{n+1}))$ is also designated, meaning that $\varphi$ preserves designated elements.

Now, a $\Sigma_{\textbf{C}}$-homomorphism $\nu:\textbf{F}(\Sigma_{\textbf{C}}, \mathcal{V})\rightarrow\mathcal{A}_{C_{n}}^{\mathcal{B}_{1}}$ is in $\mathcal{F}_{C_{n}}^{\mathcal{B}_{1}}$ when, for all formulas $\alpha$ and $\beta$:
\begin{enumerate}
\item $\nu(\alpha\wedge\neg\alpha)_{2}=\nu(\alpha)_{3}$;
\item $\nu(\alpha^{1})=(\nu(\alpha)_{3}, \nu(\alpha)_{1}\wedge\nu(\alpha)_{2}, \nu(\alpha)_{4}, \dotsc  , \nu(\alpha)_{n+1}, {\sim}(\bigwedge_{i=1}^{n+1}\nu(\alpha)_{i}))$;
\item and $\nu((\alpha^{(n)}\wedge\beta^{(n)})\rightarrow(\alpha\#\beta)^{(n)})\in D_{n}^{\mathcal{B}_{1}}$, for $\#\in\{\vee, \wedge, \rightarrow\}$.
\end{enumerate}

Quite clearly $\varphi\circ\nu:\textbf{F}(\Sigma_{\textbf{C}}, \mathcal{V})\rightarrow\mathcal{A}_{C_{n}}^{\mathcal{B}_{2}}$ is still a $\Sigma_{\textbf{C}}$-homomorphism.

\begin{enumerate}
\item We have $\varphi(\nu(\alpha\wedge\neg\alpha))_{2}=\psi(\nu(\alpha\wedge\neg\alpha)_{2})$, which equals, given $\nu$ is in $\mathcal{F}_{C_{n}}^{\mathcal{B}_{1}}$, $\psi(\nu(\alpha)_{3})=\varphi(\nu(\alpha))_{3}$.
\item See that $\varphi(\nu(\alpha^{1}))$ equals
\[\big(\psi(\nu(\alpha)_{3}), \psi(\nu(\alpha)_{1}\wedge\nu(\alpha)_{2}), \psi(\nu(\alpha)_{4}), \dotsc  , \psi(\nu(\alpha)_{n+1}), \psi({\sim}(\bigwedge_{i=1}^{n+1}\nu(\alpha)_{i}))\big)=\]
\[\big(\psi(\nu(\alpha)_{3}), \psi(\nu(\alpha)_{1})\wedge\psi(\nu(\alpha)_{2}), \psi(\nu(\alpha)_{4}), \dotsc  , \psi(\nu(\alpha)_{n+1}), {\sim}(\bigwedge_{i=1}^{n+1}\psi(\nu(\alpha)_{i}))\big)=\]
\[\big(\varphi(\nu(\alpha))_{3}, \varphi(\nu(\alpha))_{1}\wedge\varphi(\nu(\alpha))_{2}, \varphi(\nu(\alpha))_{4}, \dotsc  , \varphi(\nu(\alpha))_{n+1}, {\sim}(\bigwedge_{i=1}^{n+1}\varphi(\nu(\alpha))_{i})\big).\]
\item Since, as we already proved, $\varphi$ preserves designated elements, $\nu((\alpha^{(n)}\wedge\beta^{(n)})\rightarrow(\alpha\#\beta)^{(n)})\in D_{n}^{\mathcal{B}_{1}}$, for any $\#\in\{\vee, \wedge, \rightarrow\}$, implies that $\varphi(\nu((\alpha^{(n)}\wedge\beta^{(n)})\rightarrow(\alpha\#\beta)^{(n)}))\in D_{n}^{\mathcal{B}_{2}}$ for any $\#\in\{\vee, \wedge, \rightarrow\}$.
\end{enumerate}
This finishes proving that $\varphi\circ\nu$ lies in $\mathcal{F}_{C_{n}}^{\mathcal{B}_{2}}$.

\end{proof}

Notice that, in the spirit of the last proposition, the identity morphism $Id_{\mathcal{RM}_{C_{n}}^{\mathcal{B}}}$ of $\textbf{RSwap}_{C_{n}}$ on $\mathcal{RM}_{C_{n}}^{\mathcal{B}}$ can be written, for an arbitrary $z=(z_{1}, \dotsc  , z_{n+1})$ in $B_{n}^{\mathcal{B}}$, as
\[Id_{\mathcal{A}_{C_{n}}^{\mathcal{B}}}(z)=(z_{1}, \dotsc  , z_{n+1})=(Id_{\mathcal{B}}(z_{1}), \dotsc  , Id_{\mathcal{B}}(z_{n+1})),\]
for $Id_{\mathcal{B}}:\mathcal{B}\rightarrow\mathcal{B}$ the identity homomorphism on $\mathcal{B}$, and therefore the identity morphisms in $\textbf{RSwap}_{C_{n}}$ are, too, of the form described in the previous proposition.

A second question regarding the definition of $\textbf{RSwap}_{C_{n}}$ is whether there are morphisms in it other than those of the form described in Proposition \ref{homomorphisms of boolean algebras lead to morphisms}. The answer, in this case, is no, and in the following section we prove exactly that.

\subsection{Morphisms of $\textbf{RSwap}_{C_{n}}$}

For any function $\varphi:B_{n}^{\mathcal{B}_{1}}\rightarrow B_{n}^{\mathcal{B}_{2}}$ we will write, for any $z\in B_{n}^{\mathcal{B}_{1}}$, $\varphi(z)$ as $(\varphi_{1}(z), \dotsc  , \varphi_{n+1}(z))$: here, $\varphi_{i}$, for any $1\leq i\leq n+1$, is a function from $B_{n}^{\mathcal{B}_{1}}$ to $\mathcal{B}_{2}$ obtained by composing the $i$-th projection $\pi_{i}$\label{pii} of $B_{n}^{\mathcal{B}_{1}}$ with $\varphi$.

With this, we can say that $\varphi$ is a $\Sigma_{\textbf{C}}$-homomorphism from $\mathcal{A}_{C_{n}}^{\mathcal{B}_{1}}$ to $\mathcal{A}_{C_{n}}^{\mathcal{B}_{2}}$ if, and only if, for all snapshots $w, z\in B_{n}^{\mathcal{B}_{1}}$, and $u\in w\tilde{\#}z$ (for any $\#\in\{\vee, \wedge, \rightarrow\}$) and $v\in\tilde{\neg}z$, 
\[\varphi_{1}(u)=\varphi_{1}(w)\#\varphi_{1}(z),\quad\text{and}\quad\varphi_{1}(v)=\varphi_{2}(z)\quad\text{and}\quad\varphi_{2}(v)\leq\varphi_{1}(z),\]
what is equivalent to $\varphi(u)\in \varphi(w)\tilde{\#}\varphi(z)$ and $\varphi(v)\in\tilde{\neg}\varphi(z)$ once one considers the definitions of $\tilde{\#}$ and $\tilde{\neg}$.

So, we would like to study the function $\psi:\mathcal{B}_{1}\rightarrow\mathcal{B}_{2}$, given by $\psi(a)=$\\$\varphi((a, {\sim}a, 1, \dotsc  , 1))$ for any $a$ in $\mathcal{B}_{1}$, in the case that $\varphi$ is actually a $\Sigma_{\textbf{C}}$-homomorphism. Take a snapshot $z=(z_{1}, \dotsc  , z_{n+1})$ of $B_{n}^{\mathcal{B}_{1}}$ and consider $z^{*}=(z_{1}, {\sim}z_{1}, 1, \dotsc  , 1)$: then $\psi(z_{1})=\varphi(z^{*})$ by definition of $\psi$. Furthermore, recalling that $t^{n}_{0}=(1, 1, 0, 1, \dotsc  , 1)$ is a non-Boolean snapshot of $B_{n}^{\mathcal{B}}$, for any non-trivial Boolean algebra $\mathcal{B}$, we find that
\[z\tilde{\wedge}t^{n}_{0}=z^{*}\tilde{\wedge}t^{n}_{0}=\{w=(w_{1}, \dotsc  , w_{n+1})\in B_{n}^{\mathcal{B}_{1}}: w_{1}=z_{1}\},\]
and thus $z, z^{*}\in z\tilde{\wedge}t^{n}_{0}=z^{*}\tilde{\wedge}t^{n}_{0}$, since the first coordinates of both are precisely $z_{1}$; given $\varphi$ is a homomorphism, $z\in z\tilde{\wedge}t^{n}_{0}$ implies that $\varphi_{1}(z)=\varphi_{1}(z)\wedge\varphi_{1}(t^{n}_{0})$, while $z^{*}\in z\tilde{\wedge}t^{n}_{0}$ implies that $\varphi_{1}(z^{*})=\varphi_{1}(z)\wedge\varphi_{1}(t^{n}_{0})$, leading us to $\varphi_{1}(z)=\varphi_{1}(z^{*})$ and then to $\varphi_{1}(z)=\psi(z_{1})$. In other words, the function $\varphi_{1}$ depends exclusively on the first coordinate of a snapshot.

Furthermore, for any $a, b\in\mathcal{B}_{1}$, it is true for the snapshots $(a, {\sim}a, 1, \dotsc  , 1)$ and $(b, {\sim}b, 1, \dotsc  , 1)$ in $B_{n}^{\mathcal{B}_{1}}$ that 
\[(a, {\sim}a, 1, \dotsc  , 1)\tilde{\#}(b, {\sim}b, 1, \dotsc  , 1)=\{(a\#b, {\sim}a\#b, 1, \dotsc  , 1)\},\]
for any $\#\in\{\vee, \wedge, \rightarrow\}$. This way, taking into consideration that $\varphi$ is a homomorphism and $\varphi_{1}(z)=\psi(z_{1})$, $\psi(a\#b)=\psi(a)\#\psi(b)$. This means that $\psi$ is almost a homomorphism, and we will prove further ahead that it actually is a homomorphism of Boolean algebras.

For now,  let us also define the function $\theta:\mathcal{B}_{1}\rightarrow\mathcal{B}_{2}$ by $\theta(a)=\varphi_{2}(({\sim}a, a, 1, \dotsc  , 1))$. For an arbitrary snapshot $z=(z_{1}, \dotsc  , z_{n+1})$ in $B_{n}^{\mathcal{B}_{1}}$, consider 
\[z^{*}=({\sim}z_{2}, z_{2}, 1, \dotsc  , 1)\quad\text{and}\quad z^{\prime}=(z_{2}, {\sim}z_{2}, 1, \dotsc  , 1),\]
and by definition of $\theta$ we have that $\varphi_{2}(z^{*})=\theta(z_{2})$. Since $z$ is in $B_{n}^{\mathcal{B}_{1}}$, we have that $z_{1}\vee z_{2}=1$ and so ${\sim}z_{2}\leq z_{1}$, implying that $z^{\prime}\in\tilde{\neg}z$ and $z^{\prime}\in\tilde{\neg}z^{*}$. \footnote{Actually $\tilde{\neg}z^{*}=\{z^{\prime}\}$ and vice-versa, since both are Boolean snapshots with complementary first-coordinates.} Using, once again, that $\varphi$ is a homomorphism, $\varphi(z^{\prime})\in \tilde{\neg}\varphi(z)$ and $\varphi(z^{\prime})\in\tilde{\neg}\varphi(z^{*})$, leading to $\varphi_{1}(z^{\prime})=\varphi_{2}(z)$ and $\varphi_{1}(z^{\prime})=\varphi_{2}(z^{*})$, \textit{id est} $\varphi_{2}(z)=\varphi_{2}(z^{*})$. This leads us to $\varphi_{2}(z)=\theta(z_{2})$, meaning $\varphi_{2}$ also depends exclusively of one coordinate, in this case the second.

Finally, we may prove, with our only hypothesis being that $\varphi$ is a homomorphism, that $\psi=\theta$. To see this, take any $a\in\mathcal{B}_{1}$, and the snapshots $z=(a, {\sim}a, 1, \dotsc  , 1)$ and $z^{*}=({\sim}a, a, 1, \dotsc  , 1)$ of $B_{n}^{\mathcal{B}_{1}}$, with the property that $\tilde{\neg}z=\{z^{*}\}$ and $\tilde{\neg}z^{*}=\{z\}$. Then $\varphi(z)\in \tilde{\neg}\varphi(z^{*})$, so $\varphi_{1}(z)=\varphi_{2}(z^{*})$ (and also $\varphi_{2}(z)\leq \varphi_{2}(z^{*})$) and therefore $\psi(a)=\theta(a)$, using that $\varphi_{1}(z)=\psi(a)$ and $\varphi_{2}(z^{*})=\theta(a)$. We may summarize the results so far in the following theorem.

\begin{theorem}
If $\varphi:\mathcal{A}_{C_{n}}^{\mathcal{B}_{1}}\rightarrow\mathcal{A}_{C_{n}}^{\mathcal{B}_{2}}$ is a $\Sigma_{\textbf{C}}$-homomorphism, there exists a function $\psi:\mathcal{B}_{1}\rightarrow\mathcal{B}_{2}$ such that 
\[\varphi_{1}(z)=\psi(z_{1})\quad\text{and}\quad\varphi_{2}(z)=\psi(z_{2}),\]
for all snapshots $z=(z_{1}, \dotsc  , z_{n+1})$ of $B_{n}^{\mathcal{B}_{1}}$, and which satisfies, for all $a, b\in \mathcal{B}_{1}$ and $\#\in\{\vee, \wedge, \rightarrow\}$, $\psi(a\#b)=\psi(a)\#\psi(b)$.
\end{theorem}

Now we will demand that $\varphi:B_{n}^{\mathcal{B}_{1}}\rightarrow B_{n}^{\mathcal{B}_{2}}$ be, not only a homomorphism from $\mathcal{A}_{C_{n}}^{\mathcal{B}_{1}}$ to $\mathcal{A}_{C_{n}}^{\mathcal{B}_{1}}$, but a morphism of the category $\textbf{RSwap}_{C_{n}}$ from $\mathcal{RM}_{C_{n}}^{\mathcal{B}_{1}}$ to $\mathcal{RM}_{C_{n}}^{\mathcal{B}_{2}}$, so it also preserves designated elements and is absorbed by restricted valuations.

So, take a snapshot $z=(1, z_{2}, \dotsc  , z_{n+1})\in D_{n}^{\mathcal{B}_{1}}$: since we must have, by hypothesis over $\varphi$, $\varphi(z)\in D_{n}^{\mathcal{B}_{2}}$, and 
\[\varphi(z)=(\varphi_{1}(z), \dotsc  , \varphi_{n+1}(z))=(\psi(1), \psi(z_{2}), \varphi_{3}(z), \dotsc  , \varphi_{n+1}(z)),\]
we obtain that $\psi(1)=1$. Furthermore, for any formula $\alpha$ in the language of $C_{n}$ and valuation $\nu$ in $\mathcal{F}_{C_{n}}^{\mathcal{B}}$, Proposition \ref{Value of strong negation} shows that $\nu(\alpha\wedge\neg\alpha\wedge\alpha^{(n)})=F_{n}$; since $\varphi$ is absorbed by restricted valuations, for any $\nu\in\mathcal{F}_{C_{n}}^{\mathcal{B}_{1}}$ we have $\varphi\circ\nu\in\mathcal{F}_{C_{n}}^{\mathcal{B}_{2}}$ and so, for any formula $\alpha$ of $C_{n}$
\[F_{n}=\varphi\circ\nu(\alpha\wedge\neg\alpha\wedge\alpha^{(n)})=\varphi(F_{n})=(\psi(0), \psi(1), \varphi_{3}(F_{n}), \dotsc  , \varphi_{n+1}(F_{n})),\]
meaning that we also have $\psi(0)=0$. Since, for any $a, b\in\mathcal{B}_{1}$ and $\#\in\{\vee, \wedge, \rightarrow\}$, $\psi(a\#b)=\psi(a)\#\psi(b)$, $\psi(1)=1$ and $\psi(0)=0$, we actually have that $\psi$ is a homomorphism of Boolean algebras, as we had previously promised: given that ${\sim}a=a\rightarrow 0$, we have 
\[\psi({\sim}a)=\psi(a\rightarrow 0)=\psi(a)\rightarrow\psi(0)={\sim}\psi(a).\]

But the fact that $\varphi$ is absorbed by restricted valuations allows to prove something even stronger: remember that, for any $\nu\in\mathcal{F}^{\mathcal{B}}_{C_{n}}$ and formula $\alpha$ of $C_{n}$, $\nu(\alpha^{k})_{1}=\nu(\alpha)_{k+2}$ for $1\leq k\leq n-1$, what we proved back in Lemma \ref{Finding the coordinates} for $2\leq k\leq n-1$ and is obvious by the definition of $\mathcal{F}_{C_{n}}^{\mathcal{B}}$ for $k=1$. Given a $z$ in $B_{n}^{\mathcal{B}_{1}}$, for any propositional variable $p$ we can find a restricted valuation in $\mathcal{F}_{C_{n}}^{\mathcal{B}_{1}}$ with $\nu(p)=z$, by proceeding as in Section \ref{RM is not trivial} if necessary; so $\nu(p^{k})_{1}=z_{k+2}$.

Since $\varphi\circ\nu$ must be in $\mathcal{F}_{C_{n}}^{\mathcal{B}_{2}}$, we have that
\[\varphi\circ\nu(p^{k})_{1}=\varphi\circ\nu(p)_{k+2}=\varphi(z)_{k+2}=\varphi_{k+2}(z),\]
and at the same time,
\[\varphi\circ\nu(p^{k})_{1}=\psi(\nu(p^{k})_{1})=\psi(z_{k+2}),\]
implying that $\varphi_{k+2}(z)=\psi(z_{k+2})$, for every $1\leq k\leq n-1$. We get then the following theorem.

\begin{theorem}\label{Classifying morphisms in RSwap}
If $\varphi:\mathcal{RM}_{C_{n}}^{\mathcal{B}_{1}}\rightarrow\mathcal{RM}_{C_{n}}^{\mathcal{B}_{2}}$ is a morphism of $\textbf{RSwap}_{C_{n}}$, there exists a homomorphism of Boolean algebras $\psi:\mathcal{B}_{1}\rightarrow\mathcal{B}_{2}$ such that, for all snapshots $z\in B_{n}^{\mathcal{B}_{1}}$,
\[\varphi_{i}(z)=\psi(z_{i}),\quad\text{for all}\quad1\leq i\leq n+1.\]
\end{theorem}

\subsection{$\textbf{BA}$ and $\textbf{RSwap}_{C_{n}}$ are isomorphic}\label{BA and RSwap}

\begin{proposition}
The set $Boo_{n}^{\mathcal{B}}$ of Boolean elements of $\mathcal{A}_{C_{n}}^{\mathcal{B}}$ is a Boolean algebra isomorphic to $\mathcal{B}$ when we define, for $\#\in\{\vee, \wedge, \rightarrow\}$, and $a$ and $b$ elements of $\mathcal{B}$,
\[(a, {\sim}a, 1, \dotsc  , 1)\#(b, {\sim}b, 1, \dotsc  , 1)=(a\#b, {\sim}(a\#b), 1, \dotsc  , 1),\]
\[{\sim}(a, {\sim}a, 1, \dotsc  , 1)=({\sim}a, a, 1, \dotsc  , 1),\]
$\top=(1, 0, 1, \dotsc  , 1)$ and $\bot=(0, 1, 1, \dotsc  , 1)$
\end{proposition}

\begin{proof}
Consider the map $\rho:B\rightarrow Boo_{n}^{\mathcal{B}}$, for $B$ the universe of $\mathcal{B}$, given by $\rho(a)=$\\$(a, {\sim}a, 1, \dotsc  , 1)$. Then, for $a, b\in B$:
\begin{enumerate}
\item for $\#\in\{\vee, \wedge, \rightarrow\}$, 
\[\rho(a\#b)=(a\#b, {\sim}(a\#b), 1, \dotsc  , 1)=(a, {\sim}a, 1, \dotsc  , 1)\#(b, {\sim}b, 1, \dotsc  , 1)=\rho(a)\#\rho(b);\]
\item $\rho({\sim}a)=({\sim}a, {\sim}{\sim}a, 1, \dotsc  , 1)=({\sim}a, a, 1, \dotsc  , 1)={\sim}(a, {\sim}a, 1, \dotsc  , 1)={\sim}\rho(a)$;
\item $\rho(\top)=(1, 0, 1, \dotsc  , 1)=\top$;
\item $\rho(\bot)=(0, 1, 1, \dotsc  , 1)=\bot$.
\end{enumerate}
Furthermore, $\rho$ is clearly both injective and surjective, being therefore an isomorphism.
\end{proof}

Notice that the operations we give to $Boo_{n}^{\mathcal{B}}$ are the natural choice of Boolean algebra structure for this set, taking into consideration the multioperations on $\mathcal{A}_{C_{n}}^{\mathcal{B}}$, since for elements $a$ and $b$ of $\mathcal{B}$, and $\#\in\{\vee, \wedge, \rightarrow\}$, 
\[\rho(a)\tilde{\#}\rho(b)=\{\rho(a)\#\rho(b)\}\quad\text{and}\quad\tilde{\neg}\rho(a)=\{{\sim}\rho(a)\}.\]

Consider now the category $\textbf{BA}$\label{BA} of non-degenerate Boolean algebras\index{Boolean algebra, Non-degenerate} (that is, Boolean algebras where $0\neq 1$), with homomorphisms of Boolean algebras as morphisms, and the functors $\mathcal{A}_{n}:\textbf{BA}\rightarrow\textbf{RSwap}_{C_{n}}$\label{An}, which takes
\begin{enumerate}
\item a Boolean algebra $\mathcal{B}$ to $\mathcal{A}_{n}\mathcal{B}=\mathcal{A}_{C_{n}}^{\mathcal{B}}$;
\item a homomorphism of Boolean algebras $\psi:\mathcal{B}_{1}\rightarrow\mathcal{B}_{2}$ to the morphism $\mathcal{A}_{n}\psi:\mathcal{A}_{C_{n}}^{\mathcal{B}_{1}}\rightarrow\mathcal{A}_{C_{n}}^{\mathcal{B}_{2}}$ such that, for $z=(z_{1}, \dotsc  , z_{n+1})\in B_{n}^{\mathcal{B}_{1}}$, $\mathcal{A}_{n}\psi(z)_{i}=\psi(z_{i})$, for every $i\in\{1, \dotsc  , n+1\}$;
\end{enumerate}
and $\textbf{Boo}_{n}$\label{tBoon}, which takes
\begin{enumerate}
\item the restricted swap structure $\mathcal{A}_{C_{n}}^{\mathcal{B}}$ to the Boolean algebra $\mathcal{B}$ (or, what we saw to be equivalent, $Boo_{n}^{\mathcal{B}}$ with its natural Boolean algebra structure);
\item a morphism $\varphi:\mathcal{A}_{C_{n}}^{\mathcal{B}_{1}}\rightarrow\mathcal{A}_{C_{n}}^{\mathcal{B}_{2}}$ to the homomorphism of Boolean algebras $\textbf{Boo}_{n}\varphi:\mathcal{B}_{1}\rightarrow\mathcal{B}_{2}$ such that, for $a$ an element of $\mathcal{B}_{1}$, $\textbf{Boo}_{n}\varphi(a)=\varphi((a, {\sim}a, 1, \dotsc  , 1))_{1}$.
\end{enumerate}

\begin{proposition}
As described, $\mathcal{A}_{n}$ is, indeed, a functor.
\end{proposition}

\begin{proof}
As we proved in Theorem \ref{Classifying morphisms in RSwap}, for a homomorphism of Boolean algebras $\psi:\mathcal{B}_{1}\rightarrow\mathcal{B}_{2}$, $\mathcal{A}_{n}\psi:B_{n}^{\mathcal{B}_{1}}\rightarrow B_{n}^{\mathcal{B}_{2}}$ such that, for an arbitrary $z=(z_{1}, \dotsc  , z_{n+1})\in B_{n}^{\mathcal{B}_{1}}$, $\mathcal{A}_{n}\psi(z)_{i}=\psi(z_{i})$, for $1\leq i\leq n+1$, is indeed a morphism in $\textbf{RSwap}_{C_{n}}$.

Now, for a second homomorphism of Boolean algebras $\theta:\mathcal{B}_{2}\rightarrow\mathcal{B}_{3}$, one has
\[\mathcal{A}_{n}(\theta\circ\psi)(z)=(\theta\circ\psi(z_{1}), \dotsc  , \theta\circ\psi(z_{n+1}))=\mathcal{A}_{n}\theta((\psi(z_{1}), \dotsc  , \psi(z_{n+1}))=\]
\[\mathcal{A}_{n}\theta\circ\mathcal{A}_{n}\psi(z);\]
furthermore, for the identity $Id_{\mathcal{B}}:\mathcal{B}\rightarrow\mathcal{B}$ of $\mathcal{B}$, $\mathcal{A}_{n}Id_{\mathcal{B}}:\mathcal{A}_{C_{n}}^{\mathcal{B}}\rightarrow\mathcal{A}_{C_{n}}^{\mathcal{B}}$ satisfies that, for $z\in B_{n}^{\mathcal{B}}$, 
\[\mathcal{A}_{n}Id_{\mathcal{B}}(z)=(Id_{\mathcal{B}}(z_{1}), \dotsc  , Id_{\mathcal{B}}(z_{n+1}))=(z_{1}, \dotsc  , z_{n+1})=z,\]
and therefore is precisely the identity on $\mathcal{A}_{C_{n}}^{\mathcal{B}}$.
\end{proof}

\begin{proposition}
As described, $\textbf{Boo}_{n}$ is, indeed, a functor.
\end{proposition}

\begin{proof}
Since, for any morphism $\varphi:\mathcal{A}_{C_{n}}^{\mathcal{B}_{1}}\rightarrow\mathcal{A}_{C_{n}}^{\mathcal{B}_{2}}$, there exists a homomorphism $\psi:\mathcal{B}_{1}\rightarrow\mathcal{B}_{2}$ such that, for any $z=(z_{1}, \dotsc  , z_{n+1})\in B_{n}^{\mathcal{B}_{1}}$, $\varphi(z)_{i}=\psi(z_{i})$, for any $1\leq i\leq n+1$, one sees that 
\[\textbf{Boo}_{n}\varphi(a)=\varphi((a, {\sim}a, 1, \dotsc  , 1))_{1}=\psi(a),\]
and so $\textbf{Boo}_{n}\varphi$ is indeed a homomorphism of Boolean algebras.

For a second morphism $\eta:\mathcal{A}_{C_{n}}^{\mathcal{B}_{2}}\rightarrow \mathcal{A}_{C_{n}}^{\mathcal{B}_{3}}$ and a homomorphism $\theta:\mathcal{B}_{2}\rightarrow\mathcal{B}_{3}$ such that, for every $w\in B_{n}^{\mathcal{B}_{2}}$, $\eta(w)_{i}=\theta(w_{i})$, consider an element $a$ of $\mathcal{B}_{1}$: then
\[\textbf{Boo}_{n}(\eta\circ\varphi)(a)=\eta(\varphi((a, {\sim}a, 1, \dotsc  , 1)))_{1}=\eta_{1}(\varphi((a, {\sim}a, 1, \dotsc  , 1)))=\]
\[\theta(\varphi((a, {\sim}a, 1, \dotsc  , 1))_{1})=\theta(\varphi_{1}((a, {\sim}a, 1, \dotsc  , 1)))=\theta(\psi(a))=\theta(\textbf{Boo}_{n}\varphi(a))=\]
\[\textbf{Boo}_{n}\eta\circ\textbf{Boo}_{n}\varphi(a).\]

Finally, consider the identity homomorphism $Id_{\mathcal{A}_{C_{n}}^{\mathcal{B}}}:\mathcal{A}_{C_{n}}^{\mathcal{B}}\rightarrow\mathcal{A}_{C_{n}}^{\mathcal{B}}$: for an $a$ in $\mathcal{B}$,
\[\textbf{Boo}_{n}Id_{\mathcal{A}_{C_{n}}^{\mathcal{B}}}(a)=Id_{\mathcal{A}_{C_{n}}^{\mathcal{B}}}((a, {\sim}a, 1, \dotsc  , 1))_{1}=a,\]
proving $\textbf{Boo}_{n}Id_{\mathcal{A}_{C_{n}}^{\mathcal{B}}}$ is the identity on $\mathcal{B}$.
\end{proof}

\begin{theorem}
$\textbf{Boo}_{n}\circ\mathcal{A}_{n}=Id_{\textbf{BA}}$.
\end{theorem}

\begin{proof}
A Boolean algebra $\mathcal{B}$ is taken by $\mathcal{A}_{n}$ into $\mathcal{A}_{C_{n}}^{\mathcal{B}}$, and $\mathcal{A}_{C_{n}}^{\mathcal{B}}$ is taken to $\mathcal{B}$ again by $\textbf{Boo}_{n}$. This means $\textbf{Boo}_{n}\circ\mathcal{A}_{n}$ is the identity on objects.

Now, for Boolean algebras $\mathcal{B}_{1}$ and $\mathcal{B}_{2}$, a homomorphism $\psi:\mathcal{B}_{1}\rightarrow\mathcal{B}_{2}$, and an element $a$ of $\mathcal{B}_{1}$, let $\varphi=\mathcal{A}_{n}\psi$:
\[(\textbf{Boo}_{n}\circ\mathcal{A}_{n})\psi(a)=\textbf{Boo}_{n}\varphi(a)=\varphi((a, {\sim}a, 1, \dotsc  , 1))_{1}=\psi(a),\]
implying $\textbf{Boo}_{n}\circ\mathcal{A}_{n}$ is also identical when applied to morphisms.
\end{proof}

\begin{theorem}
$\mathcal{A}_{n}\circ\textbf{Boo}_{n}=Id_{\textbf{RSwap}_{C_{n}}}$.
\end{theorem}

\begin{proof}
$\mathcal{A}_{C_{n}}^{\mathcal{B}}$ is taken by $\textbf{Boo}_{n}$ to $\mathcal{B}$, which is taken back to $\mathcal{A}_{C_{n}}^{\mathcal{B}}$ by $\mathcal{A}_{n}$, giving us the identity on objects.

Now, for restricted swap structures $\mathcal{A}_{C_{n}}^{\mathcal{B}_{1}}$ and $\mathcal{A}_{C_{n}}^{\mathcal{B}_{2}}$, a morphism $\varphi:\mathcal{A}_{C_{n}}^{\mathcal{B}_{1}}\rightarrow\mathcal{A}_{C_{n}}^{\mathcal{B}_{2}}$ and an element $z=(z_{1}, \dotsc  , z_{n+1})\in B_{n}^{\mathcal{B}_{1}}$, let us denote $\textbf{Boo}_{n}\varphi$ by $\psi$: 
\[\mathcal{A}_{n}\circ\textbf{Boo}_{n}\varphi(z)=\mathcal{A}_{n}\psi(z)=(\psi(z_{1}), \dotsc  , \psi(z_{n+1}));\]
since $\textbf{Boo}_{n}\varphi=\psi$, for any $a$ in $\mathcal{B}_{1}$, $\varphi((a, {\sim}a, 1, \dotsc  , 1))_{1}=\psi(a)$. We also know that there exists a homomorphism $\theta:\mathcal{B}_{1}\rightarrow\mathcal{B}_{2}$ such that, for any $z\in B_{n}^{\mathcal{B}_{1}}$, $\varphi(z)_{i}=\theta(z_{i})$ for $1\leq i\leq n+1$.

From the fact that $\theta(a)=\varphi((a, {\sim}a, 1, \dotsc  , 1))_{1}=\psi(a)$, we obtain $\theta=\psi$, and so
\[\varphi(z)=(\psi(z_{1}), \dotsc  , \psi(z_{n+1}))=\mathcal{A}_{n}\circ\textbf{Boo}_{n}\varphi(z),\]
implying finally that $\mathcal{A}_{n}\circ\textbf{Boo}_{n}$ maintains morphisms fixed as well.
\end{proof}

Now that we have $\textbf{RSwap}_{C_{n}}$ is isomorphic to $\textbf{BA}$, we may use this isomorphism to bring some results from the latter category into the former: of course, $\textbf{BA}$ is one of the better known categories, so there are many such results. Let us start by reminding that every atomic, complete Boolean algebra is isomorphic to some power of the two-valued Boolean algebra $\textbf{2}$; in particular, an atomic, complete Boolean algebra with $\kappa$ atoms is isomorphic to precisely $\textbf{2}^{\kappa}$, where $\kappa$ may be a finite or infinite cardinal.

Since every finite Boolean algebra is forcibly both atomic and complete, and from the results of Section \ref{Counting snapshots} $\mathcal{RM}_{C_{n}}^{\mathcal{B}}$ is finite if and only if $\mathcal{B}$ is finite, we find a first representation result in $\textbf{RSwap}_{C_{n}}$.

\begin{corollary}
Every finite $\mathcal{RM}_{C_{n}}^{\mathcal{B}}$ is isomorphic to a power of $\mathcal{RM}_{C_{n}}$.
\end{corollary}

Another representation theorem on Boolean algebras, a stronger one, states that every Boolean algebra is isomorphic to a field of sets, that is, a Boolean subalgebra of the powerset Boolean algebra of a given set; since every powerset is a complete and atomic Boolean algebra, and is therefore isomorphic to some power $\textbf{2}^{\kappa}$ of $\textbf{2}$, it follows that every Boolean algebra is isomorphic to a subalgebra of a power of $\textbf{2}$.

Now, we still haven't defined what it means for an RNmatrix to be a "subRNmatrix" of another, although we could merely translate what this means from $\textbf{BA}$ into $\textbf{RSwap}_{C_{n}}$; we chose to give a more general approach instead. So, given RNmatrices $\mathcal{M}=(\mathcal{A}, D, \mathcal{F})$ and $\mathcal{M}^{\prime}=(\mathcal{A}^{\prime}, D^{\prime}, \mathcal{F}^{\prime})$ over the same signature, we say that $\mathcal{M}^{\prime}$ is a subRNmatrix\index{SubRNmatrix} of $\mathcal{M}$ if:
\begin{enumerate}
\item $\mathcal{A}^{\prime}$ is a submultialgebra of $\mathcal{A}$;
\item $D^{\prime}$ is a subset of $D$, and
\item $\{j\circ\nu: \nu\in\mathcal{F}^{\prime}\}\subseteq \mathcal{F}$, for $j$ the inclusion of the universe of $\mathcal{A}^{\prime}$ into the universe of $\mathcal{A}$.
\end{enumerate}
In other words, $\mathcal{M}^{\prime}$ is a subRNmatrix of $\mathcal{M}$ if $j$ is a morphism on the underlying category of RNmatrices.

Now, we state that in $\textbf{RSwap}_{C_{n}}$ the concept just described of a subRNmatrix corresponds exactly to what one would obtain by translating the notion of being a subalgebra from $\textbf{BA}$, that is, $\mathcal{RM}_{C_{n}}^{\mathcal{B}_{1}}$ is a subRNmatrix of $\mathcal{RM}_{C_{n}}^{\mathcal{B}_{2}}$ if and only if $\mathcal{B}_{1}$ is a subalgebra of $\mathcal{B}_{2}$.

In one direction, by assuming that $\mathcal{RM}_{C_{n}}^{\mathcal{B}_{1}}$ is a subRNmatrix of $\mathcal{RM}_{C_{n}}^{\mathcal{B}_{2}}$, we have that $\mathcal{A}_{C_{n}}^{\mathcal{B}_{1}}$ is a submultialgebra of $\mathcal{A}_{C_{n}}^{\mathcal{B}_{2}}$, so for every $a\in\mathcal{B}_{1}$, given that $(a, {\sim}a, 1, \dotsc  , 1)$ is an element of $B_{n}^{\mathcal{B}_{1}}$, we have that it must also be an element of $B_{n}^{\mathcal{B}_{2}}$, and therefore $a\in\mathcal{B}_{2}$. Since, for any $a, b\in\mathcal{B}_{1}$ and $\#\in\{\vee, \wedge, \rightarrow\}$, 
\[(a, {\sim}a, 1, \dotsc  , 1)\tilde{\#}(b, {\sim}b, 1, \dotsc  , 1)=\{(a\#b, {\sim}(a\#b), 1, \dotsc  , 1)\]
and $\tilde{\neg}(a, {\sim}a, 1, \dotsc  , 1)=\{({\sim}a, a, 1, \dotsc  , 1)\}$, and given that the operations $\tilde{\#}$ and $\tilde{\neg}$ in $\mathcal{A}_{C_{n}}^{\mathcal{B}_{1}}$ are the same as those in $\mathcal{A}_{C_{n}}^{\mathcal{B}_{2}}$, we obtain that the operations in $\mathcal{B}_{1}$ are the same as those in $\mathcal{B}_{2}$, proving the former is a Boolean subalgebra of the latter.

Reciprocally, suppose $\mathcal{B}_{1}$ is a Boolean subalgebra of $\mathcal{B}_{2}$: given a snapshot $z=(z_{1}, \dotsc  , z_{n+1})\in B_{n}^{\mathcal{B}_{1}}$, we have that, fist of all, $z_{1}, \dotsc  , z_{n+1}\in \mathcal{B}_{1}$, implying that $z_{1}, \dotsc  , z_{n+1}\in \mathcal{B}_{2}$; second, since $(\bigwedge_{i=1}^{k}z_{i})\vee z_{k+1}=1$, for all $1\leq k\leq n$, in $\mathcal{B}_{1}$, and since the operations in $\mathcal{B}_{2}$, when restricted to the universe of $\mathcal{B}_{1}$, coincide to the operations in $\mathcal{B}_{2}$, we obtain that $(\bigwedge_{i=1}^{k}z_{i})\vee z_{k+1}=1$, now in $\mathcal{B}_{2}$. With this, $z$ is a snapshot of $B_{n}^{\mathcal{B}_{2}}$.

So we can consider the inclusion $j:B_{n}^{\mathcal{B}_{1}}\rightarrow B_{n}^{\mathcal{B}_{2}}$, and it is easy to prove that it is a morphism of $\textbf{RSwap}_{C_{n}}$: after all, it can be written as $j(z)=(i(z_{1}), \dotsc  , i(z_{n+1}))$ for $i:\mathcal{B}_{1}\rightarrow\mathcal{B}_{2}$ the inclusion homomorphism, and $z=(z_{1}, \dotsc  , z_{n+1})$ an arbitrary snapshot in $B_{n}^{\mathcal{B}_{1}}$.

\begin{corollary}
Every $\mathcal{RM}_{C_{n}}^{\mathcal{B}}$ is a subRNmatrix of a power of $\mathcal{RM}_{C_{n}}$.
\end{corollary}

\newpage
\printbibliography[segment=\therefsegment,heading=subbibliography]
\end{refsegment}

\begin{refsegment}
\defbibfilter{notother}{not segment=\therefsegment}
\setcounter{chapter}{6}
\chapter{Logics of formal incompatibility}\label{Chapter7}\label{Chapter 7}

In classical logics, a formula and its negation are not compatible, in the sense that having both $\alpha$ and ${\sim} \alpha$ to be true trivialize whatever argument we are working over.

When dealing with logics of formal inconsistency, this is no longer true: we can have $\alpha$ and its negation $\neg\alpha$ without trivializing our logic, as long as $\alpha$ is inconsistent, that is, $\circ\alpha$ is not true, when we have at our disposal the consistency connective "$\circ$".

To formalize such a notion of incompatibility, we will consider a binary connective that, when connecting formulas $\alpha$ and $\beta$, will stand intuitively for "$\alpha$ is incompatible with $\beta$". When choosing a symbol for such connective, one somewhat natural choice is the Sheffer's stroke\index{Sheffer's stroke}: in classical propositional logic, the Sheffer's stroke may be defined from the usual connectives as
\[\alpha\Uparrow\beta={\sim}(\alpha\wedge\beta),\]
and of course having $\alpha\Uparrow\beta$ to hold, along with $\alpha$ and $\beta$, trivializes an argument. The basic axiom we will expect a system for incompatibility to satisfy will be\index{Incompatibility}\label{uparrow}
\[(\alpha\uparrow\beta)\rightarrow(\alpha\rightarrow(\beta\rightarrow\gamma)),\]
for any formula $\gamma$, or more generally, if we don't have a deduction meta-theorem,
\[\alpha\uparrow\beta, \alpha, \beta\vdash_{\mathscr{L}}\gamma;\]
in words, that means that having $\alpha$ and $\beta$ to be true while having $\alpha$ and $\beta$ to be incompatible imply that our logic is trivial.

Referring back to paraconsistency, one sees consistency may be characterized as incompatibility: $\alpha$ is consistent if, and only if, is incompatible with $\neg\alpha$. But, given a logic dealing with incompatibility, we are also tempted to say that any formula $\beta$ which is incompatible with a given $\alpha$ is a negation of $\alpha$, which leads to a notion of consistency and back again to logics having inconsistency in their scope.

Before we go any further, it is important to formally define with which structures we are working. We define the signature $\Sigma_{\bI}$\label{SigmabI} by: $(\Sigma_{\bI})_{2}=\{\vee, \wedge, \rightarrow, \uparrow\}$, and $(\Sigma_{\bI})_{n}=\emptyset$ for $n\neq 2$.

\begin{definition}\label{Logics of formal incompatibility}
A logic $\mathcal{L}$ is said to be a logic of formal incompatibility\index{Logic of formal incompatibility}, or in short $\textbf{LIp}$\label{LIp}, if it contains a family of formulas $\Uuparrow(p,q)$, exactly on the variables $p$ and $q$, such that there exist formulas $\varphi$, $\phi$ and $\psi$ for which
\[\varphi, \phi\not\vdash\psi;\]
formulas $\alpha$, $\beta$ and $\gamma$ such that:
\begin{enumerate}
\item $\Uuparrow(\alpha, \beta), \alpha\not\vdash\gamma$ and
\item $\Uuparrow(\alpha, \beta), \beta\not\vdash\gamma$,
\end{enumerate}
where $\Uuparrow(\alpha, \beta)$ is the set obtained by replacing, in each formula of $\Uuparrow(p,q)$, $p$ by $\alpha$ and $q$ by $\beta$; and, for all formulas $\gamma$, $\theta$ and $\omega$,
\[\Uuparrow(\gamma, \theta), \gamma, \theta\vdash\omega.\]
\end{definition}

More often than not, we want $\Uuparrow(p,q)$ to be composed of only one formula, given by a primitive, binary connective evaluated on $p$ and $q$ for which we will use the infix notation, that is, $p\uparrow q$.

Now, one can ask oneself when incompatibility has an interpretation in natural logics, and perhaps the most canonical example would be the one of limited resources\index{Limited resources}. Suppose we are given a basis of true statements, $\Gamma$, and formulas $\varphi$ that must be tested against $\Gamma$ in the most resource-efficient way possible, in regards to both time and amount of information recorded.

Now, if $\Gamma\vdash {\sim} \varphi$, we must discard $\varphi$: but, if the proof of this implication is unbelievably long, we waste precious time. One could, to avoid spending time unnecessarily, have a vast table pre-establishing if $\varphi$ follows from $\Gamma$ or not, but this would, of course, be memory-consuming, and in a way also time-consuming: without a strong algorithm to check the list, this could prove to be a very long search to identify $\varphi$ or ${\sim} \varphi$ on said list.

The solution could lie halfway between those two approaches: one could have a small list of statements that do not follow from $\Gamma$, which are incompatible with $\Gamma$, and, given a $\varphi$, attempt to prove or disprove the given formula, knowing a few shortcuts given by our list of incompatible statements.

This is very common not only in computer science, but also in science in general, mathematics, and logic: one does not completely prove a statement starting from the axioms, but rather accept a few results, derived from previous work, as true, and therefore accept their negations as incompatible with whatever result is being searched for, and proceed from there; it is also useful, in this context, to consider compatibility\index{Compatibility} of two formulas $\alpha$ and $\beta$ rather than incompatibility, defined trivially as $\alpha\downarrow\beta={\sim}(\alpha\uparrow\beta)$\label{downarrow}.

Another useful application of incompatibility is when dealing with partial information\index{Partial information}: suppose one must test $\varphi$ against a basis of true statements $\Delta$, without having direct access to $\Delta$ but rather $\Gamma\subset \Delta$; in this case, knowing some key incompatibilities between $\varphi$ and $\Delta$ or $\varphi$ and the complement of $\Delta$ may help deriving $\Delta\vdash \varphi$ or $\Delta\vdash{\sim} \varphi$ while only using directly $\Gamma$.

Perhaps more importantly, incompatibility has an obvious connection to probability, specifically independence (for what is to come, any reference in probability logic is sufficient, such as \cite{ProbabilityLogic}): consider a sample space $X$, that is, the set of all things that may occur in ones analysis; if, to give one example, we were to look at the outcome of flipping a coin, our sample space could be $X=\{heads, tails\}$. The next step, when dealing with probabilities, is taking a $\sigma$-algebra\footnote{A $\sigma$-algebra $\mathcal{A}$ on a set $X$ is a collection of subsets of $X$ containing $X$ itself and closed under complements and countable unions, meaning that if $A\in \mathcal{A}$ and $\{A_{n}\}_{n\in\mathbb{N}}\subseteq\mathcal{A}$, then $X\setminus A$ and $\bigcup_{n\in\mathbb{N}}A_{n}$ are both in $\mathcal{A}$.} $\mathcal{A}$ of subsets of $X$, most commonly the whole powerset of $X$. Finally, we also need a probability measure $P$ on $\mathcal{A}$ to obtain a probability space $(X, \mathcal{A}, P)$, that is, a function from $\mathcal{A}$ into the interval $[0,1]$ satisfying $P(X)=1$ and 
\[P(\bigcup_{n\in\mathbb{N}}A_{n})=\sum_{n\in\mathbb{N}}P(A_{n}),\]
for $\{A_{n}\}_{n\in\mathbb{N}}$ a family of pairwise disjoint sets of $\mathcal{A}$. Intuitively, $P$ tells us the probability of an event (\textit{i. e.} element of $\mathcal{A}$) happening, and $\mathcal{A}$ tells us to which events a probability can actually be assigned. We say two events $A$ and $B$ are mutually exclusive whenever 
\[\text{$P(A\cap B)=0$ or, what is equivalent, $P(A\cup B)=P(A)+P(B)$;}\]
intuitively, mutually exclusive events are those events that cannot both occur in an experiment. Then, if a propositional variable $p$ is associated to an event $A$, and a second propositional variable is associated to an event $B$, it is very natural to interpret the incompatibility $p\uparrow q$ of $p$ and $q$ as the mutual exclusivity of the events $A$ and $B$, in which case the axiom $(\alpha\uparrow\beta)\rightarrow(\alpha\rightarrow (\beta\rightarrow\gamma))$ clearly models the fact that both events $A$ and $B$ can not simultaneously occur if the two are independent.

This seems like a specially fruitful approach to apply to probability logics, systems used in the formal study of probability theory and statistics, as well as in Bayesian perspectives on epistemology of science. But, still on the field of probability theory, we can think of an alternative interpretation of incompatibility: two events $A$ and $B$ are independent if 
\[P(A\cap B)=P(A)P(B);\]
independent events naturally arise  from the fact that one can usually only derive more information about the probability of a complex event by knowing that its constituent events are pairwise independent, sometimes an even stronger independence condition being necessary. One could, then, for variables $p$ and $q$ standing for, respectively, events $A$ and $B$, interpret $p\uparrow q$ as the independence of $A$ and $B$, what makes sense at first glance due to the binary nature of the two concepts. However, in this interpretation, the axiom $(\alpha\uparrow\beta)\rightarrow(\alpha\rightarrow (\beta\rightarrow\gamma))$ is no longer the most desirable one, at least in a naive interpretation of implication, given that two independent events can, and often do, simultaneously happen.

Most of the research developed in this chapter can be found as a preprint in \cite{Frominconsistency}.


\section{The logic $\bI$}\label{Defining bI}

Our simplest $\textbf{LIp}$, which we shall denote by $\bI$\label{bI}, satisfies the axiom schemata of the positive fragment of classical propositional logic,

\begin{enumerate}
\item[\textbf{Ax\: 1}] $\alpha\rightarrow(\beta\rightarrow\alpha)$;
\item[\textbf{Ax\: 2}] $\big(\alpha\rightarrow (\beta\rightarrow \gamma)\big)\rightarrow\big((\alpha\rightarrow\beta)\rightarrow(\alpha\rightarrow\gamma)\big)$;
\item[\textbf{Ax\: 3}] $\alpha\rightarrow\big(\beta\rightarrow(\alpha\wedge\beta)\big)$;
\item[\textbf{Ax\: 4}] $(\alpha\wedge\beta)\rightarrow \alpha$;
\item[\textbf{Ax\: 5}] $(\alpha\wedge\beta)\rightarrow \beta$;
\item[\textbf{Ax\: 6}] $\alpha\rightarrow(\alpha\vee\beta)$;
\item[\textbf{Ax\: 7}] $\beta\rightarrow(\alpha\vee\beta)$;
\item[\textbf{Ax\: 8}] $(\alpha\rightarrow\gamma)\rightarrow\Big((\beta\rightarrow\gamma)\rightarrow \big((\alpha\vee\beta)\rightarrow\gamma\big)\Big)$;
\item[$\textbf{Ax\: 9}^{*}$] $(\alpha\rightarrow \beta)\vee\alpha$,
\end{enumerate}
plus\label{Ip}
\[\tag{\textbf{Ip}}(\alpha\uparrow\beta)\rightarrow(\alpha\rightarrow(\beta\rightarrow \gamma))\]
and\label{Comm}
\[\tag{\textbf{Comm}}(\alpha\uparrow\beta)\rightarrow(\beta\uparrow\alpha),\]
and follows the inference rule of Modus Ponens.

Compare it with the classical definition of $\textbf{mbC}$ and its extensions: they are logics over $\Sigma_{\textbf{LFI}}$, which is simply $\Sigma_{\bI}$ once we exchange $\uparrow$ for $\circ$ and add a negation, whose Hilbert calculus consists at least of the axiom schemata for the positive fragment of propositional logic, excluded middle and the schema
\[\tag{\textbf{bc1}}\circ\alpha\rightarrow(\alpha\rightarrow({\sim} \alpha\rightarrow \beta)),\]
with Modus Ponens as inference rule. The main difference here is the presence of $\textbf{Comm}$, standing for the commutativity\index{Commutativity of incompatibility} of the connective $\uparrow$. For convenience, we will denote the logic $\bI$ without $\textbf{Comm}$ by $\bI^{-}$\label{bI-}.

\begin{example}
We now present a model of $\bI$: take the signature $\Sigma_{\bI}$ and let $\textbf{CPL}$ be the classical propositional logic; we make $\textbf{CPL}$ into a model $\textbf{CPL}_{\uparrow}$ of $\bI$ by defining 
\[\alpha\uparrow\beta\quad\text{to be}\quad\alpha\Uparrow\beta={\sim}(\alpha\wedge\beta),\]
that is, $\uparrow$ is the classically defined Sheffer's stroke. Clearly $\textbf{CPL}_{\uparrow}$ satisfies the axioms of the positive fragment of classical propositional logic and Modus Ponens, remaining for us to show that it also validates $\textbf{Ip}$ and $\textbf{Comm}$.

This is easy since, if we have $\alpha\uparrow\beta$, $\alpha$ and $\beta$, and $\alpha\uparrow\beta$ means that ${\sim}(\alpha\wedge\beta)$, we have, by use of the deduction meta-theorem for classical propositional logic, that $(\alpha\wedge\beta)\wedge{\sim}(\alpha\wedge\beta)$, which by the explosivity of negation in $\textbf{CPL}$ implies any $\gamma$.

And since the conjunction is commutative in $\textbf{CPL}$, ${\sim}(\alpha\wedge\beta)\rightarrow{\sim}(\beta\wedge\alpha)$ or, what is equivalent, $(\alpha\uparrow\beta)\rightarrow(\beta\uparrow\alpha)$, proving $\textbf{Comm}$ is valid in $\textbf{CPL}_{\uparrow}$.
\end{example}


\subsection{Bottom and top elements, and classical negation}

Before anything else, we wish to show that $\bI$ has a formula equivalent to a bottom, and how we can use this to define a classical negation on it. So, for any two formulas $\alpha$ and $\beta$ in the language of $\bI$, consider
\[\bot_{\alpha\beta}=\alpha\wedge(\beta\wedge(\alpha\uparrow\beta)).\]

\begin{lemma}\label{Deduction meta-theorem for bI}
\begin{enumerate}
\item For a set $\Gamma\cup\{\alpha, \beta\}$ of formulas in $\bI$, we have that $\Gamma, \alpha\vdash_{\bI}\beta$ if and only if $\Gamma\vdash_{\bI}\alpha\rightarrow\beta$ (this is known as the deduction meta-theorem).
\item If $\Gamma\vdash_{\bI}\alpha\rightarrow\beta$ and $\Gamma\vdash_{\bI}\beta\rightarrow\gamma$, then $\Gamma\vdash_{\bI}\alpha\rightarrow\gamma$.
\item For a set $\Gamma\cup\{\alpha, \beta, \varphi\}$ of formulas in $\bI$, we have that, if $\Gamma, \alpha\vdash_{\bI}\varphi$ and $\Gamma, \beta\vdash_{\bI}\varphi$, then $\Gamma, \alpha\vee\beta\vdash_{\bI}\varphi$ (this is know as a proof by cases).
\end{enumerate}
\end{lemma}

\begin{proof}
\begin{enumerate}
\item The proof of this result for the positive fragment of $\textbf{CPL}$ can be found in \cite{Men-IntLog} (notice the proof does not use the negation), Proposition $1.9$, and the result extends to $\bI$ since this last logic has the axiom schemata of the previous one; but we make a point of proving this theorem given its importance.

Suppose first that $\Gamma, \alpha\vdash_{\bI}\beta$ and let $\varphi_{1}, \dotsc , \varphi_{n}=\beta$ be a demonstration of $\beta$ from $\Gamma\cup\{\alpha\}$; we show by induction that, for every $1\leq i\leq n$, $\Gamma\vdash_{\bI}\alpha\rightarrow\varphi_{i}$, and therefore $\Gamma\vdash_{\bI}\alpha\rightarrow\beta$ when we take $i=n$. 

The case $i=1$ follows from the general case, being an instance of an axiom or a premise, so let us assume that the result holds for $\varphi_{1}$ through $\varphi_{i}$, and prove it for $\varphi_{i+1}$; then $\varphi_{i+1}$ can be one of the following.

\begin{enumerate}
\item An instance of an axiom: in this case, since $\varphi_{i+1}$ and $\varphi_{i+1}\rightarrow(\alpha\rightarrow \varphi_{i+1})$ are instances of axioms, the second formula of $\textbf{Ax\: 1}$, we obtain by Modus Ponens that $\vdash_{\bI}\alpha\rightarrow\varphi_{i+1}$, and so $\Gamma\vdash_{\bI}\alpha\rightarrow\varphi_{i+1}$ since $\bI$ is tarskian.
\item A premise: $\varphi_{i+1}\rightarrow(\alpha\rightarrow \varphi_{i+1})$ certainly remains an instance of $\textbf{Ax\: 1}$, and through Modus Ponens $\varphi_{i+1}\vdash_{\bI}\alpha\rightarrow\varphi_{i+1}$, thus leaving us with $\Gamma\vdash_{\bI}\alpha\rightarrow\varphi_{i+1}$.
\item A conclusion of an inference rule: since we have only Modus Ponens as inference rule, this means that there exist $1\leq j<k\leq i$ with either $\varphi_{j}=\varphi_{k}\rightarrow\varphi_{i+1}$ or $\varphi_{k}=\varphi_{j}\rightarrow\varphi_{i+1}$, and we may assume without loss of generality the first case. 

From the induction hypothesis, $\Gamma\vdash_{\bI}\alpha\rightarrow\varphi_{j}$ (and so $\Gamma\vdash_{\bI}\alpha\rightarrow(\varphi_{k}\rightarrow\varphi_{i+1})$) and $\Gamma\vdash_{\bI}\alpha\rightarrow\varphi_{k}$, and from the instance
\[\big(\alpha\rightarrow(\varphi_{k}\rightarrow\varphi_{i+1})\big)\rightarrow\big((\alpha\rightarrow\varphi_{k})\rightarrow(\alpha\rightarrow\varphi_{i+1})\big)\]
of $\textbf{Ax\: 2}$ and two applications of Modus Ponens, we find $\Gamma\vdash_{\bI}\alpha\rightarrow\varphi_{i+1}$, as we needed to show.
\end{enumerate}

Reciprocally, if $\Gamma\vdash_{\bI}\alpha\rightarrow\beta$, take a demonstration $\psi_{1}, \dotsc , \psi_{m}=\alpha\rightarrow\beta$ of $\alpha\rightarrow\beta$ from $\Gamma$, and we state that $\psi_{1},\dotsc , \psi_{m}, \psi_{m+1}, \psi_{m+2}$ is a demonstration of $\beta$ from $\Gamma\cup\{\alpha\}$, where $\psi_{m+1}=\alpha$ and $\psi_{m+2}=\beta$. Of course, we have that $\psi_{1}, \dotsc , \psi_{m}$ remains a demonstration (of $\alpha\rightarrow\beta$) from $\Gamma\cup\{\alpha\}$, since this is a larger set of premises, and then:
\begin{enumerate}
\item $\psi_{m+1}=\alpha$ is a premise of $\Gamma\cup\{\alpha\}$;
\item and $\psi_{m+2}=\beta$ is the conclusion of an inference rule, with $\psi_{m}=\alpha\rightarrow\beta=$\\$\psi_{m+1}\rightarrow\psi_{m+2}$, meaning we are done.
\end{enumerate}

\item From $\Gamma\vdash_{\bI}\alpha\rightarrow\beta$ we get that $\Gamma, \alpha\vdash_{\bI}\beta$, and from $\Gamma\vdash_{\bI}\beta\rightarrow\gamma$, the fact that $\bI$ is tarskian and one application of Modus Ponens, we obtain that $\Gamma, \alpha\vdash_{\bI}\gamma$. So $\Gamma\vdash_{\bI}\alpha\rightarrow\gamma$.
\item If $\Gamma, \alpha\vdash_{\bI}\varphi$ and $\Gamma, \beta\vdash_{\bI}\varphi$, by the previous results of the lemma $\Gamma\vdash_{\bI}\alpha\rightarrow\varphi$ and $\Gamma\vdash_{\bI}\beta\rightarrow\varphi$; from the instance
\[(\alpha\rightarrow\varphi)\rightarrow((\beta\rightarrow\varphi)\rightarrow(\alpha\vee\beta)\rightarrow\varphi))\]
of $\textbf{Ax\: 8}$ and two consecutive applications of Modus Ponens, $\Gamma\vdash_{\bI}(\alpha\vee\beta)\rightarrow\varphi$. Again by the above results, this means $\Gamma, \alpha\vee\beta\vdash_{\bI}\varphi$.
\end{enumerate}
\end{proof}

\begin{proposition}
For any formulas $\alpha$, $\beta$, $\varphi$ and $\psi$ in $\bI$, $\bot_{\alpha\beta}$ and $\bot_{\varphi\psi}$ are equivalent, meaning $\vdash_{\bI}\bot_{\alpha\beta}\rightarrow\bot_{\varphi\psi}$ and $\vdash_{\bI}\bot_{\varphi\psi}\rightarrow\bot_{\alpha\beta}$.
\end{proposition}

\begin{proof}
From the deduction meta-theorem, we have both statements are equivalent to, respectively, $\bot_{\alpha\beta}\vdash_{\bI}\bot_{\varphi\psi}$ and $\bot_{\varphi\psi}\vdash_{\bI}\bot_{\alpha\beta}$, and from $\textbf{Ax\: 4}$, $\textbf{Ax\: 5}$ and $\textbf{Ax\: 3}$, those are in turn equivalent to $\alpha, \beta, \alpha\uparrow\beta\vdash_{\bI}\bot_{\varphi\psi}$ and $\varphi, \psi, \varphi\uparrow\psi\vdash_{\bI}\top_{\alpha\beta}$, both obviously true from $\textbf{Ip}$.
\end{proof}

It is also clear from this proof that each $\bot_{\alpha\beta}$ behaves like a bottom, since for any three formulas $\alpha$, $\beta$ and $\gamma$ we have that $\vdash_{\bI}\bot_{\alpha\beta}\rightarrow\gamma$.

Now, we define ${\sim}\alpha$ as the formula
\[\alpha\rightarrow\bot_{\alpha\alpha},\]
and state that is satisfies all the three axioms commonly assigned to classical negation, at least in $\textbf{CPL}$.

\begin{enumerate}
\item \mbox{}\vspace{-\baselineskip}\[(\alpha\rightarrow\beta)\rightarrow((\alpha\rightarrow{\sim} \beta)\rightarrow{\sim} \alpha)\]
Using several applications of the deduction meta-theorem, it is clear proving the validity of this statement is equivalent to proving that $\alpha\rightarrow\beta, \alpha\rightarrow{\sim} \beta\vdash_{\bI}{\sim} \alpha$. Now, $\alpha\rightarrow{\sim} \beta$ is simply $\alpha\rightarrow(\beta\rightarrow\bot_{\beta\beta})$, which by $\textbf{Ax\: 2}$ implies $(\alpha\rightarrow\beta)\rightarrow(\alpha\rightarrow\bot_{\beta\beta})$, and therefore 
\[\alpha\rightarrow\beta, \alpha\rightarrow{\sim} \beta\vdash_{\bI}\alpha\rightarrow\bot_{\beta\beta};\]
since $\vdash_{\textbf{mbC}}\bot_{\beta\beta}\rightarrow\bot_{\alpha\alpha}$, and the fact implication is transitive in $\bI$ from Lemma \ref{Deduction meta-theorem for bI}, we arrive at the desired result.

\item \mbox{}\vspace{-\baselineskip} \[\alpha\rightarrow({\sim} \alpha\rightarrow\beta)\]
Again by applying the deduction meta-theorem several times, we know it is enough to prove that $\alpha, {\sim} \alpha\vdash_{\textbf{mbC}}\beta$, and since ${\sim} \alpha=\alpha\rightarrow\bot_{\alpha\alpha}$, this is equivalent, by Modus Ponens, to $\bot_{\alpha\alpha}\vdash_{\textbf{mbC}}\beta$, which is true given all elements of the form $\bot_{\alpha\alpha}$ behave like bottoms.

\item \mbox{}\vspace{-\baselineskip} \[\alpha\vee{\sim} \alpha\]
Since ${\sim} \alpha=\alpha\rightarrow\bot_{\alpha\alpha}$, this is merely an instance of $\textbf{Ax\: 9}^{*}$.
\end{enumerate}

By defining, for a formula $\alpha$ in $\bI$,
\[\top_{\alpha}=\alpha\rightarrow\alpha,\]
we also have top elements, all equivalent to one another. To see those behave like top elements, for any two formulas $\alpha$ and $\beta$ we must show $\beta\rightarrow\top_{\alpha}$, or what is the same, that $\beta\rightarrow(\alpha\rightarrow\alpha)$. By the deduction meta-theorem, this is equivalent to $\beta\vdash_{\bI}\alpha\rightarrow\alpha$: it is clear how $\alpha\vdash_{\bI}\alpha\rightarrow\alpha$, from the instance 
\[\alpha\rightarrow(\alpha\rightarrow\alpha)\]
of $\textbf{Ax\: 1}$ and Modus Ponens; and trivially $\alpha\rightarrow\alpha\vdash_{\bI}\alpha\rightarrow\alpha$, which implies by a proof by cases that $(\alpha\rightarrow\alpha)\vee\alpha\vdash_{\bI}\alpha\rightarrow\alpha$. Since $(\alpha\rightarrow\alpha)\vee\alpha$ is an instance of $\textbf{Ax\: 9}^{*}$, we find that $\vdash_{\bI}\alpha\rightarrow\alpha$, and since $\bI$ is tarskian, $\beta\vdash_{\bI}\alpha\rightarrow\alpha$ as we wanted to prove.

And from this remark, the fact that all such top elements are equivalent is almost trivial: for any two formulas $\alpha$ and $\beta$ we have both that $\vdash_{\bI}\top_{\alpha}\rightarrow\top_{\beta}$ and $\vdash_{\bI}\top_{\beta}\rightarrow\top_{\alpha}$, since both $\top_{\alpha}$ and $\top_{\beta}$ behave like top elements, what means $\top_{\alpha}$ and $\top_{\beta}$ are equivalent.


\subsection{Bivaluations}\label{bivaluations for bI}

A \index{Bivaluation for $\bI$}bivaluation for $\bI$ is a map $\nu:F(\Sigma_{\bI}, \mathcal{V})\rightarrow \{0,1\}$ such that:
\begin{enumerate}
\item $\nu(\alpha\vee\beta)=1$ if and only if $\nu(\alpha)=1$ or $\nu(\beta)=1$;
\item $\nu(\alpha\wedge\beta)=1$ if and only if $\nu(\alpha)=\nu(\beta)=1$;
\item $\nu(\alpha\rightarrow\beta)=1$ if and only if $\nu(\alpha)=0$ or $\nu(\beta)=1$;
\item if $\nu(\alpha\uparrow\beta)=1$ and $\nu(\alpha)=1$, then $\nu(\beta)=0$;
\item $\nu(\alpha\uparrow\beta)=\nu(\beta\uparrow\alpha)$.
\end{enumerate}

Bivaluations for $\bI^{-}$ are simply bivaluations for $\bI$ that do not necessarily satisfy the condition that $\nu(\alpha\uparrow\beta)=\nu(\beta\uparrow\alpha)$, and the result that $\Gamma\vdash_{\bI^{-}}\varphi$ if and only if $\Gamma\vDash_{\bI^{-}}\varphi$ follows quite the same steps as those of the same result for $\bI$ instead.

Given a set of formulas $\Gamma\cup\{\varphi\}$ of $\bI$, we say that $\Gamma$ proves $\varphi$ semantically, and write $\Gamma\vDash_{\bI}\varphi$\label{vDashbI}, if for every bivaluation $\nu$ for $\bI$, $\nu(\Gamma)\subseteq\{1\}$ implies that $\nu(\varphi)=1$.


\subsubsection{Soundness}

We notice, first of all, that "$\vDash_{\bI}$" emulates Modus Ponens, meaning that if $\Gamma\vDash_{\bI}\alpha$ and $\Gamma\vDash_{\bI}\alpha\rightarrow\beta$, then $\Gamma\vDash_{\bI}\beta$: this is because, if $\nu$ is a bivaluation for $\bI$ such that $\nu(\Gamma)\subseteq\{1\}$, by hypothesis $\nu(\alpha)=\nu(\alpha\rightarrow\beta)=1$; since $\nu(\alpha\rightarrow\beta)=1$ implies that either $\nu(\alpha)=0$ or $\nu(\beta)=1$ and we already know $\nu(\alpha)=1$, we find that $\nu(\beta)=1$.

Then, it is possible to see that, if $\varphi$ is an instance of an axiom of classical propositional logic, then $\vDash_{\bI}\varphi$ (meaning that $\emptyset\vDash_{\bI}\varphi$); to give one example of how a proof of that would go, take the instance of axiom $\varphi=\alpha\rightarrow(\beta\rightarrow\alpha)$: if $\nu(\varphi)=0$, we must have that $\nu(\alpha)=1$ and $\nu(\beta\rightarrow\alpha)=0$, which in turn implies that $\nu(\beta)=1$ and $\nu(\alpha)=0$, making us reach a contradiction. To give another example, take an instance $\psi=(\alpha\rightarrow\beta)\vee\alpha$ of $\textbf{Ax\: 9}^{*}$: if $\nu(\psi)=0$, we must have $\nu(\alpha\rightarrow\beta)=\nu(\alpha)=0$; from $\nu(\alpha\rightarrow\beta)=0$ one gets in turn that $\nu(\alpha)=1$ and $\nu(\beta)=0$, what constitutes a contradiction. The proof for the other axioms follows analogously.

Finally, we state that if $\varphi$ is an instance of $\textbf{Ip}$ or $\textbf{Comm}$, then once again we have that $\vDash_{\bI}\varphi$: to see that, suppose $\nu$ is a bivaluation for which, by contradiction, $\nu(\varphi)=0$, implying that $\nu(\alpha\uparrow\beta)=1$ but
\[\nu(\alpha\rightarrow(\beta\rightarrow\gamma))=0;\]
this last equality implies that $\nu(\alpha)=1$ and $\nu(\beta\rightarrow\gamma)=0$, and thus $\nu(\beta)=1$ and $\nu(\gamma)=0$. But, at the same time, $\nu(\alpha\uparrow\beta)=1$ and $\nu(\alpha)=1$ imply together that $\nu(\beta)=0$, which contradicts our previous observations. One must have then that $\nu(\varphi)=1$, for any bivaluation $\nu$, and therefore $\vDash_{\bI}\varphi$.

Now, for $\textbf{Comm}$, if $\nu(\alpha\uparrow\beta)=0$, clearly we already have that $\nu((\alpha\uparrow\beta)\rightarrow(\beta\uparrow\alpha))=1$; if, otherwise, $\nu(\alpha\uparrow \beta)=1$, then by the fact $\nu$ is a bivaluation for $\bI$ we get $\nu(\beta\uparrow\alpha)=1$, and once again $\nu((\alpha\uparrow\beta)\rightarrow(\beta\uparrow\alpha))=1$.

\begin{theorem}\label{Proof of soundness for bI}
Given formulas $\Gamma\cup\{\varphi\}$ of $\bI$, if $\Gamma\vdash_{\bI}\varphi$ then $\Gamma\vDash_{\bI}\varphi$.
\end{theorem}

\begin{proof}
If $\Gamma\vdash_{\bI}\varphi$, there exists a demonstration $\alpha_{1}, \dotsc , \alpha_{n}$ of $\varphi$ from $\Gamma$, with $\alpha_{n}=\varphi$.

Let $\nu$ be a bivaluation for $\bI$ such that $\nu(\Gamma)\subseteq\{1\}$: we want to prove that, in this case, $\nu(\varphi)=1$; so we prove, by induction, that $\alpha_{1}$ through $\alpha_{n}$ have image $1$ under $\nu$, and therefore $\nu(\varphi)=\nu(\alpha_{n})=1$.

The formula $\alpha_{1}$ is either an axiom, when $\nu(\alpha_{1})=1$ since all instances of axioms have image $1$ through any bivaluations, or $\alpha_{1}$ is a premise, that is, an element of $\Gamma$, and since $\nu(\Gamma)\subseteq\{1\}$ we have that $\nu(\alpha_{1})=1$.

Suppose then that $\nu(\alpha_{1})=\cdots=\nu(\alpha_{i-1})=1$, and we have three cases to consider:
\begin{enumerate}
\item if $\alpha_{i}$ is an instance of an axiom, as we commented above $\nu(\alpha_{i})=1$;
\item if $\alpha_{i}$ is a premise, $\nu(\alpha_{i})=1$ since $\alpha_{i}\in\Gamma$ and $\nu(\Gamma)\subseteq\{1\}$;
\item if there are $\alpha_{j}$ and $\alpha_{k}$ with $1<j<k<i$ such that $\alpha_{j}=\alpha_{k}\rightarrow\alpha_{i}$ or $\alpha_{k}=\alpha_{j}\rightarrow\alpha_{i}$, since $\nu(\alpha_{j})=\nu(\alpha_{k})=1$ we find in both cases that $\nu(\alpha_{i})=1$, what ends the proof.
\end{enumerate}
\end{proof}


\subsubsection{Completeness}

Now, for a non-trivial, closed set of formulas $\Gamma$ maximal with respect to not proving $\varphi$, we want to prove that the function $\nu$, from the formulas of $\bI$ to $\{0,1\}$ and such that $\nu(\gamma)=1$ if and only if $\gamma\in\Gamma$, is a bivaluation.

\begin{enumerate}
\item If $\nu(\alpha\vee\beta)=1$, $\alpha\vee\beta\in\Gamma$: suppose, by contradiction, that $\alpha, \beta\notin\Gamma$; then $\Gamma, \alpha\vdash_{\bI}\varphi$ and $\Gamma, \beta\vdash_{\bI}\varphi$. So, by a proof by cases, $\Gamma, \alpha\vee\beta\vdash_{\bI}\varphi$, and since $\alpha\vee\beta\in\Gamma$ we find that $\Gamma\vdash_{\bI}\varphi$, which is absurd. This means either $\alpha\in\Gamma$, and then $\nu(\alpha)=1$, or $\beta\in\Gamma$, when $\nu(\beta)=1$.

Reciprocally, if $\nu(\alpha)=1$, since $\alpha\rightarrow(\alpha\vee\beta)$ is an instance of an axiom and is therefore in $\Gamma$, we find that $\alpha\vee\beta\in\Gamma$, since this set is closed; the same occurs if $\beta\in\Gamma$, forcing us to conclude $\nu(\alpha\vee\beta)=1$.

\item Now, assume $\nu(\alpha\wedge\beta)=1$; being $(\alpha\wedge\beta)\rightarrow\alpha$ an instance of an axiom, and $\Gamma$ closed, $(\alpha\wedge\beta)\rightarrow\alpha\in\Gamma$ and then $\alpha\in\Gamma$, and the same line of thought informs us that $\beta\in\Gamma$, meaning $\nu(\alpha)=\nu(\beta)=1$.

Reciprocally, if $\nu(\alpha)=\nu(\beta)=1$, $\alpha, \beta\in\Gamma$ and, given the instance $\alpha\rightarrow(\beta\rightarrow(\alpha\wedge\beta))$ of axiom $\textbf{Ax\: 3}$, by two consecutive applications of Modus Ponens we get that $\alpha\wedge\beta\in\Gamma$, meaning $\nu(\alpha\wedge\beta)=1$.

\item For implication, assume $\nu(\alpha\rightarrow\beta)=1$: then, either $\nu(\alpha)=1$, meaning $\alpha\in\Gamma$ and, since $\Gamma$ is closed and also contains $\alpha\rightarrow\beta$, $\beta\in\Gamma$ and therefore $\nu(\beta)=1$; or $\nu(\alpha)=0$.

Reciprocally, if $\nu(\beta)=1$ we have that $\beta\in\Gamma$, and since $\beta\rightarrow(\alpha\rightarrow\beta)$ is an instance of axiom $\textbf{Ax\: 1}$ and $\Gamma$ is closed, we obtain that $\alpha\rightarrow\beta\in\Gamma$.

If $\nu(\alpha)=0$, from the maximality of $\Gamma$ with respect to not proving $\varphi$ we find $\Gamma, \alpha\vdash_{\bI}\varphi$; suppose, by contradiction, that $\alpha\rightarrow\beta\notin\Gamma$, and again by the maximality of $\Gamma$ we discover $\Gamma, \alpha\rightarrow\beta\vdash_{\bI}\varphi$, and from a proof by cases
\[\Gamma, (\alpha\rightarrow\beta)\vee\alpha\vdash_{\bI}\varphi.\]
However, by $\textbf{Ax\: 9}^{*}$, $(\alpha\rightarrow\beta)\vee\alpha$ is an instance of an axiom, and we get that $\Gamma\vdash_{\bI}\varphi$, which is absurd. It follows that, if $\nu(\alpha)=0$, then $\nu(\alpha\rightarrow\beta)=1$.

\item Suppose $\nu(\alpha\uparrow\beta)=1$ and $\nu(\alpha)=1$, meaning that $\alpha\uparrow\beta, \alpha\in \Gamma$; by Modus Ponens and the correct instance of the axiom $\textbf{Ip}$,
\[(\alpha\uparrow\beta)\rightarrow(\alpha\rightarrow(\beta\rightarrow\varphi)),\]
we derive that $\alpha\rightarrow(\beta\rightarrow \varphi)\in\Gamma$ since this set is closed; using the closedness of $\Gamma$ and Modus Ponens again, $\beta\rightarrow\varphi\in\Gamma$, meaning, from what we saw before, that either $\beta\notin \Gamma$ or $\varphi\in\Gamma$. Since $\Gamma$ does not prove $\varphi$, we cannot have $\varphi\in\Gamma$, and therefore $\beta\notin\Gamma$, that is, $\nu(\beta)=0$.

\item Finally, if $\nu(\alpha\uparrow\beta)=0$, $\alpha\uparrow\beta\notin \Gamma$; if we had $\beta\uparrow\alpha\in \Gamma$, given that $\Gamma$ is closed and 
\[(\beta\uparrow\alpha)\rightarrow(\alpha\uparrow\beta)\]
is an instance of $\textbf{Comm}$, we would have $\alpha\uparrow\beta\in \Gamma$, which is a contradiction; so we must have $\beta\uparrow\alpha\notin\Gamma$ and therefore $\nu(\beta\uparrow\alpha)=0$.

If we have instead that $\nu(\alpha\uparrow\beta)=1$, $\alpha\uparrow\beta\in\Gamma$, and since $\Gamma$ is closed and $(\alpha\uparrow\beta)\rightarrow(\beta\uparrow\alpha)$ is an instance of axiom $\textbf{Comm}$, we get that $\beta\uparrow\alpha\in \Gamma$, meaning $\nu(\beta\uparrow\alpha)=1$.
\end{enumerate}

Therefore, the function $\nu$ is a bivaluation for $\bI$.

\begin{theorem}
Given formulas $\Gamma\cup\{\varphi\}$ in $\bI$, 
\[\Gamma\vdash_{\bI}\varphi\quad\text{if and only if}\quad\Gamma\vDash_{\bI}\varphi.\]
\end{theorem}

\begin{proof}
We have already proved in Theorem \ref{Proof of soundness for bI} that, if $\Gamma\vdash_{\bI}\varphi$, then $\Gamma\vDash_{\bI}\varphi$.

Reciprocally, suppose by contradiction $\Gamma\not\vdash_{\bI}\varphi$: then there exists a closed, non-trivial extension $\Delta$ of $\Gamma$ maximal with respect to not proving $\varphi$, and the function $\nu$ from the formulas of $\bI$ to $\{0,1\}$ such that $\nu(\delta)=1$ if and only if $\delta\in\Delta$ is a bivaluation.

But then $\nu$ is a bivaluation such that $\nu(\Gamma)\subseteq\{1\}$ and $\nu(\varphi)=0$, contradicting our assumption that $\Gamma\vDash_{\bI}\varphi$.

\end{proof}

Now, with the aid of bivaluations, we can prove that $\bI$ is a logic of formal incompatibility, or $\textbf{LIp}$, as in Definition \ref{Logics of formal incompatibility}, with respect to the set $\Uuparrow(p, q)=\{p\uparrow q\}$. Take any three propositional variables $p$, $q$ and $r$.
\begin{enumerate}
\item Taking a bivaluation $\nu$ with $\nu(p)=\nu(q)=1$ and $\nu(r)=0$, we then have that $p, q\not\vdash_{\bI}r$.
\item Taking a bivaluation $\nu$ with $\nu(p)=\nu(p\uparrow q)=1$ (and so $\nu(q)=0$ and $\nu(q\uparrow p)=1$) and $\nu(r)=0$, it becomes clear that $p\uparrow q, p\not\vdash_{\bI}r$.
\item Taking a bivaluation $\nu$ with $\nu(q)=\nu(p\uparrow q)=1$ (and so $\nu(p)=0$ and $\nu(q\uparrow p)=1$) and $\nu(r)=0$, it becomes clear that $p\uparrow q, q\not\vdash_{\bI}r$.
\item For any formula $\alpha$, it is true that $p\uparrow q, p, q\vdash_{\bI}\alpha$ from $\textbf{Ip}$ and Modus Ponens.
\end{enumerate}


\subsection{Fidel structures}\label{Fidel structures for bI}

Fidel structures are usually seem as generalized Boolean algebras, generalization originally achieved by adding relations to these algebraic structures. Designed to treat non-classical logics, their first occurrences in literature appear in the work of Fidel, as in \cite{Fidel3, Fidel, Fidel2}. They were first created to deal with certain paraconsistent logics, specifically da Costa's $C_{n}$, and in such a context they are probably more intuitively defined, but they will prove themselves very adequate for dealing semantically with logics of incompatibility, when correctly presented as multialgebras. When we focus our study on logics of formal inconsistency we will present a more classical approach to Fidel structures, as it can be found in Definition $6.1.2$ of \cite{ParLog}, but for now we will simply present the Fidel structures best tailored for $\bI$. 

Here, it becomes necessary to clarify the use of the nomenclature ``Fidel structures'': why are we using this name? These semantics were originally defined as certain implicative lattices\footnote{That, in the case of $C_{n}$, were proven to actually be Boolean algebras.} equipped with some unary relations to treat da Costa's hierarchy plus $C_{\omega}$. Then, as more non-classical logics needed semantical characterizations, all Boolean-like algebras (including implicative lattices, Heyting algebras and others), equipped with unary relations (usually standing for negation and consistency), could be called Fidel structures; the earliest examples of Fidel structures over algebras that are not Boolean are also found in Fidel's work, for da Costa's $C_{\omega}$ and Nelson's logics, in respectively \cite{Fidel3} and \cite{FidelNelson}. Here, our reasoning is the following: Fidel structures should be seem as RNmatrices with an underlying Boolean-like (deterministic) component, and multioperations (now of any arity) arising from what was previously presented as relations. Perhaps it is too much to classify under the title of ``Fidel structures'', but since the logics at hand are still quite similar, in spirit, to da Costa's hierarchy, the name seemed a fitting tribute to Fidel's breakthrough semantics.

Firstly, we expand the signature $\Sigma_{\bI}$ by adding a bottom, a top and a classical negation, making the signature $\Sigma_{\bI}^{\textbf{CPL}}$\label{SigmabICPL} such that: $(\Sigma_{\bI}^{\textbf{CPL}})_{0}=\{\bot, \top\}$, $(\Sigma_{\bI}^{\textbf{CPL}})_{1}=\{{\sim}\}$, $(\Sigma_{\bI}^{\textbf{CPL}})_{2}=\{\vee, \wedge, \rightarrow, \uparrow\}$, and $(\Sigma_{\bI}^{\textbf{CPL}})_{n}=\emptyset$ for $n>2$. This allows us to define a \index{Fidel structure for $\bI$}Fidel structure, presented as a $\Sigma_{\bI}^{\textbf{CPL}}$-multialgebra, for $\bI$ to be any $\Sigma_{\bI}^{\textbf{CPL}}$-multialgebra $\mathcal{A}=(A, \{\sigma_{\mathcal{A}}\}_{\sigma\in\Sigma_{\bI}^{\textbf{CPL}}})$ such that:
\begin{enumerate}
\item $(A, \{\sigma_{\mathcal{A}}\}_{\sigma\in\Sigma^{\textbf{CPL}}})$ is a Boolean algebra;\footnote{Here, a distinction is important: the operations corresponding to the symbols $\bot$, $\top$, $\sim$, $\vee$, $\wedge$ and $\rightarrow$ are deterministic, and will be treated as such; meanwhile, $\uparrow$ will be, at most times, non-deterministic.}
\item for all $a, b\in A$ and $c\in\uparrow_{\mathcal{A}}(a,b)$, 
\[\wedge_{\mathcal{A}}(a, \wedge_{\mathcal{A}}(b,c))=\bot_{\mathcal{A}};\]
\item for all $a, b\in A$, $\uparrow_{\mathcal{A}}(a, b)=\uparrow_{\mathcal{A}}(b,a)$.\footnote{Technically, this last condition could be replaced by, for all $a, b\in A$, $\uparrow_{\mathcal{A}}(a,b)\cap\uparrow_{\mathcal{A}}(b,a)\neq\emptyset$, or simply erased from our definition of Fidel structures for $\bI$; it only guarantees that, for every such Fidel structure, any function from the variables $\mathcal{V}$ into $A$ may be extended to a restricted valuation. However, given the condition $\nu(\alpha\uparrow\beta)=\nu(\beta\uparrow\alpha)$ we will ask of these restricted valuations, dropping $\uparrow_{\mathcal{A}}(a, b)=\uparrow_{\mathcal{A}}(b,a)$ from the definition of a Fidel structure would only mean that some Fidel structure would not help in establishing the validity of an argument; the class of all Fidel structures, however, would still be sound and complete.}
\end{enumerate}

For simplicity, we will drop the indexes $\mathcal{A}$ and use the standard infix notation.

Given a Fidel structure $\mathcal{A}$, presented as a $\Sigma_{\bI}^{\textbf{CPL}}$-multialgebra, for $\bI$, a valuation over $\mathcal{A}$ is a function $\nu:F(\Sigma_{\bI}, \mathcal{V})\rightarrow A$ such that:
\begin{enumerate}
\item $\nu(\alpha\#\beta)=\nu(\alpha)\#\nu(\beta)$, for $\#\in\{\vee, \wedge, \rightarrow\}$;
\item $\nu(\alpha\uparrow\beta)\in\nu(\alpha)\uparrow\nu(\beta)$.
\end{enumerate}

Notice that $\nu$ is a valuation for $\mathcal{A}$ if and only if it is a $\Sigma_{\bI}$-homomorphism, between $\textbf{F}(\Sigma_{\bI}, \mathcal{V})$ and $(A, \{\sigma_{\mathcal{A}}\}_{\sigma\in\Sigma_{\bI}})$. 

For every Fidel structure $\mathcal{A}$ for $\bI$, we will consider the restricted Nmatrix\\ $(\mathcal{A}, \{\top\}, \mathcal{F}_{\mathcal{A}})$, where $\mathcal{F}_{\mathcal{A}}$ is the set of valuations $\nu:\textbf{F}(\Sigma_{\bI}, \mathcal{V})\rightarrow\mathcal{A}$ such that 
\[\nu(\alpha\uparrow\beta)=\nu(\beta\uparrow\alpha),\]
for any two formulas $\alpha$ and $\beta$ in $F(\Sigma_{\bI}, \mathcal{V})$. If $\Gamma$ proves $\varphi$ according to such restricted Nmatrices, we will write $\Gamma\Vdash_{\mathcal{F}}^{\bI}\varphi$\label{VdashFbI}: as one could suspect, $\Gamma\vdash_{\bI}\varphi$ if and only if $\Gamma\Vdash_{\mathcal{F}}^{\bI}\varphi$.

\begin{proposition}
Each restricted Nmatrix $(\mathcal{A}, \{\top\}, \mathcal{F}_{\mathcal{A}})$, as described above, is structural.
\end{proposition}

\begin{proof}
Take any $\Sigma_{\bI}$-homomorphism $\sigma:\textbf{F}(\Sigma_{\bI}, \mathcal{V})\rightarrow\textbf{F}(\Sigma_{\bI}, \mathcal{V})$ and an element $\nu\in\mathcal{F}_{\mathcal{A}}$, which is by definition of $\mathcal{F}_{\mathcal{A}}$ a $\Sigma_{\bI}$-homomorphism $\nu:\textbf{F}(\Sigma_{\bI}, \mathcal{V})\rightarrow \mathcal{A}$ such that, for any two formulas $\alpha$ and $\beta$, $\nu(\alpha\uparrow\beta)=\nu(\beta\uparrow\alpha)$.

Clearly $\nu\circ\sigma:\textbf{F}(\Sigma_{\bI}, \mathcal{V})\rightarrow \mathcal{A}$ is a $\Sigma_{\bI}$-homomorphism, so we only need to show that, for any two formulas $\alpha$ and $\beta$, $\nu\circ\sigma(\alpha\uparrow\beta)=\nu\circ\sigma(\beta\uparrow\alpha)$; and since 
\[\nu\circ\sigma(\alpha\uparrow\beta)=\nu(\sigma(\alpha\uparrow\beta))=\nu(\sigma(\alpha)\uparrow\sigma(\beta))=\nu(\sigma(\beta)\uparrow\sigma(\alpha))=\nu(\sigma(\beta\uparrow\alpha))=\nu\circ\sigma(\beta\uparrow\alpha),\]
we finish the proof.
\end{proof}

If we simply consider the class $\mathbb{M}$ of Nmatrices $(\mathcal{A},\{\top\})$, where $\mathcal{A}$ is a Fidel structure for $\bI$, it is not hard to prove that $\Gamma\vDash_{\mathbb{M}}\varphi$ if and only if $\Gamma\vdash_{\bI^{-}}\varphi$.


\subsubsection{Soundness}

First of all, we state that for any instance of an axiom $\alpha$ of the positive fragment of classical propositional logic in $\bI$, $\Vdash_{\mathcal{F}}^{\bI}\alpha$ (meaning that $\emptyset\Vdash_{\mathcal{F}}^{\bI}\alpha$): this is true because Boolean algebras model classical propositional logic. To see one example, take an instance $\alpha\rightarrow(\beta\rightarrow\alpha)$ of $\textbf{Ax\: 1}$: we have that
\[\nu(\alpha\rightarrow(\beta\rightarrow\alpha))=\nu(\alpha)\rightarrow(\nu(\beta)\rightarrow\nu(\alpha)),\]
and remembering that in a Boolean algebra, $x\rightarrow y={\sim}x\vee y$, this equals ${\sim}\nu(\alpha)\vee({\sim}\nu(\beta)\vee\nu(\alpha))$; using that $\vee$ is commutative and associative, in this order, this equals
\[({\sim}\nu(\alpha)\vee\nu(\alpha))\vee{\sim}\nu(\beta)=\top\vee{\sim}\nu(\beta)=\top,\]
and since for any valuation $\nu$ we have that $\nu(\alpha\rightarrow(\beta\rightarrow\alpha))=\top$, we find $\Vdash_{\mathcal{F}}^{\bI}\alpha\rightarrow(\beta\rightarrow\alpha)$.

Now we take an instance $(\alpha\uparrow\beta)\rightarrow(\alpha\rightarrow(\beta\rightarrow\gamma))$ of axiom $\textbf{Ip}$; using again that $x\rightarrow y={\sim}x\vee y$, the image under a valuation $\nu$ of this formula is 
\[{\sim}\nu(\alpha\uparrow\beta)\vee({\sim}\nu(\alpha)\vee({\sim}\nu(\beta)\vee\nu(\gamma)))=({\sim}\nu(\alpha\uparrow\beta)\vee({\sim}\nu(\alpha)\vee{\sim}\nu(\beta)))\vee\nu(\gamma);\]
using that ${\sim}x\vee{\sim}y={\sim}(x\wedge y)$ in a Boolean algebra, one of the De Morgan's laws, we get this expression equals
\[{\sim}(\nu(\alpha\uparrow\beta)\wedge(\nu(\alpha)\wedge\nu(\beta)))\vee\nu(\gamma),\]
and from the requirements made on the multioperation $\uparrow$ in the definition of Fidel structures for $\bI$, we get that $\nu(\alpha\uparrow\beta)\wedge(\nu(\alpha)\wedge\nu(\beta))=\bot$, and therefore the whole expression simplifies to
\[{\sim}\bot\vee\nu(\gamma)=\top\vee\nu(\gamma)=\top.\]

Finally, take an instance $(\alpha\uparrow\beta)\rightarrow(\beta\uparrow\alpha)$ of $\textbf{Comm}$, and for any restricted Nmatrix $(\mathcal{A},\{\top\},\mathcal{F}_{\mathcal{A}})$ for $\bI$, its image under a $\nu\in\mathcal{F}_{\mathcal{A}}$ is $\nu(\alpha\uparrow\beta)\rightarrow\nu(\beta\uparrow\alpha)$. Since $\nu\in\mathcal{F}_{\mathcal{A}}$, we have that $\nu(\alpha\uparrow\beta)=\nu(\beta\uparrow\alpha)$, and therefore $\nu((\alpha\uparrow\beta)\rightarrow(\beta\uparrow\alpha))$ equals $\nu(\alpha\uparrow\beta)\rightarrow\nu(\alpha\uparrow\beta)$, which is always $\top$, meaning "$\Vdash_{\mathcal{F}}^{\bI}$" models all axioms of $\bI$.

\begin{theorem}\label{Soundness of Fidel structures for bI}
Given formulas $\Gamma\cup\{\varphi\}$ of $\bI$, if $\Gamma\vdash_{\bI}\varphi$ then $\Gamma\Vdash_{\mathcal{F}}^{\bI}\varphi$.
\end{theorem}

\begin{proof}
Let $\alpha_{1}, \dotsc , \alpha_{n}$ be a demonstration of $\varphi$ from $\Gamma$, with $\alpha_{n}=\varphi$, and let $\mathcal{A}$ be a Fidel structure, represented as a $\Sigma_{\bI}$-multialgebra, for $\bI$, with $\nu$ a valuation for $\bI$ over $\mathcal{A}$ in $\mathcal{F}_{\mathcal{A}}$.

We aim to show that, if $\nu(\Gamma)\subseteq\{\top\}$, then $\nu(\alpha_{i})=\top$ for every $i\in\{1, \dotsc , n\}$. We will proceed by induction, assuming that the result already holds for all formulas prior to $\alpha_{i}$ (notice that $\alpha_{0}$ is either an axiom or a premise). Then we have that $\alpha_{i}$ is either:
\begin{enumerate}
\item an instance of an axiom, when we already proved $\Vdash_{\mathcal{F}}^{\bI}\alpha_{i}$, and therefore $\nu(\alpha_{i})=\top$;
\item a premise, and by our hypothesis that $\nu(\Gamma)\subseteq\{\top\}$ we find $\nu(\alpha_{i})=\top$;
\item there exist $j,k<i$ such that either $\alpha_{j}=\alpha_{k}\rightarrow\alpha_{i}$ or $\alpha_{k}=\alpha_{j}\rightarrow\alpha_{i}$, where we will assume, without loss of generality, that the first case holds;

by induction hypothesis, $\nu(\alpha_{j})=\nu(\alpha_{k}\rightarrow\alpha_{i})=\top$, and therefore ${\sim}\nu(\alpha_{k})\vee\nu(\alpha_{i})=\top$, and since $\nu(\alpha_{k})=\top$, meaning ${\sim}\nu(\alpha_{k})=\bot$, we must have that $\nu(\alpha_{i})=\top$, what finishes the proof.
\end{enumerate}
\end{proof}


\subsubsection{Completeness}\label{Completeness fo Fidel structures for bI}

Now, we want to show the reciprocal of Theorem \ref{Soundness of Fidel structures for bI} is also true: if $\Gamma\Vdash_{\mathcal{F}}^{\bI}\varphi$, then $\Gamma\vdash_{\bI}\varphi$.

For a set of formulas $\Gamma$ of $\bI$, we define a relation "$\equiv_{\Gamma}^{\bI}$"\label{equivbI} between formulas of $\bI$ such that
\[\alpha\equiv_{\Gamma}^{\bI}\beta\quad\text{if and only if}\quad\Gamma\vdash_{\bI}\alpha\rightarrow\beta\quad\text{and}\quad\Gamma\vdash_{\bI}\beta\rightarrow\alpha.\]

\begin{proposition}
For any $\Gamma$, $\equiv_{\Gamma}^{\bI}$ is an equivalence relation.
\end{proposition}

\begin{proof}
\begin{enumerate}
\item For any formula $\alpha$, we trivially have that $\alpha\vdash_{\bI}\alpha$, and therefore $\vdash_{\bI}\alpha\rightarrow\alpha$; since $\bI$ is tarskian, we find $\Gamma\vdash_{\bI}\alpha\rightarrow\alpha$, and therefore $\alpha\equiv_{\Gamma}^{\bI}\alpha$.

\item If $\alpha\equiv_{\Gamma}^{\bI}\beta$, then $\Gamma\vdash_{\bI}\alpha\rightarrow\beta$ and $\Gamma\vdash_{\bI}\beta\rightarrow\alpha$; therefore $\Gamma\vdash_{\bI}\beta\rightarrow\alpha$ and $\Gamma\vdash_{\bI}\alpha\rightarrow\beta$ or, what is equivalent, $\beta\equiv_{\Gamma}^{\bI}\alpha$.

\item If $\alpha\equiv_{\Gamma}^{\bI}\beta$ and $\beta\equiv_{\Gamma}^{\bI}\gamma$, we have that $\Gamma$ demonstrates $\alpha\rightarrow\beta$, $\beta\rightarrow\alpha$, $\beta\rightarrow\gamma$ and $\gamma\rightarrow\beta$ in $\bI$; from $\Gamma\vdash_{\bI}\alpha\rightarrow\beta$ we get that $\Gamma, \alpha\vdash_{\bI}\beta$, and then from $\Gamma\vdash_{\bI}\beta\rightarrow\gamma$ and Modus Ponens we obtain that $\Gamma, \alpha\vdash_{\bI}\gamma$, or what is equivalent, that $\Gamma_{\bI}\alpha\rightarrow\gamma$. From $\Gamma\vdash_{\bI}\beta\rightarrow\alpha$ and $\Gamma\vdash_{\bI}\gamma\rightarrow\beta$ we find that $\Gamma\vdash_{\bI}\gamma\rightarrow\alpha$, meaning $\alpha\equiv_{\Gamma}^{\bI}\gamma$.
\end{enumerate}
\end{proof}

\begin{theorem}
For any $\Gamma$, $\equiv_{\Gamma}^{\bI}$ is a congruence with relation to the operators in $\{\vee, \wedge, \rightarrow\}$, meaning that for any operator $\#$ in this set, and any formulas $\alpha_{1}$, $\alpha_{2}$, $\beta_{1}$ and $\beta_{2}$ in $\bI$ such that $\alpha_{1}\equiv_{\Gamma}^{\bI}\beta_{1}$ and $\alpha_{2}\equiv_{\Gamma}^{\bI}\beta_{2}$,
\[\alpha_{1}\#\beta_{1}\equiv_{\Gamma}^{\bI}\alpha_{2}\#\beta_{2}.\]
\end{theorem}

\begin{proof}
Suppose $\alpha\equiv_{\Gamma}^{\bI}\beta$ and $\varphi\equiv_{\Gamma}^{\bI}\psi$, meaning $\Gamma$ proves $\alpha\rightarrow\beta$, $\beta\rightarrow\alpha$, $\varphi\rightarrow\psi$ and $\psi\rightarrow\varphi$. 

\begin{enumerate}
\item We have that $\Gamma, \alpha\vdash_{\bI}\beta$, and by $\textbf{Ax\: 6}$ we find $\Gamma, \alpha\vdash_{\bI}\beta\vee\psi$; analogously, we use $\Gamma, \varphi\vdash_{\bI}\psi$ to show that $\Gamma, \varphi\vdash_{\bI}\beta\vee\psi$, meaning through a proof by cases that
\[\Gamma, \alpha\vee\varphi\vdash_{\bI}\beta\vee\psi\]
and that, therefore, $\Gamma\vdash_{\bI}(\alpha\vee\varphi)\rightarrow(\beta\vee\psi)$. Using that $\Gamma$ proves $\beta\rightarrow\alpha$ and $\psi\rightarrow\varphi$, we get that $\Gamma\vdash_{\bI}(\beta\vee\psi)\rightarrow(\alpha\vee\varphi)$, and therefore $\alpha\vee\varphi\equiv_{\Gamma}^{\bI}\beta\vee\psi$.

\item From $\Gamma, \alpha\vdash_{\bI}\beta$ and $\Gamma, \varphi\vdash_{\bI}\psi$, the instance $\beta\rightarrow(\psi\rightarrow(\beta\wedge \psi))$ of $\textbf{Ax\: 3}$ and two applications of Modus Ponens, we obtain that $\Gamma, \alpha, \varphi\vdash_{\bI}\beta\wedge\psi$.

Since $\alpha\wedge\varphi\vdash_{\bI}\alpha, \varphi$, and by the fact that $\bI$ is tarskian, we obtain that $\Gamma, \alpha\wedge\varphi\vdash_{\bI}\beta\wedge\psi$, that is, $\Gamma\vdash_{\bI}(\alpha\wedge\varphi)\rightarrow(\beta\wedge\psi)$. Using $\Gamma\vdash_{\bI}\beta\rightarrow\alpha$ and $\Gamma\vdash_{\bI}\psi\rightarrow\varphi$, we obtain analogously $\Gamma\vdash_{\bI}(\beta\wedge\psi)\rightarrow(\alpha\wedge\varphi)$, and that implies $\alpha\wedge\varphi\equiv_{\Gamma}^{\bI}\beta\wedge\psi$.

\item From the fact that $\Gamma\vdash_{\bI}\beta\rightarrow\alpha$ and $\Gamma\vdash_{\bI}\varphi\rightarrow\psi$, applying twice that implication is transitive we find
\[\Gamma, \alpha\rightarrow\varphi\vdash_{\bI}\beta\rightarrow\psi,\]
and so $\Gamma\vdash_{\bI}(\alpha\rightarrow\varphi)\rightarrow(\beta\rightarrow\psi)$. From the fact that $\Gamma\vdash_{\bI}\alpha\rightarrow\beta, \psi\rightarrow\varphi$ we obtain, in a similar way, that $\Gamma\vdash_{\bI}(\beta\rightarrow\psi)\rightarrow(\alpha\rightarrow\varphi)$, meaning $\alpha\rightarrow\varphi\equiv_{\Gamma}^{\bI}\beta\rightarrow\psi$.

\end{enumerate}
\end{proof}

So the quotient set $A^{\bI}_{\Gamma}=F(\Sigma_{\bI}, \mathcal{V})/\equiv_{\Gamma}^{\bI}$\label{AbIGamma} is a well defined set, where we will denote the quotient class of the formula $\alpha$ by $[\alpha]$. More than that, from the fact that $\equiv_{\Gamma}^{\bI}$ is a congruence with respect to the connectives in $\{\vee, \wedge, \rightarrow\}$ we get that, by making
\[[\alpha]\#_{\mathcal{A}}[\beta]=[\alpha\#\beta],\quad\text{for}\quad\#\in\{\vee, \wedge, \rightarrow\},\]
such operations are well-defined. To see that, suppose both $\alpha$ and $\varphi$ are in the class $[\alpha]$ (meaning $\alpha\equiv_{\Gamma}^{\bI}\varphi$) and both $\beta$ and $\psi$ are in $[\beta]$: we get that $[\alpha]\#_{\mathcal{A}}[\beta]=[\alpha\#\beta]$ and $[\varphi]\#_{\mathcal{A}}[\psi]=[\varphi\#\psi]$, which equals $[\alpha\#\beta]$ since $\equiv_{\Gamma}^{\bI}$ is a congruence for $\#\in\{\vee, \wedge, \rightarrow\}$.

We then define, as it would be expected from our commentaries on the definable bottoms, tops and classical negation in $\bI$, $\bot_{\mathcal{A}}=[\bot_{\alpha\beta}]$ and $\top_{\mathcal{A}}=[\top_{\alpha}]$ for any formulas $\alpha$ in $\beta$ in $\bI$, what is possible since those bottom and top elements are equivalent to one another, and 
\[{\sim}_{\mathcal{A}}([\alpha])=[\alpha\rightarrow\bot_{\alpha\alpha}].\]
To see that the classical negation is well-defined, take $\beta\in[\alpha]$, and we wish to show that ${\sim}_{\mathcal{A}}([\alpha])={\sim}_{\mathcal{A}}([\beta])$, meaning $\alpha\rightarrow\bot_{\alpha\alpha}\equiv^{\bI}_{\Gamma}\beta\rightarrow\bot_{\beta\beta}$. This is equivalent to showing that, by applying the deduction meta-theorem,
\[\Gamma, \alpha\rightarrow\bot_{\alpha\alpha}\vdash_{\bI}\beta\rightarrow\bot_{\beta\beta}\quad\text{and}\quad\Gamma, \beta\rightarrow\bot_{\beta\beta}\vdash_{\bI}\alpha\rightarrow\bot_{\alpha\alpha}.\]
But since $\beta\in[\alpha]$, we have $\Gamma\vdash_{\bI}\alpha\rightarrow\beta$ and $\Gamma\vdash_{\bI}\beta\rightarrow\alpha$, while $\bot_{\alpha\alpha}\rightarrow\bot_{\beta\beta}$ and $\bot_{\beta\beta}\rightarrow\bot_{\alpha\alpha}$ are tautologies, which implies the desired result through the transitive property of implication.

\begin{theorem}
$\mathcal{A}=(A^{\bI}_{\Gamma}, \{\sigma_{\mathcal{A}}\}_{\sigma\in\Sigma^{\textbf{CPL}}})$ is a Boolean algebra.
\end{theorem}

\begin{proof}
For simplicity, we will drop the index $\mathcal{A}$. So, it must be proven that, for any formulas $\alpha$, $\beta$ and $\gamma$ in $\bI$,
\begin{enumerate}
\item $[\alpha]\vee\bot=[\alpha]$ and $[\alpha]\wedge\top=[\alpha]$;
\item $[\alpha]\vee[\beta]=[\beta]\vee[\alpha]$ and $[\alpha]\wedge[\beta]=[\beta]\wedge[\alpha]$;
\item $[\alpha]\vee([\beta]\wedge[\gamma])=([\alpha]\vee[\beta])\wedge([\alpha]\vee[\gamma])$ and $[\alpha]\wedge([\beta]\vee[\gamma])=([\alpha]\wedge[\beta])\vee([\alpha]\wedge[\gamma])$;
\item $[\alpha]\vee{\sim} [\alpha]=\top$ and $[\alpha]\wedge{\sim} [\alpha]=\bot$;
\item $[\alpha]\rightarrow[\beta]={\sim} [\alpha]\vee[\beta]$.
\end{enumerate}
For those items containing more than one assertion, we shall only prove the first, being the second completely analogous.
\begin{enumerate}
\item We must prove that $\alpha\vee\bot_{\alpha\alpha}\equiv_{\Gamma}^{\bI}\alpha$. Clearly $\Gamma\vdash_{\bI}\alpha\rightarrow\alpha\vee\bot_{\alpha\alpha}$, since this formula is an instance of $\textbf{Ax\: 6}$.

We also know that $\bot_{\alpha\alpha}$ behaves like a bottom element, i.e., for any formula $\beta$, $\bot_{\alpha\alpha}\vdash_{\bI}\beta$, and therefore $\Gamma, \bot_{\alpha\alpha}\vdash_{\bI}\alpha$; since $\alpha\rightarrow\alpha$ is a tautology, $\Gamma, \alpha\vdash_{\bI}\alpha$, and therefore $\Gamma, \alpha\vee\bot_{\alpha\alpha}\vdash_{\bI}\alpha$, implying $\Gamma\vdash_{\bI}\alpha\vee\bot_{\alpha\alpha}\rightarrow\alpha$.

\item The statement is equivalent to $\alpha\vee\beta\equiv_{\Gamma}^{\bI}\beta\vee\alpha$: from the instance $\beta\rightarrow(\alpha\vee\beta)$ of $\textbf{Ax\: 7}$ and the instance $\alpha\rightarrow(\alpha\vee\beta)$ of $\textbf{Ax\: 6}$, we get that $\Gamma,\beta\vdash_{\bI}\alpha\vee\beta$ and $\Gamma, \alpha\vdash_{\bI}\alpha\vee\beta$, and therefore $\Gamma, \beta\vee\alpha\vdash_{\bI}\alpha\vee\beta$, meaning $\Gamma\vdash_{\bI}(\beta\vee\alpha)\rightarrow(\alpha\vee\beta)$.

Using the instances of axioms $\alpha\rightarrow(\beta\vee\alpha)$ and $\beta\rightarrow(\beta\vee\alpha)$ we get the reciprocal.

\item We must prove that $\alpha\vee(\beta\wedge\gamma)\equiv_{\Gamma}^{\bI}(\alpha\vee\beta)\wedge(\alpha\vee\gamma)$: $(\beta\wedge\gamma)\rightarrow\beta$ and $(\beta\wedge\gamma)\rightarrow\gamma$ are instances of $\textbf{Ax\: 4}$ and $\textbf{Ax\: 5}$, and $\beta\rightarrow(\alpha\vee\beta)$ and $\gamma\rightarrow(\alpha\vee\gamma)$ are instances of $\textbf{Ax\: 7}$, meaning $\Gamma, (\beta\wedge\gamma)\vdash_{\bI}\alpha\vee\beta$ and $\Gamma, (\beta\wedge\gamma)\vdash_{\bI}\alpha\vee\gamma$.

Since $\alpha\rightarrow(\alpha\vee\beta)$ and $\alpha\rightarrow(\alpha\vee\gamma)$ are instances of $\textbf{Ax\: 6}$, $\Gamma, \alpha\vee(\beta\wedge\gamma)\vdash_{\bI}\alpha\vee\beta$ and $\Gamma, \alpha\vee(\beta\wedge\gamma)\vdash_{\bI}\alpha\vee\gamma$, and from the instance
\[\alpha\vee\beta\rightarrow(\alpha\vee\gamma\rightarrow((\alpha\vee\beta)\wedge(\alpha\vee\gamma)))\]
of $\textbf{Ax\: 3}$ and Modus Ponens, we get that $\Gamma\vdash_{\bI}[\alpha\vee(\beta\wedge\gamma)]\rightarrow[(\alpha\vee\beta)\wedge(\alpha\vee\gamma)]$.

Reciprocally, is enough to show, by $\textbf{Ax\: 4}$ and $\textbf{Ax\: 5}$, that $\Gamma, \alpha\vee\beta, \alpha\vee\gamma\vdash_{\bI}\alpha\vee(\beta\wedge\gamma)$. Then we have four cases:
\begin{enumerate}
\item $\Gamma, \alpha\vdash_{\bI}\alpha\vee(\beta\wedge\gamma)$, from the instance $\alpha\rightarrow[\alpha\vee(\beta\wedge\gamma)]$ of $\textbf{Ax\: 6}$;
\item $\Gamma, \alpha, \beta\vdash_{\bI}\alpha\vee(\beta\wedge\gamma)$ for the same reason as above;
\item $\Gamma, \alpha, \gamma\vdash_{\bI}\alpha\vee(\beta\wedge\gamma)$, again for the same reason;
\item $\Gamma, \beta, \gamma\vdash_{\bI}\alpha\vee(\beta\wedge\gamma)$, from the instance $\beta\rightarrow(\gamma\rightarrow(\beta\wedge\gamma))$ of $\textbf{Ax\: 3}$ and two applications of Modus Ponens that give us $\Gamma, \beta, \gamma\vdash_{\bI}\beta\wedge\gamma$, and then the use of $\textbf{Ax\: 7}$.
\end{enumerate}
So, from the first two items, $\Gamma, \alpha, \alpha\vee\beta\vdash_{\bI}\alpha\vee(\beta\wedge\gamma)$, and from the last two $\Gamma, \alpha\vee\beta, \gamma\vdash_{\bI}\alpha\vee(\beta\wedge\gamma)$, meaning
\[\Gamma, \alpha\vee\beta, \alpha\vee\gamma\vdash_{\bI}\alpha\vee(\beta\wedge\gamma).\]

\item To prove $\alpha\vee(\alpha\rightarrow\bot_{\alpha\alpha})\equiv_{\Gamma}^{\bI}\top_{\alpha}$, we notice first that $[\alpha\vee(\alpha\rightarrow\bot_{\alpha\alpha})]\rightarrow\top_{\alpha}$ follows from the fact that $\top_{\alpha}$ behaves like a top element, meaning $\Gamma\vdash_{\bI}[\alpha\vee(\alpha\rightarrow\bot_{\alpha\alpha})]\rightarrow\top_{\alpha}$.

Reciprocally, $\alpha\vee(\alpha\rightarrow\bot_{\alpha\alpha})$ is an instance of an axiom, and for this reason $\vdash_{\bI}\alpha\vee(\alpha\rightarrow\bot_{\alpha\alpha})$, meaning $\Gamma\vdash_{\bI}\top_{\alpha}\rightarrow[\alpha\vee(\alpha\rightarrow\bot_{\alpha\alpha})]$.

\item Finally, we need now only to prove that $\alpha\rightarrow\beta\equiv_{\Gamma}^{\bI}(\alpha\rightarrow\bot_{\alpha\alpha})\vee\beta$, that is, by applying the deduction meta-theorem, $\Gamma, \alpha\rightarrow\beta\vdash_{\bI}(\alpha\rightarrow\bot_{\alpha\alpha})\vee\beta$ and $\Gamma, (\alpha\rightarrow\bot_{\alpha\alpha})\vee\beta\vdash_{\bI}\alpha\rightarrow\beta$. 

Clearly $\alpha, \alpha\rightarrow\beta\vdash_{\bI}(\alpha\rightarrow\bot_{\alpha\alpha})\vee\beta$, by using Modus Ponens and the instance $\beta\rightarrow(\alpha\rightarrow\bot_{\alpha\alpha})\vee\beta$ of $\textbf{Ax\: 7}$; meanwhile, $\alpha\rightarrow\bot_{\alpha\alpha}, \alpha\rightarrow\beta\vdash_{\bI}(\alpha\rightarrow\bot_{\alpha\alpha})\vee\beta$ by the instance 
\[(\alpha\rightarrow\bot_{\alpha\alpha})\rightarrow(\alpha\rightarrow\bot_{\alpha\alpha})\vee\beta\]
of $\textbf{Ax\: 6}$. From a proof by cases, $(\alpha\rightarrow\bot_{\alpha\alpha})\vee\alpha, \alpha\rightarrow\beta\vdash_{\bI}(\alpha\rightarrow\bot_{\alpha\alpha})\vee\beta$, and since $(\alpha\rightarrow\bot_{\alpha\alpha})\vee\alpha$ is an instance of $\textbf{Ax\: 9}^{*}$ and $\bI$ is tarskian, we obtain the desired result.

Reciprocally, $\beta\vdash_{\bI}\alpha\rightarrow\beta$ by $\textbf{Ax\: 1}$, and since $\bot_{\alpha\alpha}\rightarrow\beta$ is a tautology, $\alpha\rightarrow\bot_{\alpha\alpha}\vdash_{\bI}\alpha\rightarrow\beta$, by use of the transitivity of implication. By a proof by cases, $(\alpha\rightarrow\bot_{\alpha\alpha})\vee\beta\vdash_{\bI}\alpha\rightarrow\beta$.
\end{enumerate}

\end{proof}

Now, for a given $\Gamma$, we can make $\mathcal{A}=(A^{\bI}_{\Gamma}, \{\sigma_{\mathcal{A}}\}_{\sigma\in\Sigma^{\textbf{CPL}}})$ into a Fidel structure, presented as a $\Sigma_{\bI}^{\textbf{CPL}}$-multialgebra, for $\bI$: we simply define
\[[\alpha]\uparrow[\beta]=\{[\varphi\uparrow\psi]\ :\  \varphi\in [\alpha], \psi\in[\beta]\};\]
for $[\gamma]\in[\alpha]\uparrow[\beta]$, we must show that $[\alpha]\wedge([\beta]\wedge[\gamma])=\bot$; taking $\varphi\in[\alpha]$ and $\psi\in[\beta]$ such that $\gamma\in[\varphi\uparrow\psi]$, this means
\[\varphi\wedge(\psi\wedge\gamma)\equiv_{\Gamma}^{\bI}\bot_{\varphi\psi}.\]
Clearly $\Gamma\vdash_{\bI}\bot_{\varphi\psi}\rightarrow[\varphi\wedge(\psi\wedge\gamma)]$, so it remains to be shown that $\Gamma\vdash_{\bI}[\varphi\wedge(\psi\wedge\gamma)]\rightarrow\bot_{\varphi\psi}$. This is equivalent to $\Gamma, \varphi, \psi, \gamma\vdash_{\bI}\bot_{\varphi\psi}$, from $\textbf{Ax\: 4}$, and since $\Gamma\vdash_{\bI}\gamma\rightarrow(\varphi\uparrow\psi)$, $\Gamma\vdash_{\bI}(\varphi\uparrow\psi)\rightarrow\gamma$ and $\bot_{\varphi\psi}=\varphi\wedge(\psi\wedge(\varphi\uparrow\psi))$, by Modus Ponens this equals
\[\Gamma, \varphi, \psi, \varphi\uparrow\psi\vdash_{\bI}\varphi\wedge(\psi\wedge(\varphi\uparrow\psi)),\]
but this is obviously true. We can also prove that if $\varphi\in[\alpha]$ and $\psi\in[\beta]$, $[\psi\uparrow\varphi]\in[\alpha]\uparrow[\beta]$. Very clearly $[\varphi\uparrow\psi]\in [\alpha]\uparrow[\beta]$, so it is enough to prove that $[\psi\uparrow\varphi]=[\varphi\uparrow\psi]$, or in other terms, that $\psi\uparrow\varphi\equiv_{\Gamma}^{\bI}\varphi\uparrow\psi$, meaning
\[\Gamma\vdash_{\bI}(\psi\uparrow\varphi)\rightarrow(\varphi\uparrow\psi)\quad\text{and}\quad\Gamma\vdash_{\bI}(\varphi\uparrow\psi)\rightarrow(\psi\uparrow\varphi);\]
and this is clearly true, since both implications are instances of $\textbf{Comm}$. This shows that $[\alpha]\uparrow[\beta]=[\beta]\uparrow[\alpha]$, since the elements of both are the same.

The $\Sigma_{\bI}^{\textbf{CPL}}$-multialgebra $\mathcal{A}^{\bI}_{\Gamma}$\label{mAbIGamma} just described, with universe $A^{\bI}_{\Gamma}$, is called a Lindenbaum-Tarski multialgebra of $\bI$, specifically the one associated to $\Gamma$.

\begin{theorem}
Given formulas $\Gamma\cup\{\varphi\}$ of $\bI$, if $\Gamma\Vdash^{\bI}_{\mathcal{F}}\varphi$ then $\Gamma\vdash_{\bI}\varphi$.
\end{theorem}

\begin{proof}
Suppose that $\Gamma\not\vdash_{\bI}\varphi$, and we know there exists a closed set of formulas $\Delta$ of $\bI$ that contains $\Gamma$ and is maximal with respect to not proving $\varphi$. We take the Lindenbaum-Tarski multialgebra $\mathcal{A}^{\bI}_{\Delta}$ of $\bI$ associated to $\Delta$ and its associated RNmatrix, and consider the map $\nu:F(\Sigma_{\bI}, \mathcal{V})\rightarrow A^{\bI}_{\Delta}$ such that $\nu(\alpha)=[\alpha]$. We have that:
\begin{enumerate}
\item for any $\#\in\{\vee, \wedge, \rightarrow\}$, $\nu(\alpha\#\beta)=[\alpha\#\beta]=[\alpha]\#[\beta]=\nu(\alpha)\#\nu(\beta)$;
\item $\nu(\alpha\uparrow\beta)=[\alpha\uparrow\beta]$, which is in $\{[\varphi\uparrow\psi]\ :\  \varphi\in [\alpha], \psi\in[\beta]\}=[\alpha]\uparrow[\beta]$;
\item $\nu(\alpha\uparrow\beta)=[\alpha\uparrow\beta]=[\beta\uparrow\alpha]=\nu(\beta\uparrow\alpha)$.
\end{enumerate}
So $\nu$ is a valuation over $\mathcal{A}^{\bI}_{\Delta}$ in $\mathcal{F}_{\mathcal{A}^{\bI}_{\Delta}}$. Furthermore, $\nu(\alpha)=\top$ if and only if $\Delta\vdash_{\bI}\alpha$; given that $\Delta\not\vdash_{\bI}\varphi$, we obtain $\nu(\Delta)\subseteq \{\top\}$ and $\nu(\varphi)\neq\top$, meaning that $\Delta\not\Vdash^{\bI}_{\mathcal{F}}\varphi$ and therefore $\Gamma\not\Vdash^{\bI}_{\mathcal{F}}\varphi$.
\end{proof}


\subsection{Decision method}\label{Decision Method for bI}

First of all, we remember the set $\{0,1\}$ can be made into a Boolean algebra $\textbf{2}$ when we define $\top=1$, $\bot=0$, and the other operations according to the following tables.

\begin{figure}[H]
\centering
\begin{minipage}[t]{2cm}
\centering
\begin{tabular}{|l|r|}
\hline
& ${\sim} $\\ \hline
$0$ & $1$\\ \hline
$1$ & $0$\\\hline
\end{tabular}
\caption*{Negation}
\end{minipage}
\hspace{1cm}
\centering
\begin{minipage}[t]{2.5cm}
\centering
\begin{tabular}{|l|c|r|}
\hline
$\vee$ & $0$ & $1$\\ \hline
$0$ & $0$ & $1$\\ \hline
$1$ & $1$ & $1$\\\hline
\end{tabular}
\caption*{Disjunction}
\end{minipage}
\hspace{1cm}
\centering
\begin{minipage}[t]{2.5cm}
\centering
\begin{tabular}{|l|c|r|}
\hline
$\wedge$ & $0$ & $1$\\ \hline
         $0$ & $0$ & $0$\\ \hline
         $1$ & $0$ & $1$\\\hline
\end{tabular}
\caption*{Conjunction}
\end{minipage}
\hspace{1cm}
\centering
\begin{minipage}[t]{2.5cm}
\centering
\begin{tabular}{|l|c|r|}
\hline
$\rightarrow$ & $0$ & $1$\\ \hline
             $0$ & $1$ & $1$\\ \hline
             $1$ & $0$ & $1$\\\hline
\end{tabular}
\caption*{Implication}
\end{minipage}
\end{figure}

Notice that, to give one example, technically $\bot=\{0\}$, since we are in the environment of multialgebras; however, given the operations associated to the elements in $\Sigma^{\textbf{CPL}}$ are supposed to be single-valued, since $(A, \{\sigma_{\mathcal{A}}\}_{\sigma\in\Sigma^{\textbf{CPL}}})$ is a Boolean algebra, we will identify a singleton to the element it contains, as is usual.

We make $\textbf{2}$ into a Fidel structure $\textbf{2}_{\bI}$\label{2bI}, presented as a $\Sigma_{\bI}^{\textbf{CPL}}$-multialgebra, for $\bI$ by defining $\uparrow_{\textbf{2}_{\bI}}(1,1)=\{0\}$ and $\uparrow_{\textbf{2}_{\bI}}(x,y)=\{0,1\}$ otherwise. For simplicity, we will denote $\uparrow_{\textbf{2}_{\bI}}$ simply by $\uparrow$ and use the infix notation.

\begin{figure}[H]
\centering
\captionsetup{justification=centering}
\begin{tabular}{|l|c|r|}
\hline
$\uparrow$ & $0$ & $1$\\ \hline
$0$ & $\{0,1\}$ & $\{0,1\}$\\ \hline
$1$ & $\{0,1\}$ & $\{0\}$\\\hline
\end{tabular}
\caption*{Table for $\uparrow$ in $\textbf{2}_{\bI}$}
\end{figure}

This is very clearly a Fidel structure for $\bI$ since, for $z\in x\uparrow y$, we have that $x\wedge(y\wedge z)=0$: in the case where neither $x$ nor $y$ equals $0$, the only possible value for $z$ is exactly $0$. Furthermore, one always has $x\uparrow y=y\uparrow x$.

\begin{theorem}
$\nu:F(\Sigma_{\bI}, \mathcal{V})\rightarrow\{0,1\}$ is a bivaluation for $\bI$ if, and only if, it is a $\Sigma_{\bI}$-homomorphism from $\textbf{F}(\Sigma_{\bI}, \mathcal{V})$ to $\textbf{2}_{\bI}$ in $\mathcal{F}_{\textbf{2}_{\bI}}$.
\end{theorem}

\begin{proof}
Assume first that $\nu$ is a bivaluation. Then:
\begin{enumerate}
\item if $\nu(\alpha)=1$ or $\nu(\beta)=1$, we have that $\nu(\alpha)\vee\nu(\beta)=1=\nu(\alpha\vee\beta)$; otherwise, when both $\nu(\alpha)$ and $\nu(\beta)$ are $0$, we have that $\nu(\alpha)\vee\nu(\beta)=0=\nu(\alpha\vee\beta)$;
\item if both $\nu(\alpha)$ and $\nu(\beta)$ are $1$, $\nu(\alpha)\wedge\nu(\beta)=1=\nu(\alpha\wedge\beta)$; otherwise, $\nu(\alpha)\wedge\nu(\beta)=0$, which is also the value of $\nu(\alpha\wedge\beta)$;
\item if either $\nu(\alpha)=0$ or $\nu(\beta)=1$, $\nu(\alpha)\rightarrow\nu(\beta)=1=\nu(\alpha\rightarrow\beta)$; in the case that $\nu(\alpha)=1$ and $\nu(\beta)=0$, $\nu(\alpha)\rightarrow\nu(\beta)=0=\nu(\alpha\rightarrow\beta)$;
\item for any value of $\nu(\alpha\uparrow\beta)$ such that either $\nu(\alpha)$ or $\nu(\beta)$ is not $1$, it is clear $\nu(\alpha\uparrow\beta)\in\nu(\alpha)\uparrow\nu(\beta)$, since in those cases $\nu(\alpha)\uparrow\nu(\beta)$ is the whole universe of $\textbf{2}_{\bI}$, so assume $\nu(\alpha)=\nu(\beta)=1$; if we had $\nu(\alpha\uparrow\beta)=1$, by the fact $\nu(\alpha\uparrow\beta)=1$ and $\nu(\alpha)=1$ we would get that $\nu(\beta)=0$, which is a contradiction, forcing us to have
\[\nu(\alpha\uparrow\beta)=0\in\{0\}=1\uparrow1=\nu(\alpha)\uparrow\nu(\beta);\]
\item since $\nu(\alpha\uparrow\beta)=\nu(\beta\uparrow\alpha)$ and, from the previous items, $\nu$ is a $\Sigma_{\bI}$-homomorphism, we have that $\nu$ lies in $\mathcal{F}_{\textbf{2}_{\bI}}$
\end{enumerate}
To summarize, $\nu$ is a $\Sigma_{\bI}$-homomorphism in $\mathcal{F}_{\textbf{2}_{\bI}}$.

Reciprocally, assume $\nu$ is a $\Sigma_{\bI}$-homomorphism which lies in $\mathcal{F}_{\textbf{2}_{\bI}}$.
\begin{enumerate}
\item since $\nu(\alpha\vee\beta)=\nu(\alpha)\vee\nu(\beta)$, looking at the table for disjunction we find that $\nu(\alpha\vee\beta)=1$ if and only if $\nu(\alpha)=1$ or $\nu(\beta)=1$;
\item since $\nu(\alpha\wedge\beta)=\nu(\alpha)\wedge\nu(\beta)$, $\nu(\alpha\wedge\beta)=1$ if and only if $\nu(\alpha)=\nu(\beta)=1$;
\item given $\nu(\alpha\rightarrow\beta)=\nu(\alpha)\rightarrow\nu(\beta)$, $\nu(\alpha\rightarrow\beta)=1$ if and only if $\nu(\alpha)=0$ or $\nu(\beta)=1$;
\item if $\nu(\alpha\uparrow\beta)=1$, from the table for $\uparrow$ there are three possible cases: $\nu(\alpha)=\nu(\beta)=0$, $\nu(\alpha)=0$ and $\nu(\beta)=1$, and $\nu(\alpha)=1$ and $\nu(\beta)=0$; so if $\nu(\alpha\uparrow\beta)=1$ and $\nu(\alpha)=1$, we have that $\nu(\beta)=0$;
\item since $\nu\in\mathcal{F}_{\textbf{2}_{\bI}}$, $\nu(\alpha\uparrow\beta)=\nu(\beta\uparrow\alpha)$.
\end{enumerate}
This finishes proving that $\nu$ is then a bivaluation.
\end{proof}

We will denote the restricted Nmatrix $(\textbf{2}_{\bI}, \{1\}, \mathcal{F}_{\textbf{2}_{\bI}})$ by $\mathbb{2}_{\bI}$.

\begin{theorem}\label{Decision method for BI proof}
Given formulas $\Gamma\cup\{\varphi\}$ of $\bI$, $\Gamma\vDash_{\bI}\varphi$ if and only if $\Gamma\vDash_{\mathbb{2}_{\bI}}\varphi$.
\end{theorem}

\begin{proof}
The proof is quite straightforward: suppose $\Gamma\vDash_{\bI}\varphi$; then, for any $\Sigma_{\bI}$-homomorphism $\nu:\textbf{F}(\Sigma_{\bI}, \mathcal{V})\rightarrow\textbf{2}_{\bI}$ in $\mathcal{F}_{\textbf{2}_{\bI}}$ such that $\nu(\Gamma)\subseteq\{1\}$, from the previous theorem we have that $\nu$ is a bivaluation, and since $\nu(\Gamma)\subseteq\{1\}$ and $\Gamma\vDash_{\bI}\varphi$, we find that $\nu(\varphi)=1$. Therefore, $\Gamma\vDash_{\textbf{2}_{\bI}}\varphi$.

The reciprocal is analogous.
\end{proof}

As it was discussed in depth in Section \ref{Decision methods}, not all finite RNmatrices lead straightforwardly to decision methods by row-branching, row-eliminating truth tables, since there may not exist an algorithm to decide whether a given row is, or is not, a restricted valuation. As luck would have it, in the case of $\bI$ the RNmatrix $\mathbb{2}_{\bI}$ does indeed offer a decision method; we will, however, just briefly discuss how it works, the complete proof being quite similar to the one found in Section \ref{Row-branching,row-eliminating}.

We start as it is usual: for a given formula $\varphi$, produce a list of its subformulas in non-decreasing order of complexity $\varphi_{1}, \dotsc , \varphi_{n}=\varphi$; the complexity, also known as order, of a formula $\varphi$ of $\bI$ is defined as $0$ if $\varphi$ is a variable, and if $\alpha$ and $\beta$ have complexity $p$ and $q$, respectively, and $\varphi=\alpha\#\beta$ for a $\#\in\{\vee, \wedge, \rightarrow, \uparrow\}$, then the complexity of $\varphi$ is $\max\{p, q\}+1$.

We produce a table with columns headed by the formulas $\varphi_{i}$; we have some $1\leq m<n$ such that all $\varphi_{1}, \dotsc , \varphi_{m}$ are propositional variables, and so we write down $2^{m}$ initial rows, each with a possible combination of zeros and ones for these variables. Then we start filling in the rows for $\varphi_{m+1},\dotsc , \varphi_{n}$ in the following way: if $\varphi_{l}=\varphi_{i}\#\varphi_{j}$ (where necessarily $i, j< l$), $\#\in\{\vee, \wedge, \rightarrow\}$ and, on a certain row, $\varphi_{i}$ is given the value $a$ and $\varphi_{j}$ is given the value $b$, then on this very row $\varphi_{l}$ is given the value $a\#b$. If $\varphi_{l}=\varphi_{i}\uparrow\varphi_{j}$ and both $\varphi_{i}$ and $\varphi_{j}$ are $1$ on a row, then $\varphi_{l}$ is $0$; if, however, on a specific row $\varphi_{i}$ or $\varphi_{j}$ takes the value $0$, then this row is split in two, one where $\varphi_{l}$ assumes the value $0$, the other where it assumes the value $1$.

The caveat of this last part is that undesired valuations need to be erased: if there is a formula $\varphi_{k}=\varphi_{j}\uparrow\varphi_{i}$, with $i, j<k<l$, that takes the value $0$ in the present row, then the row where $\varphi_{l}=\varphi_{i}\uparrow\varphi_{l}$ takes the value $1$ is erased, in order to $\varphi_{i}\uparrow\varphi_{j}$ and $\varphi_{j}\uparrow\varphi_{i}$ be given the same value; analogously, if in this row $\varphi_{k}$ is $1$, the row where $\varphi_{l}$ is $0$ is eliminated.\footnote{Of course, one could write down the row-branching truth table for $\varphi$, corresponding to the Nmatrix subjacent to $\mathbb{2}_{\bI}$, without simultaneously erasing the undesired rows, and only at the end select those rows where $\varphi_{i}\uparrow\varphi_{j}$ and $\varphi_{j}\uparrow\varphi_{i}$ are given different values to eliminate; this may be less efficient, but the result is the same.} At the end, the column headed by $\varphi$ is filled with nothing but ones if, and only if, the formula is a tautology for $\bI$.

For a finite set $\Gamma=\{\gamma_{1}, \dotsc , \gamma_{n}\}$, we can test whether $\Gamma\vdash_{\bI}\varphi$ by simply testing if $\gamma_{1}\rightarrow(\cdots\rightarrow(\gamma_{n}\rightarrow\varphi)\cdots)$ or $\bigwedge_{i=1}^{n}\gamma_{i}\rightarrow\varphi$ is a tautology in the same logic; of course, we can be more efficient and simply merge the tables for $\varphi$ and all $\gamma_{i}$ and check if, in all rows where all elements of $\Gamma$ are designated, so is $\varphi$.

To give one example of how our row-branching, row-eliminating truth tables work, consider the deduction $p, p\uparrow q\vdash_{\bI}q\rightarrow (q\uparrow p)$; basically, it states that if $p$ is true and incompatible with $q$, then $q$ must not be true (notice that $q\rightarrow(q\uparrow p)$, in the presence of $p$, is equivalent to ${\sim}q$).

\begin{center}
\begin{tabular}{|c|c|c|c|c|}
\hline
$p$ & $q$ & $p\uparrow q$ & $q\uparrow p$ & $q\rightarrow(q\uparrow p)$\\ \hline
$1$ & $1$ & $0$ & $0$ & $0$\\\hline
\multirow{2}{*}{$1$} & \multirow{2}{*}{$0$} & $1$ & $1$ & $1$\\ \cline{3-5}
& & $0$ & $0$ & $1$\\\hline
\multirow{2}{*}{$0$} & \multirow{2}{*}{$1$} & $1$ & $1$ & $1$\\ \cline{3-5}
& & $0$ & $0$ & $0$\\\hline
\multirow{2}{*}{$0$} & \multirow{2}{*}{$0$} & $1$ & $1$ & $1$\\ \cline{3-5}
& & $0$ & $0$ & $1$\\\hline
\end{tabular}
\end{center}

Notice that, since in every row where both $p$ and $p\uparrow q$ are $1$ so is $q\rightarrow(q\uparrow p)$, the deduction is valid.

\begin{theorem}
Let $\mathbb{2}_{\textbf{bI}^{-}}$ be the Nmatrix $(\textbf{2}_{\bI}, \{1\})$; for any formulas $\Gamma\cup\{\varphi\}$ of $\textbf{bI}^{-}$, $\Gamma\vDash_{\textbf{bI}^{-}}\varphi$ if and only if $\Gamma\vDash_{\textbf{bI}^{-}}\varphi$.
\end{theorem}

This result, whose proof is very similar to the one of Theorem \ref{Decision method for BI proof}, shows a strange property of logics of incompatibility and their semantics which deserves to be further studied: while none of the logics $\textbf{bI}^{-}$ and $\textbf{CPL}$, respectively weaker and stronger than $\textbf{bI}$, needs RNmatrices to be characterized, $\textbf{bI}$, lying squarely in the middle, does indeed require RNmatrices, what is formally proved in Section \ref{bI is not characterized by Nmatrices section}.

\subsection{Another decision method}\label{Tableaux for bI}

In Section \ref{Tableaux semantics for Cn} we showed how it is sometimes possible to extract from the tables of the operations of a finite RNmatrix tableau rules that, when paired with adequate closure conditions for branches, produce a tableau semantics for the logic characterized by said RNmatrix. This can again be done in $\bI$, but the full proof that this works is left as an exercise: it follows the lines of the same result for $C_{n}$, being however much simpler.

\begin{definition}
The following are the tableau rules with labeled formulas, where $\varphi$ and $\psi$ are formulas of $\bI$, and $\textsf{0}$ and $\textsf{1}$ are all possible labels, for the tableau calculus we denote by $\mathbb{T}_{\bI}$.

$$
\begin{array}{cp{1.5cm}cp{1.5cm}c}
\displaystyle \frac{\textsf{0}(\varphi\vee\psi)}{\begin{array}{c}\textsf{0}(\varphi) \\ \textsf{0}(\psi)\end{array}} & & \displaystyle \frac{\textsf{0}(\varphi\wedge\psi)}{\textsf{0}(\varphi)\mid\textsf{0}(\psi)} & & \displaystyle \frac{\textsf{0}(\varphi\rightarrow\psi)}{\begin{array}{c}\textsf{1}(\varphi) \\ \textsf{0}(\psi)\end{array}}  \\[2mm]
&&&&\\[2mm]
 \displaystyle \frac{\textsf{1}(\varphi\vee\psi)}{\textsf{1}(\varphi)\mid\textsf{1}(\psi)} & & \displaystyle \frac{\textsf{1}(\varphi\wedge\psi)}{\begin{array}{c}\textsf{1}(\varphi) \\ \textsf{1}(\psi)\end{array}} & & \displaystyle \frac{\textsf{1}(\varphi\rightarrow\psi)}{\textsf{0}(\varphi)\mid\textsf{1}(\psi)}\\[2mm]
&&&&\\[2mm]
\end{array}
$$

\[\frac{\textsf{1}(\varphi\uparrow\psi)}{\textsf{0}(\varphi)\mid\textsf{0}(\psi)}\]

Branches of tableaux in $\mathbb{T}_{\bI}$ are closed iff:
\begin{enumerate}
\item they contain labeled formulas $\textsf{L}(\varphi)$ and $\textsf{L}^{*}(\varphi)$ with $\textsf{L}\neq\textsf{L}^{*}$;
\item or they contain labeled formulas $\textsf{L}(\varphi\uparrow\psi)$ and $\textsf{L}^{*}(\psi\uparrow\varphi)$ with $\textsf{L}\neq\textsf{L}^{*}$.
\end{enumerate}

A branch $\theta$ is complete whenever it contains, for each $\textsf{L}(\gamma)$ in $\theta$ not of the form $\textsf{0}(\varphi\uparrow\psi)$ and with $\gamma$ not a variable, all the labeled formulas of one of the branches in the rule headed by $\textsf{L}(\gamma)$; complete branches are open if they are not closed.

A tableau in $\mathbb{T}_{\bI}$ is: closed if all of its branches are closed; complete if all of its branches are either complete or closed; and open if it is complete but not closed.
\end{definition}

Of course, the labels $\textsf{0}$ and $\textsf{1}$ correspond, respectively, to the values $0$ and $1$ of $\mathbb{2}_{\bI}$; notice, furthermore, that if we had included in the signature of $\bI$ the classical negation that is definable in this logic, it wouldn't be necessary to use labels, but this is in no way mandatory for our tableaux to properly work.

Notice that all the rules found in $\mathbb{T}_{\bI}$ are analytic, in the sense that the complexities of the formulas obtained from applying a rule are strictly smaller than the complexity of the formula that motivated the application of the rule. This means that all tableaux in $\mathbb{T}_{\bI}$ can be completed in a finite number of steps. A formula $\varphi$ of $\bI$ is provable according to tableaux in $\mathbb{T}_{\bI}$, when we write $\vdash_{\mathbb{T}_{\bI}}\varphi$, if there is a closed tableau in $\mathbb{T}_{\bI}$ starting from $\textsf{0}(\varphi)$; for a finite set of formulas $\Gamma=\{\gamma_{1}, \dotsc , \gamma_{n}\}$ of $\bI$, we say $\Gamma$ proves $\varphi$ according to $\mathbb{T}_{\bI}$ if 
\[\vdash_{\mathbb{T}_{\bI}}\bigwedge_{i=1}^{n}\gamma_{n}\rightarrow\varphi,\]
or alternatively 
\[\vdash_{\mathbb{T}_{\bI}}\gamma_{1}\rightarrow (\gamma_{2}\rightarrow\cdots (\gamma_{n}\rightarrow\varphi)\cdots).\]
The following theorem is proved as in the case of da Costa's hierarchy, found in Section \ref{Tableaux semantics for Cn}, and shows we indeed have a decision method.

\begin{theorem}
For any finite set of formulas $\Gamma\cup\{\varphi\}$ of $\bI$, $\Gamma\vdash_{\bI}\varphi$ if, and only if, $\Gamma\vdash_{\mathbb{T}_{\bI}}\varphi$.
\end{theorem}


\section{Other logics}

Of course, when exploring a new logical environment, there are many possibilities that arise according to your philosophical interpretation of what the logic at hand should accomplish.

The most basic logic of formal incompatibility, $\bI$, is just that, the most basic one: there are many, very natural, different logics of formal incompatibility that occur when dealing with this subject, from which we present a few in the next sections.


\subsection{The logic $\bIpr$}

The more one analyzes the concept of incompatibility, the more it seems it should have some sort of interplay, of relationship with other logical objects such as disjunction and conjunction. In the logic we shall call $\bIpr$\label{bIpr} we add axioms relating $\uparrow$ with $\vee$ and $\wedge$.

It seems that, in natural language, when $\gamma$ is incompatible with $\alpha$ and incompatible with $\beta$, then it is incompatible with the disjunction of $\alpha$ and $\beta$ ; alternatively, when $\gamma$ is incompatible with $\alpha$ or incompatible with $\beta$, then it is incompatible with the conjunction of $\alpha$ and $\beta$.

To see how $\textbf{CPL}$ would approach this, we remember the usual interpretation of the Sheffer's stroke is $\alpha\Uparrow\beta={\sim}(\alpha\wedge\beta)$; then, we see that
\[(\alpha\vee\beta)\Uparrow\gamma={\sim}\big[(\alpha\vee\beta)\wedge\gamma\big],\]
and by using the distributivity of conjunction over disjunction this is equivalent to ${\sim}[(\alpha\wedge\gamma)\vee(\beta\wedge\gamma)]$, which by De Morgan's law is equivalent to
\[{\sim}(\alpha\wedge\gamma)\wedge{\sim}(\beta\wedge\gamma)=(\alpha\Uparrow\gamma)\wedge(\beta\Uparrow\gamma),\]
as we suggested before. Quite analogously, $(\alpha\wedge\beta)\Uparrow\gamma$ is equivalent to $(\alpha\Uparrow\gamma)\vee(\beta\Uparrow\gamma)$ in $\textbf{CPL}$.

To make $\bI$ more similar to the classical account of incompatibility, we add to it axiom schemata regarding the propagation of incompatibility\label{prwedge}\label{prvee}
\[\tag{$\textbf{pr}_{\wedge}$}\big[(\alpha\uparrow\gamma)\wedge(\beta\uparrow\gamma)\big]\rightarrow\big[(\alpha\vee\beta)\uparrow\gamma\big]\]
and
\[\tag{$\textbf{pr}_{\vee}$}\big[(\alpha\uparrow\gamma)\vee(\beta\uparrow\gamma)\big]\rightarrow\big[(\alpha\wedge\beta)\uparrow\gamma\big],\]
what gives us the logic $\bIpr$.

We can prove that $\bIpr$ is strictly stronger than $\bI$: take the instance $[(p_{1}\uparrow p_{3})\wedge(p_{2}\uparrow p_{3})]\rightarrow[(p_{1}\vee p_{2})\uparrow p_{3}]$ of $\textbf{pr}_{\wedge}$, which is clearly true in $\bIpr$, and take a bivaluation $\nu:F(\Sigma_{\bI}, \mathcal{V})\rightarrow\{0,1\}$ for $\bI$ such that $\nu(p_{1})=\nu(p_{2})=\nu(p_{3})=0$,
\[\nu((p_{1}\vee p_{2})\uparrow p_{3})=0,\quad \nu(p_{1}\uparrow p_{3})=1\quad\text{and}\quad\nu(p_{2}\uparrow p_{3})=1.\]
Then we have that 
\[\nu([(p_{1}\uparrow p_{3})\wedge(p_{2}\uparrow p_{3})]\rightarrow[(p_{1}\vee p_{2})\uparrow p_{3}])=\nu((p_{1}\uparrow p_{3})\wedge(p_{2}\uparrow p_{3}))\rightarrow\nu((p_{1}\vee p_{2})\uparrow p_{3})=\]
\[[\nu(p_{1}\uparrow p_{3})\wedge\nu(p_{2}\uparrow p_{3})]\rightarrow0=(1\wedge1)\rightarrow0=1\rightarrow0=0,\]
meaning $\bI$ can not actually prove some instances of $\textbf{pr}_{\wedge}$.


\subsubsection{Bivaluations}

\begin{definition}\label{Bivaluation for bIpr}
A bivaluation for $\bIpr$ is a map $\nu:F(\Sigma_{\bI}, \mathcal{V})\rightarrow\{0,1\}$ that is a bivaluation for $\bI$ and also satisfies that, for any formulas $\alpha$, $\beta$ and $\gamma$,
\[\text{if}\quad \nu(\alpha\uparrow\gamma)=1\quad\text{and}\quad\nu(\beta\uparrow\gamma)=1,\quad\text{then}\quad\nu((\alpha\vee\beta)\uparrow\gamma)=1,\]
and
\[\text{if}\quad\nu(\alpha\uparrow\gamma)=1\quad\text{or}\quad\nu(\beta\uparrow \gamma)=1,\quad\text{then}\quad \nu((\alpha\wedge\beta)\uparrow\gamma)=1.\]
\end{definition}

\begin{definition}
Given a set of formulas $\Gamma\cup\{\varphi\}$ of $\bIpr$, we say $\Gamma$ proves $\varphi$ according to bivaluations for $\bIpr$, and write $\Gamma\vDash_{\bIpr}\varphi$\label{vDashbIpr}, if for every bivaluation $\nu$ for $\bIpr$ we have that, if $\nu(\Gamma)\subseteq\{1\}$, then $\nu(\varphi)=1$.
\end{definition}

As it would be expected, given instances $\phi=[(\alpha\uparrow\gamma)\wedge(\beta\uparrow\gamma)]\rightarrow[(\alpha\vee\beta)\uparrow\gamma]$ of $\textbf{pr}_{\wedge}$ and $\psi=[(\alpha\uparrow\gamma)\vee(\beta\uparrow\gamma)]\rightarrow[(\alpha\wedge\beta)\uparrow\gamma]$ of $\textbf{pr}_{\vee}$, we have that $\vDash_{\bIpr}\phi$ and $\vDash_{\bIpr}\psi$, meaning that for any bivaluation $\nu$ for $\bIpr$, $\nu(\phi)=\nu(\psi)=1$.

To prove that, we will start with $\phi$: we have $\nu(\phi)=1$ if and only if $\nu((\alpha\vee\beta)\uparrow\gamma)=1$, or 
\[\nu((\alpha\uparrow\gamma)\wedge(\beta\uparrow\gamma))=0.\]
If $\nu((\alpha\uparrow\gamma)\wedge(\beta\uparrow\gamma))=0$ there is nothing to be done, so let us assume $\nu((\alpha\uparrow\gamma)\wedge(\beta\uparrow\gamma))=1$: in this case, since $\nu$ is a bivaluation for $\bI$, $\nu(\alpha\uparrow\gamma)=1$ and $\nu(\beta\uparrow\gamma)=1$, meaning $\nu((\alpha\vee\beta)\uparrow\gamma)=1$. Either way, $\nu(\phi)=1$.

Now $\nu(\psi)=1$ if and only if $\nu((\alpha\wedge\beta)\uparrow\gamma)=1$ or 
\[\nu((\alpha\uparrow\gamma)\vee(\beta\uparrow\gamma))=0;\]
as before, if $\nu((\alpha\uparrow\gamma)\vee(\beta\uparrow\gamma))=0$ there is nothing to prove, so let us assume $\nu((\alpha\uparrow\gamma)\vee(\beta\uparrow\gamma))=1$. In this case, $\nu(\alpha\uparrow\gamma)=1$ or $\nu(\beta\uparrow\gamma)=1$, implying $\nu((\alpha\wedge\beta)\uparrow\gamma)=1$. Once again, we derive that in any case one has $\nu(\psi)=1$.

\begin{theorem}
Given formulas $\Gamma\cup\{\varphi\}$ of $\bIpr$, if $\Gamma\vdash_{\bIpr}\varphi$ then $\Gamma\vDash_{\bIpr}\varphi$.
\end{theorem}

To prove the reciprocal, we use the standard method of defining $\nu:F(\Sigma_{\bI}, \mathcal{V})\rightarrow\{0,1\}$, for a closed, non-trivial set of formulas $\Gamma$ of $\bIpr$ maximal with respect to not proving a given $\varphi$, as $\nu(\gamma)=1$ if and only if $\gamma\in \Gamma$.

We already know that $\nu$ is a bivaluation for $\bI$, since a set of formulas that is closed under $\bIpr$ must be closed under $\bI$, given $\bIpr$ extends $\bI$. So assume $\nu(\alpha\uparrow\gamma)=1$ and $\nu(\beta\uparrow\gamma)=1$, meaning that $\alpha\uparrow\gamma, \beta\uparrow\gamma\in \Gamma$; from the instance 
\[\big[(\alpha\uparrow\gamma)\wedge(\beta\uparrow\gamma)\big]\rightarrow\big[(\alpha\vee\beta)\uparrow\gamma\big]\]
of $\textbf{pr}_{\wedge}$ and the fact that $\Gamma$ is closed, we derive $(\alpha\vee\beta)\uparrow\gamma\in\Gamma$, and therefore $\nu((\alpha\vee\beta)\uparrow\gamma)=1$, proving the first extra condition of Definition \ref{Bivaluation for bIpr}. To prove the remaining property, assume $\nu(\alpha\uparrow\gamma)=1$ or $\nu(\beta\uparrow\gamma)=1$, meaning either $\alpha\uparrow\gamma\in\Gamma$ or $\beta\uparrow\gamma\in\Gamma$. From the instance 
\[\big[(\alpha\uparrow\gamma)\vee(\beta\uparrow\gamma)\big]\rightarrow\big[(\alpha\wedge\beta)\uparrow\gamma\big]\]
of $\textbf{pr}_{\vee}$ and the fact $\Gamma$ is closed, we get that $(\alpha\wedge\beta)\uparrow\gamma\in \Gamma$, and so $\nu((\alpha\wedge\beta)\uparrow\gamma)=1$; this finishes the proof.

\begin{theorem}
Given formulas $\Gamma\cup\{\varphi\}$ of $\bIpr$, if $\Gamma\vDash_{\bIpr}\varphi$ then $\Gamma\vdash_{\bIpr}\varphi$.
\end{theorem}

Notice, and the importance of this fact shall become clearer when we arrive at Section \ref{Collapsing axioms part 1}, that $\textbf{bIpr}$ is not classical propositional logic, with $\alpha\uparrow\beta$ standing for $\alpha\Uparrow\beta={\sim}(\alpha\wedge\beta)$. To see this, take propositional variables $p$ and $q$ and a bivaluation $\nu$ for $\bIpr$ with $\nu(p)=\nu(q)=\nu(p\uparrow q)=0$: it is clear such a bivaluation is possible, and yet $\nu(p\Uparrow q)=1$, implying that $p\Uparrow q\not\vdash_{\bIpr}p\uparrow q$.


\subsubsection{Fidel Structures}

The Fidel structures, as $\Sigma_{\bI}^{\textbf{CPL}}$-multialgebras, for $\bIpr$ will simply be the same as those for $\bI$. For every such Fidel structure, presented as a $\Sigma_{\bI}^{\textbf{CPL}}$-multialgebra, we consider the restricted Nmatrix $(\mathcal{A}, \{\top\}, \mathcal{F}_{\mathcal{A}})$, where $\mathcal{F}_{\mathcal{A}}$ is the set of homomorphisms $\nu:\textbf{F}(\Sigma_{\bI}, \mathcal{V})\rightarrow\mathcal{A}$ such that:
\begin{enumerate}
\item $\nu(\alpha\uparrow\beta)=\nu(\beta\uparrow\alpha)$, for any formulas $\alpha$ and $\beta$;
\item for any instance $\phi$ of $\textbf{pr}_{\wedge}$, $\nu(\phi)=\top$;
\item for any instance $\psi$ of $\textbf{pr}_{\vee}$, $\nu(\psi)=\top$.
\end{enumerate}

It is clear how such conditions keep the generalized Nmatrices in question structural: if $\sigma$ is an endomorphism of $\textbf{F}(\Sigma_{\bI}, \mathcal{V})$, and $\phi$ is an instance of $\textbf{pr}_{\wedge}$ and $\psi$ is an instance of $\textbf{pr}_{\vee}$, then $\sigma(\phi)$ and $\sigma(\psi)$ remain instances of, respectively, $\textbf{pr}_{\wedge}$ and $\textbf{pr}_{\vee}$.

The second and third conditions in the previous definition can be modified, if we remember a few simple properties of Boolean algebras: to give one example, take an instance 
\[\phi=\big[(\alpha\uparrow\gamma)\wedge(\beta\uparrow\gamma)\big]\rightarrow\big[(\alpha\vee\beta)\uparrow\gamma\big]\]
of $\textbf{pr}_{\wedge}$; applying to $\phi$ a valuation $\nu$ in $\mathcal{F}_{\mathcal{A}}$ gives us
\[{\sim}\nu(\alpha\uparrow\gamma)\vee{\sim}\nu(\beta\uparrow\gamma)\vee\nu((\alpha\vee\beta)\uparrow\gamma),\]
so, asking that $\nu(\phi)=\top$ is equivalent to asking that $x\vee{\sim}y_{1}\vee{\sim}y_{2}=\top$, where $\nu((\alpha\vee\beta)\uparrow\gamma)=x$, $\nu(\alpha\uparrow\gamma)=y_{1}$ and $\nu(\beta\uparrow\gamma)=y_{2}$.

If $\Gamma$ proves $\varphi$ according to the class of these restricted Nmatrices, we write $\Gamma\Vdash_{\mathcal{F}}^{\bIpr}\varphi$\label{VdashFbIpr}.

We already know that the axioms and rules of inference of $\bI$ are modeled by "$\Vdash_{\mathcal{F}}^{\bIpr}$", remaining for us to show that this is also true for those axioms particular to $\bIpr$, that is, $\textbf{pr}_{\wedge}$ and $\textbf{pr}_{\vee}$, meaning that for any instance $\phi$ of $\textbf{pr}_{\wedge}$ and any instance $\psi$ of $\textbf{pr}_{\vee}$, we wish to prove $\Vdash_{\mathcal{F}}^{\bIpr}\phi$ and $\Vdash_{\mathcal{F}}^{\bIpr}\psi$.

This translates to proving that, for any Fidel structure $\mathcal{A}$ for $\bI$ and any homomorphism $\nu\in\mathcal{F}_{\mathcal{A}}$, $\nu(\phi)=\top$ and $\nu(\psi)=\top$, which is trivially true given our definition of $\mathcal{F}_{\mathcal{A}}$.

\begin{theorem}
Given formulas $\Gamma\cup\{\varphi\}$ of $\bIpr$, if $\Gamma\vdash_{\bIpr}\varphi$ then $\Gamma\Vdash_{\mathcal{F}}^{\bIpr}\varphi$.
\end{theorem}

For the reciprocal, once again we define an equivalence relation on the formulas of $\bIpr$ such that $\alpha\equiv_{\Gamma}^{\bIpr}\beta$ if and only if $\Gamma\vdash_{\bIpr}\alpha\rightarrow\beta$ and $\Gamma\vdash_{\bIpr}\beta\rightarrow\alpha$. As mentioned, this is indeed an equivalence relation, but also a congruence for the connectives in $\{\vee,\wedge,\rightarrow\}$, allowing us to make $A^{\bIpr}_{\Gamma}=F(\Sigma_{\bI}, \mathcal{V})/\equiv_{\Gamma}^{\bIpr}$ into a Boolean algebra as we did when studying $\bI$ in Section \ref{Fidel structures for bI}.

We define, for formulas $\alpha$ and $\beta$,
\[\uparrow_{\mathcal{A}^{\bIpr}_{\Gamma}}([\alpha],[\beta])=\{[\varphi\uparrow\psi]\ :\  \varphi\in[\alpha], \psi\in[\beta]\}.\]

For simplicity, let us drop the index $\mathcal{A}^{\bIpr}_{\Gamma}$ and use the infix notation; first of all, it is clear how, with such a definition, $\mathcal{A}^{\bIpr}_{\Gamma}=(A^{\bIpr}_{\Gamma}, \{\sigma_{\mathcal{A}}\}_{\sigma\in\Sigma_{\uparrow}})$ is indeed a Fidel structure for $\bI$, presented as a $\Sigma_{\bI}^{\textbf{CPL}}$-multialgebra. We will call $\mathcal{A}^{\bIpr}_{\Gamma}$ the Lindenbaum-Tarski Fidel structure of $\bIpr$ associated to $\Gamma$.

\begin{theorem}
Given formulas $\Gamma\cup\{\varphi\}$ of $\bIpr$, if $\Gamma\Vdash_{\mathcal{F}}^{\bIpr}\varphi$, then $\Gamma\vdash_{\bIpr}\varphi$.
\end{theorem}

\begin{proof}
Suppose $\Gamma\not\vdash_{\bIpr}\varphi$ and take a maximal, with respect to not proving $\varphi$, closed extension $\Delta$ of $\Gamma$ and consider: the Lindenbaum-Tarski Fidel structure $\mathcal{A}^{\bIpr}_{\Delta}$ of $\bIpr$ associated to $\Delta$, and the map $\nu:F(\Sigma_{\bI}, \mathcal{V})\rightarrow A^{\bIpr}_{\Delta}$ such that $\nu(\alpha)=[\alpha]$.

Clearly $\nu$ is a $\Sigma_{\bI}$-homomorphism, given $\nu(\alpha\#\beta)=\nu(\alpha)\#\nu(\beta)$ for $\#\in\{\vee, \wedge, \rightarrow\}$ and $\nu(\alpha\uparrow\beta)\in[\alpha]\uparrow[\beta]$. It remains for us to show that $\nu\in\mathcal{F}_{\mathcal{A}}$:
\begin{enumerate}
\item since $(\alpha\uparrow\beta)\rightarrow(\beta\uparrow\alpha)$ and $(\beta\uparrow\alpha)\rightarrow(\alpha\uparrow\beta)$ are instances of $\textbf{Comm}$, we have that $[\alpha\uparrow\beta]=[\beta\uparrow\alpha]$, and therefore $\nu(\alpha\uparrow\beta)=\nu(\beta\uparrow\alpha)$;
\item given an instance $\phi$ of $\textbf{pr}_{\vee}$, we have that $\Gamma\vdash_{\bIpr}\phi\rightarrow\top_{\phi}$ and $\Gamma\vdash_{\bIpr}\top_{\phi}\rightarrow\phi$, meaning $[\phi]=\top$ and therefore $\nu(\phi)=\top$;
\item for the case of $\textbf{pr}_{\wedge}$, the proof follows that of $\textbf{pr}_{\vee}$.
\end{enumerate}
We see that $\nu(\alpha)=\top$ if and only if $\Delta\vdash_{\bIpr}\alpha$, and so $\Delta\not\Vdash_{\mathcal{F}}^{\bIpr}\varphi$, meaning that $\Gamma\not\Vdash_{\mathcal{F}}^{\bIpr}\varphi$.
\end{proof}

\subsubsection{Decision method}

Take the Fidel structure $\textbf{2}_{\bI}$ for $\bI$: we comment now on how $(\textbf{2}_{\bI}, \{1\}, \mathcal{F}_{\textbf{2}_{\bIpr}})$, which we shall denote by $\mathbb{2}_{\bIpr}$\label{2bIpr},  is a decision method for $\bIpr$, where $\mathcal{F}_{\textbf{2}_{\bIpr}}$ is the set of homomorphisms $\nu:\textbf{F}(\Sigma_{\bI}, \mathcal{V})\rightarrow\textbf{2}_{\bI}$ such that:
\begin{enumerate}
\item $\nu(\alpha\uparrow\beta)=\nu(\beta\uparrow\alpha)$, for any formulas $\alpha$ and $\beta$;
\item for any instance $\phi$ of $\textbf{pr}_{\vee}$, $\nu(\phi)=1$;
\item for any instance $\psi$ of $\textbf{pr}_{\wedge}$, $\nu(\psi)=1$.
\end{enumerate}

\begin{theorem}
$\nu:F(\Sigma_{\bI},\mathcal{V})\rightarrow\{0,1\}$ is a bivaluation for $\bIpr$ if, and only if, it is a $\Sigma_{\bI}$-homomorphism from $\textbf{F}(\Sigma_{\bI},\mathcal{V})$ to $\textbf{2}_{\bI}$ that lies in $\mathcal{F}_{\textbf{2}_{\bIpr}}$.
\end{theorem}

\begin{theorem}
Given formulas $\Gamma\cup\{\varphi\}$ of $\bIpr$, $\Gamma\vDash_{\bIpr}\varphi$ if and only if $\Gamma\vDash_{\mathbb{2}_{\bIpr}}\varphi$.
\end{theorem}

The row-branching, row-eliminating truth tables based on $\mathbb{2}_{\bIpr}$ are then constructed more or less like the ones for $\bI$ found in Section \ref{Decision Method for bI}, the key difference being that rows where $\alpha\uparrow\gamma$ and $\beta\uparrow\gamma$ are both $1$ but $(\alpha\vee\beta)\uparrow\gamma$ is $0$, or where $\alpha\uparrow\gamma$ or $\beta\uparrow\gamma$ is $1$ but $(\alpha\wedge\beta)\uparrow\gamma$ is $0$, must also be disconsidered.


\section{Collapsing axioms}\label{Collapsing axioms part 1}


\subsection{$\textbf{Ex}$}

A question that naturally arises when dealing with incompatibility is the partial equivalence between $\alpha$ and $\beta$ being incompatible, and $\alpha$ and $\beta$, together, trivializing a logic.

That is, we know that if $\alpha$, $\beta$, and $\alpha\uparrow\beta$ are simultaneously valid, then our logic becomes trivial; but if $\alpha$ and $\beta$, just the two of them, can trivialize our logic, does that mean $\alpha$ and $\beta$ are incompatible? We believe such a question must reference the philosophical notion hoped for the incompatibility to encompass, and for this end we consider the axiom\label{Ex}
\[\tag{\textbf{Ex}}(\alpha\wedge\beta\rightarrow\bot_{\alpha\beta})\rightarrow(\alpha\uparrow\beta).\]
The logic obtained from $\bI$ by addition of $\textbf{Ex}$ will be called $\bI\textbf{Ex}$\label{bIEx}, but that is not really a new logic: we shall prove that, in this case, $\alpha\uparrow\beta$ is equivalent to $\alpha\wedge\beta\rightarrow\bot_{\alpha\beta}$ or, what is the same, $\alpha\Uparrow\beta={\sim}(\alpha\wedge\beta)$ when we consider the classical negation definable in $\bI$, and therefore $\bI\textbf{Ex}$ becomes again classical propositional logic with a shorthand notation for $\alpha\wedge\beta\rightarrow\bot_{\alpha\beta}$, that is, $\alpha\uparrow\beta$.

So we wish to prove that $\alpha\uparrow\beta$ is equivalent, in this $\bI\textbf{Ex}$, to $\alpha\wedge\beta\rightarrow\bot_{\alpha\beta}$, meaning that $(\alpha\wedge\beta\rightarrow\bot_{\alpha\beta})\rightarrow(\alpha\uparrow\beta)$ but also that
\[(\alpha\uparrow\beta)\rightarrow(\alpha\wedge\beta\rightarrow\bot_{\alpha\beta}).\]
Of course, the first direction of the bi-implication follows straightforwardly from $\textbf{Ex}$, being an instance of the axiom. To prove the second direction, take the instance $(\alpha\uparrow\beta)\rightarrow(\alpha\rightarrow(\beta\rightarrow\bot_{\alpha\beta}))$ of $\textbf{Ip}$, and by applying the meta-deduction theorem three times we find that 
\[\alpha\uparrow\beta, \alpha, \beta\vdash_{\bI\textbf{Ex}}\bot_{\alpha\beta}.\]
By applying the deduction meta-theorem to the instances $\alpha\wedge\beta\rightarrow\alpha$ and $\alpha\wedge\beta\rightarrow\beta$ of, respectively, $\textbf{Ax\: 4}$ and $\textbf{Ax\: 5}$, we obtain that $\alpha\wedge\beta\vdash_{\bI\textbf{Ex}}\alpha$ and $\alpha\wedge\beta\vdash_{\bI\textbf{Ex}}\beta$, and therefore 
\[\alpha\uparrow\beta, \alpha\wedge\beta\vdash_{\bI\textbf{Ex}}\bot_{\alpha\beta},\]
implying that $\alpha\uparrow\beta\vdash_{\bI\textbf{Ex}}\alpha\wedge\beta\rightarrow\bot_{\alpha\beta}$ and, therefore, that $(\alpha\uparrow\beta)\rightarrow(\alpha\wedge\beta\rightarrow\bot_{\alpha\beta})$ is a tautology in $\bI\textbf{Ex}$.


\subsection{$\textbf{ciw}^\uparrow$}

The axiom $\textbf{ciw}$, given by
\[\circ\alpha\vee(\alpha\wedge\neg \alpha),\]
is a very important one when dealing with paraconsistency, and adding it to $\textbf{mbC}$ gives us the logic $\textbf{mbCciw}$. When dealing with incompatibility, we would like to define an equivalent axiom, namely\label{ciwuparrow}
\[\tag{$\textbf{ciw}^\uparrow$}(\alpha\uparrow\beta)\vee(\alpha\wedge\beta).\]
The logic obtained from $\bI$ by adding $\textbf{ciw}^\uparrow$ will be called $\bI\textbf{ciw}^\uparrow$\label{bIciw}, and we will prove that, as it happens with $\bI\textbf{Ex}$, $\bI\textbf{ciw}^\uparrow$ is again equivalent to classical propositional logic, with $\alpha\uparrow\beta$ corresponding to $\alpha\wedge\beta\rightarrow\bot_{\alpha\beta}$.

The proof that $\vdash_{\bI\textbf{ciw}^\uparrow}(\alpha\uparrow\beta)\rightarrow(\alpha\wedge\beta\rightarrow\bot_{\alpha\beta})$ is the same as the one for the same fact in $\bI\textbf{Ex}$, remaining for us to show that $(\alpha\wedge\beta\rightarrow\bot_{\alpha\beta})\rightarrow(\alpha\uparrow\beta)$ is a tautology in $\bI\textbf{ciw}^\uparrow$.

By an application of the deduction meta-theorem, the desired result is equivalent to $\alpha\wedge\beta\rightarrow\bot_{\alpha\beta}\vdash_{\bI\textbf{ciw}^{\uparrow}}\alpha\uparrow\beta$. It is clear how $\alpha\uparrow\beta, \alpha\wedge\beta\rightarrow\bot_{\alpha\beta}\vdash_{\bI\textbf{ciw}^{\uparrow}}\alpha\uparrow\beta$, and by an application of Modus Ponens and the fact $\bot_{\alpha\beta}$ behaves like a bottom element, it is obvious that
\[\alpha\wedge\beta, \alpha\wedge\beta\rightarrow\bot_{\alpha\beta}\vdash_{\bI\textbf{ciw}^{\uparrow}}\alpha\uparrow\beta.\]
By a proof by cases, we get that 
\[(\alpha\uparrow\beta)\vee(\alpha\wedge\beta), \alpha\wedge\beta\rightarrow\bot_{\alpha\beta}\vdash_{\bI\textbf{ciw}^{\uparrow}}\alpha\uparrow\beta,\]
and since $(\alpha\uparrow\beta)\vee(\alpha\wedge\beta)$ is an instance of an axiom of $\bI\textbf{ciw}^{\uparrow}$, we get the desired result.

\section{Brandom's notion of incompatibility}

Here, we offer a brief overview of Robert B. Brandom\index{Brandom} and Alp Aker's work with incompatibility, echoed in Jaroslav Peregrin's\index{Peregrin} research. The differences between their own methodology and ours are many: they focus on incompatibility between sets of formulas,; they aim to define consequence from incompatibility, while we assume both to coexist; and, perhaps most importantly, although many systems can be retrieved from their methods, those of a paraconsistent behavior are not among them, as the only negations considered by them are classical, or at most intuitionistic. But the connection between both works is there, and quite clear: the attempt to control, in a way or another, explosion by mediating it through incompatibility.

\subsection{Brandom and Aker's ``Between saying and doing''}

We start by pointing out that, although the author of \cite{Brandom} appears as Robert B. Brandom, he makes it clear that much of the semantics for his logic of incompatibility was developed by his Ph.D. student Alp Aker; given our focus in his book is mostly restricted to the chapters about incompatibility, we will refer to both Brandom and Aker as authors.

In \cite{Brandom}, Brandom and Aker propose an approach to logic through incompatibility, instead of consequence; they defend that a modal understanding of incompatibility fits better with the argumentative nature of epistemology, in his neo-pragmatic program. They proceed to redefine consequence, and the usual connectives, starting from incompatibility, or rather incoherence\index{Incoherence}, as a primitive notion.

Initially, Brandom and Aker start with the notions of commitment\index{Commitment} and entitlement\index{Entitlement}, which are very modal in nature. Intuitively, an agent is committed to an statement when stating such statement is necessary for the agent; and the same agent is entitled to a statement when stating it is possible for the agent. Building from this, two statements are incompatible when being committed to one implies not being entitled to the other. A notion of consequence arises from this concept when we define that a statement $p$ implies $q$ whenever every statement incompatible with $q$ is incompatible with $p$; clearly we must reach for second-order logic and quantification over formulas for this procedure to make sense.

Our approach to incompatibility is quite close, in a sense, to this first notion of Brandom and Aker, but they proceed to extend their definition: they sustain that incompatibility must deal with pairs of sets of statements, instead of pairs of statements; their argumentation is that a claim may be incompatible with a set of claims without being incompatible with any one claim of the set. First of all, we believe this problem may be circumvented by a careful distinction between incompatibility in natural discourse and as expressed formally, although this distinction may be blurred given Brandom and Aker want incompatibility to originate logic, what in turn may explain their position; second, we do not aim to encompass every reasoning involving incompatibility as Brandom  and Aker do, but rather prefer to focus on its interplay with paraconsistency and its possible semantics.

Here, a digression seems useful. In the beginning of this chapter we offered interpretations of incompatibility, two of them being related to probability: mutual exclusivity and independence. Now, if $\{A_{i}\}_{i\in I}$ is a collection of pairwise mutually exclusive events, for any $J\subseteq I$ with at least two elements $j_{1}$ and $j_{2}$ one has that $P(\bigcap_{j\in J}A_{j})\leq P(A_{j_{1}}\cap A_{j_{2}})=0$, and so $P(\bigcap_{j\in J}A_{j})=0$; this means that a set of pairwise mutually exclusive events also has a generalized mutual exclusivity. Meanwhile, if $\{A_{i}\}_{i\in I}$ is instead a collection of pairwise independent events, the fact that, for any two $i_{1}, i_{2}\in I$, one has $P(A_{i_{1}}\cap A_{i_{2}})=P(A_{i_{1}})P(A_{i_{2}})$ does not imply that, for any finite $J\subseteq I$, one also has $P(\bigcap_{j\in J}A_{j})=\prod_{j\in J}P(A_{j})$; because of that, we define instead that a finite collection of events $\{A_{i}\}_{i\in I}$ is mutually independent, something much stronger than merely pairwise independent, if, for any $J\subseteq I$, 
\[P(\bigcap_{j\in J}A_{j})=\prod_{j\in J}P(A_{j})\]
holds. This difference between these two notions carries a lot of similarities with the difference between ours, and Brandom and Aker`s approach: in our understanding of incompatibility, a set of propositions is incompatible if, and only if, any pair of its elements is, itself, incompatible, exactly as in the case of mutually exclusive events; meanwhile, in Brandom and Aker, the more general notion of incompatible set can not be reduced to pairwise incompatibility, as in the problem of mutual independence as opposed to pairwise independence. This suggests interesting connections between distinct incompatibilities and several concepts in probability theory.

Then, \cite{Brandom} demands only two properties of incompatibility, derived from its natural interpretation and intuition:
\begin{enumerate}
\item if $X$ is incompatible with $Y$, $Y$ is incompatible with $X$ (symmetry)\index{Symmetry};
\item if $X$ is incompatible with $Y$, and $Z$ contains $Y$, then $X$ is incompatible with $Z$ (persistence)\index{Persistence}.
\end{enumerate}

In our logics of incompatibility, symmetry is analogous to the commutative axiom $\textbf{Comm}$, where $\alpha\uparrow\beta\rightarrow\beta\uparrow\alpha$; regarding persistence, ignoring the obvious interpretation in second-order logic, one very natural take would be, in our language, the axiom
\[(\alpha\uparrow\beta)\rightarrow\big((\gamma\rightarrow\beta)\rightarrow(\alpha\uparrow\gamma)\big).\]
But we also may, for sets of formulas $X$ and $Y$, write a generalized incompatibility operator $X\uparrow Y$ whenever there exist formulas $\alpha$ and $\beta$ such that $X\vdash \alpha$, $Y\vdash\beta$ and $\alpha\uparrow\beta$. Then, if $Z$ contains $Y$ and $X\uparrow Y$, given formulas $\alpha$ and $\beta$ as above, since $Y\subseteq Z$ one finds $Z\vdash \beta$, and therefore $X\uparrow Z$ and persistence is reobtained on this environment. Of course, this is not to say that our logics of incompatibility naturally model Brandom and Aker's approach, but rather to show they are plastic enough to do so.

Brandom and Aker stress that their incompatibility should not be limited to truth values, as they intend to define the latter starting from the former, but they still mention \textit{en passant} that the obvious interpretation of, restricting ourselves to single statements once again, $p$ being incompatible to $q$ should be that it is impossible to have both $p$ and $q$ simultaneously true; this, in turn, is of course very distant from our own take on incompatibility, as it limits the concept to a modal version of Sheffer's stroke $\square{\sim}(p\wedge q)$.

To further distance its incompatibility from ours, \cite{Brandom} takes as a tacit starting point that a statement should be incompatible with its negation, appealing to classical logic to do so; this may be justified if they intend to replicate merely classical negation from incompatibility, but it seems an inadmissible loss not to consider more modern negations, objects around which our logics of incompatibility are actually built. It seems that, to define a negation from incompatibility, it makes more sense, and is more interesting, to start from the assumption there is no prior connection between incompatibility and negation; then one gains the interesting notion of ``well-behavedness'' among statements, a resource well established amid logics of formal inconsistency. 

Nevertheless, Brandom and Aker then define the negation of a claim $p$ as the minimum claim incompatible to $p$, meaning a claim ${\sim}p$ incompatible with $p$ such that, whenever $X$ is a set of statements incompatible with $p$, one has $X\vdash {\sim}p$. Notice that, using the classical negation that may be defined in $\bI$, we find $\alpha\uparrow\beta\rightarrow(\alpha\rightarrow{\sim}\beta)$, meaning that if $p$ and $q$ are incompatible, then $p$ implies the negation of $q$; this may seem awfully close to Brandom and Aker's account, but notice ${\sim}q$ is not, necessarily, incompatible with $q$, and we do not intend to define negation from this relationship, but rather derive the relationship from the definition of the strong negation and the fact $\alpha\vee(\alpha\rightarrow\beta)$ is an instance of an axiom. It becomes clear \cite{Brandom} has no interest in dealing with paraconsistency when it defines an inconsistent set of formulas as any set which derives both a claim and its negation; this, of course, is not interesting to us at all, since paraconsistency intends to precisely avoid this.

Brandom and Aker define the conjunction of statements $p$ and $q$ as the minimal statement incompatible with every set $X$ incompatible with $\{p,q\}$: they then go on to show that the logic obtained from both these connectives they have defined is classical propositional logic, what shows one may recover classical accounts of logic from the notion of incompatibility; although not in line with what we hope to accomplish with our own logics of incompatibility, Brandom and Aker's result is no less important, and establishes that incompatibility, at least in classical logic, is a notion no less important than that of deduction.

We have not studied modalities and their logical counterpart very deeply, but this is one of Brandom and Aker's concern, and they define that something is incompatible with the necessity of a statement $p$, $\square p$, whenever it is compatible with something that does not imply $p$; meaning, $X$ is incompatible with $\square p$ if and only if it is compatible with a $Y$ such that $Y\not\vdash p$, and one may define necessity as the minimal statement which can replace $\square p$ in the previous discussion. Very interestingly, this definition makes of the system at hand precisely $S5$. Nevertheless, one may wonder if the definitions \cite{Brandom} takes are the ideal ones: if we are recovering only the most standard systems, perhaps a little more plasticity should be allowed when redefining connectives from incompatibility.

In one last stretch, approaching their own account again to ours, Brandom and Aker notice that, in his logic of incompatibility, connectives do not have the semantic sub-formula property, meaning they are not functional: one can not derive the interpretation of a formula by merely looking at the interpretation of its sub-formulas. This non-deterministic behavior is very clear in our logics of incompatibility, and specially clear when one analyses the easiness with which one treats logics such as $\bI$ and $\nbI$ (see Chapter \ref{Chapter8}) with the aid of restricted non-deterministic matrices and, more generally, multialgebras.

In a more formal treatment of their system, Brandom and Aker wish to analyze when a consequence relation $\vdash$ may be characterized by incompatibility, in their own terms of ``to be characterized'' as previously described. They reach two necessary and sufficient conditions:
\begin{enumerate}
\item if $X\vdash Y$ and $Y\cup W\vdash Z$, then $X\cup W\vdash Z$ (general transitivity);
\item if $X\not\vdash Y$, there exist a $W$ and a $Z$ such that, $X\cup Z\not\vdash W$ but, for all $V$, $Y\cup Z\vdash V$ (defeasibility)\index{Defeasibility}.
\end{enumerate}

General transitivity is, of course, a reformulation of the cut rule; defeasibility states that, if $Y$ is not a consequence of $X$, then there exists something that when added to $Y$ will collapse the deduction, but not when added to $X$. Although somewhat natural, the way this conditions are used to achieve incompatibility is diametrically opposed to what we look for, as two sets are deemed incompatible whenever their union is inconsistent; disregarding our problem with inconsistency, the axiom schema $\textbf{Ip}$ invites one to produce instead the concept that, if two sets are incompatible, then their union should be trivializing, and not the other way around. Brandom and Aker go on to give several interesting definitions, but all based on these basic concepts, and therefore distant from our ideal understanding of incompatibility. Concisely, given a fixed set of sets of incoherent (intuitively, inconsistent) formulas, they demand of it only that if $X$ is incoherent and $X\subseteq Y$, then $Y$ is also incoherent; then, $X$ is incompatible with $Y$ if, and only if, $X\cup Y$ is incoherent. Then, $X$ derives a statement $p$, regarding their incompatibilities, when every set incompatible to $\{p\}$ is incompatible to $X$. Several strong theorems are then provided for those systems, but to showcase how they are not aligned to our notion of incompatibility, one of these theorems is that $p\wedge{\sim}p$ derives any set $Y$, meaning their negation is necessarily classically behaved.

\subsection{Peregrin}

\subsubsection{``Brandom's incompatibility semantics''}

In this first article \cite{Peregrin1}, the author Jaroslav Peregrin studies the notion of incompatibility, and their corresponding semantics, of Brandom and Aker's \cite{Brandom} from a more philosophical standpoint. More precisely, his main concerns are tied to the problem of whether formal semantics are truly compatible with pragmatic and inferentialist views and of how should such semantics look like; his argument, a very interesting one, may be summarized as stating that, as merely a model for natural processes, formal semantics can indeed be used by the pragmatist, as long as the distinction between model and what is modeled is not ignored.

More to the point we are interested in, Peregrin asks how should incompatibility come into play in logic: he first stresses how the usual approach is to say that sets of formulas $X$ and $Y$ are incompatible whenever their union can deduce anything; much in line with our own views of the problem, Peregrin continues to state that, however standard this approach may be, reducing incompatibility to inference is, first of all, wasteful, as it disregards the possible intricacies the concept may carry; second, incompatibility, as a byproduct of inference, becomes dependent on how strong is the logic we work over, meaning how expressive it is.

\cite{Peregrin1} proceeds to point some problematic aspects of Brandom and Aker's interpretation of incompatibility: first of all, when one defines that $X\vdash p$ whenever everything that is incompatible with $p$ is incompatible with $X$, this may be read as that everything compatible with $X$ is compatible with $p$, implying compatibility is preserved by consequence, and since so is truth, the two concepts may become at moments indistinguishable.  Of course, this is not of much concern for us, since we do not intend to rebuild consequence from incompatibility, but rather have incompatibility to exist in an environment with a predefined consequence operator. A second problem, still not very troublesome for us, is that consequence is usually taken to have a finitary character, while incompatibility, at least in Brandom and Aker's take, often involves quantification over all statements: so, by defining a consequence from incompatibility, one may end with a consequence operator that is more model-theoretic than proof-theoretic, that is, that has a clear non-finitary propensity.

To us, one of the most important developments found in this article is the connection established between incompatibility and Kripke\index{Kripke} semantics, \textit{id est}, semantics of possible worlds. Peregrin defines a possible world, once a concept of coherence is given, as a maximal coherent set of formulas; here, a set is a coherent one once it is not incompatible with itself, a definition that goes back to the work of Brandom and Aker. The truth of a statement on a given world is then taken to be the belonging of this statement to the wold, which is, again, a set of formulas; here, Peregrin notices once again how coherence behaves dangerously alike truth in Brandom and Aker's \cite{Brandom}. Notice how, given a statement $p$ and a possible world $w$, $w$ proves $p$ in this semantic of possible worlds if and only if $w$ derives $p$, or what is equivalent, when $p$ is compatible with $w$. Reciprocally, he derives a notion of incoherence from a Krikpkean semantics of possible worlds by saying a set of formulas $X$ is incoherent if no formula of this set is validated in any world.

Peregrin also points out how Brandom and Aker's definition of the necessity of a statement $p$ is equivalent to, in his semantic of possible words, validity: so a possible world $w$ verifies $\square p$ if, and only if, the set $\{p\}$ is not incoherent; this of course means that Brandom and Aker's take on modal logic trough incompatibility only derives the most basic interpretations of necessity and possibility. He offers an interesting alternative, in order to add richness to those modal logics: a second level of incompatibility, or rather incoherence, a meta-coherence\index{Meta-coherence} if you will; Brandom and Aker only take into account a set of sets of formulas, said to be the set of incoherent sets of formulas, while Peregrin goes further and defines a set of sets of sets of formulas, corresponding to a set of incoherent sets of sets of formulas, that offers a second level of control over incompatibility. With this, he is able to characterize modal logics more complex than $S5$, what of course suggests the power such generalizations may carry, although this level of abstractness is not within our scope of research.

\subsubsection{``Logic as based on incompatibility''}

In this second of Peregrin's article of great interest to us, he lays clear what he believes are the reasons for Brandom and Aker's incompatibility program: first, to study what the possible minimal foundations of logic may be; and second, to offer a working systems for those philosophers who have based their reasoning on incompatibility alone. While still reasoning about Brandom and Aker's program, the author defends that the most natural logic to arise from defining inference through incompatibility, in the opposite direction of Brandom and Aker's work, is intuitionistic; of course, this is not of great interest to us, since we hope for incompatibility to control \textit{ex contradictione quolibet} and not the other way around. Peregrin points out that, by defining inference via incompatibility, ``Between saying and doing'' reaches a logic of, instead of intuitionistic, classical character, but goes on to argue that this due not to the nature of incompatibility, but rather to Brandom and Aker's method.

\cite{Peregrin2} defines then an environment that should deal, simultaneously, with inference and incompatibility: a triple $(S, \bot, \vdash)$ is called a generalized inferential structure ($\textit{gis}$) when $S$ is a set, $\bot\subseteq\mathcal{P}(S)$ and $\vdash\subseteq\mathcal{P}(S)\times S$; for simplicity, if $X\in \bot$ we write $\bot X$, and if $X\cup Y\cup\{p\}$ is in $\bot$, we may simply write $\bot X,Y,p$. The basic conditions demanded of these objects are the following:
\begin{enumerate}
\item[$(\bot)$] if $\bot X$ and $X\subseteq Y$, $\bot Y$;
\item[$(\vdash)$]\begin{enumerate}\item $X, p\vdash p$;
\item if $X, p\vdash q$ and $Y\vdash p$, then $X,Y\vdash q$.
\end{enumerate}
\end{enumerate}

Peregrin suggests a possible interplay between the two concepts:
\begin{enumerate}
\item[$(\bot\vdash)$] if $\bot X$, then $X\vdash p$ for every $p$;
\item[$(\vdash\bot)$] if $X\vdash p$, then $\bot Y, p$ implies $\bot X,Y$ for every $Y$.
\end{enumerate}
Of course, the first condition is very in line with what we expect of incompatibility: that it derives anything and trivializes an argument, or, in other words, a controlled explosion. The second condition has been presented before, in a much more complicated way: it corresponds to Brandom and Aker's defeasibility. The author then points out how assuming $(\bot\vdash)$ and its converse amounts to reducing incompatibility to inference, while assuming $(\vdash\bot)$ and its converse reduce inference to incompatibility. \cite{Peregrin2} also has definitions of an incompatibility $\vartriangle$ defined trough inference and of an inference $\vartriangleright$ defined trough incompatibility:
\begin{enumerate}
\item $\vartriangle X$ if $X\vdash p$, for every $p$;
\item $X\vartriangleright p$ if, for every $Y$, $\bot Y, p$ implies $\bot X,Y$.
\end{enumerate}
Peregrin then offers interesting conditions, related to $(\bot\vdash)$ and $(\vdash\bot)$, under which $\vartriangle$ and $\bot$, and $\vartriangleright$ and $\vdash$ are equivalent. But, more importantly, the author states what the minimal requirements for an intuitionistic negation based on incompatibility should be:
\begin{enumerate}
\item $\bot p, {\sim}p$;
\item if $\bot X, p$, then $X\vdash{\sim}p$.
\end{enumerate}
Of course, the first condition is not desirable for our systems, while the second is derivable in systems as simple as $\nbI$ (which we define in Chapter \ref{Chapter8}); with this, we see that the plasticity that the author is hoping to obtain by modifying Brandom and Aker's stipulations does not encompass paraconsistency, nor is this his objective, as far as we can tell. He seems, instead, more concerned with modal systems, allowing our logics to fill a gap in his approach.

The additional requirement that $\bot X, {\sim}p$ implies $X\vdash p$ is then equivalent to stating that the negation in question is of classical behavior; Peregrin also shows that the definition Brandom and Aker give of negation is equivalent to his own negation with this last additional requirement, proving that the methodology found ``Between saying and doing'' can only lead to classical negation. Conjunction is defined as in Brandom and Aker, meaning 
\begin{enumerate}
\item if $\bot X, p\wedge q$, then $\bot X, p, q$;
\item if $\bot X, p, q$, then $\bot X, p\wedge q$.
\end{enumerate}
After comparing Brandom and Aker's, and his own, definitions with those necessary to define classical propositional logic trough inference alone, the author reaches the conclusion that while Brandom and Aker's approach inexorably leads to classical logic, varying the techniques found in his article is significantly more far-reaching. For an example, by removing the very basic condition that
\[\bot X\quad\text{and}\quad X\subseteq Y\quad\text{imply}\quad \bot Y,\]
which is equivalent, given the correct assumptions, to the fact that $X\vdash p$ implies $X, q\vdash p$, one finds a system of relevant logic; without going into further details, Peregrin suggests that changing incompatibility to a property between multisets of formulas, instead of sets, leads to linear logic.

\newpage
\printbibliography[segment=\therefsegment,heading=subbibliography]
\end{refsegment}

\begin{refsegment}
\defbibfilter{notother}{not segment=\therefsegment}
\setcounter{chapter}{7}
\chapter{Adding a negation to logics of incompatibility}\label{Chapter8}\label{Chapter 8}

When studying paraconsistent logics, our main focus is in the properties of negations, as is the case when working with paracomplete, i.e. intuitionistic, logics. Given such a prominent role negation plays in non-classical logics, is quite natural to shift our focus in the logics of incompatibility from $\uparrow$ to a non-classical negation, or better yet, to the possible interplay between $\uparrow$ and such a negation.

Our first step is adding such a negation: we have seen that in any logic extending $\bI$, is always possible to define a somewhat quite classical negation; for any formulas $\alpha$ and $\beta$, the formulas
\[\bot_{\alpha\beta}=\alpha\wedge(\beta\wedge(\alpha\uparrow\beta))\]
are all equivalent to each other, and furthermore play the expected role for the bottom, meaning that in the presence of $\textbf{Ip}$, for any formula $\gamma$ it is true that $\bot_{\alpha\beta}\rightarrow\gamma$. By having a bottom element, we can define a negation of classical behavior ${\sim} \alpha$ of a formula $\alpha$ by $\alpha\rightarrow \bot_{\alpha\alpha}$, or any other bottom since all of those are equivalent. 

We, therefore, need a new negation, weaker than that negation intrinsic to $\bI$. So we need a symbol for it: we define the signature $\Sigma_{\nbI}$\label{SigmanbI} as the signature obtained from $\Sigma_{\bI}$ by addition of an unary symbol $\neg$, that is, $(\Sigma_{\nbI})_{1}=\{\neg\}$, $(\Sigma_{\nbI})_{2}=\{\vee, \wedge, \rightarrow, \uparrow\}$ and $(\Sigma_{\nbI})_{n}=\emptyset$ for $n\notin\{1,2\}$.

Most of the research in this chapter was submitted as the preprint \cite{Frominconsistency}.

\section{The logic $\nbI$}

We start by adding to the simplest $\textbf{LIp}$, that is, $\bI$, a paraconsistent negation: so, to the axiom schemata and rules of inference of $\bI$, we add
\begin{enumerate}
\item[$\textbf{Ax\: 11}^{*}$] $\alpha\vee\neg \alpha$.
\end{enumerate}
We shall call this new logic $\nbI$\label{nbI}.


\subsection{Bivaluations}\label{Bivaluations for nbI}

A \index{Bivaluation for $\nbI$} bivaluation for $\nbI$ is a map $\nu:F(\Sigma_{\nbI},\mathcal{V})\rightarrow\{0,1\}$ such that:
\begin{enumerate}
\item if $\nu(\neg \alpha)=0$, $\nu(\alpha)=1$;
\item $\nu(\alpha\vee\beta)=1$ if and only if $\nu(\alpha)=1$ or $\nu(\beta)=1$;
\item $\nu(\alpha\wedge\beta)=1$ if and only if $\nu(\alpha)=\nu(\beta)=1$;
\item $\nu(\alpha\rightarrow\beta)=1$ if and only if $\nu(\alpha)=0$ or $\nu(\beta)=1$;
\item if $\nu(\alpha\uparrow\beta)=1$ and $\nu(\alpha)=1$, $\nu(\beta)=0$;
\item $\nu(\alpha\uparrow\beta)=\nu(\beta\uparrow\alpha)$.
\end{enumerate}

Given a set of formulas $\Gamma\cup\{\varphi\}$ of $\nbI$, we say $\Gamma$ proves $\varphi$ according to bivaluations, and write $\Gamma\vDash_{\nbI}\varphi$\label{vDashnbI}, if for every bivaluation $\nu$ for $\nbI$ we have that, if $\nu(\Gamma)\subseteq\{1\}$, then $\nu(\varphi)=1$.

Once most conditions demanded of a bivaluation for $\nbI$ are essentially the same as those demanded of a bivaluation for $\bI$ in Section \ref{bivaluations for bI}, we can easily find that "$\vDash_{\nbI}$" models those axioms and rules of inference of $\bI$: to see that it also models $\textbf{Ax\: 11}^{*}$, take a formula $\alpha$ and a bivaluation $\nu$ for $\nbI$
\begin{enumerate}
\item If $\nu(\neg \alpha)=1$, $\nu(\alpha\vee\neg \alpha)=1$ from the behavior of bivaluations for $\nbI$ regarding disjunctions.
\item If $\nu(\neg \alpha)=0$, we find $\nu(\alpha)=1$ and again $\nu(\alpha\vee\neg \alpha)=1$.
\end{enumerate}
Either way, one sees that "$\vDash_{\nbI}$" also models $\textbf{Ax\: 11}^{*}$.

\begin{theorem}
Given formulas $\Gamma\cup\{\varphi\}$ of $\nbI$, if $\Gamma\vdash_{\nbI}\varphi$ then $\Gamma\vDash_{\nbI}\varphi$.
\end{theorem}

The reciprocal result goes as it would be expected: for a closed, non-trivial set of formulas $\Gamma$ of $\nbI$ maximal with respect to not proving $\varphi$, we define $\nu:F(\Sigma_{\nbI}, \mathcal{V})\rightarrow\{0,1\}$ by means of $\nu(\gamma)=1$ if and only if $\gamma\in\Gamma$; $\nu$ is then a bivaluation for $\nbI$.

That it is a bivaluation for $\bI$ is quite obvious, remaining for us to show that, if $\nu(\neg \alpha)=0$, then $\nu(\alpha)=1$. So, suppose $\nu(\neg \alpha)=0$, meaning $\neg \alpha\not\in\Gamma$: since $\alpha\vee\neg \alpha$ is an instance of $\textbf{Ax\: 11}^{*}$ and $\Gamma$ is closed, $\alpha\vee\neg \alpha\in\Gamma$. 

If we had $\nu(\alpha)=0$, then we would have $\nu(\alpha\vee\neg \alpha)=0$, implying $\alpha\vee\neg \alpha\not\in\Gamma$, contradicting our previous remark. We must have then that $\nu(\alpha)=1$, what ends the proof that $\nu$ is a bivaluation for $\nbI$.

\begin{theorem}
Given formulas $\Gamma\cup\{\varphi\}$ of $\nbI$, if $\Gamma\vDash_{\nbI}\varphi$ then $\Gamma\vdash_{\nbI}\varphi$.
\end{theorem}


\subsection{Fidel structures}

Once again we find ourselves in need of bottom and top elements, and a classical negation, so that we must define the signature\label{SigmanbICPL} $\Sigma_{\nbI}^{\textbf{CPL}}$ with $(\Sigma_{\nbI}^{\textbf{CPL}})_{0}=\{\bot, \top\}$, $(\Sigma_{\nbI}^{\textbf{CPL}})_{1}=\{{\sim},\neg\}$, $(\Sigma_{\nbI}^{\textbf{CPL}})_{2}=\{\vee, \wedge, \rightarrow, \uparrow\}$ and $(\Sigma_{\nbI}^{\textbf{CPL}})_{n}=\emptyset$ for $n>2$.

So, it becomes easy to define Fidel structures for $\nbI$: \index{Fidel structure for $\nbI$}a Fidel structure for $\nbI$ is any $\Sigma_{\nbI}^{\textbf{CPL}}$-multialgebra $\mathcal{A}=(A, \{\sigma_{\mathcal{A}}\}_{\sigma\in\Sigma_{\nbI}^{\textbf{CPL}}})$ such that:
\begin{enumerate}
\item $(A,\{\sigma_{\mathcal{A}}\}_{\sigma\in\Sigma^{\textbf{CPL}}})$ is a Boolean algebra;
\item for all $a\in A$ and $b\in \neg _{\mathcal{A}}(a)$, $a\vee b=\top$;
\item for all $a,b\in A$ and $c\in\uparrow_{\mathcal{A}}(a,b)$, $a\wedge(b\wedge c)=\bot$;
\item for all $a, b\in A$, $a\uparrow b=b\uparrow a$.
\end{enumerate}

As we did before, we shall drop the index $\mathcal{A}$ and use the infix notation.
 
Given a Fidel structure $\mathcal{A}$, presented as a $\Sigma_{\nbI}^{\textbf{CPL}}$-multialgebra, for $\nbI$, a valuation for $\mathcal{A}$ is any $\Sigma_{\nbI}$-homomorphism $\nu:\textbf{F}(\Sigma_{\nbI}, \mathcal{V})\rightarrow \mathcal{A}$; and, for every Fidel structure $\mathcal{A}$ for $\nbI$, we will consider the restricted Nmatrix $(\mathcal{A}, \{\top\}, \mathcal{F}_{\mathcal{A}})$, where $\mathcal{F}_{\mathcal{A}}$ is the set of valuations $\nu:\textbf{F}(\Sigma_{\nbI}, \mathcal{V})\rightarrow\mathcal{A}$ such that 
\[\nu(\alpha\uparrow\beta)=\nu(\beta\uparrow\alpha),\]
for any two formulas $\alpha$ and $\beta$ in $F(\Sigma_{\nbI}, \mathcal{V})$; it is clear how these RNmatrices are structural. If $\Gamma$ proves $\varphi$ according to such restricted Nmatrices, we will write\label{VdashFnbI} $\Gamma\Vdash_{\mathcal{F}}^{\nbI}\varphi$.

It is easy to see that if $\phi$ is an instance of an axiom of $\bI$, $\Vdash_{\mathcal{F}}^{\nbI}\phi$, and if $\Vdash_{\mathcal{F}}^{\nbI}\alpha$ and $\Vdash_{\mathcal{F}}^{\nbI}\alpha\rightarrow\beta$ then $\Vdash_{\mathcal{F}}^{\nbI}\beta$: we would like to prove that, for any formula $\alpha$ of $\nbI$, we also have 
\[\Vdash_{\mathcal{F}}^{\nbI}\alpha\vee\neg\alpha,\]
implying that "$\Vdash_{\mathcal{F}}^{\nbI}$" models the axiom schemata and rules of inference of $\nbI$.

So, for any Fidel structure $\mathcal{A}$ for $\nbI$ and valuation $\nu$ over $\mathcal{A}$, we have that $\nu(\alpha\vee\neg \alpha)=\nu(\alpha)\vee\nu(\neg \alpha)$; now, $\nu(\neg \alpha)\in\neg\nu(\alpha)$, and therefore, from the conditions demanded of a Fidel structure for $\nbI$, $\nu(\alpha)\vee\nu(\neg \alpha)=\top$, what finishes the desired proof.

\begin{theorem}
Given formulas $\Gamma\cup\{\varphi\}$ of $\nbI$, if $\Gamma\vdash_{\nbI}\varphi$ then $\Gamma\Vdash_{\mathcal{F}}^{\nbI}\varphi$.
\end{theorem}

To prove the reciprocal result, we define the equivalence relation between formulas of $\nbI$ such that, for a fixed set of formulas $\Gamma$, \label{equivnbI}$\alpha\equiv_{\Gamma}^{\nbI}\beta$ if and only if $\Gamma\vdash_{\nbI}\alpha\rightarrow\beta$ and $\Gamma\vdash_{\nbI}\beta\rightarrow\alpha$; this relation is also a congruence with respect to the connectives in $\{\vee, \wedge,\rightarrow\}$, allowing us to make $A^{\nbI}_{\Gamma}=F(\Sigma_{\nbI},\mathcal{V})/\equiv_{\Gamma}^{\nbI}$ into the universe of a Boolean algebra. By defining, for classes of formulas $[\alpha]$ and $[\beta]$,
\[\neg[\alpha]=\{[\neg \varphi]\ :\  \varphi\in[\alpha]\}\]
and
\[[\alpha]\uparrow[\beta]=\{[\varphi\uparrow\psi]\ :\  \varphi\in[\alpha], \psi\in[\beta]\},\]
we can prove \label{AnbIGamma}$\mathcal{A}^{\nbI}_{\Gamma}=(A^{\nbI}_{\Gamma},\{\sigma_{\mathcal{A}}\}_{\sigma\in\Sigma_{\nbI}^{\textbf{CPL}}})$ is a Fidel structure, presented as a $\Sigma_{\nbI}^{\textbf{CPL}}$-multialgebra, for $\nbI$, which we shall call the Lindenbaum-Tarski Fidel structure of $\nbI$ associated to $\Gamma$.

The proof that, for $[\gamma]\in [\alpha]\uparrow[\beta]$, we have $[\alpha]\wedge([\beta]\wedge[\gamma])=\bot$, and that $[\alpha]\uparrow[\beta]=[\beta]\uparrow[\alpha]$ are exactly the same as those for the same facts in $\bI$ (found in Section \ref{Completeness fo Fidel structures for bI}), so we will rather focus on proving that, for any class of formulas $[\alpha]$ and any $[\beta]\in \neg[\alpha]$, $[\alpha]\vee[\beta]=\top$: by definition, there exists $\varphi\in\alpha$ such that $[\beta]=[\neg\varphi]$.

Then $[\alpha]\vee[\beta]=[\alpha\vee\neg \varphi]$: very clearly $\Gamma\vdash_{\nbI}\alpha\vee\neg \varphi\rightarrow\top_{\alpha}$, remaining for us to show that $\Gamma\vdash_{\nbI}\top_{\alpha}\rightarrow\alpha\vee\neg \varphi$ or, better yet, that $\Gamma\vdash_{\nbI}\alpha\vee\neg \varphi$, by one application of the deduction meta-theorem and the fact that a top element is always a tautology.

We know $\varphi\vee\neg \varphi$ is an instance of $\textbf{Ax\: 11}^{*}$ and $\Gamma\vdash_{\nbI}\varphi\rightarrow\alpha$:
\begin{enumerate}
\item $\Gamma,\varphi\vdash_{\nbI}\alpha$ by the deduction meta-theorem, and from the instance $\alpha\rightarrow\alpha\vee\neg \varphi$ of $\textbf{Ax\: 6}$ and one application of Modus Ponens we get that $\Gamma,\varphi\vdash_{\nbI}\alpha\vee\neg \varphi$;
\item $\Gamma,\neg \varphi\vdash_{\nbI}\neg \varphi$, and from the instance $\neg \varphi\rightarrow\alpha\vee\neg \varphi$ of $\textbf{Ax\: 7}$ and the deduction meta-theorem we get that $\Gamma,\neg \varphi\vdash_{\nbI}\alpha\vee\neg \varphi$;
\end{enumerate}
from this, we get that $\Gamma, \varphi\vee\neg \varphi\vdash_{\nbI}\alpha\vee\neg \varphi$ by a proof by cases, and since $\varphi\vee\neg \varphi$ is an instance of an axiom, $\Gamma\vdash_{\nbI}\alpha\vee\neg \varphi$, as we wanted to prove.

\begin{theorem}
Given formulas $\Gamma\cup\{\varphi\}$ of $\nbI$, if $\Gamma\Vdash_{\mathcal{F}}^{\nbI}\varphi$ then $\Gamma\vdash_{\nbI}\varphi$.
\end{theorem}


\subsection{Decision method}

We take the Boolean algebra $\textbf{2}$ once again and define the unary multioperation $\neg:\{0,1\}\rightarrow\mathcal{P}(\{0,1\})\setminus\{\emptyset\}$ trough $\neg 0=\{1\}$ and $\neg 1=\{0,1\}$,

\begin{figure}[H]
\centering
\begin{tabular}{|l|r|}
\hline
& $\neg $\\ \hline
$0$ & $\{1\}$\\ \hline
$1$ & $\{0,1\}$\\\hline
\end{tabular}
\caption*{Table for our Paraconsistent Negation}
\end{figure}

and, as we did in Section \ref{Decision Method for bI}, we define $1\uparrow 1=\{0\}$ and $x\uparrow y=\{0,1\}$ for any other values of $x$ and $y$. By adding the previous operations to $\textbf{2}$, we transform it into a $\Sigma_{\nbI}^{\textbf{CPL}}$-multialgebra\label{2nbI} $\textbf{2}_{\nbI}$ with a finite universe $\{0,1\}$. It is easy to prove that:
\begin{enumerate}
\item $(\{0,1\}, \{\sigma_{\textbf{2}_{\nbI}}\}_{\sigma\in\Sigma^{\textbf{CPL}}})$ is a Boolean algebra, since it equals $\textbf{2}$;
\item for any $x,y\in\{0,1\}$ and $z\in x\uparrow y$, $x\wedge(y\wedge z)=0$;
\item for any $x, y\in\{0,1\}$, $x\uparrow y=y\uparrow x$.
\end{enumerate}

It remains to be shown, to prove that $\textbf{2}_{\nbI}$ is a Fidel structure presented as a $\Sigma_{\nbI}^{\textbf{CPL}}$-multialgebra for $\nbI$, that for any $x\in \{0,1\}$ and $y\in \neg x$, $x\vee y=1$. If $x=0$, we must have $y=1$, and then $x\vee y=0\vee 1=1$; if $x=1$, we have that, either $y=0$, when $x\vee y=1\vee 0=1$, or $y=1$, when $x\vee y=1\vee1=1$, so we are done

\begin{theorem}
$\nu:F(\Sigma_{\nbI},\mathcal{V})\rightarrow\{0,1\}$ is a bivaluation for $\nbI$ if, and only if, it is a $\Sigma_{\nbI}$-homomorphism from $\textbf{F}(\Sigma_{\nbI}, \mathcal{V})$ to $\textbf{2}_{\nbI}$ which lies in $\mathcal{F}_{\textbf{2}_{\nbI}}$.
\end{theorem}

As expected, we will denote the restricted Nmatrix $(\textbf{2}_{\nbI}, \{1\}, \mathcal{F}_{\textbf{2}_{\nbI}})$ simply by $\mathbb{2}_{\nbI}$.

\begin{theorem}
Given formulas $\Gamma\cup\{\varphi\}$ of $\nbI$, $\Gamma\vDash_{\nbI}\varphi$ if and only if $\Gamma\vDash_{\mathbb{2}_{\nbI}}\varphi$.
\end{theorem}

So, take the Nmatrix $(\textbf{2}_{\nbI}, \{1\})$ subjacent to the RNmatrix $(\textbf{2}_{\nbI}, \{1\}, \mathcal{F}_{\textbf{2}_{\nbI}})$: if one writes the row-branching truth table for the said Nmatrix to test a formula $\varphi$ of $\nbI$ and simply erases the rows where $\alpha\uparrow\beta$ and $\beta\uparrow\alpha$ are given different values, we get a row-branching, row-eliminating truth table that decides the validity, in a finite number of steps, of $\varphi$ in $\nbI$. Essentially, this is the same decision method as the one presented for $\bI$ in Section \ref{Decision Method for bI}, with an extra multi-operation standing for negation.

\subsection{Another decision method}

Given the similitude between $\bI$ and $\nbI$, one should hope that a tableau calculus for $\nbI$, that we will name $\mathbb{T}_{\nbI}$, could be easily obtained from $\mathbb{T}_{\bI}$ of Section \ref{Tableaux for bI}: this is, in fact, the case. It is sufficient to add one new rule governing the case in which $\neg\varphi$ is false, 
\[\frac{\textsf{0}(\neg\varphi)}{\textsf{1}(\varphi)}\]
and change the conditions for a branch to be complete: a branch $\theta$ is complete, now in $\mathbb{T}_{\nbI}$, when it contains, for every $\textsf{L}(\gamma)$ in $\theta$ neither of the form $\textsf{0}(\varphi\uparrow\psi)$ or $\textsf{1}(\neg\varphi)$ and with $\gamma$ not a variable, all the labeled formulas of one of the branches in the rule headed by $\textsf{L}(\gamma)$.


\section{Collapsing axioms}


\subsection{$\textbf{ci}^{\uparrow}$}\label{ciuparrow}

Now that we have in our signature a paraconsistent negation, we are capable of doing to the axioms schemata $\textbf{ci}$ and $\textbf{cl}$ much the same we did to $\textbf{ciw}$ when we transformed it into $\textbf{ciw}^{\uparrow}$: consider\label{ciuparrow}
\[\tag{$\textbf{ci}^{\uparrow}$} \neg(\alpha\uparrow\beta)\rightarrow(\alpha\wedge\beta)\]
and the logic $\nbI\textbf{ci}^{\uparrow}$ obtained from $\nbI$ by addition of $\textbf{ci}^{\uparrow}$; we shall prove, once again, that $\nbI\textbf{ci}^{\uparrow}$ collapses back to $\textbf{CPL}$, and $\alpha\uparrow\beta$ is equivalent to $\alpha\wedge\beta\rightarrow\bot_{\alpha\beta}$.

It is clear how $\alpha\uparrow\beta\vdash_{\nbI\textbf{ci}^{\uparrow}}\alpha\wedge\beta\rightarrow\bot_{\alpha\beta}$, since such a deduction could be made in $\bI$ itself, so let us concentrate in proving that $\alpha\wedge\beta\rightarrow\bot_{\alpha\beta}\vdash_{\nbI\textbf{ci}^{\uparrow}}\alpha\uparrow\beta$.

It is also clear that 
\[\alpha\uparrow\beta, \alpha\wedge\beta\rightarrow\bot_{\alpha\beta}\vdash_{\nbI\textbf{ci}^{\uparrow}}\alpha\uparrow\beta,\]
but from the instance $\neg(\alpha\uparrow\beta)\rightarrow(\alpha\wedge\beta)$ of $\textbf{ci}^{\uparrow}$ plus Modus Ponens we can similarly get 
\[\neg(\alpha\uparrow\beta), \alpha\wedge\beta\rightarrow\bot_{\alpha\beta}\vdash_{\nbI\textbf{ci}^{\uparrow}}\alpha\uparrow\beta,\]
and using a proof by cases we get
\[(\alpha\uparrow\beta)\vee\neg(\alpha\uparrow\beta), \alpha\wedge\beta\rightarrow\bot_{\alpha\beta}\vdash_{\nbI\textbf{ci}^{\uparrow}}\alpha\uparrow\beta.\]
Since $(\alpha\uparrow\beta)\vee\neg(\alpha\uparrow\beta)$ is an instance of $\textbf{Ax\: 11}^{*}$, we get the desired result.

But we can prove an even nicer result, that implies the previous: an instance of $\textbf{ci}^{\uparrow}$ actually implies its corresponding instance of $\textbf{ciw}^{\uparrow}$ in $\nbI$, that is,
\[\neg(\alpha\uparrow\beta)\rightarrow(\alpha\wedge\beta)\vdash_{\nbI}(\alpha\uparrow\beta)\vee(\alpha\wedge\beta).\]
We have that
\[\neg(\alpha\uparrow\beta),\neg(\alpha\uparrow\beta)\rightarrow(\alpha\wedge\beta)\vdash_{\nbI}(\alpha\uparrow\beta)\vee(\alpha\wedge\beta)\]
by an application of Modus Ponens and remembering $(\alpha\wedge\beta)\rightarrow[(\alpha\uparrow\beta)\vee(\alpha\wedge\beta)]$ is an instance of $\textbf{Ax\: 7}$; furthermore, we trivially find
\[\alpha\uparrow\beta,\neg(\alpha\uparrow\beta)\rightarrow(\alpha\wedge\beta)\vdash_{\nbI}(\alpha\uparrow\beta)\vee(\alpha\wedge\beta)\]
by remembering $(\alpha\uparrow\beta)\rightarrow[(\alpha\uparrow\beta)\vee(\alpha\wedge\beta)]$ is an instance of $\textbf{Ax\: 6}$. By a proof by cases, 
\[(\alpha\uparrow\beta)\vee\neg(\alpha\uparrow\beta),\neg(\alpha\uparrow\beta)\rightarrow(\alpha\wedge\beta)\vdash_{\nbI}(\alpha\uparrow\beta)\vee(\alpha\wedge\beta),\]
and since $(\alpha\uparrow\beta)\vee\neg(\alpha\uparrow\beta)$ is an instance of $\textbf{Ax\: 11}^{*}$, we derive the aforementioned result.


\subsection{$\textbf{cl}^{\uparrow}$}\label{cluparrow}

As one should expect, the axiom schema\label{cluparrow}
\[\tag{$\textbf{cl}^{\uparrow}$} \neg(\alpha\wedge\beta)\rightarrow(\alpha\uparrow\beta)\]
will collapse its corresponding logic $\nbI\textbf{cl}^{\uparrow}$ back to $\textbf{CPL}$, with $\alpha\uparrow\beta$ standing for $\alpha\wedge\beta\rightarrow\bot_{\alpha\beta}$.

To see that $\alpha\wedge\beta\rightarrow\bot_{\alpha\beta}\vdash_{\nbI\textbf{cl}^{\uparrow}}\alpha\uparrow\beta$, we begin by noticing
\[\alpha\wedge\beta, \alpha\wedge\beta\rightarrow\bot_{\alpha\beta}\vdash_{\nbI\textbf{cl}^{\uparrow}}\alpha\uparrow\beta,\]
by use of Modus Ponens and the fact $\bot_{\alpha\beta}$ behaves like a bottom element; at the same time,
\[\neg(\alpha\wedge\beta), \alpha\wedge\beta\rightarrow\bot_{\alpha\beta}\vdash_{\nbI\textbf{cl}^{\uparrow}}\alpha\uparrow\beta,\]
given the instance $ \neg(\alpha\wedge\beta)\rightarrow(\alpha\uparrow\beta)$ of $\textbf{cl}^{\uparrow}$. With a proof by cases,
\[(\alpha\wedge\beta)\vee\neg(\alpha\wedge\beta), \alpha\wedge\beta\rightarrow\bot_{\alpha\beta}\vdash_{\nbI\textbf{cl}^{\uparrow}}\alpha\uparrow\beta,\]
and since $(\alpha\wedge\beta)\vee\neg(\alpha\wedge\beta)$ is an instance of $\textbf{Ax\: 11}^{*}$ we obtain the desired result.

However, we can again prove a stronger result: an instance of $\textbf{cl}^{\uparrow}$ can prove, in $\nbI$, the corresponding instance of $\textbf{ciw}^{\uparrow}$, what proves, as a special case, the previous collapse. To summarize, we wish to prove
\[\neg(\alpha\wedge\beta)\rightarrow(\alpha\uparrow\beta)\vdash_{\nbI}(\alpha\uparrow\beta)\vee(\alpha\wedge\beta).\]
This is pretty simple:
\[\alpha\wedge\beta, \neg(\alpha\wedge\beta)\rightarrow(\alpha\uparrow\beta)\vdash_{\nbI}(\alpha\uparrow\beta)\vee(\alpha\wedge\beta),\]
by the instance $(\alpha\wedge\beta)\rightarrow[(\alpha\uparrow\beta)\vee(\alpha\wedge\beta)]$ of $\textbf{Ax\: 7}$, and
\[\neg(\alpha\wedge\beta), \neg(\alpha\wedge\beta)\rightarrow(\alpha\uparrow\beta)\vdash_{\nbI}(\alpha\uparrow\beta)\vee(\alpha\wedge\beta)\]
by Modus Ponens and the instance $(\alpha\uparrow\beta)\rightarrow[(\alpha\uparrow\beta)\vee(\alpha\wedge\beta)]$ of $\textbf{Ax\: 6}$; with a proof by cases, 
\[(\alpha\wedge\beta)\vee\neg(\alpha\wedge\beta), \neg(\alpha\wedge\beta)\rightarrow(\alpha\uparrow\beta)\vdash_{\nbI}(\alpha\uparrow\beta)\vee(\alpha\wedge\beta),\]
and since $(\alpha\wedge\beta)\vee\neg(\alpha\wedge\beta)$ is an instance of $\textbf{Ax\: 11}^{*}$, the proof is done.


\section{The logics $\nbIciw$, $\nbIci$ and $\nbIcl$}

One notices that one of the basic axiom schema, $\textbf{Ip}$, of $\bI$ is
\[(\alpha\uparrow\beta)\rightarrow(\alpha\rightarrow(\beta\rightarrow\gamma)),\]
which can be seem as the basic axiom of $\textbf{mbC}$, $\textbf{bc1}$,
\[\circ\alpha\rightarrow(\alpha\rightarrow(\neg\alpha\rightarrow\beta)),\]
with $\neg\alpha$ replaced by $\beta$ and $\circ\alpha$, that is, that $\alpha$ is consistent, replaced by $\alpha\uparrow\beta$, that is, that $\alpha$ is incompatible with $\beta$.

An attempt to make something similar to the axiom $\textbf{ciw}$,
\[\circ\alpha\vee(\alpha\wedge\neg\alpha),\]
is not successful, as the axiom $\textbf{ciw}^{\uparrow}$
\[(\alpha\uparrow\beta)\vee(\alpha\wedge\beta)\]
collapses the incompatibility operator back to its classical interpretation. However, in $\nbI$, where we find at our disposal a non-classical negation, we can take an adaptation of the axiom $\textbf{ciw}$ to the language $\Sigma_{\nbI}$, instead of the generalization $\textbf{ciw}^{\uparrow}$, and study the resulting logic. We will do the same with the axioms $\textbf{ci}$ and $\textbf{cl}$, producing three logics that mix paraconsistency and incompatibility.

So, over the signature $\Sigma_{\nbI}$, we consider the systems:
\begin{enumerate}
\item $\nbIciw$\label{nbIciw}, obtained from $\nbI$ by addition of the axiom schema\label{ciw*}
\[\tag{$\textbf{ciw}^{*}$} (\alpha\uparrow\neg\alpha)\vee(\alpha\wedge\neg\alpha);\]
\item $\nbIci$\label{nbIci}, obtained from $\nbI$ by adding the axiom schema\label{ci*}
\[\tag{$\textbf{ci}^{*}$}\neg(\alpha\uparrow\neg\alpha)\rightarrow(\alpha\wedge\neg\alpha);\]
\item $\nbIcl$\label{nbIcl}, obtained from $\nbI$ by adding\label{cl*}
\[\tag{$\textbf{cl}^{*}$}\neg(\alpha\wedge\neg\alpha)\rightarrow(\alpha\uparrow\neg\alpha).\]
\end{enumerate}


\subsection{Bivaluations}

\begin{definition}\label{definition of Bivaluation for nbIciw, nbIci, ...}
A \index{Bivaluation for $\nbIciw$}\index{Bivaluation for $\nbIci$}\index{Bivaluation for $\nbIcl$}bivaluation for $\mathcal{L}\in\{\nbIciw, \nbIci, \nbIcl\}$ is a valuation for $\nbI$ satisfying in addition that:
\begin{enumerate}
\item $\nu(\alpha\uparrow\neg\alpha)=1$ if and only if $\nu(\alpha)=0$ or $\nu(\neg\alpha)=0$;
\item \begin{enumerate}
\item if $\mathcal{L}=\nbIci$ and $\nu(\neg(\alpha\uparrow\neg\alpha))=1$, then $\nu(\alpha)=\nu(\neg\alpha)=1$;
\item if $\mathcal{L}=\nbIcl$ and $\nu(\neg(\alpha\wedge\neg\alpha))=1$, then $\nu(\alpha\uparrow\neg\alpha)=1$.
\end{enumerate}
\end{enumerate}
\end{definition}

As usual, given a set of formulas $\Gamma\cup\{\varphi\}$ on the signature $\Sigma_{\nbI}$, we say $\Gamma$ proves $\varphi$ according to bivaluations for $\mathcal{L}\in\{\nbIciw, \nbIci, \nbIcl\}$, and write \label{vDashnbIciw}\label{vDashnbIci}\label{vDashnbIcl}$\Gamma\vDash_{\mathcal{L}}\varphi$, if, for every bivaluation $\nu$ for $\mathcal{L}$ such that $\nu(\Gamma)\subseteq\{1\}$, one has $\nu(\varphi)=1$.

Since bivaluations for $\mathcal{L}\in\{\nbIciw, \nbIci, \nbIcl\}$ are bivaluations for $\nbI$ satisfying some additional property, it is easy to see that "$\vDash_{\mathcal{L}}$" models all those axiom schemata of $\nbI$, plus Modus Ponens, its only inference rule; we wish now to prove that "$\vDash_{\nbIciw}$" also models $\textbf{ciw}^{*}$, and an analogous result holds for the other logics.

\begin{enumerate}
\item Take an instance $(\alpha\uparrow\neg\alpha)\vee(\alpha\wedge\neg\alpha)$ of $\textbf{ciw}^{*}$ and a bivaluation $\nu$ for $\nbIciw$, and we see that 
\[\nu((\alpha\uparrow\neg\alpha)\vee(\alpha\wedge\neg\alpha))=0\]
if and only if $\nu(\alpha\uparrow\neg\alpha)=0$ and $\nu(\alpha\wedge\neg\alpha)=0$; the first of these equalities holds if and only $\nu(\alpha)=1$ and $\nu(\neg\alpha)=1$, what would imply that $\nu(\alpha\wedge\neg\alpha)=1$, reaching a contradiction.

\item Take a bivaluation $\nu$ for $\nbIci$ and suppose 
\[\nu(\neg(\alpha\uparrow\neg\alpha)\rightarrow(\alpha\wedge\neg\alpha))=0,\]
what happens if and only if $\nu(\neg(\alpha\uparrow\neg\alpha))=1$ but $\nu(\alpha\wedge\neg\alpha)=0$; from the first equality, $\nu(\alpha)=\nu(\neg\alpha)=1$, what again contradicts $\nu(\alpha\wedge\neg\alpha)=0$.

\item Finally, we take a bivaluation $\nu$ for $\nbIcl$ and assume it is possible to have
\[\nu(\neg(\alpha\wedge\neg\alpha)\rightarrow(\alpha\uparrow\neg\alpha))=0.\]
This is equivalent to having $\nu(\neg(\alpha\wedge\neg\alpha))=1$ and $\nu(\alpha\uparrow\neg\alpha)=0$, but the first equality implies $\nu(\alpha\uparrow\neg\alpha)=1$, which is evidently contradictory.
\end{enumerate}

\begin{theorem}
Given formulas $\Gamma\cup\{\varphi\}$ over the signature $\Sigma_{\nbI}$ and $\mathcal{L}\in\{\nbIciw, \nbIci, \nbIcl\}$, if $\Gamma\vdash_{\mathcal{L}}\varphi$ then $\Gamma\vDash_{\mathcal{L}}\varphi$.
\end{theorem}

To prove the reciprocal, take a set of formulas $\Gamma$ maximal with respect to not proving $\varphi$ in $\mathcal{L}\in\{\nbIciw, \nbIci, \nbIcl\}$ which is also closed (again, with respect to $\mathcal{L}$) and non-trivial. We define $\nu:F(\Sigma_{\nbI},\mathcal{V})\rightarrow\{0,1\}$ such that $\nu(\gamma)=1$ if and only if $\gamma\in\Gamma$, and wish to prove that $\nu$ is a bivaluation for $\mathcal{L}$.

It is clear $\Gamma$ does not prove $\varphi$ in $\nbI$, since $\mathcal{L}$ extends $\nbI$, and also that $\Gamma$ is closed in $\nbI$, given it is closed according to $\mathcal{L}$ and this logic is stronger than the previous; this implies $\nu$ is at least a bivaluation for $\nbI$, according to the definition found in Section \ref{Bivaluations for nbI}.

\begin{enumerate}
\item For any of the logics $\mathcal{L}\in\{\nbIciw, \nbIci, \nbIcl\}$, if $\nu(\alpha\uparrow\neg\alpha)=1$, suppose $\nu(\alpha)=1$ and $\nu(\neg\alpha)=1$: by definition of $\nu$, this means $\alpha\uparrow\neg\alpha, \alpha, \neg\alpha\in\Gamma$, and since $\Gamma$ is closed and 
\[(\alpha\uparrow\neg\alpha)\rightarrow(\alpha\rightarrow(\neg\alpha\rightarrow\varphi))\]
is an instance of $\textbf{Ip}$, we get that $\varphi\in\Gamma$, what is a contradiction, We must then have either $\nu(\alpha)=0$ or $\nu(\neg\alpha)=0$.

Suppose, reciprocally, that $\nu(\alpha)=0$ or $\nu(\neg\alpha)=0$, and for a proof by contradiction, let us take $\nu(\alpha\uparrow\neg\alpha)=0$, meaning $\alpha\uparrow\neg\alpha\notin\Gamma$ and either $\alpha\notin\Gamma$ or $\neg\alpha\notin\Gamma$, implying $\alpha\wedge\neg\alpha\notin\Gamma$. By the maximality of $\Gamma$, one has $\Gamma, \alpha\uparrow\neg\alpha\vdash_{\mathcal{L}}\varphi$ and $\Gamma,\alpha\wedge\neg\alpha\vdash_{\mathcal{L}}\varphi$, and by a proof by cases
\[\Gamma, (\alpha\uparrow\neg\alpha)\vee(\alpha\wedge\neg\alpha)\vdash_{\mathcal{L}}\varphi.\]
Since $(\alpha\uparrow\neg\alpha)\vee(\alpha\wedge\neg\alpha)$ is an instance of $\textbf{ciw}^{*}$, if $\mathcal{L}=\nbIciw$ we get $\Gamma\vdash_{\nbIciw}\varphi$, a contradiction; for the other two logics, given that, from Sections \ref{ciuparrow} and \ref{cluparrow},
\[\neg(\alpha\uparrow\beta)\rightarrow(\alpha\wedge\beta)\vdash_{\nbI}(\alpha\uparrow\beta)\vee(\alpha\wedge\beta)\]
and
\[\neg(\alpha\wedge\beta)\rightarrow(\alpha\uparrow\beta)\vdash_{\nbI}(\alpha\uparrow\beta)\vee(\alpha\wedge\beta),\]
by replacing $\beta$ with $\neg\alpha$ we obtain that we still have $\Gamma\vdash_{\mathcal{L}}\varphi$, so still a contradiction. This means we must have $\nu(\alpha\uparrow\neg\alpha)=1$.

\item If $\mathcal{L}=\nbIci$ and $\nu(\neg(\alpha\uparrow\neg\alpha))=1$, this means $\neg(\alpha\uparrow\neg\alpha)\in\Gamma$; since $\neg(\alpha\uparrow\neg\alpha)\rightarrow(\alpha\wedge\neg\alpha)$ is an instance of $\textbf{ci}^{*}$ and $\Gamma$ is closed, $\alpha\wedge\neg\alpha\in\Gamma$, implying that $\alpha, \neg\alpha\in\Gamma$ and therefore $\nu(\alpha)=1$ and $\nu(\neg\alpha)=1$.

\item If $\mathcal{L}=\nbIcl$ and $\nu(\neg(\alpha\wedge\neg\alpha))=1$, we have $\neg(\alpha\wedge\neg\alpha)\in\Gamma$; since $\neg(\alpha\uparrow\neg\alpha)\rightarrow(\alpha\wedge\neg\alpha)$ is an instance of $\textbf{cl}^{*}$ and $\Gamma$ is closed, $\alpha\uparrow\neg\alpha\in\Gamma$, and therefore $\nu(\alpha\uparrow\neg\alpha)=1$.
\end{enumerate}

\begin{theorem}
Given formulas $\Gamma\cup\{\varphi\}$ over the signature $\Sigma_{\nbI}$ and $\mathcal{L}\in\{\nbIciw, \nbIci, \nbIcl\}$, if $\Gamma\vDash_{\mathcal{L}}\varphi$ then $\Gamma\vdash_{\mathcal{L}}\varphi$.
\end{theorem}

In the spirit of the arguments just used, since
\[\neg(\alpha\uparrow\beta)\rightarrow(\alpha\wedge\beta)\vdash_{\nbI}(\alpha\uparrow\beta)\vee(\alpha\wedge\beta)\]
and
\[\neg(\alpha\wedge\beta)\rightarrow(\alpha\uparrow\beta)\vdash_{\nbI}(\alpha\uparrow\beta)\vee(\alpha\wedge\beta),\]
by replacing $\beta$ with $\neg\alpha$ we discover both $\nbIci$ and $\nbIcl$ are extensions of $\nbIciw$, and with the aid of bivaluations, we can prove even more.

\begin{proposition}
\begin{enumerate}
\item $\nbIciw$ can not prove $\textbf{ci}^{*}$, and therefore $\nbIci$ is strictly stronger than $\nbIciw$;
\item $\nbIciw$ can not prove $\textbf{cl}^{*}$, and therefore $\nbIcl$ is strictly stronger than $\nbIciw$.
\end{enumerate}
\end{proposition}

\begin{proof}
\begin{enumerate}
\item Take a bivaluation $\nu$ for $\nbIciw$ such that $\nu(\alpha)=0$, $\nu(\neg\alpha)=0$ (and therefore $\nu(\alpha\wedge\neg\alpha)=0$), $\nu(\alpha\uparrow\neg\alpha)=1$ and $\nu(\neg(\alpha\uparrow\neg\alpha))=1$: then
\[\nu(\neg(\alpha\uparrow\neg\alpha)\rightarrow(\alpha\wedge\neg\alpha))=0.\]
\item Take a bivaluation $\nu$ for $\nbIciw$ such that $\nu(\alpha)=1$, $\nu(\neg\alpha)=1$ (and therefore $\nu(\alpha\wedge\neg\alpha)=1$), $\nu(\alpha\uparrow\neg\alpha)=0$ and $\nu(\neg(\alpha\wedge\neg\alpha))=1$: then 
\[\nu(\neg(\alpha\wedge\neg\alpha)\rightarrow(\alpha\uparrow\neg\alpha))=0.\]
\end{enumerate}
\end{proof}

\begin{theorem}
\begin{enumerate}
\item The formula $(\alpha\uparrow\neg\alpha)\uparrow\neg(\alpha\uparrow\neg\alpha)$ is a tautology of $\nbIci$;
\item $\nbIci$ is obtained from $\nbIciw$ by adding the axiom schema\label{cc*}
\[\tag{$\textbf{cc}^{*}$} (\alpha\uparrow\neg\alpha)\uparrow\neg(\alpha\uparrow\neg\alpha).\]
\end{enumerate}
\end{theorem}

\begin{proof}
\begin{enumerate}
\item It is clear that
\[(\alpha\uparrow\neg\alpha)\uparrow\neg(\alpha\uparrow\neg\alpha)\vdash_{\nbIci}(\alpha\uparrow\neg\alpha)\uparrow\neg(\alpha\uparrow\neg\alpha),\]
but from the instance
\[\neg[(\alpha\uparrow\neg\alpha)\uparrow\neg(\alpha\uparrow\neg\alpha)]\rightarrow[(\alpha\uparrow\neg\alpha)\wedge\neg(\alpha\uparrow\neg\alpha)]\]
of $\textbf{ci}^{*}$ we see that
\[\neg[(\alpha\uparrow\neg\alpha)\uparrow\neg(\alpha\uparrow\neg\alpha)]\vdash_{\nbIci}(\alpha\uparrow\neg\alpha)\wedge\neg(\alpha\uparrow\neg\alpha).\]
$(\alpha\uparrow\neg\alpha)\wedge\neg(\alpha\uparrow\neg\alpha)$ easily deduces both $\alpha\uparrow\neg\alpha$ and $\neg(\alpha\uparrow\neg\alpha)$, and from this last formula and the instance $\neg(\alpha\uparrow\neg\alpha)\rightarrow(\alpha\wedge\neg\alpha)$ of $\textbf{ci}^{*}$, we get $\neg[(\alpha\uparrow\neg\alpha)\uparrow\neg(\alpha\uparrow\neg\alpha)]$ deduces $\alpha\uparrow\neg\alpha$, $\alpha$ and $\neg\alpha$: from the instance of $\textbf{Ip}$
\[(\alpha\uparrow\neg\alpha)\rightarrow\big[\alpha\rightarrow\big[\neg\alpha\rightarrow\big((\alpha\uparrow\neg\alpha)\uparrow\neg(\alpha\uparrow\neg\alpha)\big)\big]\big],\]
this means 
\[\neg[(\alpha\uparrow\neg\alpha)\uparrow\neg(\alpha\uparrow\neg\alpha)]\vdash_{\nbIci}(\alpha\uparrow\neg\alpha)\uparrow\neg(\alpha\uparrow\neg\alpha),\]
and by a proof by cases 
\[[(\alpha\uparrow\neg\alpha)\uparrow\neg(\alpha\uparrow\neg\alpha)]\vee\neg[(\alpha\uparrow\neg\alpha)\uparrow\neg(\alpha\uparrow\neg\alpha)]\vdash_{\nbIci}(\alpha\uparrow\neg\alpha)\uparrow\neg(\alpha\uparrow\neg\alpha).\]
Since the antecedent in this argument is an instance of $\textbf{Ax\: 11}^{*}$, we obtain that $(\alpha\uparrow\neg\alpha)\uparrow\neg(\alpha\uparrow\neg\alpha)$ is a tautology.

\item So we have that $\nbIci$ proves any instance of $\textbf{ciw}^{*}$ and $\textbf{cc}^{*}$, being therefore stronger than the logic obtained from $\nbIciw$ by addition of $\textbf{cc}^{*}$, which we will briefly call $\nbIciw+$. To prove that $\nbIciw+$ is as strong as $\nbIci$, and therefore that both are equal, we only need to prove that any instance of $\textbf{ci}^{*}$ is a tautology in $\nbIciw+$, or what is equivalent, that
\[\neg(\alpha\uparrow\neg\alpha)\vdash_{\nbIciw+}\alpha\wedge\neg\alpha.\]
Obviously
\[\alpha\wedge\neg\alpha, \neg(\alpha\uparrow\neg\alpha)\vdash_{\nbIciw+}\alpha\wedge\neg\alpha;\]
now, from the instance
\[[(\alpha\uparrow\neg\alpha)\uparrow\neg(\alpha\uparrow\neg\alpha)]\rightarrow\big[(\alpha\uparrow\neg\alpha)\rightarrow\big[\neg(\alpha\uparrow\neg\alpha)\rightarrow(\alpha\wedge\neg\alpha)\big]\big]\]
of $\textbf{Ip}$ and the fact $(\alpha\uparrow\neg\alpha)\uparrow\neg(\alpha\uparrow\neg\alpha)$ is an instance of $\textbf{cc}^{*}$, we can see by applying the deduction meta-theorem as needed that
\[\alpha\uparrow\neg\alpha, \neg(\alpha\uparrow\neg\alpha)\vdash_{\nbIciw+}\alpha\wedge\neg\alpha.\]

By a proof by cases,
\[(\alpha\uparrow\neg\alpha)\vee(\alpha\wedge\neg\alpha), \neg(\alpha\uparrow\neg\alpha)\vdash_{\nbIciw+}\alpha\wedge\neg\alpha,\]
and since $(\alpha\uparrow\neg\alpha)\vee(\alpha\wedge\neg\alpha)$ is an instance of $\textbf{ciw}^{*}$, we discover that $\nbIci$ is $\nbIciw+$, as we wanted to show.
\end{enumerate}
\end{proof}


\subsection{Fidel structures}

A \index{Fidel structure for $\nbIciw$}\index{Fidel structure for $\nbIci$}\index{Fidel structure for $\nbIcl$}Fidel structure, presented as a $\Sigma_{\nbI}^{\textbf{CPL}}$-multialgebra, for $\mathcal{L}\in\{\nbIciw, \nbIci, \nbIcl\}$ is any $\Sigma_{\nbI}^{\textbf{CPL}}$-multialgebra $\mathcal{A}=(A,\{\sigma_{\mathcal{A}}\}_{\sigma\in\Sigma_{\nbI}^{\textbf{CPL}}})$ that is a Fidel structure for $\nbI$ and satisfies, additionally:
\begin{enumerate}
\item if $b\in\neg a$, ${\sim}(a\wedge b)\in a\uparrow b$;
\item\begin{enumerate}
\item if $\mathcal{L}=\nbIci$ and $b\in\neg a$, $a\wedge b\in\neg{\sim}(a\wedge b)$;
\item if $\mathcal{L}=\nbIcl$ and $b\in\neg a$, ${\sim}(a\wedge b)\in\neg(a\wedge b)$.
\end{enumerate}
\end{enumerate}

A valuation for a Fidel structure $\mathcal{A}$, presented as a $\Sigma_{\nbI}^{\textbf{CPL}}$-multialgebra, for $\mathcal{L}\in\{\nbIciw, \nbIci, \nbIcl\}$, will still be any $\Sigma_{\nbI}$-homomorphism $\nu:\textbf{F}(\Sigma_{\nbI},\mathcal{V})\rightarrow\mathcal{A}$; and now we will consider the restricted Nmatrices $(\mathcal{A},\{\top\},\mathcal{F}_{\mathcal{A}}^{\mathcal{L}})$, where $\mathcal{A}$ is a Fidel structure for $\mathcal{L}$ and $\mathcal{F}_{\mathcal{A}}^{\mathcal{L}}$ is the set of valuations $\nu:\textbf{F}(\Sigma_{\nbI},\mathcal{V})\rightarrow\mathcal{A}$ such that, for any formulas $\alpha$ and $\beta$ over the signature $\Sigma_{\nbI}$:
\begin{enumerate}
\item $\nu(\alpha\uparrow\beta)=\nu(\beta\uparrow\alpha)$;
\item $\nu(\alpha\uparrow\neg\alpha)={\sim}(\nu(\alpha)\wedge\nu(\neg\alpha))$;
\item \begin{enumerate}
\item if $\mathcal{L}=\nbIci$, $\nu(\neg(\alpha\uparrow\neg\alpha))=\nu(\alpha)\wedge\nu(\neg\alpha)$;
\item if $\mathcal{L}=\nbIcl$, $\nu(\alpha\uparrow\neg\alpha)=\nu(\neg(\alpha\wedge\neg\alpha))$.
\end{enumerate}
\end{enumerate}

\begin{proposition}
The RNmatrices $(\mathcal{A}, \{\top\}, \mathcal{F}_{\mathcal{A}}^{\mathcal{L}})$ as described above, for $\mathcal{L}\in$\\$\{\nbIciw, \nbIci, \nbIcl\}$, are structural.
\end{proposition}

\begin{proof}
Let $\sigma:\textbf{F}(\Sigma_{\nbI},\mathcal{V})\rightarrow\textbf{F}(\Sigma_{\nbI},\mathcal{V})$ be a $\Sigma_{\nbI}$-homomorphism and $\nu\in\mathcal{F}_{\mathcal{A}}^{\mathcal{L}}$: it is clear that $\nu\circ\sigma: \textbf{F}(\Sigma_{\nbI},\mathcal{V})\rightarrow\mathcal{A}$ is always a $\Sigma_{\nbI}$-homomorphism, satisfying additionally that for all formulas $\alpha$ and $\beta$ on the signature $\Sigma_{\nbI}$, $\nu\circ\sigma(\alpha\uparrow\beta)=\nu\circ\sigma(\beta\uparrow\alpha)$.

\begin{enumerate}
\item We must prove that $\nu\circ\sigma(\alpha\uparrow\neg\alpha)={\sim}(\nu\circ\sigma(\alpha)\wedge\nu\circ\sigma(\neg\alpha))$, what is easy since
\[\nu\circ\sigma(\alpha\uparrow\neg\alpha)=\nu(\sigma(\alpha\uparrow\neg\alpha))=\nu(\sigma(\alpha)\uparrow\sigma(\neg\alpha))=\nu(\sigma(\alpha)\uparrow\neg\sigma(\alpha))=\]
\[{\sim}\big(\nu(\sigma(\alpha))\wedge\neg\nu(\sigma(\alpha))\big)={\sim}(\nu\circ\sigma(\alpha)\wedge\nu\circ\sigma(\neg\alpha)).\]

\item\begin{enumerate}
\item If $\mathcal{L}=\nbIci$, it remains for us to show that $\nu\circ\sigma(\neg(\alpha\uparrow\neg\alpha))=\nu\circ\sigma(\alpha)\wedge\nu\circ\sigma(\neg\alpha)$, and we see that
\[\nu\circ\sigma(\neg(\alpha\uparrow\neg\alpha))=\nu\big(\neg(\sigma(\alpha)\uparrow\neg\sigma(\alpha))\big)=\nu(\sigma(\alpha))\wedge\nu(\neg\sigma(\alpha))=\nu\circ\sigma(\alpha)\wedge\nu\circ\sigma(\neg\alpha).\]

\item If $\mathcal{L}=\nbIcl$, we only need to prove that $\nu\circ\sigma(\alpha\uparrow\neg\alpha)=\nu\circ\sigma(\neg(\alpha\wedge\neg\alpha))$:
\[\nu\circ\sigma(\alpha\uparrow\neg\alpha)=\nu(\sigma(\alpha)\uparrow\neg\sigma(\alpha))=\nu\big(\neg(\sigma(\alpha)\wedge\neg\sigma(\alpha))\big)=\nu\circ\sigma(\neg(\alpha\wedge\neg\alpha)).\]
\end{enumerate}
\end{enumerate}
\end{proof}

Given formulas $\Gamma\cup\{\varphi\}$ over the signature $\Sigma_{\nbI}$, if $\Gamma$ proves $\varphi$ according to such restricted Nmatrices we will write \label{VdashFnbIciw}\label{VdashFnbIci}\label{VdashFnbIcl}$\Gamma\Vdash_{\mathcal{F}}^{\mathcal{L}}\varphi$, and say that $\Gamma$ proves $\varphi$ according to Fidel structures for $\mathcal{L}$.

It is easy to see how "$\Vdash_{\mathcal{F}}^{\mathcal{L}}$" models the axiom schemata of $\nbI$ and its rules of inference, but we can also prove "$\Vdash_{\mathcal{F}}^{\nbIciw}$" proves any instance of $\textbf{ciw}^{*}$, "$\Vdash_{\mathcal{F}}^{\nbIci}$" proves any instance of $\textbf{ci}^{*}$ and "$\Vdash_{\mathcal{F}}^{\nbIcl}$" proves any instance of $\textbf{cl}^{*}$, so take a formula $\alpha$ over the signature $\Sigma_{\nbI}$.

\begin{enumerate}
\item Given an instance $(\alpha\uparrow\neg\alpha)\vee(\alpha\wedge\neg\alpha)$ of $\textbf{ciw}^{*}$, we have for a Fidel structure $\mathcal{A}$ for $\nbIciw$ and a $\nu\in \mathcal{F}_{\mathcal{A}}^{\nbIciw}$ that
\[\nu((\alpha\uparrow\neg\alpha)\vee(\alpha\wedge\neg\alpha))=\nu(\alpha\uparrow\neg\alpha)\vee\nu(\alpha\wedge\neg\alpha)=\]
\[{\sim}(\nu(\alpha)\wedge\nu(\neg\alpha))\vee(\nu(\alpha)\wedge\nu(\neg\alpha))=({\sim}\nu(\alpha)\vee{\sim}\nu(\neg\alpha))\vee(\nu(\alpha)\wedge\nu(\neg\alpha))=\]
\[\big[({\sim}\nu(\alpha)\vee{\sim}\nu(\neg\alpha))\vee\nu(\alpha)\big]\wedge\big[({\sim}\nu(\alpha)\vee{\sim}\nu(\neg\alpha))\vee\nu(\neg\alpha)\big]=\]
\[[{\sim}\nu(\neg\alpha)\vee\top]\wedge[{\sim}\nu(\alpha)\vee\top]=\top\wedge\top=\top.\]

\item Take an instance $\neg(\alpha\uparrow\neg\alpha)\rightarrow(\alpha\wedge\neg\alpha)$ of $\textbf{ci}^{*}$, a Fidel structure $\mathcal{A}$ for $\nbIci$ and a $\nu\in\mathcal{F}_{\mathcal{A}}^{\nbIci}$:
\[\nu(\neg(\alpha\uparrow\neg\alpha)\rightarrow(\alpha\wedge\neg\alpha))=\nu(\neg(\alpha\uparrow\neg\alpha))\rightarrow\nu(\alpha\wedge\neg\alpha)=\]
\[\big[\nu(\alpha)\wedge\nu(\neg\alpha)\big]\rightarrow\big[\nu(\alpha)\wedge\nu(\neg\alpha)\big]=\top.\]

\item For $\neg(\alpha\wedge\neg\alpha)\rightarrow(\alpha\uparrow\neg\alpha)$ an instance of $\textbf{cl}^{*}$, a Fidel structure $\mathcal{A}$ for $\nbIcl$ and a $\nu\in\mathcal{F}_{\mathcal{A}}^{\nbIcl}$,
\[\nu(\neg(\alpha\wedge\neg\alpha)\rightarrow(\alpha\uparrow\neg\alpha))=\nu(\neg(\alpha\wedge\neg\alpha))\rightarrow\nu(\alpha\uparrow\neg\alpha)=\nu(\alpha\uparrow\neg\alpha)\rightarrow\nu(\alpha\uparrow\neg\alpha)=\top.\]
\end{enumerate}

\begin{theorem}
Given formulas $\Gamma\cup\{\varphi\}$ over the signature $\Sigma_{\nbI}$, for any $\mathcal{L}\in$\\$\{\nbIciw, \nbIci, \nbIcl\}$ we have that, if $\Gamma\vdash_{\mathcal{L}}\varphi$, then $\Gamma\vDash_{\mathcal{F}}^{\mathcal{L}}\varphi$.
\end{theorem}

Reciprocally, we define, for a logic $\mathcal{L}\in\{\nbIciw,\nbIci,\nbIcl\}$ and a set of formulas $\Gamma$ over the signature $\Sigma_{\nbI}$, the equivalence relation such that, for formulas $\alpha$ and $\beta$ still over the signature $\Sigma_{\nbI}$, $\alpha\equiv^{\mathcal{L}}_{\Gamma}\beta$ if and only if $\Gamma\vdash_{\mathcal{L}}\alpha\rightarrow\beta$ and $\Gamma\vdash_{\mathcal{L}}\beta\rightarrow\alpha$.

As we did several times before, the well-defined quotient $A^{\mathcal{L}}_{\Gamma}=F(\Sigma_{\nbI}, \mathcal{V})/\equiv^{\mathcal{L}}_{\Gamma}$ is made into a Boolean algebra, and by defining, for classes of formulas $[\alpha]$ and $[\beta]$, $\neg\alpha=$\\$\{[\neg\varphi]\ :\  \varphi\in[\alpha]\}$ and
\begin{enumerate}
\item if $[\beta]\in \neg[\alpha]$,
\[[\alpha]\uparrow[\beta]=[{\sim}(\alpha\wedge\beta)];\]
\item otherwise,
\[[\alpha]\uparrow[\beta]=\{[\varphi\uparrow\psi]\ :\  \varphi\in[\alpha], \psi\in[\beta]\};\]
\end{enumerate}
we shall prove $\mathcal{A}^{\mathcal{L}}_{\Gamma}=(A^{\mathcal{L}}_{\Gamma},\{\sigma_{\mathcal{A}}\}_{\sigma\in \Sigma_{\nbI}^{\textbf{CPL}}})$ is a Fidel structure for $\mathcal{L}$, which we shall call the Lindenbaum-Tarski Fidel structure of $\mathcal{L}$ associated to $\Gamma$.

Quite clearly $(A^{\mathcal{L}}_{\Gamma},\{\sigma_{\mathcal{A}}\}_{\sigma\in\Sigma^{\textbf{CPL}}})$ is a Boolean algebra; concerning the paraconsistent negation, for every $[\beta]\in\neg[\alpha]$, there exists $\varphi\equiv^{\mathcal{L}}_{\Gamma}\alpha$ such that $\beta\equiv^{\mathcal{L}}_{\Gamma}\neg\varphi$ and then $[\alpha]\vee[\beta]=[\alpha\vee\beta]=[\varphi\vee\neg\varphi]=\top$; regarding the incompatibility connective, we begin by proving the operation is well-defined in the case that $[\beta]\in\neg[\alpha]$.

If $\varphi\in[\alpha]$ and $\psi\in[\beta]$, we have that $\alpha\equiv^{\mathcal{L}}_{\Gamma}\varphi$ and $\beta\equiv^{\mathcal{L}}_{\Gamma}\psi$, implying that ${\sim}(\alpha\wedge\beta)\equiv^{\mathcal{L}}_{\Gamma}{\sim}(\varphi\wedge\psi)$ given "$\equiv^{\mathcal{L}}_{\Gamma}$" is a congruence for the connectives in $\{\vee, \wedge, \rightarrow, {\sim}\}$, and therefore $[\alpha]\uparrow[\beta]=[\varphi]\uparrow[\psi]$. In the case of $[\alpha]\uparrow[\beta]$ with $[\beta]\not\in\neg[\alpha$, we use the same reasoning that worked for $\bI$ in Section \ref{Completeness fo Fidel structures for bI}.

To prove $\mathcal{A}^{\mathcal{L}}_{\Gamma}$ is a Fidel structure for $\nbI$, we are yet to prove that, if $[\beta]\in\neg[\alpha]$, for every value $[\gamma]$ in $[\alpha]\uparrow[\beta]$ (in this case, there is only one), we have $[\alpha]\wedge([\beta]\wedge[\gamma])=\bot$, the result being clear in the case that $[\beta]\notin\neg[\alpha]$. Since $\gamma\equiv^{\mathcal{L}}_{\Gamma}{\sim}(\alpha\wedge\beta)$,
\[\alpha\wedge(\beta\wedge\gamma)\equiv^{\mathcal{L}}_{\Gamma}\alpha\wedge(\beta\wedge{\sim}(\alpha\wedge\beta))\equiv^{\mathcal{L}}_{\Gamma}\bot,\]
and the result holds; the case in which $[\alpha]\in\neg[\beta]$ is analogous.

It only remains to be shown that $\mathcal{A}^{\mathcal{L}}_{\Gamma}$ is a Fidel structure for $\mathcal{L}$.
\begin{lemma}\label{some equivalences in extensions of nbI}
\begin{enumerate}
\item For any logic $\mathcal{L}\in\{\nbIciw, \nbIci, \nbIcl\}$, $\alpha\uparrow\neg\alpha$ and ${\sim}(\alpha\wedge\neg\alpha)$ are equivalent.
\item $\neg(\alpha\uparrow\neg\alpha)$ and $\alpha\wedge\neg\alpha$ are equivalent in $\nbIci$.
\item $\neg(\alpha\wedge\neg\alpha)$ and $\alpha\uparrow\neg\alpha$ are equivalent in $\nbIcl$.
\end{enumerate}
\end{lemma}

\begin{proof}
\begin{enumerate}
\item For simplicity, we denote $\bot_{\alpha\wedge\neg\alpha,\alpha\wedge\neg\alpha}$ simply by $\bot$. Applying the deduction meta-theorem, we must prove $\alpha\uparrow\neg\alpha, \alpha\wedge\neg\alpha\vdash_{\mathcal{L}}\bot$ and $(\alpha\wedge\neg\alpha)\rightarrow\bot\vdash_{\mathcal{L}}\alpha\uparrow\neg\alpha$.

The first implication is obvious, since $\alpha\uparrow\neg\alpha$ together with $\alpha\wedge\neg\alpha$ is exactly $\bot_{\alpha,\neg\alpha}$, and all bottom elements are equivalent to each other. Reciprocally, quite obviously
\[(\alpha\wedge\neg\alpha)\rightarrow\bot, \alpha\uparrow\neg\alpha\vdash_{\mathcal{L}}\alpha\uparrow\neg\alpha,\]
and
\[(\alpha\wedge\neg\alpha)\rightarrow\bot, \alpha\wedge\neg\alpha\vdash_{\mathcal{L}}\alpha\uparrow\neg\alpha\]
by Modus Ponens and the fact a bottom element implies any formula. With a proof by cases and the fact that $(\alpha\uparrow\neg\alpha)\vee(\alpha\wedge\neg\alpha)$ is an instance of $\textbf{ciw}^{*}$, which is a tautology in any of the logics $\mathcal{L}$, we get the desired result.

\item We must prove $\neg(\alpha\uparrow\neg\alpha)\vdash_{\nbIci}\alpha\wedge\neg\alpha$ and vice-versa, being the first direction a clear application of $\textbf{ci}^{*}$. Reciprocally, we have 
\[\alpha\wedge\neg\alpha, \neg(\alpha\uparrow\neg\alpha)\vdash_{\nbIci}\neg(\alpha\uparrow\neg\alpha)\]
with easy, while 
\[\alpha\wedge\neg\alpha, \alpha\uparrow\neg\alpha\vdash_{\nbIci}\neg(\alpha\uparrow\neg\alpha)\]
follows from the fact $\alpha\wedge\neg\alpha$ implies both $\alpha$ and $\neg\alpha$, followed by an application of $\textbf{Ip}$. With a proof by cases and the fact $(\alpha\uparrow\neg\alpha)\vee\neg(\alpha\uparrow\neg\alpha)$ is an instance of $\textbf{Ax\: 11}^{*}$, the result follows.

\item Now, we must prove $\neg(\alpha\wedge\neg\alpha)\vdash_{\nbIcl}\alpha\uparrow\neg\alpha$ and its reciprocal, being the first implication a direct application of $\textbf{cl}^{*}$. Reciprocally, 
\[\alpha\uparrow\neg\alpha, \neg(\alpha\wedge\neg\alpha)\vdash_{\nbIcl}\neg(\alpha\wedge\neg\alpha)\]
and, from the fact $\alpha\wedge\neg\alpha$ implies both $\alpha$ and $\neg\alpha$, and by applying $\textbf{Ip}$,
\[\alpha\uparrow\neg\alpha, \alpha\wedge\neg\alpha\vdash_{\nbIcl}\neg(\alpha\wedge\neg\alpha);\]
through a proof by cases and the fact $(\alpha\wedge\neg\alpha)\vee\neg(\alpha\wedge\neg\alpha)$ is an instance of $\textbf{Ax\: 11}^{*}$, we conclude the proof.
\end{enumerate}
\end{proof}

\begin{enumerate}
\item By definition, if $[\beta]\in\neg[\alpha]$, then $[\alpha]\uparrow[\beta]$ is single-valued and equal to ${\sim}([\alpha]\wedge[\beta])$.
\item \begin{enumerate}
\item If $\mathcal{L}=\nbIci$ and $[\beta]\in\neg[\alpha]$, implying that there exists $\phi\equiv^{\nbIci}_{\Gamma}\alpha$ such that $\beta\equiv^{\nbIci}_{\Gamma}\neg\phi$, by use of Lemma \ref{some equivalences in extensions of nbI} we find
\[[{\sim}(\alpha\wedge\beta)]=[{\sim}(\phi\wedge\neg\phi)]=[\phi\uparrow\neg\phi],\]
and therefore $\neg[{\sim}(\alpha\wedge\beta)]$ contains $[\neg(\phi\uparrow\neg\phi)]=[\phi\wedge\neg\phi]$, which equals $[\alpha]\wedge[\beta]$, as we wanted to show.
\item If $\mathcal{L}=\nbIcl$ and $[\beta]\in\neg[\alpha]$, let $\phi\equiv^{\nbIcl}_{\Gamma}\alpha$ be an element of $[\alpha]$ such that $\beta\equiv^{\nbIcl}_{\Gamma}\neg\phi$, and then $[\alpha\wedge\beta]=[\phi\wedge\neg\phi]$, while $[\phi\uparrow\neg\phi]=[\neg(\phi\wedge\neg\phi)]\in \neg[\phi\wedge\neg\phi]$. 

Since $[\phi\uparrow\neg\phi]=[{\sim}(\phi\wedge\neg\phi)]$ and ${\sim}([\phi]\wedge[\neg\phi])={\sim}([\alpha]\wedge[\beta])$, the result is done: ${\sim}([\alpha]\wedge[\beta])\in \neg([\alpha]\wedge[\beta])$.
\end{enumerate}
\end{enumerate}

\begin{theorem}
Given formulas $\Gamma\cup\{\varphi\}$ over the signature $\Sigma_{\nbI}$, for any logic $\mathcal{L}\in$\\$\{\nbIciw, \nbIci, \nbIcl\}$ we have that, if $\Gamma\Vdash_{\mathcal{F}}^{\mathcal{L}}\varphi$, then $\Gamma\vdash_{\mathcal{L}}\varphi$.
\end{theorem}

\subsection{Decision method}\label{Decision method for nbIciw, ...}

We take the Boolean algebra $\textbf{2}$ and extend it to the $\Sigma^{\textbf{CPL}}_{\textbf{nbI}}$-multialgebra $\textbf{2}_{\textbf{nbIciw}}$ with operations given by the tables below.

\begin{figure}[H]
\centering
\begin{minipage}[t]{4cm}
\centering
\begin{tabular}{|l|r|}
\hline
& $\neg$\\ \hline
$0$ & $\{1\}$\\ \hline
$1$ & $\{0,1\}$\\\hline
\end{tabular}
\caption*{Negation}
\end{minipage}
\hspace{3cm}
\centering
\begin{minipage}[t]{4cm}
\centering
\begin{tabular}{|l|c|r|}
\hline
$\uparrow$ & $0$ & $1$\\ \hline
$0$ & $\{0,1\}$ & $\{0, 1\}$\\ \hline
$1$ & $\{0, 1\}$ & $\{0\}$\\\hline
\end{tabular}
\caption*{Incompatibility}
\end{minipage}
\end{figure}

It is easy to see that:
\begin{enumerate}
\item $(\{0,1\}, \{\sigma_{\textbf{2}_{\textbf{nbIciw}}}\}_{\sigma\in\Sigma^{\textbf{CPL}}})$ is a Boolean algebra (that is, $\textbf{2}$);
\item for any $x,y\in\{0,1\}$ and $z\in x\uparrow y$, $x\wedge(y\wedge z)=0$;
\item for any $x\in \{0,1\}$ and $y\in\neg x$, $x\vee y=1$, since in the case that $x=0$ we must have $y=1$.
\end{enumerate}

So $\textbf{2}_{\textbf{nbIciw}}$ is a Fidel structure for $\textbf{nbI}$. But we can prove even more, that it is a Fidel structure for $\mathcal{L}\in\{\textbf{nbIciw}, \textbf{nbIci}, \textbf{nbIcl}\}$.

\begin{enumerate}
\item We see that, for any values of $x$ and $y$ in $\{0,1\}$, ${\sim}(x\wedge y)\in x\uparrow y$ even if $y\not\in\neg x$;
\item
\begin{enumerate}
\item suppose that, for some values $a$ and $b$ in $\{0,1\}$ such that $b\in\neg a$, $a\wedge b$ is not in $\neg{\sim}(a\wedge b)$, and since $\neg 1=\{0,1\}$, we must then have ${\sim}(a\wedge b)=0$, that is, $a\wedge b=1$ and therefore $a=b=1$; but then $\neg{\sim}(a\wedge b)=\{1\}$, which contains exactly the value of $a\wedge b=1$, leading to a contradiction; so, to summarize, for any values $a, b\in\{0,1\}$ we have $a\wedge b\in\neg{\sim}(a\wedge b)$;
\item suppose that, for some values $a$ and $b$ in $\{0,1\}$ such that $b\in\neg a$, ${\sim}(a\wedge b)$ is not in $\neg(a\wedge b)$, which means $a\wedge b=0$ (since $\neg 1=\{0,1\}$); but then ${\sim}(a\wedge b)=1$ and $\neg(a\wedge b)=\{1\}$, what is a contradiction; so, for any values $a,b\in\{0,1\}$, ${\sim}(a\wedge b)\in\neg(a\wedge b)$.
\end{enumerate}
\end{enumerate}

Finally, we can define, for $\mathcal{L}\in\{\textbf{nbIciw}, \textbf{nbIci}, \textbf{nbIcl}\}$, the restricted Nmatrix 
\[\mathbb{2}_{\mathcal{L}}=(\textbf{2}_{\textbf{nbIciw}}, \{1\}, \mathcal{F}_{\textbf{2}_{\textbf{nbIciw}}}^{\mathcal{L}})\]
such that $\mathcal{F}_{\textbf{2}_{\textbf{nbIciw}}}^{\mathcal{L}}$ is the set of homomorphisms $\nu:\textbf{F}(\Sigma_{\textbf{nbI}}, \mathcal{V})\rightarrow \textbf{2}_{\textbf{nbIciw}}$ satisfying that, for any formulas $\alpha$ and $\beta$:
\begin{enumerate}
\item $\nu(\alpha\uparrow\beta)=\nu(\beta\uparrow\alpha)$;
\item $\nu(\alpha\uparrow\neg\alpha)={\sim}(\nu(\alpha)\wedge\nu(\neg\alpha))$;
\item\begin{enumerate}
\item if $\mathcal{L}=\textbf{nbIci}$, $\nu(\neg(\alpha\uparrow\neg\alpha))=\nu(\alpha)\wedge\nu(\neg\alpha)$;
\item if $\mathcal{L}=\textbf{nbIcl}$, $\nu(\alpha\uparrow\neg\alpha)=\nu(\neg(\alpha\wedge\neg\alpha))$.
\end{enumerate}
\end{enumerate}
Clearly such RNmatrices are structural.

\begin{theorem}
$\nu: F(\Sigma_{\textbf{nbI}},\mathcal{V})\rightarrow\{0,1\}$ is a bivaluation for $\mathcal{L}\in\{\textbf{nbIciw}, \textbf{nbIci}, \textbf{nbIcl}\}$ if, and only if, it is a $\Sigma_{\textbf{nbI}}$-homomorphism from $\textbf{F}(\Sigma_{\textbf{nbI}}, \mathcal{V})$ to $\textbf{2}_{\textbf{nbIciw}}$ which lies in $\mathcal{F}_{\textbf{2}_{\textbf{nbIciw}}}^{\mathcal{L}}$.
\end{theorem}

\begin{theorem}
Given formulas $\Gamma\cup\{\varphi\}$ of $\mathcal{L}\in\{\textbf{nbIciw}, \textbf{nbIci}, \textbf{nbIcl}\}$, $\Gamma\vDash_{\mathcal{L}}\varphi$ if and only if $\Gamma\vDash_{\mathbb{2}_{\mathcal{L}}}\varphi$.
\end{theorem}

Our finite RNmatrices again lead to decision methods through row-branching, row-eliminating truth tables, and for all three of the logics $\nbIciw$, $\nbIci$ and $\nbIcl$. We believe the subjacent Nmatrices are explicit enough, so let us focus on the conditions for a row to be erased.
\begin{enumerate}
\item If $\alpha\uparrow\beta$ and $\beta\uparrow\alpha$ have different values.
\item If $\alpha\uparrow\neg\alpha$ or $\neg\alpha\uparrow\alpha$ is $0$, and either $\alpha$ or $\neg\alpha$ is also $0$.
\item \begin{enumerate}
\item If the logic is $\nbIci$: if $\neg(\alpha\uparrow\neg\alpha)$ is $1$ but either $\alpha$ or $\neg\alpha$ is $0$.
\item If the logic is $\nbIcl$, $\neg(\alpha\wedge\neg\alpha)$ is $1$ and: either $\alpha\uparrow\neg\alpha$ or $\neg\alpha\uparrow\alpha$ is $0$, or $\alpha$ and $\neg\alpha$ are both $1$.
\end{enumerate}
\end{enumerate}

\subsection{Another decision method}

We can add rules to the tableau calculus $\mathbb{T}_{\nbI}$ inspired by the tables found in Section \ref{Decision method for nbIciw, ...} to obtain tableau calculi $\mathbb{T}_{\nbIciw}$, $\mathbb{T}_{\nbIci}$ and $\mathbb{T}_{\nbIcl}$ capable of characterizing their respective logics. So, consider the following tableau rules:

$$
\begin{array}{cp{1cm}cp{1cm}cp{1cm}c}
\displaystyle \frac{\textsf{0}(\varphi\uparrow\neg\varphi)}{\begin{array}{c}\textsf{1}(\varphi) \\ \textsf{1}(\neg\varphi)\end{array}} & &
\displaystyle \frac{\textsf{0}(\neg\varphi\uparrow\varphi)}{\begin{array}{c}\textsf{1}(\varphi) \\ \textsf{1}(\neg\varphi)\end{array}} & & 
\displaystyle \frac{\textsf{1}(\neg(\varphi\uparrow\neg\varphi))}{\begin{array}{c}\textsf{1}(\varphi) \\ \textsf{1}(\neg\varphi)\end{array}} & & 
\displaystyle \frac{\textsf{1}(\neg(\varphi\wedge\neg\varphi))}{\textsf{0}(\varphi)\mid\textsf{0}(\neg\varphi)}\\[2mm]
\end{array}
$$

Then, by addition to $\mathbb{T}_{\nbI}$ of:
\begin{enumerate}
\item the two first rules, one obtains $\mathbb{T}_{\nbIciw}$;
\item the three first rules, one obtains $\mathbb{T}_{\nbIci}$;
\item the two first rules plus the fourth, one obtains $\mathbb{T}_{\nbIcl}$.
\end{enumerate}

\section{Our logics are not algebraizable by Blok and Pigozzi}

\subsection{$\bI$ is not algebraizable by Blok and Pigozzi}\label{bI is not algebraizable}

We wish to prove $\bI$ is not algebraizable even by some quite expressive standards: specifically, we wish to show it is not \index{Algebraizable}algebraizable according to Blok and Pigozzi \cite{BlokPigozzi}. To accomplish that we will use \index{Lewin}Lewin, Mikenberg and Schwarze's construction found in\cite{Lewin}, and prove as they do that there is a model for $\bI$ for which the Leibniz operator $\Omega_{\bI}$, that sends a $\bI$-filter to its largest compatible congruence, is not a bijection; from Theorem $5.1$ of \cite{BlokPigozzi}, this proves $\bI$ is not algebraizable.

\begin{definition}
Given a signature $\Sigma$ and a $\Sigma$-algebra $\mathcal{A}=(A, \{\sigma_{\mathcal{A}}\}_{\sigma\in\Sigma})$, a \index{Congruence}congruence in $\mathcal{A}$ is a relation $\theta$ on $A\times A$ such that, if $\sigma\in\Sigma_{n}$ and $a_{1}, b_{1}, \dotsc  , a_{n}, b_{n}\in A$ are elements such that $a_{1}\theta b_{1}, \dotsc , a_{n}\theta b_{n}$, then
\[\sigma_{\mathcal{A}}(a_{1}, \dotsc  , a_{n})\theta\sigma_{\mathcal{A}}(b_{1}, \dotsc  , b_{n}).\]
\end{definition}

\begin{definition}
Given a signature $\Sigma$, a logic $\mathcal{L}$ and a $\Sigma$-algebra $\mathcal{A}=(A, \{\sigma_{\mathcal{A}}\}_{\sigma\in\Sigma})$, an \index{Filter, $\mathcal{L}$-}$\mathcal{L}$-filter in $\mathcal{A}$ is a subset $F\subseteq A$ such that 
\[\Gamma\vdash_{\mathcal{L}}\varphi\quad\text{implies}\quad \Gamma\vDash_{(\mathcal{A}, F)}\varphi,\]
for every set of formulas $\Gamma\cup\{\varphi\}$ over $\Sigma$.
\end{definition}

The \index{Congruence, Largest compatible}largest compatible congruence $\theta$ to a filter $F$ is the largest congruence such that, if $a\theta b$ and $a\in F$, then $b\in F$.

So, over the signature $\Sigma_{\bI}$, for which $(\Sigma_{\bI})_{2}=\{\vee, \wedge, \rightarrow, \uparrow\}$ and $(\Sigma_{\bI})_{n}=\emptyset$ for every $n\neq2$, we consider the $\Sigma_{\bI}$-algebra $\mathfrak{L}$ with universe $L=\{u, 1, a, b, 0\}$ and operations given by the tables bellow (all, with the exception of incompatibility, extracted from \cite{Lewin}).

\begin{figure}[H]
\centering
\begin{minipage}[t]{5cm}
\centering
\begin{tabular}{l|ccccr}
$\vee$ & $u$ & $1$ & $a$ & $b$ & $0$ \\\hline
$u$ & $u$ & $u$ & $u$ & $u$ & $u$\\
$1$ & $u$ & $1$ & $1$ & $1$ & $1$\\
$a$ & $u$ & $1$ & $a$ & $1$ & $a$\\
$b$ & $u$ & $1$ & $1$ & $b$ & $b$\\
$0$ & $u$ & $1$ & $a$ & $b$ & $0$
\end{tabular}
\caption*{Disjunction}
\end{minipage}
\centering
\begin{minipage}[t]{5cm}
\centering
\begin{tabular}{l|ccccr}
$\wedge$ & $u$ & $1$ & $a$ & $b$ & $0$ \\\hline
$u$ & $u$ & $1$ & $a$ & $b$ & $0$\\
$1$ & $1$ & $1$ & $a$ & $b$ & $0$\\
$a$ & $a$ & $a$ & $a$ & $0$ & $0$\\
$b$ & $b$ & $b$ & $0$ & $b$ & $0$\\
$0$ & $0$ & $0$ & $0$ & $0$ & $0$
\end{tabular}
\caption*{Conjunction}
\end{minipage}
\end{figure}

\begin{figure}[H]
\centering
\begin{minipage}[t]{5cm}
\centering
\begin{tabular}{l|ccccr}
$\rightarrow$ & $u$ & $1$ & $a$ & $b$ & $0$ \\\hline
$u$ & $u$ & $u$ & $a$ & $b$ & $0$\\
$1$ & $u$ & $1$ & $a$ & $b$ & $0$\\
$a$ & $u$ & $1$ & $1$ & $b$ & $b$\\
$b$ & $u$ & $1$ & $a$ & $1$ & $a$\\
$0$ & $u$ & $1$ & $1$ & $1$ & $1$
\end{tabular}
\caption*{Implication}
\end{minipage}
\begin{minipage}[t]{5cm}
\centering
\begin{tabular}{l|ccccr}
$\uparrow$ & $u$ & $1$ & $a$ & $b$ & $0$ \\\hline
$u$ & $0$ & $0$ & $0$ & $0$ & $1$\\
$1$ & $0$ & $0$ & $b$ & $a$ & $1$\\
$a$ & $0$ & $b$ & $b$ & $1$ & $1$\\
$b$ & $0$ & $a$ & $1$ & $a$ & $1$\\
$0$ & $1$ & $1$ & $1$ & $1$ & $1$
\end{tabular}
\caption*{Incompatibility}
\end{minipage}
\end{figure}

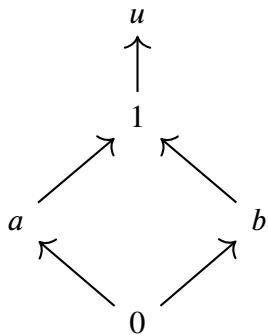
\begin{figure}[H]
\centering
\begin{tikzcd}
 & u & \\
 & 1 \arrow[u] & \\
 a \arrow[ru] & & b \arrow[lu]\\
 & 0 \arrow[ru]\arrow[lu] & \\
  \end{tikzcd}
\caption*{The lattice $(L, \vee, \wedge)$}
\end{figure}

We then consider the logical matrix $\mathfrak{M}=(\mathfrak{L}, D)$, with $D=\{u,1\}$. 

\subsubsection{$\mathfrak{M}$ is a model of $\bI$}

Notice that $\textbf{Ax\: 1}$, $\textbf{Ax\: 3}$, $\textbf{Ax\: 4}$, $\textbf{Ax\: 5}$, $\textbf{Ax\: 6}$, $\textbf{Ax\: 7}$, $\textbf{Ax\: 8}$ and Modus Ponens of the definition of $\bI$ ate the beginning of Section \ref{Defining bI} correspond, respectively, to the axiom schemata and rules $1$, $6$, $4$, $5$, $7$, $8$, $9$ and $10$ of $\textbf{C}_{1}$ as defined in \cite{Lewin}; since the operations $\vee$, $\wedge$ and $\rightarrow$ in $\mathfrak{L}$ are exactly the same as the ones in the algebra of Lewin, Mikenberg and Schwarze's article, and since the logical matrix formed by that algebra and $D$ models $\textbf{C}_{1}$, we obtain that $\mathfrak{M}$ models at least these axiom schemata and rules of inference. It remains to be shown it also models $\textbf{Ax\: 2}$, $\textbf{Ax\: 9}^{*}$, $\textbf{Ip}$ and $\textbf{Comm}$.

\begin{enumerate}
\item Concerning $\textbf{Ax\: 2}$ of $\bI$, $(\alpha\rightarrow(\beta\rightarrow \gamma)\big)\rightarrow\big((\alpha\rightarrow\beta)\rightarrow(\alpha\rightarrow\gamma))$, notice that by two applications of the deduction meta-theorem we have that the validity of the axiom schema $2$ of \cite{Lewin} in $\bI$ is equivalent to stating that $\alpha\rightarrow \beta, \alpha\rightarrow(\beta\rightarrow \gamma)\vdash_{\bI}\alpha\rightarrow\gamma$, which by two new applications of the deduction meta-theorem is now equivalent to the validity of $\textbf{Ax\: 2}$. Since the matrix of Lewin, Mikenberg and Schwarze's algebra validates the axiom schema $2$, we have $\mathfrak{M}$ validates $\textbf{Ax\: 2}$.

\item For $\textbf{Ax\: 9}^{*}$, $(\alpha\rightarrow\beta)\vee\alpha$, we see from the tables below that the image of $\textbf{Ax\: 9}^{*}$ under any homomorphism is in $D$, and therefore this axiom schema is validated by $\mathfrak{M}$.

\begin{figure}[H]
\centering
\begin{tabular}{l|ccccc|ccccr}
& \multicolumn{5}{c|}{$x\rightarrow y$} & \multicolumn{5}{c}{$(x\rightarrow y)\vee x$}\\ \hline
\diag{.1em}{0.2cm}{$x$}{$y$} & $u$ & $1$ & $a$ & $b$ & $0$ & $u$ & $1$ & $a$ & $b$ & $0$ \\\hline
$u$ & $u$ & $u$ & $a$ & $b$ & $0$ & $u$ & $u$ & $u$ & $u$ & $u$\\
$1$ & $u$ & $1$ & $a$ & $b$ & $0$ & $u$ & $1$ & $1$ & $1$ & $1$\\
$a$ & $u$ & $1$ & $1$ & $b$ & $b$ & $u$ & $1$ & $1$ & $1$ & $1$\\
$b$ & $u$ & $1$ & $a$ & $1$ & $a$ & $u$ & $1$ & $1$ & $1$ & $1$\\
$0$ & $u$ & $1$ & $1$ & $1$ & $1$ & $u$ & $1$ & $1$ & $1$ & $1$
\end{tabular}
\caption*{Table for $\textbf{Ax\: 9}^{*}$}
\end{figure}

\item To prove $\mathfrak{M}$ models $\textbf{Ip}$, given by $(\alpha\uparrow\beta)\rightarrow(\alpha\rightarrow(\beta\rightarrow\gamma))$, we begin with the two tables below.

\begin{figure}[H]
\centering
\begin{minipage}[t]{5cm}
\centering
\begin{tabular}{l|ccccr}
\diag{.1em}{0.2cm}{$y$}{$z$} & $u$ & $1$ & $a$ & $b$ & $0$ \\\hline
$u$ & $u$ & $u$ & $a$ & $b$ & $0$\\
$1$ & $u$ & $1$ & $a$ & $b$ & $0$\\
$a$ & $u$ & $1$ & $1$ & $b$ & $b$\\
$b$ & $u$ & $1$ & $a$ & $1$ & $a$\\
$0$ & $u$ & $1$ & $1$ & $1$ & $1$
\end{tabular}
\caption*{Table for $y\rightarrow z$}
\end{minipage}
\centering
\begin{minipage}[t]{5cm}
\centering
\begin{tabular}{l|ccccr}
\diag{.1em}{0.2cm}{$x$}{$y$} & $u$ & $1$ & $a$ & $b$ & $0$ \\\hline
$u$ & $0$ & $0$ & $0$ & $0$ & $1$\\
$1$ & $0$ & $0$ & $b$ & $a$ & $1$\\
$a$ & $0$ & $b$ & $b$ & $1$ & $1$\\
$b$ & $0$ & $a$ & $1$ & $a$ & $1$\\
$0$ & $1$ & $1$ & $1$ & $1$ & $1$
\end{tabular}
\caption*{Table for $x\uparrow y$}
\end{minipage}
\end{figure}

By looking at the table for $\textbf{Ip}$ below, we see that, since the image of $\textbf{Ip}$ under any homomorphism is always in $D$, we have $\mathfrak{M}$ validates this axiom schema.

\begin{figure}[H]
\centering
\begin{tabular}{l|c|c|ccccc|C{1.2em}C{1.2em}C{1.2em}C{1.2em}C{1.2em}}
& & $x\uparrow y$ & \multicolumn{5}{c|}{$x\rightarrow(y\rightarrow z)$} & \multicolumn{5}{c}{$(x\uparrow y)\rightarrow(x\rightarrow(y\rightarrow z))$}\\ \hline
$x$ & \diag{.1em}{0.2cm}{$y$}{$z$} & & $u$ & $1$ & $a$ & $b$ & $0$ & $u$ & $1$ & $a$ & $b$ & $0$ \\ \hline
\multirow{5}{*}{$u$} & $u$ & $0$ & $u$ & $u$ & $a$ & $b$ & $0$ & $u$ & $u$ & $1$ & $1$ & $1$\\
& $1$ & $0$ & $u$ & $u$ & $a$ & $b$ & $0$ & $u$ & $u$ & $1$ & $1$ & $1$\\
& $a$ & $0$ & $u$ & $u$ & $u$ & $b$ & $b$ & $u$ & $u$ & $u$ & $1$ & $1$\\
& $b$ & $0$ & $u$ & $u$ & $a$ & $u$ & $a$ & $u$ & $u$ & $1$ & $u$ & $1$\\
& $0$ & $1$ & $u$ & $u$ & $u$ & $u$ & $u$ & $u$ & $u$ & $u$ & $u$ & $u$\\ \hline
\multirow{5}{*}{$1$} & $u$ & $0$ & $u$ & $u$ & $a$ & $b$ & $0$ & $u$ & $u$ & $1$ & $1$ & $1$\\
& $1$ & $0$ & $u$ & $1$ & $a$ & $b$ & $0$ & $u$ & $1$ & $1$ & $1$ & $1$\\
& $a$ & $b$ & $u$ & $1$ & $1$ & $b$ & $b$ & $u$ & $1$ & $1$ & $1$ & $1$\\
& $b$ & $a$ & $u$ & $1$ & $a$ & $1$ & $a$ & $u$ & $1$ & $1$ & $1$ & $1$\\
& $0$ & $1$ & $u$ & $1$ & $1$ & $1$ & $1$ & $u$ & $1$ & $1$ & $1$ & $1$\\ \hline
\multirow{5}{*}{$a$} & $u$ & $0$ & $u$ & $u$ & $1$ & $b$ & $b$ & $u$ & $u$ & $1$ & $1$ & $1$\\
& $1$ & $b$ & $u$ & $1$ & $1$ & $b$ & $b$ & $u$ & $1$ & $1$ & $1$ & $1$\\
& $a$ & $b$ & $u$ & $1$ & $1$ & $b$ & $b$ & $u$ & $1$ & $1$ & $1$ & $1$\\
& $b$ & $1$ & $u$ & $1$ & $1$ & $1$ & $1$ & $u$ & $1$ & $1$ & $1$ & $1$\\
& $0$ & $1$ & $u$ & $1$ & $1$ & $1$ & $1$ & $u$ & $1$ & $1$ & $1$ & $1$\\ \hline
\multirow{5}{*}{$b$} & $u$ & $0$ & $u$ & $u$ & $a$ & $1$ & $a$ & $u$ & $u$ & $1$ & $1$ & $1$\\
& $1$ & $a$ & $u$ & $1$ & $a$ & $1$ & $a$ & $u$ & $1$ & $1$ & $1$ & $1$\\
& $a$ & $1$ & $u$ & $1$ & $1$ & $1$ & $1$ & $u$ & $1$ & $1$ & $1$ & $1$\\
& $b$ & $a$ & $u$ & $1$ & $a$ & $1$ & $a$ & $u$ & $1$ & $1$ & $1$ & $1$\\
& $0$ & $1$ & $u$ & $1$ & $1$ & $1$ & $1$ & $u$ & $1$ & $1$ & $1$ & $1$\\ \hline
\multirow{5}{*}{$0$} & $u$ & $1$ & $u$ & $u$ & $1$ & $1$ & $1$ & $u$ & $u$ & $1$ & $1$ & $1$\\
& $1$ & $1$ & $u$ & $1$ & $1$ & $1$ & $1$ & $u$ & $1$ & $1$ & $1$ & $1$\\
& $a$ & $1$ & $u$ & $1$ & $1$ & $1$ & $1$ & $u$ & $1$ & $1$ & $1$ & $1$\\
& $b$ & $1$ & $u$ & $1$ & $1$ & $1$ & $1$ & $u$ & $1$ & $1$ & $1$ & $1$\\
& $0$ & $1$ & $u$ & $1$ & $1$ & $1$ & $1$ & $u$ & $1$ & $1$ & $1$ & $1$\\
\end{tabular}
\caption*{Table for $\textbf{Ip}$}
\end{figure}

\item It only remains to show $\mathfrak{M}$ validates $\textbf{Comm}$, an axiom schema given by $(\alpha\uparrow\beta)\rightarrow(\beta\uparrow\alpha)$, what is done on the following table; notice that the tables for $x\uparrow y$ and $y\uparrow x$ are the same since incompatibility is commutative in $\mathfrak{M}$.

\begin{figure}[H]
\centering
\begin{tabular}{l|ccccc|ccccr}
& \multicolumn{5}{c|}{$x\uparrow y=y\uparrow x$} & \multicolumn{5}{c}{$(x\uparrow y)\rightarrow(y\uparrow x)$}\\\hline
\diag{.1em}{0.2cm}{$x$}{$y$} & $u$ & $1$ & $a$ & $b$ & $0$ & $u$ & $1$ & $a$ & $b$ & $0$ \\\hline
$u$ & $0$ & $0$ & $0$ & $0$ & $1$ & $1$ & $1$ & $1$ & $1$ & $1$\\
$1$ & $0$ & $0$ & $b$ & $a$ & $1$ & $1$ & $1$ & $1$ & $1$ & $1$\\
$a$ & $0$ & $b$ & $b$ & $1$ & $1$ & $1$ & $1$ & $1$ & $1$ & $1$\\
$b$ & $0$ & $a$ & $1$ & $a$ & $1$ & $1$ & $1$ & $1$ & $1$ & $1$\\
$0$ & $1$ & $1$ & $1$ & $1$ & $1$ & $1$ & $1$ & $1$ & $1$ & $1$
\end{tabular}
\caption*{Table for $\textbf{Comm}$}
\end{figure}

\end{enumerate}

With this, we have proved that $\mathfrak{M}$ is a model for $\bI$, meaning that, for any formulas $\Gamma\cup\{\varphi\}$ on the signature $\Sigma_{\bI}$, $\Gamma\vdash_{\bI}\varphi$ implies $\Gamma\vDash_{\mathfrak{M}}\varphi$.

\subsubsection{There are no non-trivial congruences on $\mathfrak{L}$}

Now, we wish to prove that the $\Sigma_{\bI}$-algebra $\mathfrak{L}$ has only two congruences, the ones we call trivial: \label{nabla}$\nabla$, equal to $L\times L$, and \label{Delta}$\Delta$, given by $\{(x,x)\ :\  x\in L\}$. Here, we will perpetrate an abuse of notation: take a signature $\Sigma$, a $\Sigma$-algebra $\mathcal{A}$ with universe $A$, a congruence $\theta$ on $\mathcal{A}$, $\sigma\in\Sigma_{n}$ and $\omega\in\Sigma_{m}$, and $a_{1}, \dotsc  , a_{n}, b_{1}, \dotsc  , b_{m}\in A$: then, if $\sigma_{\mathcal{A}}(a_{1}, \dotsc  , a_{n})=a$ and $\omega_{\mathcal{A}}(b_{1}, \dotsc  , b_{m})=b$, and $\sigma_{\mathcal{A}}(a_{1}, \dotsc  , a_{n})\theta\omega_{\mathcal{A}}(b_{1}, \dotsc  , b_{m})$, we may simply write $a=\sigma_{\mathcal{A}}(a_{1}, \dotsc  , a_{n})\theta\omega_{\mathcal{A}}(b_{1}, \dotsc  , b_{m})=b$ to make it clear that $a\theta b$.

\begin{enumerate}
\item 
\begin{enumerate}
\item If $u\theta 1$, $0=(u\uparrow a)\theta(1\uparrow a)=b$, $0=(u\uparrow b)\theta(1\uparrow b)=a$ and $1=(a\rightarrow a)\theta (b\rightarrow a)=a$, and therefore $\theta=\nabla$.
\item If $u\theta a$, then $u=(u\vee 1)\theta (a\vee 1)=1$ (meaning $a\theta 1$), $b=(u\wedge b)\theta (a\wedge b)=0$ and $0=(u\uparrow b)\theta(1\uparrow b)=a$, implying $\theta=\nabla$.
\item If $u\theta b$, then $u=(u\vee 1)\theta(b\vee 1)=1$, $a=(u\wedge a)\theta(b\wedge a)=0$ and $0=(u\uparrow a)\theta(1\uparrow a)=b$, and again $\theta=\nabla$.
\item If $u\theta 0$, then for every $x$ we have $u=(u\vee x)\theta (0\vee x)=x$, meaning $\theta=\nabla$.
\end{enumerate}
With this, we have proved that, if any pair $(u,x)$, with $x\in L\setminus\{u\}$, is in $\theta$, then $\theta=\nabla$.
\item We will ignore the case $1\theta u$, since it is equivalent to $u\theta 1$, a case considered before.
\begin{enumerate}
\item If $1\theta a$, $b=(1\wedge b)\theta (a\wedge b)=0$, $u=(u\rightarrow 1)\theta(u\rightarrow a)=a$ (and therefore $u\theta 1$) and $0=(u\uparrow b)\theta(1\uparrow b)=a$, meaning $\theta=\nabla$.
\item If $1\theta b$, $a=(1\wedge a)\theta(b\wedge a)=0$, $u=(u\rightarrow1)\theta(u\rightarrow b)=b$ (and therefore $u\theta 1$) and $0=(u\uparrow a)\theta(1\uparrow a)=b$, and so $\theta=\nabla$.
\item If $1\theta 0$, $a=(1\rightarrow a)\theta(0\rightarrow a)=1$, $b=(1\rightarrow b)\theta(0\rightarrow b)=1$ and $u=(u\rightarrow 1)\theta (u\rightarrow 0)=0$, and therefore $\theta=\nabla$.
\end{enumerate}
Once again, we have that if any pair $(1,x)$, with $x\in L\setminus\{1\}$, is in $\theta$, then $\theta=\nabla$.
\item The cases $a\theta u$ and $a\theta 1$ correspond to the cases $u\theta a$ and $1\theta a$ above.
\begin{enumerate}
\item If $a\theta 0$, $1=(a\vee b)\theta(0\vee b)=b$, $u=(u\rightarrow1)\theta(u\rightarrow b)=b$ (and therefore $u\theta 1$) and $0=(u\uparrow a)\theta(1\uparrow a)=b$, what implies $\theta=\nabla$.
\item If $a\theta b$, $b=(a\rightarrow b)\theta (b\rightarrow b)=1$ (implying $a\theta1$), $0=(a\wedge b)\theta(b\wedge b)=b$ and $u=(u\rightarrow 1)\theta(u\rightarrow a)=a$, what means $\theta=\nabla$.
\end{enumerate}
So, if $(a,x)\in\theta$, for $x\in L\setminus\{a\}$, $\theta=\nabla$.
\item The cases $b\theta u$, $b\theta 1$ and $b\theta a$ equal the previous cases, respectively, $u\theta b$, $1\theta b$ and $a\theta b$. 
\begin{enumerate}
\item If $b\theta 0$, $1=(a\vee b)\theta(a\vee 0)=a$, $u=(u\rightarrow 1)\theta(u\rightarrow a)=a$ (and therefore $u\theta 1$) and $0=(u\uparrow b)\theta(1\uparrow b)=a$, and therefore $\theta=\nabla$.
\end{enumerate}
Again, if $(b,x)\in \theta$, for $b\in L\setminus\{b\}$, then $\theta=\nabla$.
\item Since the cases $0\theta u$, $0\theta1$, $0\theta a$ and $0\theta b$ correspond, respectively, to $u\theta 0$, $1\theta0$, $a\theta0$ and $b\theta0$, we have that, if $(0,x)\in\theta$, for $x\in L\setminus\{0\}$, then $\theta=\nabla$.
\end{enumerate}

With all of this, we discover that for any congruence $\theta$ in $\mathfrak{L}$, if $(x,y)\in \theta$, for $x\neq y$, then $\theta=\nabla$, implying as we had mentioned that the only two congruences in $\mathfrak{L}$ are $\nabla$ and $\Delta$.

\subsubsection{Two filters whose largest compatible congruence is $\nabla$}

\begin{lemma}\label{Classifying filters}
Given a signature $\Sigma$, a $\Sigma$-algebra $\mathcal{A}$ and a logic $\mathcal{L}$ over $\Sigma$, $F$ is a $\mathcal{L}$-filter in $\mathcal{A}$ if, and only if, the following are satisfied:
\begin{enumerate}
\item for every $\Sigma$-homomorphism $\sigma:\textbf{F}(\Sigma, \mathcal{V})\rightarrow\mathcal{A}$ and instance of axiom $\psi$ of $\mathcal{L}$, $\sigma(\psi)\in F$;
\item for every $\Sigma$-homomorphism $\sigma:\textbf{F}(\Sigma, \mathcal{V})\rightarrow\mathcal{A}$ and instance of a rule of inference\\ $\psi_{1}, \dotsc  , \psi_{n}|\psi$ of $\mathcal{L}$, if $\sigma(\psi_{1}), \dotsc  , \sigma(\psi_{n})\in F$ then $\sigma(\psi)\in F$.
\end{enumerate}
\end{lemma}

\begin{proof}
Suppose $F$ is a $\mathcal{L}$-filter in $\mathcal{A}$ and $\sigma:\textbf{F}(\Sigma, \mathcal{V})\rightarrow\mathcal{A}$ is a $\Sigma$-homomorphism.
\begin{enumerate}
\item If $\psi$ is an instance of an axiom, $\vdash_{\mathcal{L}}\psi$, meaning $\vDash_{(\mathcal{A}, F)}\psi$ and, therefore, $\sigma(\psi)\in F$.
\item If $\psi_{1}, \dotsc  , \psi_{n}|\psi$ is an instance of an inference rule, $\psi_{1}, \dotsc  , \psi_{n}\vdash_{\mathcal{L}}\psi$, meaning $\psi_{1}, \dotsc  , \psi_{n}$\\$\vDash_{(\mathcal{A}, F)}\psi$ and, therefore, if $\sigma(\psi_{1}), \dotsc  , \sigma(\psi_{n})\in F$, then $\sigma(\psi)\in F$.
\end{enumerate}

Reciprocally, suppose $\Gamma\vdash_{\mathcal{L}}\varphi$ and that $\alpha_{1}, \dotsc  , \alpha_{n}$ is a proof of $\varphi$ from $\Gamma$ with $\alpha_{n}=\varphi$; we shall prove that, for any $\Sigma$-homomorphism $\sigma:\textbf{F}(\Sigma, \mathcal{V})\rightarrow\mathcal{A}$ such that $\sigma(\Gamma)\subseteq F$, $\sigma(\alpha_{1}), \dotsc  , \sigma(\alpha_{n})\in F$, and therefore $\sigma(\varphi)\in F$ and $\Gamma\vDash_{(\mathcal{A}, F)}\varphi$.

For induction hypothesis, being $\alpha_{1}$ either an element of $\Gamma$ or an instance of an axiom, assume $\alpha_{1}$ through $\alpha_{i-1}$ are mapped, by $\sigma$, into $F$.
\begin{enumerate}
\item If $\alpha_{i}$ is in $\Gamma$, by hypothesis on $\sigma$ we have $\sigma(\alpha_{i})\in F$.
\item If $\alpha_{i}$ is an instance of an axiom, by our hypothesis on $F$ we have $\sigma(\alpha_{i})\in F$.
\item Finally, if there are $\alpha_{i_{1}}, \dotsc  , \alpha_{i_{n}}$, with $i_{1}, \dotsc  , i_{n}\in\{1, \dotsc  , i-1\}$, such that $\alpha_{i_{1}}, \dotsc  , \alpha_{i_{n}}|\alpha_{i}$ is an instance of a rule of inference, by induction hypothesis $\sigma(\alpha_{i_{1}}), \dotsc  , \sigma(\alpha_{i_{n}})\in F$, and again by our hypothesis on $F$ we have $\sigma(\alpha_{i})\in F$.
\end{enumerate}
\end{proof}

So, we take the subsets $F_{a}=\{u, 1, a\}$ and $F_{b}=\{u, 1, b\}$ of $L$, and state that both are $\bI$-filters. First of all, for any $\Sigma_{\bI}$-homomorphism $\sigma:\textbf{F}(\Sigma_{\bI}, \mathcal{V})\rightarrow\mathfrak{L}$ and any instance of an axiom schema $\psi$ of $\bI$, since $\mathfrak{M}=(\mathfrak{L}, D)$ models $\bI$, $\sigma(\psi)\in D=\{u,1\}\subseteq F_{a}$ and in much the same way $\sigma(\psi)\in F_{b}$, implying that both $F_{a}$ and $F_{b}$ satisfy the first condition of Lemma \ref{Classifying filters} for being an $\bI$-filter.

Furthermore, there is only one rule of inference to analyze, that of Modus Ponens. From the table for implication, we see that if $x\rightarrow y$ is in $F_{a}$, then either $y\in F_{a}$ or $x\in L\setminus F_{a}$; so, if both $x$ and $x\rightarrow y$ are in $F_{a}$, then $y\in F_{a}$. Similarly, if both $x$ and $x\rightarrow y$ are in $F_{b}$, $y$ is also in $F_{b}$, what implies that both $F_{a}$ and $F_{b}$ are $\bI$-filters by Lemma \ref{Classifying filters}.

But we state that the largest compatible congruence to both $F_{a}$ and $F_{b}$ is $\Delta$, meaning the Leibniz operator $\Omega_{\bI}$ of $\bI$ is not injective and therefore $\bI$ is not algebraizable by Blok and Pigozzi. This is easy to prove: $\nabla$ is not compatible to neither $F_{a}$ nor $F_{b}$, since $u\Delta 0$ and $u\in F_{a}\cap F_{b}$, but $0$ is not in $F_{a}$ nor in $F_{b}$; clearly $\Delta$ is compatible to both $F_{a}$ and $F_{b}$, and since there are no congruences larger than $\Delta$ different from $\nabla$, we obtain the aforementioned result.

\subsection{$\nbI$ is not algebraizable by Blok and Pigozzi}

Now, we will follow the reasoning found in Section \ref{bI is not algebraizable} to show that $\nbI$ is also not algebraizable by Blok And Pigozzi, showing that its Leibniz operator is not bijective, the only difference being that we must add a negation to $\mathfrak{L}$. So, consider the $\Sigma_{\nbI}$-algebra $\mathfrak{L}_{\nbI}$ with universe $L=\{u,1,a,b,0\}$, $\sigma_{\mathfrak{L}_{\nbI}}$ equal to $\sigma_{\mathfrak{L}}$ for any $\sigma\in\{\vee, \wedge, \rightarrow, \uparrow\}$ and negation defined by the table below; we drop $\mathfrak{L}_{\nbI}$ from indexing the operations for simplicity.

\begin{figure}[H]
\centering
\begin{tabular}{l|ccccr}
$x$ & $u$ & $1$ & $a$ & $b$ & $0$\\ \hline
$\neg x$ & $1$ & $0$ & $b$ & $a$ & $1$
\end{tabular}
\caption*{Table for Negation}
\end{figure}

We then define the logical matrix $\mathfrak{M}_{\nbI}=(\mathfrak{L}_{\nbI}, D)$, with $D=\{u,1\}$. Since the operations associated to the symbols in $\{\vee, \wedge, \rightarrow, \uparrow\}$ are the same as those of $\mathfrak{L}$, we see that $\mathfrak{M}_{\nbI}$ models the axiom schemata $\textbf{Ax\: 1}$ through $\textbf{Ax\: 8}$ of $\bI$, plus $\textbf{Ax\: 9}^{*}$, $\textbf{Ip}$ and $\textbf{Comm}$, since those axiom schemata involve only the connectives on $\Sigma_{\bI}$.

The proof that $\mathfrak{M}_{\nbI}$ models $\textbf{Ax\: 11}^{*}$ is not necessary since it corresponds to the axiom scheme number $3$ for $\textbf{C}_{1}$ in the axiomatization found in \cite{Lewin}, and we use the same negation.

The proof that there are only two congruences on $\mathfrak{L}_{\nbI}$ goes as expected, beginning with the supposition that a pair $(x,y)$, with $x\neq y$, is in $\theta$, what then implies $\theta=\nabla$.

Finally, we trivially find both $F_{a}=\{u,1,a\}$ and $F_{b}=\{u,1,b\}$ are $\nbI$-filters whose largest compatible congruence is $\Delta$, what proves $\nbI$ is not-algebraizable according to Blok and Pigozzi.

\subsection{$\nbIciw$, $\nbIci$ and $\nbIcl$ are not algebraizable by Blok and Pigozzi}

We state that $\mathfrak{M}_{\nbI}$ also models $\nbIciw$, $\nbIci$ and $\nbIcl$, and since $F_{a}$ and $F_{b}$ are still, respectively, $\mathcal{L}$-filters, for $\mathcal{L}\in\{\nbIciw, \nbIci, \nbIcl\}$, we prove that none of these three is algebraizable by Blok and Pigozzi. The proof that $\mathfrak{M}$ models $\textbf{ciw}^{*}$, $\textbf{ci}^{*}$ and $\textbf{cl}^{*}$ is in the following tables.

\begin{figure}[H]
\centering
\begin{tabular}{l|cccc}
$x$ & $\neg x$ & $x\wedge\neg x$ & $x\uparrow\neg x$ & $(x\uparrow\neg x)\vee(x\wedge\neg x)$ \\\hline
$u$ &   $1$ &              $1$ &                     $0$ &                                     $1$ \\
$1$ &   $0$ &              $0$ &                     $1$ &                                     $1$ \\
$a$ &   $b$ &              $0$ &                     $1$ &                                     $1$ \\
$b$ &   $a$ &              $0$ &                     $1$ &                                     $1$ \\
$0$ &   $1$ &              $0$ &                     $1$ &                                     $1$  
\end{tabular}
\caption*{Table for $\textbf{ciw}^{*}$}
\end{figure}

\begin{figure}[H]
\centering
\begin{tabular}{l|ccccc}
$x$ & $\neg x$ & $x\wedge\neg x$ & $x\uparrow\neg x$ & $\neg(x\uparrow\neg x)$ & $\neg(x\uparrow\neg x)\rightarrow(x\wedge\neg x)$ \\\hline
$u$ &   $1$ &              $1$ &                     $0$ &                               $1$ &                                                      $1$ \\
$1$ &   $0$ &              $0$ &                     $1$ &                               $0$ &                                                      $1$ \\
$a$ &   $b$ &              $0$ &                     $1$ &                               $0$ &                                                      $1$\\
$b$ &   $a$ &              $0$ &                     $1$ &                               $0$ &                                                      $1$\\
$0$ &   $1$ &              $0$ &                     $1$ &                               $0$  &                                                      $1$
\end{tabular}
\caption*{Table for $\textbf{ci}^{*}$}
\end{figure}

\begin{figure}[H]
\centering
\begin{tabular}{l|ccccc}
$x$ & $\neg x$ & $x\wedge\neg x$ & $x\uparrow\neg x$ & $\neg(x\wedge\neg x)$ & $\neg(x\wedge\neg x)\rightarrow(x\uparrow\neg x)$ \\\hline
$u$ &   $1$ &              $1$ &                     $0$ &                               $0$ &                                                      $1$ \\
$1$ &   $0$ &              $0$ &                     $1$ &                               $1$ &                                                      $1$ \\
$a$ &   $b$ &              $0$ &                     $1$ &                               $1$ &                                                      $1$\\
$b$ &   $a$ &              $0$ &                     $1$ &                               $1$ &                                                      $1$\\
$0$ &   $1$ &              $0$ &                     $1$ &                               $1$  &                                                      $1$
\end{tabular}
\caption*{Table for $\textbf{cl}^{*}$}
\end{figure}

\section{$\bI$ and $\nbI$ are not characterizable by finite Nmatrices}\label{bI is not characterized by Nmatrices section}

It is well known that da Costa's hierarchy is both not algebraizable according to Blok and Pigozzi (\cite{Lewin}, \cite{Mortensen80}) and not characterizable by a finite Nmatrix (\cite{Avron}). We have already shown that $\bI$, and the other systems of incompatibility we have here defined, are not algebraizable according to Blok and Pigozzi, giving them a difficulty to be approached close to that of the logics $C_{n}$; but here, we show furthermore that $\bI$ and $\nbI$ are not characterizable by finite Nmatrices, and are, therefore, specially hard systems from a semantical standpoint. Interestingly, most proofs found so far in this chapter, such as those for the facts that bivaluations and Fidel structures both characterize our logics, or that these systems are not algebraizable according to Blok and Pigozzi, are incredibly similar to the proofs of these results for logics of formal inconsistency, the only difference being that some amount of care must be taking while dealing with incompatibility; however, the proofs in this section require formulas $\phi_{ij}$ specially composed to prove the non-characterizability (by either finite Nmatrices or finite Rmatrices) of systems of incompatibility, and are therefore intrinsically different from any demonstrations produced on the field of $\textbf{LFI}$'s.

We start by supposing there exists an Nmatrix $\mathcal{M}=(\mathcal{A}, D)$ that characterizes $\bI$, with $A$ the universe of $\mathcal{A}$ and $U=A\setminus D$ the set of undesignated elements; we also remember we have defined $\bot_{\alpha\beta}$ as $\alpha\wedge(\beta\wedge(\alpha\uparrow\beta))$ and ${\sim}\alpha$ as $\alpha\rightarrow\bot_{\alpha\alpha}$, since both of these will play an important role in the proof to come. For simplicity, we will drop the indexes from the operations on $\mathcal{A}$ and use the infix notation.

\begin{lemma}\label{basic facts about Nmatrix semantics}
Suppose $d_{1}, d_{2}\in D$ and $u_{1}, u_{2}\in U$:
\begin{enumerate}
\item $d_{1}\wedge d_{2}\subseteq D$, while $d_{1}\wedge u_{1}$, $u_{1}\wedge d_{1}$ and $u_{1}\wedge u_{2}$ are subsets of $U$;
\item$d_{1}\rightarrow d_{2}$, $u_{1}\rightarrow d_{1}$ and $u_{1}\rightarrow u_{2}$ are subsets of $D$, while $d_{1}\rightarrow u_{1}\subseteq U$;
\item $d_{1}\vee d_{2}$, $d_{1}\vee u_{1}$ and $u_{1}\vee d_{1}$ are subsets of $D$, while $u_{1}\vee u_{2}\subseteq U$;
\item for any formulas $\alpha$ and $\beta$ and valuation $\nu:\textbf{F}(\Sigma_{\bI}, \mathcal{V})\rightarrow \mathcal{A}$ for $\mathcal{A}$, $\nu(\bot_{\alpha\beta})\in U$;
\item for any formula $\alpha$ and valuation $\nu$, if $\nu(\alpha)\in D$, then $\nu({\sim}\alpha)\in U$; and if $\nu(\alpha)\in U$, then $\nu({\sim}\alpha)\in D$.
\end{enumerate}
\end{lemma}

\begin{proof}
\begin{enumerate}
\item Given variables $p$ and $q$, from $\textbf{Ax\: 4}$ and $\textbf{Ax\: 5}$ one gets $p\wedge q\vdash_{\bI}p$ and $p\wedge q\vdash_{\bI}q$; given a valuation $\nu$ for $\mathcal{A}$, if $\nu(p)\in U$ or $\nu(q)\in U$, then $\nu(p\wedge q)\in U$, and therefore $\nu(p)\wedge\nu(q)\subseteq U$: this is the case since, otherwise, one could define a valuation $\nu^{*}: \textbf{F}(\Sigma_{\bI}, \mathcal{V})\rightarrow \mathcal{A}$ such that $\nu^{*}(p)=\nu(p)$, $\nu^{*}(q)=\nu(q)$ but $\nu^{*}(p\wedge q)\in D$, meaning that either $p\wedge q\vdash_{\bI}p$ or $p\wedge q\vdash_{\bI}q$ is not validated by $\mathcal{M}$.

From $\textbf{Ax\: 3}$ and the deduction meta-theorem, we obtain $p, q\vdash_{\bI}p\wedge q$, meaning that if $\nu(p), \nu(q)\in D$, then $\nu(p\wedge q)\in D$, and so $\nu(p)\wedge\nu(q)\subseteq D$.

\item We have that, for variables $p$ and $q$, $q\rightarrow (p\rightarrow q)$ is an instance of an axiom schema of $\bI$, and therefore $q\vdash_{\bI}p\rightarrow q$; taking a homomorphism $\nu:\textbf{F}(\Sigma_{\bI}, \mathcal{V})\rightarrow \mathcal{A}$ such that $\nu(q)\in D$, we must have $\nu(p\rightarrow q)\in D$ and therefore $\nu(p)\rightarrow \nu(q)\subseteq D$. This, of course, corresponds to both $d_{1}\rightarrow d_{2}$ and $u_{1}\rightarrow d_{1}$ being subsets of $D$.

Now, $p, p\rightarrow q\vdash_{\bI} q$, and therefore if $\nu(q)\in U$, we must have either $\nu(p)\in U$ or $\nu(p\rightarrow q)\in U$; so, if $\nu(q)\in U$ and $\nu(p)\in D$, one must necessarily have $\nu(p\rightarrow q)\in U$, meaning $\nu(p)\rightarrow\nu(q)\subseteq U$, corresponding to $d_{1}\rightarrow u_{1}\subseteq U$. 

Finally, suppose $\nu(p), \nu(q)\in U$ and, without loss of generality, $\nu(r)\in D$; from $\textbf{Ax\: 2}$ and the deduction meta-theorem, $p\rightarrow (q\rightarrow r)\vdash_{\bI}(p\rightarrow q)\rightarrow(p\rightarrow r)$ and, from the previous investigations, one finds that $\nu(q\rightarrow r)\in D$ and $\nu(p\rightarrow(q\rightarrow r))\in D$, meaning $\nu((p\rightarrow q)\rightarrow(p\rightarrow r))\in D$. Since $\nu(p\rightarrow r)\in D$, one necessarily obtains that $\nu(p\rightarrow q)\in D$, forcibly implying that $\nu(p)\rightarrow\nu(q)\subseteq D$, what corresponds to $u_{1}\rightarrow u_{2}\subseteq D$.

\item Given variables $p$ and $q$, we have that $p\vdash_{\bI}p\vee q$ and $q\vdash_{\bI}p\vee q$ from axiom schemata $\textbf{Ax\: 6}$ and $\textbf{Ax\: 7}$ and the deduction meta-theorem. So, if $\nu(p)\in D$ or $\nu(q)\in D$, one necessarily finds $\nu(p\vee q)\in D$, and therefore $\nu(p)\rightarrow\nu(q)\subseteq D$.

Now, from axiom schema $\textbf{Ax\: 8}$ and the deduction meta-theorem, $p\rightarrow r, q\rightarrow r\vdash_{\bI}(p\vee q)\rightarrow r$; if $\nu(p)\in U$ and $\nu(q)\in U$, suppose, without loss of generality, that $\nu(r)\in U$; this means $\nu(p\rightarrow r)$ and $\nu(q\rightarrow r)$ are both in $D$, and so must be $\nu((p\vee q)\rightarrow r)$, what only happens if $\nu(p\vee q)\in U$. This, of course, implies that $\nu(p)\vee\nu(q)\subseteq U$.

\item Let $p$, $q$ and $r$ be propositional variables: from $\textbf{Ip}$ and the deduction meta-theorem one finds that $p, q, p\uparrow q\vdash_{\bI}r$, yet $p, q\not\vdash_{\bI}r$.\footnote{To see that, take a bivaluation $\nu$ for $\bI$ with $\nu(p)=\nu(q)=1$ and $\nu(r)=0$} If we suppose $\nu$ is a valuation for $\mathcal{A}$ such that $\nu(p), \nu(q)\in D$ and $\nu(r)\in U$, it becomes clear that one must have $\nu(p\uparrow q)\in U$, since otherwise one would forcibly have $\nu(r)\in D$. This means we can never have all three $\nu(\alpha)$, $\nu(\beta)$ and $\nu(\alpha\uparrow\beta)$ in $D$, and therefore 
\[\nu(\bot_{pq})=(\nu(p)\wedge\nu(q))\wedge\nu(p\uparrow q)\in U.\]
Of course, for arbitrary formulas $\alpha$ and $\beta$, $\bot_{\alpha\beta}\vdash_{\bI}\bot_{pq}$ and $\bot_{pq}\vdash_{\bI}\bot_{\alpha\beta}$, meaning that, for any $\nu$, $\nu(\bot_{\alpha\beta})\in U$.

\item If $\nu(\alpha)\in D$, since $\nu(\bot_{\alpha\alpha})\in U$ we obtain $\nu({\sim}\alpha)=\nu(\alpha\rightarrow\bot_{\alpha\alpha})\in U$. Reciprocally, if $\nu(\alpha)\in U$, $\nu({\sim}\alpha)=\nu(\alpha\rightarrow\bot_{\alpha\alpha})\in D$.
\end{enumerate}
\end{proof}

\begin{lemma}
\begin{enumerate}
\item For any two elements $a, b\in A$, either $a\uparrow b\subseteq D$ or $a\uparrow b\subseteq U$;
\item for any two elements $a, b\in A$, either both $a\uparrow b$ and $b\uparrow a$ are subsets of $D$, or both are subsets of $U$.
\end{enumerate}
\end{lemma}

\begin{proof}
Suppose that there are values $d, u\in a\uparrow b$ such that $d\in D$ and $u\in U$, and let $p$ and $q$ be propositional variables throughout the proof. Since $b\uparrow a$ is necessarily not empty, it must contain either an element $d^{*}\in D$, or an $u^{*}\in U$.
\begin{enumerate}
\item In the first case, take a valuation $\nu:\textbf{F}(\Sigma_{\bI}, \mathcal{V})\rightarrow\mathcal{A}$ satisfying $\nu(p)=a$, $\nu(q)=b$, $\nu(p\uparrow q)=u$ and $\nu(q\uparrow p)=d^{*}$, and then we have that $\nu((q\uparrow p)\rightarrow(p\uparrow q))\in d^{*}\rightarrow u$, which by Lemma \ref{basic facts about Nmatrix semantics} is contained in $U$. 

This shows $\nu$ does not validate $\textbf{Comm}$, what is absurd given that $\mathcal{M}$ characterizes $\bI$.
\item In the second case, take now a valuation $\nu$ with $\nu(p)=a$, $\nu(q)=b$, $\nu(p\uparrow q)=d$ and $\nu(q\uparrow p)=u^{*}$. Then $\nu((p\uparrow q)\rightarrow(q\uparrow p))\in d\rightarrow u^{*}$, again contained in $U$ according to Lemma \ref{basic facts about Nmatrix semantics}.
\end{enumerate}

Either way we reach a contradiction, an the conclusion must be that either $a\uparrow b$ is contained in $D$, or it is contained in $U$.

Finally, suppose that there are values $a, b\in A$ with $a\uparrow b\subseteq D$ and $b\uparrow a\subseteq U$, and then, for any valuation $\nu$ with $\nu(p)=a$ and $\nu(q)=b$ (and there are many of them), one necessarily finds that $\nu(p\uparrow q)\in D$ and $\nu(q\uparrow p)\in U$, and from Lemma \ref{basic facts about Nmatrix semantics} we have $\nu((p\uparrow q)\rightarrow(q\uparrow p))\in U$; this again contradicts the fact that $\mathcal{M}$ characterizes $\bI$, since it implies that $\nu$ does not model $\textbf{Comm}$.
\end{proof}

The following fact is trivial, but important in the following discussion so we make a point of proving it: suppose $\Gamma$ is a set of formulas of $\bI$ such that there exists $\varphi$ with the property that $\Gamma\not\vdash_{\bI}\varphi$; then we state that there exists a valuation $\nu$ for $\mathcal{A}$ with $\nu(\Gamma)\subseteq D$ but $\nu(\varphi)\in U$. This is true since, otherwise, if we had that for every valuation $\nu$ for $\mathcal{A}$ satisfying $\nu(\Gamma)\subseteq D$ one also had $\nu(\varphi)\in D$, this would imply $\Gamma\vDash_{\mathcal{M}}\varphi$. Since $\mathcal{M}$ characterizes $\bI$, this would mean $\Gamma\vdash_{\bI}\varphi$, against our suppositions.

Consider then two disjoint sets of distinct variables $\{p_{n}\ :\ n\in\mathbb{N}\}$ and $\{q_{n}\ :\ n\in\mathbb{N}\}$ and the formulas, for $i, j\in\mathbb{N}$,
\[\phi_{ij}=\begin{cases*}
      \quad p_{i}\uparrow q_{j} & if $i<j$ \\
      {\sim}(p_{i}\uparrow q_{j})   & otherwise
    \end{cases*};\]
we also define, for any $n\in\mathbb{N}$, 
\[\Gamma_{n}=\{\phi_{ij}\ :\  0\leq i, j\leq n\}.\]

\begin{figure}[H]
\centering
\begin{tabular}{cccc}
${\sim}(p_{0}\uparrow q_{0})$ & $p_{0}\uparrow q_{1}$ & $\cdots$ & $p_{0}\uparrow q_{n}$\\
${\sim}(p_{1}\uparrow q_{0})$ & ${\sim}(p_{1}\uparrow q_{1})$ & $\cdots$ & $p_{1}\uparrow q_{n}$\\
$\vdots$ & $\vdots$ & $\ddots$ & $\vdots$\\
${\sim}(p_{n}\uparrow q_{0})$ & ${\sim}(p_{n}\uparrow q_{1})$ & $\cdots$ & ${\sim}(p_{n}\uparrow q_{n})$
\end{tabular}
\caption*{The formulas in $\Gamma_{n}$}
\end{figure}

We state that, for all $n$, $\Gamma_{n}\not\vdash_{\bI}p_{0}$: to see that, we consider the bivaluation $\nu$ for which $\nu(p_{i})=\nu(p_{j})=0$, for all $i, j\in\mathbb{N}$, and 
\[\nu(p_{i}\uparrow q_{j})=\begin{cases*}
      1 & if $i<j$ \\
      0   & otherwise
    \end{cases*}.\]
Such a bivaluation is possible since: first of all, given all $p_{i}$ and $q_{j}$ are mapped into $0$, $p_{i}\uparrow q_{j}$ may assume any value; second, assigning a value to $p_{i}\uparrow q_{j}$ does not interfere with the value assigned to $p_{k}\uparrow q_{l}$, for $(i,j)\neq(k,l)$, given that the only restriction imposed by the fact $\nu(p_{i}\uparrow q_{j})$ has a certain value is that $\nu(q_{j}\uparrow p_{i})$ has the same value. Furthermore, Lemma \ref{basic facts about Nmatrix semantics} tells us that $\nu({\sim}\alpha)=1$ if, and only if, $\nu(\alpha)=0$, giving us $\nu(\phi_{ij})=1$ for all $0\leq i, j\leq n$; so $\nu(\Gamma_{n})=1$ but $\nu(p_{0})=0$, and since bivaluations characterize $\bI$, we have proved $\Gamma_{n}\not\vdash_{\bI}p_{0}$.

This means that, for any $n\in\mathbb{N}$, there exists a valuation $\nu$ for $\mathcal{A}$ such that $\nu(\Gamma_{n})\subseteq D$ but $\nu(p_{0})\in U$; however, we will prove that, assuming $\mathcal{A}$ has a finite universe with $n$ elements, there can not exist a valuation $\nu$ satisfying $\nu(\Gamma_{n})\subseteq D$. This is rather simple: suppose there exists such a valuation $\nu$; given there are $n+1$ elements among $p_{0}, p_{1}, \dotsc  , p_{n}$ and $\mathcal{A}$ has only $n$ elements, by the pigeonhole principle one finds there exist $0\leq i<j\leq n$ such that $\nu(p_{i})=\nu(p_{j})$.

Since $\nu(\Gamma_{n})\subseteq D$, $\nu(\phi_{ij})=\nu(p_{i}\uparrow q_{j})\in D$, implying $\nu(p_{i})\uparrow\nu(q_{j})\subseteq D$. Nevertheless, $\nu(p_{i})=\nu(p_{j})$ gives us that 
\[\nu(p_{j}\uparrow q_{j})\in \nu(p_{j})\uparrow\nu(q_{j})=\nu(p_{i})\uparrow\nu(q_{j})\subseteq D,\]
and therefore $\nu(\phi_{jj})=\nu({\sim}(p_{j}\uparrow q_{j}))\in U$, contradicting the supposition that $\nu(\Gamma_{n})\subseteq D$.

\begin{theorem}
There exists no finite Nmatrix which characterizes $\bI$.
\end{theorem}

Of course, this also implies $\nbI$ is not characterizable by finite Nmatrices: if it were characterizable by some $\mathcal{M}=(\mathcal{A}, D)$, for $\mathcal{A}$ a $\Sigma_{\nbI}$-multialgebra with universe $A$, by ignoring the paraconsistent negation and defining $\mathcal{A}_{-}=(A, \{\sigma_{\mathcal{A}}\}_{\sigma\in\Sigma_{\bI}})$ one would find that $\mathcal{M}_{-}=(\mathcal{A}_{-}, D)$ characterizes $\bI$, contradicting our previous theorem.

\subsection{$\bI$ is not characterizable by a finite Rmatrix}

Although Rmatrices haven't been as well studied as logical matrices, or even Nmatrices, one could ask themselves whether Rmatrix semantics wouldn't offer easier decision methods for $\bI$ and other logics of incompatibility than RNmatrices; so we stop for a moment and show that, already for $\bI$, it is not possible to find a characterizing finite Rmatrix.

Assume $\mathcal{M}=(\mathcal{A}, D, \mathcal{F})$ is a finite Rmatrix (meaning the universe $A$ of $\mathcal{A}$ is finite) which characterizes $\bI$, and consider again: two disjoint sets of distinct variables $\{p_{n}\ :\ n\in\mathbb{N}\}$ and $\{q_{n}\ :\ n\in\mathbb{N}\}$; the formulas, for $i, j\in\mathbb{N}$,
\[\phi_{ij}=\begin{cases*}
      \quad p_{i}\uparrow q_{j} & if $i<j$ \\
      {\sim}(p_{i}\uparrow q_{j})   & otherwise
    \end{cases*};\]
and the sets of formulas, for any $n\in\mathbb{N}$, $\Gamma_{n}=\{\phi_{ij}\ :\  0\leq i, j\leq n\}$.

\begin{lemma}\label{basic facts about Rmatrix semantics}
If the Rmatrix $\mathcal{M}$ characterizes $\bI$ and the image of $\nu\in\mathcal{F}$ is not contained in $D$, then $\nu(\alpha)$ and $\nu({\sim}\alpha)$ can not both belong to $D$.
\end{lemma}

\begin{proof}
We know that, since "$\sim$" behaves classically, it satisfies 
\[\alpha\rightarrow({\sim}\alpha\rightarrow\beta)\]
for any formula $\beta$, implying by the deduction meta-theorem that $\alpha, {\sim}\alpha\vdash_{\bI}\beta$. Let $\gamma$ be a formula such that $\nu(\gamma)\notin D$: if $\nu(\alpha), \nu({\sim}\alpha)\in D$, then $\nu$ would not validate the deduction $\alpha, {\sim}\alpha\vdash_{\bI}\gamma$, and therefore at most one among $\nu(\alpha)$ and $\nu({\sim}\alpha)$ belongs to $D$.
\end{proof}

As we know, $\Gamma_{n}\not\vdash_{\bI}p_{0}$ for any $n\in\mathbb{N}$, which means there must exist a valuation $\nu\in\mathcal{F}$ with $\nu(\Gamma_{n})\subseteq D$ but $\nu(\varphi)\not\in D$, since otherwise one would have, for all $\nu\in\mathcal{F}$, satisfying that $\nu(\Gamma_{n})\subseteq D$, that $\nu(\varphi)\in D$; this, of course, would show that $\Gamma_{n}\vDash_{\mathcal{M}}p_{0}$ and, since $\mathcal{M}$ characterizes $\bI$, that $\Gamma_{n}\vdash_{\bI}p_{0}$, which we know not to be true. So there exists a $\nu\in\mathcal{F}$ with undesignated elements in its image and $\nu(\Gamma_{n})\subseteq D$. If $A$ has cardinality $n$, by the pigeonhole principle there must exist two elements $p_{i}$ and $p_{j}$ among $\{p_{n}\ :\  n\in\mathbb{N}\}$ such that $\nu(p_{i})=\nu(p_{j})$, what leads us to the following problem: assume, without loss of generality, $i<j$; then 
\[\nu(p_{i}\uparrow q_{j})=\nu(p_{i})\uparrow\nu(q_{j})=\nu(p_{j})\uparrow\nu(q_{j})=\nu(p_{j}\uparrow q_{j}).\]
Since $\nu(\Gamma_{n})\subseteq D$, $\nu(p_{i}\uparrow q_{j}), \nu({\sim}(p_{j}\uparrow q_{j}))\in D$, which means $\nu$ is a valuation in $\mathcal{F}$ (with image not contained in $D$) satisfying that $\nu(p_{j}\uparrow q_{j})$ and $\nu({\sim}(p_{j}\uparrow q_{j}))$ are both in $D$, a contradiction given Lemma \ref{basic facts about Rmatrix semantics}.

The necessary conclusion is that no finite Rmatrix can characterize $\bI$. A similar argument applies to $\nbI$.

\section{Comparing matrix semantics}

So, once we introduced the notion of an RNmatrix, or restricted Nmatrix, extending the notions of Rmatrix and Nmatrix, it becomes clear that we are considering generalizations of a logical matrix in two distinct, and independent, directions: we extend matrices by allowing their underlying algebras to become multialgebras, hence exchanging deterministic operations for multioperations, non-deterministic operations; and in a second direction, we restrict ourselves to a subset of all valuations, selecting which homomorphisms to take into consideration.

This gives the semantics at hand, if one interprets independent directions as it would be done in linear algebra, two dimensions of generality. Many of such dimensions are already quite standard in the literature, if we take as starting point that a logical matrix should merely be a pair $\mathcal{M}=(\mathcal{A}, D)$, with $\mathcal{A}$ a finite algebra; we list a few ways below to extend such a matrix.
\begin{enumerate}
\item Allowing $\mathcal{A}$ to become an infinite algebra, that is, an algebra with an infinite universe (and, reciprocally, if the underlying algebraic structure of the matrix is actually a multialgebra, allowing it to have an infinite universe). This generalization turns a logical matrix from a rather strict notion of a decision method to a broader semantical approach having universes of any desired cardinality.

\item Considering, instead of an unique matrix $\mathcal{M}$, a class $\mathbb{M}$ of matrices over which we define a semantical deduction operator by $\Gamma\vDash_{\mathbb{M}}\varphi$ if and only if $\Gamma\vDash_{\mathcal{M}}\varphi$ for all $\mathcal{M}$ in $\mathbb{M}$; these matrices may be of any desired type, but it was already established by W\'ojcicki \cite{Woj, Woj2}, who envisioned this generalization, that every tarskian logic is characterizable by a class of potentially infinite logical matrices.

\item Permitting that multiple different sets of distinguished elements coexist: meaning a pair $\mathcal{M}=(\mathcal{A}, \{D_{\lambda}\}_{\lambda\in\Lambda})$, with $\bigcup_{\lambda\in\Lambda}D_{\lambda}$ contained in $\mathcal{A}$'s universe, when we define $\Gamma\vDash_{\mathcal{M}}\varphi$ if, and only if, for every valuation $\nu$ and index $\lambda$ one has that $\nu(\Gamma)\subseteq D_{\lambda}$ implies $\nu(\varphi)\in D_{\lambda}$. Again a development by W\'ojcicki, usually referred to as generalized matrices\index{Matrix, Generalized} or Gmatrices\index{Gmatrix}, he proved in \cite{Woj, Woj2} that every tarskian logic can also be characterized by a potentially infinite matrix with multiple sets of distinguished elements.
\end{enumerate}

Of course, the ideal scenario remains that of a simple, old-fashioned finite logical matrix, and the topic of how we are able to navigate the space of possible generalizations is a rather interesting one. Look at the example of the logic of incompatibility $\bI$: we know that a restricted Nmatrix $\textbf{2}(\bI)$ (Theorem \ref{2-valued RNmatrix}) with two elements may characterize it, as well as the previously defined $\mathbb{2}_{\bI}$ (Section \ref{Decision Method for bI}) based on bivaluations, also with only two elements. But the differences between the two are rather obvious: $\mathbb{2}_{\bI}$ is clearly more desirable as a semantical tool and proves to be a very efficient decision method; additionally, $\mathbb{2}_{\bI}$ has far more restrictive operations, and in exchange is far more lenient on the subject of which homomorphisms to consider.

\begin{figure}[H]
\centering
\begin{tikzpicture}
\begin{axis}[axis lines = left, 
xlabel={Restrictions on the Operations}, 
ylabel={Restrictions on the homomorphisms}, 
axis equal, 
xmin=0, 
xmax=5, 
ymin=0, 
ymax=5, 
xtick={4},
xticklabels={Deterministic},
ytick={1},
yticklabels={All},
scatter/@pre marker code/.code={},
scatter/@post marker code/.code={\node [above right] {\pgfplotspointmeta};}]
\draw (1,1) rectangle (4,4);
\draw [line width=1pt, line cap=round, dash pattern=on 0pt off 2\pgflinewidth] (0,1) -- (5,1);
\draw [line width=1pt, line cap=round, dash pattern=on 0pt off 2\pgflinewidth] (4,0) -- (4,5);
\addplot[scatter,only marks]
    plot[scatter src=explicit symbolic]
    coordinates {
(1,1) [Nmatrices]
(1,4) [RNmatrices]
(4,1) [logical matrices]
(4,4) [Rmatrices]
};
\draw [->, line width=1pt, line cap=round, dash pattern=on 0pt off 2\pgflinewidth] (1.25,3.75) -- (3.75, 1.25) node[midway,above right] {Ideal};
\end{axis}

\end{tikzpicture}
\end{figure}

Using the previous diagram, we can see that, in $\bI$'s case, we are able to start at some point in the far top-left, where $\textbf{2}(\bI)$ lies, and navigate to the bottom-right without ever leaving the, at least on the diagram, loosely-defined environment of restricted non-deterministic matrices, arriving at $\mathbb{2}_{\bI}$; notice this is very much in line with what we previously observed, that going from the first of these RNmatrices to the second involves restricting the underlying operations while expanding the accepted homomorphisms. Worthy of notice is that, while working in $\bI$, we are never able to leave RNmatrices since $\bI$ can not be characterized by either (finite) logical matrices, Nmatrices or Rmatrices.

In a different case, consider $\bI^{-}$: we know that it admits a finite RNmatrix, namely $\textbf{2}(\bI^{-})$, located somewhere on the top-left region of the schematic diagram; but, in this case, we also have a finite Nmatrix which characterizes $\bI^{-}$, although no (finite) logical matrices are available as far as we know. So, at least in this logic's case, navigating ever bottom-right indeed leads outside the environment of RNmatrices and into the methodology of Nmatrices.

Now, as we have assigned numbers to the other generalizations of a finite logical matrix, let us take $(4)$ as exchanging the underlying algebra of a matrix for a multialgebra (making a logical matrix into an Nmatrix, an Rmatrix into an RNmatrix and so on), and $(5)$ as restricting the valuations to be taken into consideration (turning a logical matrix into an Rmatrix and so on).

\begin{problem}
Which combinations of generalizations of a (finite) logical matrix have the expressive power to characterize every tarskian logic?
\end{problem}

W\'ojciki proved, in \cite{Woj} and \cite{Woj2}, that the combinations of $(1)$ and $(2)$, and $(1)$ and $(3)$, corresponding to classes of potentially infinite logical matrices and potentially infinite matrices with multiple sets of distinguished elements, are enough to characterize all tarskian logics; Piochi, in \cite{Piochi}, proved that $(1)$ plus $(5)$ is also enough, and as we have proved that every tarskian logic is characterizable by a two-elements RNmatrix (again Theorem \ref{2-valued RNmatrix}), $(4)$ plus $(5)$ is also enough. From this, it is clear that combining any four of these conditions characterizes all such logics, although it seems the same can not yet be said about combinations of three of them, as, for example, we were not able to find a reference concerning $(2)$ plus $(3)$ and $(5)$ or any subset thereof. So, we have the following table to summarize these observations.

\begin{figure}[H]
\centering
\begin{tabular}{cc|ccccc}
& & \multicolumn{5}{c}{Generalizations} \\[5pt]
     & & (1)    &  (2)   & (3)    & (4)    & (5)   \\[5pt] \cline{2-7}
\multirow{5}{*}{\rotatebox[origin=c]{90}{Generalizations}}
&  (1)  &  &  \checkmark &  \checkmark &  \checkmark &  ? \\[5pt]
 & (2)  & \checkmark &   &  ? &  ? &  ? \\[5pt]
  & (3)  & \checkmark &  ? &  &  ? &  ? \\[5pt]
&   (4) & \checkmark &  ? &  ? &   &  \checkmark \\[5pt]
  & (5) & ? &  ? &  ? &  \checkmark &  \\[5pt]
\end{tabular}
\end{figure}

We have left the diagonal of the table empty, as it does not deal with combinations of generalizations \textit{per se}, but rather the generalizations themselves. Of course, there is nothing stopping us from considering which, if any, of the generalizations can characterize all tarskian logics, leading to our second problem.

\begin{problem}
Are there any generalizations which, alone, can characterize all tarskian logics?
\end{problem}

As we have shown, $\bI$ is not characterizable by either finite Nmatrices or finite Rmatrices, so we can already place some limitatitive results in our table.

\begin{figure}[H]
\centering
\begin{tabular}{cc|ccccc}
& & \multicolumn{5}{c}{Generalizations} \\[5pt]
     & & (1)    &  (2)   & (3)    & (4)    & (5)   \\[5pt] \cline{2-7}
\multirow{5}{*}{\rotatebox[origin=c]{90}{Generalizations}}
&  (1)  & ? &  \checkmark &  \checkmark &  \checkmark &  ? \\[5pt]
 & (2)  & \checkmark & ?  &  ? &  ? &  ? \\[5pt]
  & (3)  & \checkmark &  ? & ? &  ? &  ? \\[5pt]
&   (4) & \checkmark &  ? &  ? & \ding{55}   &  \checkmark \\[5pt]
  & (5) & ? &  ? &  ? &  \checkmark & \ding{55} \\[5pt]
\end{tabular}
\end{figure}

Of course, the first course of action one probably thinks of, when finding the problems mentioned above and looking at the tables, is attempting to fill in the missing results; this does not seem impossible, but it does not seem trivial either. Proving that certain semantics can not characterize all tarskian logics involves most probably presenting counter-examples, preferably examples among already known and studied systems.

But one can always increase the rows and columns in our little illustrative tables, by considering other generalizations of logical matrices. Here, for one, we have not included making the operations of a matrix partial, instead of non-deterministic: when applied to Nmatrices, this procedure returns the semantical objects known as PNmatrices. These play a unique role indeed, as they are weaker than RNmatrices, yet we do not know if strictly so or if both semantics characterize the same logics. We are then tempted to consider a hierarchy of strength among combinations of generalizations of logical matrices, instead of the simply binary ``do they characterize all tarskian logics?'', and things start to become involved...

\newpage
\printbibliography[segment=\therefsegment,heading=subbibliography]
\end{refsegment}

\begin{refsegment}
\defbibfilter{notother}{not segment=\therefsegment}
\setcounter{chapter}{8}
\chapter{Translating paraconsistent logics}\label{Chapter9}\label{Chapter 9}

Take the signature \label{SigmaLFICPL}$\Sigma_{\textbf{LFI}}^{\textbf{CPL}}$ such that $(\Sigma_{\textbf{LFI}}^{\textbf{CPL}})_{0}=\{\bot, \top\}$, $(\Sigma_{\textbf{LFI}}^{\textbf{CPL}})_{1}=\{\neg, {\sim}, \circ\}$, $(\Sigma_{\textbf{LFI}}^{\textbf{CPL}})_{2}=\{\vee, \wedge, \rightarrow\}$ and $(\Sigma_{\textbf{LFI}}^{\textbf{CPL}})_{n}=\emptyset$ for $n>2$.

We classically define a Fidel structure to be a $\Sigma_{\textbf{LFI}}^{\textbf{CPL}}-$multialgebra $\mathcal{E}=(A, \{\sigma_{\mathcal{E}}\}_{\sigma \in\Sigma_{\textbf{LFI}}^{\textbf{CPL}}})$ such that:
\begin{enumerate}
\item $(A, \{\sigma_{\mathcal{E}}\}_{\sigma\in\Sigma^{\textbf{CPL}}})$ is a Boolean algebra;
\item for every $a\in A$ and $b\in \neg a$, $a\vee b=\top$;
\item for every $a\in A$ and $b\in \neg a$, there exists a non-empty subset $O_{ab}$ of $\circ a$, defined case by case, designed to capture the logical structure intended to be emulated by the Fidel structure.
\end{enumerate}

Intuitively, one looks at $\neg a$ as all possible negations of $a$, and at $O_{ab}$ as all possible consistencies for $a$ given that $b$ is its negation. To give one example, in $\textbf{mbC}$, where we usually denote $O_{ab}$ by $BC_{ab}$, we require that
\[a\wedge (b\wedge c)=\bot,\quad \forall b\in \neg a,\quad \forall c\in BC_{ab}.\]

When looking at previous instances in this text of a "Fidel structure", maybe the most important distinction was the intuitive replacement of consistency for incompatibility: so, for example, instead of the axiom schema $\textbf{bc1}$ of $\textbf{mbC}$ given by $\circ\alpha\rightarrow(\alpha\rightarrow(\neg\alpha\rightarrow\beta))$, we used a similar, but distinct, schema $\textbf{Ip}$, given by $(\alpha\uparrow\beta)\rightarrow(\alpha\rightarrow(\beta\rightarrow\gamma))$.

One clearly notices how our binary incompatibility operator, the generalized Sheffer's stroke, as studied above seems inherently distinct from the unary consistency operator "$\circ$". In fact, one must translate accordingly the axiomatization of $\textbf{LFI}$'s to that of our $\textbf{LIp}$'s to show that the latter generalize the former, and logics such as $\textbf{mbC}$, $\textbf{mbCciw}$, $\textbf{mbCci}$ and $\textbf{mbCcl}$ become sublogics of, respectively, $\nbI$, $\nbIciw$, $\nbIci$ and $\nbIcl$. 

But, and this is the important finding of this chapter, our intuition can still be validated: inconsistency, at least as found in the simpler paraconsistent logics here exhibited, can be obtained from incompatibility, well-behaved formulas being precisely those formulas incompatible with their negations. Since there are systems on which incompatibility clearly does not reduce back to inconsistency, such as $\textbf{bI}$ (which does not even have a paraconsistent negation), we must reach the conclusion that incompatibility appears to non-trivially generalize the notion of inconsistency, giving some extra validation to the work we performed so far. 

The developments found here have been submitted, as a preprint, in \cite{Frominconsistency}.

\section{Preliminaries}

So, let us define a translation from the usual signature $\Sigma_{\textbf{LFI}}$ for $\textbf{LFI}'s$, with $(\Sigma_{\textbf{LFI}})_{1}=\{\neg , \circ\}$, $(\Sigma_{\textbf{LFI}})_{2}=\{\vee, \wedge, \rightarrow\}$ and $(\Sigma_{\textbf{LFI}})_{n}=\emptyset$ for $n\notin\{1,2\}$, to the signature $\Sigma_{\nbI}$. For $\mathcal{V}$ a countable set of propositional variables, consider the function\label{T} 
\[T:F(\Sigma_{\textbf{LFI}}, \mathcal{V})\rightarrow F(\Sigma_{\nbI}, \mathcal{V})\]
such that:
\begin{enumerate}
\item $T(p)=p$ for every $p\in \mathcal{V}$;
\item $T(\neg\alpha)=\neg T(\alpha)$;
\item $T(\alpha\#\beta)=T(\alpha)\# T(\beta)$ for every $\#\in\{\vee, \wedge, \rightarrow\}$;
\item $T(\circ\alpha)=T(\alpha)\uparrow\neg T(\alpha)$.
\end{enumerate}
Essentially, $T$ changes all occurrences of the form $\circ\alpha$ to $\alpha\uparrow\neg\alpha$. What one does is then to consider for a set of axiom schemata $\Gamma$ of an $\textbf{LFI}$ its translation $T(\Gamma)=\{T(\psi)\ :\  \psi\in \Gamma\}$: to give one example, many instances of $\textbf{Ip}$ are translations of instances of $\textbf{bc1}$.

\begin{proposition}
$T$ is an injective function.
\end{proposition}

\begin{proof}
Suppose $T(\alpha)=T(\beta)$: we proceed by double induction on the orders of $\alpha$ and $\beta$, to show that in this case $\alpha=\beta$. 

If $\alpha$ is of order $0$, then $\alpha=p$ for some $p\in \mathcal{V}$ and therefore $T(\alpha)=\alpha$: if $\beta$ is of order $0$, then it is a propositional variable $q$, and then again $T(\beta)=\beta$, implying we have $\alpha=T(\alpha)=T(\beta)=\beta$.

So, assume $\beta$ is of order at least $1$, and we show that we can not actually have, in this case, $T(\alpha)=T(\beta)$:
\begin{enumerate}
\item if $\beta=\neg\beta_{0}$, $T(\beta)=\neg T(\beta_{0})=T(\alpha)$, which is absurd since $T(\alpha)$ is a propositional variable and can not contain an unary connective;
\item if $\beta=\beta_{0}\#\beta_{1}$, for $\#\in\{\vee, \wedge, \rightarrow\}$, $T(\beta)=T(\beta_{0})\# T(\beta_{1})=T(\alpha)$, which is again absurd since $T(\alpha)$ is a propositional variable;
\item if $\beta=\circ\beta_{0}$, $T(\beta)=T(\beta_{0})\uparrow\neg T(\beta_{0})$, which clearly contains even more than one connective and therefore can not equal $T(\alpha)$.
\end{enumerate}

Now, suppose that for every $\alpha$ of order at most $m$ we have that, if $T(\alpha)=T(\beta)$, then $\alpha=\beta$, and take a formula $\alpha$ of order $m+1$: we have that either $\alpha=\alpha_{0}\#\alpha_{1}$, for $\#\in\{\vee, \wedge, \rightarrow\}$, $\alpha=\neg\alpha_{0}$ or $\alpha=\circ\alpha_{0}$, with the orders of $\alpha_{0}$ and $\alpha_{1}$ being at most $m$.

If $\beta$ is of order $0$, then $\beta=q$ for some $q\in \mathcal{V}$, and in that case $T(\beta)=\beta$ is a formula of order $0$: we then can not really have $T(\alpha)=T(\beta)$, since $T(\alpha)$ is either $\neg T(\alpha_{0})$, $T(\alpha_{0})\# T(\alpha_{1})$ or $T(\alpha_{0})\uparrow\neg T(\alpha_{0})$, none of which is of order $0$; by vacuity, the lemma holds. 

Inductively, suppose that for all $\beta$ of order at most $n$, if $T(\alpha)=T(\beta)$ then $\alpha=\beta$, and take $\beta$ of order $n+1$: we have that either $\beta=\beta_{0}\ast\beta_{1}$, for $\ast\in\{\vee, \wedge, \rightarrow\}$, $\beta=\neg\beta_{0}$ or $\beta=\circ\beta_{0}$, with the orders of $\beta_{0}$ and $\beta_{1}$ being at most $n$.

If $\alpha=\alpha_{0}\#\alpha_{1}$, $T(\alpha)=T(\alpha_{0})\# T(\alpha_{1})$: in this case, $\beta$ can not equal $\neg\beta_{0}$, for in this case we would have $T(\beta)=\neg T(\beta_{0})$ which has a leading connective of arity different from that of $T(\alpha)$; for similar reasons we can not have $\beta=\circ\beta_{0}$ or $\beta=\beta_{0}\ast\beta_{1}$ for $\ast$ different of $\#$. The same can be done when $\alpha=\neg\alpha_{0}$, in which case $\beta$ must also be of the form $\neg\beta_{0}$, and when $\alpha=\circ\alpha_{0}$, when we must have $\beta=\circ\beta_{0}$.

So there remains three cases to check:
\begin{enumerate}
\item if $\alpha=\alpha_{0}\#\alpha_{1}$ and $\beta=\beta_{0}\#\beta_{1}$, $T(\alpha)=T(\beta)$ implies that $T(\alpha_{0})\# T(\alpha_{1})=T(\beta_{0})\# T(\beta_{1})$, and so $T(\alpha_{0})=T(\beta_{0})$ and $T(\alpha_{1})=T(\beta_{1})$; by our induction hypothesis, $\alpha_{0}=\beta_{0}$ and $\alpha_{1}=\beta_{1}$, so that $\alpha=\beta$;
\item if $\alpha=\neg\alpha_{0}$ and $\beta=\neg\beta_{0}$, $T(\alpha)=T(\beta)$ implies that $\neg T(\alpha_{0})=\neg T(\beta_{0})$ and so $T(\alpha_{0})=T(\beta_{0})$; by our induction hypothesis, $\alpha_{0}=\beta_{0}$ and therefore $\alpha=\beta$;
\item if $\alpha=\circ\alpha_{0}$ and $\beta=\circ\beta_{0}$, $T(\alpha)=T(\beta)$ implies that 
\[T(\alpha_{0})\uparrow\neg  T(\alpha_{0})=T(\beta_{0})\uparrow\neg  T(\beta_{0})\]
and so $T(\alpha_{0})=T(\beta_{0})$; by our induction hypothesis, $\alpha_{0}=\beta_{0}$ and therefore $\alpha=\beta$.
\end{enumerate}
This, of course, finishes the proof.
\end{proof}

One important thing to notice is that a formula $\alpha$ in $F(\Sigma_{\nbI},\mathcal{V})$ is not in $T(F(\Sigma_{\textbf{LFI}},\mathcal{V}))$ if and only if it contains a subformula $\beta_{1}\uparrow\beta_{2}$ such that $\beta_{2}\neq\neg\beta_{1}$.

One direction is clear: if $\alpha$ contains a subformula $\beta_{1}\uparrow\beta_{2}$ with $\beta_{2}\neq\neg \beta_{1}$ then it is not in $T(F(\Sigma_{\textbf{LFI}},\mathcal{V}))$, since such a formula is not a translation of anything over the signature $\Sigma_{\textbf{LFI}}$.

Reciprocally, we proceed by induction on the order of $\alpha$: if it is $0$, $\alpha$ is a propositional variable, and it is always the case that $\alpha$ is a translation; so assume that $\alpha$ is of order $n>1$ and that the result holds for formulas of order smaller than $n$. Then we have three cases to consider:
\begin{enumerate}
\item if $\alpha=\alpha_{0}\#\alpha_{1}$, for $\#\in\{\vee, \wedge, \rightarrow\}$, and $\alpha_{0}$ and $\alpha_{1}$ are translations, so is $\alpha$; therefore, if $\alpha\notin T(F(\Sigma_{\textbf{LFI}},\mathcal{V}))$, then one of $\alpha_{0}$ or $\alpha_{1}$ has, by induction hypothesis, a subformula $\beta_{1}\uparrow\beta_{2}$ with $\beta_{2}\neq\neg \beta_{1}$, and the result holds;
\item if $\alpha=\neg\alpha_{0}$ and $\alpha_{0}$ is a translation, so is $\alpha$; therefore, if $\alpha\notin T(F(\Sigma_{\textbf{LFI}},\mathcal{V}))$, then $\alpha_{0}\notin T(F(\Sigma_{\textbf{LFI}},\mathcal{V}))$ and thus has, by induction hypothesis, a subformula of the desired form, making the result hold once again;
\item finally, if $\alpha=\alpha_{0}\uparrow\alpha_{1}$ but $\alpha\notin T(F(\Sigma_{\textbf{LFI}},\mathcal{V}))$, either $\alpha_{0}$ and $\alpha_{1}$ are translations but $\alpha_{1}\neq\neg\alpha_{0}$, and the result holds; or one of $\alpha_{0}$ and $\alpha_{1}$ is not a translation, and therefore contains a subformula $\beta_{1}\uparrow\beta_{2}$ with $\beta_{2}\neq\neg\beta_{1}$, what ends the proof.
\end{enumerate}


\section{$T$ is a conservative translation}

Given logics $\mathcal{L}_{1}$ and $\mathcal{L}_{2}$ over the signatures $\Sigma_{1}$ and $\Sigma_{2}$, a function $\mathcal{T}:F(\Sigma_{1}, \mathcal{V})\rightarrow F(\Sigma_{2},\mathcal{V})$ is said to be a translation\index{Translation}, originally defined in \cite{Itala}, when, for every set of formulas $\Gamma\cup\{\varphi\}$ over the signature $\Sigma_{1}$,
\[\Gamma\vdash_{\mathcal{L}_{1}}\varphi\quad\text{implies}\quad\mathcal{T}(\Gamma)\vdash_{\mathcal{L}_{2}}\mathcal{T}(\varphi);\]
a translation $\mathcal{T}$ is said to be a conservative translation\index{Translation, Conservative} whenever, for every set of formulas $\Gamma\cup\{\varphi\}$ over the signature $\Sigma_{1}$,
\[\mathcal{T}(\Gamma)\vdash_{\mathcal{L}_{2}}\mathcal{T}(\varphi)\quad\text{implies}\quad\Gamma\vdash_{\mathcal{L}_{1}}\varphi.\]
 One may look at definition $2.4.1$ of \cite{ParLog} for a reference for our definition, or the original work concerning conservative translations, \cite{ConservativeTranslationsReference}: such a notion is recurring in a contemporary approach to logic, where systems may have been formulated in apparently distinct ways that prove to be, under translation, equivalent; coincidentally, we have already seen translations very briefly in Example \ref{PTS}.

We shall prove that the function $T$ we previously defined is a translation, and furthermore, a conservative one in many cases. To prove the following lemma, and consequently the following theorem, remember a $\Sigma$-homomorphism $\sigma:F(\Sigma, \mathcal{V})\rightarrow F(\Sigma, \mathcal{V})$ is determined by its action on $\mathcal{V}$, given $F(\Sigma, \mathcal{V})$ is deterministic.

\begin{lemma}\label{Existence of translated substitution}
Given a $\Sigma_{\LFI}$-homomorphism $\sigma:F(\Sigma_{\LFI}, \mathcal{V})\rightarrow F(\Sigma_{\LFI}, \mathcal{V})$, the $\Sigma_{\nbI}$-homomor\-phism $\overline{\sigma}:F(\Sigma_{\nbI}, \mathcal{V})\rightarrow F(\Sigma_{\nbI}, \mathcal{V})$ given by, for a propositional variable $p\in\mathcal{V}$, 
\[\overline{\sigma}(p)=T(\sigma(p))\]
satisfies that, for any formula $\alpha$ on $\Sigma_{\LFI}$, $T(\sigma(\alpha))=\overline{\sigma}(T(\alpha))$.
\end{lemma}

\begin{proof}
The result is trivially true for formulas of order $0$, since they are invariant under $T$. So, proceeding inductively, assume the result holds for the formulas $\alpha$ and $\beta$:
\begin{enumerate}
\item for $\#\in\{\vee, \wedge, \rightarrow\}$, 
\[T(\sigma(\alpha\#\beta))=T(\sigma(\alpha)\#\sigma(\beta))=T(\sigma(\alpha))\# T(\sigma(\beta))=\overline{\sigma}(T(\alpha))\#\overline{\sigma}(T(\beta))=\]
\[\overline{\sigma}(T(\alpha)\# T(\beta))=\overline{\sigma}(T(\alpha\#\beta));\]
\item $T(\sigma(\neg\alpha))=T(\neg\sigma(\alpha))=\neg T(\sigma(\alpha))=\neg \overline{\sigma}(T(\alpha))=\overline{\sigma}(\neg T(\alpha))=\overline{\sigma}(T(\neg\alpha))$;
\item finally, remembering $\overline{\sigma}$ is a homomorphism,
\[T(\sigma(\circ\alpha))=T(\circ\sigma(\alpha))=T(\sigma(\alpha))\uparrow\neg T(\sigma(\alpha))=\overline{\sigma}(T(\alpha))\uparrow\neg\overline{\sigma}(T(\alpha))=\]
\[\overline{\sigma}(T(\alpha))\uparrow\overline{\sigma}(\neg T(\alpha))=\overline{\sigma}(T(\alpha)\uparrow\neg T(\alpha))=\overline{\sigma}(T(\circ\alpha)),\]
what ends the proof.
\end{enumerate}
\end{proof}

\begin{theorem}\label{Soundness of translations}
If $\mathcal{L}$ is a logic over the signature $\Sigma_{\textbf{LFI}}$ with axiom schemata $\Psi$ and $\mathcal{L}^{*}$\label{L*} is the logic over the signature $\Sigma_{\nbI}$ with axiom schemata $T(\Psi)$, then $\Gamma\vdash_{\mathcal{L}} \varphi$ implies that $T(\Gamma)\vdash_{\mathcal{L}^{*}} T(\varphi)$.
\end{theorem}

\begin{proof}
Let $\alpha_{1}, \dotsc , \alpha_{n}$ be a demonstration of $\varphi$ from $\Gamma$, with $\alpha_{n}=\varphi$: we want to show that in this case $T(\alpha_{1}), \dotsc , T(\alpha_{n})$ is a demonstration of $T(\varphi)$; first of all, obviously $T(\alpha_{n})=T(\varphi)$. Then:
\begin{enumerate}
\item if $\alpha_{i}$ is an instance of an axiom schema $\psi$, $T(\alpha_{i})$ is an instance of the axiom schema $T(\psi)$, by an immediate application of Lemma \ref{Existence of translated substitution};
\item if $\alpha_{j}$ is in $\Gamma$, $T(\alpha_{j})\in T(\Gamma)$;
\item finally, given that our only rule of deduction is Modus Ponens, if $\alpha_{k}$ is such that there exist $\alpha_{i}$ and $\alpha_{j}$ with $i, j<k$ and $\alpha_{j}=\alpha_{i}\rightarrow \alpha_{k}$ or $\alpha_{i}=\alpha_{j}\rightarrow\alpha_{k}$, then $T(\alpha_{k})$ is such that either
\[T(\alpha_{j})=T(\alpha_{i}\rightarrow\alpha_{k})=T(\alpha_{i})\rightarrow T(\alpha_{k})\]
or 
\[T(\alpha_{i})=T(\alpha_{j})\rightarrow T(\alpha_{k}).\]
\end{enumerate}
This finishes proving that, if $\Gamma\vdash_{\mathcal{L}} \varphi$, then $T(\Gamma)\vdash_{\mathcal{L}^{*}} T(\varphi)$.
\end{proof}

This proves how $T$ is always a translation.

Notice how $\nbI$ has as axiom schemata all axiom schemata of \label{mbC*}$\textbf{mbC}^{*}$, and therefore is capable of all deductions that $\textbf{mbC}^{*}$ is capable, meaning that, if $\Gamma\vdash_{\textbf{mbC}^{*}}\varphi$, then $\Gamma\vdash_{\nbI}\varphi$. However, $\nbI$ is strictly stronger than $\textbf{mbC}^{*}$, as one can see from, to given one example, the formula
\[(\alpha\uparrow\alpha)\rightarrow(\alpha\rightarrow(\alpha\rightarrow\beta)),\]
which is an instance of $\textbf{Ip}$, but not of $T(\textbf{bc1})$, is a tautology of $\nbI$ but not of $\textbf{mbC}^{*}$. Analogously, $\nbIciw$, $\nbIci$ and $\nbIcl$ are strictly stronger than, respectively, \label{mbCciw*}$\textbf{mbCciw}^{*}$, \label{mbCci*}$\textbf{mbCci}^{*}$ and \label{mbCcl*}$\textbf{mbCcl}^{*}$. One more interesting example would be in $\nbIcl$: although it is well known that $\neg(\alpha\wedge\neg\alpha)$ does not imply $\neg(\neg\alpha\wedge\alpha)$ in this logic, and this remains true in $\mbCcl^{*}$ as we shall prove it, in $\nbIcl$, through use of $\textbf{cl}^{*}$, $\textbf{Comm}$ and $\textbf{Ip}$, one gets
\[\vdash_{\nbIcl}\neg(\alpha\wedge\neg\alpha)\rightarrow\neg(\neg\alpha\wedge\alpha).\]


\subsection{$T$ is conservative for $\textbf{mbC}$}

\begin{definition}
A function $\nu$ from the formulas of $\textbf{mbC}$ to $\{0,1\}$ is said to be a \index{Bivaluation for $\textbf{mbC}$}bivaluation for $\textbf{mbC}$ if it satisfies, for any formulas $\alpha$ and $\beta$,
\begin{enumerate}

\item $\nu(\alpha\wedge\beta)=1$ if and only if $\nu(\alpha)=1$ and $\nu(\beta)=1$;
\item $\nu(\alpha\vee\beta)=1$ if and only if $\nu(\alpha)=1$ or $\nu(\beta)=1$;
\item $\nu(\alpha\rightarrow\beta)=1$ if and only if $\nu(\alpha)=0$ or $\nu(\beta)=1$;
\item $\nu(\neg \alpha)=0$ implies $\nu(\alpha)=1$;
\item $\nu(\circ\alpha)=1$ implies $\nu(\alpha)=0$ or $\nu(\neg\alpha)=0$.
\end{enumerate}
\end{definition}

For $\Gamma\cup\{\varphi\}$ a set of formulas of $\textbf{mbC}$, we define $\Gamma\vDash_{\textbf{mbC}}\varphi$ to mean that for every bivaluation $\nu$ for $\textbf{mbC}$, if $\nu(\Gamma)\subseteq\{1\}$ then $\nu(\varphi)=1$.

\begin{theorem}
For formulas $\Gamma\cup\{\varphi\}$ of $\textbf{mbC}$,
\[\Gamma\vdash_{\textbf{mbC}}\varphi\quad\text{if and only if}\quad\Gamma\vDash_{\textbf{mbC}}\varphi.\]
\end{theorem}

This is a classical result, which can be found in section $2.2$ of \cite{ParLog}, heavily inspired by Newton da Costa's valuation semantics for $C_{1}$ found in \cite{Costa} and \cite{Costa2}, that will serve us to prove that, if $T(\Gamma)\vdash_{\textbf{mbC}^{*}}T(\varphi)$, then $\Gamma\vdash_{\textbf{mbC}}\varphi$.

By contraposition, suppose $\Gamma\not\vdash_{\textbf{mbC}}\varphi$: then, there exists a bivaluation $\nu$ for $\textbf{mbC}$ such that $\nu(\Gamma)\subseteq\{1\}$ and $\nu(\varphi)=0$. We define a function $\overline{\nu}$ from the formulas of $\nbI$ to $\{0,1\}$ by structural induction:
\begin{enumerate}
\item if $\alpha$ is of order $0$, and is therefore a propositional variable, $\overline{\nu}(\alpha)=\nu(\alpha)$;
\item $\overline{\nu}(\alpha\vee\beta)=1$ if and only if $\overline{\nu}(\alpha)=1$ or $\overline{\nu}(\beta)=1$;
\item $\overline{\nu}(\alpha\wedge\beta)=1$ if and only if $\overline{\nu}(\alpha)=1$ and $\overline{\nu}(\beta)=1$;
\item $\overline{\nu}(\alpha\rightarrow\beta)=1$ if and only if $\overline{\nu}(\alpha)=0$ or $\overline{\nu}(\beta)=1$;
\item $\overline{\nu}(\neg  T(\alpha))=\nu(\neg \alpha)$, and if $\alpha$ is not a translation, $\overline{\nu}(\neg \alpha)=1$;
\item $\overline{\nu}(T(\alpha)\uparrow\neg T(\alpha))=\nu(\circ\alpha)$, $\overline{\nu}(\neg T(\alpha)\uparrow T(\alpha))=\nu(\circ\alpha)$, and, in all other cases, we make $\overline{\nu}(\alpha\uparrow\beta)=0$.
\end{enumerate}

\begin{proposition}
If $\alpha$ is a formula of $\textbf{mbC}$, $\overline{\nu}(T(\alpha))=\nu(\alpha)$.
\end{proposition}

\begin{proof}
We proceed by induction on the order of $\alpha$: if it is $0$, $\alpha$ is a propositional variable and $T(\alpha)=\alpha$, and by definition of $\overline{\nu}$ we have that $\overline{\nu}(T(\alpha))=\nu(\alpha)$.

Now assume the result holds for formulas of order smaller than that of $\alpha$, and we have five cases to consider:
\begin{enumerate}
\item if $\alpha=\phi\vee\psi$, $\overline{\nu}(T(\alpha))=1$ if and only if $\overline{\nu}(T(\phi))=1$ or $\overline{\nu}(T(\psi))=1$, which by hypothesis is equivalent to either $\nu(\phi)$ or $\nu(\psi)$ being equal to $1$, which occurs if and only if $\nu(\alpha)=1$;
\item if $\alpha=\phi\wedge\psi$, $\overline{\nu}(T(\alpha))=1$ is equivalent to $\overline{\nu}(T(\phi))$ and $\overline{\nu}(T(\psi))$ equating $1$, which in turn is equivalent by induction hypothesis to $\nu(\phi)$ and $\nu(\psi)$ equating $1$, what happens if and only if $\nu(\alpha)=1$;
\item if $\alpha=\phi\rightarrow\psi$, $\overline{\nu}(T(\alpha))=1$ if and only if $\overline{\nu}(T(\phi))=0$ or $\overline{\nu}(T(\psi))=1$, what by hypothesis is equivalent to $\nu(\phi)$ being $0$ or $\nu(\psi)$ being $1$, which happens if and only if $\nu(\alpha)=1$;
\item if $\alpha=\neg \phi$, by definition of $\overline{\nu}$ we have that 
\[\overline{\nu}(T(\alpha))=\overline{\nu}(\neg  T(\phi))=\nu(\neg \phi)=\nu(\alpha);\]
\item if $\alpha=\circ\phi$, we have that
\[\overline{\nu}(T(\alpha))=\overline{\nu}(T(\phi)\uparrow\neg T(\phi))=\nu(\circ\phi)=\nu(\alpha),\]
what ends the proof.
\end{enumerate}
\end{proof}

It is clear that, by its definition, $\overline{\nu}$ satisfies those conditions of being a bivaluation for $\nbI$, found at the beginning of Section \ref{Bivaluations for nbI}, related to disjunction, conjunction, and implication.

If $\overline{\nu}(\neg \alpha)=0$, the only possibility is that $\alpha$ is the translation of some formula $\phi$ in $\textbf{mbC}$ and that $\nu(\neg \phi)=0$; in this case $\nu(\phi)=1$, what implies 
\[\overline{\nu}(\alpha)=\overline{\nu}(T(\phi))=\nu(\phi)=1.\]
If $\overline{\nu}(\alpha\uparrow\beta)=1$, this means that:
\begin{enumerate}
\item either $\alpha$ is the translation of some formula $\phi$, $\beta=\neg\alpha$ and $\nu(\circ\phi)=1$, which means that either $\nu(\phi)=0$, and therefore $\overline{\nu}(\alpha)=0$, or $\nu(\neg \phi)=0$, and then $\overline{\nu}(\beta)=0$;
\item or $\beta$ is the translation of some formula $\phi$, $\alpha=\neg\beta$ and $\nu(\circ\phi)=1$, which again means that either $\overline{\nu}(\alpha)=0$ or $\overline{\nu}(\beta)=0$.
\end{enumerate}

Finally, if $\alpha$ is a translation of $\phi$ and $\beta=\neg\alpha$, $\overline{\nu}(\alpha\uparrow\beta)=\nu(\circ\phi)=\overline{\nu}(\beta\uparrow\alpha)$, the same happening if $\beta$ is a translation and $\alpha=\neg\beta$. Otherwise, we have 
\[\overline{\nu}(\alpha\uparrow\beta)=0=\overline{\nu}(\beta\uparrow\alpha).\]
This finishes proving that $\overline{\nu}$ is a bivaluation for $\nbI$; of course, $\overline{\nu}(T(\Gamma))=\nu(\Gamma)\subseteq\{1\}$ but $\overline{\nu}(T(\varphi))=\nu(\varphi)=0$, what implies $T(\Gamma)\not\vDash_{\nbI}T(\varphi)$; this is equivalent to $T(\Gamma)\not\vdash_{\nbI}T(\varphi)$, what means that $T(\Gamma)\not\vdash_{\textbf{mbC}^{*}}T(\varphi)$, since $\nbI$ is strictly stronger than $\textbf{mbC}^{*}$. Since $\Gamma\not\vdash_{\textbf{mbC}}\varphi$ implies $T(\Gamma)\not\vdash_{\textbf{mbC}^{*}}T(\varphi)$, $T(\Gamma)\vdash_{\textbf{mbC}^{*}}T(\varphi)$ implies $\Gamma\vdash_{\textbf{mbC}}\varphi$, and therefore 
\[\Gamma\vdash_{\textbf{mbC}}\varphi\quad\text{if and only if}\quad T(\Gamma)\vdash_{\textbf{mbC}^{*}}T(\varphi).\]


\subsection{$T$ is conservative for $\textbf{mbCciw}$, $\textbf{mbCci}$ and $\textbf{mbCcl}$}

$\textbf{mbCciw}^{*}$, $\textbf{mbCci}^{*}$ and $\textbf{mbCcl}^{*}$ are, respectively, the translations, under $T$, of $\textbf{mbCciw}$, $\textbf{mbCci}$ and $\textbf{mbCcl}$, meaning:
\begin{enumerate}
\item $\textbf{mbCciw}^{*}$ is obtained from $\textbf{mbC}^{*}$ by addition of the axiom schema
\[\tag{$\textbf{ciw}^{*}$}(\alpha\uparrow\neg\alpha)\vee(\alpha\wedge\neg \alpha);\]
\item $\textbf{mbCci}^{*}$ is obtained from $\textbf{mbC}^{*}$ by adding
\[\tag{$\textbf{ci}^{*}$}\neg (\alpha\uparrow\neg\alpha)\rightarrow(\alpha\wedge\neg \alpha);\]
\item $\textbf{mbCcl}^{*}$ is obtained from $\textbf{mbC}^{*}$ by adding
\[\tag{$\textbf{cl}^{*}$}\neg (\alpha\wedge\neg \alpha)\rightarrow(\alpha\uparrow\neg\alpha).\]
\end{enumerate}

We already know by Theorem \ref{Soundness of translations} that, for any of these logics $\mathcal{L}$, if $\Gamma\vdash_{\mathcal{L}}\varphi$ then $T(\Gamma)\vdash_{\mathcal{L}^{*}}T(\varphi)$. To prove the reciprocal, as was done with $\textbf{mbC}$, we will use bivaluations.

\begin{definition}\label{Bivaluations for mbCciw, mbCci, ...}
A bivaluation for $\textbf{mbCciw}$, $\textbf{mbCci}$ or \index{Bivaluation for $\textbf{mbCci}$}$\textbf{mbCcl}$ (already defined, for $\textbf{mbCciw}$ and $\textbf{mbCcl}$, in Section \ref{RNmatrix for mbCcl})  is a valuation for $\textbf{mbC}$ satisfying also that, respectively:
\begin{enumerate}
\item $\nu(\circ\alpha)=1$ if and only if $\nu(\alpha)=0$ or $\nu(\neg \alpha)=0$;
\item the above condition, plus that $\nu(\neg \circ\alpha)=1$ implies $\nu(\alpha)=1$ and $\nu(\neg \alpha)=1$;
\item the first condition, plus that $\nu(\neg (\alpha\wedge\neg \alpha))=1$ implies $\nu(\circ\alpha)=1$.
\end{enumerate}

Given formulas $\Gamma\cup\{\varphi\}$ in the signature $\Sigma_{\textbf{LFI}}$, we say $\Gamma$ proves $\varphi$ according to bivaluations for $\textbf{mbCciw}$, and write $\Gamma\vDash_{\textbf{mbCciw}}\varphi$, if for all valuations $\nu$ for $\textbf{mbCciw}$ such that $\nu(\Gamma)\subseteq\{1\}$ we have that $\nu(\varphi)=1$; we define $\vDash_{\textbf{mbCci}}$ and $\vDash_{\textbf{mbCcl}}$ analogously.
\end{definition}

The following result is well established and can be found in sections $3.1$ and $3.3$ of \cite{ParLog}.

\begin{theorem}
Given formulas $\Gamma\cup\{\varphi\}$ over the signature $\Sigma$,
\begin{enumerate}
\item $\Gamma\vdash_{\textbf{mbCciw}}\varphi$ if and only if $\Gamma\vDash_{\textbf{mbCciw}}\varphi$;
\item $\Gamma\vdash_{\textbf{mbCci}}\varphi$ if and only if $\Gamma\vDash_{\textbf{mbCci}}\varphi$;
\item $\Gamma\vdash_{\textbf{mbCcl}}\varphi$ if and only if $\Gamma\vDash_{\textbf{mbCcl}}\varphi$.
\end{enumerate}
\end{theorem}

Now we have the tools necessary to prove that $\Gamma\vdash_{\mathcal{L}}\varphi$ if and only if $T(\Gamma)\vdash_{\mathcal{L}^{*}}T(\varphi)$, for any $\mathcal{L}\in\{\textbf{mbCciw}, \textbf{mbCci}, \textbf{mbCcl}\}$. We will proceed by contraposition, therefore suppose $\Gamma\not\vdash_{\mathcal{L}}\varphi$, and there exists a bivaluation $\nu$ for $\mathcal{L}$ such that $\nu(\Gamma)\subseteq\{1\}$ but $\nu(\varphi)=0$. We define a function $\overline{\nu}$ from the formulas over the signature $\Sigma_{\nbI}$ to $\{0,1\}$ by structural induction.

\begin{enumerate}
\item If $\alpha$ is of order $0$, and therefore a propositional variable, $\overline{\nu}(\alpha)=\nu(\alpha)$;
\item $\overline{\nu}(\alpha\vee\beta)=1$ if and only if $\overline{\nu}(\alpha)=1$ or $\overline{\nu}(\beta)=1$;
\item $\overline{\nu}(\alpha\wedge\beta)=1$ if and only if $\overline{\nu}(\alpha)=1$ and $\overline{\nu}(\beta)=1$;
\item $\overline{\nu}(\alpha\rightarrow\beta)=1$ if and only if $\overline{\nu}(\alpha)=0$ or $\overline{\nu}(\beta)=1$;
\item $\overline{\nu}(\neg T(\alpha))=\nu(\neg \alpha)$, and if $\alpha$ is not a translation:
\begin{enumerate}
\item for $\mathcal{L}=\textbf{mbCci}$ and $\alpha=\beta\uparrow\neg\beta$, if $\overline{\nu}(\beta)=0$ or $\overline{\nu}(\neg \beta)=0$, $\overline{\nu}(\neg \alpha)=0$;
\item for $\mathcal{L}=\textbf{mbCcl}$ and $\alpha=\beta\wedge\neg \beta$, if $\overline{\nu}(\beta\uparrow\neg\beta)=0$, $\overline{\nu}(\neg \alpha)=0$;
\end{enumerate}
if $\alpha$ is not a translation and we find ourselves in neither of the previous two described cases, $\overline{\nu}(\neg \alpha)=1$;
\item $\overline{\nu}(T(\alpha)\uparrow\neg T(\alpha))=\nu(\circ\alpha)$ and $\overline{\nu}(\neg T(\alpha)\uparrow T(\alpha))=\nu(\circ\alpha)$; if $\alpha$ is not a translation, $\overline{\nu}(\alpha\uparrow\neg\alpha)$ and $\overline{\nu}(\neg\alpha\uparrow\alpha)$ are $1$ if, and only if, $\overline{\nu}(\alpha)=0$ or $\overline{\nu}(\neg \alpha)=0$; and if $\alpha$ is not a translation or $\beta$ is not the negation of $\alpha$, $\overline{\nu}(\alpha\uparrow\beta)=0$.
\end{enumerate}

\begin{proposition}
If $\alpha$ is a formula over the signature $\Sigma_{\textbf{LFI}}$, $\overline{\nu}(T(\alpha))=\nu(\alpha)$.
\end{proposition}

It is clear that $\overline{\nu}$ satisfies the conditions for being a bivaluation for $\nbI$ regarding disjunction, conjunction and implication. The image under $\overline{\nu}$ of the negation of a formula is $0$ under the following conditions:
\begin{enumerate}
\item if the formula is the translation of an $\alpha$ and $\nu(\alpha)=0$, and since $\nu$ is a bivaluation for $\textbf{mbC}$, we find $\overline{\nu}(T(\alpha))=\nu(\alpha)=1$;
\item if the formula is not a translation and of the form $\beta\uparrow\neg\beta$, $\mathcal{L}=\textbf{mbCci}$ and $\overline{\nu}(\beta)=0$ or $\overline{\nu}(\neg\beta)=0$, and then we obtain that $\overline{\nu}(\beta\uparrow\neg\beta)=1$;
\item  if the formula is not a translation and of the form $\beta\wedge\neg\beta$, $\mathcal{L}=\textbf{mbCcl}$ and $\overline{\nu}(\beta\uparrow\neg \beta)=0$, in which case $\overline{\nu}(\beta)=\overline{\nu}(\neg\beta)=1$ and therefore $\overline{\nu}(\beta\wedge\neg\beta)=1$.
\end{enumerate}

To summarize, in all cases such that $\overline{\nu}(\neg\alpha)=0$, we have $\overline{\nu}(\alpha)=1$, as it would be expected of a bivaluation for $\nbI$ regarding the paraconsistent negation. Finally, the image under $\overline{\nu}$ of a formula with the incompatibility operator as leading connective is $1$ under the following conditions:
\begin{enumerate}
\item if the formula is $T(\alpha)\uparrow\neg T(\alpha)$ or $\neg T(\alpha)\uparrow T(\alpha)$ and $\nu(\circ\alpha)=1$, in which case $\nu(\alpha)=0$ or $\nu(\neg\alpha)=0$, implying either $\overline{\nu}(T(\alpha))=\nu(\alpha)=0$ or $\overline{\nu}(\neg T(\alpha))=\nu(\neg\alpha)=0$;
\item if the formula is not a translation and of the form $\alpha\uparrow\neg\alpha$ or $\neg\alpha\uparrow\alpha$, meaning that either $\overline{\nu}(\alpha)=0$ or $\overline{\nu}(\neg\alpha)=0$.
\end{enumerate}

Very clearly, for any two formulas over the signature $\Sigma_{\nbI}$ we have that $\overline{\nu}(\alpha\uparrow\beta)=\overline{\nu}(\beta\uparrow\alpha)$, allowing us to conclude that $\overline{\nu}$ is a bivaluation for $\nbI$, remaining for us to show that for each logic, $\textbf{mbCciw}$, $\textbf{mbCci}$ or $\textbf{mbCcl}$, $\overline{\nu}$ is a bivaluation for, respectively, $\nbIciw$, $\nbIci$ and $\nbIcl$ (refer back to Definition \ref{definition of Bivaluation for nbIciw, nbIci, ...} if necessary).

\begin{enumerate}
\item For any of the three logics $\textbf{mbCciw}$, $\textbf{mbCci}$ or $\textbf{mbCcl}$, it is clear that $\overline{\nu}(\alpha\uparrow\neg\alpha)=1$ if and only if $\overline{\nu}(\alpha)=0$ or $\overline{\nu}(\neg \alpha)=0$, by definition of $\overline{\nu}$.
\item If $\mathcal{L}=\textbf{mbCci}$, and $\overline{\nu}(\neg(\alpha\uparrow\neg\alpha))=1$, this means both $\overline{\nu}(\alpha)$ and $\overline{\nu}(\neg \alpha)$ are $1$ (and therefore $\overline{\nu}(\alpha\wedge\neg\alpha)=1$), since otherwise, by definition of $\overline{\nu}$, we would have $\overline{\nu}(\neg(\alpha\uparrow\neg\alpha))=0$.
\item Finally, when $\mathcal{L}=\textbf{mbCcl}$, if $\overline{\nu}(\neg(\alpha\wedge\neg \alpha))=1$, we have that $\overline{\nu}(\alpha\uparrow\neg \alpha)=1$, since otherwise, by definition of $\overline{\nu}$, we would have $\overline{\nu}(\neg (\alpha\wedge\neg \alpha))=0$.
\end{enumerate}
 
So $\overline{\nu}$ is a bivaluation for $\nbIciw$, $\nbIci$ or $\nbIcl$, whenever $\mathcal{L}$ is $\textbf{mbCciw}$, $\textbf{mbCci}$ or $\textbf{mbCcl}$, such that $\overline{\nu}(T(\Gamma))=\nu(\Gamma)\subseteq\{1\}$ but $\overline{\nu}(T(\varphi))=\nu(\varphi)=0$, what implies 
\[T(\Gamma)\not\vDash_{\nbIciw}T(\varphi), \quad T(\Gamma)\not\vDash_{\nbIci}T(\varphi)\quad \text{and}\quad T(\Gamma)\not\vDash_{\nbIcl}T(\varphi),\]
or what is equivalent,
\[T(\Gamma)\not\vdash_{\nbIciw}T(\varphi), \quad T(\Gamma)\not\vdash_{\nbIci}T(\varphi)\quad \text{and}\quad T(\Gamma)\not\vdash_{\nbIcl}T(\varphi);\]
since $\nbIciw$, $\nbIci$ and $\nbIcl$ are strictly stronger than, respectively, $\textbf{mbCciw}^{*}$, $\textbf{mbCci}^{*}$ and $\textbf{mbCcl}^{*}$, this implies that $T(\Gamma)\not\vdash_{\mathcal{L}^{*}}T(\varphi)$, and therefore we have that $T(\Gamma)\vdash_{\mathcal{L}^{*}}T(\varphi)$ implies $\Gamma\vdash_{\mathcal{L}}\varphi$, or what is equivalent,
\[\Gamma\vdash_{\mathcal{L}}\varphi\quad\text{if and only if}\quad T(\Gamma)\vdash_{\mathcal{L}^{*}}T(\varphi),\]
for any $\mathcal{L}\in\{\textbf{mbCciw}, \textbf{mbCci}, \textbf{mbCcl}\}$.

\newpage
\printbibliography[segment=\therefsegment,heading=subbibliography]
\end{refsegment}

\begin{refsegment}
\defbibfilter{notother}{not segment=\therefsegment}

\renewcommand{\chaptermark}[1]{\markboth{#1}{#1}}

\chapter*{Conclusion}
\chaptermark{Conclusion}
\addcontentsline{toc}{chapter}{Conclusion} 

We believe the main developments found in this work are, in its more algebraic side, weakly free multialgebras, and, in the logical one, RNmatrices and logics of incompatibility. There were other advances, of course, however if there are topics that may still offer further fruitful research, these are probably the ones.

We now summarize what we have hoped the ideas here presented can achieve, and future works planned based on them; for clarity, we divide our text between the two broad parts in which this thesis has been divided so far.

\subsection*{Multialgebras}

Our main motivation, when studying multialgebras, remains to search for a general theory of their identities, much like universal algebra is to algebras themselves. Of course, this is no easy task, but we believe weekly free multialgebras are a step in the right direction. In Chapter \ref{Chapter2} we have shown that these structures enjoy many characterizations similar to those of absolutely free algebras, and probably are the correct generalizations of objects without identities to the context of multialgebras. We have used this to offer simplified proofs of some standard, folkloric results in the theory of non-deterministic algebras, but we still are left wondering what should identities on a multialgebra be. There are many candidates, and the task ahead is to identify which serve what purposes; but the next logical step, once we have generalized absolutely free algebras, is to generalize the relatively free ones, a closely related endeavor since ideals of identities in universal algebra fully characterize relatively free objects.

Showing, in Chapter \ref{Chapter3}, the category of multialgebras is equivalent to a category of ordered, almost Boolean algebras seems to have mainly philosophical implications, in the following sense: non-deterministic semantics are seem by many, in what we consider an unfortunately regressive point of view, as an undesirable development in the problem of decision methods for logics, precisely because of the non-determinism. Here, we offer an alternative, not only for Nmatrices but also PNmatrices and closely related semantics: they may be recast as deterministic structures, as long as we consider, along their operations, also their orders. This invites one to recast already existing decision methods that rely on non-determinism as deterministic decision methods; alternatively, one can also analyze what the equivalences presented here, and others they may lead to being discovered, can say about the involved categories.

\subsection*{Paraconsistent Logic}

Chapter \ref{Chapter4} introduces the semantics of restricted non-deterministic matrices. We start by showing how previous semantics found in the literature, created in order to solve specific problems and without an underlying theory, can be reinterpreted as RNmatrices, what already establishes that this new tool indeed is useful. To gather further evidence for that point, we show many logics of difficult treatment, as those paraconsistent systems between $\mbCcl$ and $\CILA$, da Costa's hierarchy as a whole, and logics of formal incompatibility, can be relatively easily characterized by RNmatrices, with accompanying decision methods. So we hope to have convinced the reader of restricted Nmatrices' utility, but the issue of determining a theory for these semantics is still very open: can a theory like the one Blok and Pigozzi designed for the process of algebraization of logics exist for non-deterministic algebraization, specifically in our case through RNmatrices? A less ambitious objective, however, is searching for other systems of demanding characterization that can be more easily studied with RNmatrices: in addition to paraconsistent ones, modal logics seem like good candidates.

Still on the subject of RNmatrices, Chapter \ref{Chapter5} investigates how to obtain, from them, decision methods, in this case truth tables and tableaux, for da Costa's systems. We ponder if these procedures can be systematized, producing for all restricted Nmatrices, in some up to now unspecified class, row-branching, row-eliminating truth tables and tableau calculi; even more, we are interested in knowing if other varieties of decision methods, such as sequents, may be extracted from these semantics. Of course, regarding the decision methods and restricted non-deterministic matrices for da Costa's $C_{n}$ logics themselves, here suggested to be $n+2$-valued in nature, one should be curious of how they place against  existing decision methods for these systems, such as bivaluations, different tableau calculi and effective Nmatrices: are they more efficient? Less efficient? We are yet to survey these questions, although we are certain of the importance of having produced new instruments with which to approach such challenging logics.

Finally, even before it is ascertained whether a theory of non-deterministic algebraization is possible, Chapter \ref{Chapter6} starts to look into other directions in which to extend RNmatrices' analysis: by looking at swap structures, a general methodology that can be very well the most successful for finding restricted Nmatrices, we are lead to the notion of a category of RNmatrices. In the case of the systems $C_{n}$, this category is shown to be astoundingly well-behaved; in this thesis a broad procedure for creating categories of restricted non-deterministic matrices is already outlined, and we are interested in the problem of whether this construction always returns nice categories. Restricted swap structures for da Costa's hierarchy also suggest strong model theoretic connections, and the combinatorial investigation of their snapshots is a first step in this direction: a second step would be to scrutinize those results in model theory, if any, hold for arbitrary restricted swap structures.

Chapter \ref{Chapter7} broaches another of the more relevant contributions of this work, that of logics of incompatibility: yes, there have been previous attempts of formalizing the natural notion of incompatible properties, more prominently in Brandom's neo-pragmatic program for epistemology, but much as da Costa's systematization of contradiction brought new light into the subject of paraconsistency, we hope that a more controlled take on incompatibility could better develop the field. After defining some simpler systems, accounting for commutativity and propagation of incompatibility, we supply them with semantics, both based on bivaluations and RNamtrices, as well as decision methods constructed from the aforementioned RNmatrices; this gives additional evidence of the importance of restricted Nmatrices. In future developments, one hopes to define other systems, possibly modeling additional attributes of incompatible properties in natural language or including quantifiers and modalities.

One possible angle for studying logics of incompatibility is adding to them a negation, with which the incompatibility can interact, task we undertake in some depth in Chapter \ref{Chapter8}. Most of the systems we define are based in others for paraconsistency, and are dealt with by using bivaluations and RNmatrices, obtaining thus decision methods again through row-branching, row-eliminating truth tables and tableaux. And, what is formally proven in Chapter \ref{Chapter9} but heavily insinuated since their definitions, these logics extend those paraconsistent ones upon which they are built in a non-trivial way: that is, we have conservative translations from the latter into sublogics of the former; this means our logics of incompatibility are strictly, and non trivially, generalizing logics of formal inconsistency, and invites the interested to find the most appropriate translations for inconsistent systems, a line of future research. Furthermore, by providing characterizing RNmatrices for logics of incompatibility with negation, we feel compelled to prove that more classical semantics such as Blok and Pigozzi's algebraization, Nmatrices, and restricted matrices would not work. This means we have offered logics very simple to define, yet very difficult to characterize, and leads us to wonder about the interplay between competing generalizations of logical matrices: restricted matrices, to give one example, seem very promising although virtually unstudied.

\end{refsegment}

\newpage
\thispagestyle{plain}
\hspace{0pt}
\vfill
\begin{figure}[h]
\centering
\includegraphics[width=\textwidth]{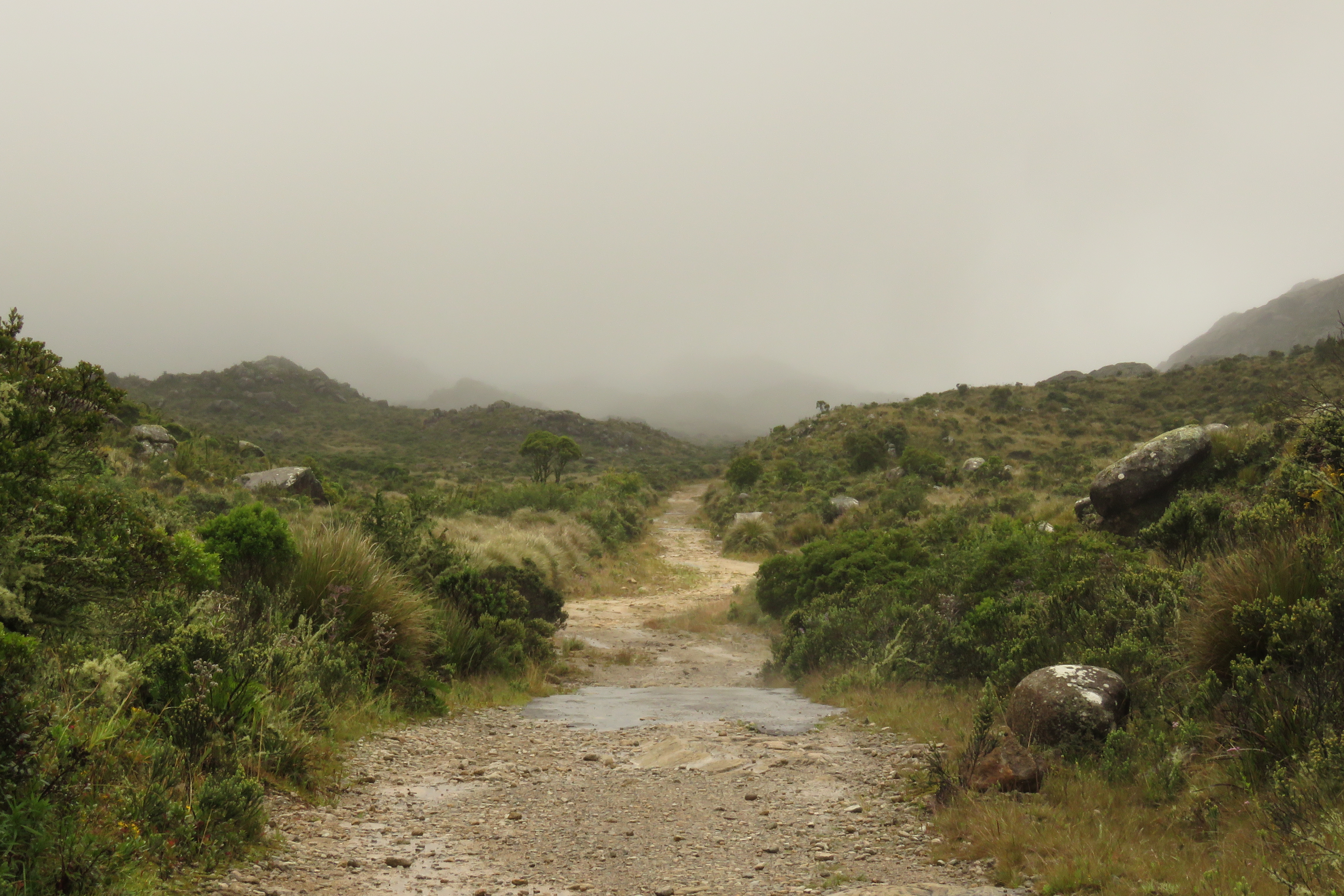}
\caption*{Mountainous region of Parque Nacional do Itatiaia, Itatiaia, Brazil.\\Photographed by Guilherme Vicentin de Toledo, all rights reserved.}
\end{figure}
\vfill
\hspace{0pt}

\cleardoublepage

\phantomsection

\begin{refsegment}
\defbibfilter{notother}{not segment=\therefsegment}

\chapter*{Lists}
\markboth{Lists}{}
\addcontentsline{toc}{chapter}{Lists}

Listed by order of appearance.

\section*{List of symbols}

\begin{longtable}{c | m{11cm} | c}

Symbols & Meaning & Page \\ \hline

$\mathcal{P}(A)$ & Powerset of the set $A$ & \pageref{powerset}\\

$\varphi:\mathcal{A}\rightarrow\mathcal{B}$ & Homomorphism from $\mathcal{A}$ to $\mathcal{B}$ & \pageref{homomorphism}\\

$F(\Sigma, \mathcal{V})$ & Formulas, on the signature $\Sigma$, over the variables $\mathcal{V}$ & \pageref{FSigmaV} \\

$\textbf{F}(\Sigma, \mathcal{V})$ & $\Sigma$-Algebra of formulas, on the signature $\Sigma$, over the variables $\mathcal{V}$ & \pageref{tFSigmaV} \\

$|\alpha|$ & Order of the formula $\alpha$ & \pageref{order}\\

$\Sigma^{\kappa}$ & Expanded signature & \pageref{expandedsignature}\\

$\textbf{mF}(\Sigma, \mathcal{V}, \kappa)$ & $\Sigma$-Multialgebra of non-deterministic formulas, on the signature $\Sigma$, over the variables $\mathcal{V}$ and bounded by $\kappa$ & \pageref{mFSigmakappaV}\\

$|X|$ & Cardinality of the set $X$ & \pageref{|X|}\\

$\textbf{cdf}$-generated & Choice-dependent freely generated & \pageref{cdfgenerated}\\

$B(\mathcal{A})$ & Build of the multialgebra $\mathcal{A}$ & \pageref{build}\\

$G(\mathcal{A})$ & Ground of the multialgebra $\mathcal{A}$ & \pageref{ground}\\

$\langle S\rangle_{m}$ & The $m$-th set obtained in the process of defining $\langle S\rangle$ & \pageref{mth generated}\\

$\langle S\rangle$ & Set generated by $S$ & \pageref{generated}\\

$M(\mathcal{A})$ & Supremum of the cardinalities of the operations on the multialgebra $\mathcal{A}$ & \pageref{M(A)}\\

$o_{B}(a)$ & $B$-order of the element $a$ of a multialgebra with strong basis $B$ & \pageref{Border}\\

$\mathcal{A}_{0}$ & Poset obtained from the poset $\mathcal{A}$ by addition of a minimum element & \pageref{A0}\\

$\mathcal{L}$ & Propositional language generated by $\Sigma$ & \pageref{language}\\

$\Gamma\vdash\varphi$ & $\Gamma$ syntactically proves $\varphi$ & \pageref{vdash}\\

$\mathscr{L}$ & An arbitrary logic & \pageref{logic}\\

$\alpha_{1}, \dotsc , \alpha_{n}|\alpha$ & Rule of inference where $\{\alpha_{1},\dotsc , \alpha_{n}\}$ deduces $\alpha$ & \pageref{ruleofinference}\\

$\mathfrak{A}$ & A set of axiom schemata & \pageref{frakA}\\

$\mathfrak{R}$ & A set of rules of inference & \pageref{frakR}\\

$\psi(p_{1}, \dotsc , p_{n})$ & A formula $\psi$ whose variables are among $\{p_{1}, \dotsc , p_{n}\}$ & \pageref{varformula}\\

$\circ\alpha$ & Formula stating the consistency of $\alpha$ & \pageref{circ}\\

$\Gamma\vDash_{\mathcal{M}}\varphi$ & $\Gamma$ semantically proves $\varphi$ according to $\mathcal{M}$ & \pageref{semant.proves}\\

$\Gamma\vDash_{\mathbb{M}}\varphi$ & $\Gamma$ semantically proves $\varphi$ according to the class $\mathbb{M}$ & \pageref{vDashM}\\

$\textbf{2}(\Sigma)$ & $\Sigma$-multialgebra with universe $\{0,1\}$ whose operations always return $\{0,1\}$ & \pageref{2Sigma}\\

$\textbf{2}(\mathfrak{L})$ & RNmatrix, with multialgebra $\textbf{2}(\Sigma)$, characterizing $\mathfrak{L}$ & \pageref{2L}\\

$K_{\mathcal{M}}(\Gamma)$ & Logical closure, according to $\mathcal{M}$, of $\Gamma$ & \pageref{KMGamma}\\

$K_{\mathbb{M}}(\Gamma)$ & Logical closure, according to the class $\mathcal{M}$, of $\Gamma$ & \pageref{KMMGamma}\\

$\mathcal{K}_{\textbf{L}}$ & Kearns RNmatrix for the logic $\textbf{L}$ & \pageref{KL}\\

$N_{a}$ & Possible negations of $a$ in a Fidel structure & \pageref{Na}\\

$\mathcal{M}^{s}$ & Static semantics for the Nmatrix $\mathcal{M}$ & \pageref{Ms}\\

$\mathcal{M}^{\emptyset}$ & RNmatrix equivalent to the PNmatrix $\mathcal{M}$ & \pageref{Mo}\\

$\mathcal{M}_{\textbf{mbCciw}}$ & Nmatrix characterizing $\textbf{mbCciw}$ & \pageref{MmbCciw}\\

$\mathcal{M}_{\textbf{mbCcl}}$ & RNmatrix characterizing $\textbf{mbCcl}$ & \pageref{MmbCcl}\\

$\mathcal{M}_{\textbf{CILA}}$ & RNmatrix characterizing $\textbf{CILA}$ & \pageref{MCILA}\\

$\alpha^{n}$ & $\alpha^{0}=\alpha$ and $\alpha^{n+1}=\neg(\alpha^{n}\wedge\neg\alpha^{n})$ & \pageref{alpha^}\\

$\alpha^{(n)}$ & $\alpha^{(1)}=\alpha^{1}$ and $\alpha^{(n+1)}=\alpha^{n+1}\wedge\alpha^{(n)}$ & \pageref{alpha^()}\\

$\alpha^{\circ}$ & Equal to $\alpha^{1}$ & \pageref{alphacirc}\\

$T_{n}$ & The snapshot $(1, 0, 1, \dotsc , 1)$ & \pageref{Tn}\\

$F_{n}$ & The snapshot $(0, 1, 1, \dotsc , 1)$ & \pageref{Fn}\\

$B_{n}$ & Set of snapshots for $C_{n}$ & \pageref{Bn}\\

$D_{n}$ & Set of designated snapshots for $C_{n}$ & \pageref{Dn}\\

$Boo_{n}$ & Set of Boolean snapshots for $C_{n}$ & \pageref{Boon}\\

$\mathcal{A}_{C_{n}}$ & Swap structure for $C_{n}$ & \pageref{ACn}\\

$\mathcal{F}_{C_{n}}$ & Set of restricted valuations for $C_{n}$ & \pageref{FCn}\\

$\mathcal{RM}_{C_{n}}$ & RNmatrix characterizing $C_{n}$ & \pageref{RMCn}\\

$\mathbb{T}_{n}$ & Tableau calculus for $C_{n}$ & \pageref{mathbbTn}\\

$\Gamma\vdash_{\mathbb{T}_{n}}\varphi$ & $\Gamma$ proves $\varphi$ according to the tableau calculus $\mathbb{T}_{n}$ & \pageref{vdashmathbbTn}\\

$\Gamma\vDash_{C_{n}}^{\mathcal{B}}\varphi$ & $\Gamma$ semantically proves $\varphi$, according to $\mathcal{B}$-valuations & \pageref{vDashBCn}\\

$B_{n}^{\mathcal{B}}$ & Set of snapshots, over $\mathcal{B}$, for $C_{n}$ & \pageref{BnB}\\

$D_{n}^{\mathcal{B}}$ & Set of designated snapshots, over $\mathcal{B}$, for $C_{n}$ & \pageref{DnB}\\

$Boo_{n}^{\mathcal{B}}$ & Set of Boolean snapshots, over $\mathcal{B}$, for $C_{n}$ & \pageref{BoonB}\\

$\mathcal{A}_{C_{n}}^{\mathcal{B}}$ & Swap structure, over $\mathcal{B}$, for $C_{n}$ & \pageref{ACnB}\\

$\mathcal{RM}_{C_{n}}^{\mathcal{B}}$ & RNmatrix, over $\mathcal{B}$, characterizing $C_{n}$ & \pageref{RMCnB}\\

$\mathcal{F}_{C_{n}}^{\mathcal{B}}$ & Set of restricted valuations, over $\mathcal{B}$, for $C_{n}$ & \pageref{FCnB}\\

$\Gamma\vDash_{\mathcal{RM}_{C_{n}}^{\mathcal{B}}}^{RN}\varphi$ & $\Gamma$ semantically proves $\varphi$, according to $\mathcal{RM}_{C_{n}}^{\mathcal{B}}$ & \pageref{vDashRNRM}\\

$\Gamma\vDash_{\mathcal{RS}_{C_{n}}}^{RN}\varphi$ & For every Boolean algebra $\mathcal{B}$, $\Gamma\vDash_{\mathcal{RM}_{C_{n}}^{\mathcal{B}}}^{RN}\varphi$ & \pageref{vDashRNRS}\\

$|\mathcal{B}|$ & Universe of the Boolean algebra $\mathcal{B}$ & \pageref{|B|}\\

$X_{m}$ & Generic set with $m$ elements & \pageref{Xm}\\

$\pi_{i}$ & Projection from a product to its $i$-th coordinate & \pageref{pii}\\

$\alpha\uparrow\beta$ & $\alpha$ and $\beta$ are incompatible & \pageref{uparrow}\\

$\alpha\downarrow\beta$ & $\alpha$ and $\beta$ are compatible & \pageref{downarrow}\\

$\Gamma\vDash_{\bI}\varphi$ & $\Gamma$ semantically proves $\varphi$, according to bivaluations for $\textbf{bI}$ & \pageref{vDashbI}\\

$\Gamma\Vdash_{\mathcal{F}}^{\bI}\varphi$ & $\Gamma$ semantically proves $\varphi$, according to Fidel structures for $\textbf{bI}$ & \pageref{VdashFbI}\\

$\alpha\equiv_{\Gamma}^{\bI}\beta$ & $\alpha$ and $\beta$ are equivalent, according to $\Gamma$, in $\bI$ & \pageref{equivbI}\\

$A^{\bI}_{\Gamma}$ & The quotient set of $F(\Sigma_{\bI}, \mathcal{V})$ by $\equiv_{\Gamma}^{\bI}$ & \pageref{AbIGamma}\\

$\mathcal{A}^{\bI}_{\Gamma}$ & The Lindenbaum-Tarski multialgebra of $\bI$ associated to $\Gamma$ & \pageref{mAbIGamma}\\

$\textbf{2}_{\bI}$ & Two-valued Fidel structure for $\bI$ & \pageref{2bI}\\

$\Gamma\vDash_{\bIpr}\varphi$ & $\Gamma$ semantically proves $\varphi$, according to bivaluations for $\textbf{bIpr}$ & \pageref{vDashbIpr}\\

$\Gamma\Vdash_{\mathcal{F}}^{\bIpr}\varphi$ & $\Gamma$ semantically proves $\varphi$, according to Fidel structures for $\textbf{bIpr}$ & \pageref{VdashFbIpr}\\

$\textbf{2}_{\bIpr}$ & Two-valued Fidel structure for $\bIpr$ & \pageref{2bIpr}\\

$\Gamma\vDash_{\nbI}\varphi$ & $\Gamma$ proves $\varphi$, according to bivaluations for $\nbI$ & \pageref{vDashnbI}\\

$\Gamma\Vdash_{\mathcal{F}}^{\nbI}\varphi$ & $\Gamma$ proves $\varphi$, according to Fidel structures for $\nbI$ & \pageref{VdashFnbI}\\

$\alpha\equiv_{\Gamma}^{\nbI}\beta$ & $\alpha$ and $\beta$ are equivalent, according to $\Gamma$, in $\nbI$ & \pageref{equivnbI}\\

$\mathcal{A}^{\nbI}_{\Gamma}$ & The Lindebaum-Tarski multialgebra of $\nbI$ associated to $\Gamma$ & \pageref{AnbIGamma}\\

$\textbf{2}_{\nbI}$ & Two-valued Fidel structure for $\nbI$ & \pageref{2nbI}\\

$\Gamma\vDash_{\nbIciw}\varphi$ & $\Gamma$ proves $\varphi$, according to bivaluations for $\nbIciw$ & \pageref{vDashnbIciw}\\

$\Gamma\vDash_{\nbIci}\varphi$ & $\Gamma$ proves $\varphi$, according to bivaluations for $\nbIci$ & \pageref{vDashnbIci}\\

$\Gamma\vDash_{\nbIcl}\varphi$ & $\Gamma$ proves $\varphi$, according to bivaluations for $\nbIcl$ & \pageref{vDashnbIcl}\\

$\Gamma\Vdash_{\mathcal{F}}^{\nbIciw}\varphi$ & $\Gamma$ proves $\varphi$, according to Fidel structures for $\nbIciw$ & \pageref{VdashFnbIciw}\\

$\Gamma\Vdash_{\mathcal{F}}^{\nbIci}\varphi$ & $\Gamma$ proves $\varphi$, according to Fidel structures for $\nbIci$ & \pageref{VdashFnbIci}\\

$\Gamma\Vdash_{\mathcal{F}}^{\nbIcl}\varphi$ & $\Gamma$ proves $\varphi$, according to Fidel structures for $\nbIcl$ & \pageref{VdashFnbIcl}\\

$\nabla$ & Trivial congruence on a multialgebra $\mathcal{A}$, equal to $A\times A$ & \pageref{nabla}\\

$\Delta$ & Trivial congruence on a multialgebra $\mathcal{A}$, equal to $\{(a,a) : a\in A\}$ & \pageref{Delta}\\

$T$ & Translation from the $\textbf{LFI}$'s to logics of incompatibility & \pageref{T}\\

\end{longtable}

\newpage

\section*{List of logics and axioms}

\begin{longtable}{c | m{11cm} | c}

Symbols & Meaning & Page \\ \hline

$\textbf{CPL}$ & Classical propositional logic & \pageref{CPL}\\

$\textbf{LFI}$ & Logics of formal inconsistency & \pageref{LFI}\\

$\textbf{mbC}$ & One of the simplest $\textbf{LFI}$'s & \pageref{mbC}\\

$\textbf{mbCciw}$ & Logic obtained from $\textbf{mbC}$ by addition of $\textbf{ciw}$ & \pageref{mbCciw}\\

$\textbf{ciw}$ & Axiom scheme $\circ\alpha\vee(\alpha\wedge\neg\alpha)$ & \pageref{ciw}\\

$\textbf{mbCci}$ & Logic obtained from $\textbf{mbC}$ by addition of $\textbf{ci}$ & \pageref{mbCci}\\

$\textbf{ci}$ & Axiom scheme $\neg \circ\alpha\rightarrow(\alpha\wedge\neg \alpha)$ & \pageref{ci}\\

$\textbf{mbCcl}$ & Logic obtained from $\textbf{mbC}$ by addition of $\textbf{cl}$ & \pageref{mbCcl}\\

$\textbf{cl}$ & Axiom scheme $\neg(\alpha\wedge\neg \alpha)\rightarrow\circ\alpha$ & \pageref{cl}\\

$\textbf{B}[\{\textbf{i1}, \textbf{i2}\}]$ & System, defined by Avron, equivalent to $\textbf{mbCci}$ & \pageref{Bi1i2}\\

$\textbf{Bi}$ & System, defined by Avron, equivalent to $\textbf{mbCci}$ & \pageref{Bi}\\

$\textbf{Bl}$ & System, defined by Avron, equivalent to $\textbf{mbCcl}$ & \pageref{Bl}\\

$\textbf{CILA}$ & An $\textbf{LFI}$ equivalent to da Costa's $C_{1}$ & \pageref{CILA}\\

$\textbf{cf}$ & Axiom scheme $\neg\neg\alpha\rightarrow\alpha$ & \pageref{cf}\\

$\textbf{Ci}$ & Logic obtained from $\textbf{mbCci}$ by addition of $\textbf{cf}$ & \pageref{Ci}\\

$C_{\omega}$ & An intuitionistic and paraconsistent system by da Costa & \pageref{Comega}\\

$C_{n}$ & The $n$-th logic in da Costa's hierarchy & \pageref{Cn}\\

$\textbf{bc}_{n}$ & Axiom scheme $\alpha^{(n)}\rightarrow(\alpha\rightarrow(\neg\alpha\rightarrow\beta))$ & \pageref{bcn}\\

$\textbf{p}\#_{n}$ & Axiom scheme $(\alpha^{(n)}\#\beta^{(n)})\rightarrow(\alpha\vee\beta)^{(n)}$ & \pageref{phashn}\\

$\textbf{LIp}$ & Logic of formal incompatibility & \pageref{LIp}\\

$\textbf{bI}$ & A basic logic of incompatibility, with commutativity & \pageref{bI}\\

$\textbf{Ip}$ & Axiom scheme $(\alpha\uparrow\beta)\rightarrow(\alpha\rightarrow(\beta\rightarrow \gamma))$ & \pageref{Ip}\\

$\textbf{Comm}$ & Axiom scheme $(\alpha\uparrow\beta)\rightarrow(\beta\uparrow\alpha)$ & \pageref{Comm}\\

$\textbf{bI}^{-}$ & The logic $\textbf{bI}$ without commutativity & \pageref{bI-}\\

$\textbf{bIpr}$ & Logic obtained from $\bI$ by addition of $\textbf{pr}_{\vee}$ and $\textbf{pr}_{\wedge}$ & \pageref{bIpr}\\

$\textbf{pr}_{\wedge}$ & Axiom scheme $[(\alpha\uparrow\gamma)\wedge(\beta\uparrow\gamma)]\rightarrow[(\alpha\vee\beta)\uparrow\gamma]$ & \pageref{prwedge}\\

$\textbf{pr}_{\vee}$ & Axiom scheme $[(\alpha\uparrow\gamma)\wedge(\beta\uparrow\gamma)]\rightarrow[(\alpha\vee\beta)\uparrow\gamma]$ & \pageref{prvee}\\

$\textbf{Ex}$ & Axiom scheme $(\alpha\wedge\beta\rightarrow\bot_{\alpha\beta})\rightarrow(\alpha\uparrow\beta)$ & \pageref{Ex}\\

$\bI\textbf{Ex}$ & Logic obtained from $\bI$ by addition of $\textbf{Ex}$ & \pageref{bIEx}\\

$\textbf{ciw}^{\uparrow}$ & Axiom scheme $(\alpha\uparrow\beta)\vee(\alpha\wedge\beta)$ & \pageref{ciwuparrow}\\

$\bI\textbf{ciw}^{\uparrow}$ & Logic obtained from $\bI$ by addition of $\textbf{ciw}^{\uparrow}$ & \pageref{bIciw}\\

$\nbI$ & Logic obtained from $\bI$ by addition of the axiom scheme $\alpha\vee\neg\alpha$ & \pageref{nbI}\\

$\textbf{ci}^{\uparrow}$ & Axiom scheme $\neg(\alpha\uparrow\beta)\rightarrow(\alpha\wedge\beta)$ & \pageref{ciuparrow}\\

$\textbf{cl}^{\uparrow}$ & Axiom scheme $\neg(\alpha\wedge\beta)\rightarrow(\alpha\uparrow\beta)$ & \pageref{cluparrow}\\

$\nbIciw$ & Logic obtained from $\nbI$ by addition of the axiom scheme $\textbf{ciw}^{*}$ & \pageref{nbIciw}\\

$\textbf{ciw}^{*}$ & Axiom scheme $(\alpha\uparrow\neg\alpha)\vee(\alpha\wedge\neg\alpha)$ & \pageref{ciw*}\\

$\nbIci$ & Logic obtained from $\nbI$ by addition of the axiom scheme $\textbf{ci}^{*}$ & \pageref{nbIci}\\

$\textbf{ci}^{*}$ & Axiom scheme $\neg(\alpha\uparrow\neg\alpha)\rightarrow(\alpha\wedge\neg\alpha)$ & \pageref{ci*}\\

$\nbIcl$ & Logic obtained from $\nbI$ by addition of the axiom scheme $\textbf{cl}^{*}$ & \pageref{nbIcl}\\

$\textbf{cl}^{*}$ & Axiom scheme $\neg(\alpha\wedge\neg\alpha)\rightarrow(\alpha\uparrow\neg\alpha)$ & \pageref{cl*}\\

$\textbf{cc}^{*}$ & Axiom scheme $(\alpha\uparrow\neg\alpha)\uparrow\neg(\alpha\uparrow\neg\alpha)$ & \pageref{cc*}\\

$\mathcal{L}^{*}$ & Translation of the $\textbf{LFI}$ logic $\mathcal{L}$ to an incompatibility logic & \pageref{L*}\\

$\textbf{mbC}^{*}$ & Translation of $\textbf{mbC}$ & \pageref{mbC*}\\

$\textbf{mbCciw}^{*}$ & Translation of $\textbf{mbCciw}$ & \pageref{mbCciw*}\\

$\textbf{mbCci}^{*}$ & Translation of $\textbf{mbCci}$ & \pageref{mbCci*}\\

$\textbf{mbCcl}^{*}$ & Translation of $\textbf{mbCcl}$ & \pageref{mbCcl*}\\

\end{longtable}

\newpage

\section*{List of categories and functors}

\begin{longtable}{c | m{11cm} | c}

Symbols & Meaning & Page \\ \hline

$\textbf{MAlg}$ & Category of $\Sigma$-multialgebras equipped with homomorphisms & \pageref{MAlg}\\

$\textbf{MAlg}_{=}$ & Category of $\Sigma$-multialgebras equipped with full homomorphisms & \pageref{MAlg=}\\

$\textbf{MMAlg}$ & Category of $\Sigma$-multialgebras equipped with multihomomorphisms & \pageref{MMAlg}\\

$\textbf{MMAlg}_{=}$ & Category of $\Sigma$-multialgebras equipped with full multihomomorphisms & \pageref{MMAlg=}\\

$\textbf{MSet}$ & Category of sets equipped with multifunctions & \pageref{MSet}\\

$\textbf{MG}(\Sigma)$ & Category of $\Sigma$-multialgebras equipped with classes of equivalence of ground-preserving homomorphisms & \pageref{MGSigma}\\

$\textbf{CABA}$ & Category of complete, atomic Boolean algebras, equipped with continuous homomorphisms & \pageref{CABA}\\

$\textbf{Alg}(\Sigma)$ & Category of $\Sigma$-algebras equipped with homomorphisms & \pageref{AlgSigma}\\

$\mathcal{P}_{\mathcal{B}}$ & Functor from $\textbf{MAlg}(\Sigma)$ to $\textbf{Alg}_{\mathcal{B}}(\Sigma)$ & \pageref{Functor PB}\\

$\textbf{Alg}_{\mathcal{B}}(\Sigma)$ & Category of $(\Sigma, \leq)$-algebras equipped with $(\Sigma, \leq)$-homomorphisms & \pageref{AlgBSigma}\\

$P$ & Endofunctor on $\textbf{MAlg}(\Sigma)$ that is part of a monad &\pageref{P}\\

$\mathbb{A}$ & Functor from $\textbf{Alg}_{\mathcal{B}}(\Sigma)$ to $\textbf{MAlg}(\Sigma)$ & \pageref{A}\\

$\textbf{Set}^{*}$ & Category of non-empty sets equipped with functions & \pageref{Set*}\\

$\textbf{CABo}$ & Category of complete, atomic and bottomless Boolean algebras equipped with continuous, atoms-preserving functions & \pageref{CABo}\\

$\textbf{PMAlg}(\Sigma)$ & Category of partial multialgebras equipped with homomorphisms & \pageref{PMAlgSigma}\\

$\textbf{Alg}_{CABA}(\Sigma)$ & Category of $\Sigma$-algebras, with a Boolean algebra structure, equipped with homomorphisms that are continuous and atoms-preserving & \pageref{AlgCABASigma}\\

$\textbf{MAlg}_{\mathcal{B}}(\Sigma)$ & Category of $(\Sigma, \leq)$-algebras equipped with continuous almost-homomorphisms & \pageref{MAlgBSigma}\\

$\textbf{RSwap}_{C_{n}}$ & Category of swap structures for $C_{n}$, equipped with homomorphisms commuting with $\mathcal{F}_{C_{n}}^{\mathcal{B}}$ & \pageref{RSwapCn}\\

$\textbf{BA}$ & Category of non-degenerate Boolean algebras, equipped with homomorphisms & \pageref{BA}\\

$\mathcal{A}_{n}$ & Functor from $\textbf{BA}$ to $\textbf{RSwap}_{C_{n}}$ & \pageref{An}\\

$\textbf{Boo}_{n}$ & Functor, inverse of $\mathcal{A}_{n}$ & \pageref{tBoon}\\

\end{longtable}

\newpage

\section*{List of signatures}

\begin{longtable}{c | m{11cm} | c}

Symbols & Meaning & Page \\ \hline

$\Sigma$ & Arbitrary signature & \pageref{signature}\\

$\Sigma_{\textbf{Lat}}$ & Signature with $\vee$ and $\wedge$ as symbols & \pageref{SigmaLat}\\

$\Sigma_{\textbf{Lat}}^{0}$ & Signature with $0, \vee$ and $\wedge$ as symbols & \pageref{SigmaLat0}\\

$\Sigma_{\textbf{Lat}}^{1}$ & Signature with $1, \vee$ and $\wedge$ as symbols & \pageref{SigmaLat1}\\

$\Sigma_{\textbf{Lat}}^{0,1}$ & Signature with $0, 1, \vee$ and $\wedge$ as symbols & \pageref{SigmaLat01}\\

$\Sigma_{\textbf{Imp}}$ & Signature with $1, \vee, \wedge$ and $\rightarrow$ as symbols & \pageref{SigmaImp}\\

$\Sigma_{\textbf{Hey}}$ & Signature with $0, 1, \vee, \wedge$ and $\rightarrow$ as symbols & \pageref{SigmaHey}\\

$\Sigma_{\textbf{Boo}}$ & Signature with $0, 1, \neg, \vee, \wedge$ and $\rightarrow$ as symbols & \pageref{SigmaBoo}\\

$\Sigma_{s}$ & Signature with an unary $s$ as only symbol & \pageref{Sigmas}\\

$\Sigma_{\textbf{CPL}}$ & Signature with $\bot, \top, \sim, \vee, \wedge$ and $\rightarrow$ as connectives & \pageref{SigmaCPL}\\

$\Sigma_{\textbf{LFI}}$ & Signature with $\neg, \circ, \vee, \wedge$ and $\rightarrow$ as connectives & \pageref{SigmaLFI}\\

$\Sigma^{\top}$ & Signature obtained from $\Sigma$ by addition of the $0$-ary connective $\top$ & \pageref{Sigmatop}\\

$\Sigma_{\textbf{C}}$ & Signature with $\neg, \vee, \wedge$ and $\rightarrow$ as connectives & \pageref{SigmaC}\\

$\Sigma_{\textbf{bI}}$ & Signature with $\vee, \wedge, \rightarrow$ and $\uparrow$ as connectives & \pageref{SigmabI}\\

$\Sigma_{\bI}^{\textbf{CPL}}$ & Signature with $\bot, \top, \sim, \vee, \wedge, \rightarrow$ and $\uparrow$ as connectives &\pageref{SigmabICPL}\\

$\Sigma_{\nbI}$ & Signature with $\neg, \vee, \wedge, \rightarrow$ and $\uparrow$ as connectives & \pageref{SigmanbI}\\

$\Sigma_{\nbI}^{\textbf{CPL}}$ & Signature with $\bot, \top, \sim, \neg, \vee, \wedge, \rightarrow$ and $\uparrow$ as connectives & \pageref{SigmanbICPL}\\

$\Sigma_{\textbf{LFI}}^{\textbf{CPL}}$ & Signature with $\bot, \top, \circ, \sim, \neg, \vee, \wedge$ and $\rightarrow$ as connectives & \pageref{SigmaLFICPL}\\

\end{longtable}

\end{refsegment}

\printindex

\addcontentsline{toc}{chapter}{Bibliography}
\nocite{*}
\printbibliography

@article{Avron,
	title = {Non-deterministic semantics for logics with a consistency operator},
	volume = {45},
	doi = {10.1016/j.ijar.2006.06.011},
	issue= {2},
	journal = {International Journal of Approximate Reasoning},
	author = {Avron, Arnon},
	year = {2007},
	pages = {271--287}
}

@book{Men-IntLog,
	title = {Introduction to mathematical logic},
	doi={10.1007/978-1-4615-7288-6},
	isbn = {978-1-46-157290-9},
	publisher = {Springer US}, 
	author = {Mendelsohn, Eliot},
	year = {1987}
}

@book{ParLog,
	series = {Logic, epistemology, and the unity of science},
	title = {Paraconsistent Logic: Consistency, Contradiction and Negation},
	isbn = {978-3-31-933203-1},
	doi={10.1007/978-3-319-33205-5},
	volume = {40},
	publisher = {Springer International Publishing},
	author = {Carnielli, Walter A. and Coniglio, Marcelo E.},
	year = {2016}
}

@article{Piochi2,
	title = {Logical matrices and non-structural consequence operators},
	volume = {42},
	doi = {10.1007/BF01418757},
	journal = {Studia Logica},
	author = {Piochi, Brunetto},
	year = {1983},
	pages = {33--42}
}

@article{Rescher,
	title = {Quasi-truth-functional systems of propositional logic},
	volume = {27},
	doi = {10.2307/2963674},
	issue = {1},
	journal = {The Journal of Symbolic Logic},
	author = {Rescher, Nicholas},
	year = {1962},
	pages = {1--10}
}

@book{Burris,
	series = {Graduate Texts in Mathematics},
	title = {A Course in Universal Algebra},
	volume = {78},
	isbn = {978-1-46-138132-7},
	publisher = {Springer-Verlag New York},
	author = {Burris, Stanley and Sankappanavar, Hanamantagouda P. },
	year = {1981},
}

@book{Woj,
	title = {Lectures on propositional calculi},
	isbn = {978-8-30-401775-7},
	publisher = {Publishing House of the Polish Academy of Sciences},
	author = {W{\'o}jcicki, Ryszard},
	year = {1984},
}

@incollection{AvronLev,
	title = {Canonical Propositional Gentzen-Type Systems},
	volume = {2083},
	isbn = {978-3-54-042254-9},
	booktitle = {Automated Reasoning},
	publisher = {Springer, Berlin, Heidelberg},
	author = {Avron, Arnon and Lev, Iddo},
	editor = {Gor\'e, Rajeev and Leitsch, Alexander and Nipkow, Tobias},
	year = {2001},
	doi = {10.1007/3-540-45744-5_45},
	pages = {529--544}
}

@article{Costa2,
	title = {A semantical analysis of the calculi $\textbf{C}_{n}$},
	volume = {18},
	doi = {10.1305/ndjfl/1093888132},
	issue = {4},
	journal = {Notre Dame Journal of Formal Logic},
	author = {da Costa, Newton C. A. and Alves, Elias H.},
	year = {1977},
	pages = {621--630}
}

@phdthesis{Fidel,
	title={Nuevos enfoques en L\'ogica Algebraica},
	author={Fidel, Manuel M.},
	year={2003},
	school={Universidad Nacional del Sur},
	address={Bahia Blanca, Argentina}
}

@unpublished{Fidel2,
	title={Una nueva sem\^antica de tipo algebraico para las l\'ogicas paraconsistentes $\textbf{C}_{n}$ de da Costa},
	author={Fidel, Manuel M.},
	year={2009},
	howpublished={Presented at the Reuni\'on Anual de la Uni\'on Matem\'atica Argentina},
	address={Mar del Plata, Argentina}
}

@article{Lewin,
	title = {$\textbf{C}_{1}$ is not algebraizable},
	volume = {32},
	doi = {10.1305/ndjfl/1093635932},
	issue = {4},
	journal = {Notre Dame Journal of Formal Logic},
	author = {Lewin, Renato A. and Mikenberg, Irene F. and Schwarze, Maria G.},
	year = {1991},
	pages = {609--611}
}

@book{BlokPigozzi,
	title = {Algebraizable Logics},
	isbn = {978-1-47-040816-9},
	series={Memoirs of the American Mathematical Society},
	author = {Blok, Willem J. and Pigozzi, Don},
	year = {1989},
	volume={77},
	issue={396},
}

@thesis{Costa3,
	title={Sistemas formais inconsistentes},
	author={da Costa, Newton C. A.},
	year={1963},
	type={Habilitation thesis},
	college={Universidade Federal do Paran{\'a}},
	address={Curitiba, Brazil},
	howpublished={Republished by Editora Universidade Federal do Paran{\'a}, 1993}
}

@article{Woj2,
	title={Some remarks on the consequence operation in sentential logics},
	doi={10.4064/fm-68-3-269-279},
	author={W{\'o}jciki, Ryszard},
	year={1970},
	journal={Fundamenta Mathematicae},
	volume={68},
	pages={269--279}
}

@article{Piochi,
	title={Matrici adequate per calcoli generali predicativi},
	author={Piochi, Brunetto},
	journal={Bolletino della Unione Matematica Italiana},
	volume={15},
	year={1978},
	pages={66--76}
}

@article{Piochi3,
	title={Nota su matrici adeguate per calcoli generali predicativi},
	author={Piochi, Brunetto},
	year={1980},
	journal={Bolletino della Unione Matematica Italiana},
	volume={1},
	pages={271--273}
}

@inproceedings{Marty,
	title={Sur une g\'en\'eralisation de la notion de groupe},
	author={Marty, Fr\'ed\'eric},
	year={1934},
	booktitle={Comptes rendus du huiti\`eme congr\`es des math\'ematiciens scandinaves},
	address={Stockholm, Sweden},
	pages={45--49}
}

@article{Costa,
	title={Une s\'emantique pour le calcul $\textbf{C}_{1}$},
	author={da Costa, Newton C. A. and Alves, Elias H.},
	year={1976},
	journal={Comptes Rendus de l'Acad\'emie de Sciences de Paris (A-B)},
	volume={283},
	pages={729--731}
}

@inproceedings{Itala,
	title={Translations between logics},
	author={da Silva, Jairo J. and D'Ottaviano, Itala M. L. and Sette, Antonio M. A.},
	year={1999},
	booktitle={Models, algebras and proofs; Selected papers of the X Latin American symposium on Mathematical Logic held in Bogot{\'a}},
	isbn={978-0-82-471970-8},
	series={Lecture Notes in Pure and Applied Mathematics},
	volume={203},
	editors={Caicedo, Xavier and Montenegro, Carlos H.},
	pages={435--448},
	publisher={Marcel Dekker, New York},
}

@inbook{Handbook,
	title={Logics of Formal Inconsistency},
	author={Carnielli, Walter A. and Coniglio, Marcelo E. and Marcos, Jo\~ao},
	year={2007},
	booktitle={Handbook of Philosophical Logic},
	edition={Second},
	editors={Gabbay, Dov M. and Guenthner, Franz},
	publisher={Springer, Dordrecht},
	volume={14},
	pages={1--93},
	doi={10.1007/978-1-4020-6324-4_1},
	isbn={978-1-40-206323-7}
}

@inproceedings{Loparic,
	title={The semantics of the systems $\textbf{C}_{n}$ of da Costa},
	author={Lopari\'c, Andr\'ea and Alves, Elias H.},
	booktitle={Proceedings of the Third Brazilian Conference on Mathematical Logic},
	year={1980},
	editors={Arruda, Ayda I. and da Costa, Newton C. A. and Sette, Antonio M. A.},
	publisher={Sociedade Brasileira de L\'ogica},
	pages={161--172}
}

@misc{Oosten,
	title={Basic Category Theory},
	author={van Oosten, Jaap},
	year={2002},
	publisher={Ultrecht University},
	address={Ultrecht, Netherlands},
	note={Lecture notes}
}

@book{CWM,
	title={Categories for the Working Mathematician},
	author={Mac Lane, Saunders},
	year={1978},
	series={Graduate Texts in Mathematics},
	volume={5},
	doi={10.1007/978-1-4757-4721-8}, 
	isbn={978-0-38-798403-2},
	publisher={Springer-Verlag, New York, NY}
}

@article{Stone,
	title={The Theory of Representations for Boolean Algebras},
	author={Stone, Marshall H.},
	year={1936},
	journal={Transactions of the American Mathematical Society},
	volume={40},
	issue={1},
	pages={37--111},
	publisher={American Mathematical Society},
	doi={10.2307/1989664}
}

@article{CFG,
	title={Non-de\-ter\-mi\-nis\-tic algebraization of logics by swap structures},
	journal={Logic Journal of the IGPL},
	author={Coniglio, Marcelo E. and Figallo-Orellano, Aldo and Golzio, Ana C.},
	year={2020},
	volume={28},
	issue={5},
	pages={1021--1059}, 
	doi={10.1093/jigpal/jzy072}
}

@article{ConiglioSernadas,
	title={A Graph-theoretic Account of Logics},
	author={Coniglio, Marcelo E. and Sernadas, Amilcar and Sernadas, Cristina and Rasga, Jo{\~a}o},
	year={2009},
	journal={Journal of Logic and Computation},
	volume={19},
	issue={6},
	pages={1281--1320},
	doi={10.1093/logcom/exp023}
}

@article{CuponaMadarasz,
	author={{\v C}upona, {\'G}orgi and  Madar{\'a}sz, Roz{\'a}lia Sz.},
	title={Free Poly-Algebras},
	year={1993},
	series={Mathematics Series of the University of Novi Sad},
	volume={23},
	issue={2},
	pages={245--261}
}

@inproceedings{CM,
	title={A taxonomy of {C}-systems},
	author={Carnielli, Walter A. and Marcos, Jo{\~a}o},
	year={2002},
	booktitle={Paraconsistency: The Logical Way to the Inconsistent},
	editors={Carnielli, Walter A. and Coniglio, Marcelo E. and D'Ottaviano, Itala M. L.},
	volume={228},
	pages={1--94},
	series={Lecture Notes in Pure and Applied Mathematics},
	isbn={978-1-13-846690-6},
	publisher={Marcel Dekker}
}

@article{Mortensen80,
	title={Every quotient algebra for $\textbf{C}_1$ is trivial},
	author={Mortensen, Chris},
	journal={Notre Dame\\ Journal of Formal Logic},
	volume={21},
	issue={4},
	year={1980},
	pages={694--700},
	doi={10.1305/ndjfl/1093883254}
}

@inbook{Mortensen89,
	title={Paraconsistency and $\textbf{C}_1$},
	author={Mortensen, Chris},
	year={1989},
	booktitle={Paraconsistent Logic: Essays on the Inconsistent},
	editors={Priest, Graham and Routley, Richard and Norman, Jean},
	pages={289--305},
	publisher={Philosophia Verlag Gmbh},
	isbn={978-3-88-405058-3}
}

@book{Woj3,
	title={Theory of Logical Calculi},
	subtitle={Basic Theory of Consequence Operations},
	author={W{\'o}jcicki, Ryszard},
	volume={199},
	year={1988},
	doi={10.1007/978-94-015-6942-2},
	publisher={Springer Netherlands},
	isbn={978-9-02-772785-5}
}

@article{CCP,
	title={Finite non-de\-ter\-mi\-nis\-tic semantics for some modal systems},
	author={Marcelo E. Coniglio and Luis F. del Cerro and Newton M. Peron},
	journal={Journal of Applied Non-Classical Logics},
	volume={25},
	issue={1},
	year={2015},
	pages={20--45},
	doi={10.1080/11663081.2015.1011543},
}

@article{CCPErr,
	title={Errata and addenda to `{F}inite non-de\-ter\-mi\-nis\-tic semantics for some modal systems'},
	author={Coniglio, Marcelo E. and del Cerro, Luis F. and Peron, Newton M.},
	journal={Journal of Applied Non-Classical Logics},
	volume={26},
	issue={4},
	year={2016},
	pages={336--345},
	doi={10.1080/11663081.2017.1300436},
}

@article{OmoriSkurt,
	title={More modal semantics without possible worlds},
	author={Omori, Hitoshi and Skurt, Daniel},
	journal={IfCoLog Journal of Logics and their Applications},
	volume={3},
	issue={5},
	year={2016},
	pages={815--846},
}

@article{Kearns,
	title={Modal semantics without possible worlds},
	author={Kearns, John T.},
	journal={The Journal of Symbolic Logic},
	volume={46},
	issue={1},
	year={1981},
	pages={77--86},
}

@article{Fidel3,
	title={The decidability of the calculi ${C}_n$},
	author={Fidel, Manuel M.},
	journal={Reports on Mathematical Logic},
	volume={8},
	year={1977},
	pages={31--40},
}

@inbook{CCCM,
	title={Two’s company: `The humbug of many logical values'},
	author={Caleiro, Carlos and Carnielli, Walter and Coniglio, Marcelo E. and Marcos, Jo{\~a}o},
	year={2005},
	booktitle={Logica Universalis},
	editors={Beziau, Jean-Yves},
	pages={169--189},
	publisher={Birkh\"auser Basel},
	doi={10.1007/3-7643-7304-0_10},
	isbn={978-3-76-437259-0}
}

@article{Dugundji,
	title={Note on a property of matrices for Lewis and Langford's calculi of propositions},
	author={Dugundji, James},
	journal={The Journal of Symbolic Logic},
	volume={5},
	issue={4},
	year={1940},
	pages={150--151},
	doi={10.2307/2268175}
}

@article{Godel,
	title={Zum intuitionistischen Aussagenkalk\"ul},
	author={G\"odel, Kurt},
	journal={Anzeiger der Akademie\-der Wissenschaften in Wien, Mathematisch-naturwissenschaftliche Klasse},
	volume={69\\},
	year={1932},
	pages={65--66}
}

@inproceedings{Avron2,
	title={Non-deterministic matrices and modular semantics of rules},
	author={Avron, Arnon},
	year={2005},
	booktitle={Logica Universalis},
	editors={Beziau, Jean-Yves},
	publisher={Birkh{\"a}user Basel},
	pages = {149--167},
	doi={10.1007/3-7643-7304-0_9}
}

@inproceedings{Avron3,
	title={Non-deterministic semantics for paraconsistent {C}-systems},
	author={Avron, Arnon},
	year={2005},
	booktitle={Symbolic and Quantitative Approaches to Reasoning with Uncertainty\\ (ECSQARU 2005)},
	editors={Godo, Llu{\'i}s},
	volume={3571},
	pages={625--637},
	series={Lecture Notes in Computer Science},
	publisher={Springer, Berlin, Heidelberg},
	isbn={978-3-54-027326-4},
	doi={10.1007/\\11518655_53}
}

@ARTICLE{AvronLev2,
	author={Avron, Arnon and Lev, Iddo},
	year={2005},
	title={Non-deterministic multi-valued structures},
	journal={Journal of Logic and Computation},
	volume={15},
	issue={3},
	pages={241--261},
	doi={10.1093/logcom/exi001}
}

@ARTICLE{Pawlowski,
	author={Pawlowski, Pawel},
	year={2020},
	title={Tree-Like Proof Systems for Finitely-Many Valued Non-de\-ter\-mi\-nis\-tic Consequence Relations},
	journal={Logica Universalis},
	volume={14},
	pages={407--420},
	note={10.1007/\\s11787-020-00263-0},
}

@ARTICLE{PawlowskUrbaniak,
	author={Pawlowski, Pawel and Urbaniak, Rafal},
	year={2018},
	title={Many-valued logic of informal provability: a non-de\-ter\-mi\-nis\-tic strategy},
	journal={The Review of Symbolic Logic},
	volume={11},
	issue={2},
	pages={207--223},
	doi={10.1017/S1755020317000363}
}

@INPROCEEDINGS{CrawEther,
	author = {Crawford, James M. and Etherington, David W.},
	title = {A non-deterministic semantics for tractable inference},
	booktitle = {Proceedings of the Fifteenth National Conference on Artificial Intelligence and Tenth Innovative Applications of Artificial Intelligence Conference, (AAAI 98, IAAI 98)},
	editor = {Mostow, Jack and Rich, Chuck},
	pages = {286--291},
	year = {1998},
	publisher = {AAAI Press/The MIT Press},
	isbn={02-6251-098-7}
}

@ARTICLE{Lev,
	author =       {Ivlev, Ju. V.},
	title =        {A semantics for modal calculi},
	journal =      {Bulletin of the Section of Logic},
	year =         {1988},
	volume =       {17},
	issue =       {3/4},
	pages =        {114--121},
}

@ARTICLE{Lev2,
	author={Ivlev, Ju. V.},
	title={Generalization of {K}almar's method for quasi-matrix logic},
	journal={Logical Investigations},
	year={2013},
	volume={19},
	pages ={281--307},
	doi={10.21146/2074-1472-2013-19-0-281-307},
}

@ARTICLE{Lev3,
	author =       {Ivlev, Ju. V.},
	title =        {Tablitznoe postrojenie propozicionalnoj modalnoj logiki},
	journal =      {Vestnik Moskovskogo Universiteta},
	series={Filosofia},
	year =         {1973},
	volume =       {6},
}

@BOOK{Lev4,
	author =       {Ivlev, Ju. V.},
	title =        {Sodierzatelnaja semantika modalnoj logiki},
	year =         {1985},
	publisher =       {Moscow},
}

@inproceedings{CF,
	title={A model-theoretic analysis of Fidel-structures for \textbf{mbC}},
	author={Coniglio, Marcelo E. and Figallo-Orellano, Aldo},
	year={2019},
	booktitle={Graham Priest on Dialetheism and Paraconsistency},
	editors={Baskent, Can and Ferguson, Thomas M.},
	series={Outstanding Contributions to Logic},
	publisher={Springer, Cham},
	volume={18},
	pages={189--216},
	isbn={978-3-03-025364-6},
	doi={10.1007/978-3-030-25365-3_10}
}

@ARTICLE{Jaskowski,
	author =       {Ja\'{s}kowski, Stanis{\l}aw},
	title =        {Rachunek zda\'{n} dla system\'{o}w dedukcyjnych sprzecznych},
	journal =      {Studia Societaiis Scientiarum Torunensis},
	year =         {1948},
	volume =       {1},
	pages =        {55--77},
}

@ARTICLE{Jaskowski2,
	author =       {Ja\'{s}kowski, Stanis{\l}aw},
	title =        {O koniunkcji dyskusyjnej w rachunku zda\'{n} dla system\'{o}w dedukcyjnych sprzecznych},
	journal =      {Studia Societaiis Scientiarum Torunensis},
	year =         {1949},
	volume =       {1},
	pages =        {171--172},
}

@inproceedings{CalMar,
	author={Caleiro, Carlos and Marcelino, S{\'e}rgio},
	title={Analytic Calculi for Monadic {PN}matrices},
	year={2019},
	booktitle={Logic, Language, Information, and Computation},
	editors={Iemhoff, Rosalie and Moortgat, Michael and de Queiroz, Ruy},
	series={Lecture Notes in Computer Science},
	volume={11541},
	isbn={978-3-66-259532-9},
	publisher={Springer, Berlin, Heidelberg},
	doi={10.1007/978-3-662-59533-6_6},
}

@article{Baaz,
	author={Baaz, Matthias and Lahav, Ori and Zamansky, Anna},
	title={A Finite-valued Semantics\\ for Canonical Labelled Calculi},
	year={2013},
	journal={Journal of Automated Reasoning},
	volume={51},
	pages={401--430},
	DOI={10.1007/s10817-013-9273-x},
}

@misc{AbsFreeHyp,
      title={Absolutely Free Hyperalgebras}, 
      author={Coniglio, Marcelo E. and Toledo, Guilherme V.},
      year={2021},
      eprint={2101.03647},
      archivePrefix={arXiv},
      primaryClass={math.LO}
}

@misc{CostaRNmatrix,
      title={A simple decision procedure for da Costa's Cn logics by Restricted Nmatrix semantics}, 
      author={Coniglio, Marcelo E. and Toledo, Guilherme V.},
      year={2020},
      eprint={2011.10151},
      archivePrefix={arXiv},
      primaryClass={math.LO}
}

@article{Peregrin1,
	title={Brandom's Incompatibility Semantics},
	journal={Philosophical Topics},
	volume={36},
	issue={2},	
	author={Peregrin, Jaroslav},
	year={2008},
	pages={99--121},
	doi={10.5840/philtopics200836221}
}

@inbook{Peregrin2,
	title={Logic as based in incompatibility},
	author={Peregrin, Jaroslav},
	year={2011},
	booktitle={The Logica Yearbook 2010},
	publisher={College Publications},
	editors={Peli{\v s}, Michal and Pun{\v}coch{\'a}{\v r} V{\'i}t},
	isbn={978-1-84-890038-7}
}

@book{Brandom,
	title={Between Saying and Doing},
	subtitle={Towards an analytic pragmatism},
	author={Brandom, Robert B.},
	publisher={Oxford University Press},
	isbn={978-0-19-954287-1},
	year={2008},
}

@article{Rivieccio,
	title={Implicative twist-structures},
	author={Rivieccio, Umberto},
	journal={Algebra universalis},
	volume={71\\},
	year={2014},
	pages={155--186},
	doi={10.1007/s00012-014-0272-5}
}

@article{OnoRivieccio,
	title={Modal twist-structures over residuated lattices},
	author={Ono, Hiroakira and Rivieccio, Umberto},
	journal={Logic Journal of the IGPL},
	volume={22},
	issue={3},
	year={2014},
	pages={440--457},
	doi={10.1093/jigpal/jzt043}
}

@article{Cignoli,
	title={The class of Kleene algebras satisfying an interpolation property and Nelson algebras},
	author={Cignoli, Roberto},
	journal={Algebra universalis},
	volume={23},
	year={1986},
	pages={262--292},
	doi={10.1007/BF01230621}
}

@inproceedings{OS:20,
	author={Omori, Hitoshi and Skurt, Daniel},
	title={A Semantics for a Failed Axiomatization of $K$},
	editor={Olivietti, N. and Verbrugge, R. and Negri, S. and Sandu, G.},
	booktitle={Advances in Modal Logic},
	volume={13},
	pages={481--501},
	publisher={College Publications},
	year={2020}
}

@article{AK:05,
	title={Multi-valued Calculi for Logics Based on Non-determinism},
	author={Avron, Arnon and Konikowska, Beata},
	journal={Logic Journal of the IGPL},
	volume={13},
	issue={4},
	year={2005},
	pages={365--387},
	doi={10.1093/jigpal/jzi030}
}

@book{Smullyan,
	title = {First-Order Logic},
	publisher = {Dover Publications},
	address = {Mineola, N.Y. USA},
	author = {Smullyan, Raymond M.},
	year = {1995},
	note = {Corrected republication of the Springer-Verlag, New York, 1968 edition}
}

@Unpublished{con:far:per:21,
	title={Tableaux systems for some {I}vlev-like (quantified) modal logics},
	author={Coniglio, Marcelo E. and del Cerro, Luis F. and Peron, Newton M.},
	note={To appear},
	year={2021},
}

@article{ItalaCastro,
	author={D'Ottaviano, Itala M. L. and de Castro, Milton A.},
	year={2006},
	title={Analytical Tableaux for da Costa's Hierarchy of Paraconsistent Logics},
	journal={Electroniv Notes in Theoretical Computer Science},
	volume={143},
	pages={27--44},
	doi={10.1016/j.entcs.2005.06.033}
}

@article{WeaklyFreeMultialgebras,
	author={Coniglio, Marcelo E. and Toledo, Guilherme V.},
	title={Weakly Free Multialgebras},
	year={2022},
	volume={51},
	number={1},
	pages={109--141},
	journal={Bulletin of the Section of Logic},
	doi={10.18778/0138-0680.2021.19}
}

@phdthesis{Golzio,
	title={Non-deterministic matrices: theory and applications to algebraic semantics},
	author={Golzio, Ana C.},
	year={2017},
	school={University of Campinas, Brazil},
	url={http://repositorio.unicamp.br/jspui/handle/REPOSIP/322436}
}

@masthersthesis{Alves-Thesis,
	author={Alves, Elias H.},
	year={1976},
	title={L\'ogica e inconsist\^encia: um estudo dos c\'alculos $C_{n}$, $1\leq n<\omega$},
	school={University of S\~ao Paulo, Brazil},
	url={https://repositorio.usp.br/item/000722069}
}

@incollection{Gratz:21,
	title={Analytic Tableaux for Non-deterministic Semantics},
	author={Gr\"atz, Lukas},
	booktitle={Automated Reasoning with Analytic Tableaux and Related Methods},
	publisher={Springer International Publishing},
	editor={Das, A. and Negri, S.},
	series={Lecture Notes in Computer Science},
	year={2021},
	pages={38--55}
}

@article{Avron:19,
	title={Paraconsistency and the need for infinite semantics},
	author={Avron, Arnon},
	journal={Soft Computing},
	volume={23},
	number={7},
	year={2019},
	pages={2167--2175},
	doi={https://doi.org/10.1007/s00500-018-3272-0}
}

@inproceedings{Avron:16,
	author={Avron, Arnon},
	title={{RM} and its nice properties},
	year={2016},
	booktitle={J. Michael Dunn on Information Based Logics},
	editor={Bimb\'{o}, K.},
	series={Outstanding Contributions to Logic},
	volume={8},
	pages={15--43},
	publisher={Sprin\-ger},
	DOI={https://doi.org/10.1007/978-3-319-29300-4_2},
}

@book{Avron:Arieli:Zamansky:18,
	title={Theory of Effective Propositional Paraconsistent Logics},
	volume={75},
	series={Studies in Logic (Mathematical Logic and Foundations)},
	author={Avron, Arnon and Arieli, Ofer and Zamansky, Anna},
	publisher={College Publications},
	year={2018}
}

@article{OConnor,
	title={Incompatible Properties},
	author={O'Connor, Daniel J.},
	journal={Analysis},
	volume={15},
	number={5},
	year={1955},
	pages={109--117},
	url={http://www.jstor.org/stable/3326361},
	doi={10.2307/3326361}
}

@book{HeytingAlgebras,
	booktitle={Heyting Algebras},
	author={Esakia, Leo},
	subtitle={Duality Theory},
	editor={Bezhanishvili, Guram and Holliday, Wesley H.},
	publisher={Springer},
	series={Trends in Logic},
	volume={50},
	year={2019},
	isbn={978-3-03-012095-5},
	doi={https://doi.org/10.1007/978-3-030-12096-2},
	translator={Anton Evseev}
}

@misc{RestrictedSwap,
      title={Restricted swap structures for da Costa's $C_{n}$ and their category}, 
      author={Coniglio, Marcelo E. and Toledo, Guilherme V.},
      year={2021},
      eprint={2112.13281},
      archivePrefix={arXiv},
      primaryClass={math.LO}
}

@article{TwoDecisionProcedures,
	title={Two Decision Procedures for\\ da Costa’s $C_{n}$ Logics Based on Restricted Nmatrix Semantics},
	author={Coniglio, Marcelo E. and Toledo, Guilherme V.},
	volume={110},
	issue={3},
	year={2022},
	pages={601--642},
	journal={Studia Logica},
	doi={https://doi.org/10.1007/s11225-021-09972-z}
}

@misc{Frominconsistency,
      title={From inconsistency to incompatibility}, 
      author={Coniglio, Marcelo E. and Toledo, Guilherme V.},
      year={2022},
      eprint={2202.10540},
      archivePrefix={arXiv},
      primaryClass={math.LO}
}

@inbook{ProbabilityLogic,
	title={Probability, Logic and Probability Logic},
	booktitle={The Blackwell Companion to Logic},
	author={H{\'a}jek, Alan},
	editor={Goble, Lou},
	publisher={Blackwell},
	pages={362--384},
	year={2001}
}

@article{Suszko,
	title={The Fregean Axiom and Polish mathematical logic in the $1920^{s}$},
	author={Suszko, Roman},
	journal={Studia Logica},
	volume={36},
	pages={377--380},
	year={1977},
	doi={https://doi.org/10.1007/BF02120672}
}

@inproceedings{PTS-defined,
	title={Possible-translations semantics},
	author={Marcos, Jo{\~a}o},
	editors={Carnielli, WalterA. and Dion{\'i}sio, Francisco M. and Mateus, Paulo C.},
	booktitle={Proceedings of the Workshop on Combination of Logics: Theory and applications},
	pages={119--128},
	year={2004},
	url={http://wslc.math.ist.utl.pt/ftp/pub/MarcosJ/04-M-pts.pdf}
}

@inproceedings{PTS-first,
	title={Many-valued logics and plausible reasoning},
	author={Carnielli, Walter A.},
	booktitle={Proceedings of the Twentieth International Symposium on Multiple-Valued Logic},
	editor={Epstein, George},
	pages={328--335},
	publisher={The IEEE Computer Society Press},
	year={1990}
}

@article{Raftery,
	title={Order algebraizable logics},
	author={Raftery, James G.},
	journal={Annals of Pure and Applied Logic},
	volume={164},
	issue={3},
	year={2013},
	pages={251--283},
	doi={https://doi.org/10.1016/j.apal.2012.10.013}
}

@article{ConservativeTranslationsReference,
	title={Conservative translations},
	author={Feitosa, H{\'e}rcules A. and D'Ottaviano, Itala M. Loffredo},
	journal={Annals of Pure and Applied Logic},
	volume={108},
	pages={205--227},
	doi={https://doi.org/10.1016/S0168-0072(00)00046-4},
	year={2001},
	issue={1--3}
}

@inproceedings{FidelNelson,
	title={An algebraic study of logic with constructible negation},
	author={Fidel, Manuel M.},
	booktitle={Proceedings of the third Brazilian conference on mathematical logic},
	editor={Arruda, Ayda I. and da Costa, Newton C. A. and Sette, Antonio M.},
	year={1980},
	pages={119--129}
}

\end{document}